\documentclass[leqno,sectrefs]{svmono}
\usepackage{type1cm}   
\usepackage[varg,cmintegrals,cmbraces]{newtxmath}
\usepackage[p,osf]{newtxtext}
\usepackage{multind}
\makeindex{word}
\makeindex{notat}
\usepackage[all,2cell,emtex,tips]{xy}
\SelectTips{cm}{10 scaled 700}
\usepackage{amsmath,amssymb,amsfonts,bbm,xspace,paralist}
\usepackage{mathrsfs}
\usepackage[counterclockwise]{rotating}
\usepackage{scrextend}
\usepackage{xcolor}
\PassOptionsToPackage{dvipsnames}{xcolor}
    \RequirePackage{xcolor}
\definecolor{Maroon}{cmyk}{0, 0.87, 0.68, 0.32}
\definecolor{RoyalBlue}{cmyk}{1, 0.50, 0, 0}
\usepackage[pdftex,breaklinks,colorlinks,
linkcolor=Maroon,
citecolor=RoyalBlue,
urlcolor=black]{hyperref}  
\newcommand{\HSpit}{H\ZZ} %Spitzweck sprectrum
\renewcommand{\to}{{\hskip -2.5pt\xymatrixcolsep{1pc}\xymatrix{\ar[r]&}\hskip -2.5pt}}
\renewcommand{\rightarrow}{{\hskip -2.5pt\xymatrixcolsep{1pc}\xymatrix{\ar[r]&}\hskip -2.5pt}}
\renewcommand{\leftarrow}{{\hskip -2.5pt\xymatrixcolsep{1pc}\xymatrix{&\ar[l]}\hskip -2.5pt}}
\renewcommand{\longrightarrow}{{\hskip -2.5pt\xymatrixcolsep{1.3pc}\xymatrix{\ar[r]&}\hskip -2.5pt}}

\renewcommand{\longmapsto}{{\hskip -2.5pt\xymatrixcolsep{1.3pc}\xymatrix{\ar@{|->}[r]&}\hskip -2.5pt}}
\renewcommand{\mapsto}{{\hskip -2.5pt\xymatrixcolsep{.9pc}\xymatrix{\ar@{|->}[r]&}\hskip -2.5pt}}
\renewcommand{\hookrightarrow}{{\hskip -1.5pt\raise 1.5pt\vbox{\xymatrixcolsep{.9pc}\xymatrix{\ar@{^{(}->}[r]&}}\hskip -3.5pt}}
\renewcommand{\xleftarrow}[1]{{\hskip -2.5pt\xymatrixcolsep{1.7pc}\xymatrix{&\ar[l]_{#1}}\hskip -2.5pt}}
\renewcommand{\xrightarrow}[1]{{\hskip -2.5pt\xymatrixcolsep{1.7pc}\xymatrix{\ar[r]^{#1}&}\hskip -2.5pt}}
\renewcommand{\rightleftarrows}{{
\hskip -2.5pt\xymatrixcolsep{1pc}
\xymatrix{ \ar@< 2.1pt>[r]&\ar@<2.1pt>[l] }\hskip -2.5pt}}
\renewcommand{\leftrightarrows}{\rightleftarrows}
\def\limproj{\mathop{\oalign{\rm lim\cr
\hidewidth$\leftarrow$\hidewidth\cr}}}%
\def\limind{\mathop{\oalign{\rm lim\cr
\hidewidth$\rightarrow$\hidewidth\cr}}}%
\renewcommand{\varinjlim}{\limind}%
\renewcommand{\varprojlim}{\limproj}%
\input cyracc.def
\DeclareFontFamily{U}{russian}{}
\DeclareFontShape{U}{russian}{m}{n}
        { <5><6> wncyr5
        <7><8><9> wncyr7
        <10><10.95><12><14.4><17.28><20.74><24.88> wncyr10 }{}
\DeclareSymbolFont{Russian}{U}{russian}{m}{n}
\DeclareSymbolFontAlphabet{\mathcyr}{Russian}
\makeatletter
\let\@math@cyr\mathcyr
\renewcommand{\mathcyr}[1]{\@math@cyr{\cyracc #1}}
\makeatother
\renewcommand{\mathcal}{\mathscr}
\newcommand{\cat}{{\mathscr C\mspace{-2.mu}at}}
\newcommand{\set}{{\mathscr S\mspace{-2.mu}et}}
\newcommand{\ab} {{\mathscr A\mspace{-2.mu}b}}
\newcommand{\tri}{{\mathscr T\mspace{-2.mu}ri}}

\newcommand{\catm}{{\mathscr C\mspace{-2.mu}at^{\otimes}}}
\newcommand{\abm} {{\mathscr A\mspace{-2.mu}b^{\otimes}}}
\newcommand{\trim}{{\mathscr T\mspace{-2.mu}ri^{\otimes}}}

\newcommand{\coin}{{\hbox{
\begin{picture}(6,6)(2,0)
\put(0,0){\line(1,0){6}}
\put(6,6){\line(0,-1){6}}
\end{picture}}}
}
\newcommand{\cocoin}{{\hbox{
\begin{picture}(6,6)(2,0)
\put(0,6){\line(1,0){6}}
\put(0,6){\line(0,-1){6}}
\end{picture}}}
}

\newcommand{\Calg}{\mathit{Comm}}
\newcommand{\Alg}{\mathit{Mon}}
\newcommand{\Mod}[1] {#1\mbox{-}\operatorname{mod}}

\newcommand{\type}{\mathscr C}

\newcommand{\site}{\mathscr S}
\newcommand{\sch}{\mathscr S}

\newcommand{\schsep}{\mathscr S^{sep}}
\newcommand{\schouv}{\mathscr S^{open}}
\newcommand{\schprop}{\mathscr S^{prop}}

\newcommand{\Pmor}{\mathscr P}
\newcommand{\Pmorx}[1]{\Pmor/{#1}}
\newcommand{\sft}{\mathscr{S}^{ft}}
\newcommand{\sm}{\mathit{Sm}}
\newcommand{\ssft}{\mathscr S^{ft}}
\newcommand{\ssm}{\mathscr S\mspace{-2.mu}m}

\newcommand{\Pmorc}{\Pmor^{cor}_{\Rc}}
\newcommand{\sftc}{\mathscr S^{ft,cor}_{\Rc}}
\newcommand{\smc}{\mathscr S\mspace{-2.mu}m^{cor}_{\Rc}}
\newcommand{\Pmorcx}[1]{\Pmor^{cor}_{\Rc,#1}}
\newcommand{\sftcx}[1]{\mathscr S^{ft,cor}_{\Rc,#1}}
\newcommand{\smcx}[1]{\mathscr S\mspace{-2.mu}m^{cor}_{\Rc,#1}}

\newcommand{\Hpt}{\mathscr H_\bullet}
\newcommand{\SH}{\operatorname{SH}}
\newcommand{\eff}{\textit{eff}}
\newcommand{\DMte}{\Der_{\AA^1,\Rc}^{\eff}}
\newcommand{\DMt}{\Der_{\AA^1,\Rc}}

\newcommand{\DMtex}[1]{\Der_{\AA^1}^{\eff}(#1,\Rc)}
\newcommand{\DMtx}[1]{\Der_{\AA^1}(#1,\Rc)}

\newcommand{\DMtgmx}[1]{\Der_{\AA^1,c}(#1,\Rc)}

\newcommand{\uDMte}{\underline{\Der\!}\,_{\AA^1,\Rc}^{\eff}}
\newcommand{\uDMt}{\underline{\Der\!}\,_{\AA^1,\Rc}}
\newcommand{\uDMtx}[1]{\underline{\Der\!}\,_{\AA^1}(#1,\Rc)}

\newcommand{\DMV}{\operatorname{DM}}
\newcommand{\DM}{\DMV_{\Rc}}
\newcommand{\DMc}{\DMV_{c,\Rc}}
\newcommand{\DMgm}{\DMV_{gm,\Rc}}
\newcommand{\DMe}{\DM^{\eff}}
\newcommand{\DMce}{\DMc^{\eff}}
\newcommand{\DMgme}{\DMgm^{\eff}}

\newcommand{\DMx}[1]{\DMV(#1,\Rc)}
\newcommand{\DMcx}[1]{\DMV_{c}(#1,\Rc)}
\newcommand{\DMgmx}[1]{\DMV_{gm}(#1,\Rc)}
\newcommand{\DMex}[1]{\DMV^{\eff}(#1,\Rc)}
\newcommand{\DMcex}[1]{\DMV_{c}^{\eff}(#1,\Rc)}
\newcommand{\DMgmex}[1]{\DMV_{gm}^{\eff}(#1,\Rc)}

\newcommand{\uDMV}{\underline{\DMV}}
\newcommand{\uDM}{\underline{\DMV}_{\Rc}}

\newcommand{\uDMe}{\underline{\DMV}^\eff_\Rc}

\newcommand{\uDMx}[1]{\underline{\DMV}(#1,\Rc)}

\newcommand{\uDMex}[1]{\underline{\DMV}^\eff(#1,\Rc)}

\newcommand{\DMB}{\DMV_\mathcyr B}
\newcommand{\DMBc}{\DMV_{\mathcyr B,c}}

\newcommand{\DMcdh}{\DMV_\cdh}

\newcommand{\DMue}{\Der^{\eff}_{\AA^1\!}} %catgorie $\AA^1$-drive associe
\newcommand{\DMuex}[1]{\Der^{\eff}_{\AA^1\!,#1}} %catgorie $\AA^1$-drive associe avec parametre
\newcommand{\uDMu}{\underline{\Der}_{\AA^1\!}} %catgorie $\AA^1$-drive associe
\DeclareMathOperator{\Spt}{Sp} % spectres de Tate d'une cat. abalienne prem.
\newcommand{\DMu}{\Der_{\AA^1\!}} %catgorie $\AA^1$-drive stable associe
\newcommand{\DMux}[1]{\Der_{\AA^1\!,#1}} %catgorie $\AA^1$-drive stableassocie avec paramtre
\newcommand{\uDMux}[1]{\underline{\Der}_{\AA^1\!,#1}} %catgorie $\AA^1$-drive associe

\newcommand{\psh}[1]{\operatorname{PSh}(#1,\Rc)}%prfaisceaux abliens
\newcommand{\sh}[2]{\operatorname{Sh}_{#1\!}\left(#2,\Rc\right)}%faisceaux abeliens
\newcommand{\ush}[2]{\operatorname{\underline{Sh}}_{#1}\!\left(#2,\Rc\right)}%faisceaux abliens
\newcommand{\shqfh}[1]{\operatorname{Sh}_{\qfh}\!\left(#1,\Rc\right)}%faisceaux qfh
\newcommand{\pshg}[1]{\operatorname{PSh}\left(#1\right)}%prfaisceaux
\newcommand{\shg}[2]{\operatorname{Sh}_{#1}\!\left(#2\right)}%faisceaux abliens
\newcommand{\ptr}[1] {\operatorname{PSh^{\mathit{tr}}}(#1,\Rc)}
\newcommand{\ftr}[1] {\operatorname{Sh^{\mathit{tr}}}(#1,\Rc)}
\newcommand{\ftrx}[2] {\operatorname{Sh^{\mathit{tr}}_{#1}}(#2,\Rc)}

\newcommand{\uftr}[1] {\operatorname{\underline{Sh}^{\mathit{tr}}}(#1,\Rc)}

\newcommand{\uftrx}[2] {\operatorname{\underline{Sh\!}\,^{\mathit{tr}}_{#1}}(#2,\Rc)}

\newcommand{\Hom}{\operatorname{Hom}}
\newcommand{\End}{\operatorname{End}}

\newcommand{\uHom}{\textit{Hom}} %Hom interne
\newcommand{\sHom}{\uHom}

\DeclareMathOperator{\Sym}{Sym}

\DeclareMathOperator{\Comp}{C}
\DeclareMathOperator{\Der}{D}
\newcommand{\Htp} {K}
\DeclareMathOperator{\K}{K}

\newcommand{\Totp}{\operatorname{Tot}^\pi}

\newcommand{\sing}{R_{\AA^1}}
\newcommand{\suslin}{\underline{C}^{*}}

\newcommand{\sus}{\Sigma^\infty}
\newcommand{\lop}{\Omega^\infty}

\newcommand{\spec}[1] {\operatorname{Spec}\left( #1 \right)}
\newcommand{\car}[1] {\operatorname{char}\!\left( #1 \right)}
\DeclareMathOperator{\pic}{Pic}

\newcommand{\mot}[2]{M_{#1}(#2)} %motifs
\newcommand{\umot}[2]{\underline{\phantom{\enskip \,}}\!\!\!\!M_{#1}(#2)} %motif gnralis

\newcommand{\rep}[2]{\Rc_{#1}(#2)}
\newcommand{\urep}[2]{\underline{\Rc\!}\,_{#1}(#2)}
\newcommand{\repx}[3]{\Rc_{#2}^{#1}(#3)}
\newcommand{\urepx}[3]{\underline{\Rc\!}_{#2}^{#1}(#3)}
\newcommand{\uRc}{\underline{\Rc}}
\newcommand{\lrep}[2]{\Rc^{tr}_{#1}(#2)}
\newcommand{\lrepNP}{\Rc^{tr}}
\newcommand{\plrepNP}{{\tilde \Rc}^{tr}}
\newcommand{\ulrep}[2]{\underline{\Rc\!}\,^\tr_{#1}(#2)}

\newcommand{\Mab}[3]{M_{#1}(#2,#3)} %motif pour une catgorie ablienne
\newcommand{\mab}[2]{M_{#1}(#2)} %motif pour une catgorie ablienne
\newcommand{\umab}[2]{\underline{M\!}\,_{#1}(#2)}
\newcommand{\MabNX}[2]{\unit_{#1}}
\newcommand{\uMabNX}[2]{\underline{\unit}_{#1}}
\newcommand{\motNP}{M}
\newcommand{\Th}{\textit Th}
\newcommand{\piso}{\mathfrak p}

\newcommand{\un}{\mathbbm 1} %unit structure monodale
\newcommand{\unit}{\un}
\renewcommand{\emptyset}{\varnothing}
\newcommand{\Hmx}[1]{{H_{\cM,#1}}}
\newcommand{\Hmrx}[1]{{H_{\!\cM,#1}^\QQ}}

\newcommand{\HB}{H_\mathcyr{B}}
\newcommand{\HBx}[1]{{H_{\mathcyr{B},#1}}}

\newcommand{\BGL}{\mathit{KGL}}
\newcommand{\BGLr}{\BGL_\QQ}
\newcommand{\BGLrx}[1]{\BGL_{\QQ,#1}}
\newcommand{\BGLrn}[1]{\BGL^{(#1)}}
\newcommand{\MGL}{\mathit{MGL}}

\newcommand{\base}{\mathcal{S}}

\newcommand{\bc}{($\Pmor$\mbox{-}BC)\xspace}
\newcommand{\pf}{($\Pmor$\mbox{-}PF)\xspace}
\newcommand{\htp}{(Htp)\xspace}
\newcommand{\stab}{(Stab)\xspace}
\newcommand{\sepx}[1]{($#1$-sep)\xspace} % for a topology t
\newcommand{\sep}{(Sep)\xspace}
\newcommand{\ssep}{(sSep)\xspace}
\newcommand{\supp}{(Supp)\xspace}
\newcommand{\suppx}[1]{(Supp$_{#1}$)\xspace}
\newcommand{\loc}{(Loc)\xspace}
\newcommand{\wloc}{(wLoc)\xspace} % restricted localization
\newcommand{\locx}[1]{(Loc$_{#1}$)\xspace} % for a closed immersion i

\newcommand{\adj}{(Adj)\xspace}
\newcommand{\adjx}[1]{(Adj$_{#1}$)\xspace} % for a proper morphism
\newcommand{\BC}{(BC)\xspace}
\newcommand{\BCx}[1]{(BC$_{#1}$)\xspace} % for a proper morphism
\newcommand{\PF}{(PF)\xspace}
\newcommand{\PFx}[1]{(PF$_{#1}$)\xspace} % for a proper morphism
\newcommand{\pur}{(Pur)\xspace} % restricted localization
\newcommand{\wpur}{(wPur)\xspace} % restricted localization
\newcommand{\bg}{(BG)\xspace}
\newcommand{\bgcdh}{(cdh)\xspace}

\DeclareMathOperator{\cycl}{Cycl}
\newcommand{\corrg}[3]{c_{#1}\left(#2,#3\right)_\Rc} %  coefficients
\newcommand{\corr}[3]{c_{#1}\left(#2,#3\right)} % entiers

\newcommand{\acycl}[1]{\langle #1 \rangle} % cycle associ
\newcommand{\otimesf}{\otimes^\flat} %produit tensoriel de cycles de Hilbert
\newcommand{\suppc}[1]{\operatorname{Supp}(#1)}

\newcommand{\tra}[1]{{}^t #1} %transpose d'une correspondance finie
\newcommand{\otr}{\otimes^{tr}} %produit tensoriel avec transferts

\newcommand{\doto} {\bullet\!\!\!\!\longrightarrow}
\newcommand{\ldoto} {\doto}
\newcommand{\xdoto}[1] {\bullet\!\!\!\!\xrightarrow{#1}}

\newcommand{\zar}{{\mathit{Zar}}}
\newcommand{\nis}{{\mathit{Nis}}}
\newcommand{\et}{\mathit{\acute{e}t}}
\newcommand{\cdh}{{\mathit{cdh}}}
\newcommand{\qfh}{{\mathit{qfh}}}
\newcommand{\h}  {{\mathit{h}}}

\newcommand{\derL}{\mathbf{L}}
\newcommand{\derR}{\mathbf{R}}

\newcommand{\To}{\to}
\newcommand{\ho} {\mathrm{Ho}} %categorie homotopique !
\newcommand{\op}[1]{{#1}^{\mathit{op}}}
\newcommand{\Int}{\textstyle{\int}}
\newcommand{\ilim} { \varinjlim }
\newcommand{\plim} { \varprojlim }

\newcommand{\tr}{\mathit{tr}} % (avec) transferts
\DeclareMathOperator{\Tr}{Tr} % morphisme Trace

\newcommand{\hcart}{\mathit{hcart}}
\newcommand{\cart}{\mathit{cart}}

\newcommand{\Rc}{\Lambda}

\newcommand{\A}{\mathscr A}
\newcommand{\B}{\mathscr B}
\newcommand{\C}{\mathscr C}
\renewcommand{\D}{\mathscr D}
\newcommand{\F}{\mathscr F}
\renewcommand{\E}{\mathscr E}
\renewcommand{\H}{\mathscr H}
\renewcommand{\I}{\mathscr I}

\newcommand{\M}{\mathscr M}
\newcommand{\N}{\mathscr N}

\newcommand{\U}{\mathscr U}
\newcommand{\V}{\mathscr V}
\newcommand{\T}{\mathscr T}
\renewcommand{\S}{\mathscr S}

\newcommand{\W}{\mathscr W}
\newcommand{\X}{\mathscr X}
\newcommand{\Y}{\mathscr Y}
\newcommand{\Z}{\mathscr Z}
\newcommand{\cD}{\mathcal D}
\newcommand{\cF}{\mathcal F}
\newcommand{\cG}{\mathcal G}
\newcommand{\cH}{\mathcal H}
\newcommand{\cI}{\mathcal I}

\newcommand{\cM}{\mathcal M}
\newcommand{\cN}{\mathcal N}
\newcommand{\cO}{\mathcal O}
\newcommand{\cS}{\mathcal S}
\newcommand{\cT}{\mathcal T}
\newcommand{\cV}{\mathcal V}
\newcommand{\cX}{\mathcal X}
\newcommand{\NN} {\mathbf N}
\newcommand{\ZZ} {\mathbf Z}
\newcommand{\QQ} {\mathbf Q}

\newcommand{\CC} {\mathbf C}
\newcommand{\fS}{\mathfrak S}
\renewcommand{\AA} {\mathbf A}
\newcommand{\PP} {\mathbf P}
\newcommand{\GG} {\mathbf{G}_m}
\newcommand{\GGx}[1] {\mathbf{G}_{m,#1}}
\newcommand{\fp}{\mathfrak p}
\newcommand{\fq}{\mathfrak q}
\newcommand{\fm}{\mathfrak m}
\newcommand{\HH} { H }
\newcommand{\HHm}[3]{H^{#1,#2}_{\cM}(#3,\Rc)}
\newcommand{\HHme}[3]{H^{#1,#2}_{\cM,\mathit{eff}}(#3,\Rc)}
\newcommand{\uM}{\underline{\phantom{\quad}}\!\!\!\!\!\mathscr M}
\newcommand{\umotNP}{\underline{\phantom{\enskip \,}}\!\!\!\!M} 
\newcommand{\uA}{\underline{\phantom{\quad}}\!\!\!\!\!\mathscr A}

\title{Triangulated categories\\
of mixed motives}

\author{Denis-Charles Cisinski \and Fr\'ed\'eric D\'eglise}
\institute{D.-C. Cisinski \at Fakult\"at f\"ur Mathematik,
Universit\"at Regensburg, 93040 Regensburg, Germany
\and
F. D\'eglise \at Institut de Math\'ematiques de Bourgogne,
MR~5584, Universit\'e de Bourgogne, 9 avenue Alain Savary, BP~47870,
21078 Dijon Cedex, France}

\spnewtheorem{thm}{Theorem}[subsection]{\bfseries}{\itshape}
\spnewtheorem{prop}[thm]{Proposition}{\bfseries}{\itshape}
\spnewtheorem{lm}[thm]{Lemma}{\bfseries}{\itshape}
\spnewtheorem{cor}[thm]{Corollary}{\bfseries}{\itshape}
\spnewtheorem{thmi}{Theorem}{\bfseries}{\itshape}
\spnewtheorem*{cori}{Corollary}{\bfseries}{\itshape}
\spnewtheorem{conj}[thm]{Conjecture}{\bfseries}{\itshape}

\spnewtheorem{rem}[thm]{Remark}{\itshape}{}
\spnewtheorem{notation}[thm]{Notations}{\itshape}{}
\spnewtheorem{ex}[thm]{Example}{\itshape}{}
\spnewtheorem{conv}[thm]{Convention}{\itshape}{}
\spnewtheorem*{remi}{Remark}{\itshape}{}

\spnewtheorem{df}[thm]{Definition}{\bfseries}{}
\spnewtheorem{num}[thm]{\nocaption}{\bfseries}{}
\spnewtheorem{paragr}[thm]{\nocaption}{\bfseries}{}
\spnewtheorem{assumption}[thm]{\nocaption}{\bfseries}{}
\numberwithin{equation}{thm}
\renewcommand{\thesubsubsection}{\thesubsection.\alph{subsubsection}}
\makeindex{word}
\makeindex{notat}
\begin{document}
\frontmatter

\maketitle

\abstract{This book discusses the construction of
triangulated categories of mixed motives
over a noetherian scheme of finite dimension, extending
Voevodsky's definition of motives over a field.
In particular, it is shown that motives with rational coefficients
satisfy the formalism of the six operations of Grothendieck.
This is achieved by studying descent properties
of motives, as well as by comparing different
presentations of these categories, following 
and extending insights and constructions
of Deligne, Beilinson, Bloch, Thomason, Gabber, Levine,
Morel, Voevodsky, Ayoub,
Spitzweck, R\"ondigs, \O stv\ae r and others.
In particular, the relation of motives with $K$-theory is addressed
in full, and we prove the absolute purity theorem with rational
coefficients, using Quillen's localization theorem in algebraic $K$-theory
together with a variation on the Grothendieck-Riemann-Roch theorem.
Using resolution of singularities via alterations
of de~Jong-Gabber, this leads to a version of Grothendieck-Verdier
duality for constructible motivic sheaves with rational coefficients
over rather general base schemes.
We also study versions with integral coefficients, constructed
via sheaves with transfers, for which we obtain partial results.
Finally, we associate to any mixed Weil cohomology a system of categories
of coefficients and well behaved realization functors,
establishing a correspondence between mixed Weil cohomologies
and suitable systems of coefficients. The results of this book
have already served as ground reference in  many subsequent works
on motivic sheaves and their realizations, and pointers to the most
recent developments of the theory are given in the introduction.}

\setcounter{tocdepth}{3}
\tableofcontents

\renewcommand{\thesection}{\Alph{section}}
\addcontentsline{toc}{part}{Introduction}
\extrachap{Introduction}

\section{Historical background}

\subsection{The conjectural theory described by Beilinson}

In a landmark paper, \cite{Bei},
 A.~Beilinson stated a series of conjectures
 which offers a complete renewal of the traditional theory of pure
 motives invented by A.~Grothendieck.
Namely, he proposes to extend the notion of pure motives
 to that of mixed motives with two models in mind: mixed Hodge structures
 defined by P.~Deligne on the one hand \cite{DelH3,DelH2},
 perverse sheaves on the other hand defined in \cite{BBD}.
One of Beilinson's main innovations, motivated by the second model, 
 is to consider a triangulated version 
 of mixed motives in which one could hope to find the more involved theory
 of abelian mixed motives through the concept of t-structures.
This hoped for theory was conjecturally described by Beilinson
 in \cite[5.10]{Bei}
 under the name of \emph{motivic complexes}.

It was modeled (see \emph{loc. cit.}, paragraph A) on the
 theory of \'etale $\ell$-adic sheaves and
 their derived category as introduced fifty years ago
 by Grothendieck and M.~Artin. The major achievement of Grothendieck
 and his collaborators in \cite{SGA4} was to define a theory of
 coefficients systems relative to any scheme with a collection of 
 operations, $f_*,f^*, f_!,f^!, \otimes, \uHom$,
 satisfying a set of formulas now called the
 \emph{Grothendieck six functors formalism}
  (see section \ref{sec:Gsixfunctors_intro}
	  in this introduction for more details)\footnote{The full derived
 formalism of $\ell$-adic complexes was fully established
 much later after \cite{SGA4} though, by Ekedahl in \cite{Ekedahl}.}. 
 This formalism, formulated in the language of triangulated categories,
 ultimately encodes a very general duality theory. Note however that
 the complete duality theory for $\ell$-adic sheaves was 
 completed only recently by the work of Gabber \cite{gabber3}.

The theory was also conjectured to be deeply linked with
 Quillen algebraic K-theory (see \cite[5.10, \textsection B]{Bei}).
 In fact, up to torsion and for a regular scheme $S$,
 the ext-groups between two Tate motives over $S$ should coincide with
 Adams graded parts of Quillen algebraic K-theory.\footnote{See
 page~\pageref{equation:cohKtheory} below
 for the precise statement.} 

The ideas of Beilinson were very fecund and
 not long after the publication of \cite{Bei},
 three candidates for a triangulated category of mixed 
 motives were proposed, respectively by
 M.~Hanamura \cite{hanamura1,hanamura2,hanamura3},
 M.~Levine \cite{LevineMM}, and
 V.~Voevodsky \cite{V1,V3,FSV}. As a matter of fact, Voevodsky
 introduced two variants: using the $\h$-topology
 (obtained by allowing proper surjective maps as coverings
 together with Zariski coverings), he defined a candidate
 for a theory of \'etale motivic sheaves \cite{V1}. Inspired
 by his knowledge of these and by his work on rigidity
 results with Suslin \cite{SV0}, he introduced a more Zariski local
 version \cite{FSV} which is the one fitting in his approach to
 the proof of the Milnor conjecture and of the Bloch-Kato conjecture.
 In this book, we will focus on Voevodsky's theories.

\subsection{Voevodsky's motivic complexes}

As briefly alluded to above,
the first attempt of Voevodsky in defining the category of motivic
 complexes, in his 1992 Harvard's thesis, introduces the fundamental
 process of $\AA^1$-localization, which amounts to make the
 affine line contractible in the category of mixed motives,
 by analogy with the topological case.
It also involves the use of the $h$-topology which was to
 become fundamental in the area of motives and cohomology.
 These two ingredients given, Voevodsky defined the triangulated
 category of (effective) $h$-motives over any base in \cite{V1}.

However, Voevodsky was aware that his definition will
 give the correct answer to Beilinson's conjectural construction
 only with rational coefficients (he was aware that the torsion
 part of the theory of $\h$-motives would be closely
 related to Grothendieck's\'etale cohomology\footnote{This is
 made precise by Suslin and Voevodsky \cite{SV0} over a field,
 and we developped this idea in full generality in \cite{CD4}, using
 the main results and constructions of the present book.}).
In \cite[chap. 5]{FSV}, he introduces another definition
 of motivic complexes over a perfect field with integral
 coefficients, still using the $\AA^1$-localization process but, this
 time, introducing the notion of Nisnevich sheaves with transfers
 and their derived category (see \cite{MVW} for a detailed exposition).
At this stage, all the properties foreseen by Beilinson
 are established for this integral category over a perfect field,
 except for the construction of the motivic $t$-structure.\footnote{This hoped for $t$-structure is described in \cite[Hyp. 0.0.21]{VBletter}.
 Moreover, Voevodsky proved in \cite{FSV} that such a $t$-structure
 does not exist with integral coefficients; however, it should exist
 with rational coefficients, and, more generally,
 for $\h$-motives with integral coefficients.}
 It remained to extend this definition to arbitrary bases
 and to establish the Grothendieck six functors formalism.

The path in this direction was laid down by Voevodsky
 in \cite{cancel} where he uses the theory of relative cycles
 introduced by Suslin and Voevodsky in \cite{FSV} to extend the definition
 of sheaf with transfers. This definition was also exploited by Ivorra
 in \cite{Ivo} to extend Voevodsky's construction of geometric motivic 
 complexes over any base,
 avoiding the use of sheaves with transfers. Nevertheless,
 constructing the Grothendieck six functors
 formalism for this definition remained untouched at this point.

\subsection{Morel and Voevodsky's homotopy theory}

Soon after the introduction of Voevodsky's motivic complexes,
 F.~Morel and Voevodsky introduced the more general
 theory of $\AA^1$-homotopy of schemes \cite{MV}
 whose design is to extend
 the framework of algebraic topology to algebraic geometry
 and is built around the $\AA^1$-localization tool.
 It is within this theory that another important
 tool in motivic homotopy theory was introduced: the $\PP^1$-stabilization
 process\footnote{At that time, with the impulse of Voevodsky's
 theory, the general process of $\otimes$-inverting an
 object such as a topological circle of $\PP^1$ in an
 homotopy-theoretic way quickly was fully documented; see \cite{Hov}.}.
 From the purely motivic point of view, this amounts to
 invert the Tate motive $\ZZ(1)$ for the tensor product.
 From the homotopical point of view, this operation is much more
 involved and reveals the theory of spectra, objects which
 incarnate cohomology theories in algebraic topology.
 These two processes, of $\AA^1$-localization and $\PP^1$-stabilization,
 applied to the category of simplicial Nisnevich sheaves,
 led to the stable $\AA^1$-homotopy category of schemes
  (see \cite{Jar}, or the last chapter of \cite{JarLoc}, for instance, or
  \cite{Robalo,hoyois1} for more modern approaches)
  a triangulated category with integral coefficients,
 defined over any base, which generalizes the category
 of motivic complexes.\footnote{Heuristically, the essential
difference between stable $\AA^1$-homotopy and motivic complexes
is the presence of transfers in the later case.}

Over a perfect field, and with rational coefficients,
the relation between $\AA^1$-homotopy invariant sheaves and motives was clarified in an
unpublished paper of Morel \cite{SHQ} (with precise statements
but without proofs):
 the rational stable $\AA^1$-homotopy
 category contains the stable (\emph{i.e.} $\PP^1$-stable) version
 of the category of motivic complexes as an explicit direct
 factor, called the \emph{+-part} of
 the stable homotopy category (that is the part where the algebraic
 Hopf fibration acts as in oriented cohomology theories).\footnote{One of the goals of
this book is to provide a proof of the generalization to abitrary base schemes
of Morel's claim; see Theorem \ref{thmi:Morel} in this introduction
 and its corollary.} Then Morel introduces this $+$-part
 as a good candidate for the rational version of the
 triangulated category of motives
  (\cite[paragraph at the end of p.2]{SHQ}).
	We will dub the objects of this category \emph{Morel motives}.

In the language of motivic stable homotopy theory,
as initiated by Spitzweck in \cite{spit},
a natural candidate for the category of
motivic sheaves is the homotopy category of modules over the motivic
ring spectrum which represents motivic cohomology.
With integral coefficients, O.~R\"ondigs and
 P.A. \O stv\ae r showed that, over a field of characteristic
 $0$, this category of modules
is equivalent to the $\PP^1$-stable category of motivic complexes
 (see \cite{RNOST}).\footnote{See also Theorem \ref{thmi:DM_modules}
 in this introduction for an extension of their result to arbitrary base,
 at the price of working with rational coefficients.
 For fields of characteristic $p>0$, this has been extended
 to $\ZZ[1/p]$-linear coefficients by Hoyois, Kelly and
 {\O}stv{\ae}r in \cite{Steenrod}. Finally, using
 the results of the present book,
 this has been extended to regular schemes of
 equal characteristic in \cite{CD5}.}
 This ring spectrum was introduced by Voevodsky \cite[\textsection 6.1]{V3},
  using the theory of relative cycles. It is defined over any
  base and one is led to consider the category of modules over 
 this ring spectrum as a possible definition of the integral
 triangulated category of motives.

\subsection{Voevodsky's cross functors and Ayoub's thesis}

The definitive step towards the six functors formalism
 in motivic homotopy theory was taken up by Voevodsky in
 a series of lectures were he laid down the theory of
 \emph{cross functors}. The main theorem of this theory
 consists in giving a criterion on a system of triangulated
 categories indexed by schemes,
 equipped with a basic functoriality,
 to be able to construct exceptional functors $(f_!,f^!)$
 satisfying the properties required by Grothendieck
 six functors formalism. In particular, the system of triangulated
 categories must satisfy three notable properties:
 the \emph{$\AA^1$-localization property},
 the \emph{$\PP^1$-stability property}
 and the \emph{localization property}.
Unfortunately, only an introductory part
 on this theory was released (see \cite{Delnotes}) in which
 the basic setup is established but which does not contain
 the proof of the main result. 

The writing of this theory was accomplished by
 J.~Ayoub in the first half of his thesis (see \cite{ayoub}).
 Ayoub uses the axioms laid down by Voevodsky:
 he calls a system of triangulated categories satisfying
 the properties alluded to above a \emph{homotopy stable functor}.
 Moreover, he goes beyond the original result of Voevodsky:
 apart from the complete theory of cross functors (concerned with $f_!,f^!$),
 he also studied monoidal structures,
 \emph{constructibility} properties and their stability under the six
 operations, \emph{homotopy t-structures}
 and \emph{specialization functors}
 such as the vanishing cycle functor.
 The main example of a stable homotopy functor is the
 stable $\AA^1$-homotopy category. This was established
 independently by R\"ondigs \cite{roendigs} and Ayoub \cite{ayoub},
 who both also derived two fundamental properties: the one of relative purity
 and the proper base change isomorphism.
 One readily deduces that the category of Morel motives
 is also a homotopy stable functor. Furthermore, Ayoub's axiomatic
 approach allows a uniform treatment which also applies to
 many natural variations of the stable $\AA^1$-homotopy category
 (as recalled in \cite{ayoub2}, we may vary the topology
 as well as the coefficients in which sheaves take their values).
 However, despite its already great level of generality, Ayoub's
 work does not allow us to reach all the constructions of interest.
 For instance, it only provides the construction of the functors $f_!$
  and $f^!$ when $f$ is quasi-projective, and the finiteness
  and duality theorems only apply under hypothesises (such
  as absolute purity)
  which are far from being obvious in practice
  (Ayoub only discusses this issue for schemes of finite type
  over a perfect field, in which case this follows from the property
  of relative purity). Moreover, the techniques recalled
  in the third chapter of Ayoub's thesis do not explain how
  to construct examples out of sheaves equipped with
extra structures, such as transfers, which are fundamental tools
to understrand how algebraic cycles play a role in $\AA^1$-homotopy
theory. The extra technicalities related to this problem
(such as having a derived tensor product as well as derived pull-back
functors for suitable notions of $\PP^1$-stable sheaves with transfers)
where addressed in two approaches: the first one, by R\"ondigs and
\O stv\ae r \cite{RNOST}, uses homotopy theory together
with enriched category theory, while
the second one, due to the authors \cite{CD1}, uses abstract methods
of homotopical algebras applied to cochain complexes in Grothendieck
abelian categories. A second kind of problems which is not addressed
in the early work of Ayoub is representability
  for $K$-theory, or for Chow groups, according to Beilinson's vision.
  
  This is why, in order to discuss the original approach of Voevodsky to
  motivic sheaves alluded to above (using $\h$-sheaves or
  Nisnevich sheaves with transfers), as well as to prove
  the absolute purity theorem, we had to take over from scratch
  many of the basic constructions. This also lead us to reach a greater level
  of generality (avoiding unnecessary quasi-projectivity hypothesises)
  as well as more precise computations. To be fair,
  we should mention that, after a first version of the present text
  has been made public in 2009, Ayoub~\cite{ayoub5} reproved
  some of the representability results as well as the
  absolute purity theorem of this book in the particular case of
  the \'etale version of the motivic
  stable homotopy theory with rational coefficients.
  We should also mention right away that the integral version
  of Voevodsky's $\h$-motives
  is only fully understood in a sequel of the present book \cite{CD4}.

Finally, we would like to end this paragraph by recalling that
the problem of constructing triangulated categories
of motives related to Chow groups
with integral coefficients and which define a homotopy stable functor
is still an open problem. For insance, it is by no means obvious that 
 the category of modules over the motivic homotopy
 ring spectrum does meet the requirements of a
 homotopy stable functor. In fact, this latter property can be seen to
be equivalent to Conjecture~17 of Voevodsky in \cite{V_OpenPB}, which
states the stability of the motivic homotopy
 ring spectrum by pull-backs;
 this is made precise in this book in Prop.~\ref{prop:reformulation of
the conjecture of voevodsky}, as an application of our main constructions.

\subsection{Grothendieck six functors formalism}
\label{sec:Gsixfunctors_intro}

\begin{num} \label{num:G6FF}
We now give the precise formulation of the
 \emph{Grothendieck six functors formalism}
 (although we do not describe all the coherences yet).
 As presented here, it is extracted from the properties
 of the derived category of $\ell$-adic sheaves.\footnote{It also 
 coincides with formulas gathered by Deligne in an unpublished note
 which he graciously shared with us.}

A triangulated category $\T$, fibred over the category of schemes,
 satisfies the \emph{Grothendieck six functors formalism} if the following
 conditions hold:
\begin{enumerate}
\item There exists 3 pairs of adjoint functors as follows:
\begin{align*}
f^*:\T(X) \rightleftarrows \T(Y):f_* \, , & \ f \text{ any morphism,} \\
f_!:\T(Y) \rightleftarrows \T(X):f^! \, , & \ f \text{ any separated morphism of finite type,} \\
(\otimes,\uHom)\, , &\ \text{ symmetric closed monoidal structure on } \T(X).
\end{align*}
The functors of type $f^*$ are monoidal.
\item There exists a structure of a covariant (resp. contravariant) 
 $2$-functors on $f \mapsto f_*$, $f \mapsto f_!$ (resp. $f \mapsto f^*$,
  $f \mapsto f^!$).
\item There exists a natural transformation 
$$\alpha_f:f_! \rightarrow f_*$$
 which is an isomorphism when $f$ is proper.
 Moreover, $\alpha$ is a morphism of $2$-functors.
\item For any smooth separated morphism $f:X \rightarrow S$ in $\sch$
 of relative dimension $d$,
 there exists a canonical natural isomorphism
\begin{align*}
\piso'_f:f^* & \xrightarrow{\sim} f^!(-d)[-2d]
\end{align*}
where $?(-d)$ denotes the inverse of the Tate twist iterated $d$-times. Moreover $\piso'$ is an isomorphism of $2$-functors.
\item For any cartesian square in $\sch$:
$$
\xymatrix@=16pt{
Y'\ar^{f'}[r]\ar_{g'}[d]\ar@{}|\Delta[rd] & X'\ar^g[d] \\
Y\ar_f[r] & X,
}
$$
such that $f$ is separated of finite type,
there exist natural isomorphisms
\begin{align*}
%%Ex(\Delta^*_!)&:
g^*f_! \xrightarrow\sim f'_!{g'}^*\, , \\
%%Ex(\Delta_*^!)&:
g'_*{f'}^! \xrightarrow\sim  f^!g_*\, .
\end{align*}
\item For any separated morphism of finite type $f:Y \rightarrow X$,
 there exist natural isomorphisms
\begin{align*}
(f_!K) \otimes_X L &\xrightarrow{\ \sim\ } f_!(K \otimes_Y f^*L)\, ,\ \\
  \uHom_X(f_!(L),K) & \xrightarrow{\ \sim\ } f_* \uHom_Y(L,f^!(K))\, ,\ \\
  f^! \uHom_X(L,M)& \xrightarrow{\ \sim\ } \uHom_Y(f^*(L),f^!(M))\, .
\end{align*}
\item[\loc] For any closed immersion $i:Z \rightarrow S$
 with complementary open immersion $j$,
 there exists a distinguished triangle of natural transformations as follows:
\begin{align*}
j_!j^! &\xrightarrow{\ \alpha'_j\ } 1 \xrightarrow{\ \alpha_i\ } i_*i^*
 \xrightarrow{\ \partial_i\ } j_!j^![1]
\end{align*}
where $\alpha'_?$ (resp. $\alpha_?$) denotes the co-unit (resp. unit)
 of the relevant adjunction.
\end{enumerate}
\end{num}
%
%\begin{rem}
%These are the properties proved in \cite[tome 3]{SGA4}
 %in the context of torsion \'etale sheaves (of exponent prime
 %to the residual characteristic).
 %It agrees with formulas gathered by Deligne in an unpublished note
 %which he graciously support us with.
%If we restrict these properties to the case of the 4 functors
 %$(f^*,f_*,f_!,f^!)$ and we modify slightly the formulation of Point (4),
 %this is the outcome of the main theorem of the theory of cross functors
 %(see \cite[1.4.2]{ayoub}) -- one can also extract the remaining formulas
  %from \cite[sec. 2.3]{ayoub}.
%In this book, we will give a variation on the cross-functor theory
 %and obtain the formalism above as a consequence of simpler
 %axioms in Theorem \ref{thm:cor3_Ayoub}.
%
%In the \'etale case, (\loc) is easily checked.
 %In the motivic case considered here, it is the crucial difficulty.
%\end{rem}

\begin{num}\label{num:G6FF_duality}
The next part of Grothendieck six functors formalism
 is concerned with duality. Historically, this is the initial motivation
 behind Grothendieck six functors formalism,
 as it appears in the first account of this formalism,
 Hartshorne's notes of Grothendieck 1963/64 seminar, \cite{Har}.
 It is considered more axiomatically,
 in the case of \'etale sheaves, in
 \cite[Exp. I]{SGA5}.\footnote{The duality properties are stated
   in unpublished notes of Deligne, as part of the complete formalism.}
 In \emph{loc. cit.},
  Grothendieck states the fundamental property
  of absolute purity and indicates its fundamental link
	with duality. We state these properties
	as natural extensions of the properties given in the preceding
	paragraph; assume $\T$ satisfies the properties listed above:
\begin{enumerate}
\item[(7)] \textit{Absolute purity}.-- 
 For any closed immersion $i:Z \rightarrow S$
  of regular schemes of (constant) codimension $c$,
	there exists a canonical isomorphism:
$$
\un_Z(-c)[-2c] \xrightarrow{\sim} i^!(\un_S)
$$
where $\un$ denotes the unit object for the tensor product.
\item[(8)] \textit{Duality}.-- Let $S$ be regular scheme
 and $K_S$ be any invertible object of $\T(S)$.
 For any separated morphism $f:X \rightarrow S$ of finite type,
 put $K_X=f^!(K_S)$.
 For any object $M$ of $\T(X)$, put $D_X(M)=\uHom(M,K_X)$.
\begin{enumerate}
\item For any $X/S$ as above, $K_X$ is a dualizing object
 of $\T(X)$: the canonical map
$$
M \rightarrow D_X(D_X(M))
$$
is an isomorphism.
\item For any $X/S$ as above, and any objects $M,N$ of $\T(X)$,
 we have a canonical isomorphism
$$D_X(M\otimes D_X(N))\simeq \uHom_X(M,N)\, . $$
\item For any morphism between separated $S$-schemes of finite type
$f:Y\To X$, we have natural isomorphisms
\begin{align*}
D_Y(f^*(M))& \simeq f^!(D_X(M))\\
f^*(D_X(M))& \simeq D_Y(f^!(M))\\
D_X(f_!(N))& \simeq f_*(D_Y(N))\\
f_!(D_Y(N))& \simeq D_X(f_*(N))\ .\\
\end{align*}
\end{enumerate}
\end{enumerate}
\end{num}

%\begin{rem}
%Dualizing objects (complexes) were introduced by Grothendieck.
%Existence conditions for them are considered axiomatically
 %in \cite[Exp. I]{SGA5} (in the \'etale context),
 %and extended in \cite[V]{SGA4.5}.
 %In the \'etale context, the optimal result is proved in \cite{gabber3}.
 %In the motivic context,
  %existence conditions are  given in \cite[2.3.10]{ayoub}
  %(even without the symmetry axiom on the tensor structure).
%In this book, we give a criterion of existence in the spirit
 %of \cite{gabber3} which extends that of Ayoub
 %in Corollary \ref{cor:localduality2}.
%\end{rem}

\begin{num}\label{num:G6FF_continuity}
The last property we want to exhibit as
 a natural extension of Grothendieck six functors formalism 
 is the compatibility with projective limits of schemes.
 The basis for the next statement is \cite[Exp. VI]{SGA4} though
 it does not appear explicitly.
 As in the case of the duality property, it should involve some
 finiteness assumption (constructibility) on the objects of $\T$.
 Note the formulation below is valid for
 an arbitrary triangulated monoidal category $\T$ fibred over schemes.
\begin{enumerate}
\item[(9)] \textit{Continuity}.-- 
Let $(S_\alpha)_{\alpha\in A}$ be an essentially affine projective system
 of schemes.
\index{word}{projective system, of schemes}
 Put $S=\varprojlim_{\alpha\in A}S_\alpha$.

Then the canonical functor
$$
2\mbox{-}\varinjlim_\alpha \T(S_\alpha)\To \T(S)
$$
is an equivalence of monoidal triangulated categories.
\end{enumerate}
\end{num}

The purpose of this book is to discuss such a formalism in
various contexts of motivic sheaves.

\section{Voevodsky's motivic complexes}
\label{seci:motivic_complexes}

The primary goal of this treatise is to develop the theory
 of Voevodsky motives,
 integrally over any base scheme\footnote{In this introduction,
  all schemes will be assumed to be noetherian of finite dimension.},
 within the framework of sheaves with transfers.
 Actually, we can define Voevodsky's motives with coefficients
 in an arbitrary ring $\Rc$ and prove all the results stated below
 in that case, but we restrict this presentation to integral coefficients
 for simplicity.

After refining and completing Suslin-Voevodsky's theory
 of relative cycles, we introduce the category 
$\renewcommand{\Rc}{\ZZ} \smcx S$
 of integral finite correspondences over smooth $S$-schemes
 and the related notion of (Nisnevich) sheaves 
 with transfers over a base scheme $S$ (Def. \ref{df:premotivic_ftr})
 as in the usual case of a perfect base field.
 Following the idea of stable homotopy,
 we define the triangulated category $\DMV(X)$ of
 \emph{stable motivic complexes} (see Def. \ref{df:Nis_DMe&DM})
 as the $\PP^1$-stabilization 
 of the $\AA^1$-localization of the derived category of
 the (Grothendieck) abelian category of sheaves with transfers
 over $S$.

One easily gets that
 the fibred category $\DMV$ is equipped with the basic functoriality
 needed by the cross-functor formalism. The main difficulty
 is the localization property,
 labelled \loc in Paragraph \ref{num:G6FF}.
Unfortunately, though all the functors involved in the formulation
 of \loc are well-defined for $\DMV$,
 we can only prove this property when $S$ and $Z$ are smooth
 over some base scheme (see Prop. \ref{prop:DMV_supp_wloc}).
In particular, the formalism of stable homotopy functors does not apply.
However, we are able to construct the six operations for $\DMV$
 using the method of Deligne, used in \cite[XVII]{SGA4},
 and partially get the Grothendieck six functors 
 formalism:
\begin{thmi}[see Th. \ref{thm:DM_six_functors}]
The triangulated category $\DMV$, fibred over the category of schemes,
 satisfies the following part of the properties stated in
 Paragraph \ref{num:G6FF}:
\begin{itemize}
\item properties (1), (2), (3) (i.e. the construction of $f_!$ and $f^!$
in $\DM$ for any separated morphism of finite type $f$),
\item property (4) when $f$ is an open immersion or $f$ is projective
 and smooth,
\item property (5) when $g$ is smooth or $f$ is projective and smooth,
\item property (6) when $f$ is projective and smooth,
\item Property (Loc) when $S$ and $Z$ are smooth over some common base scheme.
\end{itemize}
\end{thmi}
One of the applications of this theory is that we get a well-defined
 integral motivic cohomology theory for any scheme $X$:
$$
\renewcommand{\Rc}{\ZZ}
\HHm n m X=\Hom_{\DMV(X)}\big(\un_X,\un_X(m)[n]\big)
$$
which enjoys the following properties
 (see section \ref{sec:motivic_coh}):
\begin{itemize}
\item it admits a ring structure,
 pullback maps associated with any morphism of schemes compatible
 with the ring structure,
\item it admits push-forward maps with respect to projective morphisms
 between schemes smooth over some common base,
 or with respect to some finite morphisms (for example finite flat;
 see Paragraph \ref{num:corr_cohm}),
\item it coincides with Voevodsky's motivic cohomology groups
 when $X$ is smooth over a perfect field
 (see Example \ref{ex:motivic_coh&Chow}); 
 in particular one gets the following identification
 with higher Chow groups:
$$
\renewcommand{\Rc}{\ZZ}
\HHm n m X=CH^{m}(X,2m-n),
$$
\item it admits Chern classes and satisfies the projective bundle formula,
\item it admits a localization long exact sequence associated
 with a closed immersion of schemes which are smooth over some common base.
\end{itemize}

As in the classical case, any smooth $S$-scheme $X$
 admits a motive $M_S(X)$ in $\DMV(S)$.
 Moreover, one defines the Tate motive $\un_S(1)$ as the reduced
 motive of $\PP^1_S$. 
 We define the category of \emph{constructible motives}
 $\DMV_c(S)$
 as the thick triangulated subcategory of $\DMV$ generated
 by the objects of the form $M_S(X)(n)$ for a smooth $S$-scheme $X$
 and an integer $n \in \ZZ$, where $?(n)$ refers to the $n$-th
 Tate twist.
One gets the following generalization of the classical result
 obtained by Voevodsky over a perfect field:
\begin{thmi}[see Th. \ref{thm:DM_constructible_geometric}]
\label{thmi:geometric_motives}
A motive $M$ in $\DMV(S)$ is constructible if and only if
 it is compact.\footnote{Recall that $M$ is
 compact if the functor $\Hom(M,-)$ commutes with arbitrary direct
 sums.}

The category $\DMV_c(S)$ is equivalent to the
 category obtained from the bounded homotopy category 
%$\renewcommand{\Rc}{\ZZ} K^b(\smcx S)$
 of the additive category $\renewcommand{\Rc}{\ZZ} \smcx S$
 in the following way:
\begin{itemize}
\item take the Verdier quotient modulo the thick 
 triangulated subcategory generated by:
\begin{itemize}
\item
for any Nisnevich distinguished square 
$\xymatrix@=10pt{W\ar^k[r]\ar_g[d] & V\ar^f[d] \\ U\ar^j[r] & X}$ 
of smooth $S$-schemes:
$$
[W]
 \xrightarrow{g_*-k_*} [U] \oplus [V]
 \xrightarrow{j^*+f^*} [X]
$$
\item for any smooth $S$-scheme $X$, 
$p:\AA^1_X \rightarrow X$ the canonical projection:
$$
[\AA^1_X] \xrightarrow{p_*} [X],
$$
\end{itemize}
\item invert the Tate twist,
\item take the pseudo-abelian envelope.
\end{itemize}
\end{thmi}

The triangulated category $\DMV_c(X)$
 is stable under the operations $f^*$ for all $f$,
and $f_*$, when $f$ is smooth projective, as well as
 $\otimes$, but we cannot prove the stability for the
 other operations of $\DMV$ and \emph{a fortiori} do not get
 the duality properties (7) and (8) of the Grothendieck six functors
 formalism.

However, we are able to prove the continuity property (9)
 for the category $\DMV_c$:
$$
2\mbox{-}\varinjlim_\alpha \DMV_c(S_\alpha)\simeq \DMV_c(S),
$$
with the restriction that the transition morphisms
 of $(X_\alpha)$ are affine and dominant
 (see Theorem \ref{thm:DM_continuity}).
Note this result allows us to extend
 the comparison of motivic cohomology with higher
 Chow groups
 to arbitrary regular schemes of equal characteristics.

Finally, although we could not prove all the expected properties
of the six operations in $\DM$, we prove that the six operations
behave as expected in $\DM$ if and only if
Conjecture~17 of Voevodsky in \cite{V_OpenPB} is true; see
Prop.~\ref{prop:reformulation of the conjecture of voevodsky}.
\section{Beilinson motives}

\subsection{Definition and fundamental properties}

As anticipated by Morel,
 the theory of mixed motives with rational coefficients is much simpler
 and we succeed in establishing a complete formalism for them.
However, there are many candidates for $\QQ$-linear motivic sheaves
over a scheme $X$:
there are Voevodsky's $\h$-motives $\DMV_{\h,\QQ}(X)$,
Voevodsky's motivic sheaves constructed out $\QQ$-linear
sheaves with transfers $\DM(X,\QQ)$, Morel motives $\SH_{\QQ}(X)_+$,
$\QQ$-linear \'etale motives
$\Der_{\AA^1,\et}(X,\QQ)\simeq\SH_{\et,\QQ}(X)$
(also introduced by Morel, and used as length by Ayoub).
Our goal is not only to prove that the six operations act on
 each of these candidates, but also to compare all these
 versions of motivic sheaves with one another. In fact, our strategy
 consists in producing yet another candidate for $\QQ$-linear motivic sheaves,
 namely the one of \emph{Beilinson motives} $\DMB(S)$, for which
 we can prove all the features we want for it, and use them to compare
Beilinson motives with all the other versions of $\QQ$-linear motivic
sheaves mentionned above.
 
 More precisely, we construct, out of the rational motivic
 stable homotopy category and the ring spectrum associated
 with rational Quillen K-theory
 a $\QQ$-linear triangulated category $\DMB(X)$, 
 which we call the \emph{triangulated category of Beilinson motives}
 (see Def. \ref{df:Beilinson_motives}).
Essentially by construction, in the case where $X$ is
 regular, we have a natural identification\label{equation:cohKtheory}
 $$\Hom_{\DMB(X)}(\QQ_X,\QQ_X(p)[q])\simeq
 \mathit{Gr}^p_\gamma K_{2p-q}(X)_\QQ\, ,$$
where the right-hand side is the graded part of the algebraic K-theory of $X$
 with respect to the $\gamma$-filtration.
 
 These groups were first regarded by Beilinson
 as the rational motivic cohomology groups. We call them
 the \emph{Beilinson motivic cohomology groups}.

Part of the interest of our definition is that the localization
 property \loc can be easily deduced from its validity
 for the stable homotopy category.
Therefore, the cross-functor formalism and more generally, our generalization
of the results of Ayoub can be applied to $\DMB$.
\begin{thmi}[see Cor. \ref{cor:DMBmotivic} and Th. \ref{thm:cor3_Ayoub}]
All the standard Grothendieck six functors formalism
 (see Paragraph \ref{num:G6FF}) is verified by the
 fibred triangulated category $\DMB$.
\end{thmi}

% The six operations of Grothendieck act naturally on $\DMB$: we have
% all the functors $\otimes$, $\sHom$, $f^*$, $f_*$, $f_!$, $f^!$
% (with the same properties as above). The following
% result is a direct consequence of its analog in the stable homotopy
% category $\SH$.
% 
% \begin{thmi}[localization]
% For any closed immersion $i:Z\to X$ with complement open immersion $j:U\To X$,
% we have the six gluing functors:
% $$\xymatrix{
% \DMB(U)\ar@<5pt>[r]^{j_!} \ar@<-5pt>[r]_{j_*}& \DMB(X)\ar[l]|{j^*} \ar[r]|{i_*}
% & \DMB(Z) \ar@<5pt>[l]^{i^!}\ar@<-5pt>[l]_{i^*}\, .
% }$$
% \end{thmi}
% 
% The result above allows us to use Ayoub's results on cross functors, which
% gives the following properties.
% We recall that, for $f:Y\To X$ proper, we have $f_!\simeq f_*$.
% Moreover, if $f$ is smooth and quasi-projective of relative dimension $d$, we have,
% for any object $M$ of $\DMB(X)$, a purity isomorphism
% $$f^*(M)(d)[2d]\simeq f^!(M)$$
% These identifications allow one to interpret the following theorem as the
% proper base change formula and the smooth base change formula in $\DMB$.
% 
% \begin{thmi}[base change isomorphisms]
%  For any Cartesian square
% $$\xymatrix{
% Y'\ar[r]^{f'}\ar[d]_{g'}& X'\ar[d]^g\\
% Y\ar[r]_f& X,
% }$$
% with $f$ separated of finite type,
% the natural transformation $g^*\, f_!\To f'_!\, g^{\prime*}$ as well as its
% transposed transformation $g^\prime_*\, f^{\prime!}\, \To f^!\, g_* $ are invertible.
% \end{thmi}

Concerning duality for Beilinson motives,
 we first deduce from Quillen's localization theorem in algebraic K-theory
 the absolute purity theorem:

\begin{thmi}[see Th. \ref{DMBpurity}] \label{thmi:purity}
The absolute purity property (see \ref{num:G6FF_duality}(7))
 holds for $\DMB$.
\end{thmi}

As said before, this result is not enough to establish duality
 for Beilinson motives. We first have to use descent theory
 and resolution of singularities (as first explained by Grothendieck
 in \cite[I.3]{SGA5}).
 Using the existence of trace maps in algebraic K-theory, we prove
 the following result:

\begin{thmi}[$\h$-descent, see Th. \ref{thm:DMBseparated} 
 and Th. \ref{excisionGaloisalteration}]
\label{thmi:qfh-descent}
Consider a finite group $G$ and a pullback square of schemes
\begin{equation*}
\xymatrix{
T\ar[r]^h\ar[d]_g & Y\ar[d]^f \\
Z\ar[r]_i&X
}
\end{equation*}
in which $Y$ is endowed with an action of $G$ over $X$.
 Put $U=X-Z$ and assume the following three conditions are satisfied.
\begin{itemize}
\item[(a)] The morphism $f$ is proper and surjective.
\item[(b)] The induced morphism $f^{-1}(U)\To U$ is finite.
\item[(c)] The morphism $ f^{-1}(U)/G\To U$ is generically radicial.
\end{itemize}
Put $a=f \circ h=i \circ g$.
 Then, for any object $M$ of $\DMB(X)$,
 we get a canonical distinguished
 triangle in $\DMB(X)$:
$$
M \longrightarrow  i_*\,  i^*(M)\oplus  f_*\,  f^*(M)^{G}
  \longrightarrow  a_*\,  a^*(M)^{G} \longrightarrow M[1]
$$
where $?^G$ means the invariants under the action of $G$,
 and the first (resp. second) map of the triangle
 is induced by the difference (resp. sum) of the
 obvious adjunction morphisms.
% Consider a Galois alteration $p:X'\To X$ of group $G$ (i.e. $p$
% is an alteration, while $G$ is a finite group acting on $X'$ over $X$,
% such that, generically, $X'/G\To X$ is finite surjective and radicial) , as well as
% a closed subscheme $Z\subset X$, such that $U=X-Z$ is normal, and
% such that the induced map $p_U:U'=p^{-1}(U)\To U$ is a finite morphism.
% Then the pullback square
% $$\xymatrix{
% Z'\ar[r]^{i'}\ar[d]_q&X'\ar[d]^p\\
% Z\ar[r]^i&X
% }$$
% induces a canonical distinguished triangle
% $$M \To  i_*\,  i^*(M)\oplus  p_*\,  p^*(M)^{G} \To  i_*  q_*\,  q^* \,  i^*(M)^{G} \To M[1]$$
% for any object $M$ of $\DMB(X)$.
\end{thmi}
In fact, we show that this apparently simple result
 implies a much stronger descent property for the fibred triangulated
 category $\DMB$: descent for the $\h$-topology,
 thus, in particular, \'etale descent
 flat descent, as well as proper descent. The general fact that, in the
 presence of the six operations, the property of $\QQ$-linear
 $\h$-descent is essentially equivalent to the presence
 of a suitable theory trace maps is a key feature of this text;
 this is developped systemattically in
 Chapter \ref{sec:theorie_descente} of this book. This will be at the heart
 of our main comparison results explained below.

\subsection{Constructible Beilinson motives}

The next step towards duality for Beilinson motives
 is the definition of a suitable finiteness condition.
 As in the case of Voevodsky motives,
 we define the category of \emph{Beilinson constructible motives},
 denoted by $\DMV_{\mathcyr{B},c}(X)$,
 as the thick subcategory of $\DMB(X)$ generated by the motives of the form
 $M_X(Y)(p):=f_!\, f^!(\QQ_X)(p)$
 for $f:Y\To X$ separated smooth of finite type, and $p\in \ZZ$.
This category coincides with the full subcategory
 of compact objects in $\DMB(X)$.\footnote{Note the striking analogy
 with perfect complexes.}

The usefulness of this definition comes from the following result,
 which is the analog of Gabber's finiteness theorem in the  $\ell$-adic
 setting.
Analogously, its proof relies on absolute purity,
 (a weak form of) proper descent
 as well as Gabber's weak uniformization theorem.\footnote{i.e. that,
locally for the $\h$-topology, any excellent scheme is regular, and any
closed immersion between excellent schemes is
the embedding of a strict normal crossing divisor
 into a regular scheme.}
 \begin{thmi}[finiteness, see Th. \ref{thm:DMBc_6operations}]
The subcategory $\DMV_{\mathcyr{B},c}$ is stable under
 the six operations of Grothendieck when restricted
 to excellent schemes.
\end{thmi}

The final statement concerning Grothendieck six functors
 formalism in the setting of Beilinson motives is that,
 when one restricts to constructible Beilinson motives
 and separated $B$-schemes of finite type for an
 excellent scheme $B$ of dimension less than $2$,
 the complete formalism is available:\footnote{There is
 a way to avoid this extra hypothesis to get duality theorems
 (i.e. to work with
 quasi-excellent schemes over a regular base in full generality).
 However, this comes at the price of higher coherence results (i.e.
 of promoting the construction $f\mapsto f^!$ to $\infty$-categories).
 See \cite{Cis5}.\label{footnotegeneralduality}}
 
\begin{thmi}[see Th. \ref{thm:DMBc_duality} and Prop. \ref{DMBc_continuity}]
The fibred category $\DMBc$ over the category of schemes described
 above satisfies the complete Grothendieck six functors formalism
 described in section \ref{sec:Gsixfunctors_intro},
 in particular the \emph{duality property} \ref{num:G6FF_duality}(8)
 and the \emph{continuity property} \ref{num:G6FF_continuity}(9).
\end{thmi}

%
%\begin{thmi}[continuity]
%Let $S$ be a scheme which
%is the limit of an essentially affine projective system of schemes
%$\{S_\alpha\}$. Then there is a canonical
%equivalence of triangulated categories
%$$2\mbox{-}\varinjlim_\alpha\DMV_{\mathcyr{B},c}(S_\alpha)\simeq \DMV_{\mathcyr{B},c}(S)\, .$$
%\end{thmi}
%
%The following theorem uses absolute purity, (a weak form of) proper descent,
%as well as Gabber's weak uniformization theorem (i.e. that,
%locally for the $\h$-topology, any excellent scheme is regular, and any
%closed immersion between excellent schemes is
%the embedding of a strict normal crossing divisor into a regular scheme).
%Its proof relies on a mix of arguments of Gabber and Ayoub.
%
%
%
%The absolute purity theorem, proper descent, de~Jong's resolution
%of singularities by Galois alterations, and Ayoub's methods lead to:
%
%\begin{thmi}[duality]
%Let $B$ be an excellent scheme of dimension $\leq 2$.
%For any separated $B$-morphism $f:X\To S$ between $B$-schemes
%of finite type, with $S$ regular, $f^!(\QQ_S)$ is a dualizing object in $\DMB(X)$.
%\end{thmi}

\subsection {Comparison theorems} \label{seci:comparison}

In the historical part of this introduction,
 we saw many approaches for the triangulated category of (rational) motives.
 We succeed in comparing them all with our definition of Beilinson
 motives.

Denote by $\BGL_S$ the algebraic K-theory spectrum in
Morel and Voevodsky's stable homotopy category $\SH(S)$.
By virtue of a result of Riou, the $\gamma$-filtration on K-theory
induces a decomposition of $\BGL_{S,\QQ}$:
$$\BGL_{S,\QQ}\simeq \bigoplus_{n\in\ZZ}\HBx S(n)[2n]\, .$$
The ring spectrum $\HBx S$ represents Beilinson motivic cohomology.
Almost by construction,
 the category $\DMB(S)$ is the full subcategory of $\SH_\QQ(S)$
 which consists of objects $E$ such that the unit map $E\To \HBx S\otimes E$
 is an isomorphism. In fact, our first comparison result
 relates the theory of Beilinson motives with the approach
 of Spitzweck, R\"ondigs and \O stv\ae r through modules over a ring spectrum:

\begin{thmi}[see Th. \ref{thm:DMB&HB-mod}] \label{thmi:DM_modules}
For any scheme $S$, there is a canonical equivalence of categories
$$\DMB(S)\simeq \ho(\Mod{\HBx S})$$
where the right hand side denotes the homotopy category
 of modules over the ring spectrum $\HBx S$.
\end{thmi}

The next comparison involves
 the $\h$-topology: we recall that
 this is the Gro\-then\-dieck topology
 on the category of schemes, generated by \'etale surjective morphisms
 and proper surjective morphisms.
 The first published work of Voevodsky
 on triangulated categories of mixed motives \cite{V1},
 introduces the $\AA^1$-homotopy category 
 of the derived category of $\h$-sheaves. 
We consider a $\QQ$-linear and $\PP^1$-stable version of it,
 which we denote by $\uDMV_{\h,\QQ}(S)$.
By construction, for any $S$-scheme of finite type $X$,
 there is a $\h$-motive $M_S(X)$ in $\uDMV_{\h,\QQ}(S)$.
We define $\DMV_{\h,\QQ}(S)$ 
 as the smallest triangulated full subcategory
 of $\uDMV_{\h,\QQ}(S)$ which is stable under (infinite) direct sums,
 and which contains the objects $M_S(X)(p)$,
 for $X/S$ \emph{smooth of finite type}, and $p\in \ZZ$.
Using $\h$-descent in $\DMB$,
 we get the following comparison result.

\begin{thmi}[see Th. \ref{plongement01}]
If $S$ is excellent, then we have canonical equivalences of categories
$$\DMB(S) \simeq \DMV_{\h,\QQ}(S)\, .$$
\end{thmi}
In fact, we first prove this result for the variant of $\DMV_{\h,\QQ}(S)$
 obtained by replacing everywhere the $\h$-topology
 by the $\qfh$-topology -- in the later, also introduced by
 Voevodsky, coverings are generated 
 by \'etale covers and finite surjective morphisms.
 In particular, we get an equivalence of categories:
  $\DMV_{\h,\QQ}(S)\simeq \DMV_{\qfh,\QQ}(S)$.
This result allows us to link Beilinson motives with Voevodsky's motivic
 complexes. Let us denote by $\DMV_\QQ$ the rationalization of the fibred
 category of stable motivic complexes
 alluded to in Paragraph \ref{seci:motivic_complexes}. %\footnote{Beware
% that in general, we do not know if $\DMV_\QQ(S)$ is the rationalization
% of the integral category $\DMV(S)$. In fact, there is a well
% defined triangulated functor
%$$
%\DMV(S) \otimes \QQ \rightarrow \DMV_\QQ(S)
%$$
%but we can only prove it is an isomorphism when $S$ is either 
% a $\QQ$-scheme or a regular scheme
% (see Proposition \ref{prop:changeofcoef_DM}).}
 Using the preceding result in the case of the $\qfh$-topology,
  we prove:

\begin{thmi}[see Th. \ref{comparisonDMBDMV}]
If $S$ is excellent and geometrically unibranch, then there is a canonical
equivalence of categories
$$\DMB(S) \simeq \DMV_\QQ(S) \, .$$
\end{thmi}
In particular, given such a scheme $S$,
 we get a description of $\DMBc(S)$
 as in Theorem~\ref{thmi:geometric_motives} cited above.
Voevodsky's integral (resp. rational) motivic cohomology
 is represented in $\SH(S)$ by a ring spectrum $\Hmx S$
 (resp. $\Hmrx S$). The preceding theorem immediately gives
 an isomorphism of ring spectra:\footnote{Note in particular that, 
when $S$ is regular,
 we get an isomorphism:
$$
\renewcommand{\Rc}{\ZZ}
\HHm p q S \otimes \QQ
 \simeq \mathit{Gr}^p_\gamma K_{2p-q}(S)_\QQ
$$
which extends the known isomorphism when $S$
 has equal characteristics. 
 It is natural with respect to pullbacks, Gysin morphisms,
 as well as compatible with products and Chern classes.}
$$\HBx S\simeq \Hmrx S\, .$$
As Beilinson motivic cohomology ring spectra over
 different bases are compatible with pullbacks,
 we easily deduce the following corollary
 which solves affirmatively conjecture 17 of
 \cite{V_OpenPB} in some cases, and up to torsion:
\begin{cori}
For any morphism $f:T \rightarrow S$ of
 excellent geometrically unibranch schemes,
 the canonical map
$$
 f^*\Hmrx S \rightarrow\Hmrx T
$$
is an isomorphism of ring spectra.
\end{cori}

The next comparison statement is concerned with the approach of Morel,
according to whom the category $\SH_\QQ(S)$ can be decomposed into two factors,
 one of them being $\SH_{\QQ}(S)_+$,
 that is the part of $\SH_\QQ(S)$ on which
 the map $\epsilon:S^0_\QQ\To S^0_\QQ$,
 induced by the permutation of the factors in $\GG\wedge\GG$,
 acts as $-1$. Let $S^0_{\QQ+}$ be the unit object of $\SH_\QQ(S)_+$.

Using the presentation of Beilinson motives in terms
 of $\HB$-modules (Theorem \ref{thmi:DM_modules} cited above)
 as well as Morel's computation of the motivic sphere spectrum 
 in terms of Milnor-Witt K-theory, we obtain a proof
 of a statement, which, in the case where $S$ is the spectrum of a field,
was claimed by Morel in \cite{SHQ}:

\begin{thmi}[see Th. \ref{thm:Morel}] For any scheme $S$,
 the canonical map $S^0_{\QQ+}\To \HBx S$ is an isomorphism.
 \label{thmi:Morel}
\end{thmi}

In fact, we even get the following corollary:
\begin{cori}
For any scheme $S$, there is a canonical equivalence of categories
$$\SH_\QQ(S)_+\simeq \DMB(S)\, .$$
\end{cori}
Recall from Morel theory that,
 when $-1$ is a sum of squares in all the residue fields of $S$,
 $\epsilon$ is equal to $-Id$ on the whole of $\SH_\QQ(S)$.
 Thus in that particular case
 (e.g. $S$ is a scheme over an algebraically closed field),
 the category of Beilinson motives coincide with the rational
 stable homotopy category.
In general, we can introduce according to Morel 
 the \'etale variant of $\SH_\QQ(S)$ denoted by
 $\Der_{\AA^1,\et}(S,\QQ)$.\footnote{In brief,
 this is the $\PP^1$-stabilization of the $\AA^1$-localization
 of the derived category of sheaves of $\QQ$-vector spaces
 over the \emph{lisse-\'etale site} of $S$.}
 As locally for the \'etale topology, $-1$ is always a square,
 and because $\DMB$ satisfies \'etale descent, we get the following
 final illuminating comparison
 statement.\footnote{In particular, the finiteness
theorem as well as the duality property
also hold for $\Der_{\AA^1,\et}(-,\QQ)$.
The finiteness theorem and the duality theorem may be deduced
from \cite{ayoub} (Scholie 2.2.34 and Theorem 2.3.73 respectively)
when one restricts to quasi-projective schemes over a field or
over a discrete valuation ring. Nevertheless, even if
one is eager to accept such restrictions, over a discrete
valuation ring, the proof relies
in an essential way on the absolute purity property
(Theorem \ref{thmi:purity} stated above) which is
proved in the present text.}

\begin{cori}
For any scheme $S$, there is a canonical equivalence of
categories
$$\Der_{\AA^1,\et}(S,\QQ)\simeq \DMB(S)\, .$$
\end{cori}

Let us draw a conclusive picture which summarize 
 most of the comparison results we obtained:
\begin{cori}
Given any scheme $S$, the category $\DMB(S)$ is a full
 subcategory of the rational stable homotopy category
 $\SH_\QQ(S)$.
Given a rational spectrum $E$ over $S$, the following conditions
 are equivalent:
\begin{enumerate}
\item[(i)] $E$ is a Beilinson motive,
\item[(ii)] $E$ is an $\HBx S$-module,
\item[(iii)] $E$ satisfies \'etale descent,
\item[(iii')] ($S$ excellent) $E$ satisfies $\qfh$-descent,
\item[(iii'')] $(S$ excellent) $E$ satisfies $\h$-descent,
\item[(iv)] ($S$ excellent geometrically unibranch)
 $E$ admits transfers,
\item[(v)] the endomorphism $\epsilon \in \End(S^0_\QQ)$
 acts by $-Id$ on $E$ \emph{i.e.} $\epsilon \otimes 1_E=-1_E$.
\end{enumerate}
\end{cori}

\begin{remi} (see Corollary \ref{cor:carac_KBL-local})
 Points (iv) and (v) are related to the orientation
 theory for spectra (not only ring spectra).
 In fact, $\HBx S$ is the universal orientable rational ring spectrum
 over $S$.
\end{remi}

Let $\QQ.\sm_S$ be the $\QQ$-linear envelop of the category $\sm_S$.
 One obtains (see Example \ref{ex:A^1-Nis_Et_compacity_stable}
  in conjunction with Par. \ref{par:link_stablehomotopy_homology})
  that the full subcategory of compact objects of $\SH_\QQ(S)$
 is equivalent
 to the category obtained from the homotopy category $\K^b(\QQ.\sm_S)$
 by performing the following operations:
\begin{itemize}
\item take the Verdier quotient modulo the thick 
 triangulated subcategory generated by:
\begin{itemize}
\item
for any Nisnevich distinguished square 
$\xymatrix@=10pt{W\ar^k[r]\ar_g[d] & V\ar^f[d] \\ U\ar^j[r] & X}$ 
of smooth $S$-schemes:
$$
\QQ_S(W)
 \xrightarrow{g_*-k_*} \QQ_S(U) \oplus \QQ_S(V)
 \xrightarrow{j^*+f^*} \QQ_S(X)
$$
\item for any smooth $S$-scheme $X$, 
$p:\AA^1_X \rightarrow X$ the canonical projection:
$$
\QQ_S(\AA^1_X) \xrightarrow{p_*} \QQ_S(X).
$$
\end{itemize}
\item invert the Tate twist,
\item take the pseudo-abelian envelope.
\end{itemize}
Let us denote by $\renewcommand{\Rc}{\QQ} \DMtgmx S$ this category.
We finally obtain the following concrete description
 of Beilinson constructible motives:
\begin{cori}
Given any scheme $S$, the category $\DMBc(S)$
 is equivalent to the full subcategory of
 $\renewcommand{\Rc}{\QQ} \DMtgmx S$ spanned by the
 objects $E$ which satisfy one the following equivalent conditions:
\begin{enumerate}
\item[(i)] (\textit{Galois descent}) given any smooth $S$-scheme $X$
 and any Galois $S$-cover $f:Y \rightarrow X$ of group $G$, the canonical map
 $E \otimes \QQ_S(Y)/G \rightarrow E \otimes \QQ_S(X)$ is an isomorphism,
\item[(ii)] (\textit{Orientability}) $\epsilon$ acts by $-Id$ on $E$,
\end{enumerate}
\end{cori}
Recall again the following remarks:
\begin{enumerate}
\item When $(-1)$ is a sum of square in every residue fields
 of $S$, conditions (i), (ii) are true for any rational spectrum
 $E$ over $S$.
\item When $S$ is excellent and geometrically unibranch,
 the category $\DMBc(S)$ is equivalent to the
 category of rational geometric Voevodsky motives (same definition
 as in Theorem \ref{thmi:geometric_motives} 
 but replacing $\ZZ$ by $\QQ$).
% This follows from the fact that, with rational coefficients,
%  having an action of finite correspondences is a property
%  and not a structure. 
\end{enumerate}

\subsection{Realizations}

The last feature of Beilinson motives is that
 they are easily realizable in various cohomology theories.
To get this fact, we use the setting of modules
 over a strict ring spectrum.\footnote{i.e. we say a ring spectrum is \emph{strict} if it is a commutative monoid in the underlying model category.}
Given such a ring spectrum $\mathcal{E}$ in $\DMB(S)$,
one can define, for any $S$-scheme $X$, the triangulated category
$$\Der(X,\mathcal{E})=\ho(\Mod{\mathcal{E}_X})\, ,$$
where $\mathcal{E}_X=f^*\mathcal{E}$, for $f:X\To S$
the structural map.

We then have realization functors
$$\DMB(X)\To \Der(X,\mathcal{E}) \ , \quad M \mapsto \mathcal{E}_X\otimes_X M$$
which commute with the six operations of Grothendieck.
Using Ayoub's description of the Betti realization, we obtain:
\begin{thmi}
If $S=\spec k$ with $k$ a subfield of $\CC$, and if
$\mathcal{E}_{\mathit{Betti}}$ represents Betti cohomology
in $\DMB(S)$, then, for any $k$-scheme of finite type, the full subcategory of compact objects of
$\Der(X,\mathcal{E}_{\mathit{Betti}})$ is canonically
equivalent to the derived category of constructible sheaves of geometric origin $\Der^b_c(X(\CC),\QQ)$.
\end{thmi}

More generally, if $S$ is the spectrum of some field $k$,
given a mixed Weil cohomology $\mathcal{E}$,
with coefficient field (of characteristic zero) $\mathbf{K}$,
we get realization functors
$$\DMV_{\mathcyr{B},c}(X)\To \Der_c(X,\mathcal{E}) \ , \quad M \mapsto \mathcal{E}_X\otimes_X M$$
(where $\Der_c(X,\mathcal{E})$ stands for the category of compact objects of $\Der(X,\mathcal{E})$),
which commute with the six operations of Grothendieck (which
preserve compact objects on both sides). Moreover,
the category $\Der_c(S,\mathcal{E})$ is then canonically equivalent to the
bounded derived category of the abelian category of finite dimensional
$\mathbf{K}$-vector spaces. As a by-product, we get the following
concrete finiteness result: for any $k$-scheme of finite type $X$, and for any
objects $M$ and $N$ in $\Der_c(X,\mathcal{E})$, the
$\mathbf{K}$-vector space $\Hom_{\Der_c(X,\mathcal{E})}(M,N[n])$
is finite dimensional, and it is trivial for all but a finite number of values of $n$.

If the field $k$ is of characteristic zero, this
abstract construction gives essentially the usual categories of coefficients
(as seen above in the case of Betti cohomology),
and in a sequel of this work, we shall prove that one recovers in this
way the derived categories of constructible $\ell$-adic sheaves
(of geometric origin) in any characteristic. But something new happens in positive
characteristic:
\begin{thmi}
Let $V$ be a complete discrete valuation ring of mixed characteristic,
with field of functions $K$, and residue field $k$. Then rigid
cohomology is a $K$-linear mixed Weil cohomology, and thus
defines a ring spectrum $\mathcal{E}_{\mathit{rig}}$ in $\DMB(k)$.
We obtain a system of closed symmetric monoidal triangulated categories
$\Der_{\mathit{rig}}(X)=\Der_c(X,\mathcal{E}_{\mathit{rig}})$, for any $k$-scheme
of finite type $X$, such that
$$\Hom_{\Der_{\mathit{rig}}(X)}(\unit_X,\unit_X(p)[q])\simeq H^q_{\mathit{rig}}(X)(p)\, ,$$
as well as realization functors
$$R_{\mathit{rig}}:\DMV_{\mathcyr{B},c}(X)\To \Der_{\mathit{rig}}(X)$$
which preserve the six operations of Grothendieck.
\end{thmi}

\section{Detailed organization}

The book is organized in four parts that we now review in more details.

\subsection{Grothendieck six functors formalism (Part 1)}

The first part is concerned with the formalism
 described in section \ref{sec:Gsixfunctors_intro} above.
 It is the foundational part of this work.

We use the language of fibred categories
 (introduced in \cite[VI]{SGA1}),
 complemented by that of $2$-functors (or pseudo-functors),
 in order to describe the six functors formalism.
 We first describe axioms which allow one to derive
 the core formalism
  -- i.e. the part described in section \ref{num:G6FF} --
 from simpler axioms.
 We do not claim originality in this task:
  our main contribution is to give a synthesis of the approach
  of Deligne described in \cite[XVII]{SGA4}
  (see also \cite[Appendix]{Har})
  with that of Voevodsky developed by Ayoub in \cite{ayoub}.

Recall that a (cleaved) fibred category $\M$ over $\site$
 can be seen as a family of categories $\M(S)$ for
 every object $S$ of $\site$ together with a pullback
 functor $f^*:\M(S) \rightarrow \M(T)$ for any morphism
 $f:T \rightarrow S$ of $\site$.\footnote{These pullback
 functors are subject to the usual cocycle condition ; see
 section \ref{sec:pfibred}.}
 Given a suitable class $\Pmor$ of morphisms in $\site$,
 we set up a systematic study of a particular kind of
 fibred categories, called $\Pmor$-fibred categories
 (definition \ref{df:Pmor_fibred}):
 one where for any $f$ in $\Pmor$,
 the pullback functor $f^*$ admits a \emph{left} adjoint,
 generically denoted by $f_\sharp$. 
 The functor $f_\sharp$
 has to be thought as a variant
 of the \emph{exceptional direct image functor}.\footnote{This
 kind of situation frequently happens:
 the analytic case (open immersions),
 sheaves on the small \'etale site (\'etale morphisms),
 Nisnevich sheaves on the smooth site (smooth morphisms).}

\bigskip

In section \ref{sec:pfibred}, we study basic properties
 of $\Pmor$-fibred categories which will be the core of the six functors
 formalism, such as \emph{base change formulas}
 and \emph{projection formulas}
 when an additional monoidal structure is involved.
 These formulas are special cases of a compatibility relation
 between different types of functors expressed through a
 canonical comparison morphism. This kind of comparison morphisms
 are generically called \emph{exchange morphisms}.
 They are very versatile and appears everywhere in the theory
  (see Paragraphs \ref{num:exchanges1}, \ref{num:exchanges2},
  \ref{num:exchanges3}, \ref{num:exchanges4}, \ref{num:exchanges4bis},
  \ref{num:exchanges_Pfibred_morph}).
In fact, they appear fundamentally in Grothendieck six functors formalism:
 in the list of properties \ref{num:G6FF},
 they are the isomorphisms of (5), (6) and even (4).
 In the direction of the full Grothendieck functoriality,
  we introduce a core axiomatic for $\Pmor$-fibred categories 
  that we consider minimal: the categories satisfying
  this axiomatic are called \emph{$\Pmor$-premotivic}
  (section \ref{sec:premotivic_convention}).
 $\Pmor$-premotivic categories will form the basic setting in all this work.
They will appear in three different flavors,
 depending on which particular kind of additional structure
 we consider on categories: abelian, triangulated and model categories.

\bigskip

In Section \ref{sec:six_functors}, we restrict our attention
 to the triangulated and geometric case,
 meaning that we consider triangulated $\Pmor$-fibred categories 
 over a suitable category of schemes $\sch$.
 The aim of this section is to develop, and extend,
  Grothendieck six functors formalism in this basic setting.
 We exhibit many properties of such fibred categories which
  are indexed in the appendix.
Let us concentrate in this introduction on the two main properties
 which will correspond respectively to Deligne and Voevodsky's
 approach on the six functors formalism.

The first one, called the \emph{support property} and abbreviated by \supp,
 asserts that the adjoint functors
 of the kind $f_*$, for $f$ proper, and $j_\sharp$, for $j$ an open immersion,
 satisfy a gluing property that allows to use the argument of Deligne
 to construct the exceptional direct image functor $f_!$.\footnote{In
 the context of torsion \'etale sheaves of \cite[XVII]{SGA4},
  property \supp is a consequence of the proper base change theorem.}
 Several properties are derived from \supp and the basic axioms of
 $\Pmor$-fibred categories. Eventually, it leads to a partial version of the
 six functors formalism (see Theorem \ref{thm:support}).

The second property, most fundamental in the motivic context,
 is the \emph{localization property} abbreviated by \loc,
 which is in fact part of the six functors formalism
  (see Paragraph \ref{num:G6FF}).
 It has many interesting consequences 
  and reformulations that are derived in 
  Section \ref{num:general_notations_loc}.
 Note that \loc is also known in the literature
  as the ``gluing formalism''. Some properties that
  we prove in \emph{loc.cit.} are already classical (see \cite{BBD}).

The most interesting consequence of \loc was discovered
 by Voevodsky: together with the usual $\AA^1$-localization 
 and $\PP^1$-stabilization properties of the motivic context,
 it implies the complete basic six functors formalism
 as stated in Paragraph \ref{num:G6FF}.
 This was proved by Ayoub in \cite{ayoub}.
In section \ref{sec:Ayoub}, we revisit the proof of Ayoub
 and give some improvement of his theorems
 (see Theorem \ref{thm:cor3_Ayoub} for the precise
  statement):
\begin{itemize}
\item we remove the quasi-projectivity assumption
 for the existence of $f_!$, replacing it by
 the assumption that $f$ is separated of finite type;
\item we introduce the \emph{orientation property}
 which allows one to get a simpler, more usual, form
 of the purity isomorphism (the one actually stated
  in point (4) of \ref{num:G6FF});
\item we give another proof of the main theorem
 in the oriented case by showing that relative purity
 is equivalent to some (strong) duality property
 in the smooth projective case 
 (see Theorem \ref{thm:purity&duality_or});
\item we directly incorporate the monoidal structure
 whereas Ayoub gives a separate discussion for this.
\end{itemize}
Apart from these differences,
 the material of section \ref{sec:Ayoub}
 is very similar to that of \cite{ayoub}.
 Moreover, in the non oriented case,
  it should be clear that we rely 
  on the original argument of Ayoub
  for the proof of Theorem \ref{thm:purity&duality_or}.

Concerning terminology, we have called 
 \emph{motivic triangulated category}
 (Definition \ref{df:motivic_cat}) what Ayoub calls
 a ``monoidal stable homotopy functor'' (except that Ayoub only
 considers operations induced by quasi-projective morphisms).

\bigskip

The remaining of Part 1 is concerned with extensions 
 of Grothendieck six functors formalism.

 In Section \ref{sec:theorie_descente},
 we show how to use the setting of
 $\Pmor$-fibred model categories as a framework to formulate
 Deligne's cohomological descent theory.

Except in trivial cases, object of a derived category are
 not local.\footnote{The first example of this fact is the
 circle: any non-trivial connected open subset of $S^1$ is contractible
 whereas $S^1$ itself is not.}
 To formulate descent theory in derived categories,
 the main idea of Deligne was to extend the derived
 category of a scheme by one relative to a simplicial
 scheme, usually a hypercover with respect to a Grothendieck
 topology (see \cite[Vbis]{SGA4}). The construction consists
 in first extending the theory of sheaves to the case where
 the base is a simplicial scheme and then considering the
 associated derived category.

We generalize this construction to the case of an
 arbitrary $\Pmor$-fibred category equipped
 with a suitable model category structure.
%\footnote{Recall that model structures,
% introduced by Quillen, allow to perform all the usual constructions
% of derived categories by localizing an arbitrary category with
% respect to a given class of morphisms called weak equivalences.
% It contains in particular the usual case of complexes of an arbitrary
% abelian category with quasi-isomorphisms as weak equivalences.
% The main construction of the theory of Quillen is that of
% left (resp. right) derived functors which can be defined by
% replacing the usual notion of projective (resp. injective)
% resolution by that of cofibrant (resp. fibrant) resolution.}
 In fact, we show in Section \ref{sec:extension_to_diagrams}
 how to extend a $\Pmor$-fibred category over a category of schemes 
 to the corresponding category of simplicial schemes and even
 of arbitrary diagrams of schemes.
 Most importantly, we show how to extend the fibred model structure
 to the case of diagrams of schemes
  (see Prop. \ref{basicfunctdiag}).\footnote{By restricting the morphisms
 of diagrams of schemes to a
 certain class denoted by $\Pmor_\cart$, we also show how to get
 a $\Pmor_\cart$-fibred model category over diagrams of schemes
 (Rem. \ref{remPacrtfibredcat}) but this is not really needed
 in the descent theory.}
 Concretely, this means that we define a derived functor
 of the kind $\derL \varphi^*$ (resp. $\derR \varphi_*$)
 for an arbitrary morphism $\varphi$ of diagrams of schemes.
 Let us underline that these derived functors mingle two
 different kinds of functoriality: the usual pullback $f^*$
 (resp. direct image $f_*$) for a morphism of schemes $f$
 together with homotopy colimits (resp. limits) of arbitrary diagrams
 --- see the discussion in Paragraph \ref{computeadjointsdiagrams}
 till Proposition \ref{computeadjointsdiagrams3}.
With that extension in hands,
 we can easily formulate (cohomological) descent theory
 for arbitrary Grothendieck topologies on the category of schemes
 for the homotopy category of a $\Pmor$-fibred model category:
 see Definition \ref{defdescente}.

The end of Section \ref{sec:theorie_descente} is devoted 
 to concrete examples of descent
 in $\Pmor$-fibred model categories,
 and their relation with properties 
of the associated homotopy category,
 assuming it is triangulated,
 as introduced in Section \ref{sec:six_functors}.
 The first and most simple example corresponds to
 the case of a Grothendieck topology associated
 with a cd-structure in the sense of Voevodsky
 (as the Nisnevich and the $\cdh$-topology. See \cite{voecd1}
 or Paragraph \ref{df:cd_structures}).
 In that case, descent can be characterized as the existence
 of certain distinguished triangles (Mayer-Vietoris for Zariski
 topology, Brown-Gersten for Nisnevich topology): this is Theorem
 \ref{BGNis} which is in fact a reformulation of the results
 of Voevodsky.

We then proceed to the most fundamental case of descent in algebraic
 geometry, that for proper surjective maps which allows in principle 
 the use of resolution of singularities.
In fact, the main result of the whole of Section \ref{sec:theorie_descente}
 is a characterization of $\h$-descent which allows us to reduce it,
 for $\Pmor$-fibred homotopy triangulated categories which 
 are rational and motivic, to a simple property easily checked
 in practice\footnote{This is the \emph{separation} property
 defined in \ref{df:ppty:sep_ssep}. Let us mention here it is
 a consequence of the existence of well-behaved trace maps (see the proof
 of Theorem \ref{thm:DMBseparated}).}: this is Theorem
 \ref{charseparated3}. Along the way, we also proved the
 following results, interesting on their own:
\begin{itemize}
\item several characterizations of \'etale descent
 (Theorems \ref{ratdescent} and \ref{fibredetaledescent});
\item a characterization of $\qfh$-descent 
 (Theorem \ref{carratqfhdescent}) as if it was defined by
 a cd-structure.\footnote{It is at the origin of the formulation
 of descent that we gave for $\DMB$ in Theorem \ref{thmi:qfh-descent}(b)
 above. A systematic approach to such generalized cd-structures
 is developped by Park in \cite{dpark}.}
\end{itemize}
In fact, the last point is the heart of the proof of
 the main result of this section (Theorem~\ref{charseparated3}).
Whereas the extension of fibred homotopy categories to
 diagrams of schemes is not unprecedented (see \cite{ayoub2}),
 our study of proper and $\h$-descent seems to be completely new.
 In our opinion, it is one of the most important technical innovation 
 of this book.

\bigskip

In Section \ref{sec:constructible_motives},
 we study the extension of Grothendieck six functors formalism
 in \emph{rational} motivic categories, mainly duality and continuity.
 As already mentioned, the general principle is not new
 and follows mainly the path laid by Grothendieck in \cite{SGA5}.

In the case of an abstract motivic triangulated category
 --- which is for the purposes of descent theory the homotopy
 category of an underlying fibred model category ---
 the first task is to introduce a correct property of finiteness
 inherent to any duality theorem.
 This is done following Voevodsky, as in the work of Ayoub,
 by introducing the notion of \emph{constructibility} in
 Definition \ref{deftauconst}.
 The name is inspired by the \'etale case, but the notion
 of constructibility which we consider here is
 defined by a generation property which really corresponds to
 what Voevodsky called \emph{geometric motives}:
 constructible motives in our sense are generated by
 twists of motives of smooth schemes
 and are stable by cones, direct factors and finite sums.
 Let us mention that in good cases, the property of
 being constructible coincides with that of being compact
 in a triangulated category, resounding with the theory
 of perfect complexes (in the context of $\ell$-adic sheaves,
 this corresponds to ``constructible of geometric origin'').

The main point on constructible motives is the study
 of their stability under the six operations that we get from
 the axioms of a triangulated motivic category.
 This is done in Section \ref{sec:motivic_finiteness}.
 As in the \'etale case, the crucial point is the stability
 with respect to the operation $f_*$, when $f$ is a morphism
 of finite type between excellent schemes. 
 In Theorem \ref{thmfinitness}, we give conditions
 on a motivic triangulated category so that the stability
 for $f_*$ is guaranteed (then the stability by the other
 operations follows easily, see \ref{grothendieck6op}).
 Our proof essentially follows an argument
 of Gabber. The general principle, 
 going back to \cite[XIX, 5.1]{SGA4},
 is to use resolution of singularities to reduce to
 an absolute purity statement which is among our 
 assumptions.\footnote{Absolute purity will be proved later for Beilinson motives.}

In Section \ref{sec:continuity}, we introduce an important property
 of motivic triangulated categories, called \emph{continuity},
 which allows reasoning that involves projective limits of schemes.
 In fact, it is shown in Proposition \ref{continuityconstructible}
 that this property implies the property (9) of the
 (extended) Grothendieck six functors formalism
 (see Paragraph \ref{num:G6FF_continuity} above).
 We also give a criterion for continuity (\ref{abstractcontinuity})
 which will be applied later in concrete cases and draw
 some interesting consequences.

Finally, Section \ref{sec:motivic_duality} deals with duality itself
 for constructible motives, that is property (8) of 
 Paragraph \ref{num:G6FF_duality}.
 The main theorem \ref{thm:localduality} asserts that,
 under the same condition as Theorem \ref{thmfinitness},
 and if one restricts to schemes that are separated of finite type
 over an excellent base scheme $B$ of dimension less or equal to 2,
 then the full duality property holds
 (see also Corollary \ref{cor:localduality2}).
 The proof follows the same lines as the analogous Theorem 2.3.73 of
 \cite{ayoub}. In particular the main point is
 the fact that constructible motives are generated by some nice motives
 adapted to the use of resolution of singularities: see
 Corollary \ref{properregulargenerators}. The main difference
 with \emph{op. cit.} is that we implement De Jong's equivariant
 resolution of singularities~\cite{dejong2},
 so that our assumptions are much weaker.\footnote{See also
footnote \footref{footnotegeneralduality}
page \pageref{footnotegeneralduality}, which applies
to this more general setting as well.}
 
\subsection{The constructive part (Part 2)}

The purpose of this part is to give a method of construction
 of triangulated categories that satisfies the formalism
 described in Part 1.
We have chosen to mainly use the setting of derived categories.
 Also, we use our notion of $\Pmor$-fibred categories
 ($\Pmor$-premotivic with a good monoidal structure).
 Recall this means the pullback functor $f^*$
 admits a left adjoint $f_\sharp$ when $f \in \Pmor$.
 Essentially, $\Pmor$ will be either the class of smooth
 morphisms of finite type or the class of all
 morphisms of finite type (eventually separated). 

In Section \ref{sec:derived_premotivic},
 starting from a $\Pmor$-premotivic abelian category $\A$,
 we first show how to prove that the associated derived category $\Der(\A)$
 is also a $\Pmor$-premotivic category. This consists in deriving
 the structural functors of a $\Pmor$-premotivic category,
 which is done by building a suitable underlying
 $\Pmor$-fibred model category in
 Proposition \ref{t-model_category_complexes}.
 Actually, the proof of the axioms of a model category
  has already appeared in our previous work
 \cite{CD1}.
 Let us mention the flavor of this model structure:
 we can describe explicitly cofibrations as well
 as fibrations, by the use of an appropriate Grothendieck topology $t$.
 This model structure is linked with cohomological $t$-descent
 (as shown later in Proposition \ref{prop:descent&derived_P-premotivic}).
 The advantage of our framework is to easily obtain
 the functoriality of this construction
 (Paragraph \ref{num:derived_adjunction}), as well as other
 homotopical constructions (dg-structure: Rem. \ref{rem:DG-structure},
 extension to diagrams of schemes:
 Par. \ref{num:derived_premotivic_diagrams&descent}).
 In paragraph \ref{num:derived_constructible},
 we also describe in suitable cases
 the constructible objects of the derived
 category by a presentation similar to that of Voevodsky's geometric 
 motives over a perfect field.

In Section \ref{sec:A^1-derived}
 (resp. Section \ref{sec:P^1-stable-derived})
 we show how to describe the $\AA^1$-localization 
 (resp. $\PP^1$-stabilization) process in
 $\Pmor$-premotivic derived categories:
 to any $\Pmor$-premotivic abelian category $\A$
 is associated an $\AA^1$-derived category $\DMue(\A)$
 (resp. $\PP^1$-stable and $\AA^1$-derived category
  $\DMu(\A)$) in Definition \ref{df:effective_triangulated_premotives}
 (resp. \ref{df:triangulated_premotives}).
 From the model category obtained in Section \ref{sec:derived_premotivic},
  the construction uses the classical tools of motivic homotopy theory
 as introduced by Morel and Voevodsky.
 Again, our framework allows us to get the same homotopical constructions
 as in the simple derived case as well as some nice
 universal properties.
 We also get a description of constructible objects under suitable assumptions:
  Section \ref{sec:constructible_A1derived}
 (resp. \ref{sec:constructible_stableA1derived}).
 These sections are filled with concrete examples.

In Section \ref{sec:dmtilde},
 we focus on the main (in fact universal)
 example of motivic derived categories,
 the $\AA^1$-derived category of Morel,
 obtained by the process described above 
 from the abelian premotivic category of 
 abelian sheaves over the smooth Nisnevich site.
 The main point here is that one gets the localization property
 for this category by a theorem of Morel and Voevodsky.
 We give two new contributions on this topic.
 First we show in Section \ref{sec:generalized_derived}
  that the $\AA^1$-derived category can be embedded
  in a larger category which naturally contains objects
  that we can call motives of singular schemes.
 This is useful to state descent properties and will be essential
  to study $\h$-motives.
 Second, we show in Section \ref{sec:nearly_Nisnevich&loc}
  how one can use the $\AA^1$-derived category
  to obtain good properties of another premotivic derived category
  satisfying suitable assumptions.
 This will be applied to motivic complexes.

In Section \ref{sec:modules_et_anneaux},
 we go back to the case of an arbitrary monoidal $\Pmor$-fibred
 model category $\M$ and explain how to use the setting of
 ring spectra and modules over ring spectra in
 the premotivic context. The main construction
 associates to a suitable collection of (commutative) ring
 spectra $R$ in $\M$ a $\Pmor$-fibred monoidal category
 denoted by $\ho(\Mod R)$: Proposition \ref{abstractmotivicmodules}.
 This construction will be used several times:
\begin{itemize}
\item in the study of algebraic K-theory (Section \ref{sec:Ktheory}):
 the category of modules over K-theory is the fundamental technical
 tool to get motivic proper descent as well as motivic absolute purity;
\item in the study of Beilinson motives
 when we will relate them with modules over motivic cohomology
  (Theorem \ref{thm:DMB&HB-mod});
\item in the study of realizations associated with a mixed 
 Weil cohomology (Section \ref{sec:Weil}).
\end{itemize}

\subsection{Motivic complexes (Part 3)}

This part is concerned with the constructions described
 above, in Section \ref{seci:motivic_complexes}.
 Our aim is to extend the definition of Voevodsky's 
 integral motivic complexes to any base, 
 then study their functoriality
 and introduce their non-effective, or rather $\PP^1$-stable,
 counter-part.

Our first task, in Section \ref{sec:cycles},
 is to revisit Suslin-Voevodsky's theory 
 of \emph{relative cycles} exposed in \cite{SV1}.
 Indeed, they will be at the heart of the general construction.
 Our presentation is made to prepare the theory of
 \emph{finite correspondences}, a particular case of relative
 cycles. Especially, we want to give a meaning to the following
 picture representing the composition of finite correspondences
 $\alpha$ from $X$ to $Y$ and $\beta$ from $Y$ to $Z$:
$$
\xymatrix@=8pt{
\beta \otimes_Y \alpha\ar[r]\ar[d] & \beta\ar[r]\ar[d] & Z. \\
\alpha\ar[r]\ar[d] & Y & \\
X
}
$$
(see also \eqref{eq:compose_cor}). More precisely,
 we want to interpret this as a diagram of cycles.
 Thus we are led to consider cycles (with their support)
 as objects of a category. Concretely, a cycle is considered
 as a multi-pointed scheme, each point being endowed with
 some multiplicity (an integral or rational number).

This conceptual shift has the advantage of allowing a treatment
 of cycles analogous to that of algebraic varieties, or rather schemes,
 promoted by Grothendieck via studying morphisms.
Thus, we replace the various groups of relative cycles
 introduced by Suslin and Voevodsky in \emph{op. cit.}
 by properties of morphisms of cycles. Here is a list
 of the principal ones:
\begin{itemize}
\item pseudo-dominant (\ref{df:cat_cycles}),
 equidimensional (\ref{num:properties_morphism_cycles}
 and \ref{df:relative_equidim}),
\item pre-special (\ref{df:pre-special}),
\item special (\ref{df:be_special}),
\item $\Rc$-universal (\ref{df:universal_cycles}).
\end{itemize}
The most intriguing one, being \emph{pre-special},
 has no counter-part in \emph{op. cit.} Its idea comes from a mistake
 (fortunately insignificant)
 in the convention of Suslin and Voevodsky. Indeed, Lemma 3.2.4 of
 \emph{op. cit.} is false whenever the base $S$ is non reduced
 and irreducible: then any fat point $(x_0,x_1)$
 and any flat $S$-scheme give a counter-example.\footnote{Explicitly, 
 take $S=Z=\spec{k[t]/(t^2)}=\{\eta\}$, $R=\big(k[t]\big)_{(t)}$.
 The left-hand side of the equality of 3.2.4 is $2.\eta$ while
 the right-hand side is $\eta$.}
 The explanation is that the operation of specialization 
 along a fat point does not take into account 
 the geometric multiplicities of the base.
 On the contrary, when $X$ is flat over an irreducible scheme $S$,
 the geometric multiplicity of any irreducible component of $X$
 is a multiple of the geometric multiplicity of $S$.
 This leads us to the definition of a pre-special morphism
 of cycles $\beta/\alpha$,
 where a divisibility condition appears in the multiplicities 
 of $\beta$ with respect to that of $\alpha$.\footnote{To anticipate the
 rest of the construction, given a non reduced scheme $S$,
 this will allow for the operation of pull-back along the immersion
 $S_{red} \rightarrow S$ associated with the reduction of $S$:
 it simply corresponds to dividing by the geometric multiplicities of $S$,
 as the base change to $S_{red}$ does for flat $S$-schemes.}

\bigskip

The main achievement of Suslin and Voevodsky's theory 
 is the construction of a pullback operation for relative cycles.
 In our language, it corresponds to a kind of tensor product,
 more precisely a product of cycles relative to a common base
 cycle (as for example the cycle $\beta \otimes_Y \alpha$ of the
 preceding picture).
Despite our different presentation, the method to define this operation
 follows closely the original idea of Suslin and Voevodsky: 
 use the \emph{flatification theorem} of Gruson and Raynaud to
 reduce to the case of flat base change of cycles.
Recall that the key point is to find 
 the correct condition on cycles -- or rather morphisms of cycles
 in our language -- so that one obtains a uniquely defined operation
 independent of the chosen flatification. This is measured by
 a specialization procedure (Definition \ref{df:foncteur_sp})
 associated with \emph{fat points}  (Definition \ref{df:fat_points})
 and leads to the central notion of \emph{special morphisms} of cycles
 (Definition \ref{df:be_special}).
 An innovation that we introduce in the theory is to give,
 as soon as possible, local definitions at a point in the style of
 EGA. This is in particular the case for the property of being special.

Once this notion is in place, one defines for a base cycle $\alpha$,
 a special $\alpha$-cycle $\beta$ and any morphism 
 $\phi:\alpha' \rightarrow \alpha$ the relative product denoted by
 $\beta \otimes_\alpha \alpha'$. Equivalently, it
 corresponds to the base change of $\beta/\alpha$ along $\phi$
 (Definition \ref{df:tensor_product_cycles}).
 This notion is close to the correspondence morphisms
 of Section 3.2 of \emph{op. cit.}
 In particular, it usually involves denominators.
 The last important notion, being \emph{$\Rc$-universal}, corresponds
 to cycles $\beta/\alpha$ with coefficients in a ring $\Rc \subset \QQ$,
 which keeps their coefficients in $\Rc$ after any base change.

One sees that our language is especially convenient when it is time to
 consider the stability of certain properties of morphisms of cycles
 under composition (Cor. \ref{transitivity2}) or base change
  (Cor. \ref{cor:special_base_change}).
 Then the usual statements of intersection theory are proven
  in Section \ref{sec:intersection_theory},
  still following or extending Suslin and Voevodsky:
 commutativity, associativity, projection formulas.
 This makes our relative product a good extension of
 the classical notion of exterior product of cycles (over a field).

The focal point of intersection theory is the study of multiplicities.
 Thus we introduce \emph{Suslin-Voevodsky's multiplicities},
 as the ones appearing as a corollary of the existence of
 the relative cycle $\beta \otimes_\alpha \alpha'$ 
 (Definition \ref{SV_multiplicities}).
 A very important result in the theory, already enlightened by
 Suslin and Voevodsky, is the fact these multiplicities can be
 expressed in terms of \emph{Samuel multiplicities}.\footnote{When
 a correct regularity assumption is added,
  one reduces to the usual Serre's Tor-intersection formula:
  see \ref{th:mSV&tor} and \ref{rem:intersection&Serre}).}
 In fact, independently of Suslin and Voevodsky,
 we prove a new criterion for the property of
 being special at a point involving Samuel multiplicities
 at the branches of the point: see Corollary \ref{special&Samuel}.
 Roughly speaking, the multiplicities arising from Samuel's definition
 at each branches of the point must coincide: then this common
 value is simply the Suslin-Voevodsky's multiplicity.

Finally, still following the treatment of algebraic geometry by Grothendieck, 
 we add to the theory of Suslin and Voevodsky the study 
 of constructibility properties for morphisms of cycles (special and $\Rc$-universal).
 Here, our categorical point of view is plainly justified.
 Explicitly, we prove that given a relative cycle $\beta/\alpha$,
 when $\alpha$ is the cycle associated with a scheme $S$,
 the locus where $\beta$ is special (resp. $\Rc$-universal)
 is an ind-constructible subset of $S$ (Lemma \ref{ind-constructible}).
 This allows to prove the good behavior of these notions
 with respect to projective limits of schemes
  (see in particular \ref{prop:lifting_cycles}).
 This will be the {\it key point} when proving the continuity 
 property --- (9) of \ref{num:G6FF_continuity} ---
 of the fibred category $\DMV$.

\bigskip

The rest of Part 3, consists in extending the theory of
 sheaves with transfers introduced by Voevodsky,
 originally over a perfect field,
 to the case of an arbitrary base
 and apply to it the general procedures studied in Part 2
 to get the fibred category $\DMV$.

In Section \ref{sec:finite_cor}, we work out
 the theory of finite correspondences using the
 formalism of relative cycles.
 The construction is summarized 
  in Corollary \ref{prop:pfibred_corr}: given
  a class of morphisms $\Pmor$ contained in
  the class of separated morphisms of finite type
  and a ring of coefficients $\Rc$,
  we produce a monoidal $\Pmor$-fibred category,
  denoted by $\Pmorc$, whose fiber over a noetherian scheme $S$
 (eventually singular) is the category of $\Pmor$-schemes
  over $S$ with morphisms the finite correspondences.

In Section \ref{sec:fxtr}, we develop the theory
 of sheaves with transfers along the very same line as
 the original treatment of Voevodsky.
 This time, the outcome can be summarized by Corollaries
 \ref{cor:Ppremotivic_ftr} and \ref{cor:compatibility_tr&top}:
 given a class $\Pmor$ of morphisms as above
 and a suitable Grothendieck topology $t$,
 we construct an abelian premotivic category
 $\sh t {\Pmor}$ which is compatible with the topology
 $t$ (cf Part 2); its fiber over a scheme $S$
 is given by $t$-sheaves of $\Rc$-modules with transfers
  (in particular presheaves on $\Pmorcx S$).\footnote{The 
  most notable topologies $t$ that fit in this result
  are the Nisnevich and the $\cdh$ ones. 
 See Section \ref{sec:examples_ftr}.}
 The section is closed with an important comparison result,
  essentially due to Voevodsky, between Nisnevich sheaves with
 transfers and sheaves for the $\qfh$-topology 
  (with rational coefficients over geometrically unibranch bases):
  see Theorem \ref{enlargeabsheavestrqfhsheaves}.

Finally, Section \ref{sec:voevodsky} is devoted to
 gathering the work done previously and define the stable derived 
 category of motivic complexes $\DM$, given an arbitrary
 ring of coefficients $\Rc$. The out-come has already been
 described in Section \ref{seci:motivic_complexes} above.

\subsection{Beilinson motives (Part 4)}

This part contains the construction of Beilinson motives
 as well as the proof of all the properties stated before.
 It is based on the first and second parts but independent
 of the third one --- except in the comparison statements
 of Section \ref{sec:compare_BV}.

Section \ref{sec:SH_recall} contains a short review
 of the stable homotopy category and the notion of oriented
 ring spectra.

Section \ref{sec:Ktheory} is the heart of our construction.
 It contains a detailed study of the K-theory ring spectrum
 $\BGL$  and the associated notion of $\BGL$-modules in
 the homotopical sense
 (based on the formalism introduced in Section
  \ref{sec:modules_et_anneaux}).
 Using the works of several authors
  (most notably: Riou, Naumann, Spitzweck, \O stv\ae r),
 we show how the central results of Quillen on algebraic K-theory
 give important properties of $\BGL$-modules: absolute purity
 (Th. \ref{abspurityBGL}) and trace maps (Def. \ref{df:trace_BGL}).

In Section \ref{sec:Beilinson_motives},
 we finally introduce the definition of Beilinson motives.
 Let us describe it in detail now. It is based on the process
 of Bousfield localization of the stable homotopy category
 with respect to a cohomology. This operation is fundamental
 in modern algebraic topology. We apply it in algebraic geometry
 to the rational stable homotopy category
 (or, what amount to the same, to the rational stable 
 $\AA^1$-derived category of Morel, Section \ref{sec:dmtilde})
 and to the rational K-theory spectrum $\BGL_\QQ$:
 the Bousfield localization of $\Der_{\AA^1}(S,\QQ)$ with respect
 to $\BGLrx S$ is the category of Beilinson motives $\DMB(S)$ 
 over $S$ (Definition \ref{df:Beilinson_motives}).
 Using the preceding study of $\BGL_\QQ$ together with
 the decomposition of Riou recalled in the beginning of Section
 \ref{seci:comparison}, we get the main properties of the
 premotivic category $\DMB$: the $\h$-descent theorem 
 (\ref{DMB_etale&h-descent}) and 
 the absolute purity theorem (\ref{DMBpurity}).

Then the theoretical background laid down in Part 1
 is applied to $\DMB$, given in particular the complete
 Grothendieck six functors formalism for constructible
 Beilinson motives (Section \ref{sec:constructible_Beilinson}).
Our work closed with the two main subjects described
 above on Beilinson motives:
 the comparison statements (Section \ref{sec:comparisons_Beilison})
 and the study of motivic realizations (Section \ref{sec:Weil}).

\section{Developments since the first arXiv version}

The first version of this work has first appeared on arXiv
on December 2009.\footnote{A new version
 was uploaded in 2012, containing more or less the actual introduction
 which was written in order to clear-up the contributions
 and history on mixed
 motives and more specifically motivic homotopy theory.}
 During almost ten years, until the actual publication by Springer Edition,
 it has been used in several works, as well as
 completed by several other mathematicians,
 solving questions left open in the present text.
 For completeness, it appears to us beneficial
 to the reader to give an account of some of these developments
 which are the most directly related with the present contribution.
 Mathematics is indeed a collective work,
 each part of which is destined to be used, completed, renewed or superseded.

\subsection{Nisnevich motives with integral coefficients}

\begin{num}\textit{$\cdh$-motives}.--
One aim of the present work was to work out the theory of finite correspondences
 in the spirit of \cite{FSV}, whose original aim is to obtain an integral theory
 of motivic complexes related to Chow groups. The theory of $\cdh$-sheaves with transfers (see Proposition \ref{prop:cdh&transfers})
 was introduced with this motivation in mind.
 The theory of $\cdh$-motives and motivic complexes was successfully
 developed in the equal characteristic case in \cite{CD5}, provided one inverts the
 residue characteristic. In this latter work, the crucial property of localization for
 $\cdh$-motives is shown, as well as all the expected results: constructibility of
 the six operations, duality, continuity,
 comparison with modules over the $\cdh$-local version of Voevodsky's motivic cohomology,
 relation with higher Chow groups. To get these results, key points
 are the continuity property of motivic complexes which is proved in this
 book (Theorem \ref{thm:DM_continuity}),
 as a result of our reinforcement of Suslin-Voevodsky's
 theory of relative cycles
 (see in particular Section \ref{sec:constructibility} on
 constructibility for properties of relative cycles)
 together with Kelly's new  motivic descent
 results \cite{integraltransfers} which allow to use
 Gabber's improvements of de~Jong's alteration theorems \cite{gabber3}.
 Note also that \cite[3.6 and 5.1]{CD5}
 generalizes our result on Voevodsky's conjecture on base change of
 the motivic Eilenberg-MacLane spectrum (Corollary \ref{compHbeilHVoev}).
\end{num}

\begin{num}\textit{Spitzweck's motivic cohomology spectrum}.--
One of the problems with defining mixed motives
as modules over Voevodsky's
motivic cohomology
 spectrum is the compatibility of this spectrum with base change.\footnote{Recall again
 this compatibility was conjectured by Voevodsky. See Conjecture \ref{conj:voevodsky}
 for an explicit formulation. Note also that we prove the latter conjecture
 is actually equivalent to the localization property for (Nisnevich) motivic complexes:
 see Proposition \ref{prop:reformulation of the conjecture of voevodsky}.}
 The idea of Spitzweck's paper \cite{motivicHZ} is to build a spectrum which satisfies
 compatibility by base change; equivalently, one has to build a ring spectrum $\HSpit$
 over $S=\spec{\ZZ}$
 (or more generally over a Dedekind ring) which pullbacks to Voevodsky's motivic cohomology
 spectrum over the residue fields of $S$. This is what
 M. Spitzweck achieves with virtuosity in
 \emph{loc. cit.}, therefore obtaining a convenient category of $\HSpit$-modules which 
 coincides with Voevodsky's original triangulated category over the residue fields of $S$; in fact, it also coincides with $\DMcdh$ over any $k$-scheme, after inverting the residue
 characteristic of $k$. But the construction of Spitzweck works integrally. Moreover, by its very construction, the cohomology represented by $\HSpit$ coincides with
 Bloch's higher Chow groups for smooth $S$-schemes. A question left open is a possible
 comparison with Voevodsky's motivic cohomology spectrum (which is again equivalent
 to Voevodsky's base change conjecture).
\end{num}

\subsection{\'Etale motives with integral coefficients
and $\ell$-adic realization}

\begin{num}\textit{Voevodsky's motives in the \'etale topology and rigidity theorems}.--
 With rational coefficients, the comparison theorems obtained in this book
 (Section \ref{sec:comparisons_Beilison})
 show that varying the underlying topology
 is beneficial. In particular, with rational coefficients,
 we are not able to get the localization property
 for Nisnevich motivic complexes for all base schemes, but we do get that property
 when replacing Nisnevich topology with the $\qfh$-topology, or the $\h$-topology.
 This lets one believe that the transfers will be better behaved
 with integral coefficients with respect
 to stronger topologies. In \cite{CD4}, we do prove that the localization property holds for motivic complexes with torsion coefficients locally for the \'etale
 topology (\cite[Theorem 4.3.1]{CD4}). It follows that the same
 property holds 
 with integral coefficients for geometrically unibranch schemes.
 Moreover,  we prove in \emph{loc. cit.} that,
 locally for the $\h$-topology,
 motivic complexes with integral coefficients
 are perfectly well-behaved and satisfy all the
 expected properties, as listed in Section \ref{sec:Gsixfunctors_intro} of this
 introduction.\footnote{Based on the results of this book, we only get Grothendieck-Verdier
 duality for schemes of finite type over a regular $2$-dimensional excellent scheme, but
 this extra hypothesis has been removed
 by the first author in \cite{Cis5}.}
 It is remarkable that we were able to get the complete Grothendieck six
 functors formalism for Voevodsky's original construction of \'etale motives,
 as defined in his Ph. D. thesis \cite{V1},
 and show the visionary power of Voevodsky once more time. 
 Besides, we also show that one recovers the theory of $\ell$-adic complexes
 out of $\h$-motives by the categorical process of $\ell$-adic completion
 (see \cite[\textsection 7.2]{CD4}). This gives a new insight on $\ell$-adic realization of motives.
\end{num}

\begin{num}\textit{Rigidity theorems without transfers}.--
 Another extension of Suslin and Voevodsky's rigidity theorem to
 arbitrary bases is due to Ayoub, \cite{ayoub5}.
 In this latter work, Ayoub studies the category introduced in the present book
 under the notation $\DMux{\et}(S,\Lambda)$ (following Morel),
 while he uses the notation $\mathrm{DA}(S,\Lambda)$.
 The first result of \emph{loc. cit.}, inspired by an earlier
 work of R\"ondigs and {\O}stv{\ae}r~\cite{RO},
 is indeed a variation on the rigidity theorem,
 identifying the category $\DMux{\et}(S,\Lambda)$ for a $\Lambda=\ZZ/N\ZZ$
 with $N$ invertible on $S$ with the derived category of the category of sheaves
 of $\Lambda$-modules on the small \'etale site of $S$
 (under suitable hypothesises on $S$ and $N$).
 From there, one can extend the results proved in this book
 for $\DMux{\et}(S,\QQ)$
 (see Theorem \ref{etalemotives}) to the case of arbitrary coefficients:
 absolute purity, constructibility of the six operations, duality.
 Note however that, despite what is claimed
 in the appendix of Ayoub's article,
 the particular case of $2$-torsion
 for base schemes $S$ of mixed or positive characteristic
 is problematic in his approach (see \cite[Rem. 5.5.8]{CD4}).
 
Since then, Bachmann \cite{bachmann2} has extended by far the preceding
 rigidity theorems to torsion $\PP^1$-stable motivic \'etale sheaves of \emph{spectra}.
 This result also solves the aforementioned issues
 left open in \cite{ayoub5}.
\end{num}

\subsection{Motivic stable homotopy theory with rational coefficients}
\renewcommand{\Rc}{\QQ}
\begin{num}\textit{Witt sheaves}.--
In this book, following Morel's insights, we have splitted
the rational motivic stable homotopy category
$\SH(X)_\QQ\simeq\DMtx X$ into two factors
$$(\DMtx X_{+})\times(\DMtx X_{-})\simeq\DMtx X$$
and we have identified the oriented part $\DMtx X_{+}$
with Beilinson's motives $\DMB(X)$. On the other hand,
in the case where $X=\spec k$ is the spectrum of a field,
Ananyevskiy, Levine and Panin~\cite{AnLePa} have identified
the non-oriented part $\DMtx k_{-}$ with a suitable
category of Witt sheaves. The conjunction of their results
with ours may be seen as a motivic analog of (a trivial consequence of) a
theorem of Serre that the stable homotopy groups of spheres
are finite in degree $>0$;
see the introduction of \emph{loc. cit.}
The results of Ananyevskiy, Levine and Panin
have been improved by Bachmann~\cite{bachmann}, where
the comparison of $\SH(k)_{-}$ with Witt sheaves
is promoted to $\ZZ[1/2]$-linear coefficients.
Bachmann's results follow from a nice analog of the rigidity
theorem over a general base for the \emph{real {\'e}tale} topology.
\end{num}

\begin{num}\textit{Rational absolute purity}.--
D{\'e}glise, Fasel, Jin and Khan~\cite{DFJK}
have proved absolute purity property for the motivic sphere spectrum
with rational coefficients.
\end{num}
\renewcommand{\Rc}{\Lambda}

\subsection{Duality, weights and traces}

\begin{num}\textit{Weight complexes}.--
Bondarko's theory of weight complexes~\cite{Bondarko0}
has been showed to
be compatible with the six operations with
rational coefficients in \cite{Hebert,Bondarko}.
In the setting of $\cdh$-sheaves~\cite{CD5}, this has been
extended by Bondarko and Ivanov~\cite{cdhweights}
to $\ZZ[1/p]$-linear coefficients in equal characteristic,
where $p\geq 1$ denotes the exponent characteristic of the ground
field. Such weight complexes have been used by
Wildeshaus~\cite{Wildes1,Wildes2},
in order to give inconditional constructions of
motivic intersection complexes of certain Shimura varieties.
They also play a role, together with realization functors associated
to mixed Weil cohomologies, in geometric representation theory,
in the work of Soergel and his collaborators~\cite{Soergel,Soergel2}.
\end{num}

\begin{num}\textit{Motivic Lefschetz-Verdier trace formula}.--
An obvious application of the theory of motivic sheaves and
their realizations is the proof
of independence of $\ell$ results for a wealth of trace-like
constructions.
A $\QQ$-linear version of such kind of results
is provided by Olsson~\cite{Olsson1,Olsson2},
where some versions of the Motivic Lefschetz-Verdier trace formula
and of characteristic classes are discussed.
A slight improvement, allowing torsion, may be found
in \cite{Cis5}. But a full account on integral
formulas, including for characteristic classes, is settled
in the recent work of Jin and Yang \cite{jinyang}.
\end{num}

\subsection{Enriched realizations}

\begin{num}\textit{Structured mixed Weil cohomologies}.--
In his thesis~\cite{drew1}, Drew extends the formalism
of mixed Weil cohomologies to cohomologies with
values in a Tannakian category. He also
defines the realization functor into algebraic
$\mathcal{D}$-modules for schemes of finite type over a field of
characteristic zero and proves that, for any separated smooth scheme $X$
over a field of characteristic
zero, constructible modules over de~Rham cohomology in $\SH(X)$
embedd fully faithfully in algebraic
$\mathcal{D}$-modules. Drew deduces from this embedding a new purely algebraic
proof of the Riemann-Hilbert correspondence, using
motivic sheaves, as predicted in Example \ref{ex:deRhamanalytic}
in the present book. His work is also a way to define Hodge
realizations of mixed motivic sheaves; see \cite{drew2}.
\end{num}

\begin{num}\textit{Arakelov motivic cohomology}.--
Holmstrom and Scholbach~\cite{scholbach1} have extended
the representability of algebraic de~Rham cohomology to
the filtered de~Rham complex, and used it to define a motivic
version of Arakelov cohomology. The relation with
more classical versions of Arakelov cohomology and with
height pairings is discussed in  \cite{scholbach2}.
\end{num}

\section*{Thanks}
\addcontentsline{toc}{section}{Thanks}

The authors want to thank heartily the following people for help, motivation, corrections or suggestions during the elaboration of this text: Joseph Ayoub, Alexander Beilinson, Pierre Deligne, Brad Drew, David H{\'e}bert, Jens Hornbostel, Annette Huber-Klawitter, Bruno Kahn, Shane Kelly, Marc Levine, Georges Maltsiniotis, Fabien Morel, Paul Arne \O stv{\ae}r, Jo\"el Riou,
Oliver R\"ondigs, Valentina Sala, Markus Spitzweck, Vladimir Voevodsky,
Vadim Vologodsky, Chuck Weibel, and finally, as well as particularly,
J\"org Wildeshaus.\\
We would like to thank our editors of Springer
for simply making the publication of this book possible. \\
D.-C.~C is partially supported by the SFB 1085 ``Higher Invariants''
funded by the Deutsche Forschungsgemeinschaft (DFG), and F.~D. by
the French ``Investissements d'Avenir'' program, project ISITE-BFC
(contract ANR-lS-IDEX-OOOB).

\section*{Notations and conventions}
\addcontentsline{toc}{section}{Notations and conventions}

In every section, we will fix a category denoted by $\sch$
 which will contain our geometric objects.
 Most of the time, $\sch$ will be a category of schemes
which are suitable for our needs;
the required hypothesis on $\sch$ are given at the head of each section.
In the text, when no details are given, any scheme will be
assumed to be an object of $\sch$.

\bigskip

When $\A$ is an additive category,
we denote by $\A^\natural$ the pseudo-abelian envelope of $\A$.
We denote by $\Comp(\A)$
the category of complexes of $\A$. We consider $\K(\A)$
(resp. $\K^b(\A)$) the category of complexes (resp. bounded complexes) 
of $\A$ modulo the chain homotopy equivalences
and when $\A$ is abelian, we let $\Der(\A)$
be the derived category of $\A$.

If $\M$ is a model category, $\ho(\M)$ will denote its
homotopy category.

%%Unless otherwise stated, the $2$-functors considered here
%%as well as the morphisms of $2$-functors are strict.
We will use the notation
$$
\alpha:\C \rightleftarrows \D:\beta
$$
to mean a pair of functors such that $\alpha$ is left adjoint to
$\beta$. Similarly, when we speak of an adjoint pair of functors
$(\alpha,\beta)$, $\alpha$ will always be the left adjoint.
We will denote by
$$
ad(\alpha,\beta):1 \rightarrow \beta \alpha
\text{ (resp. } 
ad'(\alpha,\beta):\alpha \beta \rightarrow 1
\text{)}
$$
the unit (resp. counit) of the adjunction $(\alpha,\beta)$.
%We also refer to the fact the composition of natural transformations
%\begin{align*}
%\alpha \xrightarrow{\alpha.ad(\alpha,\beta)} \alpha \beta \alpha
% \xrightarrow{ad'(\alpha,\beta).\alpha} \alpha \\
%\beta \xrightarrow{ad(\alpha,\beta).\beta} \beta \alpha \beta
% \xrightarrow{\beta.ad'(\alpha,\beta)} \beta
%\end{align*}
%are the identity as the \emph{zig-zag equations}.
Considering a natural tranformation $\eta:F \rightarrow G$
of functors, we usually denote by the same letter $\eta$
--- when the context is clear --- the induced natural transformation $AFB \rightarrow AGB$
obtained when considering functors $A$ and $B$ composed
on the left and right  with $F$ and $G$ respectively.

In section \ref{sec:cycles},
 we will assume that equidimensional morphisms
 have constant relative dimension.

\renewcommand{\thesection}{\arabic{section}}
\setcounter{section}{0}

\mainmatter

\part{Fibred categories and the six functors formalism}
\markboth{Fibred categories and the six functors formalism}{}
\section{General definitions and axiomatic} \label{sec:pfibred}

\begin{assumption}\label{defclassP}
We assume that $\site$ is an arbitrary category.

We shall say that a class $\Pmor$ of morphisms of $\site$ is
 \emph{admissible}\index{word}{admissible, class of morphisms}
if it is has the following properties.
\begin{itemize}
\item[(Pa)] Any isomorphism is in $\Pmor$.
\item[(Pb)] The class $\Pmor$ is stable by composition.
\item[(Pc)] The class $\Pmor$ is stable by pullbacks: for any
morphism $f:X\rightarrow Y$ in $\Pmor$ and any morphism $Y'\rightarrow Y$,
the pullback $X'=Y'\times_Y X$ is representable in $\site$, and
the projection $f':X'\rightarrow Y'$ is in $\Pmor$.
%%\item[(d)] The class $\Pmor$ is stable by sums.
\end{itemize}
The morphisms which are in $\Pmor$ will be called
 the \emph{$\Pmor$-morphisms}.\footnote{In practice, 
 $\site$ will be an adequate subcategory of the category of
 noetherian schemes and $\Pmor$ will be the class of smooth
 morphisms (resp. \'etale morphisms, morphisms of finite type, 
 separated or not necessarily separated) in $\site$.}

In what follows, we assume
that an admissible class of morphisms $\Pmor$ is fixed.
\end{assumption}

\subsection{$\Pmor$-fibred categories}

\subsubsection{Definitions}

Let $\cat$ be the $2$-category of categories.

\begin{num}
Let $\M$ be a fibred\index{word}{fibred!fibred category} category
 over $\site$, seen
as a $2$-functor $\M:\op{\site}\to\cat$; see \cite[Exp.~VI]{SGA1}

Given a morphism $f:T\rightarrow S$ in $\site$, we shall
denote by
$$f^*:\M(S)\rightarrow \M(T)$$
the corresponding pullback functor between the
corresponding fibers. We shall always assume that
$(1_S)^*=1_{\M(S)}$,
 and that for any morphisms $W \xrightarrow g T \xrightarrow f S$ in $\site$,
 we have structural isomorphisms:
\begin{equation} \label{eq:f^*&composition}
g^* f^* \xrightarrow{\sim} (fg)^*
\end{equation}
which are subject to the usual cocycle condition 
 with respect to composition of morphisms.

Given a morphism $f:T\rightarrow S$ in $\site$, if the corresponding inverse image
functor $f^*$ has a left adjoint, we shall denote it by
$$f_\sharp:\M(T)\rightarrow\M(S)\, .$$
For any morphisms $W \xrightarrow g T \xrightarrow f S$ in $\site$
such that $f^*$ and $g^*$ have a left adjoint,
we have an isomorphism obtained by transposition
from the isomorphism \eqref{eq:f^*&composition}:
\begin{equation} \label{eq:f_sharp&composition}
 (fg)_\sharp\xrightarrow{\sim}f_\sharp g_\sharp \, .
\end{equation}
\end{num}

\begin{df}
A \emph{pre-$\Pmor$-fibred category}
\index{word}{fibred!pre-$\Pmor$-fibred category}
 $\M$ over $\site$ is
a fibred category $\M$ over $\site$ such that, for any morphism
$p:T\rightarrow S$ in $\Pmor$, the pullback functor $p^*:\M(S)\rightarrow\M(T)$
has a left adjoint $p_\sharp:\M(T) \rightarrow \M(S)$.
\end{df}

\begin{conv} \label{convention_identification1}
Usually,
 we will consider that \eqref{eq:f^*&composition}
 and \eqref{eq:f_sharp&composition} are identities.
 Similarly, we consider that for any object $S$ of $\site$,
 $(1_S)^*=1_{\M(S)}$ and $(1_S)_\sharp=1_{\M(S)}$.\footnote{
We can always strictify globally the fibred
category structure so that $g^* f^*= (fg)^*$
for any composable morphisms $f$ and $g$,
and so that $(1_S)^*=1_{\M(S)}$ for any
object $S$ of $\site$; moreover, for a morphism
$h$ of $\site$ such that a left adjoint of $h^*$ exists,
and we can choose the left adjoint functor $h_\sharp$ which we feel as the most
convenient for us, depending on the situation we deal with.
For instance, if $h=1_S$, we can choose $h_\sharp$ to be $1_{\M(S)}$, and
if $h=fg$, with $f^*$ and $g^*$ having left adjoints,
we can choose $h_\sharp$ to be $f_\sharp g_\sharp$ (with the
unit and counit naturally induced by composition).}
\end{conv}

\begin{ex} \label{ex:can_weak_P-fibred}
Let $S$ be an object of $\site$.
We let $\Pmorx S$ be the full subcategory of the comma category $\site/S$
 made of objects over $S$ whose structural morphism is in $\Pmor$.
 We will usually call the objects of $\Pmorx S$ the \emph{$\Pmor$-objects over $S$}.
 
Given a morphism $f:T \rightarrow S$ in $\site$
 and a $\Pmor$-morphism $\pi:X \rightarrow S$,
 we put $f^*(\pi)=\pi \times_S T$ using the property (Pc) of $\Pmor$
  (see Paragraph \ref{defclassP}).
This defines a functor $f^*:\Pmorx S \rightarrow \Pmorx T$.

Given two $\Pmor$-morphisms $f:T \rightarrow S$
 and $\pi:Y \rightarrow T$, we put 
 $f_\sharp(\pi)=f \circ \pi$ using the property (Pb) of $\Pmor$.
 this defines a functor $f_\sharp:\Pmorx T \rightarrow \Pmorx S$.
According to the property of pullbacks,
 $f_\sharp$ is left adjoint to $f^*$.

We thus get a pre-$\Pmor$-fibred category
 $\Pmorx ?:S \mapsto \Pmorx S$.
\end{ex}

\begin{ex} \label{ex:H_weak_P-fibred}
Assume $\sch$ is the category of noetherian schemes of finite dimension,
 and $\Pmor=\sm$.
For a scheme $S$ of $\sch$,
 let  $\Hpt(S)$
\index{notat}{HpS@$\Hpt(S)$}
 be the pointed homotopy category of schemes 
\index{word}{homotopy category}
over $S$ defined by Morel and Voevodsky
 in \cite{MV}.
Then according to \emph{op. cit.}, $\Hpt$ is a pre-$\sm$-fibred category
 over $\sch$.
\end{ex}

\begin{num} \label{num:exchanges1}
\textit{Exchange structures I}.--
Suppose given a pre-$\Pmor$-fibred category $\M$.

Consider a commutative square of $\site$
$$
\xymatrix{
Y\ar^q[r]\ar_g[d]\ar@{}|{\Delta}[rd] & X\ar^f[d] \\
T\ar_p[r] & S
}
$$
such that $p$ hence $q$ are $\Pmor$-morphisms,
 we get using the identification of convention \ref{convention_identification1}
 a canonical natural transformation
$$
Ex(\Delta_\sharp^*):q_\sharp g^*
 \xrightarrow{ad(p_\sharp,p^*)} q_\sharp g^* p^* p_\sharp
 = q_\sharp q^* f^* p_\sharp
 \xrightarrow{ad'(q_\sharp,q^*)} f^* p_\sharp
$$
called the \emph{exchange transformation}
\index{word}{exchange!transformation}
 between $q_\sharp$ and $g^*$.
\end{num}

\begin{rem} \label{rem:coherence_exchange1}
These exchange transformations satisfy a coherence\index{word}{coherence}
 condition
 with respect to the relations $(fg)^*=g^*f^*$ 
 and $(fg)_\sharp=f_\sharp g_\sharp$.
As an example, consider two commutative squares in $\site$:
$$
\xymatrix{
Z\ar^{q'}[r]\ar_h[d]\ar@{}|{\Theta}[rd]
 & Y\ar^q[r]\ar_g[d]\ar@{}|{\Delta}[rd] & X\ar^f[d] \\
W\ar_{p'}[r] & T\ar_p[r] & S
}$$
and let $\Delta \circ \Theta$ be the commutative square
 made by the exterior maps --- it is usually called the horizontal
 composition of the squares.
Then, the following diagram of $2$-morphisms is commutative:
$$
\xymatrix@C=40pt@R=20pt{
(q q')_\sharp h^*\ar^{Ex(\Delta \circ \Theta)_\sharp^*}[rr]\ar@{=}[d]
 && f^*(pp')_\sharp\ar@{=}[d] \\
q_\sharp q'_\sharp h^*\ar^-{Ex(\Theta_\sharp^*)}[r]
 & q_\sharp g^* p'_\sharp\ar^-{Ex(\Delta_\sharp^*)}[r]
 & f^* p_\sharp p'_\sharp
}$$
To see this, one proceeds as follows.
First, we observe that, since $ad'_q$ is a natural transformation,
for each object $M$ of $\M(T)$, the square
$$\xymatrix{
q'_\sharp q^{\prime *} g^*(M)\ar[rr]^-{q'_\sharp q^{\prime *}g^*(ad_q(M))}
\ar[d]_{ad'_q(g^*(M))}&&q'_\sharp q^{\prime *} g^*p^*p_\sharp(M)\ar[d]^{ad'_q(g^*p^*p_\sharp(M))}\\
g^*(M)\ar[rr]_{g^*(ad_q(M))}&&g^*p^*p_\sharp(M)
}$$
commutes. In other words, with a slight abuse of notations, we have the following
commutative square of functors.
$$\xymatrix{
q'_\sharp q^{\prime *} g^*\ar[r]^{ad_q}
\ar[d]_{ad'_q}&q'_\sharp q^{\prime *} g^*p^*p_\sharp\ar[d]^{ad'_q}\\
g^*(M)\ar[r]_{ad_q}&g^*p^*p_\sharp
}$$
We then consider the diagram below,
in which $ad_r$ (resp. $ad'_r$) indicates the morphism obtained 
from the obvious unit morphism (resp. counit morphism)
of the adjunction $(r_\sharp,r^*)$ by eventually adding functors
on the left side or on the right side, and we can check easily
that each cell below is commutative, proving our claim.
\begin{center}
\rotatebox{90}{
$
\xymatrix@C=10pt@R=60pt{
(qq')_\sharp h^*\ar@{=}[ddd]\ar^-{ad_{(pp')}}[rr]\ar@{}|{}[rddd]
 && (qq')_\sharp h^* (pp')_\sharp (pp')^*\ar@{=}[d]\ar@{=}[rr]
      \ar@{}|{}[rrd]
 && (qq')_\sharp (qq')^* f^* (pp')^*\ar@{=}[d]\ar^-{ad'_{(qq')}}[rr]
 && f^* (pp')_\sharp\ar@{=}[ddd]\ar@{}|{}[lddd] \\
 & q_\sharp q'_\sharp h^*p^{\prime *}p'_\sharp\ar_-{ad_{p}}[r]\ar@{=}[dd]
 & q_\sharp q'_\sharp h^*p^{\prime *}p^*p_\sharp p'_\sharp
    \ar@{=}[dr]\ar@{}|{}[rrd]
 & 
 & q_\sharp q'_\sharp q^{\prime *}q^*f^*p_\sharp p'_\sharp
    \ar_-{ad'_{q'}}[r]\ar@{=}[dl]
 & q_\sharp q^*f^*p_\sharp p'_\sharp\ar@{=}[dd]
 & \\
 &
 & q_\sharp q'_\sharp q^{\prime *}g^*p'_\sharp
    \ar_-{ad_{p}}[r]\ar@{=}[d]\ar@{}|{}[rrd]
 & q_\sharp q'_\sharp q^{\prime *}g^*p^*p_\sharp p'_\sharp\ar_-{ad'_{q'}}[r]
 & q_\sharp g^*p^*p_\sharp p'_\sharp\ar@{=}[d]
 & & \\
q_\sharp q'_\sharp h^*\ar_-{ad_{p'}}[r]
 & q_\sharp q'_\sharp h^*p^{\prime *}p'_\sharp\ar@{=}[r]
 & q_\sharp q'_\sharp q^{\prime *}g^*p'_\sharp\ar_-{ad'_{q'}}[r]
 & q_\sharp g^*p'_\sharp\ar_-{ad_{p}}[r]
 & q_\sharp g^*p^*p_\sharp p'_\sharp\ar@{=}[r]
 & q_\sharp q^*f^*p_\sharp p'_\sharp\ar_-{ad'_{q}}[r]
 & f^*p_\sharp p'_\sharp
}
$
}
\end{center}

Therefore, according to our abuse of notation for natural transformations,
 $Ex$ behaves as a contravariant functor with respect to the horizontal
 composition of squares. The same is true for vertical composition
 of commutative squares. 
\end{rem}

\begin{rem}
In the sequel, we will introduce several exchange transformation
 between various functor. We speak of an \emph{exchange isomorphism} when
 the transformation is an isomorphism. 
\index{word}{exchange!isomorphism}
 When only two kind of functors
 are involved, say of type a and b, we say that functors of type a
  and functors of type b commute when the exchange transformation 
  is an isomorphism.
\index{word}{functor!commutes}
\index{word}{commute|see{functor}}

 As an example (see also next definition),
  when the exchange transformation $Ex(\Delta_\sharp^*)$
 is an isomorphism,
  we simply say that $f^*$ and $p_\sharp$ commute
  --- or also that $f^*$ commutes with $p_\sharp$.
\end{rem}

\begin{num} \label{defPBC}
Under the setting of \ref{num:exchanges1},
 we will consider the following property:
\begin{itemize}
\item[\bc] \emph{$\Pmor$-base change}.--\index{word}{base change!$\Pmor$-base change}
 For any Cartesian square
$$
\xymatrix{
Y\ar^q[r]\ar_g[d]\ar@{}|\Delta[rd] & X\ar^f[d] \\
T\ar_p[r] & S
}
$$
such that $p$ is a $\Pmor$-morphism, the exchange transformation
$$
Ex(\Delta_\sharp^*):q_\sharp g^* \rightarrow f^* p_\sharp
$$
is an isomorphism.\footnote{In other words,
 $f^*$ commutes with $p_\sharp$.}
\end{itemize}
\end{num}

\begin{df} \label{df:Pmor_fibred}
A \emph{$\Pmor$-fibred category}
\index{word}{fibred!Pfibredcategory@$\Pmor$-fibred category}
 over $\site$ is a pre-$\Pmor$-fibred category $\M$ over $\site$ which satisfies
the property of $\Pmor$-base change.
\end{df}

%\conv \label{convention_identification2}
%Extending convention \ref{convention_identification1}
%in the situation of a $\Pmor$-fibred $\type$-category $\M$,
%we shall generally write the structural isomorphisms of the form
%$Ex(\Delta_\sharp^*)$ as identities (using implicitly
%the coherence results stated in remark \ref{rem:coherence_exchange1}
%-- see also remark \ref{rem:coherence_exchange2} below).

\begin{ex} \label{ex:can_weak_P-fibred&bc}
Consider the notations of Example \ref{ex:can_weak_P-fibred}.
Then the transitivity property of pullbacks of morphisms in $\Pmor$
 amounts to say that the category $\Pmorx ?$ satisfies
  the $\Pmor$-base change property.
Thus, $\Pmorx ?$ is in fact a $\Pmor$-fibred category,
 called \emph{the canonical 
 $\Pmor$-fibred category}.\index{word}{fibred!Pfibredcategory@$\Pmor$-fibred category!canonical}
\end{ex}

\begin{df} \label{df:Pfibred}
A $\Pmor$-fibred category
 $\M$ over $\site$
is \emph{complete}\index{word}{fibred!Pfibredcategory@$\Pmor$-fibred category!complete}
 if, for any morphism $f:T \rightarrow S$, the pullback
 functor $f^*:\M(S) \rightarrow \M(T)$ admits a right adjoint 
 $f_*:\M(S) \rightarrow \M(T)$.
\end{df}

\begin{rem}
In the case where $\Pmor$ is the class of isomorphisms,
a $\Pmor$-fibred category is what we usually call
a bifibred category\index{word}{bifibred category} over $\site$.
\end{rem}

\begin{ex} \label{ex:H_complete_P-fibred}
The pre-$\sm$-fibred category $\Hpt$ of Example \ref{ex:H_weak_P-fibred}
 is a complete $\sm$-fibred category according to
 \cite[p. 102-105, 108-110]{MV}.
\end{ex}

\begin{num} \label{num:exchanges2}
\textit{Exchange structures II}.-- 
Let $\M$ be a complete $\Pmor$-fibred category.
Consider a commutative square
$$
\xymatrix{
Y\ar^q[r]\ar_g[d]\ar@{}|{\Delta}[rd] & X\ar^f[d] \\
T\ar_p[r] & S.
}
$$
We obtain an exchange transformation\index{word}{exchange!transformation}:
$$
Ex(\Delta_*^*):p^*f_* \xrightarrow{ad(g^*,g_*)} g_*g^*p^*f_*=g_*q^*f^*f_*
 \xrightarrow{ad'(f^*,f_*)} g_*q^*.
$$
Assume moreover that $p$ and $q$ are $\Pmor$-morphism.
Then we can check that $Ex(\Delta_*^*)$ is the transpose
 of the exchange $Ex(\Delta_\sharp^*)$.
Thus, when $\Delta$ is Cartesian and $p$ is a $\Pmor$-morphism,
 $Ex(\Delta_*^*)$ is an isomorphism according to \bc.
 
We can also define an exchange transformation:
$$
Ex(\Delta_{\sharp*}):p_\sharp g_*
 \xrightarrow{ad(f^*,f_*)} f_*f^*p_\sharp g_*
 \xrightarrow{Ex(\Delta^*_\sharp)^{-1}} f_*q_\sharp g^* g_*
 \xrightarrow{ad'(g^*,g_*)} f_*q_\sharp.
$$
\end{num}

\begin{rem} \label{rem:coherence_exchange2}
As in remark \ref{rem:coherence_exchange1},
 we obtain coherence\index{word}{coherence} results
 for these exchange transformations.

First with respect to the identifications
 of the kind $f^*g^*=(gf)^*$, $(fg)_*=f_*g_*$, $(fg)_\sharp=f_\sharp g_\sharp$.
 Second when several exchange transformations of different kinds
 are involved. As an example, we consider the following commutative
 diagram in $\site$:
$$
\xymatrix@C=44pt@R=0pt{
& Y\ar^q[rd]\ar@{}|{\Gamma'}[dd] & \\
Z\ar_h[dddddd]\ar^{q'}[ru]\ar|{q'}[rd]\ar@{}_{\Theta}[rddddd]
 & & X\ar^f[dddddd]\ar@{}^{\Delta}[lddddd] \\
& Y\ar|q[ru]\ar|g[dddd] & \\
&& \\
&& \\
&& \\
& T\ar|p[rd]\ar@{}|\Gamma[dd] & \\
Q\ar|{p'}[ru]\ar_{p'}[rd] & & S \\
& T\ar_p[ru] & &
}
$$
Then the following diagram of natural transformations is commutative:
$$
\xymatrix@C=34pt@R=13pt{
q_\sharp g^*p'_*\ar^{Ex(\Delta_\sharp^*)}[r]\ar_{Ex(\Theta_*^*)}[d]
 & f^*p_\sharp p'_*\ar^-{Ex(\Gamma_{\sharp*})}[rd] & \\
q_\sharp q'_*h^*\ar_-{Ex(\Gamma'_{\sharp *})}[rd]
 & & f^*p_*p'_\sharp\ar^{Ex(\Delta_*^*)}[d] \\
& q_*q'_\sharp h^*\ar_{Ex(\Theta_\sharp^*)}[r]
 & q_*g^*p'_\sharp
}
$$
We leave the verification to the reader
 (it is analogous to that of Remark \ref{rem:coherence_exchange1}
  except that it involves also to the compatibility of the unit and counit
  of an adjunction).
\end{rem}
\begin{df} \label{df:P-fibred_transversality}
Let $\M$ be a complete $\Pmor$-fibred category.
Consider a commutative square in $\site$
$$
\xymatrix@=16pt{
Y\ar^q[r]\ar_g[d]\ar@{}|{\Delta}[rd] & X\ar^f[d] \\
T\ar_p[r] & S.
}
$$
We will say that $\Delta$ is \emph{$\M$-transversal}
\index{word}{transversal!$\M$-transversal square}
 if the exchange transformation
$$
Ex(\Delta_*^*):p^*f_* \rightarrow g_*q^*
$$
of \ref{num:exchanges2} is an isomorphism.

Given an admissible class of morphisms $Q$ in $\site$,
we say that $\M$ has
 the \emph{transversality}
\index{word}{transversality property}
 (resp. \emph{cotransversality}) 
\index{word}{cotransversality property}
 \emph{property with respect to $Q$-morphisms},
 if, for any Cartesian square $\Delta$ as above such
 that $f$ is in $Q$ (resp. $p$ is in $Q$),
 $\Delta$ is $\M$-transversal. 
\end{df}

\begin{rem}
Assume $\site$ is a sub-category of the category of schemes.
When $Q$ is the class of smooth morphisms (resp. proper morphisms),
the cotransversality (resp. transversality) property with respect to $Q$
is usually called the
 \emph{smooth base change property}\index{word}{base change!smooth base change}
(resp. \emph{proper base change property}\index{word}{base change!proper base change}).
 See also Definition \ref{df:adj_BC_PF_prop_tri_premotivic}.
\end{rem}

According to Paragraph \ref{num:exchanges2},
 we derive the following consequence of our axioms:
\begin{prop}\label{Pcotransversality}
Any complete $\Pmor$-fibred category has the cotransversality property
with respect to $\Pmor$.
\end{prop}

Let us note for future reference the following corollary:
\begin{cor}\label{immfullyfaithful}
If $\M$ is a $\Pmor$-fibred category, then, 
for any monomorphism $j:U\rightarrow S$ in $\Pmor$, the functor
$j_\sharp$ is fully faithful. 
If moreover $\M$ is complete, then
the functor $j_*$ is fully faithful as well.
\end{cor}
\begin{proof}
Because $j$ is a monomorphism, we get a Cartesian square
 in $\site$:
$$
\xymatrix@=16pt{
U\ar@{=}[r]\ar@{=}[d]\ar@{}|{\Delta}[rd] & U\ar^j[d] \\
U\ar_j[r] & S.
}
$$
Remark that $Ex(\Delta_\sharp^*):1 \rightarrow j^*j_\sharp$
is the unit of the adjunction $(j_\sharp,j^*)$.
Thus the $\Pmor$-base change property
 shows that $j_\sharp$ is fully faithful.

Assume $\M$ is complete. We remark similarly that
$Ex(\Delta_*^*):j^*j_* \rightarrow 1$ is the counit
 of the adjunction $(j^*,j_*)$. Thus, the above proposition
 shows readily that $j_*$ is fully faithful.
\end{proof}

%%%%%%%%%%%%%%%%%%%%%%%%%%%%%%%%%%%%%%%%%%%%%%%%%%%%%%%%%%%%%%%%%%%%%%%%%%%%%%

\subsubsection{Monoidal structures}

Let $\catm$ be the sub-$2$-category of $\cat$
 made of symmetric monoidal categories
 whose $1$-morphisms are (strong) symmetric monoidal functors
 and $2$-morphisms are symmetric monoidal transformations.
\begin{df}
A \emph{monoidal pre-$\Pmor$-fibred category 
\index{word}{fibred!monoidal pre-$\Pmor$-fibred category}
over $\site$}

 is a $2$-functor
$$\M:\site \rightarrow \catm$$
such that $\M$ is a pre-$\Pmor$-fibred category.
\end{df}
In other words, $\M$ is a pre-$\Pmor$-fibred category
such that each of its fibers $\M(S)$ is endowed with
a structure of a monoidal category, 
and any pullback morphism $f^*$ is monoidal, with
the obvious coherent structures.
For an object $S$ of $\site$,
we will usually denote by $\otimes_S$ (resp. $\un_S$)
the tensor product (resp. unit) of $\M(S)$. \\
In particular, we then have the following natural isomorphisms:
\begin{itemize}
\item for a morphism $f:T \rightarrow S$ in $\site$,
 and objects $M$, $N$ of $\M(S)$,
\begin{equation*} \label{eq:iso_pullback_monoid}
f^*(M) \otimes_T f^*(N) \xrightarrow{\sim} f^*(M \otimes_S N);
\end{equation*}
\item for a morphism $f:T \rightarrow S$ in $\site$,
\begin{equation*} \label{eq:iso_pullback_unit}
 f^*(\un_S)  \xrightarrow{\sim} \un_T\, .
\end{equation*}
\end{itemize}

\begin{conv} \label{convention_identification3}
As in convention \ref{convention_identification1},
 we will write formula as though these structural isomorphisms
 are identities.
\end{conv}
 
\begin{ex} \label{ex:can_mon_weak_P-fibred}
Consider the notations of Example \ref{ex:can_weak_P-fibred}.

Using the properties (Pb) and (Pc) of $\Pmor$ (see Paragraph \ref{defclassP}),
 for two $S$-objects $X$ and $Y$ in $\Pmorx S$,
 the Cartesian product $X \times_S Y$ is an object of $\Pmorx S$.
This defines a symmetric monoidal structure on $\Pmorx S$
 with unit the trivial $S$-object $S$.
Moreover, the functor $f^*$ defined in Example \ref{ex:can_weak_P-fibred}
 is monoidal. 
Thus, the pre-$\Pmor$-fibred category $\Pmorx ?$ is in fact monoidal.
\end{ex}

\begin{num} \label{num:exchanges3}
\textit{Monoidal exchange structures I}.
Let $\M$ be a monoidal pre-$\Pmor$-fibred category over $\site$.

Consider a $\Pmor$-morphism $f:T \rightarrow S$,
and $M$ (resp. $N$) an object of $\M(T)$ (resp. $\M(S)$).

We get a morphism in $\M(S)$
$$
Ex(f_\sharp^*,\otimes):
f_\sharp(M \otimes_T f^*(N))\longrightarrow f_\sharp(M) \otimes_S N
$$
as the composition
$$f_\sharp(M \otimes_T f^*(N))
\rightarrow f_\sharp(f^*f_\sharp(M) \otimes_T f^*(N))
\simeq f_\sharp f^*(f_\sharp(M) \otimes_S N)
\rightarrow f_\sharp(M) \otimes_S N\, .$$
This map is natural in $M$ and $N$.
It will be called the \emph{exchange transformation}
\index{word}{exchange!transformation}
between $f_\sharp$ and $\otimes_T$.

Remark also that the functor $f_\sharp$,
as a left adjoint of a symmetric monoidal functor,
 is colax symmetric monoidal: for any objects $M$ and $N$ of $\M(T)$,
there is a canonical morphism
\begin{equation} \label{num:exchanges3ter}
f_\sharp(M) \otimes_S f_\sharp(N) \rightarrow f_\sharp(M \otimes_T N)
\end{equation}
natural in $M$ and $N$, as well as a natural map
\begin{equation} \label{num:exchanges3quad}
f_\sharp(\un_T)\rightarrow \un_S\, .
\end{equation}
\end{num}

\begin{rem} \label{rem:coherence_exchange3}
As in remark \ref{rem:coherence_exchange1},
 the preceding exchange transformations satisfy
 a coherence\index{word}{coherence} condition
 for composable morphisms $W \xrightarrow g T \xrightarrow f S$.
We get in fact a commutative diagram:
$$
\xymatrix@C=35pt@R=18pt{
(fg)_\sharp\big(M \otimes_S (fg)^*(N)\big)\ar^-{Ex((fg)_\sharp^*,\otimes)}[rr]\ar@{=}[d]
 && \big((fg)_\sharp(M)\big) \otimes_W N\ar@{=}[d] \\
f_\sharp g_\sharp \big(M \otimes_S g^*f^*(N)\big)\ar^-{Ex(g_\sharp^*,\otimes)}[r]
 & f_\sharp \big(g_\sharp(M) \otimes_T f^*(N)\big)\ar^-{Ex(f_\sharp^*,\otimes)}[r]
 & \big(f_\sharp g_\sharp(M)\big) \otimes_W N
}
$$
As in remark \ref{rem:coherence_exchange2},
 there is also a coherence relation when different kinds
 of exchange transformations are involved.
Consider a commutative square in $\site$
$$
\xymatrix{
Y\ar^q[r]\ar_g[d]\ar@{}|{\Delta}[rd] & X\ar^f[d] \\
T\ar_p[r] & S
}
$$
such that $p$ and $q$ are $\Pmor$-morphisms.
Then the following diagram is commutative:
$$
\xymatrix@C=36pt@R=15pt{
q_\sharp g^*(M \otimes_T p^*N)\ar^{Ex(\Delta_\sharp^*)}[r]\ar@{=}[d]
 & f^*p_\sharp(M \otimes_T p^*N)\ar^{Ex(p_\sharp^*,\otimes)}[r]
 & f^*(p_\sharp M \otimes_S N)\ar@{=}[d] \\
q_\sharp(g^*M \otimes_Y q^*f^*N)\ar^{Ex(q_\sharp^*,\otimes)}[r]
 & (q_\sharp g^*M) \otimes_X f^*N\ar^{Ex(\Delta_\sharp^*)}[r]
 & (f^*p_\sharp M) \otimes_X f^*N
}
$$
We leave the verification to the reader.
\end{rem}

\begin{num} \label{basicprojformula}
Under the assumptions of \ref{num:exchanges3},
 we will consider the following property:
\begin{itemize}
\item[\pf] \emph{$\Pmor$-projection formula}.--
\index{word}{projection formula!$\Pmor$-projection formula}
For any $\Pmor$-morphism $f:T \rightarrow S$ the exchange transformation
$$
Ex(f_\sharp,\otimes_T):f_\sharp(M \otimes_T f^*(N))
 \rightarrow f_\sharp(M) \otimes_S N
$$
is an isomorphism for all $M$ and $N$.
\end{itemize}
\end{num}

\begin{df} \label{df:mono_Pmor_fibred}
A \emph{monoidal} $\Pmor$-fibred category
\index{word}{fibred!monoidalPfibredcategory@monoidal $\Pmor$-fibred category}
 over $\site$ is a
 monoidal pre-$\Pmor$-fibred category $\M:\op{\site}\rightarrow \catm$ over $\site$
 which satisfies the $\Pmor$-projection formula.
\end{df}

\begin{ex} \label{ex:can_mon_weak_P-fibred&bc}
Consider the canonical monoidal weak $\Pmor$-fibred category $\Pmorx ?$
 (see Example \ref{ex:can_mon_weak_P-fibred}).
The transitivity property of pullbacks
 implies readily that $\Pmorx ?$ satisfies the property \pf.
 Thus,  $\Pmorx ?$ is in fact a monoidal $\Pmor$-fibred category
 called \emph{canonical}.
\index{word}{fibred!Pfibredcategory@$\Pmor$-fibred category!canonical monoidal}
\end{ex}

\begin{df} \label{df:monoidal_Pfibred}
A monoidal $\Pmor$-fibred category 
\index{word}{fibred!Pfibredcategory@$\Pmor$-fibred category!complete monoidal}
$\M$ over $\site$ is \emph{complete} if it satisfies the following conditions:
\begin{enumerate}
\item $\M$ is complete as a $\Pmor$-fibred category.
\item For any object $S$ of $\site$, 
 the monoidal category $\M(S)$ is closed (i.e. has an internal Hom).
\end{enumerate}
\end{df}
In this case, we will usually denote by
$\uHom_S$ the internal Hom in $\M(S)$, so that we have natural
bijections
$$
\Hom_{\M(S)}(A\otimes_S B,C)\simeq\Hom_{\M(S)}(A,\uHom_S(B,C))\, .
$$

\begin{ex} \label{ex:H_monoidal_P-fibred}
The $\Pmor$-fibred category $\Hpt$ of Example \ref{ex:H_complete_P-fibred}
 is in fact a complete monoidal $\Pmor$-fibred category.
 The tensor product is given by the smash product (see \cite{MV}).
\end{ex}

\begin{num} \label{num:exchanges4}
\textit{Monoidal exchange structures II}.-- Let $\M$ be a complete monoidal $\Pmor$-fibred category.

Consider a morphism $f:T \rightarrow S$ in $\site$.
Then we obtain an exchange transformation:
 \index{word}{exchange!transformation}
\begin{align*}
Ex(f_*^*,\otimes_S):(f_*M) \otimes_S N
 & \xrightarrow{ad(f^*,f_*)} f_*f^*\big((f_*M) \otimes_S N\big) \\
 & =f_*\big((f^*f_*M) \otimes_T f^*N\big) 
    \xrightarrow{ad'(f^*,f_*)} f_*(M \otimes_T f^*N).
\end{align*}
\end{num}

\begin{rem} \label{rem:coherence_exchange4}
As in remark \ref{rem:coherence_exchange3},
 these exchange transformations are compatible with the identifications
 $(fg)_*=f_*g_*$ and $(fg)^*=g^*f^*$.
 Moreover, there is a coherence relation when composing the
 exchange transformations of the kind $Ex(f_*^*,\otimes)$
 with exchange transformations of the kind $Ex(\Delta_*^*)$
 as in \emph{loc. cit.}
 Finally, note that there is another kind of coherence relations
 involving $Ex(f_*^*,\otimes)$, $Ex(\Delta_\sharp^*)$
  (resp. $Ex(f_\sharp^*,\otimes)$) and $Ex(\Delta_{\sharp *})$.
  
We leave the formulation of these coherence relations to the
 reader, on the model of the preceding ones.
 \end{rem}

\begin{num} \label{num:exchanges4bis}
 \textit{Monoidal exchange structures III}.--
Let $\M$ be a complete monoidal $\Pmor$-fibred category
 and $f:T \rightarrow S$ be a morphism in $\site$.
 
Because $f^*$ is monoidal,
 we get by adjunction a canonical isomorphism
$$
\uHom_S(M,f_* N) \rightarrow f_*\uHom_T(f^*M,N).
$$
Assume that $f$ is a $\Pmor$-morphism.
Then from the $\Pmor$-projection formula,
 we get by adjunction two canonical isomorphisms:
\begin{align*}
f^* \uHom_S(M,N) &\rightarrow \uHom_T(f^*M,f^*N), \\
\uHom_S(f_\sharp M,N) &\rightarrow f_*\uHom_T(M,f^*N)
\end{align*}
These isomorphisms are generically called \emph{exchange isomorphisms}.
\index{word}{exchange!isomorphism}
\end{num}

\subsubsection{Geometric sections}

\begin{num} \label{num:def_premotives}
Consider a complete monoidal $\Pmor$-fibred category $\M$.

Let $S$ be a scheme. For any $\Pmor$-morphism $p:X \rightarrow S$,
 we put $M_S(X):=p_\sharp(\un_X)$.
 According to our conventions, 
 this object is identified with $p_\sharp p^*(\un_S)$.
 As the $\Pmor$-fibred category $\M$ is complete, the functor
 $p_\sharp p^*$ is left adjoint to $p_*p^*$.
 Consider a commutative diagram of schemes in $\base$:
$$
\xymatrix@=10pt{
Y\ar_q[rd]\ar^f[rr] && X\ar^p[ld] \\
 & S &
}
$$
such that $p$ and $q$ are in $\Pmor$. In other words, $f$ is a morphism
 in the category $\Pmor/S$ of Example \ref{ex:can_weak_P-fibred}.
 Then we get a natural transformation of functors:
\begin{equation}\label{eq:num:def_premotives}
p_*p^* \xrightarrow{ad(f^*,f_*)} p_*f_*f^*p^*=q_*q^*.
\end{equation}
By adjunction, one deduces a natural transformation:
$$
q_\sharp q^* \rightarrow p_\sharp p^*
$$
which gives a morphism $M_S(Y) \xrightarrow{f_*} M_S(X)$.
 One can check that the relation $f_*g_*=(fg)_*$ holds
 --- by reducing to the same assertion for the map
 \eqref{eq:num:def_premotives} which follows by a standard
 $2$-functoriality argument.
 Therefore, one has obtained a covariant functor
 $M_S:\Pmor/S \rightarrow \M$.

Consider a Cartesian square in $\site$
$$
\xymatrix{
Y\ar^{g}[r]\ar_{q}[d]\ar@{}|{\Delta}[rd] & X\ar^p[d] \\
T\ar_f[r] & S
}
$$
such that $p$ is a $\Pmor$-morphism.
With the notations of Example \ref{ex:can_weak_P-fibred}, $Y=f^*(X)$.
Then we get a natural exchange transformation
$$
Ex(M_T,f^*):M_T(f^*(X))=q_\sharp(\un_Y)=q_\sharp  g^*(\un_X)
 \xrightarrow{Ex(\Delta_\sharp^*)} f^* p_\sharp(\un_X)=f^* M_S(X).
$$
In other words,
 $M$ defines a lax natural transformation $\Pmorx ? \rightarrow \M$.
 
Consider $\Pmor$-morphisms $p:X \rightarrow S$, $q:Y \rightarrow S$.
Let $Z=X \times_S Y$ be the Cartesian product and consider the Cartesian
square:
$$
\xymatrix{
Z\ar^{p'}[r]\ar_{q'}[d]\ar@{}|\Theta[rd] & Y\ar^q[d] \\
X\ar|p[r] & S.
}
$$
Using the exchange transformations of the preceding paragraph,
 we get a canonical morphism
$$
Ex(M_S,\otimes_S):M_S(X\times_S Y)\longrightarrow M_S(X) \otimes_S M_S(Y)
$$
as the composition
\begin{align*}
M_S(X \times_S Y)=p_\sharp q'_\sharp  p'^*(\un_Y)
 &\xrightarrow{Ex(\Theta_\sharp^*)} p_\sharp p^* q_\sharp (\un_Y)
 =p_\sharp (\un_X \otimes_X p^* q_\sharp(\un_Y)) \\
 &\xrightarrow{Ex(p_\sharp,\otimes_X)} p_\sharp(\un_X) \otimes_S q_\sharp(\un_Y)
 =M_S(X) \otimes_S M_S(Y).
\end{align*}
In other words, the functor $M_S$ is symmetric colax monoidal.

Remark finally that for any $\Pmor$-morphism $p:T \rightarrow S$,
 and any $\Pmor$-object $Y$ over $T$,
 we obtain according to convention an identification
 $p_\sharp M_T(Y)=M_S(Y)$.
\end{num}
\begin{df} \label{df:geometric_sections}
Given a complete monoidal $\Pmor$-fibred category $\M$ over $\site$,
 the lax natural transformation $M:\Pmorx ? \rightarrow \M$ constructed
 above will be called the \emph{geometric sections
\index{word}{section!geometric}
 of $\M$}.
\end{df}

The following lemma is obvious from the definitions above:
\begin{lm} \label{lm:basics_M}
let $\M$ be a complete monoidal $\Pmor$-fibred category.
Let $M:\Pmorx ? \rightarrow \M$ be the geometric sections of $\M$.
Then:
\begin{enumerate}
\item[(i)] For any morphism $f:T \rightarrow S$ in $\site$,
 the exchange $Ex(M_T,f^*)$ defined above is an isomorphism.
\item[(ii)] For any scheme $S$,
 the exchange $Ex(M_S,\otimes_S)$ defined above is an isomorphism.
\end{enumerate}
In other words, $M$ is a Cartesian functor
 and $M_S$ is a \emph{(strong)} symmetric monoidal functor.
\end{lm}

\begin{num} \label{num:basics_M}
In the situation of the lemma we thus obtain the following
 isomorphisms:
\begin{itemize}
\item $f^*M_S(X)\simeq M_T(X \times_S T)$,
\item $p_\sharp M_T(Y)\simeq M_S(Y)$,
\item $M_S(X \times_S Y)\simeq M_S(X) \otimes_S M_S(Y)$,
\end{itemize}
whenever it makes sense.
\end{num}

\subsubsection{Twists}

\begin{num} Let $\M$ be a pre-$\Pmor$-fibred category of $\site$.
Recall that a Cartesian section of $\M$
 (\emph{i.e.} a Cartesian functor $A:\site \rightarrow \M$)
 is the data of an object $A_S$ of $\M(S)$ for each object $S$ of $\S$
 and of isomorphisms
$$
f^*(A_S) \xrightarrow \sim A_T
$$
for each morphism $f:T \rightarrow S$,
 subject to coherence identities; see \cite[Exp.~VI]{SGA1}.

If $\M$ is monoidal,
 the tensor product of two Cartesian sections is defined termwise.
\end{num} 
 
\begin{df} \label{df:twists}
let $\M$ be a monoidal pre-$\Pmor$-fibred category.
A set of \emph{twists} 
\index{word}{twist}
$\tau$ for $\M$ is a set of Cartesian sections 
of $\M$ which is stable by tensor product (up to isomorphism),
and contains the unit $\un$. For short, when $\M$ is endowed with a set of twists $\tau$,
we say also that $\M$ is \emph{$\tau$-twisted}.
\index{word}{twist!tautwisted@$\tau$-twisted}
\end{df}

\begin{num} \label{num:exchange_twists}
Let $\M$ be a monoidal pre-$\Pmor$-fibred category endowed with a set
of twists $\tau$.

The tensor product on $\tau$ induces a monoid structure
that we will denote by $+$
 (the unit object of $\tau$ will be written $0$).

Consider an object $i \in \tau$.
For any object $S$ of $\site$, 
 we thus obtain an object $t(i)_S$ in $\M(S)$ associated with $i$.
Given any object $M$ of $\M(S)$, we simply put: 
$$M\{i\}=M \otimes_S i_S$$
and call this object the twist
 of $M$ by $i$. We also define $M\{0\}=M$.
 
For any $i,j \in \tau$, and any object $M$ of $\M(S)$,
 we define $M\{i+j\}=(M\{i\})\{j\}$.
Given a morphism $f:T \rightarrow S$, an object $M$ of $\M(S)$
 and a twist $i \in \tau$, we also obtain $f^*(M\{i\})=(f^*M)\{i\}$.
If $f$ is a $\Pmor$-morphism, for any object $M$ of $\M(T)$,
the exchange transformation $Ex(f^*_\sharp,\otimes_T)$
of paragraph \ref{num:exchanges1} induces a canonical
morphism
$$
Ex(f_\sharp,\{i\}):f_\sharp(M\{i\}) \rightarrow (f_\sharp M)\{i\}.
$$
We will say that \emph{$f_\sharp$ commutes with $\tau$-twists}
\index{word}{twist!commutes with $\tau$-twists (\emph{or} twists)}
 (or simply twists when $\tau$ is clear) if for any $i \in \tau$,
 the natural transformation $Ex(f_\sharp,\{i\})$ is an isomorphism.
\end{num} 
 
\begin{df} \label{df:generating_twists}
Let $\M$ be a complete monoidal $\Pmor$-fibred category with a set of twists $\tau$
 and $M:\Pmorx ? \rightarrow \M$ be the geometric sections 
\index{word}{section!geometric}
of $\M$.
 
We say $\M$ is \emph{$\tau$-generated}
\index{word}{generated!$\tau$-generated}
\index{word}{fibred!Pfibredcategory@$\Pmor$-fibred category!$\tau$-generated|see{generated}}
 if for any object $S$ of $\site$,
the family of functors
$$
\Hom_{\M(S)}(M_S(X)\{i\},-):\M(S)\rightarrow \set
$$
indexed by a $\Pmor$-object $X/S$ and an element $i \in \tau$ is conservative.

Of course, we do not exclude the case where $\tau$ is trivial, but then,
 we shall simply say that $\M$ is \emph{geometrically generated}.
\index{word}{generated!geometrically generated}
\index{word}{fibred!Pfibredcategory@$\Pmor$-fibred category!geometrically generated|see{generated}}
%% We say $\M$ is \emph{compactly $\tau$-generated} if it is $\tau$-generated
%%  and for any $\Pmor$-object $X/S$, any twist $i\in \tau$,
%%  $M_S(X)\{i\}$ is compact\footnote{\emph{i.e.} The functor $\Hom_{\M(S)}(M_S(X)\{i\},-)$
%%   commutes with sums.}.
\end{df}

%% 
%% \begin{ex} \label{ex:H_P-fibred+twists}
%% The monoidal $\sm$-fibred category $\Hpt$ of Example \ref{ex:H_monoidal_P-fibred}
%%  is compactly $\NN$-generated where $1$ corresponds to the simplicial sphere $S^1$.
%% \end{ex}

We shall frequently use the following proposition to characterize
 complete monoidal $\Pmor$-fibred categories over $\site$:
\begin{prop} \label{prop:carac_monoidal_P-fibred}
Let $\M:\site \rightarrow \catm$ be a $2$-functor such that:
\begin{enumerate}
\item For any $\Pmor$-morphism $f:T \rightarrow S$, the pullback
 functor $f^*:\M(S) \rightarrow \M(T)$ is monoidal and 
 admits a left adjoint $f_\sharp$ in $\type$.
\item For any morphism $f:T \rightarrow S$, the pullback
 functor $f^*:\M(S) \rightarrow \M(T)$ admits a right adjoint 
 $f_*$ in $\type$.
\end{enumerate}
We consider $\M$ as a complete monoidal $\Pmor$-fibred category
 and denote by $M:\Pmorx ? \rightarrow \M$ its associated geometric sections.
Suppose given a set of twists $\tau$ such that $\M$ is $\tau$-generated.
\index{word}{generated!$\tau$-generated}
Then, the following assertions are equivalent:
\begin{enumerate}
\item[(i)] $\M$ satisfies properties \bc and \pf \\
(\emph{i.e.} $\M$ is a complete monoidal $\Pmor$-fibred category.)
\item[(ii)] 
\begin{enumerate}
\item $M$ is a Cartesian functor.
\item For any object $S$ of $\site$, $M_S$ is (strongly) monoidal.
\item For any $\Pmor$-morphism $f$, $f_\sharp$ commutes with $\tau$-twists.
\end{enumerate}
\end{enumerate}
\end{prop}
\begin{proof}
$(i) \Rightarrow (ii)$: This is obvious (see Lemma \ref{lm:basics_M}).

$(ii) \Rightarrow (i)$: We use the following easy lemma:
\begin{lm} \label{lm:generators&adjoints}
Let $\type_1$ and $\type_2$ be categories, 
 $F,G:\type_1 \rightarrow \type_2$ be two left adjoint functors,
 and $\eta:F \rightarrow G$ be a natural transformation.
Let $\cG$ be a class of objects of $\type_1$ which is generating in the
sense that the family of functors $\Hom_{\type_1}(X,-)$ for $X$ in $\cG$
is conservative.

Then the following conditions are equivalent:
\begin{enumerate}
\item $\eta$ is an isomorphism.
\item For all $X$ in $\cG$, $\eta_X$ is an isomorphism.
\end{enumerate}
\end{lm}
\noindent Given this lemma, to prove property \bc, 
we are reduced to check that the
exchange transformation $Ex(\Delta_\sharp^*)$ 
 is an isomorphism when
evaluated on an object $M_T(U)\{i\}$
 for an object $U$ of $\Pmorx T$ and a twist $i \in \tau$. 
Then it follows from (ii), \ref{num:exchange_twists}
 and Example \ref{ex:can_weak_P-fibred&bc}.\footnote{
The cautious reader will use remark \ref{rem:coherence_exchange1}
 to check that the corresponding map 
$$M_X(U \times_T Y)\{i\} \rightarrow M_X(U \times_T Y)\{i\}$$
 is the identity.} \\
To prove property \pf,
 we proceed in two steps first proving the case $M=M_T(U)\{i\}$
 and $N$ any object of $\M(S)$ using the same argument as above
 with the help of \ref{ex:can_mon_weak_P-fibred&bc}.
Then, we can prove the general case
 by another application of the same argument.
\end{proof}

Suppose given a complete monoidal $\Pmor$-fibred category $\M$
 with a set of twists $\tau$.
Let $f:T \rightarrow S$ be a morphism of $\site$.
Then the exchange transformation \ref{num:exchanges4}
 induces for any $i \in \tau$ an exchange transformation
$$
Ex(f_*,\{i\}):(f_* M)\{i\}\to f_*(M\{i\}) \, .
$$
\begin{df} \label{df:f_*_commutes_twists}
In the situation above, we say that $f_*$
\emph{commutes with $\tau$-twists}
\index{word}{twist!commutes with $\tau$-twists (\emph{or} twists)}
 (or simply with twists when $\tau$ is clear)
if, for any $i \in \tau$, the exchange transformation $Ex(f_*,\{i\})$
is an isomorphism.
\end{df}
It will frequently happen that twists are $\otimes$-invertible.
 Then $f_*$ commutes with twists as its right adjoint does.

%%%%%%%%%%%%%%%%%%%%%%%%%%%%%%%%%%%%%%%%%%%%%%%%%%%%%%%%%%%%%%%%%%%%%%%%%%%%%%%%%%
% morphisms

\subsection{Morphisms of $\Pmor$-fibred categories}
 
\subsubsection{General case}
 
\begin{num} \label{num:Pfibred_morph}
Consider two $\Pmor$-fibred categories $\M$ and $\M'$ over $\site$,
as well as a Cartesian functor $\varphi^*:\M \rightarrow \M'$
between the underlying fibred categories: for any object $S$ of $\site$,
we have a functor
$$\varphi^*_S:\M(S)\rightarrow \M'(S)\, ,$$
and for any map $f:T\rightarrow S$ in $\site$, we have an isomorphism of functors $c_f$
\begin{equation}\label{num:Pfibred_morphbis}\begin{split}
\xymatrixrowsep{.9pc} \xymatrixcolsep{1.2pc}
\UseAllTwocells \xymatrix{
\M(S)\ddrrcompositemap\omit{\quad c_{f}} 
\ar[rr]^{\varphi^*_S} \ar[dd]_{f^*}
&&\M'(S)\ar[dd]^{f^*} \\ &&  \\
\M(T)\ar[rr]_{\varphi^*_T}  &&\M'(T)&
}\end{split}
\qquad c_f:f^*\, \varphi^*_S\xrightarrow{\sim} \varphi^*_T\, f^*
\end{equation}
satisfying some cocycle condition with respect to composition in $\site$.

For any $\Pmor$-morphism $p:T \rightarrow S$,
we construct an exchange morphism
\index{word}{exchange!morphism|see{exchange transformation}}
\index{word}{exchange!transformation}
\begin{equation*}
Ex(p_\sharp,\varphi^*):p_\sharp\,  \varphi^*_T\longrightarrow \varphi_S^*\,  p_\sharp
\end{equation*}
as the composition
$$
p_\sharp \varphi^*_T
 \xrightarrow{ad(p_\sharp,p^*)} p_\sharp \varphi_T^* p^* p_\sharp
 \xrightarrow{c_p^{-1}} p_\sharp p^* \varphi_S^* p_\sharp
 \xrightarrow{ad'(p_\sharp,p^*)} \varphi_S^* p_\sharp.
$$
\end{num}

\begin{df} \label{df:Pfibred_adjunction}
Consider the situation above.
We say that the Cartesian functor
$$\varphi^*:\M \rightarrow \M'$$
is a \emph{morphism of $\Pmor$-fibred categories}
\index{word}{morphism!of $\Pmor$-fibred categories}
\index{word}{adjunction!of $\Pmor$-fibred categories|see{morphism of}}
 if, for any $\Pmor$-morphism $p$,
the exchange transformation $Ex(p_\sharp,\varphi^*)$ is an isomorphism.
\end{df}

%% \conv \label{convention_identification5}
%% In the same spirit as for convention \ref{convention_identification3},
%% we shall write the structural isomorphisms
%% \eqref{exchangemorphismofPfibredcat} of a morphisms of $\Pmor$-fibred $\type$-categories 
%% as identities.

\begin{ex}
If $\M$ is a monoidal $\Pmor$-fibred category, then the geometric sections
\index{word}{section!geometric}
$M:\Pmorx ?\rightarrow \M$ is a morphism of $\Pmor$-fibred categories
 (\ref{lm:basics_M}).
\end{ex}

\begin{df}
Let $\M$ and $\M'$ be two complete $\Pmor$-fibred categories.
A \emph{morphism of complete $\Pmor$-fibred categories}
\index{word}{morphism!of complete $\Pmor$-fibred categories}
is a morphism of $\Pmor$-fibred categories 
$$\varphi^*:\M\rightarrow \M'$$
such that, for any object $S$
of $\site$, the functor $\varphi^*_S:\M(S)\rightarrow \M'(S)$
has a right adjoint
$$\varphi_{*,S}:\M'(S)\rightarrow \M(S)\, .$$
\end{df}
When we want to indicate a notation for the right adjoint of a morphism
as above, we use the writing
$$
\varphi^*:\M \rightleftarrows \N:\varphi_*
$$
the left adjoint being in the left hand side.

\begin{num} \label{num:exchanges_Pfibred_morph}
\textit{Exchange structures III}.
Consider a morphism $\varphi^*:\M\rightarrow \M'$
 of complete $\Pmor$-fibred categories.

Then for any morphism $f:T \rightarrow S$ in $\site$,
 we define exchange transformations
\index{word}{exchange!transformation}
\begin{align}\label{exchangePmorphismfupperstar}
Ex(\varphi^*,f_*)&:\varphi^*_S f_*\longrightarrow f_* \varphi_T^*, \\
\label{eq:exchange_f_*_Padjunctions}
Ex(f^*,\varphi_*)&:f^*\varphi_{*,S} \longrightarrow \varphi_{*,T} f^*,
\end{align}
as the respective compositions
\begin{align*}
\varphi^*_S f_*
 \xrightarrow{ad(f^*,f_*)} f_* f^* \varphi_S^* f_*
&\simeq f_* \varphi_T^* f^* f_*
 \xrightarrow{ad'(f^*,f_*)} f_* \varphi_T^*, \\
f^*\varphi_{*,S}
 \xrightarrow{ad(f^*,f_*)} f^* \varphi_{*,S} f_* f^*
&\simeq f^* f_* \varphi_{*,T}  f^*
 \xrightarrow{ad'(f^*,f_*)} \varphi_{*,T}  f^*.
\end{align*}
\end{num}

\begin{rem}
We warn the reader that $\varphi_*:\M' \rightarrow \M$ is
 not a Cartesian functor in general,
 meaning that the exchange transformation $Ex(f^*,\varphi_*)$
 is not necessarily an isomorphism,
 even when $f$ is a $\Pmor$-morphism.
%We warn the reader that this notation is abusive
%for the following reasons:
%the functors $\varphi_{*,S}:\M'(S)\rightarrow \M(S)$
%do not necessarily induce a Cartesian functor $\M'\rightarrow \M$ over $\site$,
%but only a (lax) morphism of $2$-functors: for any map $f:T\rightarrow S$
%in $\site$, we have only an exchange map
%$$f^*\, \varphi_{*,S}\rightarrow \varphi_{*,T} \, f^*$$
%which have not any reason to be invertible in general;
%even in the case where $\varphi_*$ defines a Cartesian functor
%from $\M'$ to $\M$ over $\site$ (i.e. when the exchange
%map above is invertible for all $f$), there is no reason in general for
%$\varphi_*$ to be a morphism of $\Pmor$-fibred categories.
\end{rem}

\subsubsection{Monoidal case}

\begin{df} \label{df:morph_mon_P-fibred}
Let $\M$ and $\M'$ be monoidal $\Pmor$-fibred categories.

A \emph{morphism of monoidal $\Pmor$-fibred categories}
is a morphism $\varphi^*:\M \rightarrow \M'$ of $\Pmor$-fibred categories
such that for any object $S$ of $\site$,
the functor
$\varphi_S^*:\M(X) \rightarrow \N(S)$ has the structure of
a (strong) symmetric monoidal
functor, and such that the structural isomorphisms \eqref{num:Pfibred_morphbis}
are isomorphisms of symmetric monoidal functors.

In the case where $\M$ and $\M'$ are complete monoidal $\Pmor$-fibred
categories, we shall say that such a morphism $\varphi^*$
is a \emph{morphism of complete monoidal $\Pmor$-fibred categories}
if $\varphi^*$ is also a morphism of complete $\Pmor$-fibred categories.
\end{df}

\begin{rem} \label{convention_identification6}
If we denote by $M(-,\M)$ and $M(-,\M')$ the geometric sections
 of $\M$ and $\M'$ respectively, we have a natural identification:
$$
\varphi_S^*(M_S(X,\M))\simeq M_S(X,\M')\, .
$$
\end{rem}

\begin{num} \label{num:mon_exchange4}
\textit{Monoidal exchange structures IV}.
Consider a a morphism $\varphi^*:\M \rightarrow \M'$
of complete monoidal $\Pmor$-fibred categories.
For objects $M$ (resp. $N$) of $\M(S)$ (resp. $\M'(S)$),
 we define an exchange transformation
$$
Ex(\varphi_*,\otimes,\varphi^*):(\varphi_{*,S}M) \otimes_S N
 \rightarrow \varphi_{*,S}(M \otimes_T \varphi_S^*N),
$$
natural in $M$ and $N$, as the following composite
\begin{align*}
(\varphi_{*,S}M) &\otimes_S N
 \xrightarrow{ad(\varphi^*,\varphi_*)}
  \varphi_{*,S}\varphi_S^*((\varphi_{*,S}M) \otimes_S N) \\
  &=\varphi_{*,S}((\varphi_S^*\varphi_{*,S}M) \otimes_T \varphi_S^*N) 
    \xrightarrow{ad'(\varphi^*,\varphi_*)} \varphi_{*,S}(M \otimes_T \varphi_S^*N).
\end{align*}

As in remark \ref{rem:coherence_exchange4},
 we get coherence\index{word}{coherence} relations
 between the various exchange transformations
 associated with a morphism of monoidal $\Pmor$-fibred categories.
 We left the formulation to the reader. \\
Note also that, because $\varphi^*$ is monoidal,
we get by adjunction a canonical isomorphism:
$$
 \uHom_{\M(S)}(M,\varphi_{*,S} M')
 \xrightarrow{\sim} \varphi_{*,S}\uHom_{\M'(S)}(\varphi_S^*M,M') \, .
$$
\end{num}

\begin{num} \label{num:morph_monoidal_Pfibred}
Consider two monoidal $\Pmor$-fibred categories $\M$, $\M'$
 and a Cartesian functor $\varphi^*:\M \rightarrow \M'$
 such that, for any scheme $S$,
 $\varphi^*_S:\M(S) \rightarrow \M'(S)$ is monoidal.

Given a Cartesian section $K=(K_S)_{S \in \site}$ of $\M$,
 we obtain for any morphism $f:T \rightarrow S$ in $\site$ a canonical map
$$
f^*\varphi^*_S(K_S)=\varphi^*_T(f^*(K_S)) \rightarrow \varphi^*_T(K_T)
$$
which defines a Cartesian section of $\M'$, which
we denote by $\varphi^*(K)$.
\end{num}

\begin{df} \label{df:morphism&twists}
Let $(\M,\tau)$ and $(\M',\tau')$ be twisted monoidal $\Pmor$-fibred categories.
Let $\varphi^*:\M \rightarrow \M'$ be a Cartesian functor as above
 (resp. a morphism of monoidal $\Pmor$-fibred categories).

We say that
$\varphi^*:(\M,\tau)\rightarrow (\M',\tau')$
 is \emph{compatible with twists}
\index{word}{compatible with twists}
% (resp. is a morphism of twisted $\Pmor$-fibred categories)
 if for any $i \in \tau$, 
 the Cartesian section $\varphi^*(i)$ is in $\tau'$ (up to isomorphism in $\M'$).
\end{df}

\begin{rem}
In particular, $\varphi^*$ induces a map $\tau \rightarrow \tau'$
 (if we consider the isomorphism classes of objects).
Moreover, for any object $K$ of $\M(S)$ and any twist $i \in \tau$,
we get an identification:
$$
\varphi_S^*(K\{i\})\simeq (\varphi_S^*K)\{\varphi^*(i)\}.
$$
Moreover, the exchange transformation $Ex(\varphi_*,\otimes)$
induces an exchange:
$$
Ex(\varphi_*,\{i\})
 :\varphi_{*,S}(K)\{i\} \rightarrow \varphi_{*,S}\big(K\{\varphi^*(i)\}\big).
$$
When this transformation is an isomorphism for any twist $i \in \tau$,
 we say that \emph{$\varphi_*$ commutes with twists.}
\index{word}{twist!commutes with $\tau$-twists (\emph{or} twists)}
\end{rem}
 
Note finally that Lemma \ref{lm:generators&adjoints}
 allows to prove, as for Proposition \ref{prop:carac_monoidal_P-fibred},
  the following useful lemma:
\begin{lm} \label{lm:carac_morph_monoidal_P-fibred}
Consider two complete monoidal $\Pmor$-fibred categories $\M$, $\M'$
 and denote by $M(-,\M)$ and $M(-,\M')$ their respective geometric sections.
Let $\varphi^*:\M \rightarrow \M'$ be a Cartesian functor
such that
\begin{enumerate}
\item For any scheme $S$, $\varphi^*_S:\M(S) \rightarrow \M'(S)$ is monoidal.
\item For any scheme $S$, $\varphi^*_S$ admits a right adjoint 
 $\varphi_{*,S}$.
\end{enumerate}
Assume $\M$ (resp. $\M'$) is $\tau$-generated (resp. $\tau'$-twisted)
\index{word}{generated!$\tau$-generated}
 and that $\varphi^*$ induces a surjective map from the set of isomorphism
 classes of $\tau$-twists to the set of isomorphism
 classes of $\tau'$-twists.
Then the following conditions are equivalent:
\begin{enumerate}
\item[(i)] $\varphi^*$ is a morphism of complete monoidal $\Pmor$-fibred categories.
\item[(ii)] For any object $X$ of $\Pmorx S$, the exchange transformation
 (\textit{cf.} \ref{num:Pfibred_morph})
$$
\varphi^* M_S(X,\M) \rightarrow M_S(X,\M')
$$
is an isomorphism.
\end{enumerate}
\end{lm}

%%%%%%%%%%%%%%%%%%%%%%%%%%%%
% structures of $\Pmor$

\subsection{Structures on $\Pmor$-fibred categories}

\subsubsection{Abstract definition}

\begin{num} \label{num:C_structured}
We fix a sub-$2$-category $\C$ of $\cat$
with the following properties\footnote{See the following sections for examples.}:
\begin{itemize}
\item[(1)] the $2$-functor
$$\cat\to\cat'\ , \quad A\mapsto \op{A}$$
sends $\C$ to $\C'$, where $\C'$ denotes
the $2$-category whose objects and maps are those of $\C$
and whose $2$-morphisms are the $2$-morphisms of $\C$, put in the reverse
direction.
\item[(2)] $\C$ is closed under adjunction:
for any functor $u:A\rightarrow B$ in $\C$, if
a functor $v:B\to A$ is a right adjoint or a left adjoint to $u$, then
$v$ is in $\C$.
\item[(3)] the $2$-morphisms of $\C$ are closed by transposition:
if
$$u:A\rightleftarrows B:v\ \text{and} \ u':A\rightleftarrows B:v'$$
are two adjunctions in $\C$ (with the left adjoints on the left-hand side),
a natural transformation $u\rightarrow u'$ is in $\C$ if and only if
the corresponding natural transformation $v'\rightarrow v$ is in $\C$.
\end{itemize}

We can then define and manipulate $\C$-structured
$\Pmor$-fibred categories as follows.
\begin{df} \label{df:C-structures_P-fibred}
A \emph{$\C$-structured $\Pmor$-fibred category}
(resp. \emph{$\C$-structured complete $\Pmor$-fibred category})
$\M$ over $\site$ is simply a $\Pmor$-fibred category
(resp. a complete $\Pmor$-fibred category) whose underlying
$2$-functor $\M:\op{\site}\rightarrow\cat$ factors through $\C$.
\end{df}
If $\M$ and $\M'$ are $\C$-structured fibred categories over $\site$,
a Cartesian functor $\M\rightarrow \M'$ is \emph{$\C$-structured}
if the functors $\M(S)\rightarrow \M'(S)$ are in $\C$ for any object
$S$ of $\site$, and if all the structural $2$-morphisms \eqref{num:Pfibred_morphbis}
are in $\C$ as well.
\begin{df} \label{df:morph_C-structures_P-fibred}
A morphism of $\C$-structured $\Pmor$-fibred categories
(resp. $\C$-struct\-ured complete $\Pmor$-fibred categories) is a morphism
of $\Pmor$-fibred categories (resp. of complete $\Pmor$-fibred categories)
which is $\C$-structured as a Cartesian functor.
\end{df}
\end{num}

\begin{num} \label{num:monoidal_C_structured}
Consider a $2$-category $\C$ as in the paragraph \ref{num:C_structured}.
In order to deal with the monoidal case,
 we will consider also a sub-$2$-category $\C^\otimes$ of $\C$
 such that:
\begin{enumerate}
\item The objects of $\C^\otimes$ are objects of $\C$
 equipped with a symmetric monoidal structure;
\item the $1$-morphisms of $\C^\otimes$ are exactly the $1$-morphisms of $\C$
 which are symmetric monoidal as functors;
\item the $2$-morphisms of $\C^\otimes$ are exactly the $2$-morphisms of $\C$
 which are symmetric monoidal as natural transformations.
\end{enumerate}
Note that $\C^\otimes$ satisfies condition (1) of \ref{num:C_structured},
but it does not satisfy conditions (2) and (3) in general.
Instead, we get the following properties:
\begin{enumerate}
\item[(2$'$)] If $u:A \rightarrow B$ is a functor in $\C^\otimes$,
 a right (resp. left) adjoint $v$ is a lax\footnote{For any object
 $a$, $a'$ in $A$, $F$ is lax if there exists a structural map
 $F(a) \otimes F(a') \xrightarrow{(1)} F(a \otimes a')$ 
 satisfying coherence relations (see \cite[XI. 2]{MLa}). Colax is defined by reversing
 the arrow $(1)$.} (resp. colax) monoidal functor in $\C$.
\item[(3$'$)] Consider adjunctions
$$u:A\rightleftarrows B:v\ \text{and} \ u':A\rightleftarrows B:v'$$
in $\C$ (with the left adjoints on the left-hand side).
If  $u \rightarrow u'$ (resp. $v \rightarrow v'$) is a $2$-morphism in $\C^\otimes$ 
 then $v \rightarrow v'$ (resp. $u \rightarrow u'$) is a $2$-morphism in $\C$
 which is a symmetric monoidal transformation of lax (resp. colax) monoidal functors.
\end{enumerate}
We thus adopt the following definition:
%\footnote{The reader who suspects `un crime de
%l\`ese-Bourbaki' here is right: the definition of $\C^\otimes$ is
%in some sense a particular case of the definition of $\C$-structured
%fibred category. Indeed, symmetric monoidal categories can be described as
%the fibred categories over Segal's category $\Gamma=\op{F}$ (where $F$
%denotes the category of pointed finite sets) which satisfy the so called
%Segal condition. The $2$-category $\C^\otimes$ is thus the
%sub-$2$-category of $\catm$ made of symmetric monoidal
%categories which are $\C$-structured as fibred categories over $\Gamma$,
%whose $1$-morphisms are the monoidal functors which are $\C$-structured
%as cartesian functors over $\Gamma$, etc.}.
\end{num}
\begin{df} \label{df:C-structures_monoidal_P-fibred}
A \emph{$(\C,\C^\otimes)$-structured monoidal $\Pmor$-fibred category}
(resp. a \emph{$(\C,\C^\otimes)$-structured complete monoidal $\Pmor$-fibred category})
is simply a monoidal $\Pmor$-fibred category
 (resp. a complete monoidal $\Pmor$-fibred category) 
 whose underlying $2$-functor
$\M:\op{\site}\rightarrow \catm$ factors through $\C^\otimes$.
Morphisms of such objects are defined in the same way.
\end{df}
Note that, with the hypothesis made on $\C$,
 all the exchange natural transformations defined in the preceding paragraphs lie in $\C$
 and satisfy the appropriate coherence property with respect to the monoidal structure.

\subsubsection{The abelian case} \label{sec:P-fibred_abelian}

\begin{paragr} \label{num:P-fibred_abelian}
Let $\ab$ be the sub-$2$-category of $\cat$
 made of the abelian categories,
 with the additive functors as $1$-morphisms,
 and the natural transformations as $2$-morphisms.
Obviously, it satisfies properties of \ref{num:C_structured}.
When we will apply one of the definitions \ref{df:C-structures_P-fibred},
 \ref{df:morph_C-structures_P-fibred} 
 to the case $\C=\ab$,
 we will use the simple adjective \emph{abelian} for $\ab$-structured.
\index{word}{fibred!Pfibredcategory@$\Pmor$-fibred category!abelian}
 This allows speaking of
 \emph{morphisms of abelian $\Pmor$-fibred categories.}
\index{word}{morphism!of abelian $\Pmor$-fibred categories}

Let $\abm$ be the sub-$2$-category of $\ab$
 made of the abelian monoidal categories,
 with $1$-morphisms the symmetric monoidal additive functors
 and $2$-morphisms the symmetric monoidal natural transformations.
 It satisfies the hypothesis of paragraph \ref{num:monoidal_C_structured}.
When we will apply definition \ref{df:C-structures_monoidal_P-fibred}
 to the case of $(\ab,\abm)$,
 we will use the simple expression \emph{abelian monoidal}
\index{word}{fibred!Pfibredcategory@$\Pmor$-fibred category!abelian monoidal}
 for $(\ab,\abm)$-structured monoidal.
 This allows speaking of
 \emph{morphisms of abelian monoidal $\Pmor$-fibred categories.}
\index{word}{morphism!of abelian monoidal $\Pmor$-fibred categories}
\end{paragr}

\begin{lm} \label{lm:complete&abelian}
Consider an abelian $\Pmor$-fibred category $\A$
 such that for any object $S$ of $\site$,
 $\A(S)$ is a Grothendieck abelian category.
Then the following conditions are equivalent:
\begin{enumerate}
\item[(i)] $\A$ is complete.
\item[(ii)] For any morphism $f:T \rightarrow S$ in $\site$,
 $f^*$ commutes with sums.
\end{enumerate}
If in addition, $\A$ is monoidal, the following conditions are equivalent:
\begin{enumerate}
\item[(i$'$)] $\A$ is monoidal complete.
\item[(ii$'$)]
\begin{enumerate}
\item[(a)] For any morphism $f:T \rightarrow S$ in $\site$,
 $f^*$ is right exact.
\item[(b)] For any object $S$ of $\site$,
 the bifunctor $\otimes_S$ is right exact.
\end{enumerate}
\end{enumerate}
\end{lm}

In view of this lemma, we adopt the following definition:
\begin{df} \label{df:P-fibred_Grothendieck_abelian}
A Grothendieck abelian
\index{word}{fibred!Pfibredcategory@$\Pmor$-fibred category!Grothendieck abelian}
 (resp. Grothendieck abelian monoidal)
\index{word}{fibred!Pfibredcategory@$\Pmor$-fibred category!Grothendieck abelian monoidal}
 $\Pmor$-fibred category $\A$ over $\site$
 is an abelian $\Pmor$-fibred category 
 which is complete (resp. complete monoidal)
 and such that for any scheme $S$,
  $\A(S)$ is a Grothendieck abelian category.
\end{df}

\begin{rem} Let $\A$ be a Grothendieck abelian monoidal $\Pmor$-fibred category.
Conventionally, we will denote by $\Mab S - \A$ its geometric sections.
Note that if $\A$ is $\tau$-twisted,
 then any object of $\A$ is a quotient of a direct sum of objects
 of shape $\Mab S X \A\{i\}$ for a $\Pmor$-object $X/S$ and a twist $i\in \tau$.
 \end{rem}
 
\begin{num}\label{num:df:abelian_P-premotivic_well&compactly_generated}
Consider an abelian category $\A$ which admits small sums.
Recall the following definition: \\
An object $X$ of $\T$ is \emph{finitely presented}
\index{word}{finitely presented!object of a category}
 if the functor $\Hom_\T(X,-)$ commutes with small filtering colimits.
A essentially small $\cG$ of objects of $\A$ is called generating
 if for any object $A$ of $\A$ there exists an epimorphism
 of the form:
$$
\bigoplus_{i \in I} G_i \rightarrow A
$$
where $(G_i)_{i \in I}$ is a family of objects if $\cG$.
\end{num}

\begin{df} \label{df:abelian_P-premotivic_well&compactly_generated}
Let $\A$ be an abelian $\Pmor$-fibred category over $\site$.

Given a set of twists $\tau$ of $\A$,
 we say $\A$ is \emph{finitely $\tau$-presented}
\index{word}{finitely presented!finitely $\tau$-presented}
\index{word}{presented|seealso{finitely presented}}
\index{word}{fibred!Pfibredcategory@$\Pmor$-fibred category!finitely $\tau$-presented
 |see{finitely presented}}
 if for any object $S$ of $\site$,
 for any $\Pmor$-object $X/S$ and any twist $i \in \tau$,
 the object $M_S(X)\{i\}$ is finitely presented and the class
 of such objects form an essentially small generating family of $\A(S)$.
\end{df}

%\begin{rem}
%There are weaker finiteness conditions on an object $A$
% of an abelian category $\A$. From the weaker to the stronger:
%\begin{itemize}
%\item (Gabriel) $A$ is \emph{noetherian} if there is no infinite
% strictly increasing chain of sub-ojects of $A$.
%\item $A$ is \emph{finitely generated} if for any filtering family $(B_i)_{i \in I}$
% of sub-objects of $A$ such that $A=\sum_{i \in I} B_i$,
% there exists $i \in I$ such that $A=B_i$.
%\item $A$ is compact if the functor $\Hom_\A(A,-)$ commutes
% with small direct sums.
%\end{itemize}
%When 
%\end{rem}

\subsubsection{The triangulated case}

\begin{paragr} \label{num:P-fibred_triangulated}
Let $\tri$ be the sub-$2$-category of $\cat$
 made of the triangulated categories,
 with the triangulated functors as $1$-morphisms,
 and the triangulated natural transformations as $2$-morphisms.
Then $\tri$ satisfies the properties of \ref{num:C_structured}
 (property (2) can be found for instance in \cite[Lemma~2.1.23]{ayoub},
 and we leave property (3) as an exercise for the reader).
When we will apply one of the definitions \ref{df:C-structures_P-fibred},
 \ref{df:morph_C-structures_P-fibred}
 to the case $\C=\tri$, 
 we will use the simple adjective \emph{triangulated}
\index{word}{fibred!Pfibredcategory@$\Pmor$-fibred category!triangulated}
 for $\tri$-structured.
 This allows speaking of
 \emph{morphisms of triangulated $\Pmor$-fibred categories}.
\index{word}{morphism!of triangulated $\Pmor$-fibred categories}

Let $\trim$ be the sub-$2$-category of $\tri$
 made of the triangulated monoidal categories,
 with $1$-morphisms the symmetric monoidal triangulated functors
 and $2$-morphisms the symmetric monoidal natural transformations.
 It satisfies the hypothesis of paragraph \ref{num:monoidal_C_structured}.
When we will apply definition \ref{df:C-structures_monoidal_P-fibred}
 to the case of $(\tri,\trim)$,
 we will use the expression \emph{triangulated monoidal}
\index{word}{fibred!Pfibredcategory@$\Pmor$-fibred category!triangulated monoidal}
  for $(\tri,\trim)$-structured monoidal.
 This allows speaking of
 \emph{morphisms of triangulated monoidal $\Pmor$-fibred categories}.
\index{word}{morphism!of triangulated monoidal $\Pmor$-fibred categories}
\end{paragr}

\begin{conv} \label{conv:twists&triangulated}
The set of twists
\index{word}{twist!of a triangulated monoidal $\Pmor$-fibred category}
 of a triangulated monoidal $\Pmor$-fibred category $\T$
 will always be of the form $\ZZ \times \tau$, by which
 we mean that $\tau$ is a set of twists, while $\ZZ \times \tau$
 is the closure of $\tau$ by suspension functors $[n]$, $n\in \ZZ$.
 In the notation, we shall often make the abuse of only indicating $\tau$.
 In particular, the expression $\T$ is $\tau$-generated
 will mean conventionally that $\T$ is $(\ZZ \times \tau)$-generated
 in the sense of definition \ref{df:generating_twists}. 
\end{conv}

\begin{num}\label{num:def:compgenPfibredtriang}
Consider a triangulated category $\T$ which admits small sums.
Recall the following definitions: \\
An object $X$ of $\T$ is called \emph{compact}
\index{word}{compact}
 if the functor $\Hom_\T(X,-)$ commutes with small sums.
A class $\cG$ of objects of $\T$ is called generating
 if the family of functor $\Hom_\T(X[n],-)$, $X \in \cG$, $n\in\ZZ$,
 is conservative. \\
The triangulated category $\T$ is called \emph{compactly generated}
\index{word}{generated!compactly generated}
 if there exists a generating set $\cG$ of compact objects of $\T$.  %%\\
This property of being compact has been generalized by A. Neeman
 to the property of being \emph{$\alpha$-small}
 for some cardinal $\alpha$ (\textit{cf.} \cite[4.1.1]{Nee1}) ---
 recall compact=$\aleph_0$-small.
Then the property of being compactly generated has been
 generalized by Neeman to the property of being \emph{well generated};
\index{word}{generated!well generated}
see \cite{Krause} for a convenient characterization
of well generated triangulated categories.
\end{num}

\begin{df} \label{df:P-premotivic_well&compactly_generated}
Let $\T$ be a triangulated $\Pmor$-fibred category over $\site$.
We say that $\T$ is \emph{compactly generated}
\index{word}{generated!compactly generated!triangulated $\Pmor$-fibred}
 (resp. \emph{well generated}\index{word}{generated!well generated!triangulated $\Pmor$-fibred})
if for any object $S$ of $\site$, $\T(S)$ admits small sums
and is compactly generated (resp. well generated).

Given a set of twists $\tau$ for $\T$,
 we say $\T$ is \emph{compactly $\tau$-generated}
\index{word}{generated!compactly $\tau$-generated!triangulated $\Pmor$-fibred}
 if it is compactly generated in the above sense
 and for any $\Pmor$-object $X/S$,
  any twist $i \in \tau$, $M_S(X)\{i\}$ is compact.
\end{df}

\begin{paragr}\label{wellgeneratedloc}
For a triangulated category $\T$ which has small sums, given
a family $\cG$ of objects of $\T$, we denote by $\langle \cG \rangle$
the localizing subcategory of $\T$ generated by $\cG$, i.e. $\langle \cG \rangle$
is the smallest triangulated full subcategory of $\T$ which is stable by
small sums and which contains all the objects in $\cG$.
Recall that, in the case $\T$ is well generated (e.g. if $\T$ compactly
generated), then the family $\cG$ generates $\T$
(in the sense that the family of functors
$\{\Hom_\T(X,-)\}_{X\in \cG}$ is conservative)
if and only if $\T=\langle \cG \rangle$.
The following lemma is a consequence of \cite{Nee1}:
\end{paragr}
\begin{lm} \label{lm:well_gen_tri_Pfibred}
Let $\T$ be a triangulated monoidal $\Pmor$-fibred category over $\site$
 with geometric sections $M$. Assume $\T$ is $\tau$-generated.

If $\T$ is well generated, then for any object $S$ of $\site$,
$$
\T(S)=\langle M_S(X)\{i\} ; X/S\text{ a $\Pmor$-object}, i \in \tau \rangle
$$
Moreover, there exists a regular cardinal $\alpha$ such that
all the objects of shape $M_S(X)\{i\}$ are $\alpha$-compact.
\end{lm}

%\begin{df} \label{df:triangulated_geometric_premotives}
%Let $\T$ be a triangulated monoidal $\Pmor$-fibred category.
%Denote by $M$ its geometric sections.
%
%If $\T$ is $\tau$-generated,
% we denote by $\T_{\tau-gm}(S)$
% the smallest thick triangulated sub-category of $\T(S)$ 
% containing the premotives $M_S(X)\{i\}$ for a $\Pmor$-morphism $X \rightarrow S$
% and a twist $i \in \tau$.
%The objects of $\T_{\tau-gm}(S)$ are called the $\tau$-geometric premotives.
%\end{df}
%When the set of twists $\tau$ is clear or underlying to $\T$,
% we use the notation $\T_{gm}(S)$ and the terminology of geometric
% premotives for its object.

Note finally that the Brown representability theorem
\index{word}{Brown representability theorem}
 of Neeman (\textit{cf.} \cite{Nee1})
 gives the following lemma (analog of \ref{lm:complete&abelian}):
\begin{lm}\label{lm:complete&triangulated}
Consider a well generated triangulated  $\Pmor$-fibred category $\T$.
Then the following conditions are equivalent:
\begin{enumerate}
\item[(i)] $\T$ is complete.
\item[(ii)] For any morphism $f:T \rightarrow S$ in $\site$,
 $f^*$ commutes with sums.
\end{enumerate}
If in addition, $\T$ is monoidal, the following conditions are equivalent:
\begin{enumerate}
\item[(i$'$)] $\T$ is monoidal complete.
\item[(ii$'$)]
\begin{enumerate}
\item[(a)] For any morphism $f:T \rightarrow S$ in $\site$,
 $f^*$ is right exact.
\item[(b)] For any object $S$ of $\site$,
 the bifunctor $\otimes_S$ is right exact.
\end{enumerate}
\end{enumerate}
\end{lm}

We finish this section with a proposition which will constitute
 a useful trick:
\begin{prop}\label{prop:exist_right_adjoint}
Consider an adjunction of triangulated categories
$$
a:\T \rightleftarrows \T':b.
$$
Assume that $\T$ admits a set of compact generators $\cG$
 such that any object in $a(\cG)$ is compact in $\T'$.
Then $b$ commutes with direct sums.
If in addition $\T'$ is well generated
then $b$ admits a right adjoint.
\end{prop}
\begin{proof}
The second assertion follows from the first one
 according to a corollary of the Brown representability theorem
\index{word}{Brown representability theorem}
 of Neeman (\textit{cf.} \cite[8.4.4]{Nee1}).

For the first one, we consider a family $(X_i)_{i \in I}$
 of objects of $\T'$ and prove that the canonical morphism
$$
\oplus_{i \in I} b(X_i) \rightarrow b\left(\oplus_{i \in I} X_i\right)
$$
is an isomorphism in $\T$.
To prove this, it is sufficient to apply the functor $\Hom_{\T}(G,-)$
for any object $G$ of $\cG$.
Then the result is obvious from the assumptions.
\end{proof}

We shall often use the following standard argument to produce
equivalences of triangulated categories.

\begin{cor}\label{equivgeneratorscompactgentriang}
Let $a:\T\To\T'$ be a triangulated functor between triangulated categories.
Assume that the functor $a$ preserves small sums, and that $\T$
admits a small set of compact generators $\cG$, such that $a(\cG)$
form a family of compact objects in $\T'$.
Then $a$ is fully faithful if and only if, for any
couple of objects $G$ and $G'$ in $\cG$, the map
$$\Hom_\T(G,G'[n])\To\Hom_{\T'}(a(G),a(G')[n])$$
is bijective for any integer $n$.
If $a$ is fully faithful, then $a$ is an equivalence of categories
if and only if $a(\cG)$ is a generating family in $\T'$.
\end{cor}

\begin{proof}
Let us prove that this is a sufficient condition.
As $\T$ is in particular well generated,
by the Brown representability theorem,
\index{word}{Brown representability theorem}
 the functor $b$
admits a right adjoint $b:\T'\to \T$. By virtue of
the preceding proposition, the functor $b$
preserves small sums. Let us prove that $a$ is fully faithful.
We have to check that, for any object $M$ of $\T$, the map
$M\To b(a(M))$ is invertible. As $a$ and $b$ are triangulated and preserve small sums,
it is sufficient to check this when $M$ runs over a generating family of objects of $\T$
(e.g. $\cG$). As $\cG$ is generating, it is sufficient to prove that the map
$$\Hom_\T(G,M[n])\To\Hom_{\T'}(a(G),a(M)[n])=\Hom_{\T'}(a(G),b(a(M))[n])$$
is bijective for any integer $n$, which hold then by assumption.
The functor $a$ thus identifies $\T$ with the localizing subcategory of $\T'$
generated by $a(\cG)$; if moreover $a(\cG)$ is a generating family in $\T'$,
then $\T'=\langle a(\cG) \rangle$, which also proves the last assertion.
\end{proof}

\subsubsection{The model category case} \label{sec:P-fibred_model}

\begin{paragr} \label{num:P-fibred_model}
We shall use Hovey's book \cite{Hovey} for a general
reference to the theory of model categories. Note that, following
\emph{loc.~cit.}, all the model categories we shall consider
will have small limits and small colimits.

Let $\mathscr{Mod}$ be the sub-$2$-category of $\cat$ made of
 the model categories, 
 with $1$-morphisms the left Quillen functors
 and $2$-morphisms the natural transformations.
When we will apply definition \ref{df:C-structures_P-fibred}
 (resp. \ref{df:morph_C-structures_P-fibred}) to $\C=\mathscr{Mod}$,
 we will speak of a \emph{$\Pmor$-fibred model category}
\index{word}{fibred!Pfibredcategory@$\Pmor$-fibred category!model}
 for a $\mathscr{Mod}$-structured $\Pmor$-fibred category $\M$
 (resp. morphism of $\Pmor$-fibred model categories).
\index{word}{morphism!of $\Pmor$-fibred model categories}
Note that according to the definition of left Quillen functors,
 $\M$ is then automatically complete.

Given a property $(P)$ of model categories (like being
cofibrantly generated, left and/or right proper, combinatorial, stable, etc),
we will say that a $\Pmor$-fibred model category $\M$ over $\site$ has the property
$(P)$ if, for any object $S$ of $\site$, the model
category $\M(S)$ has the property $(P)$.

For the monoidal case, we let $\mathscr{Mod}^\otimes$ be
 the sub-$2$-categories of $\mathscr{Mod}$ made of the symmetric monoidal
 model categories (see \cite[Definition 4.2.6]{Hovey}),
 with $1$-morphisms the symmetric monoidal left Quillen functors
 and $2$-morphisms the symmetric monoidal natural transformations,
 following the conditions of \ref{num:monoidal_C_structured}.
When we will apply definition \ref{df:C-structures_monoidal_P-fibred}
 to the case of $(\mathscr{Mod},\mathscr{Mod}^\otimes)$,
 we will speak simply of a \emph{monoidal $\Pmor$-fibred model category} 
\index{word}{fibred!monoidalPfibredcategory@monoidal $\Pmor$-fibred category!model (category)}
 (resp. \emph{morphism of monoidal $\Pmor$-fibred model categories})
\index{word}{morphism!of monoidal $\Pmor$-fibred model category}
 for a (resp. morphism of) $(\mathscr{Mod},\mathscr{Mod}^\otimes)$-structured monoidal
 $\Pmor$-fibred category $\M$.
Again, $\M$ is then monoidal complete.
\end{paragr}

\begin{rem}\label{remiminvPmorphism}
Let $\M$ be a $\Pmor$-fibred model category over $\site$.
Then for any $\Pmor$-morphism $p:X\To Y$,
 the inverse image functor $p^*:\M(Y)\To\M(X)$ has very strong exactness properties: 
 it preserves small limits and colimits 
 (having both a left and a right adjoint), 
 and it preserves weak equivalences, cofibrations, and fibrations. 
 The only non (completely) trivial assertion here is 
 about the preservation of weak equivalences. 
For this, one notices first that it preserves trivial cofibrations 
 and trivial fibrations (being both a left Quillen functor and a right Quillen functor). 
 In particular, 
 by virtue of Ken~Brown Lemma \cite[Lemma 1.1.12]{Hovey}, 
 it preserves weak equivalences between cofibrant (resp. fibrant) objects.
Given a weak equivalence $u:M\To N$ in $\M(Y)$,
 we can find a commutative square
$$\xymatrix{
{M'}\ar[r]^{u'}\ar[d]&{N'}\ar[d]\\
M\ar[r]_u&N}$$
in which the two vertical maps are trivial fibrations, and where $u'$ is a weak equivalence
between cofibrant objects, from which we deduce easily that $p^*(u)$ is a weak
equivalence in $\M(X)$.
\end{rem}

\begin{num} \label{num:P-fibred_model_der_functors}
Consider a $\Pmor$-fibred model category $\M$ over $\site$.
By assumption, we get the following pairs of adjoint functors:
\begin{itemize}
\item[(a)] For any morphism $f:X\To S$ of $\site$,
$$
\derL f^*:\ho(\M(S))\rightleftarrows\ho(\M(X)):\derR f_*
$$
\item[(b)] For any $\Pmor$-morphism $p:T\To S$, the pullback functor
$$
\derL p_\sharp:\ho(\M(S))\rightleftarrows\ho(\M(T)):\derL p^*=p^*=\derR p^*
$$
\end{itemize}
Moreover, the canonical isomorphism of shape $(fg)^*\simeq g^*f^*$
 induces a canonical isomorphism $\derR(fg)^*\simeq \derR g^*\derR f^*$.
In the situation of the $\Pmor$-base change formula
 \ref{defPBC}, we obtain also that the base change map
$$
\derL q_\sharp \derL g^* \rightarrow \derL f^* \derL p_\sharp
$$
is an isomorphism from the equivalent property of $\M$.
Thus, we have defined a complete $\Pmor$-fibred category whose
 fiber over $S$ is $\ho(\M(S))$.
\end{num}
\begin{df}\label{df:hoPfibred}
Given a $\Pmor$-fibred model category $\M$ as above,
 the complete $\Pmor$-fibred category defined above
 will be denoted by $\ho(\M)$ and
 called the \emph{homotopy $\Pmor$-fibred category}
\index{word}{fibred!Pfibredcategory@$\Pmor$-fibred category!homotopy}
 associated with $\M$.
\end{df}

\begin{num} Assume that $\M$ is a monoidal $\Pmor$-fibred model category over $\site$.
Then, for any object $S$ of $\site$, $\ho(\M)(S)$ has the structure
 of a symmetric closed monoidal category;
  see \cite[Theorem 4.3.2]{Hovey}.
The (derived) tensor product of $\ho(\M)(S)$ will be denoted by $M\otimes^\derL_S N$,
 and the (derived) internal Hom will be written $\derR\uHom_S(M,N)$,
 while the unit object will be written $\unit_S$.

For any morphism $f:T \rightarrow S$ in $\site$,
 the derived functor $\derL f^*$ is symmetric monoidal as follows from
 the equivalent property of its counterpart $f^*$.

Moreover, for any $\Pmor$-morphism $p:T\To S$
and for any object $M$ in $\ho(\M)(T)$ and any object $N$ in $\ho(\M)(S)$,
 the exchange map of \ref{num:exchanges3}
\begin{equation*}
\derL p_\sharp(M\otimes^\derL p^*(N))\To \derL p_\sharp(M)\otimes^\derL N
\end{equation*}
is an isomorphism.
\end{num}

\begin{df}\label{df:monoidal_hoPfibred}
Given a monoidal $\Pmor$-fibred model category $\M$ as above,
 the complete monoidal $\Pmor$-fibred category defined above
 will be denoted by $\ho(\M)$
 and called the \emph{homotopy monoidal $\Pmor$-fibred category}
\index{word}{fibred!Pfibredcategory@$\Pmor$-fibred category!homotopy monoidal}
 associated with $\M$.
\end{df}

%%%%%%%%%%%%%%%%%%%%%%%%%%%%%%%%%%
\subsection{Premotivic categories} \label{sec:premotivic_convention}

In the present article, we will focus on a particular type of $\Pmor$-fibred category.
%We introduce the following terminology: \\
%
%\begin{df}\label{defadeqschemes}
%An \emph{adequate category of schemes} is
% a full subcategory $\sch$ of the category of noetherian schemes
% satisfying the following properties:
%\begin{itemize}
%\item[(a)] $\sch$ is closed under finite limits and finite sums;
%\item[(b)] for any scheme $S$ in $\sch$, any quasi-projective $S$-scheme
%belongs to $\sch$;
%\item[(c)] any separated morphism $f:Y\rightarrow X$ in $\sch$,
%admits a compactification in $\sch$ in the sense of \cite[3.2.5]{SGA4},
%i.e. admits a factorization of the form
%$$Y \xrightarrow j \bar Y \xrightarrow p X$$
%where $j$ is an open immersion, $p$ is proper, and $\bar Y$
%belongs to $\sch$.
%\end{itemize}
%\end{df}
%
%\begin{ex}
%The category of quasi-projective $B$-schemes is obviously
%an adequate category of $B$-schemes.
%\end{ex}
%
%\begin{ex}
%The category of all noetherian $B$-schemes of finite dimension
%is adequate. Similarly, the category of separated $B$-schemes of finite type
%is adequate. In these two examples, the only non trivial condition is (c):
%but, by virtue of Nagata's Theorem\footnote{Here we use
%the fact that we work with noetherian schemes. However,
%the hypothesis that the schemes in $\S$ are
%quasi-separated quasi-compact would be enough.}~\cite{Conrad},
%a $B$-morphism of finite type is separated if and only if it admits
%a compactification.
%\end{ex}

\begin{num} Let $\base$ be a scheme. 
Assume $\sch$ is a full subcategory of the category
 of $\base$-schemes. In most of this work,
  we will denote by $\sft$ the class of morphisms of finite type in $\S$
  and by $\sm$ be the class of smooth morphisms of finite type in $\sch$.
There is an exception to this rule: throughout Part 3,
 $\sft$ will be the class of separated morphisms of finite type in $\S$
  and $\sm$ will be the class of separated smooth morphisms of finite type
  in $\sch$.
 However, the axiomatic which we will present in the sequel
  can be applied identically in each cases so that the reader can
  freely use the restriction that all morphisms of $\sm$ and $\sft$
  are separated.

In any case, the classes $\sm$ and $\sft$
  are admissible in $\sch$ in the sense of Paragraph \ref{defclassP}
(this is automatic, for instance, if $\sch$ is stable by pullbacks).
\end{num}

\begin{df} \label{df:general_P-premotivic}
\label{def:premotiviccat}
Let $\Pmor$ be an admissible class of morphisms in $\sch$.

A \emph{$\Pmor$-premotivic category over $\sch$}
\index{word}{premotivic!Ppremotivic@$\Pmor$-premotivic!category}
 --- or simply \emph{$\Pmor$-premotivic category} when $\sch$ is clear ---
 is a complete monoidal $\Pmor$-fibred category over $\sch$.
 A \emph{morphism of $\Pmor$-premotivic categories}
\index{word}{morphism!of $\Pmor$-premotivic categories}
  is a morphism of complete monoidal $\Pmor$-fibred categories over $\sch$.

As a particular case,
 when $\C$ is the $2$-category
 $\tri$ of triangulated categories (resp. $\ab$ of abelian categories),
 a \emph{$\Pmor$-premotivic triangulated (resp. abelian)
  category over $\sch$}
\index{word}{premotivic!Ppremotivic@$\Pmor$-premotivic!triangulated category}
\index{word}{premotivic!Ppremotivic@$\Pmor$-premotivic!abelian category}
 is a $(\C,\C^\otimes)$-structured complete monoidal
 $\Pmor$-fibred category over $\sch$
  (def. \ref{df:C-structures_monoidal_P-fibred}).
 Morphisms of $\Pmor$-premotivic triangulated
 \index{word}{morphism!of triangulated $\Pmor$-premotivic categories}
  (resp. abelian\index{word}{morphism!of abelian $\Pmor$-premotivic categories})
  categories are defined accordingly.

We will also say: \emph{premotivic} for $\sm$-premotivic
\index{word}{premotivic!category}
 and \emph{generalized premotivic}
\index{word}{premotivic!generalized ---- category}
 for $\sft$-premotivic.
 
The sections of a $\Pmor$-premotivic category will be called
 \emph{premotives}.
\index{word}{premotive}
\end{df}

\begin{ex} \label{ex:H&SH_premotivic}
Let $\sch$ be the category of noetherian schemes of finite dimension.

For such a scheme $S$, recall $\Hpt(S)$ is 
the pointed homotopy category 
\index{word}{homotopy category}
of Morel and Voevodsky; \textit{cf.} examples \ref{ex:H_weak_P-fibred},
  \ref{ex:H_complete_P-fibred}, \ref{ex:H_monoidal_P-fibred}.
%%  and \ref{ex:H_P-fibred+twists}.
Then,
 according to the fact recalled in these examples
  the $2$-functor $\Hpt$ is a geometrically generated premotivic category
  (recall Definition \ref{df:generating_twists}).

For such a scheme $S$, consider the stable homotopy category
\index{word}{stable homotopy category of schemes}
 $\SH(S)$
\index{notat}{SHS@$\SH(S)$}
 of Morel and Voevodsky (see \cite{Jar,ayoub2}).
According to \cite{ayoub2}, it defines a triangulated premotivic category
denoted by $\SH$. Moreover,
 it is compactly $(\ZZ \times \ZZ)$-generated
  in the sense of definition \ref{df:generating_twists}
  where the first factor refers to the suspension
  and the second one refers to the Tate twist
\index{word}{twist!Tate}
  (\emph{i.e.} as a triangulated premotivic category,
   it is compactly generated by the Tate twists).
\end{ex}

\begin{num} \label{num:premotivic&twists}
Let $\T$ be a $\Pmor$-premotivic triangulated category
 with geometric sections $M$
 and $\tau$ be a set of twists\index{word}{twist}
 for $\T$ (Definition \ref{df:twists}).

Recall from Convention \ref{conv:twists&triangulated} 
 (resp. and Definition \ref{df:P-premotivic_well&compactly_generated}) that
 $\T$ is said to be $\tau$-generated\index{word}{generated!$\tau$-generated}
 (resp. compactly $\tau$-generated)\index{word}{generated!compactly $\tau$-generated}
 if for any scheme $S$,
 the family of isomorphism of classes of premotives
  of the form $M_S(X)\{i\}$ 
  for a $\Pmor$-scheme $X$ over $S$ and a twist $i \in \tau$
  is a set of generators (resp. compact generators)
 for the triangulated category $\T(S)$
 (in the respective case, we also assume $\T(S)$ admits small sums).

Let $E$ be a premotive over $S$ and $X$ be a $\Pmor$-scheme over $S$.
For any $(n,i) \in \ZZ \times \tau$, we define the cohomology of $X$ 
 in degree $n$ and twist $i$ with coefficients in $E$ as:
$$
H^{n,i}_\T(X,E)=\Hom_{\T(S)}\big(M_S(X),E\{i\}(n)\big).
$$
The fact $\T$ is $\tau$-generated amounts to say that
 any such premotive $E$ is determined by its cohomology.
\end{num}

\begin{ex} \label{ex:SH_premotivic_gen}
All the known triangulated premotivic categories are $\tau$-generated
 for a given set of twist $\tau$. In fact, one defines as usual
 the \emph{Tate twist} $\un_S(1)$ \index{word}{twist!Tate}
 in such a premotivic triangulated category $\T$ by the formula:
$$
M_S(\PP^1_S)=\un \oplus \un(1)[2].
$$
Then $\un(1)=(\un_S(1))_{S \in \base}$ is a cartesian section of $\T$.
 We will say that $\T$ is \emph{generated by Tate twists}
 if it is $\ZZ$-generated where $\ZZ$ refers to the set of twists
 $(\un(n))_{n \in \ZZ}$.

The premotivic triangulated category $SH$ of the previous example
 is compactly generated by Tate twists. Similarly,
 the stable $\AA^1$-derived category $\DMux \Rc$
 (cf. Example \ref{ex:stable_AA^1-derived_categories}),
 the category of Voevodsky motives $\DMV$
 (cf. Definition \ref{df:Nis_DMe&DM}),
 the category of $\BGL$-modules (cf. Definition \ref{df:KGL-mod})
 and the category of Beilinson motives $\DMB$
 (cf. Definition \ref{df:Beilinson_motives})
 are all compactly generated by Tate twists.
\end{ex}

\begin{df} \label{df:adjunction_premotivic}
Let $\M$ and $\M'$ be $\Pmor$-premotivic categories.

%An \emph{elementary premotivic morphism} from $\M$ to $\M'$
%is a morphism $\varphi^*:\M \rightarrow \M'$
%of twisted monoidal $\sm$-fibred (resp. $\sft$-fibred) categories.

A morphism of $\Pmor$-premotivic categories
\index{word}{morphism!of $\Pmor$-premotivic categories}
 (or simply a \emph{premotivic morphism})
\index{word}{premotivic!morphism|see{morphism of premotivic categories}}
is a morphism $\varphi^*:\M \rightarrow \M'$
of complete monoidal $\Pmor$-fibred categories.
We shall also say that
$$\varphi^*:\M \rightleftarrows \M':\varphi_*$$
is a \emph{premotivic adjunction}.
When moreover $\M$ and $\M'$ are $\Pmor$-premotivic triangulated
 (resp. abelian) categories, we will ask $\varphi^*$ 
 is a compatible with the triangulated (resp. additive)
 structure -- as in Definition \ref{df:morph_C-structures_P-fibred}.

If we assume that $\M$ (resp. $\M'$) is $\tau$-twisted
 (resp. $\tau'$-twisted), we will say as in Definition \ref{df:morphism&twists}
 that $\varphi^*$ is \emph{compatible with twists}
\index{word}{compatible with twists}
 if for any $i \in \tau$,
 $\varphi^*(i)$ belongs up to isomorphism to $\tau'$.
We say $\varphi^*$ is \emph{strictly compatible with twists}
 if it is compatible with twists and if any element of $\tau'$ is isomorphic to the image of an element of $\tau$.
\end{df}
Usually, premotivic categories comes equip with canonical twists
 (especially the Tate twist, see the above example)
 and premotivic morphisms are compatible with twists.

\begin{ex} \label{ex:H&SH_premotivic_adjunction}
With the hypothesis and notations of \ref{ex:H&SH_premotivic},
 we get a premotivic adjunction
$$
\sus:\Hpt \rightleftarrows \SH:\lop
$$
induced by the infinite suspension functor
\index{word}{functor!infinite suspension}
\index{word}{infinite suspension|seealso{functor}}
 according to \cite{Jar}.
\end{ex}

\begin{num}\label{num:df:Tgm}
Let $\T$ (resp. $\A$) be a triangulated $\Pmor$-premotivic category
 with geometric sections $M$ and a set of twists $\tau$. 
 For any scheme $S$, we let $\T_{\tau,c}(S)$ be the smallest triangulated 
 thick\footnote{\emph{i.e.} stable by direct factors.}
 subcategory of $\T(S)$ which contains premotives of shape $M_S(S)\{i\}$ 
  (resp. $\Mab S X \A\{i\}$) for a $\Pmor$-scheme $X/S$ and a twist $i \in \tau$. 
This subcategory is stable by the operations $f^*$, $p_\sharp$ and $\otimes$.
In particular, $\T_{\tau,c}$ defines a \emph{not necessarily complete}
 triangulated (resp. abelian) $\Pmor$-fibred category over $\sch$.
 We also obtain a morphism of triangulated (resp. abelian) monoidal $\Pmor$-fibred categories,
  fully faithful as a functor,
$$\iota:\T_{\tau,c} \rightarrow \T$$
\end{num}
\begin{df} \label{df:tau-geometric}
Consider the notations introduced above.
 We will call $\T_{\tau,c}$ the \emph{$\tau$-constructible part of $\T$}.
For any scheme $S$, the objects of $\T_{\tau,c}(S)$ will be called
 \emph{$\tau$-constructible}.
\index{word}{constructible!$\tau$-constructible}
\index{word}{constructible|seealso{$\tau$-constructible}}
\end{df}
When $\tau$ is clear from the context,
 we will put $\T_c:=\T_{\tau,c}$ and use the terminology \emph{constructible}.

\begin{rem}
The condition of $\tau$-constructibility is a good categorical notion of finiteness
 which extends the notion of \emph{geometric motives}
\index{word}{motive!geometric}
 as introduced by Voevodsky. 
 In the triangulated motivic case,
  it will be studied thoroughly in section \ref{sec:constructible_motives}.
\end{rem}

\begin{prop}\label{constructequivcompact}
Let $\T$ be a $\tau$-twisted $\Pmor$-premotivic triangulated category.
Let $S$ be a scheme such that:
\begin{enumerate}
\item The category $\T(S)$ admits small sums.
\item For any $\Pmor$-scheme $X$ over $S$,
 and any twist $i \in \tau$, the premotive $M_S(X)\{i\}$ is compact.
\end{enumerate}
Then, a premotive $M$ over $S$ is $\tau$-constructible if and only if it is compact.
\end{prop}
\begin{proof}
In any triangulated category $\mathscr D$, one easily obtains that
 the property of being compact is stable under extensions and retracts.
 In particular, the thick triangulated subcategory 
 of $\mathscr D$ generated by compact objects consists precisely
 of the compact objects of $\mathscr D$.
 Moreover, if $\mathscr D$ admits small sums and
 is generated by a family of compact objects $G$,
 then the thick triangulated subcategory of $\mathscr D$ generated by $G$
 contains all compact objects, and is therefore equal to the full subcategory
 of compact objects (see \cite[Lem. 2.2]{Nee4}).

Coming back to the definition of being $\tau$-constructible,
 this general fact finishes the proof.
\end{proof}
Thus, when the conditions of this proposition are fulfilled,
 the category $\T_{\tau,c}(S)$ does not depend on the particular choice
 of $\tau$.  This will often be the case
  in practice (see \ref{cor:Der_compact&constructible},
  \ref{cor:DMue_compact&constructible}, \ref{cor:geometric_premotives}).

\begin{rem}
The notion of compact objects in a triangulated category
 was heavily developed by A.~Neeman. Its relation with finiteness
 conditions is particularly emphasized when considering the
 derived category of complexes of quasi-coherent sheaves
 over a quasi-compact separated scheme: in this triangulated category,
 being compact is equivalent to being perfect (\cite[Cor. 4.3]{Nee3}).
\end{rem}

\begin{df} \label{df:enlargement}
Consider a $\tau$-generated premotivic category $\M$.

An \emph{enlargement}
\index{word}{premotivic!enlargement of ---- category}
\index{word}{enlargement, of premotivic categories|see{premotivic}}
 of $\M$ is the data of a $\tau'$-twisted generalized premotivic category $\uM$
 together with a premotivic adjunction
$$
\rho_\sharp:\M \longrightarrow \uM : \rho^*
$$
(where $\uM$ is considered as a premotivic category in the obvious way),
satisfying the following properties:
\begin{itemize}
\item[(a)] For any scheme $S$ in $\sch$, the functor
$
\rho_{\sharp,S}:\M(S) \rightarrow \uM(S)
$
is fully faithful and its right adjoint
$
\rho^*_S:\uM(S) \rightarrow \M(S)
$
commutes with sums.
\item[(b)] $\rho_\sharp$ is strictly compatible with twists.
\end{itemize}
\end{df}
Again, this notion is defined similarly for a $\C$-structured
 $\Pmor$-premotivic category.

Note that for any smooth $S$-scheme $X$, 
we get in the context of an enlargement as above the following identifications:
\begin{align*}
\rho_{\sharp,S}(M_S(X)) &\simeq \umotNP_S(X), \\
\rho^*_S(\umotNP_S(X)) &\simeq M_S(X)
\end{align*}
where $M$ (resp. $\umotNP$) denote the geometric sections of $\M$ (resp. $\uM$).

Remember also that
 for any morphism of schemes $f$ and any smooth morphism $p$,
 $\rho_\sharp$ commutes with $f^*$ and $p_\sharp$,
 while $\rho^*$ commutes with $f_*$ and $p^*$.
\section{Triangulated $\Pmor$-fibred categories in algebraic geometry}
\label{sec:six_functors}

\begin{assumption} \label{num:assumption1_sch}
In this entire section,
 we fix a base scheme $\base$, assumed to be noetherian,
  and a full subcategory $\sch$
 of the category of noetherian $\base$-schemes
 satisfying the following properties:
\begin{itemize}
\item[(a)] $\sch$ is closed under finite sums and pullback along
 morphisms of finite type.
\item[(b)] For any scheme $S$ in $\sch$, any quasi-projective $S$-scheme
belongs to $\sch$.
\end{itemize}
In sections \ref{sec:Deligne} and \ref{sec:Ayoub},
 we will add the following assumption on $\sch$:
\begin{itemize}
\item[(c)] Any separated morphism $f:Y\To X$ in $\sch$,
admits a compactification in $\sch$ in the sense of \cite[3.2.5]{SGA4},
i.e. admits a factorization of the form
$$Y \xrightarrow j \bar Y \xrightarrow p X$$
where $j$ is an open immersion, $p$ is proper,
 and $\bar Y$ belongs to $\sch$.
Furthermore, if $f$ is quasi-projective, then $p$
 can be chosen to be projective.
\item[(d)] Chow's lemma
\index{word}{Chow's lemma}
 holds in $\sch$ 
(i.e., for any proper morphism $Y\To X$ in $\sch$, 
 there exists a projective birational morphism $p:Y_0\To Y$ in $\sch$ 
 such that $fp$ is projective as well).
\end{itemize}
A category $\sch$ satisfying all these properties 
 will be called \emph{adequate}
\index{word}{adequate, category of schemes}
 for future references.\footnote{
 For instance, the scheme $\base$ can be the spectrum of a prime field
  or of a Dedekind domain.
 The category $\sch$ might be the category of
 all noetherian $\base$-schemes of finite dimension 
 or simply the category of quasi-projective $\base$-schemes.
 In all these cases, property (c) is ensured by Nagata's theorem
 (see \cite{Conrad}) and property (d) by Chow's lemma 
 (see \cite[5.6.1]{EGA2}).}

We also fix an admissible class $\Pmor$ of morphisms in $\sch$
 and a complete triangulated $\Pmor$-fibred category $\T$.
We will add the following assumptions:
\begin{itemize}
\item[(d)] In section \ref{sec:Deligne} and \ref{sec:localization},
 $\Pmor$ contains the open immersions.
\item[(e)] In section \ref{sec:Ayoub},
 $\Pmor$ contains the smooth morphisms of $\sch$.
\end{itemize}

In the case $\T$ is monoidal, we denote by
$$
M:\Pmorx ?\To \T
$$
its geometric sections.

According to the convention of \ref{df:general_P-premotivic},
 we will speak of the \emph{premotivic case}
\index{word}{premotivic!case}
 when $\Pmor$ is the  class of smooth morphisms of finite type\footnote{
 or smooth separated morphisms of finite type when applying this section
 in Part 3}
in $\sch$ and $\T$ is a premotivic triangulated category.
%The aim of this section
% is to introduce elementary properties which will ensure the existence and
% good behavior of direct image with compact support functors and their right
% adjoints.
\end{assumption}

\subsection{Elementary properties} \label{sec:triang_Pfibred_elem_ppty}

\begin{df} \label{df:additive_triangulated}
We say that $\T$ is additive, if for any finite family $(S_i)_{i \in \I}$
 of schemes in $\sch$, the canonical map
$$
\T\left(\coprod_i S_i\right) \rightarrow \prod_i \T(S_i)
$$
is an equivalence.
\end{df}
Recall this property implies in particular that $\T(\varnothing)=0$.

\begin{lm}
Let $S$ be a scheme, $p:\AA^1_S \rightarrow S$ be the canonical projection.
The following conditions are equivalent:
\begin{enumerate}
\item[(i)] The functor $p^*:\T(S) \rightarrow \T(\AA^1_S)$ is fully faithful.
\item[(ii)] The counit adjunction morphism $1 \rightarrow p_*p^*$ is an isomorphism.
\end{enumerate}
In the premotivic case, these conditions are equivalent to the following ones:
\begin{enumerate}
\item[(iii)] The unit adjunction morphism $p_\sharp p^* \rightarrow 1$ is an isomorphism.
\item[(iv)] The morphism $M_S(\AA^1_S) \xrightarrow{p_*} \un_S$ induced by $p$
 is an isomorphism.
\item[(iv')] For any smooth $S$-scheme $X$, the
morphism $M_S(\AA^1_X) \xrightarrow{(1_X \times p)_*} M_S(X)$ is an isomorphism.
\end{enumerate}
\end{lm}
The only thing to recall is that in the premotivic case,
 $p_\sharp p^*(M)=M_S(\AA^1_S) \otimes M$ and $p_* p^*(M)=\uHom_S(M_S(\AA^1_S),M)$.

\begin{df} \label{df:ppty:htp}
The equivalent conditions of the previous lemma will be called
 the \emph{homotopy property} for $\T$, denoted by \htp.
\end{df}

\begin{num} \label{num:sieves}
Recall that a \emph{sieve}
\index{word}{sieve}
 $R$ of a scheme $X$ 
 is a class of morphisms 
 in $\sch/X$ which is stable by composition on the right by any morphism of schemes
 (see \cite[I.4]{SGA4}).

Given such a sieve $R$, we will say that \emph{$\T$ is $R$-separated}
 if the class of functors $f^*$ for $f \in R$ is conservative.
Given two sieves $R$, $R'$ of $X$, the following properties are immediate:
\begin{itemize}
\item[(a)] If $R \subset R'$ then $\T$ is $R$-separated implies $\T$ is $R'$-separated.
\item[(b)] If $\T$ is $R$-separated and is $R'$-separated then $\T$ is $(R \cup R')$-separated.
\end{itemize}
 A family of morphisms $(f_i:X_i \rightarrow X)_{i \in I}$ of schemes
 defines a sieve $R=\langle f_i, i \in I \rangle$ such that $f$ is in $R$
 if and only if there exists $i \in I$ such that $f$ can be factored through $f_i$.
Obviously,
\begin{itemize}
\item[(c)] $\T$ is $R$-separated if and only if 
 the family of functors $(f_i^*)_{i \in I}$ is conservative.
\end{itemize}
Recall that a topology on $\sch$ is the data for any scheme $X$
 of a set of sieves of $X$ satisfying certain stability conditions
  (\textit{cf.} \cite[II, 1.1]{SGA4}), called $t$-covering sieves.
A pre-topology $t_0$ on $\sch$ is the data for any scheme $X$
 of a set of families of morphisms of shape $(f_i:X_i \rightarrow X)_{i \in I}$
 satisfying certain stability conditions (\textit{cf.} \cite[II, 1.3]{SGA4}), 
 called $t_0$-covers. A pre-topology $t_0$ generated a unique topology $t$.
\end{num}
\begin{df} \label{df:ppty:sepx t}
Let $t$ be a Grothendieck topology on $\site$. We say that
 $\T$ is \emph{$t$-separated} if the following property holds:
\begin{enumerate}
\item[\sepx{$t$}] For any $t$-covering sieve $R$,
 $\T$ is $R$-separated in the sense defined above.
\end{enumerate}
\end{df}
Obviously, given two topologies $t$ and $t'$ on $\sch$ such that $t'$ is finer than
 $t$, if $\T$ is $t$-separated then it is $t'$-separated.

If the topology $t$ on $\sch$ is generated by a pre-topology $t_0$
 then $\T$ is $t$-separated if and only if for any $t_0$-covers $(f_i)_{i \in I}$,
 the family of functors $(f_i^*)_{i \in I}$ is conservative -- use \cite[1.4]{SGA4}
  and \ref{num:sieves}(a)+(c).

\begin{num} \label{num:radicial}
Recall that
 a morphism of schemes $f:T \rightarrow S$ is radicial
\index{word}{morphism!radicial}
\index{word}{radicial|see{morphism}}
 if it is injective and for any point $t$ of $T$,
 the residual extension induced by $f$ at $t$ is radicial
  (\textit{cf.} \cite[3.5.4, 3.5.8]{EGA1})\footnote{
It is equivalent to ask that $f$ is universally injective.
When $f$ is surjective,
 this is equivalent to ask that $f$ is a universal homeomorphism.
\index{word}{homeomorphism, universal}}
 The following definition is inspired by \cite[Def. 2.1.160]{ayoub}.
%We say that $f$ is a \emph{radicial extension}
% if it is finite, radicial and surjective.
\end{num}

\begin{df} \label{df:ppty:sep_ssep}
We say that $\T$ is \emph{separated} (resp. \emph{semi-separated})
 if $\T$ is separated for the topology generated by
 surjective families of morphisms of finite type 
 (resp. finite radicial morphisms) in $\site$.
We also denote by \sep (resp. \ssep) this property.
\end{df}

\begin{rem} If $\T$ is additive, property \sep (resp. \ssep) is equivalent to ask
 that for any surjective morphism of finite type 
 (resp. finite surjective radicial morphism)
 $f:T \rightarrow S$ in $\sch$, the functor $f^*$ is conservative.
\end{rem}

\begin{prop}\label{stronglyquasiseparated}
Assume $\T$ is semi-separated
 and satisfies the transversality property with respect
 to finite surjective radicial morphisms.

Then for any finite surjective radicial morphism $f:Y \rightarrow X$,
 the functor
$$
f^*:\T(X) \rightarrow \T(Y)
$$
is an equivalence of categories.
\index{word}{equivalence, of categories}
\end{prop}

\begin{proof}
We first consider the case when $f=i$ is in addition a closed immersion.
In this case, we can consider the pullback square below.
$$
\xymatrix@=20pt{
Y\ar@{=}[r]\ar@{=}[d]&Y\ar[d]^i\\
Y\ar[r]_i& Z
}
$$
Using the transversality property with respect to $i$,
 we see that the counit 
$i^* i_* \rightarrow 1$
is an isomorphism.
It thus remains to prove that the unit map
$
1 \rightarrow i_*i^*
$
is an isomorphism. 
As $i^*$ is conservative by semi-separability,
it is sufficient to check that
$$
i^* \rightarrow i^*i_*i^*(M)
$$
is an isomorphism. But this is a section of the
map $i^*\, i_*i^*(M)\To i^*(M)$, which is
already known to be an isomorphism.

Consider now the general case of a finite radicial extension $f$.
We introduce the pullback square
$$
\xymatrix@=20pt{
Y\times_X Y\ar[r]^-p\ar[d]_q& Y\ar[d]^f\\
Y\ar[r]_f& X
}
$$
Consider the diagonal immersion $i:Y \rightarrow Y \times_X Y$.
Because $Y$ is noetherian and $p$ is separable,
 $i$ is finite (\textit{cf.} \cite[6.1.5]{EGA2}) thus a closed immersion.
As $p$ is a universal homeomorphism,
 the same is true for its section $i$.
The preceding case thus implies that $i^*$ is an equivalence of categories. 
Moreover, as $pi=qi=1_Y$, we see that $p^*$
and $q^*$ are both quasi-inverses to $i^*$, which
implies that they are isomorphic equivalences of categories.
More precisely, we get canonical isomorphisms of functors
$$i^*\simeq p_*\simeq q_*
\quad\text{and}\quad
i_*\simeq p^*\simeq q^*.$$

We check that the unit map $1 \rightarrow f_*f^*$ is an isomorphism.
Indeed, by semi-separability,
 it is sufficient to prove this after applying the functor $f^*$, and we get, using
the transversality property for $f$:
$$
f^*
\simeq i^* p^* f^*
\simeq q_* p^* f^*
\simeq f^* f_* f^*.
$$
We then check that the counit map $f^* f_* \rightarrow 1$ is an isomorphism as well.
In fact, using again the transversality property for $f$,
 we have isomorphisms
$$
f^* f_*(M) \simeq q_* p^*(M) \simeq i^* i_*(M) \simeq M.
$$
\end{proof}

\begin{num} \label{df:cd_structures}
Recall from \cite{voecd1} that a cd-structure
\index{word}{cdstructure@cd-structure}
 on $\sch$ is a collection $P$ of commutative squares of schemes
$$
\xymatrix@=14pt{
B\ar[r]\ar[d]\ar@{}|Q[rd] & Y\ar^f[d] \\
A\ar_e[r] & X
}
$$
which is closed under isomorphisms.
We will say that a square $Q$ in $P$ is $P$-distinguished.
\index{word}{square!$P$-distinguished}

Voevodsky associates to $P$ a topology $t_P$,
the smallest topology such that:
\begin{itemize}
\item for any $P$-distinguished square $Q$ as above,
 the sieve generated by $\{ f:A \rightarrow X, e:Y \rightarrow X \}$
\index{word}{sieve}
 is $t_P$-covering on $X$.
\item the empty sieve covers the empty scheme.
\end{itemize}
\end{num}

\begin{ex} \label{ex:lower&upper_cd_structures}
A \emph{Nisnevich distinguished square}
\index{word}{square!Nisnevich distinguished}
\index{word}{Nisnevich!distinguished square|see{distinguished}}
 is a square $Q$
 as above such that $Q$ is cartesian, $f$ is \'etale,
 $e$ is an open embedding with reduced complement $Z$
 and the induced map $f^{-1}(Z) \rightarrow Z$ is an isomorphism.
 The corresponding cd-structure is called the \emph{upper cd-structure}
\index{word}{cdstructure@cd-structure!upper}
 (see section 2 of \cite{voecd2}).
 Because we work with noetherian schemes,
 the corresponding topology is the \emph{Nisnevich topology}
\index{word}{topology!Nisnevich}
\index{word}{Nisnevich!topology|see{topology}}
 (see proposition 2.16 of \emph{loc.cit.}).

A \emph{proper $\cdh$-distinguished square}
\index{word}{square!proper $\cdh$-distinguished}
 is a square $Q$ as above such that $Q$ is Cartesian,
 $f$ is proper,
 $e$ is a closed embedding with open complement $U$
 and the induced map $f^{-1}(U) \rightarrow U$ is an isomorphism.
 The corresponding cd-structure is called the \emph{lower cd-structure}.
\index{word}{cdstructure@cd-structure!lower}
 The topology associated with the lower cd-structure is called the
 \emph{proper cdh-topology}.
\index{word}{topology!proper $\cdh$}

The topology generated by the lower and upper cd-structures
 is by definition (according to the preceding remark on Nisnevich topology)
 the \emph{cdh-topology}.
\index{word}{topology!cdhtopology@$\cdh$-topology}

All these three examples are complete cd-structures
 in the sense of \cite[2.3]{voecd1}.
\end{ex}
 
\begin{lm}\label{lm:cd-separation}
Let $P$ be a complete cd-structure (see \cite[def 2.3]{voecd1}) on $\sch$
% whose distinguished
% squares are stable by pullbacks
% (or more generally complete in the sense of \cite[2.3]{voecd1})
 and $t_P$ be the associated topology.
The following conditions are equivalent:
\begin{enumerate}
\item[(i)] $\T$ is $t_P$-separated.
\item[(ii)] For any distinguished square $Q$ for $P$ of the above form,
 the pair of functors $(e^*,f^*)$ is conservative.
\end{enumerate}
\end{lm}
\begin{proof}
This follows from the definition of a complete cd-structure and \ref{num:sieves}(a).
\end{proof}

\begin{rem} If we assume that $\sch$ is stable by arbitrary pullback
 then any cd-structure $P$ on $\sch$ such that $P$-distinguished squares
 are stable by pullback is complete (see \cite[2.4]{voecd1}).
 \end{rem}

\subsection{Exceptional functors, following Deligne} \label{sec:Deligne}

\subsubsection{The support axiom}

\begin{num} \label{num:forget_support}
Consider an open immersion $j:U \rightarrow S$.
Applying \ref{num:exchanges2} to
 the cartesian square
$$
\xymatrix{
U\ar@{=}[r]\ar@{=}[d]& U\ar^j[d] \\
U\ar_j[r] & S
}
$$
we get a canonical natural transformation
$$
\gamma_j:j_\sharp = j_\sharp 1_*
 \xrightarrow{Ex(\Delta_{\sharp*})} j_* 1_\sharp = j_*.
$$
Recall that the functors $j_\sharp$ and $j_*$ are fully faithful
 (see Corollary \ref{immfullyfaithful}).

Note that according to remark \ref{rem:coherence_exchange1},
 this natural transformation is compatible with identifications
 of the kind $(jk)_\sharp=j_\sharp k_\sharp$ and $(jk)_*=j_*k_*$.
\end{num} 
 
\begin{lm} \label{lm:localization&sums}
Let $S$ be a scheme, $U$ and $V$ be subschemes such that $S=U \sqcup V$.
We let $h:U \rightarrow S$ (resp. $k:V \rightarrow S$) be
the canonical open immersions.

Assume that the functor $(h^*,k^*):\T(S) \rightarrow \T(U) \times \T(V)$
 is conservative and that $\T(\varnothing)=0$.
Then the natural transformation $\gamma_h$ (resp. $\gamma_k$) is an isomorphism.
Moreover, the functor $(h^*,k^*)$ is then an equivalence of categories.
\end{lm}
\begin{proof}
As $h_\sharp$ and $h_*$ are fully faithful,
we have $h^*h_\sharp\simeq h^*h_*$.
By $\Pmor$-base change, we also get $k^*h_\sharp\simeq k^*h_*\simeq 0$.
It remains to prove the last assertion.
The functor $R=(h^*,k^*)$ has a left adjoint $L$
defined by L=$h_\sharp\oplus k_\sharp$:
$$L(M,N)=h_\sharp(M)\oplus k_\sharp(N)\, .$$
The natural transformation $LR\To 1$ is an isomorphism:
to see this, is it sufficient to evaluate at $h^*$ and $k^*$,
which gives an isomorphism in $\T(U)$ and $\T(V)$ respectively.
The natural transformation $1\To RL$ is also an isomorphism
because $h_\sharp$ and $k_\sharp$ are fully faithful.
\end{proof}

\begin{rem}
Assume $\T$ is Zariski separated (definition \ref{df:ppty:sepx t}).
Then, as a corollary of this lemma, 
 $\T$ is additive (definition \ref{df:additive_triangulated}) 
 if and only if $\T(\varnothing)=0$.
\end{rem}

\begin{num} \textit{Exchange structures V}.-- \label{num:exchanges5}
Assume $\T$ is additive.
We consider a commutative square of schemes
\begin{equation} \label{squareexchangeV1}\begin{split}
\xymatrix@=16pt{
V\ar^k[r]\ar_q[d]\ar@{}|\Delta[rd] & T\ar^p[d] \\
U\ar_/-2pt/j[r] & S
}\end{split}
\end{equation}
such that $j$, $k$ are an open immersions and $p$, $q$ are proper morphisms.

This diagram can be factored into the following commutative diagram:
$$
\xymatrix@C=22pt@R=12pt{
V\ar@/^9pt/^k[rrd]\ar@/_9pt/_q[rdd]\ar|-l[rd] & & \\
& U \times_S T\ar|-{j'}[r]\ar^{p'}[d]\ar@{}|\Theta[rd] & T\ar^p[d] \\
& U\ar|j[r] & S.
}
$$
Then $l$ is an open and closed immersion so that the previous lemma
 implies the canonical morphism
$\gamma_l:l_\sharp \rightarrow l_*$ is an isomorphism.
As a consequence, we get a natural exchange transformation
$$
Ex(\Delta_{\sharp*}):j_\sharp q_* = j_\sharp p'_* l_*
 \xrightarrow{Ex(\Theta_{\sharp*})} p_* j'_\sharp l_*
 \xrightarrow{\gamma_l^{-1}} p_* j'_\sharp l_\sharp=p_* k_\sharp
$$
using the exchange of \ref{num:exchanges2}.
Note that, with the notations introduced in \ref{num:forget_support},
the following diagram is commutative.
\begin{equation}\label{squareexchangeV4}\begin{split}
\xymatrix{
j_\sharp q_* \ar[rr]^{Ex(\Delta_{\sharp*})}\ar[d]_{\gamma_j q_*}&& p_*k_\sharp\ar[d]^{p_*\gamma_k}\\
j_*q_*\ar[r]^(.35)\sim & (jq)_*=(pk)_*& p_* k_*\ar[l]_(.35)\sim
}\end{split}
\end{equation}
Indeed one sees first that it is sufficient to treat the case where $\Delta$
is cartesian. Then, as $j_\sharp$ is a fully faithful left adjoint to $j^*$ it is sufficient
to check that \eqref{squareexchangeV4} commutes after having applied $j^*$.
Using the cotransversality property with respect to open immersions, one
sees then that this consists of verifying the commutativity of \eqref{squareexchangeV4}
when $j$ is the identity, in which case it is trivial.
\end{num}

\begin{df}\label{defsupportppty}
Let $p:T \rightarrow S$ be a proper morphism in $\sch$.

We say that the triangulated $\Pmor$-fibred category $\T$ satisfies the
\emph{support property} with respect to $p$, 
 denoted by \suppx p,
 if it is additive and for any commutative square
 of shape \eqref{squareexchangeV1}
 the exchange transformation 
 $Ex(\Delta_{\sharp*}):j_\sharp q_* \rightarrow p_* k_\sharp$
 defined above is an isomorphism.

We say that $\T$ satisfies the \emph{support property},
 also denoted by \supp, if it satisfies \suppx p
 for all proper morphism $p$ in $\sch$.
\end{df}
\noindent By definition,
 it is sufficient to check the last property of property \supp
 in the case where $\Delta$ is cartesian.

\subsubsection{Exceptional direct image}

\begin{num} 
We denote by $\schsep$ (resp. $\schouv$, $\schprop$)
 the sub-category of the category $\sch$ with the same objects
 but morphisms are separated morphisms of finite type
 (resp. open immersions, proper morphisms).
We denote by 
\begin{align*}
&\T_*:\sch \rightarrow \trim \\
\text{ resp. } &\T_\sharp:\schouv \rightarrow \trim 
\end{align*}
the $2$-functor defined respectively by morphisms of type $f_*$
 and $j_\sharp$ ($f$ any morphism of schemes).
The proposition below is essentially based on a result
of Deligne \cite[XVII, 3.3.2]{SGA4}:
\end{num}
%%%%%%%%%%%%%%%%%%%%%%%%%%%%%%%%%%%%%%%%%%%%
\begin{prop} \label{prop:support&exists_f_!}
Assume $\T$ is a monoidal $\Pmor$-fibred category
 and satisfies property \supp.

Then there exists a unique $2$-functor
$$
\T_!:\schsep \rightarrow \trim
$$
with the property that 
$$
\T_!|_{\schprop} = \T_*|_{\schprop}, \quad
\T_!|_{\schouv} = \T_\sharp
$$
and for any commutative square $\Delta$ of shape \eqref{squareexchangeV1}
with $p$ and $q$ proper, the composition of the structural isomorphisms
$$
j_\sharp q_*=j_! q_!\simeq (jq)_!=(pk)_!\simeq p_! k_!=p_* k_\sharp
$$
is equal to the exchange transformation $Ex(\Delta_{\sharp*})$.
\end{prop}

\begin{num} \label{num:notation_f_!}
Under the assumptions of the proposition,
for any separated morphism\index{word}{morphism!separated}
 of finite type $f:Y \rightarrow X$,
 we will denote by $f_!:\T(Y) \rightarrow \T(X)$ the functor $\T_!(f)$.
 The functor $f_!$ is called the
  \emph{direct image functor with compact support}
\index{word}{direct image with compact support|see{functor, left exceptional}}
  or the \emph{left exceptional functor}
\index{word}{functor!left exceptional}
\index{word}{exceptional functor|see{functor}}
 associated with $f$.
\end{num}
 
\begin{proof}
We recall the principle of the proof of Deligne.
Let $f:Y \rightarrow X$ be a separated morphism of finite type
 in $\sch$.

Let $\C_f$ be the category of compactifications
of $f$ in $\sch$, i.e. of factorizations of $f$ of the form
\begin{equation} \label{pf:morph_comp00}
Y \xrightarrow j \bar Y \xrightarrow p X
\end{equation}
where $j$ is an open immersion, $p$ is proper, and $\bar Y$
belongs to $\sch$. Morphisms of $\C_f$
are given by commutative diagrams of the form
\begin{equation} \label{pf:morph_comp}\begin{split}
\xymatrix@R=-2pt@C=28pt{
 & \bar Y'\ar^{p'}[rd]\ar^\pi[dd] & \\
Y\ar^{j'}[ru]\ar_{j}[rd] & & X. \\
 & \bar Y\ar_{p}[ru] &
}\end{split}
\end{equation}
in $\sch$. To any compactification
 of $f$ of shape \eqref{pf:morph_comp00}, we associate
the functor $p_*j_\sharp$. \\
To any morphism of compactifications \eqref{pf:morph_comp},
we associate a natural isomorphism
$$
p'_*j'_\sharp=p_*\pi_*j'_\sharp
 \xrightarrow{Ex(\Delta_{\sharp*})^{-1}} p_*j_\sharp1_*=p_*j_\sharp.
$$
where $\Delta$ stands for the commutative square made by removing $\pi$ in
the diagram \eqref{pf:morph_comp}, and $Ex(\Delta_{\sharp*})$
is the corresponding natural transformation (see \ref{num:exchanges5}).
The compatibility of $Ex(\Delta_{\sharp*})$ with composition of morphisms
of schemes shows that we have defined a functor
$$
\Gamma_f: \op{\C}_f \rightarrow \uHom(\T(Y),\T(X))
$$
which sends all the maps of $\C_f$ to isomorphisms (by the support property).

The category $\C_f$ is non-empty
 by the assumption \ref{num:assumption1_sch}(c) on $\sch$,
and it is in fact left filtering; see \cite[XVII, 3.2.6(ii)]{SGA4}.
This defines a canonical functor $f_!:\T(Y) \rightarrow \T(X)$,
independent of any choice compactification of $f$, defined in
the category of functors
$\uHom(\T(Y),\T(X))$ by the formula
$$f_!=\varinjlim_{\op{\C}_f}\Gamma_f\, .$$
If $f=p$ is proper, then the compactification
$$Y \xrightarrow = Y \xrightarrow{p} X$$
is an initial object of $\C_f$, which gives a canonical
identification $p_!= p_*$.
Similarly, if $f=j$ is an open immersion, then the compactification
$$Y \xrightarrow j  X \xrightarrow{=} X$$
is a terminal object of $\C_j$, so that we get a canonical
identification $j_! = j_\sharp$.

This construction is compatible with composition of morphisms.
Let $g:Z\To Y$ and $f:Y\To X$ be two separated morphisms
of finite type in $\sch$.
For any a couple of compactifications
$$Z \xrightarrow k \bar Z \xrightarrow q Y\ \text{and} \
Y \xrightarrow j \bar Y \xrightarrow p X\, $$
of $f$ and $g$ respectively, we can choose a compactification
$$\bar Z  \xrightarrow h T \xrightarrow r Y$$
of $jq$, and we get a canonical isomorphism
$$f_!\, g_! \simeq p_*\, j_\sharp\, q_*\, k_\sharp\simeq p_*\, r_*\, h_\sharp \, k_\sharp \simeq (pr)_*\, (hk)_\sharp\simeq (fg)_!\, .$$
The independence of these isomorphisms with respect to the choices of compactification
follows from \cite[XVII, 3.2.6(iii)]{SGA4}. The cocycle conditions
(i.e. the associativity) also follows formally from \cite[XVII, 3.2.6]{SGA4}.
The uniqueness statement is obvious.
\end{proof}

\begin{paragr}\label{univ2functorsupport}
This construction is functorial in the following sense.

Define a \emph{$2$-functor with support on $\T$} to be a triple $(\D,a,b)$, where:
\begin{itemize}
\item[(i)] $\D:\schsep\to \tri$
is a $2$-functor (we shall write the structural
coherence isomorphisms as $c_{g,f}:\D(gf)\xrightarrow\sim\D(g)\D(f)$
for composable arrows $f$ and $g$ in $\schsep$);
\item[(ii)] $a:\T_*|_{\schprop}\To \D|_{\schprop}$
 and $b:\T_\sharp\To \D|_{\schouv}$
are morphisms of $2$-functors which agree on objects, i.e. such that
for any scheme $S$ in $\sch$, we have
$$\psi_S=a_S=b_S:\T(S)\To \D(S)\, ;$$
\item[(iii)] for any commutative square of shape \eqref{squareexchangeV1}
in which $j$ and $k$ are open immersions, while $p$ and $q$
are proper morphisms, the diagram below commutes.
$$\xymatrix@=12pt@=15pt{
\psi_S \, j_\sharp q_*\ar[rr]^{\psi_S \mathit{Ex}(\Delta_{\sharp*})}\ar[d]_{b\,  q_*}
&& \psi_S \, p_*k_\sharp \ar[d]^{a k_\sharp}\\
\D(j)\psi_U q_* \ar[d]_{\D(j) a} && \D(p) \psi_T k_\sharp\ar[d]^{\D(p) b}  \\
\D(j)\D(q)\psi_V \ar[r]^(.4){c^{-1}_{j,q}} &\D(jq)=\D(pk)\psi_V & \ar[l]_(.4){c^{-1}_{p,k}} \D(p)\D(k)\psi_V
}$$
\end{itemize}
Morphisms of $2$-functors with support on $\T$
$$(\D,a,b)\To (\D',a',b')$$
are defined in the obvious way: these are morphisms of
$2$-functors $\D\To \D'$ which preserve all the structure on the nose.

Using the arguments of the proof of \ref{prop:support&exists_f_!},
one checks easily that the category of $2$-functors
with support has an initial object, which is nothing else but
the $2$-functor $\T_!$ together with the identities
of $\T_*|_{\schprop}$ 
and of $\T_\sharp$ respectively. In particular, for any
$2$-functor $\D:\schsep \to \tri$,
a morphism of $2$-functors $\T_!\To \D$ is completely
determined by its restrictions to $\schprop$
and $\schouv$, and by its compatibility
with the exchange isomorphisms of type $Ex(\Delta_{\sharp*})$
 in the sense described in condition (iii) above.
\end{paragr}

\begin{prop}\label{defoublisupportcompact}
Assume that $\T$ satisfies the support property
 and consider the notations of Proposition
 \ref{prop:support&exists_f_!}.
For any separated morphism of finite type $f$ in $\sch$, 
there exists a canonical natural transformation
$$\alpha_f:f_!\to f_* \, .$$
The collection of maps $\alpha_f$ defines a morphism of $2$-functors
$$
\alpha:\T_! \rightarrow \T_*|_{\schsep} \ , \quad
 f \mapsto (\alpha_f:f_! \rightarrow f_*)
$$
whose restrictions to $\schprop$
and $\schouv$ are respectively
the identity and the morphism of $2$-functors
$\gamma: \T_\sharp \To \T_*|_{\schouv}$ 
defined in \ref{num:forget_support}.
\end{prop}

\begin{proof}
The identities $f_*=f_*$ for $f$ proper (resp. projective) and
the exchange natural transformations of type
$Ex(\Delta_{\sharp*})$ turns $\T_*|_{\schsep}$
into a $2$-functor with support (resp. restricted support) on $\T$ (property (iii)
of \ref{univ2functorsupport} is expressed by the commutative square \eqref{squareexchangeV4}).
\end{proof}

\begin{prop}\label{functexceptional}
Let $\T'$ be another triangulated complete $\Pmor$-fibred category over $\sch$.
Assume that $\T$ and $\T'$ both have the support property,
 and consider given a triangulated morphism of $\Pmor$-fibred categories $\varphi^*:\T\To\T'$
 (recall definition \ref{df:Pfibred_adjunction}).

Then, there is a canonical family of natural transformations
$$
Ex(\varphi^*,f_!) : \varphi^*_X\, f_!\To f_! \, \varphi^*_Y
$$
for each separated morphism of finite type $f:Y\To X$ in $\sch$, 
which is functorial with respect to composition in $\sch$
(i.e. defines a morphism of $2$-functors) and such that, the following conditions are verified:
\begin{itemize}
\item[(a)] if $f$ is proper, then, under the identification $f_!=f_*$,
the map $Ex(\varphi^*,f_!)$ is the exchange transformation
$Ex(\varphi^*,f_*):\varphi^*_X\, f_*\To f_* \, \varphi^*_Y$
defined in \ref{num:exchanges_Pfibred_morph};
\item[(b)] if $f$ is an open immersion, then, under the identification
$f_!=f_\sharp$, the map $Ex(\varphi^*,f_!)$
is the inverse of the exchange isomorphism
$Ex(f_\sharp,\varphi^*):f_\sharp \, \varphi^*_Y\To \varphi^*_X\, f_\sharp$ defined in \ref{num:Pfibred_morph}.
\end{itemize}
\end{prop}

\begin{proof}
The exchange maps of type $Ex(\varphi^*,f_*)$
define a morphism of $2$-functors
$$a:\T_*|_{\schprop}\To \T'_*|_{\schprop}=\T'_!|_{\schprop}$$
while the inverse of the exchange isomorphisms of type $Ex(f_\sharp,\varphi^*)$
define a morphism of $2$-functors
$$b:\T_\sharp \To \T'_\sharp=\T'_!|_{\schouv}\, ,$$
in such a way that the triple $(\T'_!,a,b)$ is a $2$-functor with support on $\T$.
%% 
%% Given a separated morphism of finite type $f:Y\To X$,
%% as in the proof of Proposition \ref{prop:support&exists_f_!},
%% consider the left filtering category $\C_f$ of compactifications
%% of $f$. One then checks that
%% we have a two functors $A$ and $B$ from $\C_f$ to $\T'(X)$:
%% for a compactification of $f$ of shape \eqref{pf:morph_comp00},
%% we put $A(p,j)=\varphi^*_X\, p_* \, j_\sharp$ and $B(p,j)=p_* \, j_\sharp \, \varphi^*_Y$.
%% Using the exchange maps of type \eqref{exchangePmorphismfupperstar}
%% and \eqref{exchangemorphismofPfibredcat}, this defines a natural
%% transformation $A\To B$, and the map $\delta_f$ is defined as its colimit.
\end{proof}

%\begin{rem}
%Under the assumptions of the previous proposition, if moreover $\varphi$
%is a morphism of complete $\Pmor$-fibred categories, we have
%natural transformations
%$$\varphi_{*,Y}\, f^!\To f^! \, \varphi_{*,X}$$
%obtained from $\delta_f$ by transposition.
%\end{rem}

%%%%%%%%%%%%%%%%%%%%%%%%%%%%%%%%%%%%%%
\begin{cor} \label{prop:support_exists_f_!_plus}
Suppose $\T$ satisfies the support property
  and consider the notations of proposition
 \ref{prop:support&exists_f_!}.
\begin{enumerate}
\item For any cartesian square
$$\xymatrix@=16pt{
Y'\ar^{f'}[r]\ar_{g'}[d]\ar@{}|\Delta[rd] & X'\ar^g[d] \\
Y\ar[r]_f & X,
}$$
such that $f$ is separated of finite type,
there exists a canonical natural transformation
\begin{align*}
Ex(\Delta_!^*):g^*f_! \rightarrow f'_!{g'}^*
\end{align*}
compatible with horizontal and vertical compositions of squares,
 and satisfying the following identifications in $\T(X')$
$$\xymatrix@C=40pt@R=12pt{
\ar@{}|/-20pt/{(a) \text{ $f$ proper:}}[r]
 & & \ar@{}|/-4pt/{(b) \text{ $f$ open immersion:}}[r] & \\
g^*f_!\ar^-{Ex(\Delta^*_!)}[r]\ar@{=}[d]
 & f'_!{g'}^*\ar@{=}[d] 
 & g^*f_!\ar^-{Ex(\Delta^*_!)}[r]\ar@{=}[d] & f_!{g'}^*\ar@{=}[d] \\
g^*f_*\ar^-{Ex(\Delta^*_*)}[r] & f'_*{g'}^*,
 & g^*f_\sharp\ar^-{Ex(\Delta^*_\sharp)^{-1}}[r] & f'_\sharp {g'}^*.
}$$
Moreover, when $g$ is a $\Pmor$-morphism, $Ex(\Delta_!^*)$ is an isomorphism.
\item For any cartesian square $\Delta$ as in (1),
 assuming $f$ is separated of finite type and $g$ is a $\Pmor$-morphism,
 there exists a canonical natural transformation
\begin{align*}
Ex(\Delta_{\sharp\,!}):g_\sharp f'_! \rightarrow f_! g'_\sharp
\end{align*}
compatible with horizontal and vertical compositions of squares,
 and satisfying the following identifications in $\T(X')$
$$\xymatrix@C=40pt@R=12pt{
\ar@{}|/-20pt/{(a) \text{ $f$ proper:}}[r]
 & & \ar@{}|/-4pt/{(b) \text{ $f$ open immersion:}}[r] & \\
g_\sharp f'_!\ar^-{Ex(\Delta_{\sharp!})}[r]\ar@{=}[d]
 & f_! g'_\sharp\ar@{=}[d] 
 & g_\sharp f'_!\ar^-{Ex(\Delta_{\sharp!})}[r]\ar@{=}[d]
 & f_! g'_\sharp\ar@{=}[d] \\
g_\sharp f^{\prime*}\ar^-{Ex(\Delta_{\sharp*})}[r]
 & f_* g'_\sharp,
 & g_\sharp f'_\sharp\ar@{=}[r]
 & f_\sharp g'_\sharp.
}$$
\item If furthermore $\T$ is monoidal
% triangulated $\Pmor$-fibred category over $\sch$, 
then for any separated morphism of  finite type $f:Y \rightarrow X$,
 there is a natural transformation
$$
Ex(f_!^*,\otimes):(f_!K) \otimes L \rightarrow f_!(K \otimes f^*L)
$$
which is compatible with respect to composition in $\sch$,
and such that, in each of the following cases, we have the following identifications:
$$
\xymatrix@C=15pt@R=14pt{
&\ar@{}|/-14pt/{(a) \text{ $f$ proper:}}[r]
 && & \ar@{}|/-0pt/{(b) \text{ $f$ open immersion:}}[r] && \\
&(f_!K) \otimes L\ar^-{Ex(f_!^*,\otimes)}[rr]\ar@{=}[d]
 && f_!(K \otimes f^*L)\ar@{=}[d] 
 & (f_!K) \otimes L\ar^-{Ex(f_!^*,\otimes)}[rr]\ar@{=}[d]
 && f_!(K \otimes f^*L)\ar@{=}[d] \\
&(f_*K) \otimes L\ar^-{Ex(f_*^*,\otimes)}[rr] && f_*(K \otimes f^*L),
 & (f_\sharp K) \otimes L\ar^-{Ex(f_\sharp^*,\otimes)^{-1}}[rr]&& f_\sharp(K \otimes f^*L).
}
$$
\end{enumerate}
\end{cor}
As in the previous analogous cases,
 the natural transformations $Ex(\Delta_!^*)$, $Ex(\Delta_{\!\sharp,!})$
 and $Ex(f_!^*,\otimes)$
 will be called \emph{exchange transformations}.
\index{word}{exchange!transformation}
\begin{proof}
To prove (1), consider a fixed map $g:X'\To X$ in $\sch$.
We consider the triangulated $\Pmorx X$-fibred
categories $\T'$ and $\T''$ over $\sch/X$ defined by
$\T'(Y)=\T(Y)$ and $\T''(Y)=\T(Y')$
for any $X$-scheme $Y$ (in $\sch$), with
$g':Y'=Y\times_X X'\To Y$ the map obtained from $Y\To X$ by pullback
along $g$. The collection of functors
$$g^{\prime \, *}:\T(Y)\To\T(Y')$$
define an exact morphism of triangulated $\Pmorx X$-fibred categories over $\sch/X$
(by the $\Pmor$-base change formula):
$$\varphi^*:\T'\To\T''\, .$$
Applying the preceding proposition to the latter gives (1).
The fact that we get an isomorphism whenever $g$ is a $\Pmor$-morphism
 follows from the $\Pmor$-base change formula
 and from paragraph \ref{num:exchanges2}.

For point (2), we consider the notations above assuming that
 $g$ is a $\Pmor$-morphism. The collection of functors
$$
g'_\sharp:\T(Y') \rightarrow \T(Y)
$$
associated with an $X$-scheme $Y$,
 $g':Y'=Y \times_X X' \rightarrow Y$ obtained from $g$ as above,
 define an exact morphism of triangulated $\Pmorx X$-fibred categories over $\sch/X$
 (applying again the $\Pmor$-base change formula):
$$
\varphi^*:\T'' \rightarrow \T'\, .
$$
Applying the preceding proposition to the latter gives (2).

The proof of (3) is similar: fix a scheme $X$ in $\sch$, as well as an object
$L$ in $\T(X)$. Let $\T'$ be the restriction of $\T$ to $\sch/X$ as above.
We can consider $L$ as a cartesian section of $\T'$,
and by the $\Pmor$-projection formula,
we then have an exact morphism of triangulated $\Pmorx X$-fibred categories over $\sch/X$:
$$L\otimes (-):\T'\To\T'\, .$$
Here again, we can apply the preceding proposition and conclude.

%For point (2), we define $Ex(\Delta_{!*})$ as the following composite map:
%$$
%f_!g'_* \xrightarrow{ad(g^*,g_*)} g_*g^*f_!g'_*
% \xrightarrow{Ex(\Delta_!^*)} g_*f'_!g^{\prime*}g'_*
% \xrightarrow{ad'(g^{\prime*},g'_*)} g_*f'_!.
%$$
%According to this definition,
% the identifications (a) and (b) follow easily from their
% respective counterpart in point (1).
\end{proof}
%\rem The cautious reader will check that the exchange transformation appearing 
% in point (i) is compatible with composition. Moreover, the same principle
% as in remark \ref{rem:coherence_exchangeV} about the coherence of
% all the exchange transformations introduced so far is valid.

\subsubsection{Further properties}
We will be particularly interested in the following properties
 of the triangulated $\Pmor$-fibred category $\T$.
\begin{df} \label{df:adj_BC_PF_prop_tri_premotivic}
Let $f:Y \rightarrow X$ be a morphism in $\sch$.
We introduced the following properties for $\T$,
 assuming in the third case that $\T$ is monoidal:
\begin{itemize}
\item[\adjx f] The functor $f_*$ admits a right adjoint.
 Under this assumption, we denote by $f^!$ the right adjoint of $f_*$.
\item[\BCx f] Any cartesian square of $\sch$ of the form
$$
\xymatrix@=16pt{
Y'\ar^{f'}[r]\ar_{g'}[d]\ar@{}|\Delta[rd] & X'\ar^g[d] \\
Y\ar[r]_f & X,
}
$$
is $\T$-transversal (Def. \ref{df:P-fibred_transversality}) -- \emph{i.e.}
the exchange transformation
$$
Ex(\Delta_*^*):g^*f_* \rightarrow f'_*g^{\prime*}
$$ 
associated with $\Delta$ is an isomorphism.
\item[\PFx f] For any object premotive $M$ over $Y$, and $N$ over $X$,
 the exchange transformation (see paragraph \ref{num:exchanges4})
$$
Ex(f_*^*,\otimes_X):(f_*M) \otimes_X N
  \rightarrow f_*(M \otimes_Y f^*N)
$$
is an isomorphism.
\end{itemize}
We denote by \adj (resp. \BC, \PF) the property
 \adjx f (resp. \BCx f, \PFx f) for any \emph{proper} morphism $f$ in $\sch$
 and call it the \emph{adjoint property} 
 (resp. \emph{proper base change property}, \emph{projection formula}).
\end{df}

%Recall from definition \ref{df:P-fibred_transversality} that we say
% $\T$ satisfies the transversality property with respect to proper morphisms (in $\S$)
% if the exchange transformation $Ex(\Delta_*^*):p^*f_* \rightarrow g_*q^*$
% defined in \ref{num:exchanges2} is an isomorphism
% for any cartesian square $\Delta$ as soon as $f$ is a proper morphism (in $\S$).
%\begin{prop}\label{excepproperbasechange}
%Assume that $\T$ satisfies the support property. Then the following conditions are equivalent:
%\begin{itemize}
%\item[(i)] $\T$ satisfies the transversality property with respect to proper morphisms in $\S$.
%\item[(ii)] For any cartesian square of $\sch$
%$$
%\xymatrix@=16pt{
%Y'\ar^{f'}[r]\ar_{g'}[d]\ar@{}|\Delta[rd] & X'\ar^g[d] \\
%Y\ar[r]_f & X,
%}
%$$
%in which $f$ is separated and of finite type, the exchange transformation
%\begin{align*}
%Ex(\Delta_!^*):g^*f_! \rightarrow f'_!{g'}^*
%\end{align*}
%is an isomorphism.
%\end{itemize}
%\end{prop}
%\begin{proof}
%Property $(ii)$ is always verified in the case where $f$
%is an open immersion by assumption ($\Pmor$-base change formula).
%As any separated morphism factors as an open immersion
%followed by a proper morphism, this implies the equivalence between (i) and (ii).
%\end{proof}

We can summarize the construction and properties introduced in this section
 as follows:
\begin{thm} \label{thm:support}
Assume $\T$ satisfies the properties \supp and \adj.

Then for any separated morphism of finite type
 $f:Y \rightarrow X$ in $\S$,
 there exists an essentially unique pair of adjoint functors
$$
f_!:\T(Y) \rightleftarrows \T(X):f^!
$$
called the \emph{exceptional functors},
\index{word}{functor!exceptional}
such that:
\begin{enumerate}
\item There exists a structure of a covariant (resp. contravariant) 
 $2$-functor on $f \mapsto f_!$ (resp. $f \mapsto f^!$).
\item There exists a natural transformation $\alpha_f:f_! \rightarrow f_*$
 compatible with composition in $f$
 which is an isomorphism when $f$ is proper.
\item For any open immersion $j$, $j_!=j_\sharp$ and $j^!=j^*$.
\item For any cartesian square
$$
\xymatrix@=16pt{
Y'\ar^{f'}[r]\ar_{g'}[d]\ar@{}|\Delta[rd] & X'\ar^g[d] \\
Y\ar_/-2pt/f[r] & X,
}
$$
in which $f$ is separated and of finite type,
there exists natural transformations
\begin{align*}
Ex(\Delta^*_!)&:g^*f_! \rightarrow f'_!{g'}^*, \\
Ex(\Delta_*^!)&:g'_*{f'}^! \rightarrow f^!g_*
\end{align*}
which are isomorphisms in the following three cases:
\begin{itemize}
\item $f$ is an open immersion.
\item $g$ is a $\Pmor$-morphism.
\item $\T$ satisfies the proper base change property \BC.
\end{itemize}
\end{enumerate}
Assume that $\T$ is in addition monoidal.
Then the following property holds:
\begin{enumerate}
\item[(5)] For any separated morphism of finite type $f:Y \rightarrow X$ in $\S$,
 there exists natural transformations
\begin{align*}
Ex(f_!^*,\otimes):
(f_!K) \otimes_X L & \longrightarrow f_!(K \otimes_Y f^*L)\, ,\ \\
  \uHom_X(f_!(L),K) & \longrightarrow f_* \uHom_Y(L,f^!(K))\, ,\ \\
  f^! \uHom_X(L,M) & \longrightarrow \uHom_Y(f^*(L),f^!(M))\, .
\end{align*}
which are isomorphisms in the following cases:
\begin{itemize}
\item $f$ is an open immersion.
\item $\T$ satisfies the projection formula \PF.
\end{itemize}
\end{enumerate}
\end{thm}
\noindent Indeed
 the existence of $f_!$ follows from Proposition \ref{prop:support&exists_f_!}
 while that of $f^!$ follows directly from assumption \adj.
 Assertions (1) and (3) follows from the construction,
  (2) is Proposition \ref{defoublisupportcompact},
 (4) (resp. (5)) follows from Corollary \ref{prop:support_exists_f_!_plus}
 and the definition of \BC (resp. \PF). Note also that the second and third 
 isomorphisms in (5) are obtained by transposition from $Ex(f_!,\otimes)$.

\begin{paragr}
While the properties \BCx f and \PFx f are only reasonable in practice
 for proper morphisms, this is not the case for the property \adjx f.
Recall that an exact functor between
well generated triangulated categories admits a right
adjoint if and only if it commutes with small sums: this is
an immediate consequence of the \emph{Brown representability theorem}
\index{word}{Brown representability theorem}
proved by Neeman (\textit{cf.} \cite[8.4.4]{Nee1}).
\end{paragr}
\begin{prop}\label{lm:exist_adjoint_f^!}
Assume that $\T$ is a compactly $\tau$-generated triangulated
\index{word}{generated!compactly $\tau$-generated!triangulated $\Pmor$-fibred}
premotivic category over $\sch$.
Then, for any morphism of schemes $f:T \rightarrow S$,
the functor $f_*:\T(T) \rightarrow \T(S)$ admits a right adjoint.
\end{prop}
\begin{proof}
This follows directly from Proposition \ref{prop:exist_right_adjoint}.
\end{proof}

\subsection{The localization property} \label{sec:localization}

\subsubsection{Definition}

\begin{num} \label{num:general_notations_loc}
Consider a closed immersion $i:Z \rightarrow S$ in $\S$.
Let $U=S-Z$ be the complement open subscheme of $S$ and
 $j:U \rightarrow S$ the canonical immersion.
We will use the following consequence of the triangulated $\Pmor$-fibred
structure on $\T$:
\begin{enumerate}
\item[(a)] The unit $1  \rightarrow j^*j_\sharp$ is an isomorphism.
\item[(b)] The counit $j^*j_* \rightarrow 1$ is an isomorphism.
\item[(c)] $i^*j_\sharp=0$.
\item[(d)] $j^*i_*=0$.
\item[(e)] The composite map
$
j_\sharp j^* \xrightarrow{ad'(j_\sharp,j^*)} 1 \xrightarrow{ad(i^*,i_*)} i_*i^*
$
is zero.
\end{enumerate}
In fact, the first four relations all follow from the base
 change property \bc. Relation (e) is a consequence of (d) once
 we have noticed that the following square is commutative
$$
\xymatrix@=10pt{
j_\sharp j^*\ar[r]\ar[d] & 1\ar[d] \\
j_\sharp j^*i_*i^*\ar[r] & i_*i^*.
}
$$

For the closed immersion $i$ and the triangulated category $\T$,
  we introduce the property \locx i made of the following assumptions:
\begin{enumerate}
\item[(a)] The pair of functors $(j^*,i^*)$ is conservative.
\item[(b)] The counit $i^*i_* \xrightarrow{ad'(i^*,i_*)} 1$
 is an isomorphism.
\end{enumerate}
\end{num}

\begin{df} \label{df:localization}
We say that $\T$ satisfies the \emph{localization property},
 denoted by \loc, if:
\begin{enumerate}
\item $\T(\varnothing)=0$.
\item For any closed immersion $i$ in $\S$, \locx i is satisfied.
\end{enumerate}
\end{df}

The main consequence of the localization axiom is that
 it leads to the situation of the six gluing functor
 (\textit{cf.} \cite[prop. 1.4.5]{BBD}):
\begin{prop}\label{prop:loc_triangle}
Let $i:Z \rightarrow S$ be a closed immersion
 with complementary open immersion $j:U \rightarrow S$
 such that \locx i is satisfied.
\begin{enumerate}
\item The functor $i_*$ admits a right adjoint $i^!$.
\item
For any $K$ in $\T(S)$,
 there exists a unique map
 $\partial_{i,K}:i_*i^*K \rightarrow j_\sharp j^*K[1]$
 such that the triangle
$$
j_\sharp j^* K \xrightarrow{ad'(j_\sharp,j^*)} K
 \xrightarrow{ad(i^*,i_*)} i_*i^* K
 \xrightarrow{\partial_{i,K}} j_\sharp j^* K[1]
$$
is distinguished.
The map $\partial_{i,K}$ is functorial in $K$.
\item 
For any $K$ in $\T(S)$,
 there exists a unique map
 $\partial'_{i,K}:j_*j^*K \rightarrow i_* i^!K[1]$
 such that the triangle
$$
i_* i^! K \xrightarrow{ad'(i_*,i^!)} K
 \xrightarrow{ad(j^*,j_*)} j_*j^* K
 \xrightarrow{\partial'_{i,K}} i_*i^! K[1]
$$
is distinguished.
The map $\partial'_{i,K}$ is functorial in $K$.
\end{enumerate}
\end{prop}

Under the property \locx i, the canonical triangles appearing
 in (2) and (3) above are called the \emph{localization triangles}
\index{word}{localization!triangle, see{triangle}}
\index{word}{triangle!localization triangle}
 associated with $i$.

\begin{proof}
We first consider point (2).
For the existence, we consider a distinguished triangle
$$
j_\sharp j^*K \xrightarrow{ad'(j_\sharp,j^*)} K \xrightarrow{\ \pi\ } C \xrightarrow{+1}
$$
Applying \ref{num:general_notations_loc}(e),
 we obtain a factorization
$$
\xymatrix@R=5pt@C=20pt{
K\ar^{ad(i^*,i_*)}[rr]\ar_\pi[rd] & & i_*i^*K \\
& C\ar_-w[ru] &
}
$$
We prove $w$ is an isomorphism. According to the above triangle,
 $j^*C=0$. From \ref{num:general_notations_loc}(d), $j^*i_*i^*K=0$
 so that $j^*w$ is an isomorphism.
Applying $i^*$ to the above distinguished triangle, we obtain 
from \ref{num:general_notations_loc}(c) that $i^*\pi$ is an isomorphism.
Thus, applying $i^*$ to the above commutative diagram together
with \locx i(b), we obtain that $i^*w$ is an isomorphism which concludes.

Consider a map $K \xrightarrow u L$ in $\T(S)$ and suppose we have
 chosen maps $a$ and $b$ in the diagram:
$$
\xymatrix@R=18pt@C=42pt{
j_\sharp j^* K\ar^-{ad'(j_\sharp,j^*)}[r]\ar_u[d]
 & K\ar^-{ad(i^*,i_*)}[r]\ar_u[d]
 & i_*i^*K\ar^a[r] & j_\sharp j^* K[1]\ar^u[d] \\
j_\sharp j^* L\ar^-{ad'(j_\sharp,j^*)}[r]
 & L\ar^-{ad(i^*,i_*)}[r] & i_*i^*L\ar^b[r] & j_\sharp j^* L[1]
}
$$
such that the horizontal lines are distinguished triangles.
We can find a map $h:i_*i^*K \rightarrow i_*i^*L$ completing
the previous diagram into a morphism of triangles.
Then the map $w=h-i_*i^*(u)$ satisfy the relation $w \circ ad(i^*,i_*)=0$.
Thus it can be lifted to a map in $\Hom(j_\sharp j^*K[1],i_*i^*L)$.
But this is zero by adjunction and the relation 
\ref{num:general_notations_loc}(d). This proves
both the naturality of $\partial_{i,K}$ and its uniqueness.

For point (1) and (3), for any object $K$ of $\T(S)$,
 we consider a distinguished triangle
$$
D \rightarrow K \xrightarrow{ad(j^*,j_*)} j_*j^* K \xrightarrow{+1}
$$
According to \ref{num:general_notations_loc}(b), $j^*D=0$.
Thus according to the triangle of point (2) applied to $D$,
 we obtain $D=i_*i^*D$. Arguing as for point (2),
 we thus obtain that $D$ is unique and depends functorially
 on $K$ so that, if we put $i^!K=i^*D$, point (1) and (3) follows.
\end{proof}

\begin{rem} Consider the hypothesis and notations of the previous proposition.
\label{rem:i^!j_*=0}
\begin{enumerate}
\item By transposition from \ref{num:general_notations_loc}(d),
 we deduce that $i^!j_*=0$.
\item Assume that $i$ is a $\Pmor$-morphism.
Then the $\Pmor$-base change formula implies that $i^*j_*=0$.
Dually, we get that $i^!j_\sharp=0$.
By adjunction, we thus obtain $\partial_{i,K}=0$ and $\partial'_{i,K}=0$
 for any object $K$ so that both localization triangles are split.
In that case, we get that $\T(S)=\T(Z) \times \T(U)$.\footnote{This remark
 explains why the localization property is too strong for generalized premotivic
 categories.}
\end{enumerate}
The preceding proposition admits the following reciprocal statement:
\end{rem}

\begin{lm} \label{lm:loc&weak_exact_seq(beurk)}
Consider a closed immersion $i:Z \rightarrow S$ in $\S$
 with complementary open immersion $j:U \rightarrow S$.
Then the following properties are equivalent:
\begin{enumerate}
\item[(i)] $\T$ satisfies \locx i.
\item[(ii)]
\begin{enumerate}
\item[(a)] The functor $i_*$ is conservative.
\index{word}{conservative}
\item[(b)] For any object $K$ of $\T(S)$,
 there exists a map $i_*i^*(K) \rightarrow j_\sharp j^*(K)[1]$
 which fits into a distinguished triangle
$$
j_\sharp j^*(K) \xrightarrow{ad'(j_\sharp,j^*)} K
 \xrightarrow{ad(i^*,i_*)} i_*i^*(K)
 \rightarrow j_\sharp j^*(K)[1]
$$
\end{enumerate}
\end{enumerate}
\end{lm}
\begin{proof} The fact (i) implies (ii) follows from
 Proposition \ref{prop:loc_triangle}.
Conversely, (ii)(b) implies that the pair $(i^*,j^*)$ is conservative
 and it remains to prove \locx i(b).
Let $K$ be an object of $\T(S)$. Consider the distinguished triangle
 given by (ii)(b):
$$
j_\sharp j^*(K) \xrightarrow{ad'(j_\sharp,j^*)} K
 \xrightarrow{ad(i^*,i_*)} i_*i^*(K) 
 \rightarrow j_\sharp j^*(K)[1].
$$
If we apply $i_*$ on the left to this triangle, 
 we get using \ref{num:general_notations_loc}(d)
 that the morphism
$$
i_*(K) \xrightarrow{ad(i^*,i_*).i_*} i_*i^*i_*(K) 
$$
is an isomorphism. Hence, by the zig-zag equation, the morphism
$$
i_*i^*i_*(K) \xrightarrow{i_*.ad'(i^*,i_*)} i_*(K)
$$
is an isomorphism. Property (ii)(a) thus implies that $i^* i_*(K)\simeq K$.
\end{proof}

\subsubsection{First consequences of localization}

The following statement is straightforward.
\begin{prop} \label{prop:easy_csq_loc}
Assume $\T$ satisfies the localization property 
 and consider a scheme $S$ in $\S$.
\begin{enumerate}
\item Let $S_{red}$ be the reduced scheme associated with $S$.
The canonical immersion $S_{red} \xrightarrow \nu S$ induces
 an equivalence of categories:
$$
\nu^*:\T(S) \rightarrow \T(S_{red}).
$$
\item For any any partition
\item{partition}
 $(S_i \xrightarrow{\nu_i} S)_{i \in I}$
 of $S$ by locally closed subsets, the family of functors
 $(\nu_i^*)_{i \in I}$ is conservative ($S_i$ is considered with its
 canonical structure of a reduced subscheme of $S$).
\end{enumerate}
\end{prop}

\begin{lm} \label{lm:loc=>additive}
If $\T$ satisfies the localization property \loc then it is additive.
\end{lm}
\begin{proof}
Note that, by assumption, $\T(\varnothing)=0$.
Then the assertion follows directly from Lemma \ref{lm:localization&sums}.
\end{proof}

\begin{prop} \label{prop:loc=>sepx cdh}
If $\T$ satisfies the localization property then it satisfies the $\cdh$-separation property.
\end{prop}
\begin{proof}
Consider a cartesian square of schemes
$$
\xymatrix@=14pt{
B\ar[r]\ar[d]\ar@{}|Q[rd] & Y\ar^p[d] \\
A\ar_e[r] & X.
}
$$
According to Lemma \ref{lm:cd-separation}, we have only to check
 that the pair of functors $(e^*,p^*)$ is conservative when
 $Q$ is a Nisnevich (or respectively a proper cdh) distinguished square.
Let $\nu:A' \rightarrow X$ be the complementary closed (resp. open) immersion to $e$,
 where $A'$ has the induced reduced subscheme (resp. induced subscheme) structure.
Consider the cartesian square
$$
\xymatrix@=14pt{
Y\ar_p[d] & B'\ar^q[d]\ar[l] \\
X & A'\ar^\nu[l].
}
$$
By assumption on $Q$, $q$ is an isomorphism. According to \loc(ii), $(e^*,\nu^*)$
 is conservative. This concludes.
\end{proof}

The following proposition can be found in a slightly less precise 
 and general form in \cite[2.1.162]{ayoub}.\footnote{A warning: the proof
 in \emph{loc. cit.} seems to require that the schemes are excellent.}
\begin{prop} \label{prop:redseptoetale}
Assume $\T$ satisfies the localization property.

Then the following conditions are equivalent:
\begin{enumerate}
\item[(i)] $\T$ is separated.
\item[(ii)] For a morphism $f:T \rightarrow S$ in $\S$,
 $f^*:\T(S) \rightarrow \T(T)$ is conservative whenever $f$
 is:
\begin{enumerate}
\item[(a)] a finite \'etale cover;
\item[(b)] finite, faithfully flat
\index{word}{morphism!faithfully flat}
 and radicial.
\index{word}{morphism!radicial}
\end{enumerate}
\end{enumerate}
\end{prop}
\begin{proof}
Only $(ii) \Rightarrow (i)$ requires a proof.
Consider a surjective morphism of finite type $f:T \rightarrow S$ in $\S$.
According to \cite[17.16.4]{EGA4}, there exists a partition
 $(S_i)_{i \in I}$ of $S$ by (affine) subschemes and a family of
 maps of the form
$$
S''_i \xrightarrow{g_i} S'_i \xrightarrow{h_i} S_i
$$
such that $g_i$ (resp. $h_i$) satisfies assumption (a) (resp. (b)) above
and such that for any $i \in I$, $f \times_S S''_i$ admits a section.
Thus, Proposition \ref{prop:easy_csq_loc} concludes.
\end{proof}

\subsubsection{Localization and exchange properties}

\begin{num}
Consider a morphism of complete triangulated $\Pmor$-fibred categories
 over $\sch$:
$$
\varphi^*:\T \rightarrow \T'.
$$
Recall that for any morphism $f:Y \rightarrow X$,
 there is an exchange transformation \eqref{exchangePmorphismfupperstar}:
$$
Ex(\varphi^*,f_*):\varphi^*_X f_*\longrightarrow f_* \varphi_Y^*.
$$
If $\T$ and $\T'$ satisfies the support axiom
 and $f$ is separated of finite type,
 we have constructed (Proposition \ref{functexceptional})
 another exchange transformation:
$$
Ex(\varphi^*,f_!):\varphi^*_X f_! \longrightarrow f_! \varphi_Y^*.
$$
\end{num}
\begin{prop} \label{prop:morph&loc1}
Consider a morphism $\varphi^*:\T \rightarrow \T'$ as above.
\begin{enumerate}
\item Let $i:Z \rightarrow X$ be a closed immersion such
 that $\T$ and $\T'$ satisfy property \locx i. \\
 Then the exchange $Ex(\varphi^*,i_*):\varphi_X^*i_* \rightarrow i_*\varphi^*_Z$
 is an isomorphism.
\index{word}{exchange!isomorphism}
\item Assume $\T$ and $\T'$ satisfy property \loc. \\
Then the following conditions are equivalent:
\begin{enumerate}
\item[(i)] For any integer $n>0$ and any scheme $X$ in $\sch$,
 the exchange $Ex(\varphi^*,p_{n*})$ is an isomorphism where
 $p_n:\PP^n_X \rightarrow X$ is the canonical projection.
\item[(ii)] For any proper morphism $f:Y \rightarrow X$,
 the exchange $Ex(\varphi^*,f_*)$ is an isomorphism.
\end{enumerate}
\item Assume $\T$ and $\T'$ satisfy properties \loc and \supp. \\
Then conditions (i) and (ii) above are equivalent to the following one:
\begin{enumerate}
\item[(iii)] For any separated morphism $f:Y \rightarrow X$ of finite type,
 the exchange $Ex(\varphi^*,f_!)$ is an isomorphism.
\end{enumerate}
\end{enumerate}
\end{prop}
\begin{rem} 
We will simply say that $\varphi^*$ commutes with $f_!$
 when assertion (iii) is fulfilled.
 For an important case where this happens,
  see Proposition \ref{prop:motivic_adj_6_operations}.
\end{rem}
\begin{proof}
Assertion (1) follows easily from the conservativity of $(i^*,j^*)$
 where $j$ is the complementary open immersion
 and the relations of paragraph \ref{num:general_notations_loc}.
Assertion (3) is an easy consequence of the definition of $f_!$
 and the exchange $Ex(\varphi^*,f_!)$.

Concerning assertion (2), we have to prove that (i) implies (ii). 
We fix a morphism $f:Y \rightarrow X$ and prove that
 the exchange $Ex(\varphi^*,f_*):\varphi^*_Y f_* \rightarrow f_*\varphi^*_X$
 is an isomorphism.

We first treat the case where $f$ is projective.
According to Proposition \ref{prop:loc=>sepx cdh},
 $\T'$ satisfies the Zariski separation property. 
Using the \bc property,
 we see that the problem is local in $X$
 so that we can assume $X$ is affine.
Then $X$ admits an ample line bundle and there exists an integer $n>0$
 such that $f$ can be factored 
 (\cite[(5.5.4)(ii)]{EGA2}) into
 a closed immersion $i:Y \rightarrow \PP^n_X$ and the projection
 $p_n:\PP^n_X \rightarrow X$.
 Thus, assertion (1) and assumption (i) allow us to conclude.

To treat the general case, we argue by noetherian induction on $Y$,
 assuming that for any proper closed subscheme $T$ of $Y$, 
 the result is known for the restriction of $f$ to $T$.
In fact, the case $T=\varnothing$ is obvious because $\T(\varnothing)=0$.

According to Chow's lemma
\index{word}{Chow's lemma}
 \cite[5.6.2]{EGA2},
 there exists a morphism $p:Y_0 \rightarrow Y$ such that:
\begin{enumerate}
\item[(a)] $p$ and $f \circ p$ are projective morphisms.
\item[(b)] There exists a dense open subscheme $V_0$ of $Y$
 over which $p$ is an isomorphism.
\end{enumerate}
Let $T$ be the complement of $V$ in $Y$ equipped with its reduced
 subscheme structure. Let $j$ and $i$ be the respective immersion
 of $T$ and $V$ in $Y$.
According to point (3) of Proposition \ref{prop:loc_triangle},
 it is sufficient to prove that the following natural transformations
 are isomorphisms:
\begin{align}
\label{eq1:pf_morph&loc} \varphi^*_Y f_*i_* & \rightarrow f_*\varphi^*_X i_*. \\
\label{eq2:pf_morph&loc} \varphi^*_Y f_*j_* & \rightarrow f_*\varphi^*_X j_*.
\end{align}
Concerning the first one, we consider the following commutative diagram:
$$
\xymatrix@R=16pt@C=44pt{
\varphi^*_Y f_*i_*\ar^-{Ex(\varphi^*,f_*)}[r]\ar@{=}[d]
 & f_*\varphi^*_X i_*\ar^-{Ex(\varphi^*,i_*)}[r]
 & f_*i_*\varphi^*_X\ar@{=}[d] \\
\varphi^*_Y (fi)_*\ar^-{Ex(\varphi^*,(fi)_*)}[rr]
 & & (fi)_*\varphi^*_X.
}
$$
Thus the result follows from assertion (1) and the induction hypothesis.

Concerning the natural transformation \eqref{eq2:pf_morph&loc}, 
we consider the pullback square
$$
\xymatrix{
V_0\ar^-l[r]\ar_q[d] & Y_0\ar^p[d] \\
V\ar^-j[r] & Y.
}
$$
Assumption (b) above says that $q$ is an isomorphism
 which implies the relation: $j_*=p_*l_*q^*$. 
 In particular, it is sufficient to prove that 
 the natural transformation
 $\varphi^*_Y f_*p_* \rightarrow f_*\varphi^*_X p_*$
 is an isomorphism. This follows from the commutativity
 of the following diagram
$$
\xymatrix@R=16pt@C=44pt{
\varphi^*_Y f_*p_*\ar^-{Ex(\varphi^*,f_*)}[r]\ar@{=}[d]
 & f_*\varphi^*_X p_*\ar^-{Ex(\varphi^*,p_*)}[r]
 & f_*p_*\varphi^*_X\ar@{=}[d] \\
\varphi^*_Y (fp)_*\ar^-{Ex(\varphi^*,(fp)_*)}[rr]
 & & (fp)_*\varphi^*_X,
}
$$
according to the projective case treated above 
 and assumption (b). The proof is complete.
\end{proof}

\begin{cor} \label{cor:loc&Supp,BC,PF}
In the next statements, we assume $\T$ is monoidal when it is needed.
\begin{enumerate}
\item Let $i:Z \rightarrow X$ be a closed immersion such
 that $\T$ satisfies property \locx i. \\
 Then $\T$ satisfies property \suppx i (resp. \BCx i, \PFx i).
\item Assume $\T$ satisfies the localization property.
 Then the following properties of $\T$ are equivalent:
\begin{enumerate}
\item[(i)] For any integer $n>0$ and any scheme $X$ in $\sch$,
 $p_n:\PP^n_X \rightarrow X$ being the canonical projection,
 $\T$ satisfies \suppx{p_n} (resp. \BCx{p_n}, \PFx{p_n}).
\item[(ii)] $\T$ satisfies \supp (resp. \BC, \PF).
\end{enumerate}
\item Assume $\T$ is well generated and
 satisfies the localization property.
 Then the following properties of $\T$ are equivalent:
\begin{enumerate}
\item[(i')] For any integer $n>0$ and any scheme $X$ in $\sch$,
 $p_n:\PP^n_X \rightarrow X$ being the canonical projection,
 $\T$ satisfies \adjx{p_n}.
\item[(ii')] $\T$ satisfies \adj.
\end{enumerate}
\end{enumerate}
\end{cor}
\begin{proof}
As in the proof of Corollary \ref{prop:support_exists_f_!_plus},
 each respective case of assertions (1) and (2) 
 follows from the previous proposition
 applied to a particular type of morphisms
  $\varphi^*:\T' \rightarrow \T''$
 of complete $\Pmor$-fibred triangulated categories
 over a subcategory $\sch'$ of $\sch$.

For property \supp, we proceed as follows.
 We fix an open immersion $j:U \rightarrow X$
 and let $\sch'=\sch/X$.
 For any $Y/X$, we let $j_Y=Y \times_X U \rightarrow Y$ be the pullback 
 of $j$. We put $\T'(Y)=\T(Y \times_X U)$ and $\T''(Y)=\T(Y)$
 and let $\varphi^*_Y$ be the functor:
$$
j_{Y\sharp}:\T(Y \times_X U) \rightarrow \T(Y).
$$

For the property \BC (resp. \PF),
 we refer the reader to the proof of assertion (1) (resp. (2))
  in Corollary \ref{prop:support_exists_f_!_plus}.

Finally we consider assertion (3).
 It is sufficient to prove that (\emph{i'}) implies (\emph{ii'}). \\
According to the Brown representability theorem
\index{word}{Brown representability theorem}
 \cite[8.4.4]{Nee1},
 the property \adjx f for a proper morphism $f$
 is equivalent to ask that $f_*$ preserves small sum. \\
Consider an arbitrary set $I$.
For any scheme $S$,
 we put $\T^I(S)=\T(S)^I$,
 that is the category of families of object of $\T(S)$ indexed by $I$.
 Then $\T^I$ is obviously a complete triangulated $\Pmor$-fibred
 category over $\sch$ (limits and colimits are computed termwise).
 For any scheme $S$, we consider the functor:
$$
\varphi_S^*:\T^I(S) \rightarrow \T(S)\ ,\quad
 (M_i)_{i \in I} \mapsto \bigoplus_{i\in I} M_i.
$$
Then $\varphi^*:\T^I \rightarrow \T$ is obviously a morphism
 of complete $\Pmor$-fibred categories.
 Thus, given condition (\emph{i'}),
 the preceding proposition applied to $\varphi^*$ 
 shows that for any proper morphism $f$,
  $f_*$ commutes with sums indexed by $I$. 
 As this is true for any $I$, we obtain (\emph{ii'}).
\end{proof}

\subsubsection{Localization and monoidal structure}

\begin{num} \label{num:loc&mon_notation}
 Assume $\T$ is monoidal
 and let $M$ denote its geometric sections.
\index{word}{section!geometric}
Fix a closed immersion $i:Z \rightarrow S$ in $\S$
 with complementary open immersion $j:U \rightarrow S$.
We fix an object $M_S(S/S-Z)$ of $\T(S)$
 and a distinguished triangle
\begin{equation} \label{eq:proto_loc_triangle}
M_S(S-Z) \xrightarrow{j_*} \un_S
 \xrightarrow{p_i} M_S(S/S-Z) \xrightarrow{d_i} M_S(S-Z)[1].
\end{equation}
Remark that according to \ref{num:general_notations_loc}(c),
the map $i^*(p_i):\un_Z \rightarrow i^*M_S(S/S-Z)$
is an isomorphism.
Given any object $K$ in $\T(S)$, we thus obtain
an isomorphism
$$
i^*(M_S(S/S-Z) \otimes_S K)=i^*(M_S(S/S-Z)) \otimes_Z i^*(K)
 \xrightarrow{(i^*p_i)^{-1}} \un_Z \otimes_Z i^*(K)=i^*(K)
$$
which is natural in $K$. It induces by adjunction a map
\begin{equation} \label{eq:loc&premotivic_notation}
\psi_{i,K}:M_S(S/S-Z) \otimes_S K \rightarrow i_*i^*(K)
\end{equation}
which is natural in $K$. \\
For any $\Pmor$-scheme $X/S$,
 we put $M_S(X/X-X_Z)=M_S(S/S-Z) \otimes_S M_S(X)$
 so that we get from \eqref{eq:proto_loc_triangle}
 a canonical distinguished triangle:
$$
M_S(X-X_Z) \xrightarrow{j_{X*}} M_S(X)
 \rightarrow M_S(X/X-X_Z) \rightarrow M_S(X-X_Z)[1].
$$
The map \eqref{eq:loc&premotivic_notation} for $K=M_S(X)$ gives a canonical map
\begin{equation} \label{eq:loc&premotivic_notation_bis}
\psi_{i,X}:M_S(X/X-X_Z) \rightarrow i_*(M_Z(X_Z)).
\end{equation}
\end{num}

\begin{prop} \label{prop:loc&monoidal}
Consider the previous hypothesis and notations.
Then the following conditions are equivalent:
\begin{enumerate}
\item[(i)] $\T$ satisfies the property \locx i.
\item[(ii)]
\begin{enumerate}
\item[(a)] The functor $i_*$ is conservative.
\index{word}{conservative}
\item[(b)] The morphism $\psi_{i,S}:M_S(S/S-Z) \rightarrow i_*(\un_Z)$
is an isomorphism.
\item[(c)] For any object $K$ of $\T(S)$,
the exchange transformation
$$
Ex(i^*_*,\otimes):(i_*\un_Z) \otimes_S K \rightarrow i_*i^*K
$$
is an isomorphism.
\end{enumerate}
\item[(iii)] 
\begin{enumerate}
\item[(a)] The functor $i_*$ is conservative.
\index{word}{conservative}
\item[(b)] The morphism $\psi_{i,S}:M_S(S/S-Z) \rightarrow i_*(\un_Z)$
is an isomorphism.
\item[(c)] For any objects $K$ and $L$ of $\T(S)$,
the exchange transformation
$$
Ex(i^*_*,\otimes):
(i_* K) \otimes_S L \rightarrow i_*(K \otimes_Z i^* L)
$$
is an isomorphism.
\end{enumerate}
\end{enumerate}
Assume in addition that $\T$ is well generated and $\tau$-generated
\index{word}{generated!$\tau$-generated}
 as a triangulated $\Pmor$-fibred category. Then the above conditions
are equivalent to the following one:
\begin{enumerate}
\item[(iv)]
\begin{enumerate}
\item[(a)] The functor $i_*$ is conservative, commutes with direct sums
 and with $\tau$-twists.
\item[(b)] The morphism $\psi_{i,X}:M_S(X/X-X_Z) \rightarrow i_*(M_Z(X_Z))$ is an isomorphism
 for any $\Pmor$-scheme $X/S$.
\end{enumerate}
\end{enumerate}
\end{prop}
\noindent In particular, \locx i implies that for any object $K$ of $\T(S)$,
the localization triangle of \ref{prop:loc_triangle}
$$
j_\sharp j^*(K) \rightarrow K \rightarrow i_*i^*(K)
 \xrightarrow{\partial_K} j_\sharp j^*(K)[1]
$$
is canonically isomorphic (through exchange transformations) to 
 the triangle \eqref{eq:proto_loc_triangle} tensored with $K$.
\begin{proof}
\noindent $(i) \Rightarrow (iii)$~: 
According to \locx i(a),
 we need only to check that the maps in (iii)(b) and (iii)(c) are isomorphisms
 after applying $i^*$ and $j^*$.
This follows easily from \locx i(b).

\noindent $(iii) \Rightarrow (ii)$~: Obvious

\noindent $(ii) \Rightarrow (i)$~: According to (ii)(b),
 the distinguished triangle \eqref{eq:proto_loc_triangle} is isomorphic
 to a triangle of the form
$$
j_\sharp j^*(\un_S) \xrightarrow{ad'(j_\sharp,j^*)} \un_S
 \xrightarrow{ad(i^*,i_*)} i_*i^*(\un_S) \rightarrow j_\sharp j^*(\un_S).
$$
According to (ii)(c), this latter triangle tensored with $K$ is
 isomorphic through exchange transformations to a triangle of the form
$$
j_\sharp j^*(K) \xrightarrow{ad'(j_\sharp,j^*)} K
 \xrightarrow{ad(i^*,i_*)} i_*i^*(K) \rightarrow j_\sharp j^*(K).
$$
Thus Lemma \ref{lm:loc&weak_exact_seq(beurk)} allows us to conclude.

To end the proof, we remark by using the equations for the adjunction
 $(i^*,i_*)$ that for any object $M$ of $\T(S)$,
 the following diagram is commutative:
$$
\xymatrix@=12pt@C=40pt{
 & i_*i^*(\un_S) \otimes K\ar@{=}[r]
  & i_*(\un_Z) \otimes K\ar^{Ex(i_*^*,\otimes)}[dd] \\
M_S(S/S-Z) \otimes K\ar^-{\psi_{i} \otimes 1_K}[ru]\ar_-{\psi_{i,K}}[rd]
 & & \\
 & i_*i^*(K)\ar@{=}[r]
 & i_*(\un_Z \otimes i^*i^*(K)).
}
$$
Note that (i) implies that $i_*$ is conservative and commutes with direct sums
 (see \ref{prop:loc_triangle}) and (ii)(c) implies it commutes with twists.
According to the above diagram, (ii)(b) implies (iv)(b). \\
We prove that reciprocally that (iv) implies (ii).
Because (ii)(b) (resp. (ii)(a)) is a particular case of (iv)(b) (resp. (iv)(a)),
 we have only to prove (ii)(b).
In view of the previous diagram,
 we are reduced to prove that for any object $K$ of $\T(S)$,
 the map 
 $\psi_{i,K}$ is an isomorphism.
Consider the full subcategory $\U$ of $\T(S)$
 made of the objects $K$ such that $\psi_{i,K}$ is an isomorphism.
 Then $\U$ is triangulated. Using (iv)(a), $\U$ is stable by small sums
 and $\tau$-twists. By assumption, it contains the objects of the form $M_S(X)$
 for a $\Pmor$-scheme $X/S$.
 Thus, because $\T$ is well generated by assumption,
 Lemma \ref{lm:well_gen_tri_Pfibred} concludes.
\end{proof}

\begin{lm} \label{lm:generated&conservativity_of_i_*}
Consider a closed immersion $i:Z \rightarrow S$.
We assume the following conditions are satisfied
 in addition to that of \ref{num:assumption1_sch}:
\begin{itemize}
\item $\T$ is well generated, $\tau$-generated,
 and satisfies the Zariski separation property.
\item For any $\Pmor$-scheme $X_0/Z$ and any point $x_0$ of $X_0$,
 there exists an open neighborhood $U_0$ of $x_0$ in $X_0$
 and a $\Pmor$-scheme $U/S$ such that $U_0=U \times_S Z$.\footnote{
 This property is trivial when $\Pmor$ is the class of open immersions
  or the class of morphisms of finite type in $\sch$.
 It is also true when $\Pmor$ is the class of \'etale morphism or 
 $\Pmor=\sm$ (\textit{cf.} \cite[18.1.1]{EGA4}).}
\end{itemize}

Then the functor $i_*$ is conservative.
\end{lm}
\begin{proof}
Consider an object $K$ of $\T(Z)$ such that 
 $i_*(K)=0$. We prove that $K=0$. \\
Because $\T$ is $\tau$-generated,
 it is sufficient to prove that for a $\Pmor$-morphism
  $p_0:X_0 \rightarrow Z$
 and a twist $(n,m) \in \ZZ \times \tau$,
$$
\Hom_{\T(Z)}(M_Z(X_0)\{m\}[n],K)=0.
$$
Because $M_Z(X_0)=p_{0\sharp}(\un_{X_0})$,
 this equivalent to prove that 
$$\Hom_{\T(X_0)}(\un_{X_0}\{m\}[n],p_0^*(K))=0.$$
Using the Zariski separation property on $\T$,
 this latter assumption is local in $X_0$. Thus,
  according to the assumption on the class $\Pmor$,
 we can assume there exists a $\Pmor$-scheme $X/S$ such that $X_0=X \times_S Z$.
 Thus $M_Z(X_0)\{m\}[n]=i^*(M_S(X)\{m\}[n])$
 and the initial assumption on $K$ allows us to conclude.
\end{proof}

Note for future applications the following interesting corollaries:
\begin{cor} \label{cor:premotivic&i_*}
Assume $\T$ is a premotivic triangulated category
 which is compactly $\tau$-generated
\index{word}{generated!compactly $\tau$-generated!triangulated $\Pmor$-fibred}
 for a group of twists $\tau$ 
(\emph{i.e.} any twists in $\tau$ admits a tensor inverse)
 and which satisfies the Zariski separation property.

Then, for any closed immersion $i$, the functor $i_*$ is conservative,
 commutes with sums and with twists.
\end{cor}
This is a consequence of lemmas \ref{lm:generated&conservativity_of_i_*}
 and \ref{lm:exist_adjoint_f^!}. In fact, 
 under these conditions, $i_*$ commutes with arbitrary
 $\tau$-twists because it is true for its (left) adjoint $i^*$.

\begin{cor}\label{cor:premotivic&i_*_bis}
Assume $\T$ satisfies the assumptions of the preceding corollary.
Then the following conditions on a closed immersion
$i$ are equivalent:
\begin{enumerate}
\item[(i)] $\T$ satisfies the property \locx i.
\item[(ii)] For any scheme $S$ in $\sch$ and any smooth $S$-scheme $X$,
 the map \eqref{eq:loc&premotivic_notation_bis}
$$
\psi_{i,X}:M_S(X/X-X_Z) \rightarrow i_*M_Z(X_Z)
$$
is an isomorphism.
\end{enumerate}
\end{cor}

We finish this section with the following useful result:
\begin{prop} \label{prop:localization&exchange}
Assume $\T$ is $\tau$-generated
 and consider a $\tau'$-generated triangulated $\Pmor$-fibred category $\T'$
 and a morphism
$$
\varphi^*:(\T,\tau) \rightleftarrows (\T',\tau'):\varphi_*\, .
$$
We assume the following properties:
\begin{enumerate}
\item[(a)] the morphism $\varphi^*$ is strictly compatible with twists;
\item[(b)] $\T'$ is well generated.
\end{enumerate}
We consider a closed immersion $i:Z \rightarrow S$
 and further assume the following properties:
\begin{enumerate}
\item[(c)] $\T$ satisfies the property \locx i.
\item[(d)] The exchange transformation $Ex(\varphi^*,i_*):\varphi^* i_* \rightarrow i_* \varphi^*$
 is an isomorphism.
\item[(e)] The functor $i_*:\T'(Z) \rightarrow \T'(S)$ commutes with $\tau'$-twists.\footnote{This
 will be satisfied if any $\tau'$-twists is invertible because
 the left adjoint of $i_*$ commutes with $\tau'$-twists.}
\end{enumerate}
Then $\T'$ satisfies the property \locx i.
\end{prop}
\begin{proof}
Note that, under the above assumptions, $\varphi_*$ is conservative
 (in fact, for any $\Pmor$-scheme $X/S$ and any twists $i \in \tau'$,
 the premotive $M_S(X)\{i\}$ is in the essential image of $\varphi^*$).
Thus, if $i_*:\T(Z) \rightarrow \T(S)$ is conservative (resp. commute with sums),
 then $i_*:\T'(S) \rightarrow \T'(S)$ is conservative (resp. commute with sums)
 using the isomorphism $\varphi_*i_*\simeq i_*\varphi_*$. \\
Let $M$ (resp. $M'$) be the geometric sections of $\T$ (resp. $\T'$).
As in \ref{num:loc&mon_notation}, we fix a distinguished triangle
$$
M_S(S-Z) \xrightarrow{j_*} \un_S
 \xrightarrow{p_i} M_S(S/S-Z) \xrightarrow{d_i} M_S(S-Z)[1].
$$
and we put $M'_S(S/S-Z)=\varphi^*M_S(S/S-Z)$.
According to \emph{loc. cit.}, we thus get for any $\Pmor$-scheme $X/S$
 canonical maps
\begin{align*}
\psi_{i,X}:M_S(X/X-X_Z) &\rightarrow i_*M_Z(X_Z), \\
\psi'_{i,X}:M'_S(X/X-X_Z) &\rightarrow i_*M'_Z(X_Z). \\
\end{align*}
By construction, the following diagram is commutative:
$$
\xymatrix@=50pt@R=18pt{
\varphi^*M_S(X/X-X_Z)\ar@{=}[d]\ar^{\varphi^*\psi_{i,X}}[r]
 & \varphi^* i_*M_Z(X_Z)\ar^{Ex(\varphi^*,i_*)}[r]
 & i_* \varphi^*M_Z(X_Z)\ar@{=}[d] \\
M'_S(X/X-X_Z)\ar^{\psi'_{i,X}}[rr]
 && M'_Z(X_Z)
}
$$
Thus, Proposition \ref{prop:loc&monoidal} allows us to conclude.
\end{proof}

%%%%%%%%%%%%%%%%%%%%%%%%%%%%%%%%%%%%%%%%%%%%%%%%%%%%%%%%%%%%%%%%%%%%%%
% purity
%
%

\subsection{Purity and the theorem of Voevodsky-R\"ondigs-Ayoub} \label{sec:Ayoub}

Recall we assume $\Pmor=\sm$ in this section.

\subsubsection{The stability property}

The following section is directly inspired by the work of Ayoub in
 \cite[\textsection 1.5]{ayoub}.\footnote{See
  also \cite[\textsection 5]{Delnotes}.}
 We claim no originality except for a closer look on the needed axioms.
\begin{df}
A \emph{pointed smooth $S$-scheme}
\index{word}{pointed, smooth $S$-scheme}
 will be a couple $(f,s)$ of morphisms of $\sch$
 such that $f:X \rightarrow S$ is a smooth separated morphism
  of finite type and $s:S \rightarrow X$ is a section of $f$.

We associate with a pointed smooth scheme $(f,s)$
 the following endofunctor of $\T(S)$
$$
\Th(f,s):=f_\sharp s_*
$$
called the associated \emph{Thom transformation}.
\index{word}{Thom!transformation}

If $\T$ satisfies \adjx s (recall: $s_*$ admits a right adjoint denoted by $s^!$),
 we put
$$
\Th'(f,s):=s^!f^*
$$
and call it the associated \emph{adjoint Thom transformation}.
\index{word}{Thom!adjoint transformation}
\end{df}

\begin{rem}
Note that because $f$ is separated, $s$ is a closed immersion.
\end{rem}

\begin{ex} \label{ex:Thom&etale}
\begin{enumerate}
\item Let $p:E \rightarrow X$ be a vector bundle and $s_0$ be its zero section.
Following \cite{ayoub}, we put $\Th(E):=\Th(p,s_0)$ and simply call it
 the Thom transformation associated with $E/X$.
\item Consider a pointed smooth $S$-scheme $(f,s)$ such that $f$ is \'etale.
 Then $s$ is an open and closed immersion. Thus, if $\T$ is additive,
 $s_*=s_\sharp$ according to Lemma \ref{lm:localization&sums}.
 In particular, $\Th(f,s)=Id_S$.
%\item Consider an arbitrary smooth separated morphism of finite type $f:X \rightarrow S$.
%We consider the cartesian square:
%\begin{equation} \label{eq:self_intersection}
%\xymatrix@=20pt{
%X \times_S X\ar^{f''}[r]\ar_{f'}[d] & X\ar^f[d] \\
%X\ar_f[r] & S
%}
%\end{equation}
%and $\delta:X \rightarrow X \times_S X$ the diagonal of $X/S$.
%We define \emph{suspension associated with $f$}
\end{enumerate}
\end{ex}

\begin{df}\label{df:stabilityppty}
We will say that $\T$ satisfies the \emph{stability property}, denoted by \stab,
 if for any point smooth scheme $(f,s)$, 
 the Thom transformation $\Th(f,s)$ is an equivalence of categories.
\end{df}

\begin{num}
Consider a commutative diagram in $\sch$ of the form
\begin{equation} \label{eq:hyp_Thom&composition}
\begin{split}
\xymatrix@R=18pt@C=26pt{
S\ar^{t}[rd]\ar_{t'}[d] & & \\
Y'\ar^/-2pt/{s'}[r]\ar_{p'}[d]\ar@{}|\Delta[rd]
 & Y\ar^{g}[rd]\ar^/2pt/{p}[d] & \\
S \ar_{s}[r] & X\ar_{f}[r] & S
}
\end{split}
\end{equation}
such that $\Delta$ is a cartesian square,
 $(f,s)$, $(g,t)$ are smooth pointed schemes
 and $g$ is a smooth separated morphism of finite type.
 Then we get a canonical exchange morphism:
\begin{equation} \label{eq:Thom&composition}
\Th(g,t)=f_\sharp p_\sharp s'_* t'_*
 \xrightarrow{\ Ex(\Delta_{\sharp*}) \ }
  f_\sharp s_* p'_\sharp t'_*=\Th(f,s) \Th(p',t').
\end{equation}
This is an isomorphism as soon as $Ex(\Delta_{\sharp*})$ is an isomorphism.
The following lemma gives a sufficient condition for this to happen.
\end{num}
\begin{lm}
Consider the above notations.
If $\T$ satisfies \locx s
 then the natural transformations $Ex(\Delta_{\sharp*})$
 is an isomorphism for any square $\Delta$ as above.
\end{lm}
This lemma follows easily from the definition of \locx s,
 the relations of paragraph \ref{num:general_notations_loc}
 and the $\Pmor$-base change formula \bc.
It motivates the next definition:
\begin{df}\label{df:weak_localization_ppty}
We say that $\T$ satisfies the \emph{weak localization property}
 \wloc if it satisfies \locx s for any closed immersion $s$
 which admits a smooth retraction.
\end{df}

\begin{prop}\label{prop:Nis-sep&wloc}
Assume that $\T$ satisfies the Nisnevich separation property.
Then the following conditions are equivalent:
\begin{enumerate}
\item[(i)] $\T$ satisfies \wloc.
\item[(ii)] For any scheme $S$ 
 and any closed immersion $i:Z \rightarrow X$ between smooth $S$-schemes, 
 $\T$ satisfies \locx i.
\end{enumerate}
\end{prop}
\begin{proof}
Of course, (ii) implies (i). We prove the reciprocal statement.
The Nisnevich separation property says that for any Nisnevich cover
 $f:X' \rightarrow X$, the functor $f^*$ is conservative.
 We deduce from that point the properties \locx i(a) and \locx i(b) 
 are local in $X$ with respect to the Nisnevich topology
 -- for (b), one also uses the smooth projection formula.
 Thus, we can conclude as locally for the Nisnevich topology,
 $i$ admits a smooth retraction (see for example \cite[4.5.11]{Deg7}).
\end{proof}

Applying the second point of Example \ref{ex:Thom&etale},
 we easily deduce from that construction the following
 kind of excision property:
\begin{lm}\label{lm:Thom&excision}
Assume that $\T$ satisfies \wloc.

Then, given any diagram \eqref{eq:hyp_Thom&composition}
 satisfying the assumption as above and such that $p$ is \'etale,
 the natural transformation \eqref{eq:Thom&composition} gives an isomorphism:
$$
\Th(g,t) \xrightarrow{\ \sim \ } \Th(f,s).
$$
\end{lm}

\begin{num} \label{num:thom_coherence}
To any short exact sequence of vector bundles over a scheme $S$
\begin{equation} \tag{$\sigma$}
0 \rightarrow E' \xrightarrow \nu E
 \xrightarrow \pi E'' \rightarrow 0,
\end{equation}
we can associate a commutative diagram
$$
\xymatrix@=18pt{
S\ar[rd]\ar[d] & & \\
E'\ar^\nu[r]\ar[d]\ar@{}|\Delta[rd]
 & E\ar[rd]\ar^\pi[d] & \\
S \ar[r] & E''\ar[r] & S
}
$$
where the non labelled map are either the canonical projections
 or the zero sections of the relevant vector bundles,
 and $\Delta$ is cartesian.
 Using the notation of Example \ref{ex:Thom&etale},
 the exchange transformation \eqref{eq:Thom&composition}
 associated with this diagram has the following form:
$$
\Th(\sigma):\Th(E) \longrightarrow \Th(E'') \circ \Th(E').
$$
Recall from the above that this natural transformation
 is an isomorphism
 as soon as $\T$ satisfies \wloc.
\end{num}

\begin{prop} \label{prop:caracterisation_stability}
Assume $\T$ satisfies \wloc and \sepx \zar.
Then the following conditions are equivalent:
\begin{enumerate}
\item[(i)] The complete triangulated $\sm$-fibred category $\T$
 satisfies the stability property.
\item[(ii)] For any scheme $S$,
 the Thom transformation $\Th(\AA^1_S)$
 is an equivalence of categories.
\end{enumerate}
\end{prop}
\begin{proof}
We have to prove that (ii) implies (i). 
Note that according to the above paragraph,
 we already now that for any scheme $S$ and any integer $n \geq 0$,
 $\Th(\AA^n_S) \simeq \Th(\AA^1_S)^{\circ,n}$ is an equivalence.

We consider a smooth pointed scheme
 $(f:X \rightarrow S,s)$ and we prove that $\Th(f,s)$ is an equivalence.

Recall that \locx s implies \adj s
 (first point of Proposition \ref{prop:loc_triangle}).
 In particular, $\Th(f,s)$ admits a right adjoint $\Th'(f,s)$
 and we have to prove that the adjunction morphisms are isomorphisms.

Consider an open immersion $j:U \rightarrow S$ and let
 $(f_0,s_0)$ be the restriction of the smooth $S$-point $(f,s)$ over $U$.
 Property \locx s implies \BCx s (Corollary \ref{cor:loc&Supp,BC,PF}).
 Thus, using also property \bc, we obtain a canonical isomorphism:
$$
j^* \Th(f,s) \xrightarrow{\ \sim \ } \Th(f_0,s_0) j^*.
$$
Recall also that \locx s implies \suppx s
 (again Corollary \ref{cor:loc&Supp,BC,PF}).
 Thus we get a canonical isomorphism:
$$
j_\sharp \Th(f_0,s_0) \xrightarrow{\ \sim \ } \Th(f,s) j_\sharp
$$
which gives by adjunction an isomorphism:
$$
\Th'(f_0,s_0) j^* \xrightarrow{\ \sim \ } j^* \Th'(f,s).
$$
Thus, \sepx \zar shows that the property for $\Th(f,s)$ to be an equivalence
 is Zariski local in $S$.

Consider a point $a \in S$, $x=s(a)$. As $X$ is smooth over $S$,
 there exists an open subscheme $U \subset X$, an integer $n\geq 0$
 and an \'etale $S$-morphism
 $\pi:U \rightarrow \AA^n_S$ which fits into the following cartesian square:
$$
\xymatrix@=10pt{
S_0\ar[r]\ar[d] & U\ar^\pi[d] \\
S\ar^\nu[r] & \AA^n_S
}
$$
where $\nu$ is the zero section (\emph{cf.} \cite[17.12.2]{EGA4}).
Note that the scheme $S_0=s^{-1}(U)$ is an open neighborhood of $a$ in $S$.
Let us put $X_0=f^{-1}(S_0)$ and $U_0=U \cap X_0$. Then we get the following
commutative diagram:
$$
\xymatrix@R=14pt@C=36pt{
& X_0\ar^{f_0}[rd] & \\
S_0\ar^{s_0}[ru]\ar_{\nu_0}[rd]\ar|-{s'_0}[r]
 & U_0\ar@{^(->}[u]\ar^{\pi_0}[d]\ar|-{f'_0}[r] & S_0 \\
& \AA^n_{S_0}\ar[ru] &
}
$$
where $\pi_0$ is the restriction of $\pi$ above $S_0$
 and $\nu_0$ is again the zero section.
According to Lemma \ref{lm:Thom&excision}, we get isomorphisms
$$
\Th(f_0,s_0) \simeq \Th(f'_0,s'_0) \simeq \Th(\AA^n_S).
$$
Thus, according to the beginning of the proof,
 $\Th(f_0,s_0)$ is an equivalence.
 This concludes because $S_0$ is an open neighborhood of $a$ in $S$.
\end{proof}

\begin{df} \label{df:premotive_smooth_pts}
Assume that $\T$ is monoidal.
\begin{enumerate}
\item For any smooth pointed scheme $(f:X \rightarrow S,s)$,
 we put $M_S\big(X/X-s(S)\big):=f_\sharp s_*(\un_S)$.
\item For any vector bundle $E/S$ with projection $f$
 and zero section $s$, we define the \emph{Thom premotive}
\index{word}{Thom!premotive}
 associated with $E$ over $S$ as $MTh_S(E)=f_\sharp s_*(\un_S)$.
\end{enumerate}
\end{df}

\begin{num} \label{num:monoidal&stability}
We assume $\T$ is monoidal and satisfies properties \wloc and \sepx \zar.

In each case of the previous definition,
 if we apply $f_\sharp$ to the distinguished triangle
 obtained from point (2) of Proposition \ref{prop:loc_triangle}
 applied to $s$, 
 we get the following canonical distinguished triangles:
\begin{align*}
M_S\big(X-s(S)\big) \rightarrow M_S(X)
 \rightarrow &M_S\big(X/X-s(S)\big) \xrightarrow{+1} \\
M_S(E^\times) \rightarrow M_S(E)
 \rightarrow & MTh_S(E) \xrightarrow{+1}
\end{align*}
where the first map is induced by the obvious open immersion.

Moreover, property \locx s implies \PFx s
 (see Corollary \ref{cor:loc&Supp,BC,PF}).
 Thus for any premotive $K$ over $S$,
 the following composite map is an isomorphism:
\begin{equation} \label{eq:Thom&monoidal}
\begin{split}
\Th(f,s).K &=f_\sharp s_*(K)=f_\sharp s_*(\un_S \otimes_S s^*f^*(K))
 \xrightarrow{Ex(s_*^*,\otimes)^{-1}} f_\sharp(s_*(\un_S) \otimes_X f^*(K)) \\
 &\xrightarrow{Ex(f_\sharp^*,\otimes)} (f_\sharp s_*(\un_S)) \otimes_S K
 =M_S(X/X-s(S)) \otimes_S K
\end{split}
\end{equation}
Similarly, in the case of a vector bundle $E/S$,
 we get a canonical isomorphism:
\begin{equation*}
\Th(E).K  \xrightarrow{\ \sim \ } MTh_S(E) \otimes_S K.
\end{equation*}
From these isomorphisms, we deduce easily the following
 corollary of the previous proposition:
\end{num}
\begin{cor} \label{cor:monoidal&stability}
Consider the above notations and assumptions.
Then the following properties are equivalent:
\begin{enumerate}
\item[(i)] $\T$ satisfies the stability property.
\item[(ii)] For any smooth pointed scheme $(X \rightarrow S,s)$,
 the premotive $M_S(X/X-s(S))$ is $\otimes$-invertible.
\item[(iii)] For any vector bundle $E/S$
 the Thom premotive $MTh_S(E)$ is $\otimes$-invertible.
\item[(iv)] For any scheme $S$,
 the premotive $MTh_S(\AA^1_S)$ is $\otimes$-invertible.
\end{enumerate}
\end{cor}

\begin{rem} \label{rem:Thom_sp&exact_sequences}
Assume that $\T$ satisfies the assumptions
 and the equivalent conditions of the previous corollary.
 Then, under the notations of Paragraph \ref{num:thom_coherence},
 we associate with the exact sequence $(\sigma)$ a canonical isomorphism
\begin{equation} \label{eq:Thom_sp&exact_sequences}
Th_S(\sigma):MTh_S(E) \rightarrow MTh_S(E'') \otimes_S MTh_S(E').
\end{equation}
Recall that Deligne introduced in \cite[4.12]{Del_det}
 the Picard category
\index{word}{Picard category}
  $\underline K(S)$ of \emph{virtual vector bundle}
\index{word}{bundle!virtual vector bundle}
 over a scheme $S$.

Then, it follows from the above isomorphism and the universal properties
 of $\underline K(S)$ (see \cite[4.3]{Del_det}) 
 that the functor $MTh_S$ can be extended uniquely 
 to a symmetric monoidal functor:
$$
MTh_S:\underline K(S) \rightarrow \T(S).
$$
The reader is referred to \cite[th. 1.5.18]{ayoub}
 for a detailed argument.
\end{rem}

\begin{num}
Assume $\T$ is monoidal.
For any scheme $S$, the canonical projection $p:\PP^1_S \rightarrow S$
 is a split epimorphism. A splitting is given by the inclusion
 of the infinite point $\nu:S \rightarrow \PP^1_S$.
The induced map $p_*:M_S(\PP^1_S) \rightarrow \un_S$
 is a split epimorphism.
 Thus it admits a kernel $K$ in the triangulated category $\T(S)$.
\end{num}
\begin{df} \label{df:Tate_twist}
Under the above assumption and notations,
 we define the \emph{Tate premotive}
\index{word}{premotive!Tate premotive}
  over $S$
 as the object $\un_S(1)=K[-2]$ of $\T(S)$.

The monoid generated by the cartesian section $(\un_S)_S$
 defines a canonical $\NN$-twist on $\T$ called the \emph{Tate twist}.
\index{word}{twist!Tate}
\index{word}{Tate!twist|see{twist}}
 The $n$-th Tate twist of an object $K$ is denoted by $K(n)$.
\end{df}

\begin{num}
Consider again the assumption of Paragraph \ref{num:monoidal&stability}.

According to Lemma \ref{lm:Thom&excision},
 we get a canonical isomorphism
$$
MTh_S(\AA^1_S)=M_S(\AA^1_S/\AA^1_S-\{0\}) \rightarrow M_S(\PP^1_S/\PP^1_S-\{0\}).
$$

On the other hand, $\un_S(1)[2]$ is by definition the cokernel
 of the monomorphism $\nu_*:\un_S \rightarrow M_S(\PP^1_S)$.
 Thus we get a canonical morphism:
\begin{equation} \label{eq:comparison_of_Tate_twists}
\un_S(1)[2] \rightarrow M_S(\PP^1_S/\PP^1_S-\{0\})
 \xrightarrow{\ \sim \ } MTh_S(\AA^1_S).
\end{equation}
From this definition
 and Corollary \ref{cor:monoidal&stability}
 the following result is obvious:
\end{num}
\begin{cor}
Consider the above assumption and notations.
Then the following conditions are equivalent:
\begin{enumerate}
\item[(i)] $\T$ satisfies the homotopy property.
\item[(ii)] For any scheme $S$,
 the arrow \eqref{eq:comparison_of_Tate_twists} is an isomorphism.
\end{enumerate}
When these equivalent assertions are satisfied,
 the following conditions are equivalent:
\begin{enumerate}
\item[(iii)] $\T$ satisfies the stability property.
\item[(iv)] For any scheme $S$,
 the Tate premotive $\un_S(1)$ is $\otimes$-invertible.
 \end{enumerate}
\end{cor}
%If $\un_S(1)$ is $\otimes$-invertible,
% we will consider the Tate twist on $\T$ as a $\ZZ$-twist.

\subsubsection{The purity property}

\begin{num} \label{num:proper_purity}
Let $f:X \rightarrow S$ be a smooth proper morphism in $\sch$.
 We consider the following cartesian square:
\begin{equation} \label{eq:auto-co-intersection}
\begin{split}
\xymatrix@=22pt{
X \times_S X\ar^-{f''}[r]\ar_{f'}[d]\ar@{}|\Delta[rd]
 & X\ar^f[d] \\
X\ar_f[r] & S
}
\end{split}
\end{equation}
where $f'$ (resp. $f''$) is the projection on the first (resp. second) factor.
Let $\delta:X \rightarrow X \times_S X$ be the diagonal embedding.
Note that $(f',\delta)$ is a smooth pointed scheme which depends only on $f$.
We put:
$$
\Sigma_f:=\Th(f',\delta)=f'_\sharp \delta_*.
$$
We then define a canonical morphism:
\begin{equation*}
\piso_f:f_\sharp=f_\sharp f''_* \delta_*
 \xrightarrow{Ex(\Delta_{\sharp *})}
 f_* f'_\sharp \delta_*=f_* \circ \Sigma_f
\end{equation*}
using the exchange transformation
 introduced in paragraph \ref{num:exchanges2}.
\end{num}
\begin{df} \label{df:pure_proper_morphism}
We say that $f$ is \emph{$\T$-pure},
\index{word}{morphism!$\T$-pure|seealso{morphism, pure}}
\index{word}{morphism!pure (proper case)}
 or simply \emph{pure} when $\T$ is clear,
 when the following conditions are satisfied:
\begin{enumerate}
\item The natural transformation $\Sigma_f$ is an equivalence.
\item The morphism $\piso_f:f_\sharp \rightarrow f_* \circ \Sigma_f$ is an isomorphism.
\end{enumerate}
Then $\piso_f$ is called
 the \emph{purity isomorphism}\index{word}{purity!isomorphism (relative)}
 associated with $f$.
We say also that $f$ is \emph{universally $\T$-pure}
\index{word}{morphism!universally $\T$-pure (proper case)}
 if $f$ is pure after any base change along a morphism of $\sch$.

We introduce the following properties on $\T$:
\begin{itemize}
\item $\T$ satisfies the \emph{purity property} \pur
 if any proper smooth morphism is pure.
\item $\T$ satisfies the \emph{weak purity property} \wpur
 if for any scheme $S$  and any integer $n>0$,
 the canonical projection $p_n:\PP^n_S \rightarrow S$ is pure.
\end{itemize}
\end{df}

\begin{rem} \label{rem:Omega_f}
Consider the above notations and assume $f$ is pure.

Then $f_*$ admits a right adjoint $f^!$
 and we deduce by transposition from $\piso_f$ a canonical isomorphism:
$$
\piso'_f:f^* \rightarrow \Sigma_f^{-1} \circ f^!.
$$
Recall also that, when $\delta_*$ admits a right adjoint $\delta^!$,
  $\Sigma_f$ admits as a right adjoint the transformation
 $\Omega_f:=\delta^!f^*$. In particular, $\Omega_f=\Sigma_f^{-1}$.
\end{rem}

The following lemma shows the importance of the purity property.
\begin{lm} \label{lm:purity=>BC,supp,PF}
Assume that $\T$ satisfies \wloc.
Let $f:Y \rightarrow X$
 be a proper smooth morphism.
 If $f$ is universally pure then the following conditions hold:
\begin{enumerate}
\item $\T$ satisfies \suppx f and \BCx f.
\item For any cartesian square
$$
\xymatrix@=14pt{
Z\ar^{\tilde f}[r]\ar_{h}[d]\ar@{}|\Delta[rd] & Y\ar^g[d] \\
X\ar_f[r] & S
}
$$
such that $g$ is smooth, the exchange transformation:
$$
Ex(\Delta_{\sharp*}):g_\sharp \tilde f_* \rightarrow f_* h_\sharp
$$
is an isomorphism.
\item If moreover $\T$ is monoidal then $\T$ satisfies \PFx f.
\end{enumerate}
\end{lm}
\begin{proof}
We first prove condition (2).
By assumption, the natural transformation $\Sigma_{\tilde f}$
 is an equivalence.
 for $f$ and $\tilde f$:
 by assumption the natural transformations $\Sigma_f=f'_\sharp \delta_*$
 and $\Sigma_{\tilde f}=\tilde f' \tilde \delta_*$) are equivalences.
 Thus, it is sufficient to prove that the natural transformation
$$
g_\sharp \tilde f_* \Sigma_{\tilde f}
 \xrightarrow{Ex(\Delta_{\sharp*})}
 f_* h_\sharp \Sigma_{\tilde f}
$$
is an isomorphism.

For matter of notations, let us also introduce
 the following cartesian squares:
$$
\xymatrix@R=18pt@C=24pt{
Z\ar^-{\tilde \delta}[r]\ar_h[d]\ar@{}|\Gamma[rd]
 & Z \times_{Y} Z\ar^-{\tilde f'}[r]\ar_k[d]\ar@{}|-\Theta[rd]
 & Z\ar^h[d] \\
X\ar_-\delta[r] & X \times_S X\ar_-{f'}[r] & X
}
$$
using the notations of \ref{num:proper_purity}.
Thus, by definition: $\Sigma_f=f'_\sharp \delta_*$,
 $\Sigma_{\tilde f}=\tilde f' \tilde \delta_*$.
Then we consider the following diagram of exchange transformations:
$$
\xymatrix@C=46pt{
g_\sharp \tilde f_\sharp\ar@{=}[d]\ar^{\piso_{\tilde f}}[rrr]
 &&& g_\sharp \tilde f_* \tilde f'_\sharp \tilde \delta_*\ar^{Ex(\Delta_{\sharp *})}[d] \\
f_\sharp h_\sharp\ar^-{\piso_f}[r]
 & f_* f'_\sharp \delta_* h_\sharp
 & f_* f'_\sharp k_\sharp \tilde \delta_*\ar_-{Ex(\Gamma_{\sharp *})}[l]
 & f_* h_\sharp \tilde f'_\sharp \tilde \delta_*\ar@{=}[l]
}
$$
Note that it only involves exchange transformations
 of type $Ex(?_{\sharp *})$:
 it is commutative by compatibility of these exchange transformations 
 with composition.
By assumption, the transformations $\piso_f$ and $\piso_{\tilde f}$
 are isomorphisms. Moreover the property \locx \delta is satisfied
 and it implies \suppx \delta according to Corollary \ref{cor:loc&Supp,BC,PF}.
 Thus  $Ex(\Gamma_{\sharp *})$ is an isomorphism and this concludes the proof of (2).

For condition (1), we note that (2) already implies \suppx f.
Thus we have only to prove \BCx f. We consider a square of shape
 $\Delta$
 as in the statement of the lemma without assuming that $g$ is smooth.
 We have to prove that
$$
Ex(\Delta_*^*):g^*f_* \rightarrow \tilde f_* h^*
$$
is an isomorphism. We proceed as for condition (2).
It is sufficient to prove that $Ex(\Delta_*^*)$ is an isomorphism
 after composition on the right with $\Sigma_f$.
Then we consider the following commutative diagram of exchange transformations:
$$
\xymatrix@C=46pt{
g^* f_\sharp\ar_{Ex(\Delta_\sharp^*)}[d]\ar^{\piso_{f}}[rrr]
 &&& g^* f_* f'_\sharp \delta_*\ar^{Ex(\Delta^*_{*})}[d] \\
\tilde f_\sharp h^*\ar^-{\piso_{\tilde f}}[r]
 & \tilde f_* \tilde f'_\sharp \tilde \delta_* h^*
 & \tilde f_* \tilde f'_\sharp k^* \delta_*\ar_-{Ex(\Gamma^*_{*})}[l]
 & \tilde f_* h^* f'_\sharp \delta_*\ar_-{Ex(\Theta^*_{\sharp})}[l]
}
$$
According to \bc, $Ex(\Delta_\sharp^*)$ and $Ex(\Theta_\sharp^*)$
 are isomorphisms.
By assumption, $\piso_f$ and $\piso_{\tilde f}$ are isomorphisms.
Moreover, property \locx \delta is satisfied and this implies $Ex(\Gamma_*^*)$
 is an isomorphism according to Corollary \ref{cor:loc&Supp,BC,PF}.
 Condition (1) is proved.

It remains to prove (3). We consider again the notations
 of the cartesian diagram \eqref{eq:auto-co-intersection}.
For any premotives $K$ over $X$ and $L$ over $S$,
 we consider the following commutative diagram
 of exchange transformations (see Remark \ref{rem:coherence_exchange4}):
$$
\xymatrix@C=60pt@R=18pt{
f_\sharp\big(K \otimes f^*(L)\big)\ar_{Ex(f_\sharp^*,\otimes)}[ddd]\ar^-{\piso_f}[r]
 & f_* f'_\sharp\delta_*\big(K \otimes \delta^*f^{\prime*}f^*(L)\big) \\
 & f_* f'_\sharp\big(\delta_*(K) \otimes f^{\prime*}f^*(L)\big)
 \ar^{Ex(f_\sharp^{\prime*},\otimes)}[d]\ar_{Ex(\delta_*^*,\otimes)}[u]  \\
 & f_* \big(f'_\sharp\delta_*(K) \otimes f^*(L)\big) \\
f_\sharp(K) \otimes L\ar^-{\piso_f}[r]
 & f_* f'_\sharp\delta_*(K) \otimes L.\ar_{Ex(f_*^*,\otimes)}[u] \\
}
$$
By definition,
 the exchanges $Ex(f_\sharp^*,\otimes)$ and $Ex(f_\sharp^{\prime*},\otimes)$
 are isomorphisms.
By assumption, the arrows labeled $\piso_f$ are isomorphisms.
Moreover, the property \locx \delta is satisfied: 
 Corollary \ref{cor:loc&Supp,BC,PF} implies that $Ex(\delta_*^*,\otimes)$
 is an isomorphism. We deduce from this that the arrow $Ex(f_*^*,\otimes)$
 is an isomorphism.
 This concludes the proof of (3) as
 the functor $\Sigma_f=f'_\sharp \delta_*$ is an equivalence
 according to the hypothesis on $f$. 
\end{proof}

\begin{num}
Assume that $\T$ satisfies the support property \supp.
Then we can extend Definition \ref{df:pure_proper_morphism} to the
 case of a smooth separated morphism of finite type $f:X \rightarrow S$.
We still consider the cartesian square \eqref{eq:auto-co-intersection}
 and the diagonal embedding $\delta:X \rightarrow X \times_S X$.
Again, $(f',\delta)$ is a smooth pointed scheme so that we can put
$$
\Sigma_f:=\Th(f',\delta)=f'_\sharp \delta_*
$$
and we define a canonical morphism:
\begin{equation} \label{eq:relative_purity_iso1}
\piso_f:f_\sharp=f_\sharp f''_! \delta_!
 \xrightarrow{Ex(\Delta_{\sharp\,!})}
 f_! f'_\sharp \delta_!=f_! \circ \Sigma_f.
\end{equation}
using the exchange transformation of point (2)
in Corollary \ref{prop:support_exists_f_!_plus}.
\end{num}
\begin{df} \label{df:pure_morphism}
Using the notations above,
 we say that $f$ is \emph{$\T$-pure}, or simply \emph{pure}
 \index{word}{morphism!pure}
 when $\T$ is clear,
 when the following conditions are satisfied:
\begin{enumerate}
\item The natural transformation $\Sigma_f$ is an equivalence.
\item The morphism $\piso_f:f_\sharp \rightarrow f_! \circ \Sigma_f$ is an isomorphism.
\end{enumerate}
\end{df}
We can easily deduce from the construction of the exchange transformation
 $Ex(\Delta_{\sharp\,!})$ that, when $\T$ satisfies properties \stab
  and \pur, any smooth separated morphism of finite type $f$ is pure.
The following theorem is a consequence of the formalism developed previously.
\begin{thm} \label{thm:localizationandwpur}
Assume that $\T$ satisfies the localization and weak purity properties.
Then the following conditions hold:
\begin{enumerate}
\item $\T$ satisfies the stability property.
\item $\T$ satisfies the support and base change properties. \\
If moreover $\T$ is monoidal,
 it satisfies the projection formula.
\item Any smooth separated morphism of finite type is pure.
\item For any projective morphism $f$, the property \adjx f holds. \\
If moreover $\T$ is well generated,
 then the adjoint property holds in general.
\end{enumerate}
\end{thm}
\begin{proof}
We start by proving condition (1). As \loc implies \sepx \zar,
 we can apply Proposition \ref{prop:caracterisation_stability}
 and we have only to prove that for any scheme $S$,
 $\Th(\AA^1_S)$ is an equivalence.
Let $s:S \rightarrow \AA^1_S$ be the zero section
 and $j:\AA^1_S \rightarrow \PP^1_S$ be the canonical open immersion.
 Put $t=j \circ s$. According to Lemma \ref{lm:Thom&excision},
 $j$ induces an isomorphism $\Th(\AA^1_S)\simeq \Th(p_1,s)$.
Consider now the following cartesian squares:
$$
\xymatrix@R=16pt@C=26pt{
S\ar^s[r]\ar_s[d]
 & \PP^1_S\ar^{p_1}[r]\ar^{s'}[d]\ar@{}|\Delta[rd] & S\ar^s[d] \\
\PP^1_S\ar_-\delta[r] & \PP^1_S \times_S \PP^1_S \ar_-{p'_1}[r] & \PP^1_S
}
$$
where $p'_1$ (resp. $\delta$) is the projection on the first factor
 (resp. diagonal embedding).
The property \locx s implies that $s^*s_*=1$
 and that the exchange transformation $Ex(\Delta_{\sharp *})$
 is an isomorphism according to Corollary \ref{cor:loc&Supp,BC,PF}.
Thus we get an isomorphism of functors:
$$
\Th(p_1,s)=p_{1 \sharp} s_*=s^*s_*p_{1 \sharp} s_*
 \xrightarrow{Ex(\Delta_{\sharp *})^{-1}} s^*p'_{1 \sharp} s'_* s_*
 =s^*p'_{1 \sharp} \delta_* s_*=s^* \Sigma_{p_1} s_*
$$
and this proves (1) because $p_1$ is pure.

Condition (2) follows simply from Corollary \ref{cor:loc&Supp,BC,PF}.
In fact, for any scheme $S$,
 the weak purity assumption on $\T$ implies 
  that $p_n:\PP^n_S \rightarrow S$ is universally pure.
 Thus, Lemma \ref{lm:purity=>BC,supp,PF} implies
  properties \suppx{p_n} and \BCx{p_n}
  so that we can apply Corollary \ref{cor:loc&Supp,BC,PF} to get \supp and \BC.
 The same argument applies to the property \PF in the monoidal case.

For condition (3),
 we consider a smooth separated morphism of finite type $g:Y \rightarrow S$
  and we prove it is pure.
 According to (1), $\Sigma_g$ is an equivalence.
 Thus, by definition of $\piso_g$, it is sufficient to prove that
 for any cartesian square:
$$
\xymatrix@=14pt{
Z\ar^{\tilde f}[r]\ar_{h}[d]\ar@{}|\Delta[rd] & Y\ar^g[d] \\
X\ar_f[r] & S
}
$$
with $f$ separated of finite type, the exchange transformation
$$
Ex(\Delta_{\sharp\,!}):g_\sharp \tilde f_! \rightarrow f_! h_\sharp
$$
is an isomorphism. \\
To do this, we apply Proposition \ref{prop:morph&loc1},
 as in the case of Corollary \ref{cor:loc&Supp,BC,PF}.
 We consider the obvious complete $\sm$-fibred triangulated categories
 $\T'$ and $\T''$ over $\sch/S$ which to an $S$-scheme $Y$ associates:
\begin{itemize}
\item $\T'(Y)=\T(Y \times_S X)$.
\item $\T''(Y)=\T(Y)$.
\end{itemize}
We consider the morphism $\varphi^*:\T' \rightarrow \T''$
 such that for any $S$-scheme $Y$, $\varphi_Y^*=(Y \times_S p)_\sharp$.
 As for any scheme $S$, $p_n:\PP^n_S \rightarrow S$ is universally pure,
 Lemma \ref{lm:purity=>BC,supp,PF} shows that 
 $\varphi^*$ satisfies condition (i) of Proposition \ref{prop:morph&loc1}.
 According to that Proposition, (i) is equivalent to condition (iii),
  and (iii) is precisely what we want.

It remains only to prove condition (4).
According to property \pur,
 any smooth proper morphism $f$ satisfies \adjx f. 
According to \loc and Proposition \ref{prop:loc_triangle}
 any closed immersion $i$ satisfies \adjx i.
 It follows easily that any projective morphism $f$
 satisfies \adjx f.
 When $\T$ is well generated,
  we simply apply point (4) of Corollary \ref{cor:loc&Supp,BC,PF}.
\end{proof}

\begin{rem}
In particular,
 in the assumption of the previous theorem,
 if $\T$ satisfies properties \loc, \wpur and \adj\footnote{Note
 that under the assumptions of the previous theorem,
 we know that for any proper smooth morphism $f$,
  $f_*$ admits a right adjoint. The same is true for a proper morphism
  which can be factorized as a closed immersion followed by a smooth proper
  morphism according to \loc.},
  we can apply Theorem \ref{thm:support} to $\T$ so that we get
  a complete formalism of operations $(f^*,f_*,f_!,f^!)$
  satisfying all the desired formulas.

Thus the preceding theorem gives another look
 at the main result of \cite[1.4.2]{ayoub}.
In fact,
 the proof given here is simpler
  as the assumptions of our theorem are stronger.
However, we do not use the homotopy property in our theorem.
\end{rem}

We end up this section with a theorem due to Ayoub \cite[1.4.2]{ayoub}.
The particular case $\T(X)=\SH(X)$ was also established by R\"ondigs
in \cite{roendigs}, after Voevodsky,
with a proof which extends immediately to Ayoub's axiomatic setting.
It may be stated in a simpler form, according to theorem \ref{thm:localizationandwpur} above:
\begin{thm}[Voevodsky-R\"ondigs-Ayoub]
Assume $\T$ satisfies the localization,
homotopy and stability properties.
Then $\T$ is weakly pure.\label{thm:Ayoub_wpur}
\end{thm}
In fact, this theorem is stated explicitly in \cite[Theorem 1.7.9]{ayoub}.

\begin{rem}
Recall that Ayoub proves more than just this theorem:
 indeed he constructs the whole formalism of the six functors for
 quasi-projective morphisms for his \emph{monoidal homotopy stable functors}
 --- see again \cite{ayoub}. Similarly, the fact that one can deduce the
 proper base change formula from relative purity was also observed
 by R\"ondigs \cite{roendigs}.
The work we have done here is to isolate the crucial properties of purity
 and weak purity. Also, using the construction of Deligne, we see
 how to avoid the assumption of quasi-projectiveness made by Ayoub.
 Finally, the interest of Theorem \ref{thm:localizationandwpur} is to give a possible
 approach to the \emph{six functors formalism} without requiring the homotopy
 property.
\end{rem}

\subsubsection{Duality, purity and orientation}
\begin{num}\label{num:recall_duality}
This section is concerned with the relation between
 purity and duality.
 We will assume that $\T$ is premotivic.

Recall that an object $M$ of a monoidal category $\mathscr M$ 
 is called \emph{strongly dualizable}
\index{word}{dualizable, strongly}\index{word}{strongly dualizable|see{dualizable}}
  if there exists an object $M'$ such that $(M' \otimes -)$
 is both right and left adjoint to $(M \otimes -)$.
 Then, $M'$ is called the \emph{strong dual}
\index{word}{dual, strong}
  of $M$.

In case $\mathscr M$ is closed monoidal,
 we will say that a morphism of the form
$$
\mu:M \otimes M' \rightarrow \un
$$
is a \emph{perfect pairing}
\index{word}{perfect pairing}
 if the natural transformation
$$
(M \otimes -) \rightarrow \uHom(M',-)
$$
obtained from $\mu$ by adjunction is an isomorphism.
 Then $M$ is strongly dualizable with dual $M'$.
\end{num}
\begin{prop}\label{prop:purity&duality}
Let $f:X \rightarrow S$ be a smooth proper morphism.
If $f$ is pure\index{word}{morphism!pure}
 then the premotive $M_S(X)$ is strongly dualizable in $\T(S)$
 with dual:
$$
f_*(\un_X) \simeq f_\sharp\big(\Omega_f(\un_X)\big)
$$
where $\Omega_f$ denotes the inverse of $\Sigma_f$.
\end{prop}
\begin{proof}
By assumption, $\Sigma_f$ is an automorphism of the category
 $\T(X)$. Moreover, the identification \eqref{eq:Thom&monoidal} 
 can be rewritten as $\Sigma_f(M)=\Sigma_f(\un_X) \otimes_X M$
 for any premotive $M$ over $X$.
 The fact $\Sigma_f$ is an equivalence means that $\Sigma_f(\un_X)$
 is a $\otimes$-invertible object,
 whose inverse is $T:=\Omega_f(\un_S)$.
 In particular, we get: $\Omega_f(M)=T \otimes M$.

According to the $\sm$-projection formula,
 the functor $M_S(X) \otimes .$ is isomorphic to $f_\sharp f^*$.
 Thus, its right adjoint is $f_*f^*$. As $f$ is pure by assumption,
 this last functor is isomorphic to $f_\sharp \Omega_f f^*$.
Using the observation at the beginning of the proof
 and the $\sm$-projection formula again, we obtain:
$$
f_\sharp \Omega_f f^*(N)=f_\sharp(T \otimes f^*(N))
 =f_\sharp(T) \otimes N.
$$
Moreover, the right adjoint of $f_\sharp \Omega_f f^*$
 is $f_*\Sigma_f f^*$. Using again the purity isomorphism for $f$, 
 this last functor can be identified with $f_\sharp f^*$ and this concludes.
\end{proof}

\begin{num}
Assume again that the premotivic triangulated category $\T$ 
 satisfies properties \wloc and \sepx \nis.

Let $S$ be a scheme.
A \emph{smooth closed $S$-pair} will be pair $(X,Z)$
 of smooth $S$-schemes such that $Z$ is closed subscheme of $X$.
 We consider the canonical projection $p:X \rightarrow S$
 and the immersion $i:Z \rightarrow X$ associated with $(X,Z)$.
 Note that according to Proposition \ref{prop:Nis-sep&wloc},
 $\T$ satisfies property \locx i.
 Then we define the premotive of $(X,Z)$ as follows:
\begin{equation} \label{eq:M_with_support&i_*}
M_S(X/X-Z):=p_\sharp i_*(\un_Z).
\end{equation}
According to property \locx i,
 we thus get a canonical distinguished triangle:
\begin{equation} \label{eq:localization&sm_pairs}
M_S(X-Z) \xrightarrow{j_*} M_S(X) \rightarrow M_S(X/X-Z)
 \xrightarrow{+1}
\end{equation}
Note that given any smooth morphism $p:S \rightarrow S_0$,
 we get obviously:
\begin{equation} \label{eq:identification1_mot_rel}
p_\sharp M_S(X/X-Z)=M_{S_0}(X/X-Z).
\end{equation}
Moreover, given any morphism $f:T \rightarrow S$,
 we get an exchange isomorphism:
\begin{equation} \label{eq:identification2_mot_rel}
f^*M_S(X/X-Z) \xrightarrow{\sim} M_T(X_T/X_T-Z_T).
\end{equation}
A morphism of smooth closed $S$-pairs $(Y,T) \rightarrow (X,Z)$
 will be a couple $(f,g)$ which fits into a commutative diagram
$$
\xymatrix@=14pt{
T\ar^k[r]\ar_g[d]\ar@{}|\Delta[rd] & Y\ar^f[d] \\
Z\ar_i[r] & X,
}
$$
with $i$, $k$ the canonical immersions,
 and such that $T=f^{-1}(Z)$ as a set.
We can associate with $(f,g)$ a morphism of premotives:
\begin{align*}
M_S(Y/Y-T)& =q_\sharp k_* g^*(\un_Z) 
\xrightarrow{Ex(\Delta_*^*)^{-1}}  q_\sharp f^* i_*(\un_Z)
 \xrightarrow{Ex_\sharp^*} p_\sharp i_*(\un_Z)=M_S(X/X-Z).
\end{align*}
Indeed, the exchange map $Ex(\Delta_*^*)$ is an isomorphism
 according to \locx i and Corollary \ref{cor:loc&Supp,BC,PF}.

It is easy to check that the triangle \eqref{eq:localization&sm_pairs}
 is functorial with respect to morphisms of closed $S$-pairs.
 Before proving the next theorem,
  we state the following lemma.
\end{num}
\begin{lm}
Consider the assumptions and notations above.

Let $(f,g):(Y,T) \rightarrow (X,Z)$ be a morphism
 of smooth closed $S$-pairs such that $f$ is \'etale
 and $g$ is an isomorphism.
Then the induced map $M_S(Y/Y-T) \rightarrow M_S(X/X-Z)$
 is an isomorphism.
\end{lm}
\begin{proof}
According to the identification \ref{eq:identification1_mot_rel},
 it is sufficient to treat the case where $X=Z$.
 Let $U=X-Z$ and $j:U \rightarrow X$ be the obvious immersion.
 Then $(f,j)$ is a Nisnevich cover of $X$. According to \sepx \nis,
 it is sufficient to prove that the pullback of 
 $M_X(Y/Y-T) \rightarrow M_X(X/X-Z)$ along $f$ and $j$ is an isomorphism.
 This is obvious using \ref{eq:identification2_mot_rel}.
\end{proof}

\begin{num}
We consider again the assumption of the paragraph preceding
 the above lemma.

Fix a smooth closed $S$-pair $(X,Z)$.
Let $B_ZX$ (resp. $B_Z(\AA^1_X)$ be the blow-up of $X$
 (resp. $\AA^1_X$)
 with center in $Z$ (resp. $\{0\} \times Z$).
 We define the deformation space
\index{word}{deformation space}
  associated with $(X,Z)$ as the $S$-scheme $D_ZX=B_Z(\AA^1_X)-B_ZX$.
 Note also $D_ZZ=\AA^1_Z$ is a closed subscheme of $D_ZX$ ;
 the couple $(D_ZX,\AA^1_Z)$ is a smooth closed $S$-pair. \\
Let $N_ZX$ be the normal bundle
\index{word}{bundle!normal}
 of $Z$ in $X$.
The scheme $D_ZX$ is fibred over $\AA^1$.
 Moreover, the $0$-fiber of $(D_ZX,\AA^1)$ is the closed pair
 $(N_ZX,Z)$ corresponding to the zero section
 and the $1$-fiber is the closed pair $(X,Z)$.
 In particular, we get the following morphisms of closed pairs:
\begin{equation}
(X,Z) \xrightarrow{d_1} (D_ZX,\AA^1_Z) \xleftarrow{d_0} (N_ZX,Z)
\end{equation}
We are now ready to state the purity theorem
 for smooth closed pairs in our abstract formalism. 
 Though our assumptions are more general, 
 this theorem follows exactly from the method of Morel and Voevodsky
 used to prove this result in
  the homotopy category\index{word}{homotopy category} $\H$
 (see \cite[\textsection 3, 2.24]{MV}):
\end{num}
\begin{thm} \label{thm:purity_closed_pairs}
Consider the above assumptions and notations
 and suppose that $\T$ satisfies the homotopy property.
Then the morphisms
$$
M_S(X/X-Z) \xrightarrow{d_{1*}} M_S(D_ZX/D_ZX-\AA^1_Z)
 \xleftarrow{d_{0*}} M_S(N_ZX/N^\times_ZX)=:MTh_S(N_ZX).
$$
are isomorphisms.
\end{thm}
\begin{proof}
By noetherian induction and the preceding lemma,
 the statement is local in $X$ for the Nisnevich topology.
Thus, because $(X,Z)$ is a smooth closed $S$-pair,
 we can assume that there exists an \'etale map $\pi:X \rightarrow \AA^{n+c}_S$
 such that $\pi^{-1}(\AA^c_S)=Z$ -- \emph{cf.} \cite[17.12.2]{EGA4}.
Consider the pullback square
$$
\xymatrix@R=10pt@C=36pt{
X'\ar^p[r]\ar_q[d] & X\ar^\pi[d] \\
{}\AA^n \times Z\ar^-{1 \times \pi|Z}[r] & {}\AA^n \times \AA^c_S. }
$$
There is an obvious closed immersion $Z \rightarrow X'$
and its image is contained in $q^{-1}(Z)$.
As $q$ is \'etale, $Z$ is a direct factor of $q^{-1}(Z)$.
Put $W=q^{-1}(Z)-Z$ and $\Omega=X'-W$. 
Thus $\Omega$ is an open subscheme of $X'$, 
and the reader can check that
 $p$ and $q$ induces morphisms of smooth closed $S$-pairs
$$
(X,Z) \leftarrow (\Omega,Z) \rightarrow (\AA^n_Z,Z).
$$
Applying again the preceding lemma,
 these morphisms induces isomorphisms on the associated premotives.
 Thus we are reduced to the case of the closed $S$-pair
 $(\AA^n_Z,Z)$. 
 A direct computation shows that $D_Z(\AA^n_Z) \simeq \AA^1 \times \AA^n_Z$.
 Under this isomorphism $d_0$ (resp. $d_1$) corresponds
 to the $0$-section (resp. $1$-section) of $\AA^1 \times \AA^n_Z$
 corresponding to the first factor. 
 Thus, we conclude using the homotopy property.
\end{proof}

\begin{num} \label{num:purity&def_normal}
The interest of the previous theorem is to simplify the purity isomorphism.
 Let us restate the assumptions
 on the triangulated premotivic category $\T$:
\begin{itemize}
\item $\T$ satisfies properties \sepx \nis, \wloc and \htp.
\end{itemize}
Then applying the above theorem,
 we get for any smooth closed $S$-pair $(X,Z)$ a canonical isomorphism
\begin{equation} \label{eq:pur_closed_pairs}
\piso_{X,Z}:M_S(X/X-Z) \rightarrow MTh_S(N_ZX)
\end{equation}
\end{num}
\begin{cor} \label{cor:thm:purity_closed_pairs}
Consider the assumptions and notations above.
\begin{enumerate}
\item For any smooth pointed $S$-scheme $(f,s)$
 and any premotive $K$ over $S$, we get a canonical isomorphism
$$
\Th(f,s).K \simeq M_S(X/X-s(S)) \otimes_S K
 \xrightarrow{\piso_{X,S}} MTh_S(N_s) \otimes_S K.
$$
where the first isomorphism is given by the map \eqref{eq:Thom&monoidal}
 and $N_s$ is the normal bundle\index{word}{bundle!normal} of $s$.
\item For any smooth separated morphism of finite type
 $f:X \rightarrow S$ with tangent bundle\footnote{We define $T_f$ as
 the normal bundle of the diagonal immersion $\delta:X \rightarrow X \times_S X$.}
 \index{word}{bundle!tangent} $T_f$, and any premotive $K$ over $X$,
  we get a canonical isomorphism:
$$
\piso_{XX,X}:
\Sigma_f(K) \xrightarrow{\ \sim\ } MTh_X(T_f) \otimes_X K
$$
--- here, $(XX,X)$ stands for the closed pair corresponding
 to the diagonal embedding of $X/S$.
\end{enumerate}
\end{cor}
In the assumption of point (2),
 we thus get a canonical map:
\begin{equation} \label{eq:relative_purity_after_def}
f_\sharp(K) \xrightarrow{\piso_f} f_!(\Sigma_f K)
 \xrightarrow{\ \sim\ } f_!\big(MTh_X(T_f) \otimes_X K\big)
\end{equation}
that we will still denote by $\piso_f$ and call
 the \emph{purity isomorphism}\index{word}{purity!isomorphism (relative)}
 associated with $f$.

\begin{df} \label{df:premotivic_orientation}
Assume the triangulated premotivic category $\T$ satisfies \wloc.
As usual, $M(1)$ denotes the Tate twist of a premotive $M$.

An \emph{orientation}\index{word}{orientation!of a triangulated premotivic category}
 $\mathfrak t$ of $\T$ will be the data 
 for each smooth scheme $X$ and each vector bundle $E/X$
 of rank $n$ of an isomorphism
$$
\mathfrak t_E:MTh_X(E) \rightarrow \un_X(n)[2n],
$$
called the \emph{Thom isomorphism},\index{word}{Thom!isomorphism}
 satisfying the following coherence properties:
\begin{enumerate}
\item[(a)] Given a scheme $X$ and an isomorphism of
 vector bundles $\varphi:E \rightarrow F$ of ranks $n$ over $X$,
 the following diagram is commutative:
$$
\xymatrix@=16pt@C=30pt{
MTh_X(E)\ar_{\mathfrak t_E}[rd]\ar^{\varphi_*}[rr]
 && MTh_X(F).\ar^{\mathfrak t_F}[ld] \\
 & \un_X(n)[2n] &
}
$$
\item[(b)] For any morphism $f:Y \rightarrow X$ of schemes,
 and any vector bundle $E/X$ of rank $n$ with pullback $F$ over $Y$,
 the following diagram commutes:
$$
\xymatrix@=16pt@C=30pt{
f^*(MTh_X(E))\ar_\sim[d]\ar^{f^*\mathfrak t_E}[r] & f^*(\un_X(n)[2n])\ar^\sim[d] \\
MTh_Y(F)\ar^{\mathfrak t_F}[r] & \un_Y(n)[2n]
}
$$
where the vertical maps are the canonical isomorphisms.
\item[(c)] For any scheme $X$ and any exact sequence $(\sigma)$ of vector bundles over $X$
\begin{equation*}
0 \rightarrow E' \xrightarrow \nu E
 \xrightarrow \pi E'' \rightarrow 0,
\end{equation*}
if $n$ (resp. $m$) denotes the rank of the vector bundle $E'$ (resp. $E''$),
 the following diagram commutes:
$$
\xymatrix{
MTh_X(E)\ar_{\mathfrak t_E}[d]\ar^-{Th_X(\sigma)}[r]
 & MTh_X(E') \otimes MTh_X(E'')\ar^{\mathfrak t_{E'} \otimes \mathfrak t_{E''}}[d] \\
\un_X(n+m)[2n+2m]\ar[r] & \un_X(n)[2n] \otimes \un_X(m)[2m]
}
$$
where the map $Th_X(\sigma)$ is the isomorphism \eqref{eq:Thom_sp&exact_sequences}
 associated with $(\sigma)$ and the bottom vertical one is the obvious identification.
\end{enumerate}
We will also say that $\T$ is \emph{oriented}
 when the choice of one particular orientation 
 is not essential.
\end{df}
Note that the Thom isomorphism can be viewed as a cohomology class in 
$$
H^{2n,n}_\T(Th_X(E)):=\Hom_{\T(X)}\big(MTh_X(E),\un_S(n)[2n]\big)
$$
which in classical homotopy theory is called the \emph{Thom class}.
\index{word}{Thom!class}

\begin{num} \label{num:Gysin&purity_or}
Suppose the triangulated premotivic category $\T$ satisfies the following properties:
\begin{itemize}
\item $\T$ satisfies properties \sepx \nis, \wloc, \htp.
\item $\T$ admits an orientation $\mathfrak t$.
\end{itemize}
Consider a smooth closed $S$-pair $(X,Z)$ of codimension $n$.
Let $p$ (resp. $q$) be the structural morphism of $X/S$ (resp. $Z/S$)
 and $i:Z \rightarrow X$ the associated immersion.
 Then we associate with $(X,Z)$
 the following form of
  the purity isomorphism:\index{word}{purity!isomorphism (relative)}
\begin{equation} \label{eq:pur_closed_pairs_or}
\piso^\mathfrak t_{X,Z}:M_S(X/X-Z) \xrightarrow{\piso_{X,Z}} MTh_S(N_ZX)
 \xrightarrow{q_\sharp(\mathfrak t_{N_ZX})} M_S(Z)(n)[2n]
\end{equation}
where $\piso_{X,Z}$ is the isomorphism \eqref{eq:pur_closed_pairs}.
 For future reference, note that we deduce from this
  the so-called Gysin morphism:
\begin{equation} \label{eq:Gysin}
i^*:M_S(X) \xrightarrow{\pi} M_S(X/X-Z)
 \xrightarrow{\piso^\mathfrak t_{X,Z}} M_S(Z)(n)[2n]
\end{equation}
where $\pi$ is the following map:
$$
M_S(X)=p_\sharp(\un_X) \xrightarrow{ad(i^*,i_*)} p_\sharp i_* i^*(\un_X)=M_S(X/X-Z).
$$
As a particular case,
 we get using the notation of Corollary \ref{cor:thm:purity_closed_pairs}, point (2),
 an isomorphism:
$$
\piso_{XX,X}^\mathfrak t:\Sigma_f(K)
 \xrightarrow{\piso_{XX,X}} MTh_X(T_f) \otimes K
 \xrightarrow{\mathfrak t_{T_f}} K(d)[2d]
$$
In particular, when $\T$ satisfies property \supp,
 the purity comparison map associated with $f$ can be rewritten as:
\begin{equation} \label{eq:relative_purity_iso_or}
\piso_f^\mathfrak t:f_\sharp \xrightarrow{\ \piso_f\ } f_! \circ \Sigma_f
 \xrightarrow{\ \piso_{XX,X}^\mathfrak t\ } f_!(d)[2d]
\end{equation}
\end{num}

\begin{ex} \label{ex:or&GysinII}
Assume as in the above definition that $\T$ is premotivic
 and satisfies properties \wloc and \sepx \nis.

We suppose the following two additional conditions are fulfilled:
\begin{enumerate}
\item[(a')] There exists a morphism of triangulated premotivic categories:
$$
\varphi^*:\SH \leftrightarrows \T:\varphi_*
$$
where $\SH$ is the stable homotopy category of Morel and Voevodsky
 --- see Example \ref{ex:H&SH_premotivic}.
\item[(b')] For any scheme $X$, let $\pic(X)$ be the Picard
 group of $X$. We assume there exists an application
$$
c_1:\pic(X) \rightarrow H_\T^{2,1}(X):=\Hom_{\T(X)}(M(X),\un_X(1)[2])
$$
which is natural with respect to contravariant functoriality
--- we do not require $c_1$ is a morphism of abelian groups.
\end{enumerate}
Then one can apply the results of \cite{Deg8} to $\T(X)$ for any scheme $X$.
All the references which follows will be within \emph{loc. cit.}:
 according to section 2.3.2, the triangulated category $\T(X)$ satisfies
 the axioms of Paragraph 2.1.\footnote{Note in particular
  that for any smooth closed $S$-pair, we obtain a canonical isomorphism
  in $\T(S)$ of the form:
  $$
  \varphi^*(\Sigma^\infty X/X-Z)\simeq M_S(X/X-Z)
  $$
  where one the left-hand side $X/X-Z$ stands for the homotopy cofiber
  of the open immersion $(X-Z) \rightarrow X$ while the left-hand side 
  is defined by Equality \eqref{eq:M_with_support&i_*}.}
 Then the existence of the Thom isomorphism
 follows from Proposition 4.3 and, more explicitly, from Paragraph 4.4.
Property (a) and (b) of the above definition are easy --- explicitly,
 this is a consequence of 4.10 --- and Property (c) follows from Lemma 4.30.

To sum up, the assumptions (a') and (b') guarantees the existence
 of a canonical orientation\index{word}{orientation!of a triangulated premotivic category}
  of $\T$ in the sense of the above definition.
 Moreover,
  the purity isomorphism \eqref{eq:pur_closed_pairs_or}
  as well as the Gysin morphism \eqref{eq:Gysin} associated in the preceding
  paragraph for this particular orientation 
  coincide with the one defined in \cite{Deg8} (see in particular
  the uniqueness statement of \cite[Prop. 4.3]{Deg8}).

Note moreover that assuming $\T$ satisfies all the properties above
 except (b'), the data of an orientation of $\T$ is equivalent
 to the data of a map $c_1$ as in (b').
 Indeed, if $\mathfrak t$ is an orientation of $\T$,
 given any line bundle $L/X$ with zero section $s$,
 we put $c_1(L)=\rho(\mathfrak t_L)$
 where $\rho$ is the following composite map:
$$
H^{2,1}_\T(Th_X(L)) \rightarrow H^{2,1}_\T(L)
 \xrightarrow{s^*} H^{2,1}_\T(X)
$$
where the first map is induced by the canonical projection
 $M_X(L) \rightarrow MTh_X(L)$.
 Then $c_1$ depends only on the isomorphism classes of $L/X$
 --- property (a) of the above definition ---
 and it is compatible with pullbacks --- property (c) of the above definition.
\end{ex}

\begin{num}
We now assume the following conditions
 on the triangulated premotivic category $\T$:
\begin{itemize}
\item $\T$ satisfies properties \sepx \nis, \wloc, \htp and \stab.
\item $\T$ admits an orientation\index{word}{orientation} $\mathfrak t$.
\end{itemize}
Let  $f:X \rightarrow S$ be a smooth proper morphism of dimension $d$.
Note we do not need that $\T$ satisfies property \supp
 to rewrite the purity comparison map as follows:
\begin{equation} \label{thm:loc_wpur_or}
\piso_f^\mathfrak t:f_\sharp \rightarrow f_*(d)[2d]
\end{equation}
(see Paragraph \ref{num:Gysin&purity_or}).

Note also that using the Gysin morphism \eqref{eq:Gysin}
 associated with the diagonal immersion $\delta:X \rightarrow X \times_S X$,
 we get the following morphism:
\begin{equation}\label{eq:dual_pairing_or}
\mu_f^\mathfrak t:M_S(X) \otimes M_S(X)(-d)[-2d]=M_S(X \times_S X)(-d)[-2d]
 \xrightarrow{\ \delta^*\ } M_S(X) \xrightarrow{\ f_*\ } \un_S.
\end{equation}
\end{num}
\begin{thm} \label{thm:purity&duality_or}
Consider the assumptions and notations above.
Then the following conditions are equivalent:
\begin{enumerate}
\item[(i)] $f$ is pure\index{word}{morphism!pure}: $\piso_f$ is an isomorphism.
\item[(i')] The natural transformation $\piso_f.f^*$ is an isomorphism.
\item[(ii)] The premotive $M_S(X)$ is strongly dualizable\index{word}{dualizable, strongly}
 and $\mu_f^\mathfrak t$ is a perfect pairing.
\end{enumerate}
\end{thm}
\begin{proof}
In this proof, we put $\tau(K)=K(d)[2d]$.
As $\T$ satisfies property \stab,
 $f_*$ commutes with Tate twist (def. \ref{df:f_*_commutes_twists}).
 This means the following exchange transformation is an isomorphism:
%\begin{equation} \label{eq:twist&f^*,f_!}
%f^* \tau =\tau f^*, f_\sharp \tau=\tau f_\sharp.
%\end{equation}
\begin{equation} \label{eq:twist&f_*}
Ex_\tau:\tau f_* \rightarrow f_* \tau.
\end{equation}

We first prove that (i) is equivalent to (i').
 One implication is obvious so that we have only to prove that (i') implies (i).
 Guided by a method of Ayoub (see \cite[1.7.14, 1.7.15]{ayoub}),
 we will construct a right inverse $\phi_1$ and a left inverse $\phi_2$
 to the morphism $\piso^\mathfrak t_f$ as the following composite maps:
\begin{align*}
\phi_1:&f_* \tau
 \xrightarrow{ad(f^*,f_*)} f_*f^*f_* \tau 
 \xrightarrow{Ex_\tau^{-1}} f_* f^*\tau f_*=f_*\tau f^*f_*
 \xrightarrow{(\piso^\mathfrak t_f.f^*f_*)^{-1}} f_\sharp f^* f_*
 \xrightarrow{ad'(f^*,f_*)} f_\sharp \\
\phi_2:&f_* \tau
 \xrightarrow{\beta_f} f_* \tau f^*f_\sharp
 \xrightarrow{(\piso^\mathfrak t_f.f^*f_\sharp)^{-1}} f_\sharp f^* f_\sharp
 \xrightarrow{ad'(f_\sharp,f^*)} f_\sharp.
\end{align*}

Let us check that $\piso^\mathfrak t_f \circ \phi_1=1$. To prove this relation,
 we prove that the following diagram is commutative:
$$
\xymatrix@C=28pt{
f_*\tau\ar^-{ad(f^*,f_*)}[r]\ar@{=}[ddd]
 & f_*f^*f_*\tau\ar^-{Ex_\tau^{-1}}[r]\ar@{=}[dd]
 & f_*\tau f^*f_*\ar^-{(\piso^\mathfrak t_ff^*f_*)^{-1}}[r]\ar@{=}[d]
 & f_\sharp f^*f_*\ar^-{ad'(f^*,f_*)}[r]\ar@{=}[d]\ar@{}|{(1)}[rrd]
 & f_\sharp\ar^-{\piso^\mathfrak t_f}[r]
 & f_*\tau\ar@{=}[d] \\
 & 
 & f_*\tau f^*f_*\ar^-{(\piso^\mathfrak t_ff^*f_*)^{-1}}[r]\ar@{=}[d]\ar@{}|{(2)}[rrrd]
 & f_\sharp f^*f_*\ar^-{\piso^\mathfrak t_ff^*f_*}[r]
 & f_*\tau f\*f_*\ar^-{ad'(f^*,f_*)}[r]
 & f_*\tau\ar@{=}[d] \\
 & f_*f^*f_*\tau\ar^-{Ex_\tau^{-1}}[r]\ar@{=}[d]\ar@{}|{(3)}[rrrrd]
 & f_*\tau f^*f_*\ar|-{ad'(f^*,f_*)}[rrr]
 & 
 & 
 & f_*\tau\ar@{=}[d] \\
f_*\tau\ar^-{ad(f^*,f_*)}[r]
 & f_*f^*f_*\tau\ar|-{ad'(f^*,f_*)}[rrrr]
 & 
 & 
 & 
 & f_*\tau.
}
$$
The commutativity of (1) and (2) is obvious
 and the commutativity of (3) follows from Formula \eqref{eq:twist&f_*}
 defining $Ex_\tau$. 
 Then the result follows from the usual formula between the unit
 and counit of an adjunction.
 The relation $\phi_2 \circ \piso^\mathfrak t_f=1$ is proved
  using the same kind of computations.

\bigskip

It remains to prove that (i) and (i') are equivalent to (ii).
We already know from Proposition \ref{prop:purity&duality} that (i) implies 
 the premotive $M_S(X)$ is strongly dualizable.
 Saying that $\mu_f^\mathfrak t$ is a perfect pairing amounts to prove that the natural
 transformation obtained by adjunction
$$
d^\mathfrak t_f:(M_S(X) \otimes -) \rightarrow \uHom(M_S(X),-(d)[2d])
$$
is an isomorphism.
On the other hand, as we have already seen previously,
the smooth projection formula implies an identification of functors: 
\begin{equation} \label{eq:proof_purity&duality_or}
\begin{split}
f_\sharp f^* &\simeq (M_S(X) \otimes -), \\
f_*f^* &\simeq \uHom(M_S(X),-).
\end{split}
\end{equation}
Thus, to finish the proof, it will be enough to show that the map
$$
f_\sharp f^* \xrightarrow{\piso^\mathfrak t_ff^*} f_*\tau f^*=f_*f^*\tau.
$$
is equal to $d^\mathfrak t_f$ through the identifications \eqref{eq:proof_purity&duality_or}.

Let us consider the following cartesian square
$$
\xymatrix@=22pt{
X \times_S X\ar^-{f''}[r]\ar_{f'}[d]\ar@{}|\Delta[rd]
 & X\ar^f[d] \\
X\ar_f[r] & S
}
$$
and put $g=f \circ f''$.
According to the definition of $\mu_f^\mathfrak t$,
 and notably Formula \eqref{eq:Gysin} for the Gysin map $\delta^*$,
 the natural transformation of functors $(\mu_f^\mathfrak t \otimes -)$ can be described
 as the following composition:
\begin{equation*}
\begin{split}
f_\sharp f^* f_\sharp f^*
 \xrightarrow{Ex(\Delta_\sharp^*)}
  f_\sharp f'_\sharp f^{\prime \prime *} f^*
  & =g_\sharp g^*
 \xrightarrow{ad(\delta^*,\delta_*)}
  g_\sharp \delta_* \delta^* g^*\\
& =f_\sharp f'_\sharp \delta_* f^*
 \xrightarrow{\piso^\mathfrak t_{XX,X}} f_\sharp \tau f^*=f_\sharp f^*\tau
 \xrightarrow{ad'(f_\sharp,f^*)} \tau.
\end{split}
\end{equation*}
Note in particular that the base change map $Ex(\Delta_\sharp^*)$ 
 corresponds to the first identification in Formula \eqref{eq:dual_pairing_or}.
Thus we have to prove the preceding composite map is equal to
 the following one, obtained by adjunction from $\piso^\mathfrak t_f$:
\begin{equation*}
\begin{split}
f_\sharp f^* f_\sharp f^*
& =f_\sharp f^*f_\sharp f''_* \delta_*f^*
 \xrightarrow{Ex(\Delta_{\sharp*})}
 f_\sharp f^*f_* f'_\sharp \delta_*f^* \\
 & \xrightarrow{\piso^\mathfrak t_{XX,X}}
  f_\sharp f^*f_*\tau f^*=f_\sharp f^*f_*f^*\tau 
  \xrightarrow{ad'(f^*,f_*)} f_\sharp f^*\tau 
  \xrightarrow{ad'(f_\sharp,f^*)} \tau 
\end{split}
\end{equation*}
This amounts to prove,
 after some easy cancellation,
 the commutativity of the following diagram:
$$
\xymatrix@C=42pt{
f^*f_\sharp\ar_{Ex(\Delta^*_\sharp)}[d]\ar@{=}[r]
 & f^*f_\sharp f''_*\delta_*\ar^{Ex(\Delta_{\sharp *})}[r]
 & f^*f_*f''_\sharp\delta_*\ar^{ad'(f^*,f_*)}[d] \\
f'_\sharp f^{\prime\prime*}\ar^-{ad(\delta^*,\delta_*)}[r]
 & f'_\sharp \delta_* \delta^* f^{\prime\prime*}\ar@{=}[r]
 & f'_\sharp \delta_*.
}
$$
According to the definition of the exchange transformation
 $Ex(\Delta_{\sharp*})$ (cf Paragraph \ref{ex:H_complete_P-fibred}),
 we can divide this diagram into the following pieces:
$$
\xymatrix@C=32pt@R=38pt{
 f^*f_\sharp f''_*\delta_*\ar^-{ad(f^*,f_*)}[r]\ar_{Ex(\Delta^*_\sharp)}[d]
 & f^*f_*f^*f_\sharp f''_*\delta_*\ar^-{Ex(\Delta^*_\sharp)}[r]
 & f^*f_*f'_\sharp f^{\prime\prime*}f''_*\delta_*\ar^-{ad'(f^{\prime\prime*},f''_*)}[r]
    \ar|{ad'(f^*,f_*)}[d]
 & f^*f_*f''_\sharp\delta_*\ar|{ad'(f^*,f_*)}[d] \\
f'_\sharp f^{\prime\prime*}f''_*\delta_*\ar@{=}[d]\ar@{=}[rr]\ar^/-30pt/{ad(f^*,f_*)}[rru]
 && f'_\sharp f^{\prime\prime*}f''_*\delta_*\ar^{ad'(f^{\prime\prime*},f''_*)}[r]
 & f'_\sharp \delta_*\ar@{=}[d] \\
f'_\sharp f^{\prime\prime*}\ar|-{ad(\delta^*,\delta_*)}[rrr]\ar@{}|{(*)}[rrru]
 &
 &
 & f'_\sharp \delta_*.
}
$$
Every part of this diagram is obviously commutative except for part $(*)$.
As $f''\delta=1$, the axioms of a 2-functors (for $f^*$ and $f_*$ say)
 implies that the unit map
$$
\alpha:
f'_\sharp f^{\prime\prime*}
 \rightarrow  f'_\sharp f^{\prime\prime*}(f''\delta)_*(f''\delta)^*
$$
is the canonical identification that we get using $1_*=1$ and $1^*=1$.
We can consider the following diagram:
$$
\xymatrix@C=35pt@R=36pt{
f'_\sharp f^{\prime\prime*}\ar@{=}^-{\alpha}[rr]\ar@{=}[d]
 & & f'_\sharp f^{\prime\prime*}(f''\delta)_*(f''\delta)^*\ar@{=}[r]\ar@{=}[d]
 & f'_\sharp f^{\prime\prime*}f''_*\delta_*\ar_{ad'(f^{\prime\prime*},f''_*)}[dd] \\
f'_\sharp f^{\prime\prime*}\ar^-{ad(f^{\prime\prime*},f''_*)}[r]\ar@{=}[d]
 & f'_\sharp f^{\prime\prime*} f''_\sharp f^{\prime\prime*}
    \ar^-{ad(\delta^*,\delta_*)}[r]\ar_{ad'(f^{\prime\prime*},f''_*)}[d]
 & f'_\sharp f^{\prime\prime*}(f''\delta)_*(f''\delta)^*\ar_{ad'(f^{\prime\prime*},f''_*)}[d]
 & \\
f'_\sharp f^{\prime\prime*}\ar@{=}[r]
 & f'_\sharp f^{\prime\prime*}\ar^{ad(\delta^*,\delta_*)}[r]
 & f'_\sharp \delta_* \delta^* f^{\prime\prime*}\ar@{=}[r]
 & f'_\sharp \delta_*
}
$$
for which each part is obviously commutative. This concludes.
\end{proof}

As a corollary, together with the results of \cite{Deg8},
 we get the following theorem:
\begin{cor} \label{cor:orientationsandpurity}
Let us assume the following conditions on
 the triangulated premotivic category $\T$:
\begin{itemize}
\item[(a)] $\T$ satisfies properties \sepx \nis, \wloc, \htp and \stab.
\item[(b)] $\T$ admits an orientation $\mathfrak t$.
\item[(c)] There exists a morphism of triangulated premotivic categories:
$$
\varphi^*:\SH \rightleftarrows \T:\varphi_*\ .
$$
\end{itemize}
Then any smooth projective morphism is $\T$-pure.
 In particular, $\T$ is weakly pure.
\end{cor}
\begin{proof}
According to Example \ref{ex:or&GysinII},
 one can apply the results of \cite{Deg8} to the triangulated
 category $\T(X)$.
 Then it follows from \cite[5.23]{Deg8} that condition (ii)
 of the above theorem is satisfied.
\end{proof}

\begin{rem}
This theorem is to be compared with the result of Ayoub recalled
 in Theorem \ref{thm:Ayoub_wpur}.
 On the one hand, if $\T$ satisfies the localization property,
 we get another proof of this result under the additional assumption
 that $\T$ is oriented. On the other hand, the above theorem
 does not require the assumption that $\T$ satisfies \loc ; this
 is important as we can only prove \wloc for the category
 $\DM$ introduced in Definition \ref{df:Nis_DMe&DM}.
\end{rem}

\subsubsection{Motivic categories}
 
This section summarizes the main constructions of this part
 and draws a conclusive theorem.
\begin{df} \label{df:motivic_cat}
%A \emph{pregeometric category over $\sch$} is
%a triangulated complete $\sm$-fibred category over $\sch$ which
%satisfies the homotopy, stability and localization properties.
% 
%A \emph{geometric category over $\sch$} is
%a pregeometric category $\T$ over $\sch$ 
%such that, for any proper morphism $f$ in $\sch$, the functor
%$f_*$ has a right adjoint.
A \emph{motivic triangulated category over $\sch$} is
 a premotivic triangulated category over $\sch$ which 
 satisfies the homotopy, stability, localization and adjoint property.
\end{df}

\begin{rem}
Without the adjoint property,
 this definition corresponds to what Ayoub called a 
 \emph{monoidal stable homotopy $2$-functor}
\index{word}{monoidal stable homotopy $2$-functor}
  (cf \cite[def. 2.3.1]{ayoub}).
 We think our shorter terminology fits well in the spirit
 of the current theory of mixed motives.
\end{rem}

\begin{rem} \label{rem:compatly&motivic}
Assume $\T$ is a premotivic triangulated category such that:
\begin{enumerate}
\item $\T$ is well generated.
\item $\T$ satisfies the homotopy and stability properties.
\item $\T$ satisfies the localization property.
\end{enumerate}
Then $\T$ is a motivic triangulated category in the above sense.
Indeed, property \adj is proved under the above
 assumptions in point (4) of Theorem \ref{thm:localizationandwpur}.
 Note also that if $\T$ is compactly $\tau$-generated,
 we simply obtain property \adj from
  Lemma \ref{lm:exist_adjoint_f^!}.\footnote{In our examples,
  (1) will always be satisfied, (2) will be obtained by construction
  and (3) will be the hard point.}
\end{rem}

\begin{ex} \label{ex:SH_motivic}
According to the previous remark,
 the premotivic category $\SH$ of Example \ref{ex:H&SH_premotivic}
 is a motivic category.
 In fact,
 property (1) is proved in \cite[4.5.67]{ayoub},
 property (2) follows by definition
 and property (3) is proved in \cite[4.5.44]{ayoub}.
\end{ex}

\begin{num}
In the next theorem,
 we summarize what is now called the \emph{Grothendieck six functors formalism}.
\index{word}{formalism, Grothendieck six functors}
 In fact, this is a consequence of the axioms in the above definition,
  as a result of the work done in previous sections.
 More precisely:
\begin{itemize}
\item We apply Theorem \ref{thm:localizationandwpur} using
 the theorem of Ayoub
 recalled in \ref{thm:Ayoub_wpur}, and use the generalized theorem 
 of Morel and Voevodsky, Theorem \ref{thm:purity_closed_pairs},
 to get the form \eqref{eq:relative_purity_after_def} of the purity
 isomorphism.
\item In the case where $\T$ is oriented,
 we use the form \eqref{thm:loc_wpur_or} of the purity isomorphism.
 Recall that,
  when $\T$ satisfies assumption (c) of Corollary \ref{cor:orientationsandpurity},
  then we have given a different proof of the Theorem of Ayoub
  and the theorem below follows
  from \ref{thm:localizationandwpur} and \ref{cor:orientationsandpurity}. 
\end{itemize}
\end{num}
\begin{thm}\label{thm:cor3_Ayoub}
Let $\T$ be a motivic triangulated category.

Then, for any separated morphism of finite type $f:Y \rightarrow X$ in $\sch$,
 there exists a pair of adjoint functors, the \emph{exceptional functors},
\index{word}{functor!exceptional}
$$
f_!:\T(Y) \rightleftarrows \T(X):f^!
$$
such that:
\begin{enumerate}
\item There exists a structure of a covariant (resp. contravariant) 
 $2$-functor on $f \mapsto f_!$ (resp. $f \mapsto f^!$).
\item There exists a natural transformation $\alpha_f:f_! \rightarrow f_*$
 which is an isomorphism when $f$ is proper.
 Moreover, $\alpha$ is a morphism of $2$-functors.
\item For any smooth separated morphism of finite type $f:X \rightarrow S$ in $\sch$
 with tangent bundle\index{word}{bundle!tangent} $T_f$,
 there are canonical natural isomorphisms
\begin{align*}
\piso_f:f_\sharp & \longrightarrow f_!\big(MTh_X(T_f) \otimes_X .\big) \\
\piso'_f:f^* & \longrightarrow MTh_X(-T_f) \otimes_X f^!
\end{align*}
which are dual to each other
 -- the Thom premotive $MTh_X(T_f)$ is $\otimes$-invertible
 with inverse $MTh_X(-T_f)$.

If $\T$ admits an orientation $\mathfrak t$
 and $f$ has dimension $d$ then there are canonical natural isomorphisms
\begin{align*}
\piso^\mathfrak t_f:f_\sharp & \longrightarrow f_!(d)[2d] \\
\piso^{\prime\mathfrak t}_f:f^* & \longrightarrow f^!(-d)[-2d]
\end{align*}
which are dual to each other.
\item For any Cartesian square:
$$
\xymatrix@=16pt{
Y'\ar^{f'}[r]\ar_{g'}[d]\ar@{}|\Delta[rd] & X'\ar^g[d] \\
Y\ar_f[r] & X,
}
$$
such that $f$ is separated of finite type,
there exist natural isomorphisms
\begin{align*}
Ex(\Delta^*_!)&:
g^*f_! \xrightarrow\sim f'_!{g'}^*\, , \\
Ex(\Delta_*^!)&:
g'_*{f'}^! \xrightarrow\sim  f^!g_*\, .
\end{align*}
\item For any separated morphism of finite type $f:Y \rightarrow X$ in $\sch$,
 there exist natural isomorphisms
\begin{align*}
Ex(f_!^*,\otimes):
(f_!K) \otimes_X L &\xrightarrow{\ \sim\ } f_!(K \otimes_Y f^*L)\, ,\ \\
  \uHom_X(f_!(L),K) & \xrightarrow{\ \sim\ } f_* \uHom_Y(L,f^!(K))\, ,\ \\
  f^! \uHom_X(L,M)& \xrightarrow{\ \sim\ } \uHom_Y(f^*(L),f^!(M))\, .
\end{align*}
\end{enumerate}
\end{thm}

\begin{rem}\label{rem:purity&composition}
It is important to precise that in the case where the morphisms in $\sch$
 are assumed to be quasi-projective, this theorem is proved by Ayoub
 in \cite{ayoub} if we except the case where $\T$ is oriented in point (3).\footnote{
  This theorem was first announced by Voevodsky but only notes covering 
  the basic setting were to be found by the time Ayoub wrote the proof.}

Regarding this theorem,
 our contribution is to extend the result of Ayoub to the non quasi-projective
 case and to consider the oriented case --- which is crucial in the theory of motives.
 Recall also we have given another proof of this result in the case where the motivic category $\T$
 satisfies in addition the assumptions of Corollary \ref{cor:orientationsandpurity}
 --- which will always be the case for the different categories of motives
 introduced here.
\end{rem}

\begin{rem}
The purity isomorphism is compatible with composition.
Given smooth separated morphisms of finite type
$$
Y \xrightarrow g X \xrightarrow f S
$$
we obtain (\emph{cf.} \cite[17.2.3]{EGA4})
 an exact sequence of vector bundles over $Y$
\begin{equation} \tag{$\sigma$}
0 \rightarrow g^{-1}T_f \rightarrow T_{fg} \rightarrow T_g \rightarrow 0.
\end{equation}
which according to Remark \ref{rem:Thom_sp&exact_sequences}
 induces an isomorphism:
\begin{align*}
\epsilon_\sigma:MTh_Y(T_{fg})
 \xrightarrow{MTh_Y(\sigma)} & MTh_Y(T_g) \otimes_Y MTh_Y(g^{-1}T_f)\\
& \xrightarrow{\sim} g^*MTh_X(T_f) \otimes_Y MTh_Y(T_g).
\end{align*}
One can check the following diagram is commutative:
$$
\xymatrix{
(fg)_\sharp(K)\ar_{\piso_{fg}}[ddd]\ar@{=}[r]
 & f_\sharp g_\sharp(K)\ar^{\piso_f \circ \piso_g}[d] \\
 & f_!\Big(MTh_X(T_f) \otimes_X g_!\big(MTh_Y(T_g) \otimes_Y K\big)\Big)
    \ar^{Ex(g_!^*,\otimes)^{-1}}[d] \\
 & f_!g_!\big(g^*MTh_Y(T_f) \otimes_Y MTh_Y(T_g) \otimes_Y K\big)
    \ar^{\epsilon_\sigma^{-1}}[d] \\
(fg)_!(MTh(T_{fg}) \otimes K)\ar@{=}[r]
 & f_!g_!(MTh(T_{fg}) \otimes K).
}
$$
This is not an easy check.\footnote{The main point is to check
 that the isomorphism of Theorem \ref{thm:purity_closed_pairs}
 is compatible with composition (of closed immersions).
 On that particular point, see \cite[Th. 4.32, Cor. 4.33]{Deg8}.}
 In fact, this is one of the key technical point in
  the proof of the main Theorem of Ayoub (\cite[1.4.2]{ayoub}).
 We refer the reader to \cite[1.5]{ayoub} for details.

Note also that given the commutativity of the above diagram,
 if $\T$ admits an orientation $\mathfrak t$,
 it readily follows from axiom (c) of Definition \ref{df:premotivic_orientation}
 that the following diagram is commutative:
$$
\xymatrix{
(fg)_\sharp(K)\ar_{\piso^\mathfrak t_{fg}}[d]\ar@{=}[r]
 & f_\sharp g_\sharp(K)\ar^{\piso^\mathfrak t_f \circ \piso^\mathfrak t_g}[d] \\
(fg)_!(K)(n+m)[2n+2m]\ar@{=}[r]
 & f_!g_!(K)(n+m)[2n+2m]
}
$$
where $n$ (resp. $m$) is the relative dimension of $f$
 (resp. $g$).
\end{rem}

Morphisms of triangulated motivic categories
 are compatible with Grothendieck 6 operations in the following sense:
\begin{prop}\label{prop:motivic_adj_6_operations}
Let $\T$ and $\T'$ be motivic triangulated categories
 and
$$
\varphi^*:\T \rightleftarrows \T':\varphi_*
$$
be an adjunction of premotivic categories.
\index{word}{adjunction!of premotivic categories|see{morphism of}}
\index{word}{morphism!of triangulated premotivic categories}

Then $\varphi^*$ (resp. $\varphi_*$)
 commutes with the operations $f^*$ (resp. $f_*$),
 for any morphism of schemes $f$, as well as with the operation
 $p_!$ (resp. $p^!$), for any separated morphism of finite type $p$.

Moreover, $\varphi^*$ is monoidal and for any premotive
 $M \in \T(S)$, $N \in \T'(S)$, the canonical map
$$
\uHom(M,\varphi_*(N))
 \rightarrow \varphi_*\uHom(\varphi^*(M),N)
$$
is an isomorphism.
\end{prop}
\begin{proof}
The only thing to prove is that $\varphi^*$
 commutes with $p_!$ as the other statements follows
 either from the definitions or by adjunction.
 This follows from Proposition \ref{prop:morph&loc1},
  the purity property in $\T$ and $\T'$
   (property (3) in the above theorem)
   and the fact $\varphi^*$ commutes with $p_\sharp$
   when $p$ is smooth by assumption.
\end{proof}

\begin{rem}
With additional assumptions on $\T$ and $\T'$, and over a field,
we will see that $\varphi^*$ commutes with all the
six operations (see Theorem \ref{thm:realcommute6op}).
\end{rem}

\section{Descent in $\Pmor$-fibred model categories} \label{sec:theorie_descente}

\begin{assumption} \label{num:assumption_hdescent_sch}
In this section, $\sch$ is an abstract category
 and $\Pmor$ an admissible class of morphisms in $\sch$.

In section \ref{sec:hdescent} however, we will consider
 as in \ref{num:assumption1_sch} a noetherian base scheme $\base$
 and we will assume that $\sch$ is an adequate category of $\base$-schemes
 satisfying the following condition on $\sch$:
\begin{itemize}
\item[(a)] Any scheme in $\sch$ is finite dimensional.
\end{itemize}
Moreover,
 in sections \ref{sec:proper_descent&Galois} and \ref{sec:proper_descent&separation},
 we will even assume:
\begin{itemize}
\item[(a$'$)] Any scheme in $\sch$ is quasi-excellent and finite dimensional.
\end{itemize}

We fix an admissible class $\Pmor$ of morphisms in $\sch$
 which contains the class of {\'e}tale morphisms in $\sch$
  and a stable combinatorial $\Pmor$-fibred model category $\M$ over
	$\sch$ (See Paragraph \ref{num:P-fibred_model}). \\
In section \ref{sec:proper_descent&separation},
 we will assume furthermore that:
\begin{itemize}
\item[(b)] The stable model $\Pmor$-fibred category $\M$ is $\QQ$-linear
 (see \ref{defQlinear}).
\end{itemize}

\end{assumption}

\subsection{Extension of $\Pmor$-fibred categories to diagrams}
\label{sec:extension_to_diagrams}

\subsubsection{The general case}

\begin{paragr} \label{paragr:extension_to_diagrams}
Assume given a $\Pmor$-fibered category $\M$ over $\site$.
Then $\M$ can be extended to $\site$-diagrams
\index{word}{diagram!Sdiagram@$\site$-diagram}
(i.e. functors from a small category to $\site$) as follows.
Let $I$ be a small category, and $\X$ a functor from $I$ to $\site$.
For an object $i$ of $I$, we will denote by $\X_i$
the fiber of $\X$ at $i$ (i.e. the evaluation of $\X$ at $i$), and, for
a map $u:i\To j$ in $I$, we will still denote by $u:\X_i\To\X_j$ the morphism
induced by $u$. We define the category $\M(\X,I)$ as follows.

An object of $\M(\X,I)$ is a couple $(M,a)$, where $M$ is the data
of an object $M_i$ in $\M(\X_i)$ for any object $i$ of $I$,
and $a$ is the data of a morphism $a_u:u^*(M_j)\To M_i$
for any morphism $u:i\To j$ in $I$, such that, for any object $i$
of $I$, the map $a_{1_i}$ is the identity of $M_i$ (we will always
assume that $1^*_i$ is the identity functor), and,
for any composable morphisms
$u:i\To j$ and $v:j\To k$ in $I$, the following diagram commutes.
$$\xymatrix{
u^*v^*(M_k)\ar[d]_{u^*(a_v)}\ar[r]^\simeq&(vu)^*(M_k)\ar[d]^{a_{vu}}\\
u^*(M_j)\ar[r]_{a_u}& M_i
}$$
A morphism $p:(M,a)\To (N,b)$ is a collection of morphisms
$$p_i:M_i\To N_i$$
in $\M(\X_i)$, for each object $i$ in $I$, such that, for any morphism
$u:i\To j$ in $I$, the following diagram commutes.
$$\xymatrix{
u^*(M_j)\ar[r]^{u^*(p_j)}\ar[d]_{a_u}&u^*(N_j)\ar[d]^{b_u}\\
M_i\ar[r]_{p_i}&N_i
}$$
In the case where $\M$ is a monoidal $\Pmor$-fibred category,
the category $\M(\X,I)$ is naturally endowed with a symmetric monoidal
structure. Given two objects $(M,a)$ and $(N,b)$ of $\M(\X,I)$,
their tensor product
$$(M,a)\otimes (N,b)=(M\otimes N,a\otimes b)$$
is defined as follows. For any object $i$ of $I$,
$$(M\otimes N)_i=M_i\otimes N_i\, ,$$
and for any map $u:i\To j$ in $I$, the map
$(a\otimes b)_u$ is the composition of the isomorphism
$u^*(M_j\otimes N_j)\simeq u^*(M_j)\otimes u^*(N_j)$ with the morphism
$$a_u\otimes b_u:u^*(M_j)\otimes u^*(N_j)\To M_i\otimes N_i\, .$$
Note finally that if $\M$ is a complete monoidal $\Pmor$-fibred category,
 then $\M(\X,I)$ admits an internal Hom.
\end{paragr}

\begin{paragr}\label{defevaluation} \textit{Evaluation functors}.
Assume now that for any $S$, $\M(S)$ admits small sums.

For each object $i$ of $I$, we have a functor
\begin{equation}\label{evobject}
\begin{aligned}
i^*:\M(\X,I)&\To\M(\X_i)\\
(M,a)&\longmapsto M_i
\end{aligned}
\end{equation}
called the \emph{evaluation functor}\index{word}{functor!evaluation}
 associated with $i$.
This functor $i^*$ has a left adjoint
\begin{equation}\label{leftadjevobject}
i_\sharp:\M(\X_i)\To\M(\X,I)
\end{equation}
defined as follows. If $M$ is an object of $\M(\X_i)$, then
$i_\sharp(M)$ is the data $(M',a')$ such that for any object $j$ of $I$,
\begin{equation}\label{leftadjevobject2}
\left(i_\sharp(M)\right)_j=M'_j=\coprod_{u\in\Hom_I(j,i)}u^*(M)\, ,
\end{equation}
and, for any morphism $v:k\To j$ in $I$, the map $a'_v$ is the
canonical map induced by the collection of maps
\begin{equation}\label{leftadjevobject3}
v^*u^*(M)\simeq (uv)^*(M)\To\coprod_{w\in\Hom_I(k,i)}w^*(M)
\end{equation}
for $u\in\Hom_I(j,i)$.

If we assume that $\M$ is a complete $\Pmor$-fibred category
 and that $\M(S)$ admits small products for any $S$,
then $i^*$ has a right adjoint
\begin{equation}\label{rightadjevobject}
i_*:\M(\X_i)\To\M(\X,I)
\end{equation}
given, for any object $M$ of $\M(\X_i)$ by the formula
\begin{equation}\label{rightadjevobject2}
\left(i_*(M)\right)_j=\prod_{u\in\Hom_I(i,j)}u_*(M),
\end{equation}
with transition map given by the dual formula of \ref{leftadjevobject3}.
\end{paragr}

\begin{paragr} \label{defiminvvarphi} \textit{Functoriality}.
Assume that $\M$ if a $\Pmor$-fibred category such that
 for any object $S$ of $\site$, $\M(S)$ has small colimits.

Remember that, if $\X$ and $\Y$ are $\site$-diagrams, indexed
respectively by small categories $I$ and $J$,
 a morphism of $\site$-diagrams
 \index{word}{morphism!of Sdiagrams@of $\site$-diagrams}
  $\varphi:(\X,I)\To(\Y,J)$
 is a couple $\varphi=(\alpha,f)$, where $f:I\To J$ is a functor,
and $\alpha:\X\To f^*(\Y)$ is a natural transformation (where $f^*(\Y)=\Y\circ f$).
In particular, for any object $i$ of $I$, we have a morphism
$$\alpha_i:\X_i\To\Y_{f(i)}$$
in $\site$. This turns $\site$-diagrams into a strict $2$-category:
the identity of $(\X,I)$
is the couple $(1_\X,1_I)$, and, if $\varphi=(\alpha,f):(\X,I)\To(\Y,J)$
and $\psi=(\beta,g):(\Y,J)\To(\Z,K)$ are two composable morphisms,
the composed morphism $\psi\circ\varphi:(\X,I)\To(\Z,K)$ is the couple $(gf,\gamma)$,
where for each object $i$ of $I$, the map
$$\gamma_i:\X_i\To \Z_{g(f(i))}$$ is the composition
$$\X_i\xrightarrow{ \ \alpha_i \ }\Y_{f(i)}\xrightarrow{\beta_{f(i)}}\Z_{g(f(i))}\, .$$
There is also a notion of natural transformation between
morphisms of $\site$-diagrams: if $\varphi=(\alpha,f)$ and $\varphi'=(\alpha',f')$
are two morphisms from $(\X,I)$ to $(\Y,J)$, a natural transformation $t$
from $\varphi$ to $\varphi'$ is a natural transformation $t:f\To f'$ such that
the following diagram of functors commutes.
$$\xymatrix{
&\X\ar[dl]_{\alpha}\ar[dr]^{\alpha'}&\\
\Y\circ f\ar[rr]_t&&\Y\circ f'
}$$
This makes the category of $\site$-diagrams a (strict) $2$-category.

To a morphism of diagrams $\varphi=(\alpha,f):(\X,I)\To(\Y,J)$, we associate a functor
$$\varphi^*:\M(\Y,J)\To\M(\X,I)$$
as follows. For an object $(M,a)$ of $\M(\Y)$, $\varphi^*(M,a)=(\varphi^*(M),\varphi^*(a))$
is the object of $\M(\X)$ defined by
$\varphi^*(M)_i=\alpha^*_i(M_{f(i)})$ for $i$ in $I$, and by the formula
$$\varphi^*(a)_u=\alpha^*_i(a_{f(u)}):
\alpha^*_i\, f(u)^*(M_{f(j)})=u^*\, \alpha^*_j(M_{f(j)})\To \alpha^*_i(M_{f(i)})$$
for $u:i\To j$ in $I$.

We will say that a morphism $\varphi:(\X,I)\To(\Y,J)$ is a \emph{$\Pmor$-morphism} if,
for any object $i$ in $I$, the morphism $\alpha_i:\X_i\To\Y_{f(i)}$ is
a $\Pmor$-morphism.
For such a morphism $\varphi$,
 the functor $\varphi^*$ has a left adjoint which we denote by
$$\varphi_\sharp:\M(\X,I)\To\M(\Y,J)\, .$$
For instance, given a $\site$-diagram $\X$ indexed by a small category $I$,
each object $i$ of $I$ defines a $\Pmor$-morphism of diagrams $i:\X_i\To(\X,I)$
(where $\X_i$ is indexed by the terminal category), so that the corresponding
the functor $i_\sharp$ corresponds precisely to \eqref{leftadjevobject}.

Assume that $\M$ is a complete $\Pmor$-fibred category such
 that $\M(S)$ has small limits for any object $S$ of $\site$.
Then the functor $\varphi^*$ has a right adjoint which we denote by
$$\varphi_*:\M(\X,I)\To\M(\Y,J)\, .$$
In the case where $\varphi$ is the morphism $i:\X_i\To(\X,I)$ defined
 by an object $i$ of $I$, $i_*$ corresponds precisely
 to \eqref{rightadjevobject}.
\end{paragr}

\begin{rem}
This construction can be applied in particular
 to any Grothendieck abelian (monoidal) $\Pmor$-fibred category
  (\textit{cf.} definition \ref{df:P-fibred_Grothendieck_abelian}).
The triangulated case cannot be treated in general without assuming
 a thorough structure -- this is the purpose of the next section.
\end{rem}

\subsubsection{The model category case} \label{sec:fibredmodcat}

\begin{paragr}
Let $\M$ be a $\Pmor$-fibred model category over $\site$
 (\textit{cf.} \ref{num:P-fibred_model}).
Given a $\site$-diagram $\X$ indexed by a small category $I$,
we will say that a morphism of $\M(\X,I)$ is a \emph{termwise weak equivalence}
 \index{word}{equivalence!termwise weak equivalence}
 \index{word}{weak equivalence|see{equivalence}}
(resp. a \emph{termwise fibration}\index{word}{fibration!termwise},
 resp. a \emph{termwise cofibration}\index{word}{cofibration!termwise}) if,
 for any object $i$ of $I$, its image by the functor $i^*$ is
 a weak equivalence (resp. a fibration, resp. a cofibration) in $\M(\X_i)$.
\end{paragr}

\begin{prop}\label{projdiamodcat}
If $\M$ is a cofibrantly generated $\Pmor$-fibred model category over $\site$,
then, for any $\site$-diagram $\X$ indexed by a small category $I$,
the category $\M(\X,I)$ is a cofibrantly generated
model category whose weak equivalences (resp. fibrations) are the
termwise weak equivalences (resp. the termwise fibrations).
This model category structure on $\M(\X,I)$ will be called the
\emph{projective model structure.}\index{word}{model structure!projective (diagrams)}

Moreover, any cofibration of $\M(\X,I)$ is a termwise cofibration, and
the family of functors
$$i^*:\ho(\M)(\X,I)\To\ho(\M)(\X_i)\ , \quad i\in\mathrm{Ob}(I)\, ,$$
is conservative.

If $\M$ is left proper (resp. right proper, resp. combinatorial, resp.
stable), then so is the projective model category structure on $\M(\X)$.
\end{prop}

\begin{proof}
Let $\X^{\delta}$ be the $\site$-diagram indexed by the
set of objects of $I$ (seen as a discrete category), whose
fiber at $i$ is $\X_i$. Let $\varphi:(\X^\delta,\mathit{Ob}\, I)\To(\X,I)$
be the inclusion (i.e. the map which is the identity on objects
and which is the identity on each fiber). As $\varphi$ is
clearly a $\Pmor$-morphism, we have an adjunction
$$\varphi_\sharp:\M(\X^\delta,\mathit{Ob}\, I)\simeq\prod_{i}\M(\X_i)\rightleftarrows\M(\X,I):\varphi^*\, .$$
The functor $\varphi_\sharp$ can be made explicit: it sends a family of
objects $(M_i)_i$ (with $M_i$ in $\M(\X_i)$) to the sum of the $i_\sharp(M_i)$'s
indexed by the set of objects of $I$.
Note also that this proposition is trivially verified whenever $\X^\delta=\X$.
Using the explicit formula for $i_\sharp$
given in \ref{defevaluation}, it is then straightforward to check that
the adjunction $(\varphi_\sharp,\varphi^*)$ satisfies the assumptions of \cite[Theorem 3.3]{crans},
which proves the existence of the projective model structure on $\M(\X,I)$.
Furthermore, the generating cofibrations (resp. trivial cofibrations of $\M(\X,I)$)
can be described as follows. For each object $i$ of $I$, let $A_i$ (resp. $B_i$) be
a generating set of cofibrations (resp. of trivial cofibrations in
$\M(\X_i)$. The class of termwise trivial fibrations (resp. of termwise fibrations)
of $\M(\X,I)$ is the class of maps which have the right lifting property
with respect to the set $A=\cup_{i\in I}\, i_\sharp(A_i)$ (resp. to the set
$B=\cup_{i\in I}\, i_\sharp(B_i)$). Hence, the set $A$ (resp. $B$) generates the
class of cofibrations (resp. of trivial cofibrations). In particular,
as any element of $A$ is a termwise cofibration (which follows immediately
from the explicit formula for $i_\sharp$
given in \ref{defevaluation}), and as termwise cofibrations
are stable by pushouts, transfinite compositions and retracts, any cofibration
is a termwise cofibration (by the small object argument). 

As any fibration (resp. cofibration) of $\M(\X,I)$ is a termwise fibration (resp.
a termwise cofibration), it is clear that, whenever the model categories $\M(\X_i)$
are right (resp. left) proper, the model category $\M(\X,I)$ has the same property.

The functor $\varphi^*$ preserves fibrations and cofibrations, while
it also preserves and detects weak equivalences (by definition). This implies
that the induced functor
$$\varphi^*:\ho(\M)(\X,I)\To\ho(\M)(\X^\delta,\mathit{Ob}\, I)\simeq\prod_i\ho(\M)(\X_i)$$
is conservative (using the facts that the set of maps from a cofibrant
object to a fibrant object in the homotopy category of a model category
is the set of homotopy classes of maps, and that a morphism of
a model category is a weak equivalence if and only if it induces an
isomorphism in the homotopy category). As $\varphi^*$ commutes to
limits and colimits, this implies that it commutes to homotopy limits
and to homotopy colimits (up to weak equivalences). Using the conservativity
property, this implies that a commutative square of $\M(\X,I)$ is a
homotopy pushout (resp. a homotopy pullback) if and only if it is
so in $\M(\X^\delta,\mathit{Ob}\, I)$. Remember that stable model categories
are characterized as those in which a commutative square is a
homotopy pullback square if and only if it is a homotopy pushout square.
As a consequence, if all the model categories $\M(\X_i)$ are stable,
as $\M(\X^\delta,\mathit{Ob}\, I)$ is then obviously stable as well, the model category
$\M(\X,I)$ has the same property.

It remains to prove that, if $\M(X,I)$ is a combinatorial model category for any
object $X$ of $\site$, then $\M(\X,I)$ is combinatorial as well.
For each object $i$ in $I$, let $G_i$ be a set of accessible generators
of $\M(\X_i)$. Note that, for any object $i$ of $I$, the functor $i_\sharp$
has a left adjoint $i^*$ which commutes to colimits (having itself a right
adjoint $i_*$).
It is then easy to check that the set of objects of shape $i_\sharp(M)$,
for $M$ in $G_i$ and $i$ in $I$, is a small set of accessible generators of $\M(\X,I)$.
This implies that $\M(\X,I)$ is accessible and ends the proof.
\end{proof}

\begin{prop}\label{injdiamodcat}
Let $\M$ be a combinatorial $\Pmor$-fibred model category over $\site$.
Then, for any $\site$-diagram $\X$ indexed by a small category $I$,
the category $\M(\X,I)$ is a combinatorial
model category whose weak equivalences (resp. cofibrations) are the
termwise weak equivalences (resp. the termwise cofibrations).
This model category structure on $\M(\X,I)$ will be called the
\emph{injective model structure}\footnote{Quite unfortunately,
this corresponds to the `semi-projective' model structure
introduced in \cite[Def. 4.5.8]{ayoub}.}.
\index{word}{model structure!injective (diagrams)}
Moreover, any fibration of the injective model structure on $\M(\X,I)$ is a termwise fibration.

If $\M$ is left proper (resp. right proper, resp. stable), then so is the injective
model category structure on $\M(\X,I)$.
\end{prop}

\begin{proof} See \cite[Theorem 2.28]{Bar} for the existence of such
a model structure (if, for any object $X$ in $\site$,
all the cofibrations of $\M(X)$ are monomorphisms,
this can also be done following mutatis mutandis the proof of
\cite[Proposition 4.5.9]{ayoub}).
Any trivial cofibration of the projective model structure
being a termwise trivial cofibration, any fibration of the injective model structure is a fibration
of the projective model structure, hence a termwise fibration.

The assertions about properness follow from their analogs for the projective
model structure and from \cite[Corollary 1.5.21]{Cis3} (or can be proved directly;
see \cite[Proposition 2.31]{Bar}).
Similarly, the assertion on stability follows from their analogs for the
projective model structure.
\end{proof}

\begin{paragr}
From now on, we assume that a combinatorial $\Pmor$-fibred model category $\M$ over $\site$ is given.
%such that,
%for any $S$-scheme $X$, all the cofibrations of $\M(X)$ are monomorphisms.
Then, for any $\site$-diagram $(\X,I)$, we have two model category structures on $\M(\X,I)$, and
the identity defines a left Quillen equivalence from the projective model structure to the injective model
structure. This fact will be used for the understanding
 of the functorialities coming from morphisms of diagrams of $S$-schemes.
\end{paragr}

\begin{paragr}
The category of $\site$-diagrams admits small sums.
If $\{(\Y_j,I_j)\}_{j\in J}$ is a small family of $\site$-diagrams, then
their sum is the $\site$-diagram $(\X,I)$, where
$$I=\coprod_{j\in J}I_j\, , $$
and $\X$ is the functor from $I$ to $\site$ defined by
$$\X_i=\Y_{j}\quad\text{whenever $i\in I_j$.}$$
\end{paragr}

\begin{prop}\label{sumCdiag}
For any small family of $\site$-diagrams $\{(\Y_j,I_j)\}_{j\in J}$,
the canonical functor
$$\ho(\M)\Big(\coprod_{j\in J}\Y_j \Big)\To \prod_{j\in J}\ho(\M)(\Y_j)$$
is an equivalence of categories.
\end{prop}

\begin{proof}
The functor
$$ \M\Big(\coprod_{j\in J}\Y_j \Big)\To \prod_{j\in J} \M(\Y_j)$$
is an equivalence of categories. It thus remains an equivalence
after localization. To conclude, it is sufficient to see that the
homotopy category of a product of model categories is
the product of their homotopy categories, which follows
rather easily from the explicit description of the homotopy category
of a model category; see e.g. \cite[Theorem 1.2.10]{Hovey}.
\end{proof}

\begin{prop}\label{basicfunctdiag}
Let $\varphi=(\alpha,f):(\X,I)\To(\Y,J)$ be a morphism of $\site$-diagrams.
\begin{itemize}
\item[(i)] The adjunction $\varphi^*:\M(\Y,J)\rightleftarrows\M(\X,I):\varphi_*$
is a Quillen adjunction\index{word}{adjunction!Quillen adjunction}
 with respect to the injective model structures. In particular, it induces a derived adjunction
$$\derL\varphi^*:\ho(\M)(\Y,J)\rightleftarrows\ho(\M)(\X,I):\derR\varphi_*\, .$$
\item[(ii)] If $\varphi$ is a $\Pmor$-morphism, then the adjunction
$\varphi_\sharp:\M(\X,I)\rightleftarrows\M(\Y,J):\varphi^*$ is a Quillen adjunction with respect to the
projective model structures, and the functor $\varphi^*$ preserves weak equivalences. In particular,
we get a derived adjunction
$$\derL\varphi_\sharp : \ho(\M)(\X,I)\rightleftarrows\ho(\M)(\Y,J):\derL\varphi^*=\varphi^*=\derR\varphi^*\, .$$
\end{itemize}
\end{prop}

\begin{proof}
The functor $\varphi^*$ obviously preserves termwise cofibrations and termwise trivial cofibrations
(we reduce to the case of a morphism of $\site$ using the explicit description of $\varphi^*$
given in \ref{defiminvvarphi}), which proves the first assertion. Similarly, the second assertion follows from
the fact that, under the assumption that $\varphi$ is a $\Pmor$-morphism, the functor $\varphi^*$ preserves
termwise weak equivalences (see Remark \ref{remiminvPmorphism}), as well as termwise fibrations.
\end{proof}

\begin{paragr}\label{computeadjointsdiagrams}
The computation of the (derived) functors $\derR\varphi_*$
(and $\derL\varphi_\sharp$ whenever it makes sense)
given by Proposition \ref{basicfunctdiag} has to do with
 homotopy limits\index{word}{homotopy!limit}
  (and homotopy colimits\index{word}{homotopy!colimit}).
 It is easier to first understand this
  in the non derived version as follows.

Consider first the trivial case of a constant $\site$-diagram: let $X$
be an object of $\site$, and $I$ a small category. Then, seeing $X$
as the constant functor $I\To\site$ with value $X$, we have a
projection map $p^{}_I:(X,I)\To X$.
From the very definition, the category $\M(X,I)$ is simply the
category of   on $I$ with values in $\M(X)$, so that
the inverse image functor
\begin{equation}\label{limits0}
p^*_I:\M(X)\To\M(X,I)=\M(X)^{\op{I}}
\end{equation}
is the `constant diagram functor', while its right adjoint
\begin{equation}\label{limits1}
\varprojlim_{\op{I}}=p^{}_{I,*}:\M(X,I)\To\M(X)
\end{equation}
is the limit functor, and its left adjoint,
\begin{equation}\label{limits2}
\varinjlim_{\op{I}}=p^{}_{I,\sharp}:\M(X,I)\To\M(X)
\end{equation}
is the colimit functor.

Let $S$ be an object of $\site$. A \emph{$\site$-diagram over $S$} is the data of
a $\site$-diagram $(\X,I)$, together with a morphism of $\site$-diagrams
$p :(\X,I)\To S$ (i.e. its a $\site/S$-diagram). Such a map $p$
factors as
\begin{equation}\label{limits3}
(\X,I)\overset{\pi}{\To}(S,I)\overset{p^{}_I}{\To}S\, ,
\end{equation}
where $\pi=(p,1_I)$. Then one checks easily that, for any object $M$
of $\M(\X,I)$, and for any object $i$ of $I$, one has
\begin{equation}\label{limits4}
\pi_*(M)_i\simeq p^{}_{i,*}(M_i)\, ,
\end{equation}
where $p^{}_i:\X_i\To S$ is the structural map, from which we deduce the formula
\begin{equation}\label{limits5}
p_*(M)\simeq \varprojlim_{i\in \op{I} }\pi_*(M)_i\simeq \varprojlim_{i\in \op{I} } p^{}_{i,*}(M_i)\, ,
\end{equation}
Remark that, if $I$ is a small category with a terminal object $\omega$, then any
$\site$-diagram $\X$ indexed by $I$ is a $\site$-diagram over $\X_\omega$, and we deduce from the
computations above that, if $p:(\X,I)\To\X_\omega$ denotes the canonical map, then, for
any object $M$ of $\M(\X,I)$,
\begin{equation}\label{limits6}
p_*(M) \simeq M_\omega\, .
\end{equation}

Consider now a morphism of $\site$-diagrams $\varphi=(\alpha,f):(\X,I)\To(\Y,J)$.
For each object $j$, we can form the following pullback square of categories.
\begin{equation}\label{limits7}\begin{split}\xymatrix{
I/j\ar[r]^{u_j}\ar[d]_{f/j}&I\ar[d]^f\\
J/j\ar[r]_{v_j}&J
}\end{split}\end{equation}
in which $J/j$ is the category of objects of $J$ over $j$
(which has a terminal object, namely $(j,1_j)$, and $v_j$
is the canonical projection; the category $I/j$ is thus the category
of pairs $(i,a)$, where $i$ is an object of $I$, and $a:f(i)\to j$
a morphism in $J$. From this, we can form the following pullback
of $\site$-diagrams
\begin{equation}\label{limits8}\begin{split}\xymatrix{
(\X/j,I/j)\ar[r]^{\mu_j}\ar[d]_{\varphi/j}&(\X,I)\ar[d]^\varphi \\
(\Y/j,J/j)\ar[r]_{\nu_j}&(\Y,J)
}\end{split}\end{equation}
in which $\X/j=\X\circ u_j$, $\Y/j=\Y\circ v_j$, and the maps $\mu_j$
and $\nu_j$ are the one induced by $u_j$ and $v_j$ respectively.
For an object $M$ of $\M(\X,I)$ (resp. an object $N$ of $\M(\Y,J)$), we define
$M/j$ (resp. $N/j$) as the object of $\M(\X/j,I/j)$ (resp. of $\M(\Y/j,J/j)$)
obtained as $M/j=\mu^*_j(M)$ (resp. $N/j=\nu^*_j(N)$).
With these conventions, for any object $M$ of $\M(\X,I)$ and any object $j$
of the indexing category $J$, one gets the formula
\begin{equation}\label{limits9}
\varphi_*(M)_j\simeq (\varphi/j)_*(M/j)_{(j,1_j)}
\simeq \varprojlim_{(i,a)\in\op{I/j}}\alpha_{i,*}(M_i) \, .
\end{equation}
This implies that the natural map
\begin{equation}\label{limits10}
\varphi_*(M)/j=\nu^*_j\, \varphi_*(M)\To (\varphi/j)_*\, \mu^*_j(M)=(\varphi/j)_*(M/j)
\end{equation}
is an isomorphism: to prove this, it is sufficient to obtain
an isomorphism from \eqref{limits10} after evaluating by any object $(j',a:j'\To j)$
of $J/j$, which follows readily from \eqref{limits9} and from the obvious fact that
$(I/j)/(j',a)$ is canonically isomorphic to $I/j'$.

In order to deduce from the computations above their derived versions,
we need two lemmata.
\end{paragr}

\begin{lm}\label{iminvevaluationrq}
Let $\X$ be a $\site$-diagram indexed by a small category $I$, and $i$
an object of $I$. Then the evaluation functor\index{word}{functor!evaluation}
$$i^*:\M(\X,I)\To\M(\X_i)$$
is a right Quillen functor\index{word}{functor!Quillen}
 with respect to the injective model structure,
and it preserves weak equivalences.
\end{lm}

\begin{proof}
Proving that the functor $i^*$ is a right Quillen functor
is equivalent to proving that its left adjoint \eqref{leftadjevobject}
is a left Quillen functor with respect to the injective model structure,
which follows immediately from its computation \eqref{leftadjevobject2},
as, in any model category, cofibrations and trivial cofibrations
are stable by small sums. The last assertion is obvious from
the very definition of the weak equivalences in $\M(\X,I)$.
\end{proof}

\begin{lm}\label{iminvdiscfibrq}
For any pullback square of $\site$-diagrams of shape \eqref{limits8}, the functors
\begin{align*}
\mu^*_j:\M(\X,I)&\To\M(\X/j,I/j) \, , \quad M\mapsto M/j\\
\nu^*_j:\M(\Y,I)\, &\To\M(\Y/j,J/j)\ , \quad N\mapsto N/j
\end{align*}
are right Quillen functors with respect to the injective model
structure, and they preserve weak equivalences.
\end{lm}

\begin{proof}
It is sufficient to prove this for the functor $\mu^*_j$
(as $\nu^*_j$ is simply the special case where $I=J$ and $f$ is the identity).
The fact that $\mu^*_j$ preserves weak equivalences is obvious, so that it remains
to prove that it is a right Quillen functor. We thus have to prove that
left adjoint of $\mu^*_j$,
$$\mu_{j,\sharp}:\M(\X/j,I/j)\To\M(\X,I)\, ,$$
is a left Quillen functor. In other words, we have to prove that, for any
object $i$ of $I$, the functor
$$i^*\mu_{j,\sharp}:\M(\X,I)\To\M(\X)$$
is a left Quillen functor. For any object $M$ of $\M(\X,I)$,
we have a natural isomorphism
$$i^*\mu_{j,\sharp}(M)\simeq\coprod_{a\in\Hom_J(f(i),j)}(i,a)_\sharp(M_i)\, .$$
But we know that the functors $(i,a)_\sharp$ are left Quillen functors, so that
the stability of cofibrations and trivial cofibrations by small sums
and this description of the functor $i^*\mu_{j,\sharp}$ achieves the proof.
\end{proof}

\begin{prop}\label{computeadjointsdiagrams2}
Let $S$ be an object of $\site$, and $p:(\X,I)\To S$ a $\site$-diagram over $S$,
and consider the canonical factorization \eqref{limits3}.
For any object $M$ of $\ho(\M)(\X,I)$, there are canonical
isomorphisms and $\ho(\M)(S)$:
$$\derR\pi_*(M)_i\simeq \derR p^{}_{i,*}(M_i)\quad\text{and}\quad
\derR p_*(M) \simeq \underset{\quad i\in \op{I} }{\derR\varprojlim} \derR p^{}_{i,*}(M_i)\, .$$
In particular, if furthermore the category $I$ has a terminal object $\omega$, then
$$\derR p_*(M)\simeq \derR p_{\omega,*}(M_\omega)\, .$$
\end{prop}

\begin{proof}
This follows immediately from Formulas \eqref{limits4}, \eqref{limits5}
 and from the fact that deriving (right) Quillen functors is compatible
 with composition.
\end{proof}

\begin{prop}\label{computeadjointsdiagrams3}
We consider the pullback square of $\site$-diagrams \eqref{limits8}
(as well as the notations thereof).
For any object $M$ of $\ho(\M)(\X,I)$, and any object $j$ of $J$,
we have natural isomorphisms
$$\derR \varphi_*(M)_j \simeq
\underset{\phantom{R} (i,a)\in\op{I/j}}{\derR \varprojlim} \derR \alpha_{i,*}(M_i) 
\quad\text{and}\quad
\derR \varphi_*(M)/j \simeq \derR (\varphi/j)_*(M/j)$$
in $\ho(\M)(\Y_j)$ and in $\ho(\M)(\Y/j,J/j)$ respectively.
\end{prop}

\begin{proof}
Using again the fact that deriving right Quillen functors
is compatible with composition,
by virtue of Lemma \ref{iminvevaluationrq}
and Lemma \ref{iminvdiscfibrq},
this is a direct translation of \eqref{limits9} and \eqref{limits10}.
\end{proof}

\begin{prop}\label{basicPbasechange}
Let $u:T\To S$ be a $\Pmor$-morphism of $\site$, and $p:(\X,I)\To S$
a $\site$-diagram over $S$. Consider the pullback square of $\site$-diagrams
$$\xymatrix{
(\Y,I)\ar[r]^\varphi \ar[d]_q &(\X,I)\ar[d]^p \\
T\ar[r]_u & S
}$$
(i.e. $\Y_i=T \times_{S} \X_i$ for any object $i$ of $I$).
Then, for any object $M$ of $\ho(\M)(\X,I)$, the canonical map
$$\derL u^*\, \derR p_*(M)\To \derR q_*\, \derL v^*(M)$$
is an isomorphism in $\ho(\M)(T)$.
\end{prop}

\begin{proof}
By Remark \ref{remiminvPmorphism},
the functor $\nu^*$ is both a left and a right Quillen functor which preserves
weak equivalences, so that the functor $\derL \nu^*=\nu^*=\derR \nu^*$
preserves homotopy limits. Hence, by Proposition \ref{computeadjointsdiagrams2}, one
reduces to the case where $I$ is the terminal category, i.e. to
the transposition of the isomorphism given by the $\Pmor$-base change formula
\bc for the homotopy $\Pmor$-fibred category $\ho(\M)$ (see \ref{Pcotransversality}).
\end{proof}

\begin{paragr}
A morphism of $\site$-diagrams $\nu=(\alpha,f):(\Y',J')\To(\Y,J)$, is
\emph{cartesian}\index{word}{morphism!cartesian ---- of $\site$-diagrams}
\index{word}{cartesian morphism|see{morphism}}
 if, for any arrow $i\To j$ in $J'$, the induced commutative square
$$\xymatrix{
\Y'_i\ar[r]\ar[d]_{\alpha_i}&\Y'_j\ar[d]^{\alpha_j}\\
\Y_{f(i)}\ar[r]&\Y_{f(j)}
}$$
is cartesian.

A morphism of $\site$-diagrams $\nu=(\alpha,f):(\Y',J')\To(\Y,J)$ is \emph{reduced}
if $J=J'$ and $f=1_J$.
\end{paragr}

\begin{prop}\label{Pcartbasechange1}
Let $\nu:(\Y',J)\To (\Y,J)$ be a reduced cartesian $\Pmor$-morphism of $\site$-diagrams,
and $\varphi=(\alpha,f):(\X,I)\To (\Y,J)$ a morphism of $\site$-diagrams.
Consider the pullback square of $\site$-diagrams
$$\xymatrix{
(\X',I)\ar[r]^\mu \ar[d]_\psi &(\X,I)\ar[d]^\varphi \\
(\Y',J)\ar[r]_\nu & (\Y,J)
}$$
(i.e. $\X'_i=\Y'_{f(i)} \times_{\Y_{f(i)}} \X_i$ for any object $i$ of $I$).
Then, for any object $M$ of $\ho(\M)(\X,I)$, the canonical map
$$\derL \nu^*\, \derR \varphi_*(M)\To \derR \psi_*\, \derL \mu^*(M)$$
is an isomorphism in $\ho(\M)(\Y',J)$.
\end{prop}

\begin{proof}
By virtue of Proposition \ref{projdiamodcat}, it is sufficient to prove that
the map
$$j^* \derL \nu^*\, \derR \varphi_*(M)\To j^* \derR \psi_*\, \derL \mu^*(M)$$
is an isomorphism for any object $j$ of $J$.
Let $p:(\X/j,I/j)\To \Y_j$ and $q:(\X'/j,J,j)\To \Y'_j$ be the canonical maps.
As $\nu$ is cartesian, we have a pullback square of $\site$-diagrams
$$\xymatrix{
(\X'/j,I/j)\ar[r]^{\mu/j}\ar[d]_q&(\X/j,I/j)\ar[d]^p\\
\Y'_j\ar[r]_{\nu_j}&\Y_j
}$$
But $\nu_j$ being a $\Pmor$-morphism,
by virtue of Proposition \ref{basicPbasechange}, we thus have an isomorphism
$$\derL \nu^*_j\, \derR p_*(M/j)\simeq \derR q_*\, \derL(\mu/j)^*(M/j)
=\derR q_*(\derL \mu^*(M)/j)\, .$$
Applying Proposition \ref{computeadjointsdiagrams3}
and the last assertion of Proposition \ref{computeadjointsdiagrams2} twice, we also
have canonical isomorphisms
$$j^* \derR\varphi_*(M)\simeq \derR p_*(M/j)\quad\text{and}\quad
j^* \derR\psi_*\, \derL\mu^*(M)\simeq \derR q_*(\derL \mu^*(M)/j) \, .$$
The obvious identity $j^* \derL \nu^*=\derL\nu^*_j j^*$
achieves the proof.
\end{proof}

\begin{cor}\label{Pcartbasechange2}
Under the assumptions of Proposition \ref{Pcartbasechange1},
for any object $N$ of the category $\ho(\M)(\Y',j)$, 
the canonical map
$$\derL\mu_\sharp\, \derL\psi^*(N)\To \derL\varphi^*\, \derL\nu_\sharp(N)$$
is an isomorphism in $\ho(\M)(\X,I)$.
\end{cor}

\begin{rem}\label{remPacrtfibredcat}
The class of cartesian $\Pmor$-morphisms form an admissible
class of morphisms in the category of $\site$-diagrams, which we denote
by $\Pmor_\cart$.
\index{notat}{pmorcart@$\Pmor_\cart$}
 Proposition \ref{basicfunctdiag} and the preceding corollary
thus asserts that $\ho(\M)$ is a $\Pmor_\cart$-fibred category over the category
of $\site$-diagrams.
\end{rem}

\begin{paragr}\label{diagdiag}
We shall sometimes deal with diagrams of $\site$-diagrams.
Let $I$ be a small category, and $\F$ a functor from $I$ to the category
of $\site$-diagrams. For each object $i$ of $I$, we have
a $\site$-diagram $(\F(i),J_i)$, and, for each map $u:i\to i'$, we have
a functor $f_u:J_i\To J_{i'}$ as well as a natural transformation
$\alpha_u:\F(i)\To\F(i')\circ f_u$, subject to coherence identities.
In particular, the correspondence $i\mapsto J_i$ defines a
functor from $I$ to the category of small categories.
Let $I_\F$ be the cofibred category over $I$ associated to it; see \cite[Exp. VI]{SGA1}.
Explicitly, $I_\F$ is described as follows. The objects are the couples
$(i,x)$, where $i$ is an object of $I$, and $x$ is an object of $J_i$.
A morphism $(i,x)\To(i',x')$ is a couple $(u,v)$, where $u:i\To i'$ is
a morphism of $I$, and $v:f_u(x)\To x'$ is a morphism of $J_{i'}$.
The identity of $(i,x)$ is the couple $(1_i,1_x)$, and,
for two morphisms $(u,v):(i,x)\To(i',x')$ and $(u',v'):(i',x')\To(i'',x'')$,
their composition $(u'',v''):(i,x)\To(i'',x'')$ is defined by $u''=u'\circ u$,
while $v''$ is the composition of the map
$$f_{u''}(x)=f_{u'}(f_u(x))\xrightarrow{f_{u'}(v)}f_{u'}(x')\xrightarrow{\ v'\ }x''\, .$$
The functor $p:I_\F\To I$ is simply the projection $(i,x)\mapsto i$.
For each object $i$ of $I$, we get a canonical pullback square of categories
\begin{equation}\label{diagdiagfiber1}
\begin{split}
\xymatrix{J_i\ar[d]_q\ar[r]^{\ell_i}&I_\F\ar[d]^p\\
e\ar[r]_i&I}
\end{split}\end{equation}
in which $i$ is the functor from the terminal category $e$ which
corresponds to the object $i$, and $\ell_i$ is the functor defined
by $\ell_i(x)=(i,x)$.

The functor $\F$ defines a $\site$-diagram $(\Int\F,I_\F)$: for an object $(i,x)$ of $I_\F$,
$(\Int\F)_{(i,x)}=\F(i)_x$, and for a morphism $(u,v):(i,x)\To(i',x')$, the map
$$(u,v):(\Int\F)_{(i,x)}=\F(i)_x\To(\Int\F)_{(i',x')}=\F(i')_{x'}$$
is simply the morphism induced by $\alpha_u$ and $v$.
For each object $i$ of $I$, there is a natural morphism of $\site$-diagrams
\begin{equation}\label{diagdiagfiber2}
\lambda_i: (\F(i),J_i)\To(\Int\F,I_\F)\, ,
\end{equation}
given by $\lambda_i=(1_{\F(i)},\ell_i)$
%% Consider now a $\site$ diagram $\X$ indexed by $I$, and
%% denote by $\tilde\X$ the functor from $I$ to $\site$-diagrams
%% defined by $\X$ (i.e. the functor which associates to each object $i$ of $I$
%% the functor from $e$ to $\site$ with value $\X_i$). Let
%% $f:\F\To \tilde \X$ be a morphism of functors (where $\X$ is considered here
%% as a functor from $I$ to $\site$-diagrams, each $\X_i$ being seen as a
%% functor from $e$ to $\site$). Then we have a canonical
%% isomorphism of $\site$-diagrams $(\Int \tilde \X,I_{\tilde \X})=(\X,I)$, and $f$
%% induces a morphism of $\site$-diagrams
%% \begin{equation}\label{diagdiagfiber3}
%% \varphi: (\Int \F, I_\F)\To (\X,I)
%% \end{equation}
%% whose underlying functor $I_\F\To I$ is simply the projection $p$.
%% Consider now a cartesian morphism $\nu:(\X',I')\To(\X,I)$
%% with underlying functor $a:I'\To I$. We get
%% a functor $a^*(\F)=\F\circ a$ from $I'$ to the category of $\site$-diagrams, and
%% one checks easily that there is a canonical pullback square of $\site$-diagrams
%% \begin{equation}\label{diagdiagfiber4}\begin{split}\xymatrix{
%% (\Int^*(\F),I'_{a^*(\F)})\ar[r]^{\mu}\ar[d]_{\varphi'}&(\Int\F,I_\F) \ar[d]^\varphi\\
%% (\X ',I')\ar[r]_\nu &(\X ,I)
%% }\end{split}\end{equation}
\end{paragr}
%% 
%% \begin{prop}\label{directimagediagdiag0}
%% If $\nu:(\X',I')\To(\X,I)$ is a cartesian $\Pmor$-morphism, then, for any object
%% $M$ of $\ho(\M)(\F,I_\F)$, the exchange map
%% $$\nu^*\, \derR\varphi_*(M)\To \derR\varphi'_*\, \mu^*(M)$$
%% is an isomorphism in $\ho(\M)(\X',I')$.
%% \end{prop}

\begin{prop}\label{directimagediagdiag}
Let $X$ be an object of $\site$, and $f:\F\To X$ a morphism
of functors (where $X$ is considered
as the constant functor from $I$ to $\site$-diagrams
with value the functor from $e$ to $\site$ defined by $X$).
Then, for each object $i$ of $I$,
we have a canonical pullback square of $\site$-diagrams
$$\xymatrix{
(\F(i),J_i)\ar[r]^{\lambda_i}\ar[d]_{\varphi_i}&(\Int\F,I_\F) \ar[d]^\varphi\\
X\ar[r]_i&(X,I)
}$$
in which $\varphi$ and $\varphi_i$ are the obvious morphisms induced by $f$
(where, this time, $(X,I)$ is seen as the constant functor from $I$ to $\site$
with value $X$).

Moreover, for any object $M$ of $\ho(\M)(\Int\F,I_\F)$, the natural map
$$i^*\, \derR\varphi_*(M)=\derR\varphi_*(M)_i\To \derR \varphi_{i,*}\, \lambda^*_i(M)$$
is an isomorphism. In particular, if we also write by abuse of notation $f$
for the induced map of $\site$-diagrams from $(\Int\F,I_\F)$ to $X$, we have a natural isomorphism
$$\derR f_*(M)\simeq
\underset{\quad i\in \op{I} }{\derR\varprojlim} \derR \varphi_{i,*}\, \lambda^*_i(M)\, .$$
\end{prop}

\begin{proof}
This pullback square is the one induced by \eqref{diagdiagfiber1}.
We shall prove first that the map
$$i^*\, \derR\varphi_*(M)=\derR\varphi_*(M)_i\To \derR \varphi_{i,*}\, \lambda^*_i(M)$$
is an isomorphism in the particular case where $I$
has a terminal object $\omega$ and $i=\omega$.
By virtue of Propositions \ref{computeadjointsdiagrams2}
and \ref{computeadjointsdiagrams3}, we have
isomorphisms
\begin{equation}\label{directimagediagdiag1}
\omega^*\, \derR\varphi_*(M)\simeq
\underset{\quad i\in \op{I} }{\derR\varprojlim}\derR\varphi_*(M)_i\simeq
\underset{\quad (i,x)\in \op{I}_\F }{\derR\varprojlim}
\derR \varphi_{i,x,*}\, (M_{(i,x)})\, ,
\end{equation}
where $\varphi_{i,x}:\F(i)_x\To X$ denotes the map induced
by $f$.
We are thus reduced to prove that the canonical map
\begin{equation}\label{directimagediagdiag2}
\underset{\quad (i,x)\in \op{I}_\F }{\derR\varprojlim} \derR \varphi_{i,x,*}\, (M_{(i,x)})
\To \underset{\quad x\in \op{J}_\omega }{\derR\varprojlim}
\derR \varphi_{\omega ,x,*}\, (M_{(\omega,x)})
\simeq  \derR \varphi_{\omega,*}\, \lambda^*_\omega(M)
\end{equation}
is an isomorphism. As $I_\F$ is cofibred over $I$, and as $\omega$
is a terminal object of $I$, the inclusion functor $\ell_\omega : J_\omega\To I_\F$
has a left adjoint, whence is coaspherical in any weak basic localizer
(i.e. is homotopy cofinal); see \cite[1.1.9, 1.1.16 and 1.1.25]{Mal}.
As any model category defines a Grothendieck derivator
\index{word}{derivator, Grothendieck}
 (\cite[Thm. 6.11]{Cis1}), it follows from \cite[Cor. 1.15]{Cis1}
 that the map \eqref{directimagediagdiag2} is an isomorphism.

To prove the general case, we proceed as follows.
Let $\F/i$ be the functor obtained by composing $\F$
 with the canonical functor $v_i:I/i\To I$. 
 Then, keeping track of the conventions adopted 
 in \ref{computeadjointsdiagrams},
 we check easily that
$(I/i)_{\F/i}=(I_\F)/i$ and that $\Int (\F/i)=(\Int \F)/i$.
Moreover, the pullback square \eqref{diagdiagfiber1} is the composition
of the following pullback squares of categories.
$$\xymatrix{
J_i\ar[r]^{a_i}\ar[d]_q&I_\F/i \ar[r]^{u_i}\ar[d]_{p/i}& I_\F\ar[d]^p\\
e\ar[r]_{(i,1_i)}& I/i\ar[r]_{v_i}& I
}$$
The pullback square of the proposition is thus
the composition of the following pullback squares.
$$\xymatrix{
(\F(i),J_i)\ar[r]^{\alpha_i}\ar[d]_{\varphi_i}
&(\Int\F/i,I_\F/i) \ar[r]^{\mu_i}\ar[d]_{\varphi /i}&(\Int \F,I_\F)\ar[d]^\varphi \\
X\ar[r]_{(i,1_i)}&(X,I/i)\ar[r]_{v_i}& (X,I)
}$$
The natural transformations
$$(i,1_i)^*\, \derR (\varphi /i)_*\To \derR \varphi_{i,*}\, \alpha^*_i\quad\text{and}
\quad v_i^*\, \derR \varphi_* \To \derR (\varphi /i)_* \, \mu^*_i$$
are both isomorphisms: the first one comes from the fact that $(i,1_i)$
is a terminal object of $I/i$, and the second one
from Proposition \ref{computeadjointsdiagrams3}. We thus get:
\begin{align*}
i^*\, \derR \varphi_*(M)
& \simeq (i,1_i)^*\, v^*_i \, \derR\varphi_*(M)\\
& \simeq (i,1_i)^*\, \derR (\varphi /i)_*\, \mu^*_i(M) \\
& \simeq  \derR \varphi_{i,*}\, \alpha^*_i \, \mu^*_i(M)\\
& \simeq \derR \varphi_{i,*}\, \lambda^*_i(M)\, .
\end{align*}
The last assertion of the proposition is then a straightforward
application of Proposition \ref{computeadjointsdiagrams2}.
\end{proof}

\begin{prop}\label{derivedtensordiaginj}
If $\M$ is a monoidal $\Pmor$-fibred combinatorial model category over $\site$,
then, for any $\site$-diagram $\X$ indexed by a small category $I$,
the injective model structure turns
$\M(\X,I)$ into a symmetric monoidal model category.
In particular, the categories $\ho(\M)(\X,I)$ are canonically endowed with
a closed symmetric monoidal structure, in such a way that, for any morphism
of $\site$-diagrams $\varphi:(\X,I)\To(\Y,J)$, the functor
$\derL \varphi^*:\ho(\M)(\Y,J)\To\ho(\M)(\X,I)$ is symmetric monoidal.
\end{prop}

\begin{proof}
This is obvious from the definition of a symmetric
monoidal model category, as the tensor product of $\M(\X,I)$ is defined termwise,
as well as the cofibrations and the trivial cofibrations.
\end{proof}
%\begin{prop}\label{derivedtensordiagproj}
%If $\M$ is a symmetric monoidal cofibrantly generated $\Pmor$-fibred model category over $\site$,
%then, for any $\site$-diagram $\X$ indexed by a small category $I$
%which admits finite products, the projective model structure turns
%$\M(\X,I)$ into a symmetric monoidal model category. 
%\end{prop}
%\begin{proof}
%\end{proof}

\begin{prop}\label{diagprojformula}
Assume that $\M$ is a monoidal $\Pmor$-fibred combinatorial model category over $\site$,
and consider a reduced cartesian $\Pmor$-morphism $\varphi=(\alpha,f):(\X,I)\To(\Y,I)$.
Then, for any object $M$ in $\ho(\M)(\X,I)$ and any object $N$ in $\ho(\M)(\Y,I)$, the
canonical map
$$\derL \varphi_\sharp(M\otimes^\derL \varphi^*(N))\To \derL \varphi_\sharp(M)\otimes^\derL N$$
is an isomorphism.
\end{prop}

\begin{proof}
Let $i$ be an object of $I$. It is sufficient to prove that the map
$$i^* \derL \varphi_\sharp(M\otimes^\derL \varphi^*(N))\To i^* \derL \varphi_\sharp(M)\otimes^\derL N$$
is an isomorphism in $\ho(\M)(\X_i)$. Using Corollary \ref{Pcartbasechange2},
we see that this map can be identified with the map
$$\derL \varphi_{i,\sharp}(M_i\otimes^\derL \varphi^*_i(N_i))
\To \derL \varphi_{i,\sharp}(M_i)\otimes^\derL N_i\, ,$$
which is an isomorphism according to the $\Pmor$-projection formula
for the homotopy $\Pmor$-fibred category $\ho(\M)$.
\end{proof}

\begin{paragr}\label{defhomotopycartsection}
Let $(\X,I)$ be a $\site$-diagram. An object $M$ of $\M(\X,I)$
 is \emph{homotopy cartesian}
\index{word}{homotopy cartesian!object over a diagram}
 if, for any map $u:i\To j$ in $I$,
the structural map $u^*(M_j)\To M_i$ induces an isomorphism
$$\derL u^*(M_i)\simeq M_j$$
in $\ho(\M)(\X,I)$ (i.e. if there exists a weak equivalence $M'_j\To M_j$
with $M'_j$ cofibrant in $\M(\X_j)$ such that the map
$u^*(M'_j)\To M_i$ is a weak equivalence in $\M(\X_i)$).

We denote by $\ho(\M)(\X,I)_\hcart$ the full subcategory of
$\ho(\M)(\X,I)$ spanned by homotopy cartesian sections.
\end{paragr}

\begin{df}\label{deftractable}
A cofibrantly generated model category $\V$ is \emph{tractable}
\index{word}{tractable}
if there exist sets $I$ and $J$ of cofibrations between
cofibrant objects which generate the class of
cofibrations and the class of trivial cofibrations respectively.
\end{df}

\begin{rem}
If $\M$ is a combinatorial
 and tractable $\Pmor$-fibred model category over $\site$, 
 then so are the projective and the injective model structures on $\M(\X,I)$;
see \cite[Thm. 2.28 and 2.30]{Bar}.
\end{rem}

\begin{prop}\label{existencehcartsections}
If $\M$ is tractable, then the inclusion functor
$$\ho(\M)(\X,I)_\hcart\To\ho(\M)(\X,I)$$
admits a right adjoint.
\end{prop}

\begin{proof}
This follows from the fact that the cofibrant homotopy cartesian
sections are the cofibrant objects of a right Bousfield localization
of the injective model structure on $\M(\X,I)$; see
\cite[Theorem 5.25]{Bar}.
\end{proof}

\begin{df}\label{defQuillenPfibredmorphism}
Let $\M$ and $\M'$ two $\Pmor$-fibred model categories over $\site$.
A \emph{Quillen morphism}
\index{word}{morphism!Quillen ---- of $\Pmor$-fibred model categories}
 $\gamma$ from $\M$ to $\M'$
is a morphism of $\Pmor$-fibred categories $\gamma:\M\To\M'$
such that $\gamma^*:\M(X)\To\M'(X)$ is a left Quillen functor
for any object $X$ of $\site$.
\end{df}

\begin{rem}
If $\gamma:\M\To\M'$ is a Quillen morphism between $\Pmor$-fibred
combinatorial model categories, then, for any $\site$-diagram $(\X,I)$,
we get a Quillen adjunction
$$\gamma^*:\M(\X,I)\rightleftarrows\M'(\X,I):\gamma_*$$
(with the injective model structures as well as with the projective
model structures).
\end{rem}

\begin{prop}\label{derivedPfibredQuillenfunct}
For any Quillen morphism $\gamma:\M\To\M'$,
the derived adjunction
$$\derL\gamma^*:\ho(\M)(X)\rightleftarrows\ho(\M')(X):\derR\gamma_*$$
defines a morphism of $\Pmor$-fibred categories $\ho(\M)\To\ho(\M')$ over $\site$.
If moreover $\M$ and $\M'$ are combinatorial, then the
morphism $\ho(\M)\To\ho(\M')$ extends to a morphism of
$\Pmor_\cart$-fibred categories over the category of $\site$-diagrams.
\end{prop}

\begin{proof}
This follows immediately from \cite[Theorem 1.4.3]{Hovey}.
\end{proof}

\subsection{Hypercovers, descent, and derived global sections}

\begin{paragr}\label{hypothesissite}
Let $\site$ be an essentially small category,
and $\Pmor$ an admissible class of morphisms in $\site$.
We assume that a Grothendieck topology $t$ on $\site$ is given.
We shall write $\site^\amalg$ for the full subcategory of the category
of $\site$-diagrams whose objects are the small families $X=\{X_i\}_{i\in I}$ of
objects of $\site$ (seen as functors from a discrete category
to $\site$). The category $\site^\amalg$ is equivalent to the
full subcategory of the category of presheaves of sets on $\site$
spanned by sums of representable presheaves.
In particular, small sums are representable in $\site^\amalg$
(but note that the functor from $\site$ to $\site^\amalg$
does not preserve sums). Finally, we remark that the topology
$t$ extends naturally to a Grothendieck topology on $\site^\amalg$
such that the topology $t$ on $\site$ is the topology induced
from the inclusion $\site\subset\site^\amalg$. The covering maps
for this topology on $\site^\amalg$ will be called \emph{$t$-covers}\index{word}{cover}
(note that the inclusion $\site\subset\site^\amalg$ is continuous and induces
an equivalence between the topos of $t$-sheaves on $\site$ and the topos
of $t$-sheaves on $\site^\amalg$).

Let $\Delta$ be the category of non-empty finite ordinals.
Remember that a simplicial object of $\site^\amalg$ is a presheaf on $\Delta$
with values in $\site^\amalg$.
For a simplicial set $K$ and an object $X$ of $\site^\amalg$, we denote by
$K\times X$ the simplicial object of $\site^\amalg$ defined by
$$(K\times X)_n=\coprod_{x\in K_n} X\quad , \qquad n\geq 0 \, .$$
We write $\Delta^n$ for the standard combinatorial simplex of dimension $n$,
and $i_n:\partial\Delta^n\To\Delta^n$ for its boundary inclusion.

A morphism $p:\X\To\Y$ between simplicial objects of $\site^\amalg$
is a \emph{$t$-hypercover}\index{word}{hypercover}
 if, locally for the $t$-topology,
it has the right lifting property with respect to boundary inclusions of standard simplices,
which, in a more precise way, means that, for any integer $n\geq 0$,
any object $U$ of $\site^\amalg$, and any commutative square
$$\xymatrix{
&\partial\Delta^n\times U\ar[r]^{x}\ar[d]_{{i_n\times 1}}&\X\ar[d]^p&\\
&\Delta^n\times U\ar[r]_y&\Y&,
}$$
there exists a $t$-covering $q:V\To U$, and a morphism
of simplicial objects $z:\Delta^n\times V\To \X$, such that the
diagram bellow commutes.
$$\xymatrix{
\partial\Delta^n\times V\ar[r]^(.63){x(1\times q)}\ar[d]_{{i_n\times 1}}&\X\ar[d]^p\\
\Delta^n\times V\ar[r]_(.63){y(1\times q)}\ar[ur]^z&\Y
}$$
A \emph{$t$-hypercover} of an object $X$ of $\site^\amalg$ is a a $t$-hypercover
$p:\X\To X$ (where $X$ is considered as a constant simplicial object).
\end{paragr}

\begin{rem}\label{remord}
This definition of $t$-hypercover
is equivalent to the one given in \cite[Exp.~V, 7.3.1.4]{SGA4}.
\end{rem}

\begin{paragr}\label{prepdefdescente}
Let $\X$ be a simplicial object of $\site^\amalg$.
It is in particular a functor from the category $\op{\Delta}$
to the category of $\site$-diagrams, so that the constructions and
considerations of \ref{diagdiag} apply to $\X$.
In particular, there is a $\site$-diagram $\tilde\X$ associated to
$\X$, namely $\tilde\X=(\Int \X, (\op{\Delta})_\X)$.
More explicitly, for each integer $n\geq 0$, there is
a family $\{\X_{n,x}\}_{x\in K_n}$ of objects of $\site$, such that
\begin{equation}\label{presentsimplhypercover}
\X_n=\coprod_{x\in K_n}\X_{n,x}\, .
\end{equation}
In fact, the sets $K_n$ form a simplicial set $K$, and the
category $(\op{\Delta})_\X$ can be identified over $\op{\Delta}$
to the category $\op{(\Delta/K)}$, where $\Delta/K$ is the
fibred category over $\Delta$ whose fiber over $n$ is the set $K_n$
(seen as a discrete category), i.e. the category of simplices of $K$.
We shall call $K$ the \emph{underlying simplicial set of $\X$},
\index{word}{underlying simplicial set!of a simplicial object}
while the decomposition \eqref{presentsimplhypercover}
will be called the 
\emph{local presentation\index{word}{presentation!local presentation
 of a simplicial object} of $\X$}.
The construction $\X\mapsto\tilde\X$ is functorial.
If $p:\X\To\Y$ is a morphism of simplicial objects of $\site^\amalg$,
we shall still denote by $p:\tilde \X\To \tilde \Y$ the induced
morphism of $\site$-diagrams. In particular, for
a morphism of $p:\X\To X$, where $X$ is an object of $\site^\amalg$,
$p:\tilde \X\To X$ denotes the corresponding morphism of $\site$-diagrams.

Let $\M$ be a $\Pmor$-fibred combinatorial model category over $\site$.
Given a simplicial object $\X$ of $\site^\amalg$, we define the category
$\ho(\M)(\X)$ by the formula:
\begin{equation}\label{defHoMhyperrec}
\ho(\M)(\X)=\ho(\M)(\Int \X, (\op{\Delta})_\X)\, .
\end{equation}
Given an object $X$ of $\site^\amalg$ and a morphism $p:\X\to X$,
we have a derived adjunction
\begin{equation}\label{defHoMhyperrecfunct}
\derL p^* : \ho(\M)(X)\rightleftarrows \ho(\M)(\X):\derR p_* \, .
\end{equation}
\end{paragr}

\begin{prop}\label{computesimpldirectimage}
Consider an object $X$ of $\site$, a simplicial object $\X$
of $\site^\amalg$, as well as a morphism $p:\X\To X$.
Denote by $K$ the underlying simplicial set of $\X$,
and for each integer $n\geq 0$ and each simplex $x\in K_n$,
write $p_{n,x}:\X_{n,x}\To X$ for the morphism of $\site^\amalg$
induced by the local presentation of $\X$ \eqref{presentsimplhypercover}.
Then, for any object $M$ of $\ho(\M)(X)$, there are canonical
isomorphisms
$$\derR p_*\derR p^*(M)\simeq
\underset{\quad n\in \Delta}{\derR\varprojlim} \derR p_{n,*}\derL p^*_{n}(M)
\simeq \underset{\quad n\in \Delta}{\derR\varprojlim}
\Big(\prod_{x\in K_n}\derR p_{n,x,*}\derL p^*_{n,x}(M)\Big)\, .$$
\end{prop}

\begin{proof}
The first isomorphism is a direct application of the last assertion of
Proposition \ref{directimagediagdiag}
for $\F=\X$, while the second one follows from the first one by Proposition \ref{sumCdiag}.
\end{proof}

\begin{df}\label{defdescente}
Given an object $Y$ of $\site^\amalg$, an object $M$ of $\ho(\M)(Y)$
will be said to satisfy \emph{$t$-descent}\index{word}{descent!tdescent@$t$-descent}
if it has the following property:
for any morphism $f:X\To Y$ and any $t$-hypercover $p:\X\To X$,
the map
$$\derR f_*\, \derL f^*(M)\To \derR f_*\, \derR p_*\, \derL p^*\, \derL f^*(M)$$
is an isomorphism in $\ho(\M)(Y)$.

We shall say that $\M$ (or by abuse, that $\ho(\M)$) satisfies \emph{$t$-descent} if,
for any object $Y$ of $\site^\amalg$, any object of $\ho(\M)(Y)$ satisfies
$t$-descent.
\end{df}

\begin{prop}\label{triv}
If $Y=\{Y_i \}_{i\in I}$ is a small family of objects of $\site$
(seen as an object of $\site^\amalg$), then an object $M$ of
$\ho(\M)(Y)$ satisfies $t$-descent if and only if, for any $i\in I$,
any morphism $f:X\to Y_i$ of $\site$, and
any $t$-hypercover $p:\X\To X$, the map
$$\derR f_*\, \derL f^*(M_i)\To \derR f_*\, \derR p_*\, \derL p^*\, \derL f^*(M_i)$$
is an isomorphism in $\ho(\M)(Y_i)$.
\end{prop}

\begin{proof}
This follows from the definition and from
Proposition \ref{sumCdiag}.
\end{proof}

\begin{cor}\label{trivglobchartaudescent}
The $\Pmor$-fibred model category $\M$ satisfies $t$-descent if and only if,
for any object $X$ of $\site$, and
any $t$-hypercover $p:\X\To X$, the functor
$$\derL p^*:\ho(\M)(X)\To\ho(\M)(\X)$$
is fully faithful.
\end{cor}

\begin{prop}\label{coverconservative}
If $\M$ satisfies $t$-descent, then,
for any $t$-cover $f:Y\To X$, the functor
$$\derL f^*:\ho(\M)(X)\To\ho(\M)(Y)$$
is conservative.
\end{prop}

\begin{proof}
Let $f:Y\To X$ be a $t$-cover,
and $u:M\to M'$ a morphism of $\ho(\M)(X)$
whose image by $\derL f^*$ is an isomorphism.
We can  consider the \v{C}ech $t$-hypercover
associated to $f$, that is the simplicial object $\Y$ over $X$
defined by
$$\Y_n=\underset{\text{$n+1$ times}}{\underbrace{Y\times_X Y\times_X\dots \times_X Y}}\, .$$
Let $p:\Y\To X$ be the canonical map. For each $n\geq 0$, the map
$p_n:\Y_n\To X$ factor through $f$, from which we deduce that the functor
$$\derL p^*_n:\ho(\M)(X)\To\ho(\M)(\Y_n)$$
sends $u$ to an isomorphism. This implies that the functor
$$\derL p^*:\ho(\M)(X)\To\ho(\M)(\Y)$$
sends $u$ to an isomorphism as well. But, as $\Y$ is a $t$-hypercover of
$X$, the functor $\derL p^*$ is fully faithful, from which we deduce that $u$
is an isomorphism by the Yoneda Lemma.
\end{proof}

\begin{paragr}\label{enrichedpref}
Let $\V$ be a complete and cocomplete category.
For an object $X$ of $\site$, define $\pshg{\site/X,\V}$
as the category of presheaves
on $\site/X$ with values in $\V$.
Then $\pshg{C/-,\V}$
is a $\Pmor$-fibred category (where, by abuse of
notations, $\site$ denotes also the class of all maps in $\site$):
this is a special case of the
constructions explained in \ref{defevaluation}
applied to $\V$, seen as a fibred category over the
terminal category. To be more explicit,
for each object $X$ of $\site^\amalg$, we have a $\V$-enriched Yoneda embedding
\begin{equation}\label{enrichedpref1}
\site^\amalg/X\times\V\To \pshg{\site/X,\V}\quad , \qquad (U,M\}\mapsto U\otimes M\, ,
\end{equation}
where, if $U=\{U_i\}_{i\in I}$ is a small family of objects of $\site/X$, $U\otimes M$ is the presheaf
\begin{equation}\label{enrichedpref2}
V\mapsto \coprod_{i\in I} \ \coprod_{a\in \Hom_{\site/S}(V,U_i)} M\, .
\end{equation}
For a morphism $f:X\To Y$ in $\site$, the functor
$$f^*:\pshg{\site/Y,\V}\To\pshg{\site/X,\V}$$
is the functor defined by composition with the corresponding functor $\site/X\To\site/Y$.
The functor $f^*$ has always a left adjoint
$$f_\sharp: \pshg{\site/X,\V}\To\pshg{\site/Y,\V}\, ,$$
which is the unique colimit preserving functor defined by
$$f_\sharp(U\otimes M)=U\otimes M\, ,$$
where, on the left hand side $U$ is considered as an object over $X$,
while, on the right hand side, $U$ is considered as an object over $Y$
by composition with $f$. Similarly, if all the pullbacks by $f$
are representable in $\site$ (e.g. if $f$ is a $\Pmor$-morphism),
the functor $f^*$ can be described as the colimit preserving functor
defined by the formula
$$f^*(U\otimes M)=(X\times_Y U)\otimes M\, .$$

If $\V$ is a cofibrantly generated model category, then,
for each object $X$ of $\site$, the category $\pshg{\site/X,\V}$
is naturally endowed with the \emph{projective model category
structure}, i.e. with the cofibrantly generated model category
structure whose weak equivalences and fibrations are defined
termwise (this is Proposition \ref{projdiamodcat}
applied to $\V$, seen as a fibred category over the terminal category).
The cofibrations of the projective model category structure on $\pshg{\site/X,\V}$
will be called the \emph{projective cofibrations.}
If moreover $\V$ is combinatorial (resp.
left proper, resp. right proper, resp. stable), so is $\pshg{\site/X,\V}$.
In particular, if $\V$ is a combinatorial model category, then
$\pshg{\site/-,\V}$ is a $\Pmor$-fibred combinatorial model category over $\site$.

According to Definition \ref{defdescente},
it thus makes sense to speak of $t$-descent in $\pshg{\site/-,\V}$.

If $U=\{U_i\}_{i\in I}$ is a small family of objects of $\site$ over $X$,
and if $F$ is a presheaf over $\site/X$, we define
\begin{equation}\label{enrichedpref3}
F(U)=\prod_{i\in I}F(U_i)\, .
\end{equation}
the functor $F\mapsto F(U)$ is a right adjoint
to the functor $E\mapsto U\otimes E$.

We remark that a termwise fibrant presheaf $F$ on $\site/X$
satisfies $t$-descent if and only if, for any object $Y$
of $\site^\amalg$, and any $t$-hypercover
$\Y\To Y$ over $X$, the map
$$F(Y)\To \underset{\phantom{Ru}n\in \Delta}{\derR\varprojlim} F(\Y_n)$$
is an isomorphism in $\ho(\V)$.
\end{paragr}

\begin{prop}\label{taulocmodcat}
If $\V$ is combinatorial and left proper, then the category of presheaves
$\pshg{\site/X,\V}$
admits a combinatorial model category structure whose cofibrations are
the projective cofibrations, and whose fibrant objects are the termwise
fibrant objects which satisfy $t$-descent. This model
category structure will be called the \emph{$t$-local model category
structure}, and the corresponding homotopy category will be denoted by
$\ho_t(\pshg{\site/X,\V})$.

Moreover, any termwise weak equivalence is a weak equivalence for
the $t$-local model structure, and the induced functor
$$a^*:\ho(\pshg{\site/X,\V})\To\ho_t(\pshg{\site/X,\V})$$
admits a fully faithful right adjoint
$$a_*:\ho_t(\pshg{\site/X,\V})\To\ho(\pshg{\site/X,\V})$$
whose essential image consists precisely of the full subcategory of
$\ho(\pshg{\site/X,\V})$ spanned by the presheaves which satisfy
$t$-descent.
\end{prop}

\begin{proof}
Let $H$ be the class of maps of shape
\begin{equation}\label{gentaulocmodcat}
\underset{n\in \op{\Delta}}{\mathrm{hocolim}}\, \Y_n\otimes E\To Y\otimes E\, ,
\end{equation}
where $Y$ is an object of $\site^\amalg$ over $X$,
$\Y\To Y$ is a $t$-hypercover,
and $E$ is a cofibrant replacement of an object which is either
a source or a target of a generating cofibration of $\V$.
Define the $t$-local model category structure as the left Bousfield
localization of $\mathit{Pr}(\site/X,\V)$ by $H$; see \cite[Theorem 4.7]{Bar}.
We shall call \emph{$t$-local weak equivalences} the weak equivalences
of the $t$-local model category structure.
For each object $Y$ over $X$, the functor $Y\otimes (-)$ is
a left Quillen functor from $\V$ to $\mathit{Pr}(\site/X,\V)$.
We thus get a total left derived functor
$$Y\otimes^\derL(-):\ho(\V)\To\ho_t(\pshg{\site/X,\V})$$
whose right adjoint is the evaluation at $Y$.
For any object $E$ of $\V$ and any $t$-local fibrant presheaf $F$
on $\site/X$ with values in $\V$, we thus have natural bijections
\begin{equation}\label{gentaulocmodcat2}
\Hom(E,F(Y))\simeq\Hom(Y\otimes^\derL E,F)\, ,
\end{equation}
and, for any simplicial object $\Y$ of $\site/X$, identifications
\begin{equation}\label{gentaulocmodcat3}
\Hom(E, \underset{\phantom{Ru}n\in \Delta}{\derR\varprojlim} F(\Y_n))
\simeq\Hom
( \underset{\phantom{Ru}n\in \Delta}{\derL\varinjlim} \Y_n \otimes^\derL E,F)\, ,
\end{equation}
One sees easily that, for any $t$-hypercover $\Y\To Y$
and any cofibrant object $E$ of $\V$, the map
\begin{equation}\label{gentaulocmodcat4}
\underset{\phantom{Ru}n\in \Delta}{\derL\varinjlim} \Y_n \otimes^\derL E
\To Y\otimes^\derL E
\end{equation}
is an isomorphism in the $t$-local homotopy category $\ho_t(\pshg{\site/X,\V})$:
by the small object argument, the smallest full subcategory of $\ho(\pshg{\site/X,\V})$
which is stable by homotopy colimits and which contains the source and the targets
of the generating cofibrations is $\ho_t(\pshg{\site/X,\V})$ itself, and
the class of objects $E$ of $\V$ such that the map \eqref{gentaulocmodcat4}
is an isomorphism in $\ho(\V)$ is sable by homotopy colimits.
Similarly, we see that, for any object $E$, the functor $(-)\otimes^\derL E$
preserves sums. As a consequence, we get from \eqref{gentaulocmodcat2}
and \eqref{gentaulocmodcat3} that the fibrant objects
of the $t$-local model category structure are precisely
the termwise fibrant objects $F$ of the projective model structure
which satisfy $t$-descent.
The last part of the proposition follows from the general yoga of left Bousfield
localizations.
\end{proof}

%% \begin{paragr}\label{defhomotopyenrichment}
%% Let $\V$ be  a symmetric monoidal model category,
%% and $\M$ a \emph{$\V$-enriched}
%% $\Pmor$-fibred combinatorial symmetric monoidal model category over $\site$, that is
%% $\Pmor$-fibred combinatorial symmetric monoidal model category over $\site$
%% endowed with a symmetric monoidal left Quillen functor
%% $$K:\V\To \M(e)\, ,$$
%% where $e$ stands for the terminal object of $\site$
%% (for instance, any $\Pmor$-fibred combinatorial symmetric monoidal model category $\M$ over $\site$
%% is $\M(e)$-enriched). We denote by
%% $$\Gamma:\M(e)\To\V$$
%% the right adjoint of $K$. We then have a
%% total right derived functor
%% $$\derR\Gamma:\ho(\M)(e)\To\ho(\V)\, .$$
%% 
%% Let $S$ be a fixed object of $\site$, and $M$ an object of $\ho(\M)(S)$.
%% For any morphism $f:X\To S$ of $\site$, we define
%% \begin{equation}\label{expliciteglobsect1}
%% %\begin{aligned}
%% \derR\Gamma(X,M) =\derR\Gamma(\derR\sHom_X(\unit_X,\derL f^*(M))\\
%% \simeq \derR\Gamma(S,\derR f_* \derL f^*(M))\, .
%% %\end{aligned}
%% \end{equation}
%% This formula leads formally to the canonical identification
%% \begin{equation}\label{expliciteglobsect2}
%% \Hom_{\ho(\V)}(\unit,\derR\Gamma(X,M))=\Hom_{\ho(\M)(X)}(\unit_X,M)\, .
%% \end{equation}
%% This construction defines a presheaf on the category $\site/S$,
%% with values in the homotopy category $\ho(\V)$. This can be strictified a follows.
%% \end{paragr}

\begin{paragr}\label{strictifyglobalsection}
Let $\M$ be a $\Pmor$-fibred combinatorial model category over $\site$, and
$S$ an object of $\site$ . Denote by
$$\S:\site/S\To \site$$
the canonical forgetful functor. Then there is a canonical morphism
of $\site$-diagrams
\begin{equation}\label{strictifyglobalsection1}
\sigma: (\S,\site/S)\To(S,\site/S)
\end{equation}
(where $(S,\site/S)$ stands for the constant diagram with value $S$).
This defines a functor
\begin{equation}\label{strictifyglobalsection2}
\derR\sigma_*:\ho(\M)(\S,\site/S)\To\ho(\M)(S,\site/S)=\ho(\pshg{\site/S,\M(S)}) \, .
\end{equation}
For an object $M$ of $\ho(\M)(S)$, one defines the presheaf
of \emph{geometric derived global sections
\index{word}{section!geometric derived global section}
\index{word}{global section|see{section}} of $M$ over $S$}
by the formula
\begin{equation}\label{strictifyglobalsection3}
\derR\Gamma_{\mathit{geom}}(-,M)=\derR\sigma_*\, \derL\sigma^*(M)\, .
\end{equation}
This is a presheaf on $\site/S$ with values in $\M(S)$
whose evaluation on a morphism $f:X\To S$ is, by virtue of
Propositions \ref{computeadjointsdiagrams2} and \ref{computeadjointsdiagrams3},
\begin{equation}\label{strictifyglobalsection4}
\derR\Gamma_{\mathit{geom}}(X,M)\simeq \derR f_*\, \derL f^*(M)\, .
\end{equation}
%% By Proposition \ref{taulocmodcat}, we have an adjunction
%% $$a^*:\ho(\mathit{Pr}(\site/S,\M(S))\rightleftarrows\ho_t(\mathit{Pr}(\site/S,\M(S)):a_*\, ,$$
%% and, for an object $X$ of $\site$ over $S$, we define
%% \begin{equation}\label{strictifyglobalsection5}
%% \derR\Gamma_{\mathit{geom},t}(X,M_t)=(a_*\, a^*(\derR\sigma_*\, \derL\sigma^*(M)))(X)\, .
%% \end{equation}
%% The object $\derR\Gamma_{\mathit{geom},t}(X,M_t)$ is the \emph{geometric
%% $t$-hypercohomology of $X$ over $S$ with coefficients in $M$}.
%% The presheaf $\derR\Gamma_{\mathit{geom},t}(-,M_t)=a_*\, a^*(\derR\sigma_*\, \derL\sigma^*(M))$
%% satisfies $t$-descent
%% by definition, and we have, by construction, a natural morphism
%% \begin{equation}\label{strictifyglobalsection6}
%% \derR\Gamma_{\mathit{geom}}(X,M)\To\derR\Gamma_{\mathit{geom},t}(X,M_t)\, .
%% \end{equation}
\end{paragr}

\begin{prop}\label{taudescentvscohomology}
For an object $M$ of $\ho(\M)(S)$, the following conditions
are equivalent.
\begin{itemize}
\item[(a)] The object $M$ satisfies $t$-descent.
\item[(b)] The presheaf $\derR\Gamma_{\mathit{geom}}(-,M)$ satisfies $t$-descent.
%%\item[(c)] The maps \eqref{strictifyglobalsection6} are isomorphisms for any $X$ over $S$.
\end{itemize}
\end{prop}

\begin{proof}
%%The equivalence between conditions (b) and (c) is given by Proposition \ref{taulocmodcat}.
For any morphism $f:X\To S$ and any
$t$-hypercover $p:\X\To X$ over $S$, 
we have, by Proposition \ref{computesimpldirectimage}
and  formula \eqref{strictifyglobalsection4}, an isomorphism
$$\derR f_* \derR p_*\, \derL p^*\, \derL f^*(M)\simeq
 \underset{\phantom{Ru}n\in \Delta}{\derR\varprojlim}\derR\Gamma_{\mathit{geom}} (\X_n,M)\, .$$
From there, we see easily that conditions (a) and (b) are equivalent.
\end{proof}

\begin{paragr}\label{derivators}
The preceding proposition allows us to reduce descent problems in a fibred
model category to descent problems in a category of presheaves
with values in a model category. On can even go further and
reduce the problem to category of presheaves with values in
an `elementary model category' as follows.

Consider a model category $\V$. Then one can associate to $\V$
its corresponding \emph{prederivator}\index{word}{prederivator}
 $\mathbf{Ho}(\V)$, that is the
strict $2$-functor
from the $2$-category of small categories to the $2$-category
of categories, defined by
\begin{equation}\label{defderV}
\mathbf{Ho}(\V)(I)=\ho(\V^{\op{I}})=\ho(\pshg{I,\V})
\end{equation}
for any small category $I$.
More explicitly: for any functor $u : I\To J$, one gets a functor
$$u^*: \mathbf{Ho}(\V) ( J ) \To \mathbf{Ho}(\V) ( I ) $$
(induced by the composition with $u$), and for any morphism of functors
$$\UseTwocells \xymatrix{ I\rrtwocell<4>^u_v{\,\alpha}&&J} \quad , $$
one has a morphism of functors
$$\UseTwocells \xymatrix{\mathbf{Ho}(\V)(I)&&\mathbf{Ho}(\V)(J)
         \lltwocell<4>^{v^*}_{u^*}{\alpha^* \ } }\quad . $$
Moreover, the prederivator $\mathbf{Ho}(\V)$ is then a
 Grothendieck derivator\index{word}{derivator, Grothendieck};
  see \cite[Thm. 6.11]{Cis1}. This means in particular that,
for any functor between small categories $u:I\To J$, the functor $u^*$
has a left adjoint
\begin{equation}\label{defderV2}
\derL u_\sharp:\mathbf{Ho}(\V)(I)\To \mathbf{Ho}(\V)(J)
\end{equation}
as well as a right adjoint
\begin{equation}\label{defderV3}
\derR u_*:\mathbf{Ho}(\V)(I)\To \mathbf{Ho}(\V)(J)
\end{equation}
(in the case where $J=e$ is the terminal category, then
$\derL u_\sharp$ is the homotopy colimit functor, while
$\derR u_*$ is the homotopy limit functor).

If $\V$ and $\V'$ are two model categories, a \emph{morphism of
derivators}\index{word}{morphism!of derivators}
$$\varPhi: \mathbf{Ho}(\V)\To  \mathbf{Ho}(\V')$$
is simply a morphism of $2$-functors, that is the data of functors
$$\varPhi_I : \mathbf{Ho}(\V)(I)\To  \mathbf{Ho}(\V')(I)$$
together with coherent isomorphisms
$$u^*(\varPhi_J(F))\simeq \varPhi_I (u^*(F))$$
for any functor $u:I\To J$ and any presheaf $F$ on $J$
with values in $\V$ (see \cite[p.~210]{Cis1} for
a precise definition).

Such a morphism $\varPhi$ is said to be
 \emph{continuous}\emph{morphism!continuous}
if, for any functor $u:I\To J$, and any object $F$ of $\mathbf{Ho}(\V)(I)$,
the canonical map
\begin{equation}\label{defderV4}
\varPhi_J\, \derR u_*(F)\To \derR u_*\, \varPhi_I(F)
\end{equation}
is an isomorphism. One can check that a morphism of derivators
$\varPhi$ is continuous if and only if it commutes with
homotopy limits (i.e. if and only if the maps $\eqref{defderV4}$
are isomorphisms in the case where $J=e$ is the terminal category);
see \cite[Prop.~2.6]{Cis2}. For instance, the total
right derived functor of any right Quillen functor
defines a continuous morphism of derivators; see \cite[Prop.~6.12]{Cis1}.

Dually a morphism $\varPhi$ of derivators is
 \emph{cocontinuous}\index{word}{morphism!cocontinuous}
if, for any functor $u:I\To J$, and any object $F$ of $\mathbf{Ho}(\V)(I)$,
the canonical map
\begin{equation}\label{defderV5}
 \derL u_!\, \varPhi_I(F)\To \varPhi_J\, \derL u_!(F)
\end{equation}
is an isomorphism. 
\end{paragr}

\begin{paragr}\label{defQlinear}
We shall say that a stable model category $\V$ is \emph{$\QQ$-linear}
 \index{word}{linear!$\QQ$-linear (stable model category)}
if all the objects of the triangulated category $\ho(\V)$
are uniquely divisible.
\end{paragr}

\begin{thm}\label{niceenrichment}
Let $\V$ be a model category (resp. a stable model category,
resp. a $\QQ$-linear stable model category), and denote
by $\mathcal{S}$ the model category of simplicial sets
(resp.  the stable model category of $S^1$-spectra, resp.
the $\QQ$-linear stable model category of complexes of
$\QQ$-vector spaces). Denote by $\unit$ the unit object
of the closed symmetric monoidal category $\ho(\mathcal{S})$.
Then, for each object $E$ of $\ho(\V)$, there exists a unique continuous
morphism of derivators
$$\derR\Hom(E,-):\mathbf{Ho}(\V)\To \mathbf{Ho}(\mathcal{S})$$
such that, for any object $F$ of $\ho(\V)$, there is a functorial bijection
$$\Hom_{\ho(\mathcal{S})}(\unit,\derR\Hom(E,F))\simeq\Hom_{\ho(\V)}(E,F))\, .$$
\end{thm}

\begin{proof}
Note that the stable $\QQ$-linear case follows from the stable case
and from the fact that the derivator of complexes of $\QQ$-vector
spaces is (equivalent to) the full subderivator of the derivator of $S^1$-spectra
spanned by uniquely divisible objects.

It thus remains to prove the theorem in the case where
$\V$ be a model category (resp. a stable model category)
and $\mathcal{S}$ is the model category of simplicial sets
(resp.  the stable model category of $S^1$-spectra).
The existence of $\derR\Hom(E,-)$ follows then from
\cite[Prop.~6.13]{Cis1} (resp. \cite[Lemma~A.6]{CisTab}).

For the unicity, as we don't really need it here, we shall only sketch the proof
(the case of simplicial sets is done in \cite[Rem.~6.14]{Cis1}).
One uses the universal property of the derivator 
$\mathbf{Ho}(\mathcal{S})$:
by virtue of \cite[Cor.~3.26]{Cis2} (resp. of \cite[Thm.~A.5]{CisTab}),
for any model category (resp. stable model category)
$\V'$ there is a canonical equivalence of categories
between the category of cocontinous morphisms from
$\mathbf{Ho}(\mathcal{S})$ to $\mathbf{Ho}(\V')$
and the homotopy category $\ho(\V)$. As a consequence,
the derivator $\mathbf{Ho}(\mathcal{S})$
admits a unique closed symmetric monoidal structure, and any
derivator (resp. triangulated derivator) is naturally and uniquely
enriched in $\mathbf{Ho}(\mathcal{S})$; see \cite[Thm.~5.22]{Cis2}.
More concretely, this universal property gives,
for any object $E$ in $\ho(\V')$, a unique
cocontinuous morphism of derivators
$$\mathbf{Ho}(\mathcal{S})\To \mathbf{Ho}(\V')\quad ,
\qquad K\mapsto K\otimes E$$
such that $\unit\otimes E=E$.
For a fixed $K$ in $\mathbf{Ho}(\mathcal{S})(I)$,
this defines a cocontinuous morphism of derivators
$$\mathbf{Ho}(\V')\To \mathbf{Ho}(\V^{\prime\op{I}})\quad ,
\qquad E\mapsto K\otimes E$$
which has a right adjoint
$$\mathbf{Ho}(\V^{\prime\op{I}})\To \mathbf{Ho}(\V')\quad ,
\qquad F\mapsto F^K\, .$$
Let
$$\derR\Hom(E,-):\mathbf{Ho}(\V)\To \mathbf{Ho}(\mathcal{S})$$
be a continuous morphism such that, for any object $F$ of $\V$, there is a functorial bijection
$$i_{F}:\Hom_{\ho(\mathcal{S})}(\unit,\derR\Hom(E,F))\simeq\Hom_{\ho(\V)}(E,F))\, .$$
Then, for any object $K$ of $\mathbf{Ho}(\mathcal{S})(I)$,
and any object $F$ of $\mathbf{Ho}(\V)(I)$
a canonical isomorphism
$$\derR\Hom(E,F^K)\simeq\derR\Hom(E,F)^K$$
which is completely determined by being the identity for $K=\unit$
(this requires the full universal property of $\mathbf{Ho}(\mathcal{S})$ given by
by \cite[Thm.~3.24]{Cis2} (resp. by the dual version of \cite[Thm.~A.5]{CisTab})).
We thus get from the functorial bijections $i_F$ the natural bijections:
\begin{align*}
\Hom_{\mathbf{Ho}(\mathcal{S})(I)}(K,\derR\Hom(E,F))
\simeq & \Hom_{\ho(\mathcal{S})}(\unit,\derR\Hom(E,F)^K)\\
\simeq & \Hom_{\ho(\mathcal{S})}(\unit,\derR\Hom(E,F^K))\\
\simeq & \Hom_{\ho(\V)}(E, F^K)\\
\simeq & \Hom_{\mathbf{Ho}(\V)(I)}(K\otimes E,F)\, .
\end{align*}
In other words, $\derR\Hom(E,-)$ has to be a right adjoint to
$(-)\otimes E$.
\end{proof}

\begin{rem}\label{remordsderivateurs}
The preceding theorem mostly holds for abstract derivators.
The only problem is for the existence of the morphism $\derR\Hom(E,-)$
(the unicity is always clear). However, this problem disappears for derivators which
have a Quillen model (as we have seen above),
as well as for triangulated derivators (see  \cite[Lemma~A.6]{CisTab}).
Hence Theorem \ref{niceenrichment} holds in fact for any triangulated
Grothendieck derivator.\index{word}{derivator, Grothendieck}

In the case when $\V$ is a combinatorial model category
(which, in practice, will essentially always be the case),
the enrichment over simplicial sets (resp, in the stable case, over
spectra) can be constructed via Quillen functors by
Dugger's presentation theorems \cite{dug0}
(resp. \cite{dug1}).
\end{rem}

\begin{cor}\label{reductionelementarymodcatbasic}
Let $\M$ be a $\Pmor$-fibred combinatorial model category
(resp. a stable $\Pmor$-fibred combinatorial model category,
resp. a $\QQ$-linear stable $\Pmor$-fibred combinatorial model category)
over $\site$, and $\mathcal{S}$ the model category of simplicial sets
(resp.  the stable model category of $S^1$-spectra, resp.
the $\QQ$-linear stable model category of complexes of
$\QQ$-vector spaces).

Consider an object $S$ of $\site$, a morphism
$f:X\To S$, and a morphism of $\site$-diagrams
$p:(\X,I)\To X$ over $S$. Then, for an object $M$ of $\ho(\M)(S)$,
the following conditions are equivalent.
\begin{itemize}
\item[(a)] The map
$$\derR f_*\, \derL f^*(M)\To \derR f_*\, \derR p_*\, \derL p^*\, \derL f^*(M)$$
is an isomorphism in $\ho(\M)(S)$.
\item[(b)] The map 
$$\derR\Gamma_{\mathit{geom}}(X,M)\To
\underset{\phantom{Ru}i\in \op{I}}{\derR\varprojlim}\derR\Gamma_{\mathit{geom}}(\X_i,M)$$
is an isomorphism in $\ho(\M)(S)$.
\item[(c)] For any object $E$ of $\ho(\M)(S)$, the map
$$\derR\Hom(E,\derR\Gamma_{\mathit{geom}}(X,M))\To
\underset{\phantom{Ru}i\in \op{I}}{\derR\varprojlim}
\derR\Hom(E,\derR\Gamma_{\mathit{geom}}(\X_i,M))$$
is an isomorphism in $\ho(\mathcal{S})$.
\end{itemize}
\end{cor}

\begin{proof}
The equivalence between (a) and (b) follows from Propositions
\ref{computeadjointsdiagrams2} and \ref{computeadjointsdiagrams3}, which
give the formula
$$\derR f_*\, \derR p_*\, \derL p^*\, \derL f^*(M)\simeq 
\underset{\phantom{Ru}i\in \op{I}}{\derR\varprojlim}\derR\Gamma_{\mathit{geom}}(\X_i,M)\, .$$
The identification
$$\Hom_{\ho(\mathcal{S})}(\unit,\derR\Hom(E,F))\simeq\Hom_{\ho(\M)(S)}(E,F)$$
and the Yoneda Lemma show that a map in $\ho(\M)(S)$ is an isomorphism
if and only its image by $\derR\Hom(E,-)$ is an isomorphism for any
object $E$ of $\ho(\M)(S)$.
Moreover, as $\derR\Hom(E,-)$ is continuous,
for any small category $I$ and any presheaf $F$ on $I$
with values in $\M(S)$, there is a canonical isomorphism
$$\derR\Hom(E,\underset{\phantom{Ru}i\in \op{I}}{\derR\varprojlim} \, F_i))\simeq
\underset{\phantom{Ru}i\in \op{I}}{\derR\varprojlim}\derR\Hom(E,F_i))\, .$$
This proves the equivalence between conditions (b) and (c).
\end{proof}

\begin{cor}\label{reductionelementarymodcat}
Under the assumptions of Corollary \ref{reductionelementarymodcatbasic},
given an object $S$ of $\site$, an object $M$ of $\ho(\M)(S)$
satisfies $t$-descent if and only if, for any object $E$ of $\ho(\M)(S)$
the presheaf of simplicial sets (resp. of $S^1$-spectra, resp. of
complexes of $\QQ$-vector spaces) $$\derR\Hom(E,\derR\Gamma_{\mathit{geom}}(-,M))$$
satisfies $t$-descent over $\site/S$.
\end{cor}

\begin{proof}
This follows from the preceding corollary, using the formula given
by Proposition \ref{computesimpldirectimage}.
\end{proof}

\begin{rem}
We need the category $\site$ to be small in some sense
to apply the two preceding corollaries because we need to make sense of the
projective model category structure of Proposition \ref{taulocmodcat}.
However, we can use these corollaries even if the site $\site$
is not small as well: we can either use the theory of universes, or
apply these corollaries to all the adequate small subsites of $\site$.
As a consequence, we shall feel free to use
Corollaries \ref{reductionelementarymodcatbasic} and
\ref{reductionelementarymodcat} for non necessarily small sites $\site$,
leaving to the reader the task to avoid set-theoretic difficulties
according to her/his taste.
\end{rem}

\begin{df}\label{defderivedabsoluteglobal sections}
For an $S^1$-spectrum $E$ and an integer $n$, we define its $n$th
cohomology group $H^n(E)$ by the formula
$$H^n(E)=\pi_{-n}(E)\, ,$$
where $\pi_i$ stands for the $i$th stable homotopy group functor.

Let $\M$ be a monoidal $\Pmor$-fibred
stable combinatorial model category over $\site$.
Given an object $S$ of $\site$ as well as an object $M$ of $\ho(\M)(S)$,
we define the presheaf of \emph{absolute derived global sections
\index{word}{section!absolute derived global section}
 of $M$ over $S$} by the formula
$$\derR\Gamma(-,M)=\derR\Hom(\unit_S,\derR\Gamma_{\mathit{geom}}(-,M))\, .$$
For a map $X\To S$ of $\site$, we thus have the \emph{absolute cohomology
of $X$ with coefficients in $M$}, $\derR\Gamma(X,M)$, as well as the
\emph{cohomology groups of $X$ with coefficients in $M$}:
$$H^n(X,M)=H^n(\derR\Gamma(X,M))\, .$$
We have canonical isomorphisms of abelian groups
\begin{align*}
H^n(X,M)
& \simeq \Hom_{\ho(\M)(S)}(\unit_S,\derR f_*\, \derL f^*(M))\\
& \simeq \Hom_{\ho(\M)(X)}(\unit_X, \derL f^*(M))\, .
\end{align*}

%% Similarly, one defines the \emph{$t$-hypercohomology of $X$ with coefficients
%% in $M$} as
%% $$\derR\Gamma_t(X,M_t)=\derR\Hom(\unit_S,\derR\Gamma_{\mathit{geom},t}(X,M_t))\, ,$$
%% as well as the \emph{$t$-hypercohomology groups of $X$ with coefficients in $M$}
%% $$H^n_t(X,M_t)=H^n(\derR\Gamma_t(X,M_t))\, .$$

Note that, if moreover $\M$ is $\QQ$-linear, the presheaf
$\derR\Gamma(-,M)$ %% and $\derR\Gamma_t(-,M_t)$
can be considered as a presheaf of complexes of $\QQ$-vector spaces on $\site/S$.
\end{df}
\subsection{Descent over schemes} \label{sec:hdescent}

The aim of this section is to give natural sufficient conditions for
$\M$ to satisfy descent with respect to various Grothendieck
 topologies\footnote{
In fact, using
remark \ref{remordsderivateurs}, all of this section
(results and proofs) holds for an abstract algebraic prederivator
 in the sense of Ayoub~\cite[Def. 2.4.13]{ayoub} without any changes
(note that the results of \ref{sec:fibredmodcat} are in fact a proof that (stable)
combinatorial fibred model categories over $\sch$ give rise
to algebraic prederivators). The only interest
of considering a fibred model category over $\sch$ is that it allows formulating
things in a little more naive way. Of course, the optimal setting in which to
formulate descent theory is the one of $\infty$-categories. However, restricting to
presentable $\infty$-categories, using Dugger's presentation theorem~\cite{dug0},
as well as rectification theorems such as \cite[Thm.~7.5.30 and 7.9.8]{HCHA}
as well as those from \cite{balzin}, we can see that the case of model categories
remains meaningful.}

\subsubsection{Localization and Nisnevich descent}

\begin{paragr}
Recall from example \ref{ex:lower&upper_cd_structures}
 that a \emph{Nisnevich distinguished square}
 \index{word}{square!Nisnevich distinguished}
  is a pullback square of schemes
\begin{equation}\label{Nisdist}
\begin{split}
\xymatrix{
V\ar[r]^l\ar[d]_g&Y\ar[d]^f\\
U\ar[r]_j&X}
\end{split}
\end{equation}
in which $f$ is \'etale, $j$ is an open immersion with reduced complement $Z$
 and the induced morphism $f^{-1}(Z)\To Z$ is an isomorphism.
\end{paragr}

For any scheme $X$ in $\sch$,
 we denote by $X_\nis$ the small Nisnevich site of $X$.
\begin{thm}[Morel-Voevodsky]\label{BGNis}
Let $\V$ be a (combinatorial) model category and $T$
 a scheme in $\sch$.
For a presheaf $F$ on $T_\nis$ with values in $\V$, the following conditions
are equivalent.
\begin{itemize}
\item[(i)] $F(\varnothing)$ is a terminal object in $\ho(\V)$, and
for any Nisnevich distinguished square \eqref{Nisdist} in $T_\nis$,
the square
$$\xymatrix{
F(X)\ar[r]\ar[d]&F(Y)\ar[d]\\
F(U)\ar[r]&F(V)
}$$
is a homotopy pullback square in $\V$.
\item[(ii)] The presheaf $F$ satisfies Nisnevich descent
\index{word}{descent!Nisnevich}
 on $T_\nis$.
\end{itemize}
\end{thm}
\begin{proof}
By virtue of corollaries \ref{reductionelementarymodcatbasic}
and \ref{reductionelementarymodcat},
it is sufficient to prove this in the case where $\V$ is the
usual model category of simplicial sets,
in which case this is precisely Morel and Voevodsky's theorem; see \cite{MV,voecd1,voecd2}.
\end{proof}

\num Consider a Nisnevich distinguished square \eqref{Nisdist}
 and put $a=jg=fl$. According to our general assumption \ref{num:assumption_hdescent_sch},
 the maps $a$, $j$ and $f$ are $\Pmor$-morphisms.
For any object $M$ of $\M(X)$, we obtain a commutative square in $\M$
(which is well-defined as an object in the homotopy of commutative squares in $\M(X)$):
\begin{equation}\label{localizationNisnevichdescent00}
\begin{split}
\xymatrix{
\derL a_\sharp a^*M \ar[r]\ar[d]
 & \derL f_\sharp f^*(M) \ar[d]\\
\derL j_\sharp j^*(M) \ar[r]
 & M. }
\end{split}
\end{equation}
We also obtain another commutative square in $\M$ by applying the functor $\derR\sHom_X(-,\unit_X)$:
\begin{equation}\label{localizationNisnevichdescent00_dual}
\begin{split}
\xymatrix{
M \ar[r]\ar[d]& \derR f_*\,  f^*(M) \ar[d]\\
\derR j_*\, j^*(M) \ar[r]&\derR a_*\, a^*(M). }
\end{split}
\end{equation}
\begin{prop}\label{localizationNisnevichdescent0}
If the category $\ho(\M)$
 has the localization property,
 then for any Nisnevich distinguished square \eqref{Nisdist}
 and any object $M$ of $\ho(\M)(X)$,
 the squares \eqref{localizationNisnevichdescent00} and
 \eqref{localizationNisnevichdescent00_dual} are homotopy cartesians.
 \index{word}{homotopy cartesian!square}
\end{prop}
\begin{proof}
Let $i:Z\To X$ be the complement of the open immersion $j$ ($Z$
being endowed with the reduced structure)
 and $p:f^{-1}(Z)\To Z$ the map induced by $f$.

We have only to prove that one of the squares \eqref{localizationNisnevichdescent00},
 \eqref{localizationNisnevichdescent00_dual} are cartesian.
We choose the square \eqref{localizationNisnevichdescent00}.

Because the pair of functor $(\derL i^*,j^*)$ is conservative on $\ho(\M)(X)$,
 we have only to check that the pullback of \eqref{localizationNisnevichdescent00}
 along $j^*$ or $\derL i^*$ is homotopy cartesian.
But, using the $\Pmor$-base change property, 
we see that the image of \eqref{localizationNisnevichdescent00} by $j^*$
is (canonically isomorphic to) the commutative square
$$
\xymatrix{
\derL g_\sharp a^*(M) \ar[d]\ar@{=}[r]& \derL g_\sharp a^*(M) \ar[d]\\
j^*(M) \ar@{=}[r] & j^*(M) }
$$
which is obviously homotopy cartesian.

Using again the $\Pmor$-base change property,
we obtain that the image of \eqref{localizationNisnevichdescent00} by $\derL i^*$
is isomorphic in $\ho(\M)$ to the square
$$
\xymatrix{
0 \ar@{=}[d]\ar[r] & p_\sharp p^* \derL i^*(M) \ar[d]\\
0 \ar[r] & \derL i^*(M) }
$$
which is again obviously homotopy cartesian because $p$ is an isomorphism
 (note for this last reason, $p_\sharp=\derL p_\sharp)$.
\end{proof}

\begin{cor}\label{localizationNisnevichdescent}
If $\ho(\M)$ has the localization property then it satisfies Nisnevich descent.
\index{word}{descent!Nisnevich}
\end{cor}

\begin{proof}
This corollary thus follows immediately from
Corollary \ref{reductionelementarymodcatbasic},
Theorem \ref{BGNis} and Proposition \ref{localizationNisnevichdescent0}.
\end{proof}

\begin{rem} \label{rem:MV_triangle}
Note that using Theorem \ref{BGNis},
 if we assume only that $\ho(\M)$ satisfies Nisnevich descent,
 then the squares \eqref{localizationNisnevichdescent00}
 and \eqref{localizationNisnevichdescent00_dual} are homotopy cartesian
 for any Nisnevich distinguished square \eqref{Nisdist}.

Assume that $\M$ is monoidal with geometric sections $M$.
 Let $S$ be a base scheme and consider a Nisnevich distinguished square \eqref{Nisdist}
 of smooth $S$-schemes.
 Then the fact that the square \eqref{localizationNisnevichdescent00} is homotopy cartesian
 implies there exists a \emph{canonical} distinguished triangle:
\begin{equation*}
M_S(V) \xrightarrow{g_*+l_*} M_S(U) \oplus M_S(Y) \xrightarrow{f_*+j_*} M_S(X)
 \longrightarrow M_S(V)[1]
\end{equation*}
It is called the
 \emph{Mayer-Vietoris triangle}\index{word}{triangle!Mayer-Vietoris triangle}
 associated with the square \eqref{Nisdist}.
\end{rem}

% \begin{cor}\label{Niscovercons}
% If $\M$ has the localization property, then, for any Nisnevich cover of $S$-schemes $f:Y\To X$,
% the functor
% $$\derL f^*:\ho(\M)(X)\To\ho(\M)(Y)$$
% is conservative.
% \end{cor}
% 
% \begin{proof}
% This follows from the preceding corollary and from Proposition \ref{coverconservative}.
% \end{proof}

\subsubsection{Proper base change isomorphism and descent by blow-ups}

\begin{paragr}\label{defcdhdistsquare}
Recall from example \ref{ex:lower&upper_cd_structures}
 that a
 \emph{$\cdh$-distinguished square}
\index{word}{square!cdhdistinguished@$\cdh$-distinguished}
 is a pullback square of schemes
\begin{equation}\label{cdhdist}
\begin{split}
\xymatrix{
T\ar[r]^k\ar[d]_g&Y\ar[d]^f\\
Z\ar[r]_i&X
}
\end{split}
\end{equation}
in which $f$ is proper surjective, $i$ a closed immersion
 and the induced map $f^{-1}(X-Z)\To X-Z$ is an isomorphism.

Recall from Example \ref{ex:lower&upper_cd_structures}
 the \emph{$\cdh$-topology} is the Grothendieck topology
  on the category of schemes generated by 
  Nisnevich coverings and by coverings of shape
$\{Z\To X,Y\To X\}$ for any $\cdh$-distinguished square \eqref{cdhdist}.
\end{paragr}

\begin{thm}[Voevodsky] Let $\V$ be a (combinatorial) model category.
For a presheaf $F$ on $\sch$ with values in $\V$, the following conditions
are equivalent.
\begin{itemize}
\item[(i)] The presheaf $F$ satisfies $\cdh$-descent
\index{word}{descent!cdhdescent@$\cdh$-descent}
 on $\sch$.
\item[(ii)] The presheaf $F$ satisfies Nisnevich descent and,
for any $\cdh$-distinguished square \eqref{cdhdist} of $\sch$,
the square
$$\xymatrix{
F(X)\ar[r]\ar[d]&F(Y)\ar[d]\\
F(Z)\ar[r]&F(T)
}$$
is a homotopy pullback\index{word}{homotopy pullback|see{homotopy cartesian}}
\index{word}{homotopy cartesian}
 square in $\V$.
\end{itemize}\label{BGcdh}
\end{thm}

\begin{proof}
It is sufficient to prove this in the case where $\V$ is the
usual model category of simplicial sets;
see  corollaries \ref{reductionelementarymodcatbasic}
and \ref{reductionelementarymodcat}.
As the distinguished $\cdh$-squares
define a bounded regular and reduced $cd$-structure on $\sch$, the equivalence between
(i) and (ii) follows from Voevodsky's theorems on
descent with respect to topologies defined by $cd$-structures \cite{voecd1,voecd2}.
\end{proof}

\num Consider a $\cdh$-distinguished square \eqref{cdhdist}
 and put $a=ig=fk$.
For any object $M$ of $\M(X)$, we obtain a commutative square in $\M$
(which is well-defined as an object in the homotopy of commutative squares in $\M(X)$):
\begin{equation}\label{cdhdescent00}\begin{split}
\xymatrix{
M \ar[r]\ar[d]& \derR f_*\,  \derL f^*(M) \ar[d]\\
\derR i_*\, \derL i^*(M) \ar[r]&\derR a_*\, \derL a^*(M) }
\end{split}\end{equation}
\begin{prop}\label{cdhdescent}
Assume $\ho(\M)$ satisfies the localization property
 and the transversality property with respect to proper morphisms.
Then the following conditions hold:
\begin{enumerate}
\item[(i)] For any $\cdh$-distinguished square \eqref{cdhdist},
 and any object $M$ of $\ho(\M)(X)$
 the commutative square \eqref{cdhdescent00} is homotopy cartesian.
\item[(ii)] The $\Pmor$-fibred model category $\ho(\M)$ satisfies $\cdh$-descent.
\index{word}{descent!cdhdescent@$\cdh$-descent}
\end{enumerate}
\end{prop}
\begin{proof}
We first prove (i). Consider a $\cdh$-distinguished square \eqref{cdhdist}
 and let $j:U\To X$ be the complement open immersion of $i$.
As the pair of functor $(\derL i^*,j^*)$ is conservative on $\ho(\M)(X)$,
 we have only to check that the image of \eqref{cdhdescent00} under $\derL i^*$
  and $j^*$ is homotopy cartesian.

Using projective transversality, we see that the image of \eqref{cdhdescent00}
by the functor $\derL i^*$ is (isomorphic to) the homotopy pullback square
$$\xymatrix{
&\derL i^* (M)\ar[r]\ar@{=}[d] & \derR g_* \, \derL g^* \, \derL i^* (M)\ar@{=}[d]& \\
&\derL i^* (M) \ar[r]& \derR g_* \, \derL g^* \, \derL i^* (M)&.
}$$

Let $h:f^{-1}(U) \rightarrow U$ be the pullback of $f$ over $U$.
As $j$ is an open immersion, it is by assumption a $\Pmor$-morphism
 and the $\Pmor$-base change formula implies 
 that the image of \eqref{cdhdescent00} by $j^*$
 is (isomorphic to) the commutative square
$$\xymatrix{
&\derL j^* (M)\ar[r] \ar[d] & \derR h_* \derL h^* \derL j^* (M)\ar[d]&\\
&0 \ar@{=}[r]& 0&
}$$
which is obviously homotopy cartesian because $h$ is an isomorphism.

We then prove (ii).
We already know that $\M$ satisfies Nisnevich descent (Corollary \ref{localizationNisnevichdescent}).
Thus, by virtue of the equivalence between conditions (i) and (ii)
of Theorem \ref{BGcdh}, the computation above,
together with corollaries \ref{reductionelementarymodcatbasic}
and \ref{reductionelementarymodcat}
imply that $\M$ satisfies $\cdh$-descent.
\end{proof}

\begin{paragr}\label{universalcdhdescent0}
To any $\cdh$-distinguished square \eqref{cdhdist}, one
associates a diagram of schemes $\Y$ over $X$ as follows.
Let $\cocoin$ be the category freely generated by the
oriented graph
\begin{equation}\label{universalcdhdescent000}\begin{split}
\xymatrix{
a\ar[r]\ar[d]&b\\
c&
}
\end{split}\end{equation}
Then $\Y$ is the functor from $\cocoin$ to $\sch/X$ defined by the
following diagram.
\begin{equation}\label{universalcdhdescent00}\begin{split}
\xymatrix{
T\ar[r]^k\ar[d]_g&Y \\
Z &}
\end{split}\end{equation}
We then have a canonical map $\varphi : \Y\To X$, and
the second assertion of Theorem \ref{cdhdescent}
can be reformulated by saying that the adjunction map
$$M\To \derR \varphi_*\, \derL \varphi^*(M)$$
is an isomorphism for any object $M$ of $\ho(\M)(X)$: indeed, by
virtue of Proposition \ref{computeadjointsdiagrams2}, 
$\derR \varphi_*\, \derL \varphi^*(M)$ is the homotopy limit of the diagram
$$\xymatrix{
& \derR f_*\,  \derL f^*(M) \ar[d]\\
\derR i_*\, \derL i^*(M) \ar[r]&\derR a_*\, \derL a^*(M) }$$
in $\ho(\M)(X)$. In other words,
if $\M$ has the properties of localization and of projective transversality, then
the functor
$$\derL\varphi^*:\ho(\M)(X)\To\ho(\M)(\Y,\cocoin)$$
is fully faithful.
\end{paragr}

\subsubsection{Proper descent with rational coefficients I: Galois excision}
\label{sec:proper_descent&Galois}

From now on, we assume that any scheme in $\sch$ is
quasi-excellent\footnote{See \ref{num:quasi-excellent}
below for a reminder on quasi-excellent schemes.}
(in fact, we shall only use the fact that the normalization of
a quasi-excellent schemes gives rise to a finite surjective morphism, so that,
in fact, universally japanese schemes would be enough).
We fix a scheme $S$ in $\sch$, and we shall work with $S$-schemes in $\sch$
(assuming these form an essentially small category). 

\begin{paragr}\label{defqfhandhtop}
The \emph{$\h$-topology}\index{word}{topology!htopology@$\h$-topology}
 (resp. the \emph{$\qfh$-topology)}\index{word}{topology!qfhtopology@$\qfh$-topology}
is the Grothendieck topology on the category of schemes
associated to the pretopology whose coverings
are the universal topological epimorphisms (resp.
the quasi-finite universal topological epimorphisms).
This topology has been introduced and studied by
Voevodsky in \cite{V1}.

The $\h$-topology is finer than the $\cdh$-topology
and, of course, finer than the $\qfh$-topology.
The $\qfh$-topology is in turn finer than the \'etale topology.
An interesting feature of the $\h$-topology (resp. of the $\qfh$-topology)
is that any proper (resp. finite) surjective
map is an $\h$-cover. In fact, the $\h$-topology
(resp. the $\qfh$-topology) can be
described as the topology generated by
the Nisnevich coverings and by the proper (resp. finite) surjective
maps; see Lemma \ref{refinedhvoverings}
(resp. Lemma \ref{qfhcoverings2}) below for a precise statement.
\end{paragr}

\num
Consider a morphism of schemes $f:Y \rightarrow X$. 
Consider the group of automorphisms $G=\mathrm{Aut}_Y(X)$ of the $X$-scheme $Y$.

Assuming $X$ is connected, we say according to \cite[exp. V]{SGA1} that
 $f$ is a \emph{Galois cover}\index{word}{cover!Galois cover}
  if it is finite \'etale (thus surjective)
 and $G$ operates transitively and faithfully on any (or simply one) of
 the geometric fibers of $Y/X$. 
Then $G$ is called the \emph{Galois group}\index{word}{group!Galois group}
\index{word}{Galois group|see{group}}
 of $Y/X$.\footnote{The map $f$ induces a one to one correspondence 
 between the generic points of $Y$
 and that of $X$. For any generic point $y \in Y$, $x=f(y)$, the residual extension
 $\kappa_y/\kappa_x$ is a Galois extension with Galois group $G$.}

When $X$ is not connected, we will still say that $f$ is a \emph{Galois cover}
 if it is so over any connected component of $X$. Then $G$ will be called
 the \emph{Galois group} of $X$.
 If $(X_i)_{i \in I}$ is the family connected
 components of $X$, then $G$ is the product of the Galois groups $G_i$
 of $f \times_X X_i$ for each $i \in I$. The group $G_i$ is equal to the Galois
 group of any residual extension over a generic point of $X_i$. 

The following definition is an extension of the definition 5.5 of \cite{SV1}:
\begin{df}\label{df:pseudo-galois}
A \emph{pseudo-Galois cover}\index{word}{cover!pseudo-Galois cover}
 is a finite surjective morphism of schemes $f:Y\To X$
 which can be factored as
$$
Y \xrightarrow{f'} X' \xrightarrow p X
$$
where $f'$ is a Galois cover and $p$ is radicial\footnote{
 See \ref{num:radicial} for a reminder on radicial morphisms.} (such a $p$ is automatically
finite and surjective).
\end{df}
Note that the group $G$ defined by the Galois cover $f'$
 is independent of the choice of the factorization.
 In fact, if $\bar X$ denotes the semi-localization of $X$ at its
 generic points, considering the cartesian squares
$$
\xymatrix@R=10pt@C=20pt{
\bar Y\ar[r]\ar[d]
 & \bar{X'}\ar[r]\ar[d] & \bar X\ar[d] \\
Y\ar^{f'}[r] & X'\ar^p[r] & X
}
$$
then $G=\mathrm{Aut}_{\bar X}(\bar Y)$ -- for any point $y \in \bar Y$,
 $x'=f'(y)$, $x=f(y)$, $\kappa_{x'}/\kappa_x$ is the maximal radicial sub-extension
 of the normal extension $\kappa_y/\kappa_x$.
It will be called the \emph{Galois group} of $Y/X$.

Remark also that $Y$ is a $G$-torsor over $X$ locally for the $\qfh$-topology
(i.e. it is a Galois object of group $G$ in the $\qfh$-topos of $X$):
this comes from the fact that finite radicial epimorphisms
are isomorphisms locally for the $\qfh$-topology
(any universal homeomorphism has this property by \cite[prop. 3.2.5]{V1}).

Let $f:Y\To X$ be a finite morphism, and $G$ a finite group
acting on $Y$ over $X$. Note that, as $Y$ is affine on $X$, the
scheme theoretic quotient $Y/G$ exists; see \cite[Exp.~V, Cor. 1.8]{SGA1}.
Such scheme-theoretic quotients are stable by flat pullbacks;
see  \cite[Exp.~V, Prop. 1.9]{SGA1}.

\begin{df}\label{settinghdesc}
Let $G$ be finite group. A \emph{$\qfh$-distinguished square
\index{word}{square!qfhdistinguished@$\protect\qfh$-distinguished}
 of group $G$}
is a pullback square of $S$-schemes of shape
\begin{equation}\label{qfhdist}
\begin{split}
\xymatrix{
T\ar[r]^h\ar[d]_g&Y\ar[d]^f \\
Z\ar[r]_i&X
}
\end{split}
\end{equation}
in which $Y$ is endowed with an action of $G$ over $X$, and satisfying the
following three conditions.
\begin{itemize}
\item[(a)] The morphism $f$ is finite and surjective.
\item[(b)] %%The closed subscheme $Z$ is nowhere dense, and
The induced morphism $f^{-1}(X-Z)\To f^{-1}(X-Z)/G$ is flat.
\item[(c)] The morphism $ f^{-1}(X-Z)/G\To X-Z$ is radicial.
\end{itemize}
\end{df}

Immediate examples of $\qfh$-distinguished squares of trivial group are the
following. The scheme $Y$ might be the normalization of $X$, and $Z$ is a
nowhere dense closed
subscheme out of which $f$ is an isomorphism; or
$Y$ is dense open subscheme of $X$ which is the disjoint
union of its irreducible components;
or $Y$ is a closed subscheme of $X$ inducing
an isomorphism $Y_{\mathit{red}}\simeq X_{\mathit{red}}$.

A $\qfh$-distinguished square of group $G$ \eqref{qfhdist}
will be said to be \emph{pseudo-Galois}
\index{word}{square!pseudo-Galois $\qfh$-distinguished}
\index{word}{pseudo-Galois|see{cover \emph{or} distinguished}}
 if $Z$ is nowhere dense in $X$ and if
the map $f^{-1}(X-Z)\To X-Z$ is a pseudo-Galois cover of group $G$.

The main examples of pseudo-Galois $\qfh$-distinguished squares
will come from the following situation.

\begin{prop}\label{genpsdgalhreg}
Consider an irreducible normal scheme $X$, and a finite extension $L$ of its field
of functions $k(X)$. Let $K$ be the inseparable closure of $k(X)$ in $L$, and
assume that $L/K$ is a Galois extension of group $G$.
Denote by $Y$ the normalization of $X$ in $L$.
Then the action of $G$ on $k(Y)=L$ extends naturally to an
action on $Y$ over $X$. Furthermore, there exists a closed
subscheme $Z$ of $X$, such that
the pullback square
$$\xymatrix{
T\ar[r]\ar[d]&Y\ar[d]^f\\
Z\ar[r]_i&X
}$$
is a pseudo-Galois $\qfh$-distinguished square of group $G$.
\end{prop}

\begin{proof}
The action of $G$ on $L$ extends naturally to an action on $Y$ over $X$
by functoriality. Furthermore, $Y/G$ is the normalization of $X$ in $K$,
so that $Y/G\To X$ is finite radicial and surjective
(see \cite[Lemma 3.1.7]{V1} or \cite[V, \textsection 2, n\textordmasculine~3, lem. 4]{BourbakiAC}).
By construction, $Y$ is generically a Galois cover over $Y/G$,
which implies the result (see \cite[Cor. 18.2.4]{EGA4}).
\end{proof}

\begin{paragr}
For a given $S$-scheme $T$, we shall denote by $L(T)$ the corresponding
representable $\qfh$-sheaf of sets (remember that the $\qfh$-topology is not subcanonical,
so that $L(T)$ has to be distinguished from $T$ itself). Beware that, in general, there
is no reason that, given a finite group $G$ acting on $T$, the scheme-theoretic
quotient $L(T/G)$ (whenever defined) and the $\qfh$-sheaf-theoretic quotient $L(T)/G$ would coincide.
\end{paragr}

\begin{lm}
\label{lemmerecqfhdiag}
Let $f:Y\To X$ be a separated morphism, $G$ a finite
group acting on $Y$ over $X$, and $Z$ a closed subscheme of $X$
such that $f$ is finite and surjective over $X-Z$,
and such that the quotient map $f^{-1}(X-Z)\To f^{-1}(X-Z)/G$
is flat, while the map $f^{-1}(X-Z)/G\To X-Z$ is radicial.
For $g\in G$, write $g:Y\To Y$
for the corresponding automorphism of $Y$, and define $Y_g$
as the image of the diagonal $Y\To Y\times_X Y$ composed
with the automorphism $1_Y\times_X g: Y\times_X Y\To Y\times_X Y$.
Then, if $T=Z\times_X Y$, we get a $\qfh$-cover of $Y\times_X Y$
by closed subschemes:
$$Y\times_X Y=(T\times_Z T)\cup \bigcup_{g\in G} Y_g \, .$$
\end{lm}

\begin{proof}
Note that, as $f$ is separated,
the diagonal $Y\To Y\times_X Y$ is a closed embedding,
so that the $Y_g$'s are closed subschemes of $Y\times_X Y$.
As the map $Y\times_{Y/G}Y\To Y\times_X Y$
is a universal homeomorphism, we may assume that $Y/G=X$.
It is sufficient to prove that, if $y$ and $y'$ are two geometric points of $Y$ whose images coincide
in $X$ and do not belong to $Z$, there exists an element $g$ of $G$ such that $y'=gy$
(which means that the pair $(y,y')$ belongs to $Y_g$).
For this purpose, we may assume, without loss of generality, that $Z=\varnothing$.
Then, by assumption, $Y$ is flat over $X$, from which we get the identification
$(Y\times_X Y)/G\simeq Y\times_X (Y/G)\simeq Y$
(where the action of $G$ on $Y\times_X Y$ is trivial on the first factor and
is induced by the action on $Y$ on the second factor). This achieves the proof.
\end{proof}

\begin{prop}\label{prepqfhexcision}
For any $\qfh$-distinguished square of group $G$ \eqref{qfhdist},
the commutative square
$$\xymatrix{L(T)/G\ar[r]\ar[d]&L(Y)/G\ar[d]\\
L(Z)\ar[r]&L(X)}$$
is a pullback and a pushout in the category of $\qfh$-sheaves.
Moreover, if $X$ is normal and if $Z$ is nowhere dense in $X$, then
the canonical map $L(Y)/G\To L(Y/G)\simeq L(X)$ is an isomorphism of $\qfh$-sheaves
(which implies that $L(T)/G\To L(Z)$ is an isomorphism as well).
\end{prop}

\begin{proof}
Note that this commutative square
is a pullback because it was so before taking the quotients by $G$
(as colimits are universal in any topos).
As $f$ is a $\qfh$-cover, it is sufficient to prove that
$$\xymatrix{L(T)\times_{L(Z)}L(T)/G\ar[r]\ar[d]&L(Y)\times_{L(X)}L(Y)/G\ar[d]\\
L(T)\ar[r]&L(Y)}$$
is a pushout square. This latter square fits into the following
commutative diagram
$$\xymatrix{
L(T)\ar[r]\ar[d]&L(Y)\ar[d]\\
L(T)\times_{L(Z)}L(T)/G\ar[r]\ar[d]&L(Y)\times_{L(X)}L(Y)/G\ar[d]\\
L(T)\ar[r]&L(Y)}$$
in which the two vertical composed maps are identities
(the vertical maps of the upper commutative square are obtained
from the diagonals by taking the quotients under the natural action
of $G$ on the right component).
It is thus sufficient to prove that the upper square is a pushout.
As the lower square is a pullback, the upper one shares the same property;
moreover, all the maps in the upper commutative square are monomorphisms
of $\qfh$-sheaves, so that it is sufficient to prove that the
map $(L(T)\times_{L(Z)}L(T)/G)\amalg L(Y)\To L(Y)\times_{L(X)}L(Y)/G$ is an epimorphism
of $\qfh$-sheaves. According to Lemma \ref{lemmerecqfhdiag},
this follows from the commutativity of the diagram
$$\xymatrix{
L(T\times_Z T)\amalg \Big(\coprod_{g\in G} L(Y_g)\Big) \ar[r]\ar[d]& L(Y\times_X Y)\ar[d]\\
(L(T)\times_{L(Z)}L(T)/G)\amalg L(Y)\ar[r]& L(Y)\times_{L(X)}L(Y)/G
}$$
in which the vertical maps are obviously epimorphic.

Assume now that $X$ is normal and that $Z$ is nowhere dense in $X$,
and let us prove that the canonical map $L(Y)/G\To L(X)$ is an isomorphism of $\qfh$-sheaves.
This is equivalent to prove that, for any $\qfh$-sheaf of sets $F$,
the map $f^*:F(X)\To F(Y)$ induces a bijection
$$F(X)\simeq F(Y)^G\, .$$
Let $F$ be a $\qfh$-sheaf.
The map $f^*:F(X)\To F(Y)$ is injective because $f$ is a $\qfh$-cover,
and it is clear that the image of $f^*$ lies in $F(Y)^G$.

Let $a$ be a section of $F$ over $Y$
which is invariant under the action of $G$.
Denote by $\mathit{pr}_1,\mathit{pr}_2:Y\times_XY\To Y$ the two
canonical projections. With the notations introduced in Lemma \ref{lemmerecqfhdiag},
we have
$$\mathit{pr}_1^*(a)|_{Y_g}=a=a.g=\mathit{pr}_2^*(a)|_{Y_g}$$
for every element $g$ in $G$.
As $Z$ does not contain any generic point of $X$,
the scheme $T\times_Z T$ does not contain
any generic point of $Y\times_X Y$ neither: as any irreducible
component of $Y$ dominates an irreducible component of $X$, 
and, as $X$ is normal, the finite map $Y\To X$ is universally open;
in particular, the projection $\mathit{pr}_1:Y\times_XY\To Y$
is universally open, which implies that any generic point of $Y\times_X Y$
lies over a generic point of $Y$.
By virtue of \cite[prop. 3.1.4]{V1}, Lemma \ref{lemmerecqfhdiag}
thus gives a $\qfh$-cover of $Y\times_X Y$ by closed subschemes
of shape
$$Y\times_X Y=\bigcup_{g\in G}Y_g\, .$$
This implies that
$$\mathit{pr}_1^*(a)=\mathit{pr}_2^*(a)\, .$$
The morphism $Y \To X$ being a $\qfh$-cover
and $F$ a $\qfh$-sheaf, we deduce that the section $a$ lies in the image of $f^*$.
\end{proof}

\begin{cor}\label{preqfhdescente}
For any $\qfh$-distinguished square of group $G$ \eqref{qfhdist},
we get a bicartesian square of $\qfh$-sheaves of abelian groups
$$\xymatrix{\ZZ_\qfh(T)_G\ar[r]\ar[d]&\ZZ_\qfh(Y)_G\ar[d]\\
\ZZ_\qfh(Z)\ar[r]&\ZZ_\qfh(X)}$$
(where the subscript $G$ stands for the coinvariants under the action of $G$).
In other words, there is a canonical short exact sequence of sheaves of abelian groups
$$0\To\ZZ_\qfh(T)_G\To\ZZ_\qfh(Z)\oplus\ZZ_\qfh(Y)_G\To\ZZ_\qfh(X)\To 0\, .$$
\end{cor}

\begin{proof}
As the abelianization functor preserves colimits and monomorphisms,
the preceding proposition implies formally that we have
a short exact sequence of shape
$$\ZZ_\qfh(T)_G\To\ZZ_\qfh(Z)\oplus\ZZ_\qfh(Y)_G\To\ZZ_\qfh(X)\To 0\, ,$$
while the left exactness follows from the fact that $Z\To X$
being a monomorphism, the map obtained by pullback,
$L(T)/G\To L(Y)/G$, is a monomorphism as well.
\end{proof}

\begin{paragr}
Let $\V$ be a $\QQ$-linear stable model category (see \ref{defQlinear}).

Consider a finite group $G$, and an object $E$ of $\V$, endowed with
an action of $G$. By viewing $G$ as a category with one object we can see
$E$ as functor from $G$ to $\V$ and take its homotopy limit
in $\ho(\V)$, which we denote by $E^{hG}$ (in the literature,
$E^{hG}$ is called the \emph{object of homotopy fixed points}
\index{word}{homotopy!object of homotopy fixed points}
under the action of $G$ on $E$). One the other hand, the category $\ho(\V)$
is, by assumption, a $\QQ$-linear triangulated category with small sums,
and, in particular, a $\QQ$-linear pseudo-abelian category so that we can
define $E^G$ as the object of $\ho(\V)$ defined by
\begin{equation}\label{defQlinGinv1}
E^G=\mathrm{Im}\, p\, ,
\end{equation}
where $p:E\To E$ is the projector defined in $\ho(\V)$
by the formula
\begin{equation}\label{defQlinGinv2}
p(x)=\frac{1}{\# G}\sum_{g\in G} g.x\, .
\end{equation}
The inclusion $E^G\To E$ induces a canonical isomorphism
\begin{equation}\label{defQlinGinv3}
E^{G}\overset{\sim}{\To} E^{hG}
\end{equation}
in $\ho(\V)$: to see this, by virtue of Theorem \ref{niceenrichment}, we can
assume that $\V$ is the model category of complexes of $\QQ$-vector spaces,
in which case it is obvious.
\end{paragr}

\begin{cor}\label{derivedqfhgenericGaloisdescent}
Let $C$ be a presheaf of complexes of $\QQ$-vector spaces on the
category of $S$-schemes. Then, for any $\qfh$-distinguished square of group $G$ \eqref{qfhdist},
the commutative square
$$\xymatrix{
\derR\Gamma_\qfh(X,C_\qfh)\ar[r]\ar[d]&\derR\Gamma_\qfh(Y,C_\qfh)^G\ar[d]\\
\derR\Gamma_\qfh(Z,C_\qfh)\ar[r]&\derR\Gamma_\qfh(T,C_\qfh)^G
}$$
is a homotopy pullback\index{word}{homotopy cartesian}
 square in the derived category of $\QQ$-vector spaces.
In particular, we get a long exact sequence of shape
$$\cdots\To H^n_\qfh(X,C_\qfh)\To H^n_\qfh(Z,C_\qfh) \oplus H^n_\qfh(Y,C_\qfh)^G
\To H^n_\qfh(T,C_\qfh)^G \To\cdots$$
If furthermore $X$ is normal and $Z$ is nowhere dense in $X$, then the maps
$$H^n_\qfh(X,C_\qfh)\To H^n_\qfh(Y,C_\qfh)^G\quad
\text{and}\quad H^n_\qfh(Z,C_\qfh)\To H^n_\qfh(T,C_\qfh)^G$$
are isomorphisms for any integer $n$.
\end{cor}

\begin{proof}
Let $C_\qfh \To C'$ be a fibrant resolution
in the $\qfh$-local injective model category structure
on the category of $\qfh$-sheaves of complexes of $\QQ$-vector spaces;
see for instance \cite[Cor.~4.4.42]{ayoub}.
Then for $U=Y,T$, we have a natural isomorphism of complexes
$$\Hom(\QQ_\qfh(U)_G,C')=C'(U)^G$$
which gives an isomorphism
$$\derR\Hom(\QQ_\qfh(U)_G,C_\qfh)\simeq \derR\Gamma_\qfh(U,C_\qfh)^G$$
in the derived category of the abelian category of $\QQ$-vector spaces.
This corollary thus follows formally from Corollary \ref{preqfhdescente}
by evaluating at the derived functor $\derR\Hom(-,C_\qfh)$.

If furthermore $X$ is normal,
then one deduces the isomorphism $H^n_\qfh(X,C_\qfh)\simeq H^n_\qfh(Y,C_\qfh)^G$
from the fact that $L(Y)/G\simeq L(Y/G)\simeq X$
(Proposition \ref{prepqfhexcision}), which implies that
$\ZZ_\qfh(Y)_G\simeq\ZZ_\qfh(X)$.
The isomorphism $H^n_\qfh(Z,C_\qfh)\simeq H^n_\qfh(T,C_\qfh)^G$
then comes as a byproduct of the long exact sequence above.
\end{proof}

\begin{thm}\label{ratdescent}
Let $X$ be a scheme, and $C$ be a presheaf of complexes of $\QQ$-vector
spaces on the small \'etale site of $X$. Then $C$ satisfies \'etale descent
\index{word}{descent!etale@\'etale}
if and only if it has the following properties.
\begin{itemize}
\item[(a)] The complex $C$ satisfies Nisnevich descent.
\item[(b)] For any \'etale $X$-scheme $U$ and any
Galois cover $V\To U$ of group $G$, the map
$C(U)\To C(V)^G$ is a quasi-isomorphism.
\end{itemize}
\end{thm}

\begin{proof}
These are certainly necessary conditions.
To prove that they are sufficient, note that the Nisnevich cohomological dimension and the
rational \'etale cohomological dimension of a noetherian scheme are bounded by the
dimension; see \cite[proposition 1.8, page 98]{MV} and \cite[Lemma 3.4.7]{V1}.
By virtue of \cite[Theorem 0.3]{SV2},
for $\tau=\nis,\et$, we have strongly convergent spectral sequences
$$E^{p,q}_2=H^p_\tau(U,H^q(C)_\tau)\Rightarrow H^{p+q}_\tau(U,C_\tau)\, .$$
Condition (a) gives isomorphisms $H^{p+q}(C(U))\simeq H^{p+q}_\nis(U,C_\nis)$,
so that it is sufficient to prove that,
for each of the cohomology presheaves $F=H^q(C)$, we have
$$H^p_\nis(U,F_\nis)\simeq H^p_\et(U,F_\et)\, .$$
As the rational \'etale cohomology
of any henselian scheme is trivial in non-zero degrees, it is sufficient to prove that,
for any local henselian scheme $U$ (obtained as the henselisation of an
\'etale $X$-scheme at some point), $F_\nis(U)\simeq F_\et(U)$.
Let $G$ be the absolute Galois group of the closed point of $U$. Then we have
$$F_\nis(U)=F(U)\quad \text{and}  \quad F_\et(U)=\varinjlim_\alpha F(U_\alpha)^{G_\alpha}\, ,$$
where the $U_\alpha$'s run over all the Galois covers of $U$ corresponding
to the finite quotients $G\To G_\alpha$. But it follows from (b) that
$F(U)\simeq F(U_\alpha)^{G_\alpha}$ for any $\alpha$, so that
$F_\nis(U)\simeq F_\et(U)$.
\end{proof}

\begin{lm}\label{qfhcoverings}
Any $\qfh$-cover\index{word}{cover!qfhcover@$\qfh$-cover}
 admits a refinement of the form $Z\To Y\To X$,
where $Z\To Y$ is a finite surjective morphism, and $Y\To X$ is an
\'etale cover.
\end{lm}

\begin{proof}
This property being clearly local on $X$ with respect to the
\'etale topology, we can assume that $X$ is strictly henselian,
in which case this follows from \cite[Lemma 3.4.2]{V1}.
\end{proof}

\begin{thm}\label{carratqfhdescent}
A presheaf of complexes of $\QQ$-vector spaces $C$ on the category of $S$-schemes
satisfies $\qfh$-descent\index{word}{descent!qfhdescent@$\qfh$-descent}
 if and only if it has the following two properties:
\begin{itemize}
\item[(a)] the complex $C$ satisfies Nisnevich descent;
\item[(b)] for any pseudo-Galois $\qfh$-distinguished square of group $G$ \eqref{qfhdist},
the commutative square
$$\xymatrix{
C(X)\ar[r]\ar[d]&C(Y)^G\ar[d]\\
C(Z)\ar[r]&C(T)^G
}$$
is a homotopy pullback square\index{word}{homotopy cartesian}
 in the derived category of $\QQ$-vector spaces.
\end{itemize}
\end{thm}

\begin{proof}
Any complex of presheaves of $\QQ$-vector spaces satisfying $\qfh$-descent
satisfies properties (a) and (b):
property (a) follows from the fact that the $\qfh$-topology is finer
than the \'etale topology;
property (b) is Corollary \ref{derivedqfhgenericGaloisdescent}.

Assume now that $C$ satisfies these two
properties. Let $\varphi: C\To C'$ be a morphism of presheaves of complexes of $\QQ$-vector spaces
which is a quasi-isomorphism locally for the $\qfh$-topology, and such that
$C'$ satisfies $\qfh$-descent (such a morphism exists thanks to the
$\qfh$-local model category structure on the category of presheaves of
complexes of $\QQ$-vector spaces; see Proposition \ref{taulocmodcat}).
Then the cone of $\varphi$ also satisfies
conditions (a) and (b). Hence it is sufficient to prove the theorem in the case where $C$
is acyclic locally for the $\qfh$-topology.

Assume from now on that $C_\qfh$ is an acyclic complex of $\qfh$-sheaves, and
denote by $H^n(C)$ the $n\text{th}$ cohomology presheaf associated to $C$.
We know that the associated $\qfh$-sheaves vanish, and we want to deduce that
$H^n(C)=0$.

We shall prove by induction on $d$ that, for any $S$-scheme $X$ of dimension $d$
and for any integer $n$, the group $H^n(C)(X)=H^n(C(X))$ vanishes.
The case where $d<0$ follows from the fact, that by (a), the presheaves $H^n(C)$
send finite sums to finite direct sums, so that, in particular, $H^n(C)(\varnothing)=0$.
Before going further,
notice that condition (b) implies $H^n(C)(X_\mathit{red})=H^n(C)(X)$ for any $S$-scheme $X$
(consider the case where, in the diagram \eqref{qfhdist},
$Z=Y=T=X_\mathit{red}$), so that it is always
harmless to replace $X$ by its reduction. 
Assume now that $d\geq 0$, and
that the vanishing of $H^n(C)(X)$ is known whenever $X$ is of dimension $<d$ 
and for any integer $n$. Under this inductive assumption,
we have the following reduction principle.

Consider a pseudo-Galois $\qfh$-distinguished square of group $G$  \eqref{qfhdist}.
If $Z$ and $T$ are of dimension $<d$, then by condition (b), the map
$H^n(C)(X)\To H^n(C)(Y)^G$ is an isomorphism: indeed, we have an exact
sequence of shape
$$H^{n-1}(C)(T)^G\To H^n(C)(X)\To H^n(C)(Z)\oplus H^n(C)(Y)^G\To H^n(C)(T)^G\, ,$$
which implies our assertion by induction on $d$.

We shall prove now the vanishing of $H^n(C)(T)$
for normal $S$-schemes $T$ of dimension $d$.
Let $a$ be a section of $H^n(C)$ over such a $T$.
As $H^n(C)_\qfh(T)=0$, there exists a $\qfh$-cover $g:Y\To T$
such that $g^*(a)=0$. But, by virtue of Lemma \ref{qfhcoverings},
we can assume $g$ is the composition of a finite surjective morphism
$f:Y\To X$ and of an \'etale cover $e:X\To T$.
We claim that $e^*(a)=0$. To prove it, as, by (a), the presheaf
$H^n(C)$ sends finite  sums to finite direct sums, we can
assume that $X$ is normal and connected. Refining $f$ further, we can
assume that $Y$ is the normalization of $X$ in a finite extension of $k(X)$,
and that $k(Y)$ is a Galois extension of group $G$ over the inseparable closure of $k(X)$
in $k(Y)$. By virtue of Proposition \ref{genpsdgalhreg}, we get by
the reduction principle the identification
$H^n(C)(X)=H^n(C)(Y)^G$, whence $e^*(a)=0$.
As a consequence, the restriction of the presheaf of complexes $C$
to the category of normal $S$-schemes of dimension $\leq d$
is acyclic locally for the \'etale topology (note that this is
quite meaningful, as any \'etale scheme over a normal
scheme is normal; see \cite[Prop. 18.10.7]{EGA4}).
But $C$ satisfies \'etale descent\index{word}{descent!etale@\'etale}
 (by virtue of Theorem \ref{ratdescent}
this follows formally from property (a) and from property (b) for $Z=\varnothing$), so that
$H^n(C)(T)=H^n_\et(T,C_\et)=0$ for any normal $S$-scheme $T$ of dimension $\leq d$ and any
integer $n$.

Consider now a reduced $S$-scheme $X$ of dimension $\leq d$.
Let $p:T\To X$ be the normalization of $X$. As $p$ is
birational (see \cite[Cor.~6.3.8]{EGA2}) and finite surjective
(because $X$ is quasi-excellent), we can apply the reduction principle
and see that the pullback map $p^*:H^n(C)(X)\To H^n(C)(T)=0$
is an isomorphism for any integer $n$,
which achieves the induction and the proof.
\end{proof}

\begin{lm}\label{etcovering}
\'Etale coverings are finite \'etale coverings locally
for the Nisnevich topology:
any \'etale cover admits a refinement of the form
$Z\To Y\To X$, where $Z\To Y$ is a finite \'etale cover and $Y\To X$
is a Nisnevich cover.
\end{lm}

\begin{proof}
This property being local on $X$ for the Nisnevich
topology, it is sufficient to prove this in the case where $X$
is local henselian. Then, by virtue of \cite[Cor. 18.5.12 and Prop. 18.5.15]{EGA4},
we can even assume that $X$ is the spectrum of field, in which
case this is obvious.
\end{proof}

\begin{lm}\label{qfhcoverings2}
Any $\qfh$-cover\index{word}{cover!qfhcover@$\qfh$-cover}
 admits a refinement of the form
$Z\To Y\To X$, where $Z\To Y$ is a finite surjective morphism, and $Y\To X$
is a Nisnevich cover.
\end{lm}

\begin{proof}
As finite surjective morphisms are stable by pullback and
composition, this follows immediately from
lemmata \ref{qfhcoverings} and \ref{etcovering}.
\end{proof}

\begin{lm}\label{refinedhvoverings}
Any $\h$-cover\index{word}{cover!hcover@$\h$-cover}
 of an integral scheme $X$ admits a refinement of the form
$$U\To Z\To Y \To X\, ,$$
where $U\To Z$ is a finite surjective morphism,
$Z\To Y$ is a Nisnevich cover,
$Y\To X$ is a proper surjective birational map, and $Y$
is normal.
\end{lm}

\begin{proof}
By virtue of \cite[Theorem 3.1.9]{V1}, any $\h$-cover
admits a refinement of shape
$$W \To V \To X\, ,$$
where $W\To V$ is a $\qfh$-cover, and $V\To X$ is a proper surjective
birational map. By replacing $V$ by its normalization $Y$,
we get a refinement of shape
$$W\times_V Y\To Y \To X$$
where $W\times_V Y\To Y$ is a $\qfh$-cover, and $Y\To X$ is proper
surjective birational map.
We conclude by Lemma \ref{qfhcoverings2}.
\end{proof}

\begin{lm}\label{perphdesc}
Let $C$ be a presheaf of complexes of $\QQ$-vector spaces on the
category of $S$-schemes satisfying $\qfh$-descent.
Then, for any finite surjective morphism
$f:Y\To X$ with $X$ normal, the map
$f^*:H^n(C)(X)\To H^n(C)(Y)$ is a monomorphism.
\end{lm}

\begin{proof}
It is clearly sufficient to prove this when $X$ is connected.
Then, up to refinement, we can assume that
$f$ is a map as in Proposition \ref{genpsdgalhreg}. In this case, by virtue
of Corollary \ref{derivedqfhgenericGaloisdescent}, the $\QQ$-vector
space $H^n(C)(X)\simeq H^n(C)(Y)^G$ is a direct factor of
$H^n(C)(Y)$.
\end{proof}

\begin{thm}\label{hdesceqqfhcdhdesc}
A presheaf of complexes of $\QQ$-vector spaces
on the category of $S$-schemes satisfies
 $\h$-descent\index{word}{descent!hdescent@$\h$-descent}
 if and only if it satisfies
 $\qfh$-descent\index{word}{descent!qfhdescent@$\qfh$-descent}
 and $\cdh$-descent.\index{word}{descent!cdhdescent@$\cdh$-descent}
\end{thm}

\begin{proof}
This is certainly a necessary condition, as the $\h$-topology
is finer than the $\qfh$-topology and the $\cdh$-topology.
For the converse, as in the proof of Theorem \ref{carratqfhdescent},
it is sufficient to prove that any presheaf of complexes of
$\QQ$-vector spaces $C$ on the
category of $S$-schemes satisfying $\qfh$-descent
and $\cdh$-descent, and which is acyclic locally for the $\h$-topology,
is acyclic. We shall prove by noetherian induction that,
given such a complex $C$, for
any integer $n$, and any $S$-scheme $X$, for any section $a$ of $H^n(C)$
over $X$, there exists a $\cdh$-cover $X'\To X$ on which $a$ vanishes.
In other words, we shall get that $C$ is acyclic locally for the $\cdh$-topology,
and, as $C$ satisfies $\cdh$-descent, this will imply that
$H^n(C)(X)=H^n_\cdh(X,C_\cdh)=0$
for any integer $n$ and any $S$-scheme $X$.
Note that the presheaves $H^n(C)$ send finite
sums to finite direct sums (which follows, for instance, from the fact
that $C$ satisfies Nisnevich descent). In particular, $H^n(C)(\varnothing)=0$
for any integer $n$.

Let $X$ be an $S$-scheme, and $a\in H^n(C)(X)$. We have
a $\cdh$-cover of $X$ of shape $X'\amalg X''\To X$,
where $X'$ is the sum of the irreducible components of $X_{\mathit{red}}$
and $X''$ is a nowhere dense closed subscheme of $X$,
so that we can assume $X$ is integral.
Let $a$ be a section of the presheaf $H^n(C)$ over $X$.
As $H^n(C)_h=0$, by virtue of Lemma \ref{refinedhvoverings},
there exists a proper surjective birational map $p:Y\To X$
with $Y$ normal, a Nisnevich cover $q:Z\To Y$, and a surjective
finite morphism $r:U\To Z$ such that $r^*(q^*(p^*(a)))=0$ in $H^n(C)(U)$.
But then, $Z$ is normal as well (see \cite[Prop. 18.10.7]{EGA4}),
so that, by Lemma \ref{perphdesc}, we have
$q^*(p^*(a))=0$ in $H^n(C)(Z)$.
Let $T$ be a nowhere dense closed subscheme of $X$ such that $p$
is an isomorphism over $X-T$. By noetherian induction,
there exists a $\cdh$-cover $T'\To T$ such that $a|_{T'}$ vanishes. 
Hence the section $a$ vanishes on 
the $\cdh$-cover $T'\amalg Z\To X$.
\end{proof}
%% \begin{paragr}
%% An \emph{alteration} is a proper surjective morphism
%% $f:Y\To X$ which is generically finite. The \emph{center}
%% of an alteration $f$ is the reduced closed subscheme $Z$
%% such that $U=X-Z$ is the maximal open subscheme of $X$
%% such that the map $f^{-1}(U)\To U$ is finite.
%% \end{paragr}
%% 
%% \begin{paragr}
%% Let $f:Y\To X$ be an alteration between integral schemes, and suppose
%% furthermore that $Y$ is normal.
%% Remark that, there exists a finite extension $K$ of the field
%% of functions of $Y$ such that, if $Y'$ denotes the
%%\end{paragr}

\subsubsection{Proper descent with rational coefficients II: separation}
\label{sec:proper_descent&separation}
From now on, we assume that $\ho(\M)$ is $\QQ$-linear.

\begin{prop}\label{computeschemequotientaction}
Let $f:Y\To X$ be a morphism of schemes in $\sch$, and $G$ a finite group
acting on $Y$ over $X$. Denote by $\Y$ the scheme $Y$
considered a functor from $G$ to the category of $S$-schemes, and
denote by $\varphi:(\Y,G)\To X$ the morphism induced by $f$.
Then, for any object $M$ of $\ho(\M)(X)$, there are canonical
isomorphisms
$$(\derR f_*\, \derL f^*(M))^G\simeq(\derR f_*\, \derL f^*(M))^{hG}
\simeq \derR\varphi_*\, \derL \varphi^*(M)\, .$$
(where $G$ acts by functoriality of the construction
$\derR f_*\, \derL f^*$, as expressed by formulas
\eqref{strictifyglobalsection3} and \eqref{strictifyglobalsection4}).
\end{prop}

\begin{proof}
The second isomorphism comes from Proposition \ref{computeadjointsdiagrams2},
and the first, from \eqref{defQlinGinv3}.
\end{proof}

\begin{thm}\label{fibredetaledescent}
If $\ho(\M)$ satisfies Nisnevich descent, the following conditions are equivalent:
\begin{itemize}
\item[(i)] $\ho(\M)$ satisfies \'etale descent.
\item[(ii)] for any finite \'etale cover $f:Y\To X$, the functor
$$\derL f^*:\ho(\M)(X)\To\ho(\M)(Y)$$
is conservative;
\item[(iii)] for any finite Galois cover $f:Y\To X$ of group $G$, and for any
object $M$ of $\ho(\M)(X)$, the canonical map
$$M\To(\derR f_*\, \derL f^*(M))^G$$
is an isomorphism.
\end{itemize}
\end{thm}

\begin{proof}
The equivalence between (i) and (iii) follows from
Theorem \ref{ratdescent} by corollaries \ref{reductionelementarymodcatbasic}
and \ref{reductionelementarymodcat}, and Proposition \ref{coverconservative}
shows that (i) implies (ii). It is thus sufficient to prove that
(ii) implies (iii).
Let $f:Y\To X$ be a finite Galois cover of group $G$.
As the functor $f^*=\derL f^*$ is conservative
by assumption, it is sufficient to check that the map
$M\To(\derR f_*\, \derL f^*(M))^G$
becomes an isomorphism after applying $f^*$.
By virtue of Proposition \ref{basicPbasechange},
this just means that it is sufficient to prove (iii) when
$f$ has a section, i.e. when $Y$ is
isomorphic to the trivial $G$-torsor over $X$.
In this case, we have the (equivariant) identification
$ \bigoplus_{g\in G}M \simeq \derR f_*\, \derL f^*(M)$,
where $G$ acts on the left term by permuting the factors.
Hence $M\simeq (\derR f_*\, \derL f^*(M))^G$.
\end{proof}

\begin{prop}\label{prop:charseparated1}
Assume that $\ho(\M)$ has the localization property.
The following conditions are equivalent:
\begin{itemize}
\item[(i)] $\ho(\M)$ is separated.
\item[(ii)] $\ho(\M)$ is semi-separated and satisfies \'etale descent.
\end{itemize}
\end{prop}
\begin{proof}
This follows from Proposition \ref{prop:redseptoetale} and
Theorem \ref{fibredetaledescent}.
\end{proof}

\begin{cor}\label{charseparated2}
Assume that all the residue fields of $\base$ are of characteristic zero,
and that $\ho(\M)$ has the property of localization.
Then the following conditions are equivalent:
\begin{itemize}
\item[(i)] $\ho(\M)$ is separated.
\item[(ii)] $\ho(\M)$ satisfies \'etale descent.
\end{itemize}
\end{cor}

\begin{proof}
Consider a radicial finite surjective morphism $f:Y\To X$ in $\sch$.
To prove that the functor $\derL f^*$ is conservative, as
$\ho(\M)$ has the property of localization, by noetherian induction,
we may replace $X$ by any dense open subscheme $U$ (and $Y$ by
$U\times_X Y$). The residue fields of $X$ being of characteristic zero,
this means that we may assume that $f$ induces an isomorphism
after reduction $Y_{\mathit{red}}\simeq X_{\mathit{red}}$.
But it is clear that, by the localization property, such a morphism $f$
induces an equivalence of categories $\derL f^*$, so that $\ho(\M)$
is automatically semi-separated.
We conclude by Proposition \ref{prop:charseparated1}.
\end{proof}

\begin{prop}\label{Galoisdescseparated}
Assume that $\ho(\M)$ is separated, satisfies the localization property
 the proper transversality property.
Then, for any pseudo-Galois cover $f:Y\To X$ of group $G$, and for any
object $M$ of $\ho(\M)(X)$, the canonical map
$$M\To(\derR f_*\, \derL f^*(M))^G$$
is an isomorphism.
\end{prop}

\begin{proof}
By Proposition \ref{prop:charseparated1},
this is an easy consequence of Proposition \ref{stronglyquasiseparated}
and of condition (iii) of Theorem \ref{fibredetaledescent}.
\end{proof}

\begin{paragr}
From now on, we assume furthermore that any scheme in $\sch$
is quasi-excellent.
\end{paragr}

\begin{thm}\label{charseparated3}
Assume that $\ho(\M)$ satisfies the localization and proper transversality properties.
Then the following conditions are equivalent:
\begin{itemize}
\item[(i)] $\ho(\M)$ is separated;
\item[(ii)] $\ho(\M)$ satisfies
 $\h$-descent;\index{word}{descent!hdescent@$\h$-descent}
\item[(iii)] $\ho(\M)$ satisfies
 $\qfh$-descent;\index{word}{descent!qfhdescent@$\qfh$-descent}
\item[(iv)] for any  $\qfh$-distinguished
\index{word}{square!qfhdistinguished@$\qfh $-distinguished}
square \eqref{qfhdist} of group $G$, if we write $a=fh=ig:T\To X$ for
the composed map, then, for any object $M$ of $\ho(\M)(X)$,
the commutative square
\begin{equation}\label{qfhdescent00}\begin{split}
\xymatrix{
M \ar[r]\ar[d]& (\derR f_*\,  \derL f^*(M))^G \ar[d]\\
\derR i_*\, \derL i^*(M) \ar[r]&(\derR a_*\, \derL a^*(M))^G }
\end{split}\end{equation}
is homotopy cartesian;\index{word}{homotopy cartesian}
\item[(v)] the same as condition (iv), but only for
pseudo-Galois $\qfh$-distinguished squares.
\end{itemize}
\end{thm}

\begin{proof}
As $\M$ satisfies $\cdh$-descent (Theorem \ref{cdhdescent}),
the equivalence between conditions (ii) and (iii) follows
from Theorem \ref{hdesceqqfhcdhdesc} by
Corollary \ref{reductionelementarymodcat}.
Similarly, Theorem \ref{carratqfhdescent}
and corollaries \ref{derivedqfhgenericGaloisdescent},
\ref{reductionelementarymodcatbasic}
and \ref{reductionelementarymodcat} show that
conditions (iii), (iv) and (v) are equivalent.
As \'etale surjective morphisms as well as
finite radicial epimorphisms are $\qfh$-coverings,
it follows from Proposition \ref{coverconservative},
Theorem \ref{fibredetaledescent} and Proposition \ref{prop:charseparated1}, that
condition (iii) implies condition (i).
It thus remains to prove that condition (i) implies condition (v).
So let us consider a pseudo-Galois $\qfh$-distinguished square \eqref{qfhdist}
of group $G$, and prove that \eqref{qfhdescent00} is homotopy cartesian.
Using proper transversality, we see that the image of \eqref{qfhdescent00}
by the functor $\derL i^*$ is (isomorphic to) the homotopy pullback square
$$\xymatrix{
&\derL i^* (M)\ar[r]\ar@{=}[d] & (\derR g_* \, \derL g^* \, \derL i^* (M))^G\ar@{=}[d]& \\
&\derL i^* (M) \ar[r]& (\derR g_* \, \derL g^* \, \derL i^* (M))^G&.
}$$
Write $j:U\To X$ for the complement open immersion of $i$,
and $b:f^{-1}(U)\To U$ for the map induced by $f$.
As $j$ is \'etale, we see, using Proposition \ref{basicPbasechange},
that the image of \eqref{cdhdescent00} by $j^*=\derL j^*$
is (isomorphic to) the square
$$\xymatrix{
& j^* (M)\ar[r] \ar[d] & (\derR b_*\, \derL b^*\,  j^* (M))^G\ar[d]&\\
&0 \ar@{=}[r]& 0&.
}$$
in which the upper horizontal map is an isomorphism by Proposition \ref{Galoisdescseparated}.
Hence it is a homotopy pullback square.
Thus, because the pair of functors $(\derL i^*,j^*)$ is conservative
 on $\ho(\M)(X)$, the square \eqref{qfhdescent00} is homotopy cartesian.
\end{proof}

\begin{cor}\label{charseparated4}
Assume that all the residue fields of $\base$ are of characteristic zero,
and that $\ho(\M)$ has the localization and proper transversality properties.
Then $\ho(\M)$ satisfies $\h$-descent\index{word}{descent!hdescent@$\h$-descent}
 if and only if it satisfies
  \'etale descent.\index{word}{descent!etale@\'etale}
\end{cor}

\begin{proof}
This follows from Corollary \ref{charseparated2} and Theorem \ref{charseparated3}.
\end{proof}

\begin{cor}\label{stronggenericGaloisdescent}
Assume that $\ho(\M)$ is separated
 and has the localization and proper transversality properties.
Let $f:Y\To X$ be a finite surjective morphism, with $X$ normal, and $G$ a group
acting on $Y$ over $X$, such that the map $Y/G\To X$ is generically radicial
(i.e. radicial over a dense open subscheme of $X$).
Consider at last a pullback square of the following shape.
$$\xymatrix{
Y'\ar[r]\ar[d]_{f'}& Y\ar[d]^f\\
X'\ar[r]& X
}$$
Then, for any object $M$ of $\ho(\M)(X')$, the natural
map
$$M\To (\derR f'_*\, \derL f^{\prime*}(M))^G$$
is an isomorphism.
\end{cor}

\begin{proof}
For any presheaf $C$ of complexes of $\QQ$-vector spaces on $\sch/X$,
one has an isomorphism
$$\derR\Gamma_\qfh(X',C_\qfh)\simeq \derR\Gamma_\qfh(Y',C_\qfh)^G\, .$$
This follows from the fact that we have an isomorphism
of $\qfh$-sheaves of sets $L(Y)/G\simeq L(X)$ (the map $Y\To Y/G$ being generically
flat, this is Proposition \ref{prepqfhexcision}),
which implies that the map $L(Y')/G\To L(X')$ is an isomorphism
of $\qfh$-sheaves (by the universality of colimits in topoi), and implies
this assertion (as in the proof of \ref{derivedqfhgenericGaloisdescent}).

By virtue of Theorem \ref{charseparated3}, $\ho(\M)$ satisfies $\qfh$-descent,
so that the preceding computations imply the result by
corollaries \ref{reductionelementarymodcatbasic}
and \ref{reductionelementarymodcat}.
\end{proof}

\begin{cor}\label{stronggenericGaloisdescentcor}
Assume that $\ho(\M)$ is separated
 and has the localization and proper transversality properties.
Then for any finite surjective morphism $f:Y\To X$ with $X$ normal, the morphism
$$M\To \derR f_*\, \derL f^*(M)$$
is a monomorphism and admits a functorial splitting in $\ho(\M)(X)$.
Furthermore, this remains true
after base change by any map $X'\To X$.
\end{cor}

\begin{proof}
It is sufficient to treat the case where $X$ is connected.
We may replace $Y$ by a normalization of $X$ in a suitable finite extension
of its field of function, and assume that a finite group $G$ acts on $Y$ over $X$,
so that the properties described in the preceding corollary
are fulfilled (see \ref{genpsdgalhreg}).
\end{proof}

\begin{rem}
The condition (iv) of Theorem \ref{charseparated3} can be reformulated
in a more global way as follows (this won't be used in these notes, but
this might be useful for the reader who might want to formulate all this
in terms of (pre-)algebraic derivators \cite[Def.~2.4.13]{ayoub}).
Given a $\qfh$-distinguished square \eqref{qfhdist}
of group $G$, we can form a functor $\F$ from
category $I=\cocoin$ \eqref{universalcdhdescent000}
to the category of diagrams of $S$-schemes corresponding to
the diagram of diagrams of $S$-schemes
\begin{equation*}%%\label{universalqfhdescent00}
\begin{split}
\xymatrix{
(\T,G)\ar[r]^{(h,1_G)}\ar[d]_{g}&(\Y,G) \\
Z &}
\end{split}\end{equation*}
in which $\T$ and $\Y$ correspond to $T$ anf $Y$
respectively, seen as functor from $G$ to $\sch/X$.
The construction of \ref{diagdiag} gives a diagram of $X$-schemes
$(\Int \F,I_\F)$ which can be described explicitly as follows.
The category $I_\F$ is the cofibred category over $\cocoin$
associated to the functor from $\cocoin$ to the category of small
categories defined by the diagram
\begin{equation*}%%\label{universalqfhdescent01}
\begin{split}
\xymatrix{
G\ar[r]^{1_G}\ar[d]&G \\
e &}
\end{split}\end{equation*}
in which $e$ stands for the terminal category, and $G$ for the
category with one object associated to $G$.
It has thus three objects $a,b,c$ (see \eqref{universalcdhdescent000}),
and the morphisms are determined by
$$\Hom_{I_\F}(x,y)=
\begin{cases}
$*$ &\text{if $y=c$;}\\
\varnothing &\text{if $x\neq y$ and $x=b,c$;}\\
$G$&\text{otherwise.}
\end{cases}$$
The functor $\F$ sends $a,b,c$ to $T,Y,Z$ respectively, and
simply encodes the fact that the diagram
$$\xymatrix{
T\ar[r]^{h}\ar[d]_{g}&Y \\
Z &}$$
is $G$-equivariant, the action on $Z$ being trivial. Now, by propositions \ref{directimagediagdiag}
and \ref{computeschemequotientaction},
if $\varphi: (\F,I_\F)\To (X,\cocoin)$ denotes the canonical map, for any object
$M$ of $\ho(\M)(X)$, the object $\derR\varphi_*\, \derL \varphi^*(M)$
is the functor from $\coin=\op{\cocoin}$ to $\M(X)$ corresponding to the diagram
below (of course, this is well defined only in the homotopy category of
the category of functors from $\coin$ to $\M(X)$).
\begin{equation*} \begin{split}
\xymatrix{
& (\derR f_*\,  \derL f^*(M))^G \ar[d]\\
\derR i_*\, \derL i^*(M) \ar[r]&(\derR a_*\, \derL a^*(M))^G }
\end{split}\end{equation*}
As a consequence, if $\psi:(\Int\F,I_\F)\To X$ denotes the structural map,
the object $\derR \psi_*\, \derL \psi^*(M)$ is simply the 
homotopy limit of the diagram of $\M(X)$ above, so that
condition (iv) of Theorem \ref{charseparated3} can now be reformulated by saying
that the map
$$M\To \derR \psi_*\, \derL \psi^*(M)$$
is an isomorphism, i.e. that the functor
$$\derL \psi^*: \ho(\M)(X)\To\ho(\M)(\Int \F,I_\F)$$
is fully faithful.
\end{rem}
\section{Constructible motives} \label{sec:constructible_motives}

\begin{assumption} \label{assumption_constructible_motives}
Consider as in \ref{num:assumption1_sch}
 a base scheme $\base$ and a sub-category $\sch$
  of the category of $\base$-schemes.
 In section \ref{sec:motivic_duality}, and for the main theorem
 of section  \ref{sec:motivic_finiteness}, we will assume:
\begin{itemize}
\item[(a)] Any scheme in $\sch$ is quasi-excellent.\footnote{See
 Paragraph \ref{num:quasi-excellent}. The reader can safely restrict his
 attention to the more classical notion of an excellent scheme
 (\cite[IV, 7.8.5]{EGA4}).}
\end{itemize}
Apart in Definition \ref{df:continuous} and the subsequent proposition,
 where we will consider an abstract situation,
 we will be concerned with the study of a fixed premotivic triangulated category $\T$
 over $\sch$ (recall Definition \ref{df:motivic_cat}) such that:
\begin{itemize}
\item[(b)] $\T$ is motivic (see Definition \ref{df:motivic_cat}).
\item[(c)] $\T$ is endowed with a set of twists $\tau$
  (see Paragraph \ref{num:premotivic&twists})
   which is stable under Tate twists $\un(p)[q]$, for $p,q\in\ZZ$.
\item[(d)] $\T$ is the homotopy category associated with
 a stable combinatorial $\sm$-fibred model category $\M$ over $\sch$.\footnote{
 We need this assumption to apply descent theory
 as described in section \ref{sec:hdescent}.}
\end{itemize}
As usual, the geometric sections of $\T$ will be denoted by $M$.

Unless explicitly referring to the underlying model category $\M$,
 we will not indicate in the notation of the six operations 
 that the functors are derived functors.
\end{assumption}

\subsection{Resolution of singularities}

The aim of this subsection is to gather the results
 from the theory of resolution of singularities
 \index{word}{resolution of singularities}
  that will be used subsequently.

\begin{num} \label{num:quasi-excellent}
In \cite[IV, 7.8.2]{EGA4}, Grothendieck defined the notion
 of an \emph{excellent ring}.
 Matsumura introduced in \cite{Mat} the weaker notion of
  a \emph{quasi-excellent}\index{word}{quasi-excellent} ring $A$.
 Recall $A$ is quasi-excellent if the following conditions hold:
\begin{enumerate}
\item $A$ is noetherian.
\item For any prime ideal $\mathfrak p$, $\hat A_\mathfrak p$
 being the completion of $A$ at $\mathfrak p$,
 the canonical morphism $A \rightarrow \hat A_{\mathfrak p}$ is regular
 (see \ref{num:regular_morphism} below).
\item For any $A$-algebra $B$ of finite type,
 the regular locus of $\spec B$ is open.
\end{enumerate}
Then a ring $A$ is excellent if
 it is quasi-excellent and universally catenary.
Following Gabber, we say a scheme $X$ is \emph{quasi-excellent}
\index{word}{scheme!quasi-excellent}
(\emph{excellent})
\index{word}{scheme!excellent}
 if it admits an open cover by affine schemes 
 whose rings are quasi-excellent (excellent, respectively).
\end{num}

\begin{thm}[Gabber's weak local uniformization] Let $X$ be a
quasi-excellent scheme, and $Z\subset X$
a nowhere dense closed subscheme. Then there exists a finite $\h$-cover
\index{word}{cover!hcover@$\h$-cover}
$\{f_i:Y_i\To X\}_{i\in I}$
such that for all $i$ in $I$, $f_i$ is a morphism of finite type, the scheme
$Y_i$ is regular, and $f^{-1}_i(Z)$ is either empty or the support of a strict normal crossing divisor in $Y_i$.\label{Gabber1}
\end{thm}
\noindent See \cite{gabber3} for a proof.
Note that, if we are only interested in schemes of finite type over $\spec R$,
for $R$ either a field, a complete discrete valuation ring, or a Dedekind domain
whose field of functions is a global field, this is an immediate consequence of
de~Jong's resolution of singularities by alterations; see \cite{dejong}.
One can also deduce the case of schemes of finite type over
an excellent noetherian scheme of dimension lesser or equal to $2$
from \cite{dejong2}; see Theorem \ref{thm:resuptoquotient}
and Corollary \ref{cor:dejongdimleq2} below for a precise statement.

\begin{rem}
This theorem will be used in the proof of
 Lemma \ref{gabbergeomlemmathmfinitness}
 which is the key point for the proof of Theorem \ref{thmfinitness0}.
\end{rem}

\begin{num} \label{num:regular_morphism}
Recall that a morphism of rings $u:A\To B$ is \emph{regular} if
it is flat, and if, for any prime ideal $\mathfrak{p}$ in $A$, with residue field $\kappa(\mathfrak{p})$,
the $\kappa(\mathfrak{p})$-algebra $\kappa(\mathfrak{p})\otimes_A B$ is geometrically regular
(equivalently, this means that, for any prime ideal $\mathfrak{q}$ of $B$, the
$A$-algebra $B_{\mathfrak{q}}$ is formally smooth for the $\mathfrak{q}$-adic topology).
We recall the following great generalization of Neron's desingularization theorem:
\end{num}
\begin{thm}[Popescu-Spivakovsky] A morphism of noetherian
rings $u:A\To B$ is regular if and only if $B$ is a filtered colimit of
smooth $A$-algebras of finite type.\label{popescu}
\end{thm}
\begin{proof}
See \cite[theorems 1.1 and 1.2]{spivak}.
\end{proof}

\begin{paragr}
Recall that an \emph{alteration}\index{word}{alteration} is
 a proper surjective morphism
$p:X'\To X$ which is generically finite, i.e. such that there exists
a dense open subscheme $U\subset X$ over which $p$ is finite.
\end{paragr}

\begin{df}[de Jong] Let $X$ be a noetherian scheme
endowed with an action of a finite group $G$.
A \emph{Galois alteration}\index{word}{alteration!Galois alteration} of
the couple $(X,G)$ is the data of a finite group $G'$, 
of a surjective morphism
of groups $G'\To G$, of an alteration $X'\To X$, 
and of an action of $G'$ on $X'$, such that:
\begin{itemize}
\item[(i)] the map $X'\To X$ is $G'$-equivariant;
\item[(ii)] for any irreducible component $T$ of $X$, there exists
a unique irreducible component $T'$ of $X'$ over $T$, and
the corresponding finite field extension
$$k(T)^G\subset k(T')^{G'}$$
is purely inseparable.
\end{itemize}
In practice, we shall keep the morphism of groups $G'\To G$ implicit,
and we shall say that $(X'\To X,G')$ is a Galois alteration of $(X,G)$.

Given a noetherian scheme $X$, a \emph{Galois alteration} of $X$ is
a Galois alteration $(X'\To X,G)$ of $(X,e)$, where $e$ denotes the trivial group.
In this case, we shall say that $X'\To X$ is a \emph{Galois alteration of $X$ of group $G$}.\label{df:Galoisalteration}
\end{df}

\begin{rem}
If $p:X'\To X$ is a Galois alteration of group $G$ over $X$, then, if $X$
and $X'$ are normal, irreducible and quasi-excellent, $p$ can be factored as
a radicial finite surjective morphism $X''\To X$, followed by a Galois alteration
$X'\To X''$ of group $G$, such that $k(X'')=k(X')^G$
(just define $X''$ as the normalization of $X$ in $k(X')^G$).
In other words, up to a radicial finite surjective morphism, $X$ is generically the quotient of $X'$
under the action of $G$.
\end{rem}

\begin{df}\label{df:ressingupquotient}
A noetherian scheme $S$ \emph{admits canonical dominant resolution of singularities up to quotient singularities} 
\index{word}{resolution of singularities!canonical dominant ---- up to quotient singularities}
if, for any Galois alteration $S'\To S$ of group $G$, and for any
$G$-equivariant nowhere dense closed subscheme $Z'\subset S'$, there exists
a Galois alteration $(p:S''\To S',G')$ of $(S',G)$, such that $S''$ is regular
and projective over $S$, and such that the inverse image of $Z'$ in $S''$ is contained in a $G'$-equivariant
strict normal crossing divisor (i.e. a strict normal crossing divisor
whose irreducible components are stable under the action of $G'$).

A noetherian scheme $S$ admits \emph{canonical resolution
 of singularities up to quotient singularities}
\index{word}{resolution of singularities!canonical ---- up to quotient singularities}
  if any integral closed subscheme of $S$
admits canonical dominant resolution of singularities up to quotient
singularities.

A noetherian scheme $S$ admits \emph{wide resolution of singularities
 up to quotient singularities}
\index{word}{resolution of singularities!wide ---- up to quotient singularities}
  if, for any separated $S$-scheme of finite type $X$, and any nowhere
dense closed subscheme $Z\subset X$, there exists a projective Galois
alteration $p:X'\To X$ of group $G$, with $X'$ regular,
such that, in each connected component of $X'$,
$Z'=p^{-1}(Z)$ is either empty or the support of a strict normal crossing divisor.
\end{df}

\begin{thm}[de Jong] If an excellent noetherian
scheme of finite dimension $S$
%%which has an ample family of line bundles
admits canonical resolution of singularities up to quotient singularities, then any separated
$S$-scheme of finite type admits canonical resolution of singularities up to quotient singularities.\label{thm:resuptoquotient}
%% In particular, for any separated $S$-scheme of finite type $Y$,
%% for any dominant proper morphism $X\To Y$, and any nowhere
%% dense closed subscheme $Z\subset X$, there exists a Galois
%% alteration $p:X'\To X$ of group $G$, with $X'$ regular and projective over $Y$,
%% such that, in each connected component of $X'$,
%% $Z'=p^{-1}(Z)$ is either empty, either the support of a strict normal crossing divisor.
\end{thm}

\begin{proof}
Let $X$ be an integral separated $S$-scheme of finite type.
There exists a finite morphism $S'\To S$, with $S'$ integral,
an integral dominant $S'$-scheme $X'$
and a radicial extension $X'\To X$ over $S$,
such that $X'$ has a geometrically irreducible generic fiber over $S'$.
It follows then from (the proof of) \cite[theorem 5.13]{dejong2}
that $X'$ admits canonical dominant resolution of singularities up to quotient singularities,
which implies that $X$ has the same property.
%% It remains to check the second assertion.
%% Let $Y$ be a separated $S$-scheme of finite type, $X\To Y$
%% a proper morphism, and $Z$ a nowhere dense closed subscheme in $X$.
%% Replacing $X$ by the disjoint union of the reduction of its irreducible
%% component, we see that it is sufficient to consider the case where $X$
%% is integral. By Chow's lemma, we may assume that $X$ is projective over $Y$,
%% so that it is harmless to assume that $X=Y$. If $Z$ is not empty, by blowing up $Z$ in $X$,
%% we may also assume that $Z$ a divisor in $X$.
%% Then, by de~Jong's theorem, there exists a projective
%% Galois alteration $p:X'\To X$ of group $G$, with $X'$ regular,
%% such that $Z'=p^{-1}(Z)$ is contained in the support of a strict normal crossing divisor in $X'$.
%% In this situation, if $Z$ is a divisor, then $Z'$ has to be a divisor itself, which implies
%% it is the support of a strict normal crossing divisor.
\end{proof}
%% 
%% \begin{thm}[Lipman~{\cite{lipman}}]
%% Any excellent noetherian scheme of dimension lesser or equal to $2$ admits canonical
%% resolution of singularities (hence, in particular, canonical dominant resolution of singularities up to
%% quotient singularities).
%% \end{thm}

\begin{cor}[de Jong] Let $S$ be an excellent
noetherian scheme of dimension lesser or equal to $2$.
%%which has an ample family of line bundles.
Then any separated scheme of finite type over $S$
admits canonical resolution of singularities up to quotient
singularities. In particular, $S$ admits wide resolution of singularities up to quotient singularities.\nocite{lipman}\label{cor:dejongdimleq2}
\end{cor}

\begin{proof}
See \cite[corollary 5.15]{dejong2}.
\end{proof}

\subsection{Finiteness theorems} \label{sec:motivic_finiteness}
%\label{generalassumptionfiniteness}
%We fix a quasi-excellent noetherian scheme of finite dimension $\base$
%as well as an adequate category of $\base$-schemes $\sch$ (\ref{num:assumption1_sch}).
%We consider given a stable combinatorial $\sm$-fibred model
%category $\M$ over $\sch$, such that $\T$ is a motivic category over $\sch$,
%with a generating set of twists $\tau$, which is assumed to be stable under negative Tate twists.
%As usual, the geometric section of $\T$ will be denoted by $M$.

The aim of this section is to study the notion of $\tau$-constructibility
\index{word}{constructible!$\tau$-constructible|(}
\index{word}{constructibility|see{constructible}}
 in the triangulated motivic case and to study its stability properties under
 Grothendieck six operations. Recall the following particular case
  of Definition \ref{df:tau-geometric}:
\begin{df}\label{deftauconst}
For a scheme $X$ in $\sch$,
 we denote by $\T_{c}(X)$ the thick triangulated sub-category
 of $\T(X)$ generated by premotives of the form
 $M_X(Y)\{i\}$ for a smooth $X$-scheme $Y$ and a twist $i \in \tau$.
 We will say that a premotive in $\T_{c}(X)$ is \emph{$\tau$-constructible},
 or, simply, \emph{constructible}.
%
%When $\tau=\ZZ$ is the group generated by the Tate twist $\un(1)$,
% we say simply constructible and denote by $\T_c$ the corresponding category.
\end{df}

\begin{rem}\renewcommand{\Rc}{\Lambda}
Let us mention that our main examples:
\begin{itemize}
\item the stable homotopy category $\SH$
 (cf. Example \ref{ex:H&SH_premotivic}),
\item the category of Voevodsky motives $\DMV$
 (cf. Definition \ref{df:Nis_DMe&DM}),
\item the category of Beilinson motives $\DMB$
 (cf. Definition \ref{df:Beilinson_motives})
\end{itemize}
are all generated by the Tate twists (\emph{i.e.} $\tau=\ZZ$).
Recall also Proposition \ref{constructequivcompact}: it applies to all
 these examples so that constructible premotives coincides with compact
  objects.\footnote{Notice
 however this fact is not true for \'etale motivic complexes.}
\end{rem}

%\begin{prop}\label{constructequivcompact}
%Assume that, for any scheme $X$, the triangulated category
%$\T(X)$ admits finite sums and that,
%for any smooth $X$ scheme $Y$ and any $n\in \tau$, the
%object $M_X(Y)\{n\}$ is compact. Then, an object of $\T(X)$
%is constructible if and only if it is compact.
%\end{prop}
%
%\begin{proof}
%If $\T$ is any compactly generated triangulated category,
%then, for any small family $C$ of compact generators, the
%thick triangulated category of $\T$ generated by $C$ consists exactly
%of the compact objects of $\T$.
%\end{proof}
%\begin{ex} \renewcommand{\Rc}{\Lambda}
%Assume $X$ is noetherian of finite dimension.
%The preceding proposition can be applied to the category $\SH(X)$.
% It will also be applied to the categories:
% $\DMtx X$, $\DM(X)$, $\DMB(X)$.
%\end{ex}

\begin{prop}\label{soritestauconstruct1}
If $M$ and $N$ are constructible in $\T(X)$, so is
$M\otimes_X N$.
\end{prop}

\begin{proof}
For a fixed $M$, the full subcategory of $\T(X)$ spanned by objects
such that $M\otimes_X N$ is constructible is a thick
triangulated subcategory of $\T(X)$.
In the case $M$ is of shape $M_X(Y)\{n\}$ for $Y$ smooth over $X$
and $n\in \tau$, this proves that $M\otimes_X N$ is constructible
whenever $N$ is. By the same argument, using the symmetry of the
tensor product, we get to the general case.
\end{proof}

Similarly, one has the following conservation property.

\begin{prop}\label{soritestauconstruct2}
For any morphism $f:X\To Y$ of schemes, the functor
$$f^*:\T(Y)\To\T(X)$$
preserves constructible objects.
If moreover $f$ is smooth, the functor
$$f_\sharp: \T(X)\To \T(Y)$$
also preserves constructible objects.
\end{prop}

\begin{cor}
The categories $\T_{c}(X)$ form a thick triangulated monoidal
$\sm$-fibred subcategory of $\T$.
\end{cor}

\begin{prop}\label{thmfinitnessredlocZar0}
Let $X$ a scheme, and $X=\bigcup_{i\in I}U_i$ a cover of $X$
by open subschemes. An object $M$ of $\T(X)$ is
constructible if and only if its restriction to $U_i$
is constructible in $\T(U_i)$ for all $i\in I$.
\end{prop}

\begin{proof}
This is a necessary condition by \ref{soritestauconstruct2}.
For the converse, as $X$ is noetherian, it is sufficient to treat
the case where $I$ is finite. Proceeding by induction on the cardinal
of $I$ it is sufficient to treat the case of a cover by two open
subschemes $X=U\cup V$. For an open immersion $j:W\To X$, write $M_W=j_\sharp \, j^*(M)$.
If the restrictions of $M$ to $U$ and $V$ are constructible, then
so is its restriction to $U\cap V$.
According to Proposition \ref{localizationNisnevichdescent0},
 we get a distinguished triangle
$$M_{U\cap V} \To M_U \oplus M_V \To M \To M_{U\cap V}[1]$$
in which $M_W$ is constructible for $W=U,V,U\cap V$
(using \ref{soritestauconstruct2} again). 
Thus the premotive $M$ is constructible.
\end{proof}

\begin{cor}\label{Thomconstruct}
For any scheme $X$ and any vector bundle $E$ over $X$, the functors
$\Th(E)$ and $\Th(-E)$ preserve constructible objects in $\T(X)$.
\end{cor}

\begin{proof}
To prove that $\Th(E)$ and $\Th(-E)$ preserves constructible
objects, by virtue of the preceding proposition, we may assume that $E$ is trivial of rank $r$.
It is thus sufficient to prove that $M(r)$ is constructible whenever $M$ is so
for any integer $r$. For we may assume that $M=\unit_X\{n\}$
for some $n\in \tau$ (using \ref{soritestauconstruct2}),
this is true by assumption on $\tau$; see \ref{assumption_constructible_motives}(c).
\end{proof}

\begin{cor}\label{thmfinitnessredlocZar}
Let $f:X\To Y$ a morphism of finite type. The property that the functor
$$f_*:\T(X)\To \T(Y)$$
preserves constructible objects is local on $Y$ with respect to
the Zariski topology.
\end{cor}

\begin{proof}
Consider a finite Zariski cover $\{v_i:Y_i\to Y\}_{i\in I}$, and write
$f_i:X_i\To Y_i$ for the pullback of $f$ along $v_i$ for each $i$ in $I$.
Assume that the functors $f_{i,*}$ preserves constructible objects;
we shall prove that $f_*$ has the same property.
Let $M$ be a constructible object in $\T(X)$. Then
for $i\in I$, using the smooth base change isomorphism (for open immersions),
we see that the restriction of $f_*(M)$ to $Y_i$ is isomorphic to the
image by $f_{i,*}$ of the restriction of $M$ to $X_i$, hence
is constructible. The preceding proposition thus implies that $f_*(M)$ is constructible.
%% 
%% This property is also local on $X$: if $X=U\cup V$ is an open cover of $X$
%% such that $(f|_U)_*$ and $(f|_V)_*$ preserve constructible objects,
%% then $f_*$ have the same property because the family of objects of shape $M_X(T)\{n\}$,
%% with $T$ smooth over $U$ or $V$ and $n\in \tau$ form a generating family of
%% $\T_c(X)$.
\end{proof}

\begin{prop}\label{thmfinitnessclosedimmersions}
For any closed immersion $i:Z\To X$, the functor
$$i_*:\T(Z)\To \T(X)$$
preserves constructible objects.
\end{prop}

\begin{proof}
It is sufficient to prove that for any smooth $Z$-scheme $Y_0$
 and any twist $n \in \tau$, the premotive $i_*(M_Z(Y_0)\{n\})$
 is constructible in $\T(X)$.
 According to the Mayer-Vietoris triangle (see Remark \ref{rem:MV_triangle}),
 this assertion is local in $X$.
 Thus we can assume there exists a smooth $X$-scheme $Y$ such that
 $Y_0=Y \times_X Z$ (apply \cite[18.1.1]{EGA4}). 
% By virtue of Corollary \ref{thmfinitnessredlocZar}, we may assume that
%$X$ is affine. The category $\T_{c}(Z)$ is then
%the thick subcategory of $\T(Z)$ generated
%by the objects of shape $M_Z(Y_0)\{n\}$ for
%$Y_0$ affine and smooth over $Z$, and $n\in\tau$.
%By virtue of , for any affine and smooth $Z$-scheme $Y_0$,
%there exists a smooth and affine $X$-scheme $Y$
%whose pullback along $i$ is isomorphic to $Y_0$.
%Therefore, the category $\T_{c}(Z)$ is the
%the thick subcategory of $\T(Z)$ generated
%by the objects of shape $i^*(M)$, where $M$
%is a constructible object in $\T(X)$.
%To prove the proposition, it is thus sufficient to prove that
%the functor $i_*\, i^*$ preserves constructible
%objects.
Put $U=X-Z$ and let $j:U\To X$ be the obvious open immersion.
From the localization property, we get a distinguished triangle
$$
M_X(Y \times_X U)\{n\}\To M_X(Y)\{n\} \To i_*(M_Z(Y_0)\{n\})
 \To M_X(Y \times_X U)\{n\}[1]
$$
and this concludes.
\end{proof}

\begin{cor}\label{corstratconstruct}
Let $i:Z\To X$ be a closed immersion with open complement $j:U\To X$.
an object $M$ of $\T(X)$ is constructible if and only if
$j^*(M)$ and $i^*(M)$ are constructible in $\T(U)$
and $\T(Z)$ respectively.
\end{cor}

\begin{proof}
We have a distinguished triangle
$$j_\sharp \, j^*(M)\To M \To i_*\, i^*(M)\To j_\sharp \, j^*(M)[1]\, .$$
Hence this assertion follows from propositions \ref{soritestauconstruct2}
and \ref{thmfinitnessclosedimmersions}.
\end{proof}

\begin{prop}\label{thmfinitnessproper}
If $f:X\To Y$ is proper, then  the functor
$$f_*:\T(X)\To \T(Y)$$
preserves constructible objects.
\end{prop}

\begin{proof}
We shall first consider the case where $f$ is projective.
As this property is local on $Y$ (Corollary \ref{thmfinitnessredlocZar}), we may assume
that $f$ factors as a closed immersion $i:X\To \PP^n_Y$
followed by the canonical projection $p:\PP^n_Y \To Y$.
By virtue of Proposition \ref{thmfinitnessclosedimmersions}, we can assume that $f=p$.
In this case, the functor $p_*$ is isomorphic to $p_\sharp$
composed with the quasi-inverse of the Thom endofunctor associated to the cotangent bundle of $p$;
see \ref{thm:cor3_Ayoub}~(3). Therefore, the functor $p_*$
preserves constructible objects by virtue of Proposition \ref{soritestauconstruct2}
and of Corollary \ref{Thomconstruct}.
The case where $f$ is proper follows easily from
the projective case, using Chow's lemma and $\cdh$-descent
(the homotopy pullback squares \eqref{cdhdescent00}),
by induction on the dimension of $X$.
\end{proof}

\begin{cor}\label{corthmfinitnessproper}
If $f:X\To Y$ is separated of finite type, then the functor
$$f_!:\T(X)\To\T(Y)$$
preserves constructible objects.
\end{cor}

\begin{proof}
It is sufficient to treat the case where $f$ is either an open immersion, either
a proper morphism, which follows respectively from \ref{soritestauconstruct2}
and \ref{thmfinitnessproper}.
\end{proof}

\begin{prop}\label{thmfinitnessproper2}
Let $X$ be a scheme. The category of constructible
objects in $\T(X)$ is the smallest thick triangulated subcategory
which contains the objects of shape $f_*(\unit_{X'}\{n\})$,
where $f:X'\To X$ is a (strictly) projective morphism, and $n\in \tau$.
\end{prop}

\begin{proof}
Let $\T_p(X)$ be the smallest thick triangulated subcategory
which contains the objects of shape $f_*(\unit_{X'}\{n\})$,
where $f:X'\To X$ is a (strictly) projective morphism, and $n\in \tau$. Proposition \ref{thmfinitnessproper}
shows that $\T_p(X)\subset\T_c(X)$, to that it is sufficient
to prove the reverse inclusion.
Note that, for any separated smooth morphism $f$,
locally for the Zariski topology, $f_\sharp$ coincides
with $f_!$ up to a Tate twist. In other words, it is
sufficient to prove that, for any separated
morphism of finite type $f:Y\To X$, $f_!(\un_Y)$
belongs to $\T_p(X)$. If we factor $f$ into
an open immersion $j:Y\To X'$ followed by
a proper morphism $p:X'\To X$, we see that is sufficient
to prove that $j_\sharp(\un_Y)$ belongs to $\T_p(X')$.
This follows straight away from the localization property.
\end{proof}

The following lemma is the key geometrical point
 for the finiteness Theorem \ref{thmfinitness0}
\begin{lm}\label{gabbergeomlemmathmfinitness}
Let $j:U\To X$ be a dense open immersion such that $X$ is quasi-excellent.
Then, there exists the following data:
\begin{itemize}
\item[(i)] a finite $h$-cover $\{f_i:Y_i\To X\}_{i\in I}$
such that for all $i$ in $I$, $f_i$ is a morphism of finite type, the scheme
$Y_i$ is regular, and $f^{-1}_i(U)$ is either $Y_i$ itself or the
complement of a strict normal crossing divisor in $Y_i$; we shall write
$$f:Y=\coprod_{i\in I}Y_i\To X$$
for the induced global $h$-cover;
\item[(ii)] a commutative diagram
\begin{equation}\label{diagthmfinitness}\begin{split}
\xymatrix{
X'''\ar[rr]^g\ar[d]_q&&Y\ar[d]^f\\
X''\ar[r]^u&X'\ar[r]^p&  X
}\end{split}\end{equation}
in which: $p$ is a proper birational morphism, $X'$ is normal,
$u$ is a Nisnevich cover, and $q$ is a finite surjective morphism.
\end{itemize}
Let $T$ (resp. $T'$) be a closed subscheme of $X$ (resp. $X'$)
 and assume that for any irreducible component $T_0$ of $T$,
 the following inequality is satisfied:
\begin{equation} \label{eq:gabberlemma}
\mathrm{codim}_{X'}(T')\geq \mathrm{codim}_X(T_0),
\end{equation}
Then, possibly after shrinking $X$ in an open neighborhood of the generic points of $T$ in $X$,
one can replace $X''$ by an open cover and $X'''$ by its pullback along this cover,
in such a way that we have in addition the following properties:
\begin{itemize}
\item[(iii)] $p(T') \subset T$
 and the induced map $T' \rightarrow T$
 is finite and pseudo-dominant;\footnote{Recall from \ref{num:properties_morphism_cycles}
that this means that any irreducible component
  of $T'$ dominates an irreducible component of $T$.}
%$p(T' \times_X )\subset T$, the induced map $p:T'\To T$ is finite
% and any irreducible component of $T'$ 
%
%while the restriction of $p$ to each of the irreducible components of $T'$
%is dominant over an irreducible component of $T$; in particular, we have
%$$T'=\varnothing\quad\text{or}\quad \mathrm{codim}_{X'} T'=\mathrm{codim}_X T\, ;$$
\item[(iv)] if we write $T''=u^{-1}(T')$, the induced map $T''\To T'$
is an isomorphism.
\end{itemize}
\end{lm}

\begin{proof}
The existence of $f:Y\To X$ as in (i) follows from
Gabber's weak uniformization theorem (see \ref{Gabber1}),
 while the commutative diagram \eqref{diagthmfinitness}
 satisfying property (ii) is ensured by Lemma \ref{refinedhvoverings}.

Consider moreover closed subschemes $T \subset X$ and $T' \subset X'$
 satisfying \eqref{eq:gabberlemma}.

We first show that,
 by shrinking $X$ in an open neighborhood of the generic points of $T$
 and by replacing the diagram \eqref{diagthmfinitness} by its pullback
 over this neighborhood, we can assume that condition (iii) is satisfied.
 Note that shrinking $X$ in this way does not change the condition
 \eqref{eq:gabberlemma} because $\mathrm{codim}_X(T_0)$ does not change
  and $\mathrm{codim}_{X'}(T')$ can only increase.\footnote{Remember
  that for any scheme $X$, $\mathrm{codim}_X(\varnothing)=+\infty$.}

Note first that, by shrinking $X$,
 we can assume that any irreducible component $T'_0$ of $T'$
 dominates an irreducible component $T_0$ of $T$.
In fact, given an irreducible component $T'_0$ which does not satisfy
 this condition, $p(T'_0)$ is a closed subscheme of $X$ disjoint from the set
 of generic points of $T$ and
 replacing $X$ by $X-f(T'_0)$, we can throw out $T'_0$.

Further, shrinking $X$ again,
 we can assume that for any pair $(T'_0,T_0)$ as in the preceding paragraph,
 $p(T'_0) \subset T_0$. In fact, in any case, as $p(T'_0)$ is closed
 we get that $T_0 \subset p(T'_0)$. Let $Z$ be the closure of $p(T'_0)-T_0$ in $X$.
 Then $Z$ does not contain any generic point of $T$ (because $p(T'_0)$ is irreducible),
 and $p(T'_0) \cap (X-Z) \subset T_0$. Thus it is sufficient to replace $X$ by $X-Z$
 to ensure this assumption.

Consider again a pair $(T'_0,T_0)$ as in the two preceding paragraphs
 and the induced commutative square:
\begin{equation}\label{diagthmfinitnessproof1}\begin{split}
\xymatrix@=20pt{
T'_0\ar[r]\ar_{p_0}[d] & X'\ar^p[d] \\
T_0\ar[r] & X
}\end{split}\end{equation}

We show that the map $p_0$ is generically finite.
In fact, this will conclude the first step,
 because if it is true for any irreducible component $T'_0$
 of $T'$, we can shrink $X$ again so that the dominant
 morphism $p_0:T'_0 \rightarrow T_0$ becomes finite.

Denote by $c'$ (resp. $c)$ the codimension of $T_0$ in $X$
 (resp. $T'_0$ in $X'$). Note that \eqref{eq:gabberlemma} gives the inequality $c' \geq c$.
Let $t_0$ be the generic point of $T_0$, $\Omega$ the localization of $X$ at $t_0$
 and consider the pullback of \eqref{diagthmfinitnessproof1}:
\begin{equation}\label{diagthmfinitnessproof2}\begin{split}
\xymatrix@=20pt{
W'\ar[r]\ar_{q_0}[d] & \Omega'\ar^q[d] \\
\{t_0\}\ar[r] & \Omega.
}\end{split}\end{equation}
We have to prove that $\dim(W')=0$.
Consider an irreducible component $\Omega'_0$ of $\Omega'$ containing $W'$.
As $q$ is still proper birational, $\Omega'_0$ corresponds to a unique
irreducible component $\Omega_0$ of $\Omega$ such that $q$ induces
a proper birational map $\Omega_0' \rightarrow \Omega_0$.
According to \cite[5.6.6]{EGA4}, we get the inequality
$$
\dim(\Omega_0') \leq \dim(\Omega_0).
$$
Thus, we obtain the following inequalities:
\begin{align*}
\dim(W')
&\leq \dim(\Omega_0')-\mathrm{codim}_{\Omega_0'}(W')\\
&\leq \dim(\Omega_0)-\mathrm{codim}_{\Omega_0'}(W')\\
&\leq \dim(\Omega)-\mathrm{codim}_{\Omega_0'}(W').
\end{align*}
As this is true for any irreducible component $\Omega'_0$ of $\Omega'$,
 we finally obtain:
$$
\dim(W')\leq \dim(\Omega)-\mathrm{codim}_{\Omega'}(W') \leq c-c'
$$
and this concludes the first step.

Keeping $T'$ and $T$ as above, as the map from $T''$ to $T'$
is a Nisnevich cover, it is a split epimorphism in a neighborhood
of the generic points of $T'$ in $X'$. 
Hence, as the map $X'\To X$ is proper and birational, we can find a
neighborhood of the generic points of $T$ in $X$ over which the map
$T''\To T'$ admits a section $s:T'\To T''$.
Let $S$ be a closed subset of $X''$ such that
$T''=s(T')\amalg S$ (which exists because $X''\To X'$ is \'etale).
The map $(X''-T'')\amalg(X''-S)\To X'$ is then a Nisnevich cover.
Replacing $X''$ by $(X''-T'')\amalg(X''-S)$ (and $X'''$ by the pullback of
$X'''\To X''$ along $(X''-T'')\amalg(X''-S)\To X'$), we may assume that
the induced map $T''\to T'$ is an isomorphism, without modifying further
the data $f$, $p$, $T$ and $T'$. This gives property (iv) and
ends the proof the lemma.
\end{proof}

\begin{paragr}\label{assumptionsfiniteness}
Let $\T_0$ be a full $\mathit{Open}$-fibred subcategory of $\T$
(where $\mathit{Open}$ stands for the class of open immersions).
We assume that $\T_0$ has the following properties:
\begin{itemize}
\item[(a)] for any scheme $X$ in $\sch$, $\T_0(X)$ is
a thick subcategory of $\T(X)$ which contains
the objects of the form $\unit_X\{n\}$, $n\in\tau$;
%%\item[(b)] for any morphism $f:X\To Y$ in $\sch$,
%%$\T$ is stable under the operation $f^*$;
%%\item[(c)] for any smooth morphism $f:X\To Y$ in $\sch$,
%%$\T$ is stable under the operation $f_\sharp$;
\item[(b)] for any separated morphism of finite type $f:X\To Y$ in $\sch$,
$\T_0$ is stable under $f_!$;
\item[(c)] for any dense open immersion $j:U\To X$, with $X$ regular,
which is the complement of a strict normal crossing divisor,
$j_*(\unit_U\{n\})$ is in $\T_0(U)$ for any $n\in \tau$.
\end{itemize}
Properties (a) and (b) have the following consequences:
any constructible object belongs to $\T_0$;
given a closed immersion $i:Z\To X$
with complement open immersion $j:U\To X$,
an object $M$ of $\T(X)$ belongs to $\T_0(X)$ if and only
if $j^*(M)$ and $i^*(M)$ belongs to $\T_0(U)$ and $\T_0(Z)$
respectively; for any scheme $X$ in $\sch$, the condition
that an object of $\T(X)$ belongs to $\T_0(X)$ is
local on $X$ for the Zariski topology.
\end{paragr}

\begin{thm}\label{thmfinitness0}
Consider the above hypothesis
 and assume that $\T$ is $\QQ$-linear and separated.
Let $Y$ be a quasi-excellent scheme
\index{word}{scheme!quasi-excellent}
 and  $f:X\To Y$ be a morphism of finite type.
 Then for any constructible object $M$ of $\T(X)$,
  the object $f_*(M)$ belongs to $\T_0(Y)$.
\end{thm}

\begin{proof}
It is sufficient to prove that, for any dense open immersion
$j:U\To X$, and for any $n\in \tau$, the object $j_*(\unit_U\{n\})$
is in $\T_0$. Indeed, assume this is known.
We want to prove that $f_*(M)$ is in $\T_0(Y)$ whenever $M$ is constructible.
We deduce from property (b) of \ref{assumptionsfiniteness}
and from Proposition \ref{thmfinitnessproper2}
that it is sufficient to consider the case where $M=\unit_X\{n\}$, with $n\in \tau$.
Then, as this property is assumed to be known for dense
open immersions, by an easy Mayer-Vietoris argument, we see that
the condition that $f_*(\unit_X\{n\})$ belongs to $\T_0$
is local on $X$ with respect to the Zariski topology. Therefore, we may assume that
$f$ is separated. Consider a compactification of $f$, i.e. a commutative diagram
$$\xymatrix{
Y\ar[r]^j\ar[d]_f& \bar Y\ar[dl]^{\bar f}\\
X &
}$$
with $j$ a dense open immersion, and $\bar f$ proper.
By property (b) of \ref{assumptionsfiniteness}, we may assume that $f=j$ is a dense
open immersion.

Let $j:U\To X$ be a dense open immersion.
We shall prove by induction on the dimension of $X$
that , for any $n\in \tau$, the object $j_*(\unit_U\{n\})$
is in $\T_0$. The case where $X$ is of dimension $\leq 0$
follows from the fact the map $j$ is then an isomorphism, which implies that
$j_\sharp\simeq j_*$, and allows to conclude  (because $\T_0$
is $\mathit{Open}$-fibred).

Assume that $\mathrm{dim}\, X>0$.
Following an argument used by Gabber \cite{gabber3} in the context of
$\ell$-adic sheaves, we shall prove by induction on $c\ge 0$
that there exists a closed subscheme $T\subset X$
of codimension $>c$ such that, for any $n\in \tau$, the restriction of $j_*(\unit_U\{n\})$
to $X-T$ is in $\T_0(X-T)$. As $X$ is of finite dimension,
this will obviously prove Theorem \ref{thmfinitness0}.
%% 
%%  because, if
%% $\{T_c\}_{c\geq 0}$ is a family of closed subschemes of $X$, each $T_c$ being of codimension $>c$,
%% then we have a cover by open subschemes $X=\bigcup_{c\geq 0} X-T_c$,
%% and we conclude by Corollary \ref{thmfinitnessredlocZar0} (as $X$ is noetherian).

The case where $c=0$ is clear: we can choose $T$ such that $X-T=U$.
If $c>0$, we choose a closed subscheme $T$ of $X$, of codimension $>c-1$,
such that the restriction of $j_*( \unit_U \{ n \} )$
to $X-T$ is in $\T_0$. It is then sufficient to find a dense
open subscheme $V$ of $X$, which contains all the generic points of $T$, and such that
the restriction of $j_*(\unit_U\{n\})$ to $V$ is in $\T_0$: for such a $V$,
we shall obtain that
the restriction of $j_*(\unit_U\{n\})$ to $V\cup (X-T)$ is in $\T_0$,
the complement of $V\cup (X-T)$ being the support of a closed subscheme of
codimension $>c$ in $X$. In particular, using the smooth base change isomorphism
(for open immersions), we can always replace $X$ by a generic neighborhood of $T$.
It is sufficient to prove that, possibly after shrinking $X$ as above,
the pullback of $j_*(\unit_U\{n\})$ along $T\To X$
is in $\T_0$ (as we already know that its restriction to $X-T$ is in $\T_0$).

We may assume that $T$ is purely of codimension $c$.
We may assume that we have data as in points (i) and (ii) of Lemma \ref{gabbergeomlemmathmfinitness}.
We let $j':U' \rightarrow X'$ denote the pullback of $j$ along $p:X' \rightarrow X$.
Then, we can find, by induction on $c$, a closed subscheme $T'$ in $X'$, of codimension $>c-1$,
such that the restriction of $j'_*(\unit_{U'}\{n\})$
to $X'-T'$ is in $\T_0$. By shrinking $X$,
we may assume that conditions (iii) and (iv) of Lemma \ref{gabbergeomlemmathmfinitness}
are fulfilled as well.

For an $X$-scheme $w:W\To X$ and a closed subscheme $Z\subset W$, we shall write
$$\varphi(W,Z)=w_*\, i_*\, i^* \, j_{W,*}\, j^*_W(\unit_W\{n\}) \, ,$$
where $i:Z\To W$ denotes the inclusion, and $j_W:W_U\To W$
stands for the pullback of $j$ along $w$.
This construction is functorial with respect to morphisms of pairs of
$X$-schemes: if $W'\To W$ is a morphism of $X$-schemes, with $Z'$ and $Z$
two closed subschemes of $W'$ and $W$ respectively, such that $Z'$ is sent to $Z$,
then we get a natural map $\varphi(W,Z)\To \varphi(W',Z')$.
Remember that we want to prove that $\varphi(X,T)$ is in $\T_0$.
This will be done via the following lemmas (which hold assuming
all the conditions stated in Lemma \ref{gabbergeomlemmathmfinitness}
as well as our inductive assumptions).

\begin{lm}\label{thmfinitness01}
The cone of the map $\varphi(X,T)\To \varphi(X',T')$ is in $\T_0$.
\end{lm}

\noindent The map $\varphi(X,T)\To \varphi(X',T')$ factors as
$$\varphi(X,T)\To \varphi(X',p^{-1}(T))\To\varphi(X',T')\, .$$
By the octahedral axiom, it is sufficient to prove that each of these two maps
has a cone in $\T_0$.

We shall prove first that the cone of the map $\varphi(X',p^{-1}(T))\To \varphi(X',T')$
is in $\T_0$. Given an immersion $a:S\To X'$, we shall write
$$M_S=a_! \, a^*(M)\, .$$
We then have distinguished triangles
$$M_{p^{-1}(T)-T'}\To M_{p^{-1}(T)}
\To M_{T'} \To M_{p^{-1}(T)-T'}[1]\, .$$
For $M=j'_*(\unit_{U'}\{n\})$ (recall $j'$ is the pullback of $j$ along $p$)
the image of this triangle by $p_*$ gives a distinguished triangle
$$p_*(M_{p^{-1}(T)-T'})\To \varphi(X',p^{-1}(T))
\To\varphi(X',T') \To p_*(M_{p^{-1}(T)-T'})[1]\, .$$
As the restriction of $M=j'_*(\unit_{U'}\{n\})$ to $X'-T'$ is in $\T_0$
by assumption on $T'$, the object $M_{p^{-1}(T)-T'}$ is in $\T_0$
as well (by property (b) of \ref{assumptionsfiniteness} and because $\T_0$
is $\mathit{Open}$-fibred), from which
we deduce that $p_*(M_{p^{-1}(T)-T'})$ is in $\T_0$
 (using condition (iii) of Lemma \ref{gabbergeomlemmathmfinitness}
  and property (b) of \ref{assumptionsfiniteness}).

Let $V$ be a dense open subscheme of $X$ such that $p^{-1}(V)\to V$
is an isomorphism. We may assume that $V\subset U$, and write
$i: Z\To U$ for the complement closed immersion.
Let $p_U:U'=p^{-1}(U)\To U$ be the pullback of $p$ along $j$, and
let $\bar Z$ be the reduced closure of $Z$ in $X$. 
We thus get the commutative squares
of immersions below,
\begin{equation*}
\begin{split}
\xymatrix{
Z\ar[r]^k\ar[d]_i& \bar Z\ar[d]^l\\
U\ar[r]_j& X
}\end{split}
\quad\text{and}\quad
\begin{split}
\xymatrix{
Z'\ar[r]^{k'}\ar[d]_{i'}& \bar Z'\ar[d]^{l'}\\
U'\ar[r]_{j'}& X'
}\end{split}
\end{equation*}
where the square on the right is obtained from the one on the left by pulling back
along $p:X'\To X$.
As $p$ is an isomorphism over $V$,
we get by $\cdh$-descent (Proposition \ref{cdhdescent}) 
%%(and using the proper base change formula)
the homotopy pullback square below.
%% the restriction of the map
%% $$\unit_U\{n\}\To p_{U,*}(\unit_{U'}\{n\})$$
%% to $V$ is an isomorphism, from which we deduce, by the localization
%% property, that the commutative square
$$\xymatrix{
\unit_U\{n\}\ar[r]\ar[d]& p_{U,*}(\unit_{U'}\{n\})\ar[d]\\
i_*\, i^*(\unit_{Z}\{n\})\ar[r]& i_*\, i^* \, p_{U,*}(\unit_{U'}\{n\})
}$$
%%is homotopy cartesian.
If $a:T\To X$ denotes the inclusion, applying the functor
$a_*\, a^*\, j_*$ to the commutative square above,
we see from the proper base change formula and from the identification
$j_*\, i_*\simeq l_*\, k_*$
that we get a commutative square
isomorphic to the following one
$$\xymatrix{
\varphi(X,T)\ar[r]\ar[d]&\varphi(X',p^{-1}(T))\ar[d] \\
\varphi(\bar Z,\bar Z\cap T)\ar[r]&\varphi(\bar Z',p^{-1}(\bar Z\cap T))\, ,
}$$
which is thus homotopy cartesian as well. It is sufficient to prove that the two objects
$\varphi(\bar Z,\bar Z\cap T)$ and $\varphi(\bar Z',p^{-1}(\bar Z\cap T))$
are in $\T_0$. It follows from the proper base change formula
that the object $\varphi(\bar Z,\bar Z\cap T)$ is canonically isomorphic to
the restriction to $T$ of $l_*\, k_*(\unit_{Z}\{n\})$.
As $\mathrm{dim}\, \bar Z<\mathrm{dim}\, X$,
we know that the object $k_*(\unit_{Z}\{n\})$ is in $\T_0$.
By property (b) of \ref{assumptionsfiniteness},
we obtain that $\varphi(\bar Z,\bar Z\cap T)$ is in $\T_0$.
Similarly, the object $\varphi(\bar Z',p^{-1}(\bar Z\cap T))$ is canonically isomorphic to
the restriction of $p_*\, l'_*\, k'_*(\unit_{Z'}\{n\})$ to $T$, and,
as $\mathrm{dim}\, \bar Z'<\mathrm{dim}\, X'$
(because, $p$ being an isomorphism over the dense open subscheme $V$ of $X$, $\bar Z'$
does not contain any generic point of $X'$), 
$k'_*(\unit_{Z'}\{n\})$ is in $\T_0$.
We deduce again from property (b) of \ref{assumptionsfiniteness}
that $\varphi(\bar Z',p^{-1}(\bar Z\cap T))$ is in $\T_0$ as well,
which achieves the proof of the lemma.

\begin{lm}\label{thmfinitness03}
The map $\varphi(X',T')\To \varphi(X'',T'')$
is an isomorphism in $\T(X)$.
\end{lm}

\noindent Condition (iv) of Lemma \ref{gabbergeomlemmathmfinitness}
can be reformulated by saying that
we have the Nisnevich distinguished square below.
$$\xymatrix{
X''-T''\ar[r]\ar[d]& X''\ar[d]^{v}\\
X'-T'\ar[r]&X'
}$$
This lemma follows then by Nisnevich excision
 (Proposition \ref{localizationNisnevichdescent0})
and smooth base change (for \'etale maps).

\begin{lm}\label{thmfinitness02}
Let $T'''$ be the pullback of $T''$ along the finite surjective morphism
$X'''\To X''$. The map $\varphi(X'',T'')\To \varphi(X''',T''')$
is a split monomorphism in $\T(X)$.
\end{lm}

\noindent We have the following pullback squares
$$\xymatrix{
T'''\ar[r]^t\ar[d]_{r}& X'''\ar[d]^q & U'''\ar[l]_{j'''}\ar[d]^{q^{}_U}\\
T''\ar[r]^s& X'' & U'\ar[l]_{j''} 
}$$
in which $j''$ and $j'''$ denote the pullback of $j$ along $pu$
and $puq$ respectively, while $s$ and $t$ are the inclusions.
By the proper base change formula applied to the left-hand square,
we see that the map $\varphi(X'',T'')\To \varphi(X''',T''')$ is isomorphic to
the image of the map
$$j''_*(\unit_{U''}\{n\})\To q_*\, q^* \, j''_*(\unit_{U''}\{n\})
\To q_*\, j'''_*(\unit_{U'''}\{n\}) \, .$$
by $f_*\, s^*$, where $f:T''\To T$ is the map induced by $p$
(note that $f$ is proper as $T''\simeq T'$ by assumption).
As $q_*\, j'''_*\simeq j''_*\, q_{U,*}$,
we are thus reduced to prove that the unit map
$$\unit_{U''}\{n\} \To q_{U,*}(\unit_{U'''}\{n\})$$
is a split monomorphism.
As $X''$ is normal (because $X'$ is so by assumption,
while $X''\To X'$ is \'etale), this follows immediately from Corollary \ref{stronggenericGaloisdescentcor}.

Now, we can finish the proof of Theorem \ref{thmfinitness0}.
Consider the Verdier quotient
$$D=\T(X)/\T_0(X)\, .$$
We want to prove that, under the conditions stated in Lemma \ref{gabbergeomlemmathmfinitness},
we have $\varphi(X,T)\simeq 0$ in $D$.
Let $\pi: T'''\To X$ be the map induced by $puq:X'''\To X$.
If $a:T'''\To Y$ denotes the map induced by $g:X'''\To Y$,
and $j_Y:Y_U\To Y$ the pullback of $j$ by $f$,
we have the commutative diagram below.
$$\xymatrix{
\varphi(X,T)\ar[rr]\ar[dr]&& \varphi(X''',T''')\\
& \pi_* \, a^*\, j_{Y,*}(\unit_{Y_U}\{n\})\ar[ur]& 
}$$
By virtue of lemmas \ref{thmfinitness01}, \ref{thmfinitness02}, and
\ref{thmfinitness03}, the horizontal map
is a split monomorphism in $D$. It is thus sufficient to prove that
this map vanishes in $D$, for which it will be sufficient to prove that
$\pi_* \, a^*\, j_{Y,*}(\unit_{Y_U}\{n\})$ is
in $\T_0$. The morphism $\pi$ is finite (by construction,
the map $T''\To T'$ is an isomorphism, while the maps
$T'''\To T''$ and $T'\To T$ are finite). Under this condition, $\T_0$
is stable under the operations $\pi_*$ and $a^*$.
To finish the proof of the theorem, it remains to check that
$j_{Y,*}(\unit_{Y_U}\{n\})$ is in $\T_0$, which
follows from property (c) of \ref{assumptionsfiniteness} (and additivity).
\end{proof}

\begin{df}\label{df:weaklytaupure}
We shall say that $\T$ is \emph{$\tau$-compatible}
 if it satisfies the following two conditions.
\begin{itemize}
\item[(a)] For any closed immersion $i:Z\To X$ between regular schemes
 in $\sch$, 
 the image of $\unit_X\{n\}$, $n\in \tau$, 
 by the exceptional inverse image functor
$i^!:\T(X)\To \T(Z)$ is constructible.
\item[(b)] For any scheme $X$, any $n\in\tau$, and any constructible
object $M$  in $\T(X)$, the object $\Hom_X(\unit_X\{n\},M)$ is constructible.
\end{itemize}
As usual,
 when $\tau$ is the monoid generated by the Tate twist,
 we say \emph{compatible with Tate twists}.
\end{df}

\begin{rem}
Condition (b) of the definition above will come essentially for free
if the objects $\unit_X\{n\}$ are $\otimes$-invertible with
constructible $\otimes$-quasi-inverse (which will hold
in practice, essentially by definition).
\end{rem}

\begin{ex}
In practice, condition (a) of the definition above
 will be a consequence of the \emph{absolute purity theorem}.
 In particular,
  the category of Beilinson motives $\DMB$ is compatible with Tate twist
  as a corollary of the fact the Tate twist is invertible and
  Theorem \ref{DMBpurity}.
\end{ex}

\begin{lm}\label{thmfinitnessncd}
Assume that $\T$ is $\tau$-compatible.
Let $i:Z\To X$ be a closed immersion, with $X$ regular, and $Z$ the support
of a strict normal crossing divisor. Then $i^!(\unit_X\{n\})$
is constructible for any $n\in \tau$. As a consequence, if $j:U\To X$
denotes the complement open immersion, then
$j_*(\unit_U\{n\})$ is constructible for any $n\in \tau$.
\end{lm}

\begin{proof}
The first assertion follows easily by induction on the number of irreducible components of $Z$,
using Proposition \ref{thmfinitnessredlocZar0}. The second assertion follows from the
distinguished triangles
$$i_*\, i^!(M)\To M \To j_*\, j^*(M) \To i_*\, i^!(M)[1]$$
and from Lemma \ref{thmfinitnessclosedimmersions}.
\end{proof}

\begin{thm}\label{thmfinitness}
Assume that $\T$ is $\QQ$-linear, separated, and $\tau$-compatible.
Then,
for any morphism of finite type $f:X\To Y$ such that $Y$ is quasi-excellent,
 the functor
$$f_*:\T(X)\To \T(Y)$$
preserves constructible objects.
\end{thm}

\begin{proof}
By virtue of propositions \ref{soritestauconstruct2} and \ref{thmfinitnessproper}
as well as of Lemma \ref{thmfinitnessncd}, if $\T$ is $\tau$-compatible,
we can apply Theorem \ref{thmfinitness0}, where $\T$ stands for the subcategory
of constructible objects.
\end{proof}

\begin{cor}\label{constructibleinternalHom}
Under the assumptions of the above theorem, 
 for any quasi-excellent scheme $X$, and for any couple of constructible
objects $M$ and $N$ in $\T(X)$, the object
$\sHom_X(M,N)$ is constructible.
\end{cor}
\begin{proof}
It is sufficient to treat the case where $M=f_\sharp(\unit_Y\{n\})$,
for $n\in \tau$ and $f:Y\To X$ a smooth morphism.
But then, we have, by transposition of the $\sm$-projection formula,
a natural isomorphism:
$$\sHom_X(M,N)\simeq f_*\, \sHom(\unit_Y\{n\},f^*(N))\, .$$
This corollary follows then immediately from
Proposition \ref{soritestauconstruct2} and from Theorem \ref{thmfinitness}.
\end{proof}

\begin{cor}\label{contructexcepinvimclosedimm}
Under the assumptions of the above theorem, 
 for any closed immersion $i:Z\To X$ such that $X$ is quasi-excellent,
 the functor
$$i^!:\T(X)\To\T(Z)$$
preserves constructible objects.
\end{cor}

\begin{proof}
Let $j:U\to X$ be the complement open immersion.
For an object $M$ of $\T(X)$, we have the following distinguished
triangle.
$$i_*\, i^!(M)\To M \To j_*\, j^*(M) \To i_*\, i^!(M)[1] \, .$$
By virtue of Proposition \ref{thmfinitnessredlocZar0} and Theorem \ref{thmfinitness},
if $M$ is constructible, then $j_*\, j^*(M)$ have the same property,
which allows us to conclude.
\end{proof}

\begin{lm}\label{constructexcepinvimlocal}
Let $f:X\To Y$ be a separated morphism of finite type.
The condition that the functor $f^!$ preserves constructible objects in $\T$
is local over $X$ and over $Y$ for the Zariski topology.
\end{lm}

\begin{proof}
If $u:X'\To X$ is a Zariski cover, then we have, by definition, $u^*=u^!$,
so that, by Proposition \ref{thmfinitnessredlocZar0}, the condition that $f^!$
preserves $\tau$-constructibility is equivalent to the condition that
the functors
$u^*\, f^!\simeq (fu)^!$ preserve $\tau$-constructibility.
Let $v:Y'\To Y$ be a Zariski cover, and consider the following pullback square.
$$\xymatrix{
X'\ar[r]^u\ar[d]_g& X\ar[d]^f\\
Y'\ar[r]_v& Y
}$$
We then have a natural isomorphism $u^*\, f^!\simeq g^!\, v^*$, and, as $u$
is still a Zariski cover, we deduce again from Proposition \ref{thmfinitnessredlocZar0}
that, if $g^!$ preserves $\tau$-constructibility, so does $f^!$.
\end{proof}

\begin{cor}\label{constructexcepinvim}
Under the assumptions of the above theorem,
 for any separated morphism of finite type $f:X\To Y$, the functor
$$f^!:\T(Y)\To\T(X)$$
preserves constructible objects.
\end{cor}

\begin{proof}
By virtue of the preceding lemma, we may assume that $f$ is affine.
We can then factor $f$ as an immersion $i:X\To \AA^n_Y$
followed by the canonical projection $p:\AA^n_Y \To Y$.
The case of an immersion is reduced to
 the case of an open immersion (\ref{soritestauconstruct2})
 and to the case of a closed immersion (\ref{contructexcepinvimclosedimm}). 
 Thus we may assume that $f=p$, in which case $p^!\simeq p^*(-)(n)[2n]$
 (according to point (3) of Theorem \ref{thm:cor3_Ayoub}), 
 so that we conclude by
 \ref{soritestauconstruct2} and \ref{thmfinitnessclosedimmersions}.
\end{proof}

In conclusion, we have proved the following finiteness theorem:
\index{word}{finiteness theorem}
\begin{thm}\label{grothendieck6op}
Assume the motivic triangulated category $\T$ is
 $\QQ$-linear, separated and $\tau$-compatible.\footnote{Remember
 also that $\T$ is associated with a combinatorial stable
  premotivic model category.}

Then constructible objects of $\T$ are closed under
the six operations of Grothendieck when restricted
 to the subcategory $\sch'$ of $\sch$
  made of quasi-excellent schemes and morphisms of finite type.
In particular, $\T_{c}$ is a $\tau$-generated
\index{word}{generated!$\tau$-generated}
 motivic category over $\sch'$.
\end{thm}

\subsection{Continuity} \label{sec:continuity}

\begin{paragr} \label{paragr:continuous}
For the next definition,
 we consider an admissible class $\Pmor$ of morphisms in $\sch$
 and an abstract symmetric monoidal $\Pmor$-fibred category $\T$ over $\sch$.

Let $(S_\alpha)_{\alpha\in A}$ be a projective system of schemes 
\index{word}{projective system, of schemes}
 in $\sch$, with affine transition maps, and
such that $S=\varprojlim_{\alpha\in A}S_\alpha$ is representable in $\sch$
(we assume that $A$ is a partially ordered set to keep the notations simple). 
For each index $\alpha$, we denote by $p^{}_\alpha:S\To S_\alpha$
the canonical projection. Given an index $\alpha_0\in A$ and
an object $E_{\alpha_0}$ in $\T(S_{\alpha_0})$, we write
$E_\alpha$ for the pullback of $E_{\alpha_0}$ along the map $S_\alpha\To S_{\alpha_0}$,
and put $E=\derL p^*_{\alpha}(E_\alpha)$.
We will say that $(S_\alpha)_{\alpha \in A}$ is \emph{dominant} if
 each transition map is furthermore dominant.
\end{paragr}
\begin{df}\label{df:continuous}
Consider the assumptions above and let $\tau$ be a set of twists of $\T$.

We say that $\T$ is \emph{$\tau$-continuous}
 (resp. \emph{weakly $\tau$-continuous}),
or simply \emph{continuous} (resp. \emph{weakly continuous})
if $\tau$ is clearly specified by the context,
 if it is $\tau$-generated and if,
 given any projective system (resp. dominant projective system)
 of schemes $(S_\alpha)$ as above,
 for any index $\alpha_0$, any object $E_{\alpha_0}$ in $\T(S_{\alpha_0})$,
 and any twist $n \in \tau$, the canonical map
$$
\varinjlim_{\alpha\geq \alpha_0}
 \Hom_{\T(S_\alpha)}(\unit_{S_\alpha}\{n\},E_\alpha)
 \To \Hom_{\T(S)}(\unit_S\{n\},E),
$$
is bijective.
%
%We will simply say continuous instead of $\tau$-continuous
% in one of the following situations:
%\begin {itemize}
%\item $\tau=\varnothing$;
%\item  $\T$ is a motivic triangulated category
% and $\tau=\ZZ$ is the group generated by the Tate twist.
%\end {itemize}
\end{df}

\begin{ex}
The main examples of $\tau$-continuous categories
 will be seen afterwards:
\begin{itemize}
\item the $\AA^1$-derived category $\DMt$
 (Example \ref{DMtcontinuous});
\item the motivic category $\DMB$ of Beilinson motives 
 (Proposition \ref{DMBcontinuous}).
\end{itemize}
The triangulated motivic category
 of motivic complexes $\DM$,
 as well as its effective counterpart $\DMe$,
 is weakly continuous (Theorem \ref{thm:DM_continuity}).
 We are only able to prove it is continuous in some special cases
 (namely when it compares to Beilinson motives,
 see Theorem \ref{comparisonDMBDMV}). 
\end{ex}

The interest of the continuity property is to allow
 a description of constructible
objects over $S$ in terms of constructible objects over the $S_\alpha$'s.
\begin{prop}\label{continuityconstructible}
Assume, under the hypothesis of \ref{paragr:continuous},
 that $\T$ is $\tau$-continuous (resp. weakly $\tau$-continuous).
Consider a scheme $S$ in $\sch$,
 as well as a projective system of schemes $(S_\alpha)_{\alpha \in A}$ in $\sch$
 with affine (resp. affine dominant) transition maps
 and such that $S=\varprojlim_\alpha S_\alpha$. 

Then, for any index $\alpha_0$,
and for any objects $C_{\alpha_0}$ and $E_{\alpha_0}$ in $\T(S_{\alpha_0})$,
if $C_{\alpha_0}$ is constructible, then the canonical map
\begin{equation}\label{continuityconstructible1}
\varinjlim_{\alpha\geq\alpha_0}\Hom_{\T(S_\alpha)}(C_\alpha,E_\alpha)
\To \Hom_{\T(S)}(C,E)
\end{equation}
is bijective. Moreover, the canonical functor
\begin{equation}\label{continuityconstructible2}
2\mbox{-}\varinjlim_\alpha\T_c(S_\alpha)\To \T_c(S)
\end{equation}
is an equivalence of monoidal triangulated categories.
\end{prop}

\begin{proof}
To prove the first assertion, we may assume, without loss of generality, that
$C_{\alpha_0}=M_{S_{\alpha_0}}(X_{\alpha_0})\{n\}$ for some smooth
$S_{\alpha_0}$-scheme of finite type $X_{\alpha_0}$, and $n\in \tau$.
Consider an object $E_{\alpha_0}$ in $\T(S_{\alpha_0})$.
For $\alpha\geq\alpha_0$, write $X_\alpha$ (resp. $E_\alpha$)
for the pullback of $X_{\alpha_0}$ (resp. of $E_{\alpha_0}$)
along the map $S_{\alpha}\To S_{\alpha_0}$. Similarly, write $X$
(resp. $E$) for the pullback of $X_{\alpha_0}$ (resp. of $E_{\alpha_0}$)
along the map $S\To S_{\alpha_0}$. We shall also write $E'_\alpha$
(resp. $E'$) for the pullback of $E_\alpha$ (resp. $E$) along the smooth map
$X_\alpha\To S_\alpha$ (resp. $X\To S$). Then, $(X_\alpha)$ is 
a projective system of schemes in $\sch$, with affine transition maps,
such that $X=\varprojlim_\alpha X_\alpha$. Note that 
 if $(S_\alpha)$ is dominant in the sense of Paragraph \ref{paragr:continuous},
 then $(X_\alpha)$ is dominant, as dominant morphisms are stable under
 smooth base change.
Then, by continuity (resp. weak continuity), we have the following
natural isomorphism, which proves the first assertion.
$$\begin{aligned}
\varinjlim_\alpha\Hom_{\T(S_\alpha)}(M_{S_\alpha}(X_\alpha)\{n\},E_\alpha)
\simeq & \varinjlim_\alpha\Hom_{\T(X_\alpha)}(\unit_{X_\alpha}\{n\},E'_\alpha)\\
\simeq & \Hom_{\T(X)}(\unit_X\{n\},E')\\
\simeq & \Hom_{\T(S)}(M_S(X)\{n\},E)
\end{aligned}$$
Note that the first assertion implies that the functor \eqref{continuityconstructible2}
is fully faithful. Pseudo-abelian triangulated categories are stable under
filtered $2$-colimits. In particular, the source of the functor \eqref{continuityconstructible2}
can be seen as a thick subcategory of $\T(S)$.
The essential surjectivity of \eqref{continuityconstructible2} follows from the fact that,
for any smooth $S$-scheme of finite type $X$, there exists some index $\alpha$,
and some smooth $S_\alpha$-scheme $X_\alpha$, such that $X\simeq S\times_{S_\alpha}X_\alpha$;
see \cite[8.8.2 and 17.7.8]{EGA4}: this implies that the essential image of the
fully faithful functor \eqref{continuityconstructible2} contains all the
objects of shape $M_S(X)\{n\}$ for $n\in \tau$ and $X$ smooth over $S$, so that it contains
$\T_c(S)$, by definition.
\end{proof}

\begin{paragr}\label{settingcontinuity}
Before showing how the assumption of weak continuity can be used in the
 case of motivic categories, we state a proposition which later on 
 will allow us to show continuity or weak continuity in concrete cases.
 Let $\M$ be a symmetric monoidal $\Pmor$-fibred model category $\M$ over $\sch$.

We consider again the assumptions and notations of \ref{paragr:continuous},
assuming the transition maps of the pro-scheme $(S_\alpha)$ are $\Pmor$-morphisms,
with $\T=\ho(\M)$.
For each index $\alpha\in A$, we choose a small set $I_\alpha$
(resp. $J_\alpha$) of generating cofibrations (resp. of generating
trivial cofibration) in $\ho(\M)(S_\alpha)$.
We also choose a small set $I$ (resp. $J$)
of generating cofibrations (resp. of generating
trivial cofibration) in $\ho(\M)(S)$.

Consider the following assumptions:
\begin{itemize}
%%\item[(a)] For any index $\alpha\in A$ and any $\Pmor$-morphism $X_\alpha\To S_\alpha$,
%%the pullback $X=S\times_{S_\alpha}X_\alpha$ is representable in $\sch$.
%%\item[(b)] For any $\Pmor$-morphism $f:X\To S$, there exists
%%an index $\alpha \in A$ as well as a $\Pmor$-morphism $f^{}_\alpha:X_\alpha\To S_\alpha$
%%such that $X=S\times_{S_\alpha}X_\alpha$.
\item[(a)] We have
$I\subset\bigcup_{\alpha\in A}p^*_\alpha(I_\alpha)$ and
$J\subset\bigcup_{\alpha\in A}p^*_\alpha(J_\alpha)$.
\item[(b)] For any index $\alpha_0$, if $C_{\alpha_0}$ and
$E_{\alpha_0}$ are two objects of $\M(S_{\alpha_0})$, with $C_{\alpha_0}$
either a source or a target of a map in $I_{\alpha_0}\cup J_{\alpha_0}$,
the natural map
$$\varinjlim_{\alpha\in A}\Hom_{\M(S_\alpha)}(C_{\alpha},E_\alpha)
\To \Hom_{\M(S)}(C,E)$$
is bijective.
\end{itemize}
\end{paragr}

\begin{prop}\label{abstractcontinuity}
Under the assumptions of \ref{settingcontinuity}, for any index
$\alpha_0\in A$, the pullback functor
$p^*_{\alpha_0}:\M(S_{\alpha_0})\To \M(S)$ preserves fibrations
and trivial fibrations. Moreover, given an index $\alpha_0\in A$, as well as two objects
$C_{\alpha_0}$ and $E_{\alpha_0}$ in $\M(S_{\alpha_0})$, if
$C_{\alpha_0}$ belongs to the smallest full subcategory of $\T(S_{\alpha_0})$
which is closed under finite homotopy colimits and which contains the source
and targets of $I_{\alpha_0}$, then, the canonical map
$$\varinjlim_{\alpha\in A}\Hom_{\ho(\M)(S_\alpha)}(C_\alpha,E_\alpha)
\To \Hom_{\ho(\M)(S)}(C,E)$$
is bijective.
\end{prop}

\begin{proof}
We shall prove first that, for any index $\alpha_0\in A$, the pullback functor
$p^*_{\alpha_0}$ preserves fibrations and trivial fibrations.
By assumption, for any $\alpha\geq \alpha_0$, the pullback functor
along the $\Pmor$-morphism $S_\alpha\To S_{\alpha_0}$ is both a left Quillen functor and
a right Quillen functor. Let $E_{\alpha_0}\To F_{\alpha_0}$
be a trivial fibration (resp. a fibration) of $\M(S_{\alpha_0})$.
Let $i:C\To D$ a generating cofibration (resp. a generating trivial cofibration)
in $\M(S)$. By condition (a) of \ref{settingcontinuity},
we may assume that there exists $\alpha_1\in A$, a
cofibration (resp. a trivial cofibration) $i_{\alpha_1}:C_{\alpha_1}\To D_{\alpha_1}$, such that
$i=p^*_{\alpha_1}(i_{\alpha_1})$. We want to prove that the map
$$\Hom(D,E)\To \Hom(C,E)\times_{\Hom(C,F)}\Hom(D,F)$$
is surjective. But, by condition (b) of \ref{settingcontinuity}, this map
is isomorphic to the filtered colimit of the surjective maps
$$\Hom(D_\alpha,E_\alpha)\To
\Hom(C_\alpha,E_\alpha)\times_{\Hom(C_\alpha,F_\alpha)}\Hom(D_\alpha,F_\alpha)$$
with $\alpha\geq\mathrm{sup}(\alpha_0,\alpha_1)$, which proves the first
assertion. 

To prove the second assertion, we may assume that
$C_{\alpha_0}$ is cofibrant and that $E_{\alpha_0}$ if fibrant.
The set of maps from a cofibrant object to a fibrant object
in the homotopy category of a model category can be described
as homotopy classes of maps. Therefore, using
the fact that $p^*_{\alpha_0}$ preserves cofibrations and fibrations,
 as well as the trivial ones, we see it is sufficient to prove that
the map
$$\varinjlim_{\alpha\in A}\Hom_{\M(S_\alpha)}(C_{\alpha},E_\alpha)
\To \Hom_{\M(S)}(C,E)$$
is bijective for some nice cofibrant replacement of $C_{\alpha_0}$.
But the assumptions on $C_{\alpha_0}$ imply that it is weakly
equivalent to an object $C'_{\alpha_0}$ such that the map
$\varnothing \To C'_{\alpha_0}$ belongs to the smallest class of maps
in $\M(S_{\alpha_0})$, which contains $I_{\alpha_0}$, and which is closed
under pushouts and (finite) compositions. We may thus
assume that $C_{\alpha_0}=C'_{\alpha_0}$. In that case, $C_{\alpha_0}$
is in particular contained in the smallest full subcategory of $\M(S_{\alpha_0})$
which is stable by finite colimits and which contains the
source and targets of $I_{\alpha_0}$. As filtered colimits commute
with finite limits in the category of sets, we conclude by using
again condition (a) of \ref{settingcontinuity}.
\end{proof}

We now go back to the situation of a motivic triangulated category $\T$
 satisfying our general assumptions \ref{assumption_constructible_motives}
 on page \pageref{assumption_constructible_motives}
\begin{lm}\label{prePopescuRHom}
Let $a:X\To Y$ be a morphism in $\sch$. Assume that $X=\varprojlim_\alpha X_\alpha$,
where $(X_\alpha)_{\alpha \in A}$ is a projective system of smooth affine $Y$-schemes.

If $\T$ is $\tau$-continuous, then, for any objects $E$ and $F$ in $\T(Y)$,
with $E$ constructible, then the exchange morphism 
$$a^*\, \sHom_Y(E,F)\simeq \sHom_X(a^*(E),a^*(F)),$$
defined in  Paragraph \ref{num:exchanges4bis}, is an isomorphism.

The same conclusion holds if $\T$ is weakly $\tau$-continuous and
 the transition maps of $(X_\alpha)$ are dominant.
\end{lm}

\begin{proof}
We have
$$a_*\, \sHom_X(a^*(E),a^*(F))\simeq
\sHom_Y(E,a_* \, a^*(F))\, ,$$
so that the map $F\To a_* \, a^*(F))$ induces a map
$$\sHom_Y(E,F)\To a_*\, \sHom_X(a^*(E),a^*(F))\, ,$$
hence, by adjunction, a map
$$a^*\, \sHom_Y(E,F)\To \sHom_X(a^*(E),a^*(F))\, .$$
We already know that the later is an isomorphism whenever $a$ is smooth.

Let us write $a^{}_\alpha:X_\alpha\To Y$ for the structural maps.
Let $C$ be a constructible object in $\T(X)$. By Proposition \ref{continuityconstructible},
we may assume that there exists an index $\alpha_0$, and a constructible
object $C_{\alpha_0}$ in $\T(X_{\alpha_0})$, such that, if we write
$C_\alpha$ for the pullback of $C_{\alpha_0}$ along the map $X_\alpha\To X_{\alpha_0}$
for $\alpha\geq\alpha_0$, we have isomorphisms:
$$\begin{aligned}
\Hom(C,a^*\, \sHom_Y(E,F))
\simeq &\varinjlim_\alpha \Hom(C_\alpha,a_\alpha^*\, \sHom_Y(E,F))\\
\simeq & \varinjlim_\alpha \Hom(C_\alpha,\sHom_X(a_\alpha^*(E),a_\alpha^*(F)))\\
\simeq & \varinjlim_\alpha \Hom(C_\alpha\otimes_{X_\alpha}a_\alpha^*(E),a_\alpha^*(F))\\
\simeq & \Hom(C\otimes_X a^*(E),a^*(F))\\
\simeq & \Hom(C,\sHom_X(a^*(E),a^*(F)))\, .
\end{aligned}$$
As constructible objects generate $\T(X)$, this proves the first assertion.
 The second assertion obviously follows from the same argument.
\end{proof}

\begin{paragr}\label{abstractNisnevichstacks}
Let $X$ be a scheme in $\sch$.
Assume that, for any point $x$ of $X$,
the corresponding morphism $i_x:\spec{\mathscr{O}^h_{X,x}}\To X$ is in $\sch$
(where $\mathscr{O}^h_{X,x}$ denotes the henselization\index{word}{henselization}
 of $\mathscr{O}^{}_{X,x}$).
Consider at last a scheme of finite type $Y$ over $X$, and write
$$a_x:Y_x=\spec{\mathscr{O}^h_{X,x}}\times_X Y\To Y$$
for the morphism obtained by pullback. Finally, for an object $E$ of $\T(Y)$,
let us write $$E_x=a^*_x(E)\, .$$
\end{paragr}

\begin{prop}\label{Nispointscontinuity}
If the motivic category $\T$ is weakly $\tau$-continuous,
 then the family of functors
$$\T(Y)\To \T(Y_x)\ , \quad
E\longmapsto E_x \ , \quad x\in X\, ,$$
is conservative.\index{word}{conservative}
\end{prop}

\begin{proof}
Let $E$ be an object of $\T(Y)$ such that $E_x\simeq 0$ for any
point $x$ of $X$. For any constructible object $C$ of $\T(Y)$, we have
a presheaf of $S^1$-spectra on the small Nisnevich site of $X$:
$$F : U\longmapsto F(U)=\Hom(M_Y(U\times_X Y),\sHom_Y(C,E))\, .$$
It is sufficient to prove that $F(X)$ is acyclic.
As $\T$ satisfies Nisnevich descent (\ref{localizationNisnevichdescent0}),
it is sufficient to prove that $F$ is acyclic locally for the Nisnevich
topology, i.e. that, for any point $x$ of $X$, the spectrum
$F(\spec{\mathscr{O}^h_{X,x}})$ is acyclic.
Writing $\spec{\mathscr{O}^h_{X,x}}$ as the projective limit
of the Nisnevich neighborhoods of $x$ in $X$, we see easily,
using Proposition \ref{continuityconstructible} and Lemma \ref{prePopescuRHom},
that, for any integer $i$,
$\pi_i(F(\spec{\mathscr{O}^h_{X,x}})
\simeq \Hom(C_x,E_x[i])\simeq 0$.
\end{proof}

\begin{prop}\label{henseliancompletioncons}
Let $S$ be a quasi-excellent noetherian and henselian scheme.
Write $\hat S$ for its completion along its closed point, and
assume that both $S$ and $\hat S$ are in $\sch$.
Consider an $S$-scheme of finite type $X$, and write
$i:\hat S\times_S X\To X$ for the induced map.
If $\T$ is $\tau$-continuous, then the pullback functor
$$i^*: \T(X)\To \T(\hat S\times_S X)$$
is conservative.\index{word}{conservative}
\end{prop}

\begin{proof}
As $S$ is quasi-excellent, the map $\hat S\To S$ is
regular. By Popescu's theorem, we can then write
$\hat S=\varprojlim_\alpha S_\alpha$, where $\{S_\alpha\}$
is a projective system of schemes with affine transition maps, and such that
each scheme $S_\alpha$ is smooth over $S$. Moreover, as $\hat S$
and $S$ have the same residue field, and as $S$ is henselian, each map
$S_\alpha$ has a section. Write $X_\alpha=S_\alpha\times_S X$, so
that we have $X=\varprojlim_\alpha X_\alpha$.
Consider a constructible object $C$ and an object $E$ in $\T(X)$.
Then, as the maps $X_\alpha\To X$ have sections, it follows from
the first assertion of Proposition \ref{continuityconstructible} that the map
$$\Hom_{\T(X)}(C,E)\To \Hom_{\T(\hat S\times_S X)}(i^*(C),i^*(E))$$
is a monomorphism (as a filtered colimit of such things).
Hence, if $i^*(E)\simeq 0$, for any constructible object $C$ in $\T(X)$,
we have $\Hom_{\T(X)}(C,E)\simeq 0$. Therefore, as  $\tau$-cons\-truc\-tible objects
generate $\T(X)$, we get $E\simeq 0$.
\end{proof}

\begin{prop}\label{regpullbackvsRHom}
Let $a:X\To Y$ be a regular morphism in $\sch$.
If $\T$ is $\tau$-continuous, then, for any objects $E$ and $F$ in $\T(Y)$,
with $E$ constructible, there is a canonical isomorphism
$$a^*\, \sHom_Y(E,F)\simeq \sHom_X(a^*(E),a^*(F))\, .$$
\end{prop}

\begin{proof}
We want to prove that the canonical map
$$a^*\, \sHom_Y(E,F)\To \sHom_X(a^*(E),a^*(F))$$
is an isomorphism, while we already know it is so whenever $a$ is smooth.
Therefore, to prove the general case, we see that the problem is local
on $X$ and on $Y$ with respect to the Zariski topology.
In particular, we may assume that both $X$ and $Y$ are affine.
By Popescu's Theorem \ref{popescu}, we thus have $X=\varprojlim_\alpha X_\alpha$,
where $\{X_\alpha\}$ is a projective system of smooth affine $Y$-schemes.
We conclude by Lemma \ref{prePopescuRHom}.
\end{proof}

\begin{num}
Consider the following pullback square in $\sch$
$$
\xymatrix{
X'\ar[r]^a\ar[d]_g\ar@{}|\Delta[rd]
 & X\ar[d]^f\\
Y'\ar[r]_b & Y
}
$$
and assume that $f$ is separated of finite type.
 Then one gets, using the recipe that we have seen several times before,
 the following exchange transformation:
$$
Ex(\Delta^{*!}):a^*f^! \xrightarrow{ad(b^*,b_*)}
 a^*f^!b_*b^* \xrightarrow{[Ex(\Delta^!_*)]^{-1}} a^*a_*g^!b^*
 \xrightarrow{ad'(a^*,a_*)} g^!b^*
$$
where $Ex(\delta^!_*)$ is the exchange isomorphism
 of Theorem \ref{thm:cor3_Ayoub}, point (4).
\end{num}
\begin{prop}\label{regpullbacksvsexceptpullbacks}
Consider the previous notations
 and assume that $b$ is regular and $\T$ is $\tau$-continuous.
 Then the exchange transformation defined above
$$
Ex(\Delta^{*!}):a^*f^! \rightarrow g^!b^*
$$
is an isomorphism.
\end{prop}

\begin{proof}
The exchange transformation $Ex(\Delta^{*!})$ is invertible whenever $b$ is smooth: this is obvious
in the case of an open immersion, so that, by Zariski descent, it is
sufficient to treat the case where $b$ is smooth with trivial cotangent
bundle of rank $d$; in this case, by relative purity (\ref{thm:cor3_Ayoub}~(3)),
this reduces to the canonical isomorphism $a^!f^!\simeq g^!b^!$
evaluated at $E(-d)[-2d]$. To prove the general case,
as the condition is local on $X$ and on $Y$ for the
Zariski topology, we may assume that $f$ factors as an immersion
$X\To \PP^n_Y$, followed by the canonical projection $\PP^n_Y\To Y$.
We deduce from there that it is sufficient to treat the case where $f$
is either a closed immersion, either a smooth morphism of finite type.
The case where $f$ (hence also $g$) is smooth follows by relative purity (\ref{thm:cor3_Ayoub}):
we can then replace $f^!$ and $g^!$ by $f^*$ and $g^*$
respectively, and the formula follows from
the fact that $a^*f^*\simeq g^*b^*$. We may thus assume that $f$
is a closed immersion. As $g$ is a closed immersion as well,
the functor $g_!$ is conservative (it is fully faithful).
Therefore, it is sufficient to prove that the map
$$b^*\, f_! \, f^!(E)\simeq g_!\, a^*\, f^!(E)\To g_!\, g^!\, b^*(E)$$
is invertible. Then, using Proposition \ref{regpullbackvsRHom}
(which makes sense because the functor $f_!$ preserves $\tau$-constructibility by \ref{thmfinitnessproper}),
and the projection formula, we have
$$\begin{aligned}
b^*\, f_! \, f^!(E)
& \simeq b^*\sHom_Y(f_!(\unit_X),E)\\
& \simeq \sHom_{Y'}(b^* \, f_!(\unit_X),b^*(E))\\
& \simeq \sHom_{Y'}(g_!(\unit_{X'}),b^*(E))\\
& \simeq g_!\, g^!\, b^*(E)\, ,
\end{aligned}$$
which achieves the proof.
\end{proof}

\begin{lm}\label{bigdirectimages}
Let $f:X\To Y$ be a morphism in $\sch$.
Assume that $X=\varprojlim_\alpha X_\alpha$
and $Y=\varprojlim_\alpha Y_\alpha$, where $\{X_\alpha\}$
and $\{Y_\alpha\}$ are projective systems of schemes with
affine (resp. affine and dominant) transition maps, while $f$ is induced by a system
of morphisms $f_\alpha:X_\alpha\To Y_\alpha$.
Let $\alpha_0$ be some index, $C_{\alpha_0}$ a
constructible object of $\T(Y_{\alpha_0})$, and
$E_{\alpha_0}$ an object of $\T(X_{\alpha_0})$.
If $\T$ is $\tau$-continuous (resp. weakly $\tau$-continuous),
then we have a natural isomorphism
of abelian groups
$$\varinjlim_{\alpha\geq\alpha_0}
\Hom_{\T(Y_\alpha)}(C_\alpha,f_{\alpha,*}(E_\alpha))
\simeq \Hom_{\T(Y)}(C,f_*(E))\, .$$
\end{lm}

\begin{proof}
By virtue of Proposition \ref{continuityconstructible}, we have a natural isomorphism
$$\varinjlim_{\alpha\geq\alpha_0}
\Hom_{\T(X_\alpha)}(f^*_\alpha (C_\alpha), E_\alpha)
\simeq \Hom_{\T(Y)}(f^* (C),E)\, .$$
The expected formula follows by adjunction.
\end{proof}

\begin{prop}\label{regularbasechange}
Consider the following pullback square in $\sch$.
$$\xymatrix{
X'\ar[r]^a\ar[d]_g&X\ar[d]^f\\
Y'\ar[r]_b&Y
}$$
with $b$ regular. If $\T$ is $\tau$-continuous, then, for any object
$E$ in $\T(X)$, there is a canonical isomorphism
in $\T(Y')$:
$$b^*\, f_*(E)\simeq g_*\, a^*(E)\, .$$
\end{prop}

\begin{proof}
This proposition is true in the case where $b$ is smooth
(by definition of $\sm$-fibred categories), from which we deduce, by
Zariski separation, that this property is local on $Y$ and on $Y'$
for the Zariski topology. In particular, we may assume that both $Y$
and $Y'$ are affine. Then, by Popescu's Theorem \ref{popescu},
we may assume that $Y'=\varprojlim_\alpha Y'_\alpha$, where
$\{Y'_\alpha\}$ is a projective system of smooth $Y$-algebras. Then,
using the preceding lemma as well as Proposition \ref{continuityconstructible},
we reduce easily the proposition to the case where $b$ is smooth.
\end{proof}

\begin{prop}\label{insepclosure}
Assume that $\T$ is weakly $\tau$-continuous,
 $\QQ$-linear and semi-separated,
and consider a field $k$, with inseparable closure $k'$,
such that both $\spec k$ and $\spec{k'}$ are in $\sch$.
Given a $k$-scheme $X$ write $X'=k'\otimes_k X$, and
$f:X'\To X$ for the canonical projection.
Then the functor
$$f^* : \T(X)\To \T(X')$$
is an equivalence of categories.
\end{prop}

\begin{proof}
Note that $X'$ is a projective limit of $k$-schemes with affine
 and dominant (even flat) transition maps.
Thus, it follows from weak $\tau$-continuity,
 Proposition \ref{continuityconstructible}
 and Proposition \ref{stronglyquasiseparated} that the functor
$$f^* : \T_c(X)\To \T_c(X')$$
is an equivalence of categories. Similarly, for any objects $C$
and $E$ in $\T(X)$, if $C$ is  constructible, the map
$$\Hom_{\T(X)}(C,E)\To
\Hom_{\T(X)}(f^*(C),f^*(E))$$
is bijective. As constructible objects generate $\T(X)$, this
implies that the functor
$$f^* : \T(X)\To \T(X')$$
is fully faithful. As the latter is essentially surjective on a set of generators,
this implies that it is an equivalence of categories (see \ref{equivgeneratorscompactgentriang}).
\end{proof}

\begin{prop}\label{prop:pointsconservative}
Assume that $\T$ is weakly $\tau$-continuous.
Then, for any scheme $X$ in $\sch$, the family of functors
induced by its points
$$x^*:\T(X)\To\T(\spec{\kappa(x)}\ , \quad x\in X\, ,$$
is conservative.
\end{prop}
\begin{proof}
We proceed by induction on the dimension $d$ of $X$. If $d\leq 0$, this is trivial.
If $d>0$, using Proposition \ref{Nispointscontinuity}, we may assume
that $X$ is local. By induction, the proposition is true on the
complement of the closed point of $x$. Therefore, Proposition \ref{prop:easy_csq_loc}
achieves the proof.
\end{proof}

%
%Here is a slightly more general version of Proposition \ref{DMBseparated}.
%
%\begin{prop}\label{DMBseparated2}
%The motivic category $\DMB$ is separated on noetherian schemes of finite
%dimension.
%\end{prop}
%
%\begin{proof}
%As in the proof of \ref{DMBseparated}, it is sufficient to prove that,
%given a finite surjective morphism $f:T\To S$, the pullback functor
%$f^*:\DMB(S)\To \DMB(T)$ is conservative.
%By virtue of Proposition \ref{Nispointscontinuity}, we may assume that $S$
%is henselian. Using the localization property, we may even assume
%(by induction on the dimension) that $S$ is the spectrum of a field.
%Replacing $T$ by its reduction, we may thus assume that both $S$ and $T$
%are regular. We can then conclude by a trace argument, as in the proof
%of Proposition \ref{DMBseparated}.
%\end{proof}
%
%\begin{cor}\label{DMBetaledescent}
%The motivic category $\DMB$ satisfies \'etale descent on
%noetherian schemes of finite dimension.
%\end{cor}
%
%\begin{proof}
%This follows from the preceding proposition and from
%Theorem \ref{fibredetaledescent}.
%\end{proof}

\subsection{Duality} \label{sec:motivic_duality}

%We fix a noetherian base scheme $\base$ of finite dimension,
%as well as an adequate category of schemes $\sch$. We shall also assume
%that all the schemes in $\sch$ are quasi-excellent.
%
%We consider given a stable combinatorial  $\sm$-fibred model
%category $\M$ over $\sch$,
%such that $\T$ is motivic, and
%is endowed with a fixed set of twists $\tau$, which is
%assumed to be stable under negative Tate twists.
The aim of this section is to prove a local duality theorem in $\T$
(see \ref{thm:localduality} and \ref{cor:localduality2}).
%%under the hypothesis
%%that singularities of $B$ can be solved up to quotient singularities (see \ref{df:ressingupquotient}).

If we work with rational coefficients, resolution of
singularities up to quotient singularities is almost as good
as classical resolution of singularities:
we have the following replacement of the blow-up formula.

\begin{thm}\label{excisionGaloisalteration}
Assume that $\T$ is $\QQ$-linear and separated.
Let $X$ be a scheme in $\sch$. Consider a proper surjective morphism
$p:X'\To X$ and a finite group $G$ acting on $X'$ over $X$.
Assume that there is a closed subscheme $Z\subset X$ such that
$U=X-Z$ is normal, while the induced map $p_U:U'=p^{-1}(U)\To U$ is finite,
and the map $U'/G\To U$ is generically radicial (i.e. is radicial over an open dense subscheme of $U$) --- e.g. this situation occurs when $p$ is a Galois alteration\index{word}{alteration!Galois alteration.}.
Then the pullback square
\begin{equation}\label{excisionGaloisalteration1}\begin{split}
\xymatrix{
Z'\ar[r]^{i'}\ar[d]_q&X'\ar[d]^p\\
Z\ar[r]^i&X
}\end{split}\end{equation}
induces an homotopy pullback\index{word}{homotopy cartesian} square
\begin{equation}\label{excisionGaloisalteration2}
\begin{split}
\xymatrix{
M\ar[r]\ar[d]&(\derR p_*\, \derL p^*(M))^{G}\ar[d]\\
\derR i_*\, \derL i^*(M)\ar[r]&(\derR i_*\derR q_*\, \derL q^* \, \derL i^*(M))^{G}
}
\end{split}
\end{equation}
for any object $M$ of $\T(X)$.
\end{thm}

\begin{proof}
We already know that, for any object $N$ of $\T(U)$, the map
$$N\To (\derR p^{}_{U*}\, \derL p^*_U(N))^G$$
is an isomorphism (Corollary \ref{stronggenericGaloisdescent}).
The proof is then similar to the proof of condition (iv) of Theorem \ref{charseparated3}.
\end{proof}

\begin{rem}\label{exceptionnalexcisionGaloisalteration}
Under the assumptions of the preceding theorem, applying the total derived functor
$\derR\Hom_X(-,E)$ to the homotopy pullback square \eqref{excisionGaloisalteration2}
for $M=\unit_X$, we obtain the homotopy pushout square
\begin{equation}\label{excisionGaloisalteration3}
\begin{split}
\xymatrix{
(i_!\, q_!\, q^! \, i^!(E))_{G}\ar[r]\ar[d]&(p_!\, p^!(E))_{G}\ar[d]\\
i_!\, i^!(E)\ar[r]&E
}
\end{split}
\end{equation}
for any object $E$ of $\T(X)$ .
\end{rem}

\begin{cor}\label{properregulargenerators}
Assume that $\T$ is $\QQ$-linear and separated.
Let $B$ be a scheme in $\sch$, admitting wide resolution of singularities up to quotient singularities.
Consider a separated $B$-scheme of finite type $S$, endowed with a
closed subscheme $T\subset S$. The category of constructible
objects in $\T(S)$ is the smallest thick triangulated subcategory which contains
the objects of shape $\derR f_*(\unit_X\{n\})$  for $n\in\tau$, and for $f:X\To S$
a projective morphism, with $X$ regular and connected, such that
$f^{-1}(T)_{\mathit{red}}$ is either empty,
either $X$ itself or the support of a strict normal crossing divisor.
\end{cor}

\begin{proof}
Let $\T(S)'$ be the smallest thick triangulated subcategory of $\T(S)$ which contains
the objects of shape $\derR f_*(\unit_X\{n\})$ for $n\in\tau$ and $f:X\To S$
a projective morphism with $X$ regular and connected, while
$f^{-1}(T)_{\mathit{red}}$ is empty, or $X$ itself,
or the support of a strict normal crossing divisor.
We clearly have $\T(S)'\subset\T_c(S)$ (Proposition \ref{thmfinitnessproper}).
To prove the reverse inclusion, by virtue of Proposition \ref{thmfinitnessproper2},
it is sufficient to prove that,
for any $n\in\tau$, and any projective morphism $f:X\To S$, the object
$\derR f_*(\unit_X\{n\})$ belongs to $\T(S)'$.
We shall proceed by induction on the dimension of $X$.
If $X$ is of dimension $\leq 0$, we may replace it by its reduction, which
is regular. If $X$ is of dimension $>0$,
by assumption on $B$, there exists
a Galois alteration $p:X'\To X$ of group $G$, with $X'$ regular
and projective over $S$ (and in which $T$ becomes either empty,
either $X'$ itself, either the support of
a strict normal crossing divisor, in each connected component of $X'$).
Choose a closed subscheme $Z\subset X$, such that $U=X-Z$ is
a normal dense open subscheme, and
such that the induced map $r:U'=p^{-1}(U)\To U$ is a finite morphism, and consider
the pullback square \eqref{excisionGaloisalteration1}.
As $Z$ and $Z'=p^{-1}(Z)$ are of dimension smaller than the dimension of $X$,
we conclude from the homotopy pullback square obtained by applying
the functor $\derR f_*$ to \eqref{excisionGaloisalteration2}
for $M=\unit_X\{n\}$, $n\in \tau$.
\end{proof}

\begin{df}
Let $S$ be a scheme in $\sch$.
An object $R$ of $\T(S)$ is \emph{$\tau$-dualizing}
\index{word}{dualizing!$\tau$-dualizing}
 if it satisfies the following conditions.
\begin{itemize}
\item[(i)] The object $R$ is constructible.
\item[(ii)] For any constructible object $M$ of $\T(S)$, the
natural map
$$M\To \derR\sHom_S(\derR\sHom_S(M,R),R)$$
is an isomorphism.
\end{itemize}
\end{df}

\begin{rem}
If $\T$ is $\tau$-compatible, $\QQ$-linear and separated,
then, in particular, the six operations of Grothendieck preserve
$\tau$-constructibility in $\T$ (\ref{grothendieck6op}).
Under this assumption, for any scheme $X$ in $\sch$, and any
$\otimes$-invertible object $U$ in $\T(X)$ which is constructible,
its quasi-inverse is constructible: the quasi-inverse of $U$ is
simply its dual $U^\wedge=\derR\sHom(U,\unit_X)$,
which is constructible by virtue of \ref{constructibleinternalHom}.
\end{rem}

\begin{prop}\label{unicitydualizing}
Assume that $\T$ is $\tau$-compatible, $\QQ$-linear and separated,
and consider a scheme $X$ in $\sch$.
\begin{itemize}
\item[(i)] Let $R$ be a $\tau$-dualizing object, and $U$ be a constructible
$\otimes$-invertible object in $\T(X)$. Then $U\otimes^\derL_S R$ is $\tau$-dualizing.
\item[(ii)] Let $R$ and $R'$ be two $\tau$-dualizing objects in $\T(X)$.
Then the evaluation map
$$\derR \sHom_S (R,R')\otimes^\derL_S R\To R'$$
is an isomorphism.
\end{itemize}
\end{prop}

\begin{proof}
This follows immediately from \cite[2.1.139]{ayoub}.
\end{proof}

\begin{prop}\label{dualizingopenimmersions}
Consider an open immersion $j:U\To X$ in $\sch$.
If $R$ is a $\tau$-dualizing object in $\T(X)$, then $j^!(R)$
is $\tau$-dualizing in $\T(U)$.
\end{prop}

\begin{proof}
If $M$ is a constructible object in $\T(U)$, then
$j_!(M)$ is constructible, and the map
\begin{equation}\label{dualizingopenimmersions01}
j_!(M)\To\derR\sHom_X(\derR\sHom_X(j_!(M),R),R)
\end{equation}
is an isomorphism. Using the isomorphisms of type
$$M\simeq j^*\, j_!(M)=j^!\, j_!(M)\quad\text{and}\quad
j^*\derR\sHom_X(A,B)\simeq\derR\sHom_U(j^*(A),j^*(B))\, ,$$
we see that the image of the map \eqref{dualizingopenimmersions01}
by the functor $j^*=j^!$ is isomorphic to the map
\begin{equation}\label{dualizingopenimmersions02}
M\To\derR\sHom_U(\derR\sHom_U(M,j^!(R)),j^!(R))\, ,
\end{equation}
which proves the proposition.
\end{proof}

\begin{prop}\label{dualizingpptyislocal}
Let $X$ be a scheme in $\sch$, and $R$
an object in $\T(X)$.
Assume there exists an open cover $X=\bigcup_{i\in I}U_i$
such that the restriction of $R$ on each of the open subschemes $U_i$
is $\tau$-dualizing in $\T(U_i)$. Then $R$ is $\tau$-dualizing.
\end{prop}

\begin{proof}
We already know that the property of $\tau$-constructibility
is local with respect to the Zariski topology (\ref{thmfinitnessredlocZar0}).
Denote by $j_i:U_i\To X$ the corresponding open immersions, and put $R_i=j^!_i(R)$.
Let $M$ be a constructible object in $\T(X)$.
Then, for all $i\in I$, the image by $j^*_i=j^!_i$ of the map
$$M\To\derR\sHom_X(\derR\sHom_X(M,R),R)$$
is isomorphic to the map
$$j^*_i(M)\To\derR\sHom_{U_i}(\derR\sHom_{U_i}(j^*_i(M),R_i),R_i)\, .$$
This proposition thus follows from the property of separation with respect to the
Zariski topology.
\end{proof}

\begin{cor}\label{cor:dualizingpptyislocal}
Let $f:X\To Y$ be a separated morphism of finite type in $\sch$.
Given an object $R$ of $\T(Y)$, the property for $f^!(R)$ of being
a $\tau$-dualizing object in $\T(X)$ is local over $X$ and over $Y$ for the
Zariski topology.
\end{cor}

\begin{prop}\label{intimmersions}
Assume that $\T$ is $\tau$-compatible.
Let $i:Z\To X$ be a closed immersion and $R$ be a $\tau$-dualizing object in $\T(X)$.
Then $i^!(R)$ is $\tau$-dualizing in $\T(Z)$.
\end{prop}
 
\begin{proof}
As $\T$ is $\tau$-compatible, we already know that $i^!(R)$
is constructible. For any objects $M$ and $R$ of $\T(Z)$ and $\T(X)$
respectively, we have the identification:
\begin{equation*}
i_! \, \derR\sHom_Z(M,i^!(R))
\simeq \derR\sHom_X(i_!(M),R)\, .
\end{equation*}
Let $j:U\To X$ be the complement immersion.
Then we have
$$j^!\derR\sHom_X(i_!(M),R)
\simeq \derR\sHom_U(j^* \, i_!(M),j^!(R))\simeq 0\, ,$$
so that
$$\derR\sHom_X(i_!(M),R)\simeq
i_! \, \derL i^* \derR\sHom_X(i_!(M),R)\, .$$
As $i_!$ is fully faithful, this provides a canonical isomorphism
$$\derL i^* \derR\sHom_X(i_!(M),R)\simeq i^!\derR\sHom_X,(i_!(M),R)\, .$$
Under this identification, we see easily that the map
\begin{equation*}
i_!(M)\To \derR\sHom_{X}(\derR\sHom_{X}(i_!(M),R),R)
\end{equation*}
is isomorphic to the image by $i_!$ of the map
\begin{equation*}
M\To \derR\sHom_{Z}(\derR\sHom_{Z}(M,i^!(R)),i^!(R))\, .
\end{equation*}
As $i_!$ is fully faithful, it is conservative, and this ends the proof.
\end{proof}

\begin{prop}\label{prop:localduality}
Assume that $\T$ is $\tau$-compatible,
$\QQ$-linear and separated, and consider a scheme
$B$ in $\sch$ which admits wide resolution of singularities up to quotient singularities.
Consider a separated $B$-scheme of finite type $S$, and a constructible object
$R$ in $\T(S)$. The following conditions are equivalent.
\begin{itemize}
\item[(i)] For any separated morphism of finite type $f:X\To S$,
the object $f^!(R)$ is $\tau$-dualizing.
\item[(ii)] For any  projective morphism $f:X\To S$,
the object $f^!(R)$ is $\tau$-dualizing.
\item[(iii)] For any  projective morphism $f:X\To S$, with $X$
regular, the object $f^!(R)$ is $\tau$-dualizing.
\item[(iv)] For any projective morphism $f:X\To S$, with $X$
regular, and for any $n\in \tau$, the map
\begin{equation}\label{prop:localduality1}
\unit_X\{n\} \To \derR\sHom_X(\derR\sHom_X(\unit_X\{n\},f^!(R)),f^!(R))
\end{equation}
is an isomorphism in $\T(X)$.
\end{itemize}
If, furthermore, for any regular separated $B$-scheme of finite type $X$,
and for any $n\in\tau$, the object $\unit_X\{n\}$ is $\otimes$-invertible, then
these conditions are equivalent to the following one.
\begin{itemize}
\item[(v)] For any projective morphism $f:X\To S$, with $X$ regular, the map
\begin{equation}\label{prop:localduality1bis}
\unit_X \To \derR\sHom_X(f^!(R),f^!(R))
\end{equation}
is an isomorphism in $\T(X)$.
\end{itemize}
\end{prop}

\begin{proof}
It is clear that (i) implies (ii), which implies (iii), which implies (iv).
Let us check that condition (ii) also implies condition (i).
Let $f:X\To S$ be a morphism of separated $B$-schemes of finite type, with $S$
regular. We want to prove that $f^!(\unit_S)$ is $\tau$-dualizing,
while we already know it is true whenever $f$ is projective.
In the general case, by virtue of Corollary \ref{cor:dualizingpptyislocal}, we may assume
that $f$ is quasi-projective, so that $f=pj$, where $p$ is projective, and
$j$ is an open immersion. As $f^!\simeq j^!\, p^!$,
we conclude with Proposition \ref{dualizingopenimmersions}.
Under the additional assumption, the equivalence between (iv) and (v) is obvious.
It thus remains to prove that (iv) implies (ii).
It is in fact sufficient to prove that, under condition (iv),
the object $R$ itself is $\tau$-dualizing.
To prove that the map
\begin{equation}\label{prop:localduality2}
M\To \derR\sHom_X(\derR\sHom_X(M,R),R)
\end{equation}
is an isomorphism for any constructible object $M$ of $\T(S)$,
it is sufficient to consider the case where $M=\derR f_*(\unit_X\{n\})=f_!(\unit_X\{n\})$,
where $n\in\tau$ and $f:X\To S$ is a projective morphism with $X$ regular
(Corollary \ref{properregulargenerators}).
For any object $A$ of $\T(X)$, we have canonical isomorphisms
\begin{align*}
\derR\sHom_S(f_!(A),R)
& \simeq \derR f_*\, \derR \sHom_X(A,f^!(R))\\
& = f_!\, \derR \sHom_X(A,f^!(R))\, ,
\end{align*}
from which we get a natural isomorphism:
$$\derR\sHom_S(\derR\sHom_S(f_!(A),R),R)\simeq
f_!\, \derR\sHom_X(\derR\sHom_X(A,f^!(R)),f^!(R))\, .$$
Under these identifications, the map \eqref{prop:localduality2}
for $M=f_!(\unit_X\{n\})$ is the image of the map \eqref{prop:localduality1}
by the functor $f_!$. As \eqref{prop:localduality1} is invertible by assumption,
this proves that $R$ is $\tau$-dualizing.
\end{proof}

\begin{lm}\label{invertiblepptyislocal}
Let $X$ be a scheme in $\sch$, and $R$ be an object of $\T(X)$.
The property for $R$ of being $\otimes$-invertible is local over $X$ with respect to the Zariski
topology.
\end{lm}

\begin{proof}
Let $R^\wedge=\derR\sHom(R,\unit_X)$ be the dual of $R$.
The object $R$ is $\otimes$-invertible if and only if the evaluation map
$$R^\wedge\otimes^\derL_X R\To \unit_X$$
is invertible. Let $j:U\To X$ be an open immersion.
Then, for any objects $M$ and $N$ in $\T(X)$, we have
the identification
$$j^*\derR\sHom_X(M,N)\simeq \derR \sHom_U(j^*(M),j^*(N))\, .$$
In particular, we have $j^*(R^\wedge)\simeq j^*(R)^\wedge$.
As $j^*$ is monoidal, the lemma follows from the fact that $\T$
has the property of separation with respect to the Zariski topology. 
\end{proof}

\begin{df}\label{deftaupure}
We shall say that $\T$ is \emph{$\tau$-dualizable}
 if it satisfies the following conditions:
\begin{itemize}
\item[(i)] $\T$ is $\tau$-compatible (\ref{df:weaklytaupure});
\item[(ii)] for any closed immersion between regular schemes
$i:Z\To S$ in $\sch$, the object $i^!(\unit_S)$
is $\otimes$-invertible (i.e. the functor $i^!(\unit_S)\otimes^\derL_S(-)$
is an equivalence of categories);
\item[(ii)] for any regular scheme $X$ in $\sch$,
and for any $n\in\tau$, the map
$$\unit_X\{n\}\To \derR\sHom_X(\derR\sHom_X(\unit_X\{n\},\unit_X),\unit_X)$$
is an isomorphism.
\end{itemize}
\end{df}
As in other similar situations,
 we simply say \emph{dualizable with respect to Tate twist}
 when the set of twists $\tau$ is generated by the Tate twist.

\begin{ex}
In practice, the property of being dualizable with respect
 to Tate twist is a consequence of the absolute purity theorem.
 Our main example is the motivic category $\DMB$ of Beilinson
  motives over excellent noetherian schemes, as a consequence
  of Theorem \ref{DMBpurity}.
\end{ex}

\begin{rem}\label{trivialtaupurityreductions}
Note that, whenever the set of twists $\tau$ consists of rigid objects
(which will be the case in practice),
conditions (i) and (ii) of the preceding definition are equivalent to
the condition that $i^!(\unit_X)$ is constructible and $\otimes$-invertible for any closed immersion $i$
between regular separated schemes in $\sch$, while condition (iii) is then automatic.
This principle gives easily the property of $\tau$-purity when $\sch$ is made
of schemes of finite type over some perfect field:
\end{rem}

\begin{prop}\label{purityoverfields}
Assume that $\sch$ consists exactly of schemes of finite type over
a field $k$.
If the objects $\unit\{n\}$ are rigid with constructible duals
in $\T(\spec k)$ for all $n\in\tau$,
then $\T$ is $\tau$-dualizable.
\end{prop}

\begin{proof}
For any $k$-scheme of finite type $f:X\To \spec k$,
as the functor $\derL f^*$ is symmetric monoidal, the
objects $\unit_X\{n\}$ are rigid in $\T(X)$ for all $n\in\tau$.
Therefore, as stated in remark \ref{trivialtaupurityreductions}, we have only to prove that,
for any closed immersion $i:Z\To X$ between regular $k$-schemes
of finite type, the object $i^!(\unit_X)$ is $\otimes$-invertible and constructible.
Note that, as $k$ is perfect, $X$ and $Z$ are in fact smooth.
Using \ref{invertiblepptyislocal} and \ref{thmfinitnessredlocZar0},
we may also assume that $X$ is quasi-projective and that $Z$ is purely
of codimension $c$ in $X$, while the normal bundle of $i$
is trivial. This proposition is then a consequence
of relative purity (\ref{thm:cor3_Ayoub}), which gives a canonical isomorphism
$i^!(\unit_X)\simeq\unit_Z(-c)[-2c]$.
\end{proof}

\begin{prop}\label{poincarerigidfield}
Assume that $\site$ consists of schemes
of finite type over a field $k$ and that $\T$ has the following properties:
\begin{itemize}
\item[(a)] it is $\tau$-dualizable;
\item[(b)] for any $n\in\tau$, $\unit\{n\}$ is rigid;
\item[(c)] either $k$ is  perfect, either $\T$
is continuous.
\end{itemize}
Then, any constructible object of $\T(k)$ is rigid.
\end{prop}

\begin{proof}
By \ref{insepclosure}, it is sufficient to treat the case where $k$ is perfect.
It is well-known that rigid objects form a thick subcategory of $\T$.
Thus, we conclude easily from Corollary \ref{properregulargenerators}
and Proposition \ref{prop:purity&duality}.
\end{proof}

\begin{lm}\label{iminvexcepinvertible}
Assume that $\T$ is $\tau$-dualizable. Then, for any
projective morphism $f:X\To S$ between
regular schemes in $\sch$, the object $f^!(\unit_S)$
is $\otimes$-invertible and constructible.
\end{lm}

\begin{proof}
As, for any open immersion $j:U\To X$, one has $j^*=j^!$,
we deduce easily from Lemma \ref{invertiblepptyislocal}
(resp. Proposition \ref{thmfinitnessredlocZar0})
that the property for $f^!(\unit_S)$ of being $\otimes$-invertible
(resp. constructible) is local on $S$ for the Zariski topology.
Therefore, we may assume that $S$ is separated over $B$ and that
$f$ factors as a closed immersion $i:X\To \PP^n_S$
followed by the canonical projection $p:\PP^n_S\To S$.
Using relative purity for $p$, we have the following computations:
$$f^!(\unit_S)\simeq i^!\, p^!(\unit_S)\simeq i^!(\unit_{\PP^n_S}(n)[2n])
\simeq i^!(\unit_{\PP^n_S})(n)[2n]\, .$$
As $i$ is a closed immersion between regular schemes,
the object $i^!(\unit_{\PP^n_S})$ is $\otimes$-invertible
and constructible by assumption on $\T$,
which implies that $f^!(\unit_S)$ is $\otimes$-invertible
and constructible as well.
%%The property of $\tau$-constructibility is ensured by \ref{constructexcepinvim}.
\end{proof}

\begin{df}\label{def:localduality}
Let $B$ a scheme in $\sch$.
We shall say that \emph{local duality\index{word}{duality!local duality}
 holds over $B$ in $\T$}
if, for any separated morphism of finite type $f:X\To S$, with
$S$ regular and of finite type over $B$, the object $f^!(\unit_S)$ is $\tau$-dualizing
in $\T(X)$.
\end{df}

\begin{rem}
By definition, if $\T$ is $\tau$-compatible, and if
local duality holds over $B$ in $\T$, then the restriction of $\T$
to the category of $B$-schemes of finite type is $\tau$-dualizable.
A convenient sufficient condition for local duality to hold in $\T$ is the following
(in particular, using the result below as well as Proposition \ref{purityoverfields},
local duality holds almost systematically over fields).
\end{rem}

\begin{thm}\label{thm:localduality}
Assume that $\T$ is $\tau$-dualizable,
 $\QQ$-linear and separated, and consider a scheme
$B$ in $\sch$ which admits wide resolution of singularities up to quotient singularities
(e.g. $B$ might be any scheme which is separated and of finite type over an excellent noetherian scheme of
dimension lesser or equal to $2$ in $\sch$; see \ref{cor:dejongdimleq2}).
Then local duality holds over $B$ in $\T$.
\end{thm}

\begin{proof}
Let $S$ be a regular separated $B$-scheme of finite type.
Then, for any separated morphism of finite type
$f:X\To S$, the object $f^!(\unit_S)$ is $\tau$-dualizing:
Lemma \ref{iminvexcepinvertible}
implies immediately condition (iv) of Proposition \ref{prop:localduality}.
The general case (without the separation assumption on $S$)
follows easily from Corollary \ref{dualizingpptyislocal}.
\end{proof}

\begin{prop}\label{prop:axiomaticlocalduality}
Consider a scheme $B$ in $\sch$. Assume that $\T$ is $\tau$-dualizable, and
that local duality holds over $B$ in $\T$.
Consider a regular $B$-scheme of finite type $S$.
\begin{itemize}
\item[(i)] An object of $\T(S)$ is $\tau$-dualizing if and only if it is
constructible and $\otimes$-invertible.
\item[(ii)] For any separated morphism of $S$-schemes of finite
type $f:X\To Y$, and for any $\tau$-dualizing object $R$ in $\T(Y)$,
the object $f^!(R)$ is $\tau$-dualizing in $\T(X)$.
\end{itemize}
\end{prop}

\begin{proof}
As the unit of $\T(S)$ is $\tau$-dualizing by assumption, Proposition \ref{unicitydualizing}
implies that an object of $\T(S)$ is $\tau$-dualizing if and only if it is constructible
and $\otimes$-invertible.

Consider a regular $B$-scheme of finite type $S$, as well as a separated morphism of $S$-schemes
of finite type $f:X\To Y$, as well as a $\tau$-dualizing object $R$ in $\T(Y)$.
To prove that $f^!(R)$ is $\tau$-dualizing, by virtue of Corollary \ref{dualizingpptyislocal},
we may assume that $Y$ is separated over $S$.
Denote by $u$ and $v$ the structural maps from $X$ and $Y$ to $S$ respectively.
As we already know that $v^!(\unit_S)$ is $\tau$-dualizing,
by virtue of Proposition \ref{unicitydualizing}, there exists a constructible
and $\otimes$-invertible object $U$ in $\T(Y)$ such that $U\otimes^\derL_Y R\simeq v^!(\unit_S)$.
As the functor $\derL f^*$ is symmetric monoidal, it preserves $\otimes$-invertible objects
and their duals, from which we deduce the following isomorphisms:
\begin{align*}
u^!(\unit_S)& \simeq f^!\, v^! (\unit_S) \\
& \simeq f^!(U\otimes^\derL_Y R)\\
& \simeq f^!\, \derR\sHom_Y(U^\wedge,R)\\
& \simeq \derR\sHom_X(\derL f^*(U^\wedge),f^!(R))\\
& \simeq \derR\sHom_X(\derL f^*(U)^\wedge,f^!(R))\\
& \simeq \derL f^*(U)\otimes^\derL_X f^!(R) \, .
\end{align*}
The object $a^!(\unit_S)$ being $\tau$-dualizing, while $\derL f^*(U)$
is constructible and invertible, we deduce from Proposition \ref{unicitydualizing}
that $f^!(R)$ is $\tau$-dualizing as well.
\end{proof}

\begin{paragr}\label{deflocaldualityfunctors}
Assume that $\T$ is $\tau$-dualizable, $\QQ$-linear and separated.

Consider a scheme $B$ in $\sch$, 
 such that local duality holds over $B$ in $\T$
 --- this is the case if $B$ admits wide resolution 
 of singularities up to quotient singularities according to
 the above Theorem.
Consider a fixed regular $B$-scheme of finite type $S$, as well
as a constructible and $\otimes$-invertible object $R$
in $\T(S)$ (in the case $S$ is of pure dimension $d$, it might
be wise to consider $R=\unit_S(d)[2d]$, but an arbitrary $R$
as above is eligible by \ref{prop:axiomaticlocalduality}).
Then, for any separated $S$-scheme of finite type $f:X\To S$,
we define the \emph{local duality functor}
$$D_X: \op{\T(X)}\To\T(X)$$
by the formula
$$D_X(M)=\derR\sHom_X(M,f^!(R))\, .$$
This functor $D_X$ is right adjoint to itself.
\end{paragr}
\begin{cor}\label{cor:localduality2}
Under the above assumptions,
we have the following properties of the motivic triangulated
 category $\T$:
\begin{itemize}
\item[(a)] For any separated $S$-scheme of finite type $X$,
the functor $D_X$ preserves constructible objects.
\item[(b)] For any separated $S$-scheme of finite type $X$, the natural map
$$M\To D_X(D_X(M))$$
is an isomorphism for any constructible object $M$ in $\T(X)$.
\item[(c)] For any separated $S$-scheme of finite type $X$, and
for any objects $M$ and $N$ in $\T(X)$, if $N$ is constructible,
then we have a canonical isomorphism
$$D_X(M\otimes^\derL_X D_X(N))\simeq \derR\sHom_X(M,N)\, . $$
\item[(d)] For any morphism between separated $S$-schemes of finite type
$f:Y\To X$, we have natural isomorphisms
\begin{align*}
D_Y(f^*(M))& \simeq f^!(D_X(M))\\
f^*(D_X(M))& \simeq D_Y(f^!(M))\\
D_X(f_!(N))& \simeq f_*(D_Y(N))\\
f_!(D_Y(N))& \simeq D_X(f_*(N))\\
\end{align*}
for any constructible objects $M$ and $N$
in $\T(X)$ and $\T(Y)$ respectively.
\end{itemize}
\end{cor}
This corollary sums up what must be called the \emph{Grothendieck duality}
\index{word}{duality, Grothendieck}
 property for the motivic triangulated category $\T$
 with respect to the set of twists $\tau$.
\begin{proof}
Assertions (a) and (b) are only stated for the record\footnote{We have
put to a lot of assumptions here: in fact, if $\T$
is $\tau$-dualizable and if local duality holds over $B$ in $\T$, the six
Grothendieck operations preserve constructible objects
on the restriction of $\T$ to $B$-schemes of finite type; we leave this
as a formal exercise for the reader.}; see \ref{constructibleinternalHom}.
To prove (c), we see that we have an obvious isomorphism
$$D_X(M\otimes^\derL_X P)\simeq \derR\sHom_X(M,D_X(P))$$
for any objects $M$ and $P$. If $N$ is constructible,
we may replace $P$ by $D_X(N)$ and get the expected formula
using (b).
The identification $D_Y\, f^*\simeq f^!\, D_X$ is a special case of the formula
$$\derR\sHom_Y(f^*(A),f^!(B)) \simeq f^!\, \derR\sHom_X(A,B) \, .$$
Therefore, we also get:
$$f^*\, D_X\simeq D^2_Y \, f^*\, D_X \simeq D_Y\, f^! \, D^2_X \simeq D_Y \, f^!\, .$$
The two other formulas of (d) follow by adjunction.
\end{proof}
\index{word}{constructible!$\tau$-constructible|)}

\begin{thm}\label{thm:realcommute6op}
Assume that $\site$ consists of schemes
of finite type over a field $k$. We consider a
$\tau'$-generated
motivic triangulated category $\T'$ over $\site$ as well as
a premotivic morphism
$$\varphi^*:\T \To \T'\, .$$
We suppose that the following properties hold:
\begin{itemize}
\item[(a)] $\T$ is $\tau$-dualizable,
$\QQ$-linear and separated;
\item[(b)] $\T'$ is $\QQ$-linear and separated;
\item[(b)] the object $\unit\{i\}$ is rigid in $\T(k)$
for any $i\in\tau$.
\end{itemize}
Then, the premotivic morphism
$$\varphi^*:\T_c \To \T'$$
commutes with the six operations.
\end{thm}

\begin{rem}
Remark that, as a corollary, we obtain immediately,
under the assumptions of the theorem that $\T'$ is
$\varphi^*(\tau)$-dualizable and that
the functor $\varphi^*$ commutes with the duality
 functors on $\T$ and $\T'$,
 respectively obtained by applying
 the above corollary in the case $B=\spec k$.
\end{rem}

\begin{proof}
Given a morphism of finite type $f:X\to\spec k$, let us consider the
following property.
\begin{itemize}
\item[$(\ast)_f$] \emph{For any constructible object $M$ in $\T(X)$,
 the natural exchange map
$$\varphi^*\, f_*(M)\to f_*\,\varphi^*(M)$$
is invertible.}
\end{itemize}

We will first prove the theorem assuming that property $(\ast)_f$ holds for any $f$.

Let $u:X\to Y$ be a $k$-morphism of finite type.
We claim that the exchange map
$$\varphi^*\, u_*(M)\To u_*\, \varphi^*(M)$$
is invertible for any $\tau$-constructible object $M$ of $\T(X)$.

It is sufficient to prove that, for any smooth separated
$k$-morphism of finite type $g:T\to X$, any constructible
object $M$ in $\T(X)$ and any twist $i$ in $\tau'$,
the natural map
$$\Hom_{\T'(X)}(g_\sharp(\unit_T)\{i\},\varphi^*\, u_*(M))
\To\Hom_{\T'(X)}(g_\sharp(\unit_T)\{i\}, u_*\varphi^*(M))$$
is bijective.
Consider the following commutative diagram of morphisms of schemes:
$$
\xymatrix@C=6pt@R=16pt{
V\ar^v[rr]\ar_h[d] && T\ar^g[d] \\
X \ar_-a[rd]\ar^u[rr] && Y\ar^-b[ld] \\
& \spec k & 
}
$$
in which the square is cartesian. Recall that the functor
 $v_*$ preserves constructible objects by virtue of
 Theorem \ref{thmfinitness0}.
 Then we conclude by the computations below:
\begin{align*}
\Hom_{\T'(Y)}(g_\sharp(\unit_T)\{i\},\varphi^*\, u_*(M))
&= \Hom_{\T'(T)}(\unit_T\{i\},g^*\,\varphi^*\, u_*(M)) \\
&= \Hom_{\T'(T)}(\unit_T\{i\},\varphi^*\,g^*\, u_*(M)) \\
&= \Hom_{\T'(T)}(g^*b^*(\unit_k)\{i\},\varphi^*\,g^*\, u_*(M)) \\
&= \Hom_{\T'(T)}(g^*b^*(\unit_k)\{i\},\varphi^*\,v_*\,h^*(M)) \\
&= \Hom_{\T'(k)}(\unit_k\{i\},(bg)_*\,\varphi^*\,v_*\,h^*(M)) \\
&= \Hom_{\T'(k)}(\unit_k\{i\},\varphi^*\,(bg)_*\,v_*\,h^*(M)) 
 \quad \text{(by } (\ast)_{bg}) \\
&= \Hom_{\T'(k)}(\unit_k\{i\},(bgv)_*\,\varphi^*\,h^*(M)) 
 \quad \text{(by } (\ast)_{bgv}) \\
&= \Hom_{\T'(k)}(\unit_k\{i\},(bg)_*\,g^*\,u_*\,\varphi^*(M)) \\
&= \Hom_{\T'(Y)}(g_\sharp(\unit_T)\{i\},u_*\,\varphi^*(M)) \\
\end{align*}
From there, we see that, for any $k$-scheme of finite type $X$
and any $\tau$-constructible objects $M$ and $N$ of $\T(X)$, the natural
map
$$\varphi^*(\uHom_X(M,N))\To\uHom_X(\varphi^*(M),\varphi^*(N))$$
is invertible in $\T'(X)$. For this, we may assume that
$M=f_\sharp(\unit_Y\{i\})$ for a smooth morphism of finite type $f:Y\To X$
and a twist $i$, in which case we have
$$\varphi^*(\uHom_X(M,N))=\varphi^*\, f_*\, f^*(N)\simeq
f_*\, f^*\, \varphi^*(N)=\uHom_X(\varphi^*(M),\varphi^*(N)).$$
It remains to prove that for any separated $k$-morphism $f:X \rightarrow Y$
of finite type and any constructible object $N$ in $\T(X)$,
the exchange map:
$$
\varphi^*\,f^!(N) \rightarrow f^!\,\varphi^*(N)
$$
is an isomorphism. It is easy to see that this property is
local for the Zariski topology, both on $X$ and on $Y$, so that we may assume that the morphism $f$ is affine. Therefore, it is sufficient to consider the situation where $f$ i either a closed immersion or a separated smooth map. In the smooth case, as the functor $f^!$ is of the form $f^*(d)[2d]$, this is obvious. If $f=i$ is
a closed immersion with open complement $j$, as we already
know that $\varphi^*$ commutes with $u_*$
for any morphism  $u$, this property follows straight
away from the localization distinguished triangles
$$i_*\, i^!\To 1 \To j_*\, j^*\To\, .$$
%
%It is sufficient to prove that for any constructible
% object $M$ in $\T(X)$, the induced map:
%$$
%\Hom_{\T'(Y)}(\varphi^*(M),\varphi^*\,f^!(N)) \rightarrow
% \Hom_{\T'(Y)}(\varphi^*(M),f^!\,\varphi^*(N))
%$$
%is bijective, because $\varphi^*$ is essentially surjective on the 
% family of generators of $\T'_c(Y)$ of the form $g_\sharp(\unit_W\{i\})$
% for a smooth morphism $g$ and a twist $i$.
%This follows formally by applying the functor $\Hom_{\T'(Y)}(\unit_Y,-)$
% to the isomorphisms:
%\begin{align*}
%f_*\uHom_X(\varphi^*(M),\varphi^*f^!(N))
%&=\varphi^*\,f_*\,\uHom_X(M,f^!(N)) \\
%&=\varphi^*\,\uHom_X(f_!(M),N) \\
%&=\uHom_X(f_!\,\varphi^*(M),\varphi^*(N)) \\
%&=f_*\uHom_X(\varphi^*(M),f^!\,\varphi^*(N))
%\end{align*}

It remains to prove property $(\ast)_f$ for any morphism $f$ of finite type.

We claim it is sufficient to prove that,
 for any $k$-scheme of finite type $X$ with structural morphism $f$,
 the following property holds:
\begin{itemize}
\item[$(\ast\ast)_X$] \emph{For any twist $i \in \tau$,
 the natural exchange map
$$\varphi^*\, f_*(\unit_X\{i\})\to f_*\,\varphi^*(\unit_X\{i\})$$
is invertible.}
\end{itemize}
Indeed, by virtue of Theorem \ref{thmfinitnessproper2}, we may assume that
$M=w_*(\unit_W\{i\})$ for $w:W\to X$ a projective $k$-morphism, and $i\in \tau$.
As the exchange map $\varphi^*\, w_*\To w_*\, \varphi^*$ is invertible
(Proposition \ref{prop:motivic_adj_6_operations}), we see that we may assume that
$M=\unit_X\{i\}$ for some twist $i$.

Let us prove property $(\ast\ast)_X$ in the case $X$ is in addition smooth over $k$.
As $\varphi^*$ is monoidal, for any rigid object $M$ of $\T(k)$, we get
the identification:
$$
\varphi^*(M^\vee)=\varphi^*(M)^\vee.
$$
On the other hand, according to assumption (b), the object $f_\sharp(\unit_X)$
 is rigid in $\T(k)$ as well as in $\T'(k)$ (because the functor $\varphi^*$
 is symmetric monoidal and commutes with the operations of the form $f_\sharp$
 for $f$ smooth). Thus we get:
$$
f_*(\unit_X\{i\})=\uHom_k(f_\sharp(\unit_X),\un_k\{i\})
=f_\sharp(\unit_X)^\vee\{i\}.
$$
Then property $(\ast\ast)_X$ readily follows.

We finally prove property $(\ast\ast)_X$ for any algebraic  $k$-scheme $X$.
We will proceed by induction on the dimension of $X$.

In case $\dim(X)<0$, the result is obvious.
Let us assume $\dim(X) \geq 0$.
According to the localization property, we can assume that $X$ is reduced.
Let $\bar k$ be an inseparable closure of $k$
 and $\bar X=X \otimes_k \bar k$.
 According to De Jong theorem applied to $\bar X$
 (see Th. \ref{thm:resuptoquotient} for $S=\spec{\bar k}$),
 there exists a Galois alteration $\bar X' \rightarrow \bar X$
 of group $G$ such that $\bar X'$ is smooth over $\bar k$.

We can assume that such a smooth alteration exists over
 a finite inseparable extension field $E/k$.
 Because $\T$ (resp. $\T'$) is $\QQ$-linear and separated,
 the base change functor $\pi^*$ associated
 with the finite morphism $\pi:\spec E \rightarrow \spec k$
 and relative to the premotivic category $\T$ (resp. $\T'$)
 is an equivalence of categories 
 (see Proposition \ref{stronglyquasiseparated}).
 Thus we can replace $k$ by $E$ and assume
 that there exists a Galois alteration $p:X' \rightarrow X$
 of group $G$ such that $X'$ is a smooth $k$-scheme.
Using the localization property, we can assume $X$ is reduced.
 Then there exists a nowhere dense closed subscheme $\nu:Z \rightarrow X$
 such that $U=X-Z$ is regular (thus normal) and
 the induced map $p|_U:p^{-1}(U) \rightarrow U$ is finite.
 Thus we can apply Theorem \ref{excisionGaloisalteration} to
 the cartesian square:
$$
\xymatrix{
Z'\ar[r]^{\nu'}\ar[d]_q&X'\ar[d]^p\\
Z\ar[r]^\nu&X
}
$$
and we get the distinguished triangle in $\T(X)$
(thus in $\T'(X)$ as well, as the functor $\varphi^*$
is monoidal and commutes with the operations of the form $u_*$
for any proper morphism $u$) of the form:
$$
\un_X\{i\} \rightarrow p_*(\un_{X'}\{i\})^G \oplus \nu_*(\un_Z\{i\})
 \rightarrow (\nu q)_*(\un_{Z'}\{i\})^G \xrightarrow{+1}
$$
for any twist $i$.
If we consider the triangles in $\T(k)$ and $\T'(k)$
 obtained by applying the functor $f_*$,
 where $f$ is the structural morphism of $X/k$,
 we deduce that property $(\ast\ast)_X$
 follows from properties $(\ast\ast)_{X'}$, $(\ast\ast)_{Z}$, $(\ast\ast)_{Z'}$.
Thus we can conclude applying either the case of a smooth $k$-scheme
 treated above or the induction hypothesis as $\dim(Z)=\dim(Z')<\dim(X)$. 
\end{proof}

\part{Construction of fibred categories}
\markboth{Construction of fibred categories}{}
\section{Fibred derived categories}\label{sec:fibred_derived}

\begin{assumption} \label{num:assumption_sch_derived}
In this entire section,
 we fix a full subcategory $\sch$
 of the category of noetherian $\base$-schemes
 satisfying the following properties:
\begin{itemize}
\item[(a)] $\sch$ is closed under finite sums and pullback along
 morphisms of finite type.
\item[(b)] For any scheme $S$ in $\sch$, any quasi-projective $S$-scheme
belongs to $\sch$.
\end{itemize}

We fix an admissible class of morphisms $\Pmor$ of $\sch$.
All our $\Pmor$-premotivic categories
 (\textit{cf.} definition \ref{df:general_P-premotivic})
 are defined over $\sch$.
Moreover, for any abelian $\Pmor$-premotivic category $\A$ in this section,
 we assume the following:
\begin{itemize}
\item[(c)] $\A$ is a \emph{Grothendieck} abelian $\Pmor$-premotivic category
 (see definition \ref{df:P-fibred_Grothendieck_abelian}
 and the recall below).
\item[(d)] $\A$ is given with a generating set of twists $\tau$.
 We sometimes refer to it as \emph{the twists of $\A$}.
\item[(e)] We will denote by $\Mab S X \A$, or
simply by $\mab S X$, the geometric section over a $\Pmor$-scheme $X/S$.
\end{itemize}
Without precision, any scheme will be assumed to be an object of $\sch$.

In section \ref{sec:A^1-derived}, except possibly for \ref{sec:loc_derived},
 we assume further:
\begin{itemize}
\item[(f)] $\Pmor$ contains the class of smooth morphisms of finite type.
\end{itemize}
%%In section \ref{sec:P^1-stable-derived}, we assume (f) and instead of (d)
%%above.
%%\begin{itemize}
%%\item[(d')] Without precision, an abelian $\Pmor$-premotivic category $\A$
%% is geometrically generated (recall this means $\tau$-generated for $\tau=0$).
%%\end{itemize}
\end{assumption}

\begin{num} \label{num:canonical_DG-structure}
We will sometimes refer to the canonical dg-structure of
 the category of complexes $\Comp(\A)$ over an abelian category $\A$.
 Recall that to any complexes $K$ and $L$ over $\A$,
 we associate a complex of abelian groups $\Hom^\bullet_\A(K,L)$
\index{notat}{hombullet@$\Hom^\bullet(-,-)$}
  whose component
 in degree $n \in \ZZ$ is
$$
\prod_{p \in \ZZ} \Hom_\A(K^p,L^{p+n})
$$
and whose differential in degree $n \in \ZZ$ is defined by the formula:
$$
(f_p)_{p \in \ZZ}
 \mapsto \big(d_L \circ f_p-(-1)^n.f_{p+1} \circ d_K)\big)_{p \in \ZZ}.
$$
In other words, this is the image of the bicomplex
 $\Hom_\A(K,L)$ by the Tot-product functor which we denote by $\Totp$.
\index{notat}{totp@$\Totp$}
Of course, the associated homotopy category is the category $\K(\A)$
 of complexes up to chain homotopy equivalence.
 \end{num}

\subsection{From abelian premotives to triangulated premotives}
\label{sec:derived_premotivic}

\subsubsection{Abelian premotives: recall and examples}

Consider an abelian $\Pmor$-premotivic category $\A$.
According to the convention of \ref{num:assumption_sch_derived},
 for any scheme $S$, $\A_S$ is a Grothendieck abelian closed symmetric monoidal category.
Moreover, if $\tau$ denotes the twists of $\A$, the essentially small family
$$
\big(\mab S X \{i\}\big)_{X \in \Pmorx S, i \in \tau}
$$
is a family of generators of $\A_S$ in the sense of \cite{Gro1}.

%Recall that according to definition \ref{df:generating_twists},
% $\A$ is compactly $\tau$-generated if for any $\Pmor$-scheme
% $X/S$ and any twist $i \in \tau$, the premotive $\Mab S X \A\{i\}$
% is compact.

\begin{ex} \label{ex:abelian_premotivic_presheaves}
Consider a fixed ring $\Rc$.
Let $\psh{\Pmorx S}$
\index{notat}{pshP@$\psh{\Pmorx S}$}
 be the category of $\Rc$-presheaves
\index{word}{presheaf!lambdapresheaf@$\Lambda$-presheaf}
(i.e. presheaves of $\Rc$-modules) on $\Pmorx S$. 
For any $\Pmor$-scheme $X/S$,
we let $\rep S X$ be the free $\Rc$-presheaf
on $\Pmorx S$ represented by $X$. Then $\psh{\Pmorx S}$
 is a Grothendieck abelian category generated by
the essentially small family 
$\big(\rep S X\big)_{X \in \Pmorx S}$.

There is a unique symmetric closed monoidal structure on
 $\psh{\Pmorx S}$ such that
$$
\rep S X \otimes_S \rep S Y=\rep S {X \times_S Y}.
$$
Finally the existence of functors $f^*$,
$f_*$ and, in the case when $f$ is a $\Pmor$-morphism, of $f_\sharp$,
follows from general sheaf theory (\textit{cf.} \cite{SGA4}).

Thus, $\psh{\Pmor}$ defines an abelian $\Pmor$-premotivic category.
\end{ex}

\begin{num} \label{num:universality_presheaves}
\renewcommand{\Rc}{\ZZ}
Consider an abstract abelian $\Pmor$-premotivic category $\A$. 
To any premotive $M$ of $\A_S$, we can associate a presheaf of
abelian groups
$$X \mapsto \Hom_{\A_S}(\mab S X,M)$$
which we denote by $\gamma_*(M)$.  \\
This defines a functor
$\gamma_*:\A_S \rightarrow \psh{\Pmorx S}$.
It admits the following left adjoint:
$$
\gamma^*:
 \psh{\Pmorx S} \rightarrow \A_S \ , \quad F \mapsto \ilim_{X/F} \Mab S X \A
$$
where the colimit runs over the category of representable
presheaves over $F$.

It is now easy to check we have defined
 a morphism of (complete) abelian $\Pmor$-premotivic categories:
\begin{equation} \label{eq:universality_presheaves}
\gamma^*:\psh{\Pmor} \rightleftarrows \A:\gamma_*.
\end{equation}
Moreover $\psh{\Pmor}$ appears as the initial abelian
 $\Pmor$-premotivic category.

Remark that the functor $\gamma_*:\A_S \rightarrow \psh{\Pmorx S}$ is
conservative if the set of twists $\tau$ of $\A$ is trivial.
\end{num}

\begin{df}
A $\Pmor$-admissible topology
\index{word}{topology!Padmissible@$\Pmor$-admissible}
 $t$ is a Grothendieck pretopology $t$
on the category $\sch$, such that any $t$-covering family consists of
$\Pmor$-morphisms.
\end{df}
Note that, for any scheme $S$ in $\sch$, such a topology $t$ induces a
pretopology on $\Pmorx S$ (which we denote by the same letter).
For any morphism (resp. $\Pmor$-morphism) $f:T\To S$,
the functor $f^*$ (resp. $f_\sharp$) preserves $t$-covering families.

As $\Pmor$ is fixed in all this section,
 we will simply say \emph{admissible}
\index{word}{topology!admissible}
\index{word}{admissible topology|see{topology}}
  for $\Pmor$-admissible.

\begin{ex}
\label{ex:abelian_premotivic_sheaves}
Let $t$ be an admissible topology.
We denote by $\sh t {\Pmorx S}$
\index{notat}{shtP@$\sh t {\Pmorx S}$}
 the category of $t$-sheaves of $\Rc$-modules
\index{word}{sheaf!tsheaf@$t$-sheaf of $\Lambda$-modules}
  on $\Pmorx S$.
Given a $\Pmor$-scheme $X/S$, 
 we let $\repx t S X$
\index{notat}{lambdat@$\repx t S X$}
  be the free $\Rc$-linear $t$-sheaf represented by $X$.
Then, $\sh t {\Pmorx S}$ is an abelian Grothendieck category
 with generators $(\repx t S X)_{X \in \Pmorx S}$.

As in the preceding example, 
 the category $\sh t {\Pmorx S}$ admits
  a unique closed symmetric monoidal structure 
 such that $\repx t S X \otimes_S \repx t S Y=\repx t S {X \times_S Y}$. 
Finally, for any morphism $f:T \rightarrow S$ of schemes,
 the existence of functors $f^*$, $f_*$ 
 (resp. $f_\sharp$ when $f$ is a $\Pmor$-morphism) follows from the general
 theory of sheaves (see again \cite{SGA4}: according to our assumption on $t$
 and \cite[III, 1.6]{SGA4}, the functors
 $f^*:\Pmorx S \rightarrow \Pmorx T$ and $f_\sharp:\Pmorx T \rightarrow \Pmorx S$
 (for $f$ in $\Pmor$) are continuous).

Thus, $\sh t {\Pmor}$ defines an abelian $\Pmor$-premotivic category
(with trivial set of twists).

The associated $t$-sheaf functor induces a morphism
\begin{equation} \label{eq:abelian_premotivic_sheaves}
a^*_t:\psh{\Pmor} \rightleftarrows \sh t {\Pmor}:a_{t,*}.
\end{equation}
\end{ex}

\begin{rem}
\renewcommand{\Rc}{\ZZ}
Recall the abelian category $\sh t {\Pmorx S}$ is a localization
of the category $\psh S$ in the sense of Gabriel-Zisman.
In particular, 
given an abstract abelian $\Pmor$-premotivic category $\A$,
the canonical morphism
$$
\gamma^*:\psh{\Pmorx S} \rightleftarrows \A_S:\gamma_*
$$
induces a unique morphism
$$
\sh t {\Pmorx S} \rightleftarrows \A_S
$$
if and only if for any presheaf of abelian groups $F$ on $\Pmorx S$ such that
$a_t(F)=F_t=0$, one has $\gamma^*(F)=0$.

We leave to the reader the exercise which consists of formulating
the universal property of the abelian $\Pmor$-premotivic category
$\sh t {\Pmor}$.\footnote{We will formulate a derived version in the paragraph
on descent properties for derived premotives (\textit{cf.} \ref{descent_derived_premotives}).}
\end{rem}

\subsubsection{The $t$-descent model category structure}

\begin{paragr}
Consider an abelian $\Pmor$-premotivic category $\A$ with set of twists $\tau$.

We let $\Comp(\A)$ be the $\Pmor$-fibred abelian category
over $\sch$ whose fibers over a scheme $S$ is the
category $\Comp(\A_S)$ of (unbounded) complexes in $\A_S$.
For any scheme $S$, we let
$\iota_S:\A_S \rightarrow \Comp(\A_S)$ the embedding
which sends an object of $\A_S$ to the corresponding complex 
concentrated in degree zero.

If $\A$ is $\tau$-twisted, then the category $\Comp(\A_S)$ is obviously $(\ZZ \times \tau)$-twisted.
The following lemma is straightforward:
\end{paragr}
\begin{lm}
\label{lm:complexes_of_premotivic_abelian}
With the notations above,
 there is a unique structure of abelian $\Pmor$-premotivic category
on $\Comp(\A)$ such that the functor
$\iota:\A \rightarrow \Comp(\A)$ is
 a morphism of abelian $\Pmor$-premotivic categories.
\end{lm}

\begin{paragr} \label{num:motives_simplicial_sm}
For a scheme $S$, let $(\Pmorx S)^\amalg$ be the category introduced
in \ref{hypothesissite}. The functor $\mab S -$
can be extended to $(\Pmorx S)^\amalg$ by associating to a
family $(X_i)_{i \in I}$ of $\Pmor$-schemes over $S$
the premotive
$$\bigoplus_{i \in I} \mab S {X_i}.$$
If $\cX$ is a simplicial object of $(\Pmorx S)^\amalg$,
we denote by $\mab S \cX$ the complex associated 
with the simplicial object of $\A_S$ obtained
by applying degreewise the above extension of $\mab S -$.
\end{paragr}

\begin{df}
\label{df:basic_complexes&topology}
Let $\A$ be an abelian $\Pmor$-premotivic category
and $t$ be an admissible topology.

Let $S$ be a scheme and $C$ be an object of $\Comp(\A_S)$:
\begin{enumerate}
\item The complex $C$ is said to be \emph{local}
\index{word}{local}
(with respect to the geometric section)  if, for any $\Pmor$-scheme $X/S$ and
any pair $(n,i) \in \ZZ \times \tau$, the canonical morphism
$$
\Hom_{\K(\A_S)}(\mab S X\{i\}[n],C) 
 \rightarrow \Hom_{\Der(\A_S)}(\mab S X\{i\}[n],C)
$$
is an isomorphism.
\item The complex $C$ is said to be \emph{$t$-flasque}
\index{word}{flasque, $t$-flasque complex}
 if for any $t$-hypercover $\mathcal X \rightarrow X$ in $\Pmorx S$,
 for any $(n,i) \in \ZZ \times \tau$, the canonical morphism

$$
\Hom_{\K(\A_S)}(\mab S X\{i\}[n],C)
 \rightarrow \Hom_{\K(\A_S)}(\mab S \cX\{i\}[n],C)
$$
is an isomorphism.
\end{enumerate}
We say the abelian $\Pmor$-premotivic category $\A$ satisfies
\emph{cohomological $t$-descent}
\index{word}{descent!cohomological $t$-descent}
 if for any $t$-hypercover $\mathcal X \rightarrow X$ of a
$\Pmor$-scheme $X/S$, and for any $i\in\tau$, the map
$$
\mab S \cX \{i\} \rightarrow \mab S X \{i\}
$$
is a quasi-isomorphism (or equivalently, if any local complex
is $t$-flasque).

We say that $\A$ is \emph{compatible with $t$}
\index{word}{compatible with (a topology) $t$}
 if $\A$ satisfies cohomological $t$-descent, and if,
for any scheme $S$, any $t$-flasque complex of $\A_S$ is local.
\end{df}
%It follows from \cite{CD1} that 
%if $\A$ is compatible with an admissible topology $t$
% then it satisfies cohomological $t$-descent 
%(\textit{cf.} remark \ref{rem:resolutions&derivations} for a proof).

\begin{ex} \label{ex:t-local_sheaves}
Consider the notations of \ref{ex:abelian_premotivic_sheaves}.

Consider the canonical dg-structure
\index{word}{dgstructure@dg-structure}
 on $\Comp(\sh t {\Pmorx S})$
 (see \ref{ex:abelian_premotivic_presheaves}).
By definition, for any complexes $D$ and $C$ of sheaves,
 we get an equality:
\begin{align*}
\Hom_{\K(\sh t {\Pmorx S})}(D,C)
&=H^0(\Hom_{\sh t {\Pmorx S}}^\bullet(D,C))\\
&=H^0(\Totp \Hom_{\sh t {\Pmorx S}}(D,C)).
\end{align*}
In the case where $D=\Rc_S^t(X)$ (resp. $D=\Rc_S^t(\cX)$)
 for a $\Pmor$-scheme $X/S$ (resp. a simplicial $\Pmor$-scheme over $S$)
 we obtain the following identification:
\begin{align*}
& \Hom_{\K(\sh t {\Pmorx S})}(\Rc^t_S(X),C)=H^0(C(X)). \\
\text{(resp. } & \Hom_{\K(\sh t {\Pmorx S})}(\Rc^t_S(\cX),C)
 =H^0(\Totp C(\cX)) \, ).
\end{align*}
Thus, we get the following equivalences:
\begin{align*}
C \text{ is local} & \Leftrightarrow \text{for any $\Pmor$-scheme $X/S$, }
H^n_t(X,C)\simeq H^n(C(X)). \\
C \text{ is $t$-flasque} & \Leftrightarrow \text{for any $t$-hypercover
$\cX \rightarrow X$, }
H^n(C(X))\simeq H^n(\Totp C(\cX)).
\end{align*}

According to the computation of cohomology with hypercovers
(\textit{cf.} \cite{Brown}),
if the complex $C$ is $t$-flasque, it is local. In other words,
we have the expected property that
the abelian $\Pmor$-premotivic category $\sh t {\Pmor}$ is 
compatible with $t$.\index{word}{compatible with (a topology) $t$}
\end{ex}

%\begin{num} \label{Nisnevich&BG}
%Let $C$ be a complex of presheaves of $\Rc$-modules
%on $\sm/S$.
%We say that $C$ satisfies the B.G.-property if 
%the image by $C$ of any distinguished square (\textit{cf.} \ref{BG-property-triangulated})
%$$
%\xymatrix@=10pt{
%W\ar^k[r]\ar_g[d] & V\ar^f[d] \\ U\ar^j[r] & X
%}
%$$
%is homotopy cartesian in the derived category of complexes of $\Rc$-modules.
%%$
%%\xymatrix@=10pt{
%%W\ar^k[r]\ar_g[d] & V\ar^f[d] \\ U\ar^j[r] & X
%%}
%%$
%%the image square
%%$$
%%\xymatrix@=10pt{
%%C(X)\ar^{f^*}[r]\ar_{j^*}[d] & C(V)\ar^{k^*}[d] \\ C(U)\ar^{k^*}[r] & C(W)
%%}
%%$$
%%is homotopy cartesian.
%
%It follows from Morel and Voevodsky's lemma  \cite[1.17, page 100]{MV}
%that $C$ is $\nis$-local if and only if it satisfies the B.G.-property.
%\end{num}

\begin{paragr} \label{t-descent_model_structure}
Consider an abelian $\Pmor$-premotivic category  $\A$
 and an admissible topology $t$.

Fix a base scheme $S$.
A morphism $p:C \rightarrow D$ of complexes on $\A_S$
is called a \emph{$t$-fibration}
\index{word}{fibration!$t$-fibration}
 if its kernel is a $t$-flasque complex
and if for any $\Pmor$-scheme $X/S$, any $i \in \tau$ and any integer
$n \in \ZZ$, the map of abelian groups
$$\Hom_{\A_{S}}(\mab S X\{i\},C^n)\To\Hom_{\A_{S}}(\mab S X\{i\},D^n)$$
is surjective.

For any object $A$ of $\A_S$, we let $S^nA$ (resp. $D^nA$)
be the complex with only one non-trivial term (resp. two non-trivial terms)
equal to $A$ in degree $n$ (resp. in degree $n$ and $n+1$, with the identity
as only non-trivial differential). We define the class of \emph{cofibrations}
\index{word}{cofibration}
as the smallest class of morphisms of $\Comp(\A_S)$ which:
\begin{enumerate}
\item contains the map
$S^{n+1} \mab S X\{i\} \rightarrow D^n \mab S X\{i\}$ for any $\Pmor$-scheme
$X/S$, any $i \in \tau$, and any integer $n$;
\item is stable by pushout, transfinite composition and retract.
\end{enumerate}
A complex $C$ is said to be \emph{cofibrant} if the canonical map $0 \rightarrow C$
is a cofibration. For instance, for any $\Pmor$-scheme $X/S$ and any $i \in \tau$, 
the complex $\mab S X \{i\}[n]$ is cofibrant.
\end{paragr}

Let $\cG_S$ be the essentially small family made of premotives $\mab S X\{i\}$
for a $\Pmor$-scheme $X/S$ and a twist $i \in \tau$, and $\cH_S$ be the family of 
complexes of the form 
$\mathrm{Cone}(\mab S \cX\{i\} \rightarrow \mab S X\{i\})$
for any $t$-hypercover $\cX \rightarrow X$ and any twist $i \in \tau$.
By the very definition, as $\A$ is compatible with $t$ 
(definition \ref{df:basic_complexes&topology}),
$(\cG_S,\cH_S)$ is a descent structure on $\A_S$ in the sense of
 \cite[def. 2.2]{CD1}.
Moreover, it is weakly flat in the sense of \cite[par. 3.1]{CD1}.
Thus the following proposition is a particular case of 
\cite[theorem 2.5, proposition 3.2, and corollary 5.5]{CD1}:
\begin{prop}\label{t-model_category_complexes}
Let $\A$ be an abelian $\Pmor$-premotivic category, which we assume to be
compatible with an admissible topology $t$. Then for any scheme $S$, 
the category $\Comp(\A_S)$ with the preceding definition
of fibrations and cofibrations, with quasi-isomorphisms as weak equivalences
is a proper symmetric monoidal model category.
%%(\textit{cf.} \cite{Hov} for definitions). 
\end{prop}

\begin{num} \label{num:derived&generators}
We will call this model structure on $\Comp(\A_S)$
 the \emph{$t$-descent model category structure}
\index{word}{model structure!$t$-descent}
  (over $S$).
Note that, for any $\Pmor$-scheme $X/S$ and any twist $i \in \tau$,
the complex $\mab S X\{i\}$ concentrated in degree $0$ is cofibrant
by definition, as well as any of its suspensions and twists.
They form a family of generators for the triangulated category $\Der(A_S)$. \\
Observe also that the fibrant objects
for the $t$-descent model category structure are exactly the
$t$-flasque complexes in $\A_S$.
Moreover, essentially by definition,
a complex of $\A_S$ is local 
if and only if it is $t$-flasque (see \cite[2.5]{CD1}).
\end{num}

\begin{num} \label{t-descent_model_structure_bis}
Consider again the notations and hypothesis of \ref{t-descent_model_structure}.

Consider a morphism of schemes $f:T \rightarrow S$.
Then the functor 
$$
f^*:\Comp(\A_S) \rightarrow \Comp(\A_T)
$$
sends $\cG_S$ in $\cG_T$,
 and $\cH_S$ in $\cH_T$ because the topology $t$ is admissible.
This means it satisfies descent according to the definition of \cite[2.4]{CD1}.
Applying theorem 2.14 of \emph{op. cit.}, the functor $f^*$ preserves cofibrations
and trivial cofibrations, i.e. the pair of functors $(f^*,f_*)$ is a Quillen adjunction 
with respect to the $t$-descent model category structures.

Assume that $f$ is a $\Pmor$-morphism. Then, similarly, the functor
$$
f_\sharp:\Comp(\A_T) \rightarrow \Comp(\A_S)
$$
sends $\cG_S$ (resp. $\cH_S$) in $\cG_T$
 (resp. $\cH_T$) so that $f_\sharp$ also satisfies descent
in the sense of \emph{op. cit}.
Therefore, it preserves cofibrations and trivial cofibrations,
and the pair of adjoint functors $(f_\sharp,f^*)$ is a Quillen adjunction for the
$t$-descent model category structures.
 
In other words, we have obtained the following result.
\end{num}

\begin{cor}
Let $\A$ be an abelian $\Pmor$-premotivic category
 compatible with an admissible topology $t$.
The $\Pmor$-fibred category $\Comp(\A)$ with the $t$-descent
model category structure defined in \ref{t-model_category_complexes}
is a symmetric monoidal $\Pmor$-fibred model
\index{word}{fibred!Pfibredcategory@$\Pmor$-fibred category!monoidal $\Pmor$-fibred model category}
 category. 
Moreover, it is stable, proper and combinatorial.
\end{cor}

\begin{num} \label{num:resolutions&derivations}
Recall the following consequences of this corollary
 (see also \ref{num:P-fibred_model_der_functors} for the general theory).
Consider a morphism $f:T \rightarrow S$ of schemes.
Then the pair of adjoint functors $(f^*,f_*)$ admits total left/right derived
 functors
$$
\derL f^*:\Der(\A_S) \rightleftarrows \Der(\A_T):\derR f_*.
$$
More precisely,
 $f_*$ (resp. $f^*$) preserves $t$-local (resp. cofibrant) complexes.
For any complex $K$ on $\A_S$,
 $\derR f_*(K)=f_*(K')$ (resp. $\derL f^*(K)=f^*(K'')$) where $K' \rightarrow K$
 (resp. $K \rightarrow K''$)
 is a $t$-local (resp. cofibrant) resolution of $K$.\footnote{Recall also
 that fibrant/cofibrant resolutions can be made functorially,
 because our model categories are cofibrantly generated,
 so that the left or right derived functors are in fact
 defined at the level of complexes.}

When $f$ is a $\Pmor$-morphism, the functor $f^*$ is even exact
 and thus preserves quasi-isomorphisms. This implies that $\derL f^*=f^*$.
 The functor $f_\sharp$ admits a total left derived functor
$$
\derL f_\sharp:\Der(\A_T) \rightleftarrows \Der(\A_S):\derR f^*
$$
defined by the formula $\derL f_\sharp(K)=f_\sharp(K'')$ for a complex
$K$ on $\A_T$ and a cofibrant resolution $K'' \rightarrow K$.

Note also that the tensor product (resp. internal Hom) of $\Comp(\A_S)$
 admits a total left derived functor (resp. total right derived functor).
For any complexes $K$ and $L$ on $\A_S$, this derived functors are defined
 by the formula:
\begin{align*}
K \otimes_S^\derL L& =K'' \otimes_S L'' \\
\derR\! \uHom_S(K,L)&=\uHom_S(K'',L')
\end{align*}
where $K \rightarrow K''$ and $L \rightarrow L''$ are cofibrant resolutions
 and $L' \rightarrow L$ is a $t$-local resolution.

It is now easy to check that these functors define a
 triangulated $\Pmor$-premotivic category $\Der(\A)$,
 which is $\tau$-generated according to \ref{num:derived&generators}.
\end{num}
\begin{df} \label{df:derived_premotivic_category}
Let $\A$ be an abelian $\Pmor$-premotivic category
compatible with an admissible topology $t$.

The triangulated $\Pmor$-premotivic category $\Der(\A)$
 defined above is called the \emph{derived $\Pmor$-premotivic category}
\index{word}{premotivic!Ppremotivic@$\Pmor$-premotivic!derived category}
 associated with $\A$.\footnote{Indeed remark that
  $\Der(\A)$ does not depend on the topology $t$.}
\end{df}
The geometric section of a $\Pmor$-scheme $X/S$ in
 the category $\Der(\A)$ is the complex concentrated in degree
 $0$ equal to the object $\mab S X$.
The triangulated $\Pmor$-fibred category is $\tau$-generated
 and well generated in the sense of \ref{df:P-premotivic_well&compactly_generated}.
 Recall this means that $\Der(\A_S)$ is equal to
 the localizing\footnote{\emph{i.e.} triangulated and stable by sums.} 
 subcategory generated by the family
\begin{equation}\label{eq:premotivic_derived_generated}
\{ \mab S X \{i\} ; X/S\text{ $\Pmor$-scheme}, i \in \tau \}.
\end{equation}
 
\begin{ex} \label{ex:abelian_premotivic_derived_sheaves}
Given any admissible topology $t$,
the abelian $\Pmor$-premotivic category $\sh t {\Pmor}$ introduced
in example \ref{ex:abelian_premotivic_sheaves}
is compatible with $t$ (\textit{cf.} \ref{ex:t-local_sheaves})
and defines the derived $\Pmor$-premotivic category
$\Der(\sh t {\Pmor})$.

Remark also that the abelian $\Pmor$-premotivic category $\psh{\Pmor}$
 introduced in example \ref{ex:abelian_premotivic_presheaves}
 is compatible with the coarse topology
 and gives the derived $\Pmor$-premotivic category $\Der(\psh{\Pmor})$.
\end{ex}

\begin{rem} \label{rem:DG-structure}
Recall from \ref{num:canonical_DG-structure}
there exists a canonical dg-structure on $\Comp(\A_S)$.
Then we can define a derived dg-structure
by defining for any complexes $K$ and $L$ of $\A_S$,
the complex of morphisms:
$$
\derR \Hom_{\A_S}(K,L)=\Hom^\bullet_{\A_S}(Q(K),R(L))
$$
where $R$ and $Q$ are respectively some fibrant and cofibrant (functorial)
resolutions for the $t$-descent model structure. 
The homotopy category associated with this new dg-structure
on $\Comp(\A_S)$ is the derived category $\Der(\A_S)$.
Moreover, for any morphism (resp. $\Pmor$-morphism) of schemes $f$,
the pair $(\derL f^*,\derR f_*)$ (resp. $(\derL f_\sharp,f^*)$)
is a dg-adjunction. 
The same is true for the pair of bifunctors
 $(\otimes_S^\derL,\derR\! \uHom_S)$.
\end{rem}

\begin{num} \label{num:derived_premotivic_diagrams&descent}
Consider an abelian $\Pmor$-premotivic category $\A$ compatible
 with a topology $t$.
According to section \ref{sec:fibredmodcat},
 the $2$-functor $\Der(\A)$ can be extended
 to the category of $\sch$-diagrams:
\index{word}{diagram!Sdiagram@$\site$-diagram}
to any diagram of schemes $\X:I \rightarrow \sch$
 indexed by a small category $I$,
 we can associate a closed symmetric monoidal triangulated category $\Der(\A)(\X,I)$
 which coincides with $\Der(\A)(X)$ when $I=e$, $\X=X$ for a scheme $X$.

Let us be more specific.
The fibred category $\A$ admits an extension to $\sch$-diagrams:
 a section of $\A$ over a diagram of schemes $\X:I \rightarrow \sch$,
 indexed by a small category $I$, is the following data:
\begin{enumerate}
\item A family $(A_i)_{i \in I}$ such that $A_i$ is an object of $\A_{X_i}$.
\item A family $(a_u)_{u \in Fl(I)}$ such that for any arrow $u:i \rightarrow j$ in $I$,
 $a_u:u^*(A_j) \rightarrow A_i$ is a morphism in $\A_{X_i}$
 and this family of morphisms satisfies a cocyle condition
 (see paragraph \ref{paragr:extension_to_diagrams}).
\end{enumerate}
Then, $\Der(\A)(\X,I)$ is
 the derived category of the abelian category $\A(\X,I)$.
In particular, objects of $\Der(\A)(\X,I)$
 are complexes of sections of $\A$ over $(\X,I)$ (or, what amount to the same thing,
 families of complexes $(K_i)_{i \in I}$ with transition maps $(a_u)$
 as above, relative to the fibred category $\Comp(\A)$).

Recall that a morphism of $\sch$-diagrams $\varphi:(\X,I) \rightarrow (\Y,J)$
 is given by a functor $f:I \rightarrow J$ and a natural transformation
 $\varphi:\X \rightarrow \Y \circ f$.
We say that $\varphi$ is a $\Pmor$-morphism if
 for any $i \in I$, $\varphi_i:\X_i \rightarrow \Y_{f(i)}$ is a $\Pmor$-morphism.
For any morphism (resp. $\Pmor$-morphism) $\varphi$,
 we have defined in \ref{defiminvvarphi} adjunctions of (abelian) categories:
\begin{align*}
\varphi^*:\A(\Y,J)
 & \rightleftarrows \A(\X,I):\varphi_* \\
\text{resp. } \varphi_\sharp:\A(\X,I)
 & \rightleftarrows \A(\Y,J):\varphi^*
\end{align*}
which extends the adjunctions we had on trivial diagrams.

According to Proposition \ref{basicfunctdiag}, these respective adjunctions
 admits left/right derived functors as follows:
\begin{align} \label{eq:derived_premotivic_diag1}
\derL \varphi^*:\Der(\A)(\Y,J)
 & \rightleftarrows \Der(\A)(\X,I):\derR \varphi_* \\
\label{eq:derived_premotivic_diag2}
\text{resp. } \derL \varphi_\sharp:\Der(\A)(\X,I)
 & \rightleftarrows \Der(\A)(\Y,J):\derL \varphi^*=\varphi^*
\end{align}
Again, these adjunctions coincide on trivial diagrams
 with the map we already had.

Note also that the symmetric closed monoidal structure on $\Comp(\A(\X,I))$
 can be derived and induces a symmetric monoidal structure on $\Der(\A)(\X,I)$
 (see Proposition \ref{derivedtensordiaginj}).\footnote{In fact,
 $\Der(\A)$ is then a monoidal $\Pmor_\cart$-fibred category over the category
 of $\sch$-diagrams (remark \ref{remPacrtfibredcat}).}

Recall from \ref{defdescente} and \ref{trivglobchartaudescent}
that, given a topology $t'$ (not necessarily admissible) over $\sch$,
 we say that $\Der(\A)$ satisfies $t'$-descent
 if for any $t'$-hypercover $p:\X \rightarrow X$ 
 (here $\X$ is considered as a $\sch$-diagram),
 the functor
\begin{equation} \label{eq:derived_premotivic_descent}
\derL p^*:\Der(\A)(X) \rightarrow \Der(\A)(\X)
\end{equation}
is fully faithful
 (see Corollary \ref{trivglobchartaudescent}).
\end{num}
\begin{prop}\label{equivcohdescentabstractdescent}
Consider the notations and hypothesis introduced above.
Let $t'$ be an admissible topology on $\sch$.
Then the following conditions are equivalent:
\begin{enumerate}
\item[(i)] $\Der(\A)$ satisfies $t'$-descent,
\index{word}{descent!tdescent@$t$-descent}
in the sense recalled above.
\item[(ii)] $\A$ satisfies cohomological $t'$-descent.
\index{word}{descent!cohomological $t$-descent}
\end{enumerate}
\end{prop}
\begin{proof}
We prove (i) implies (ii).
Consider a $t'$-hypercover $p:\X \rightarrow X$ in $\Pmorx S$.
This is a $\Pmor$-morphism. 
Thus, by the fully faithfulness of \eqref{eq:derived_premotivic_descent},
the counit map
$\derL p_\sharp p^* \rightarrow 1$ is an isomorphism.
By applying the latter to the unit object $\un_X$ of $\Der(\A_X)$, we thus obtain that
$$
\mab X \X  \To \un_X
$$
is an isomorphism in $\Der(\A_X)$. If $\pi:X \rightarrow S$ is the structural
$\Pmor$-morphism, by applying the functor $\derL \pi_\sharp$ to this isomorphism,
we obtain that
$$
\mab S \X \To \mab S X
$$
is an isomorphism in $\Der(\A_S)$ and this concludes.

Reciprocally, to prove (i), we can restrict to $t'$-hypercovers
 $p:\X \rightarrow X$ which are $\Pmor$-morphisms because $t'$ is admissible.
Because $\derR p^*=p^*$ admits a left adjoint $\derL p_\sharp$,
 we have to prove that the counit
$$
\derL p_\sharp p^* \rightarrow 1
$$
is an isomorphism. 
This is a natural transformation between triangulated functors
 which commutes with small sums. Thus, according to \eqref{eq:premotivic_derived_generated},
 we have only to check this is an isomorphism when evaluated at a complex
 of the form $\mab X Y\{i\}$ for a $\Pmor$-scheme $Y/X$ and a twist $i \in \tau$.
But the resulting morphism is then
$\mab X {\X \times_X Y}\{i\} \rightarrow \mab X Y\{i\}$ and we can conclude
because $\X \times_X Y \rightarrow Y$ is a $t'$-hypercover in $\Pmorx S$ (again
 because $t'$ is admissible).
\end{proof}

\begin{num} \label{num:premotivic_derived_lim}.
Consider the situation of \ref{num:derived_premotivic_diagrams&descent}
Let $S$ be a scheme.
 An interesting particular case is given for constant $\sch$-diagrams 
 over $S$; for a small category $I$, we let $I_S$ be the constant $\sch$-diagram
 $I \rightarrow \sch, i \mapsto S, u \mapsto 1_S$.
Then the adjunctions \eqref{eq:derived_premotivic_diag1} for this kind of diagrams
define a \emph{Grothendieck derivator}
\index{word}{derivator, Grothendieck}
%%\footnote{
%%One can also obtain this result more directly by using the main theorem of \cite{Cis1}.
%% We refer to this article for further informations on derivators.}
 $$I \mapsto \Der(\A)(I_S).$$
Recall that,
 if $f:I \rightarrow e$ is the canonical functor to
 the terminal category and $\varphi=f_X:I_X \rightarrow X$
 the corresponding morphism of $\sch$-diagrams,
 for any $I$-diagram $K_\bullet=(K_i)_{i \in I}$ of complexes over $\A_S$,
 we get right derived limits and left derived colimits:
\begin{align*}
\derR \varphi_*(K_\bullet)&=\derR \plim_{i \in I} K_i. \\
\derL \varphi_\sharp(K_\bullet)&=\derL \ilim_{i \in I} K_i.
\end{align*}
\end{num}

%\begin{ex}
%Consider the above extension and assume $\A$ is $\tau$-generated.
%Let $K$ be an object of $\Der(\A_S)$ for a scheme $S$.
%Consider $K$ as an object in $\Comp(\A_S)$
% and put $G_K$ be the subcategory of $\Comp(\A_S)/K$ made of the objects
% of the form 
%$$
%\Mab S X \A\{i\}[n] \rightarrow K
%$$
%for a $\Pmor$-scheme $X/S$ and a couple $(i,n) \in \tau \times \ZZ$.
%Note that $G_K$ is essentially small because $\Pmorx S$ is.
%By assumption, we obtain
%$$
%K=\ilim \Mab S X \A\{i\}[n]
%$$
%where the limit is taken over $G_K$.
%
%As a complex, $K$ can be written as
%$$
%$$
%\end{ex}

\begin{num} \label{num:derived_adjunction}
The associated derived $\Pmor$-premotivic category is functorial in the following sense.

Consider an adjunction 
$$
\varphi:\A \rightleftarrows \B:\psi
$$
of abelian $\Pmor$-premotivic categories. 
Let $\tau$ (resp. $\tau'$) be the set of twists of $\A$ (resp. $\B$),
and recall that $\varphi$ induces a morphism of monoid
 $\tau \rightarrow \tau'$ still denoted by $\varphi$.
Consider two topologies $t$ and $t'$ such that $t'$ is finer than $t$.
Suppose $\A$ (resp. $\B$) is compatible with $t$ (resp. $t'$)
 and let $(\cG^\A_S,\cH^\A_S)$ (resp. $(\cG^\B_S,\cH^\B_S)$) be the descent
 structure on $\A_S$ (resp. $\B_S$) defined in \ref{t-descent_model_structure}.

For any scheme $S$, consider the evident extensions
$$
\varphi_S:\Comp(\A_S) \rightleftarrows \Comp(\B_S):\psi_S
$$
of the above adjoint functors to complexes.
Recall that for any $\Pmor$-scheme $X/S$ and any twist $i \in \tau$,
 $\varphi_S(\Mab S X \A\{i\})=\Mab S X \B\{\varphi(i)\}$ by definition.
Thus, $\varphi_S$ sends $\cG_S^\A$ to $\cG_S^\A$.
Because $t'$ is finer than $t$, it sends also $\cH_S^\A$ to $\cH_S^\B$.
In other words, it satisfies descent in the sense of \cite[par. 2.4]{CD1}
 so that the pair $(\varphi_S,\psi_S)$ is a Quillen adjunction 
with respect to the respective $t$-descent and $t'$-descent
model structure on $\Comp(\A_S)$ and $\Comp(\B_S)$.

Considering the derived functors, 
it is now easy to check we have obtained a $\Pmor$-premotivic
 adjunction\footnote{
Remark also that this adjunction extends on $\sch$-diagrams
 considering the situation described in
 \ref{num:derived_premotivic_diagrams&descent}:
for any diagram $\cX:I \rightarrow \sch$, we get an adjunction
$$
\derL \varphi_\cX:\Der(\A)(\cX) \rightleftarrows \Der(\B)(\cX):\derR \psi_\cX
$$
and this defines a morphism of triangulated monoidal $\Pmor_\cart$-fibred
 categories over the $\sch$-diagrams 
 (\textit{cf.} Proposition \ref{derivedPfibredQuillenfunct}).}
$$
\derL \varphi:\Der(\A) \rightleftarrows \Der(\B):\derR \psi.
$$
\end{num}

\begin{ex}
Let $t$ be an admissible topology.
Consider an abelian $\Pmor$-premotivic category $\A$
 compatible with $t$.
Then the morphism of abelian $\Pmor$-premotivic categories
 \eqref{eq:universality_presheaves}
 induces a morphism of triangulated $\Pmor$-premotivic categories:
\begin{equation} \renewcommand{\Rc}{\ZZ}
\label{eq:universal_derived_presheaves}
\derL \gamma^*:\Der(\psh{\Pmor})
 \rightleftarrows \Der(\A):\derR \gamma_*
\end{equation}
Similarly, the morphism \eqref{eq:abelian_premotivic_sheaves}
 induces a morphism of triangulated $\Pmor$-premotivic categories
\begin{equation}
\label{ex:associated_sheaf=adjunction}
a^*_t:\Der(\psh{\Pmor})
 \rightleftarrows
  \Der(\sh t {\Pmor}):\derR a_{t,*}.
\end{equation}
Note that $a^*_t=\derL a^*_t$ on objects, because the functor $a^*_t$ is exact.
\end{ex}

\begin{ex} \label{ex:derived_chg_coef}
\renewcommand{\Rc}{\Lambda}
Consider an admissible topology $t$.
Let $\varphi:\Lambda \rightarrow \Lambda'$ be a morphism
 of rings.
 For any scheme $S$, it induces a pair of adjoint functors:
\begin{equation} \label{eq:abelian_chg_coef}
\varphi^*:\sh t {\Pmor_S} \rightleftarrows
 \renewcommand{\Rc}{\Lambda'} 
 \sh t {\Pmor_S}:\varphi_*
\end{equation}
such that $\varphi^*$ (resp. $\varphi_*$)
 is induced by the obvious extension (resp. restriction) of scalars
 functor.
By definition, for any $\Pmor$-scheme $X/S$,
 the functor $\varphi^*$ sends the representable sheaf of
 $\Lambda$-modules $\repx t S X$ to the representable sheaf 
 of $\Lambda'$-modules 
 $\renewcommand{\Rc}{\Lambda^{\prime t}} \rep S X$.
 Thus $(\varphi^*,\varphi_*)$ defines an adjunction of
 abelian $\Pmor$-premotivic categories.
 Applying the results of Paragraph \ref{num:derived_adjunction},
 one deduces a $\Pmor$-premotivic adjunction:
$$
\derL \varphi^*:\Der(\sh t {\Pmor})
 \rightleftarrows
 \renewcommand{\Rc}{\Lambda'} 
 \Der(\sh t {\Pmor}):\derR \varphi_*.
$$
The functor $\varphi_*$ is exact so that $\derR \varphi_*=\varphi_*$.
Similarly when $\Lambda'/\Lambda$ is flat,
 $\derL \varphi^*=\varphi^*$.
\end{ex}

The following result can be used to check the compatibility
 to a given admissible topology:
\begin{prop} \label{adjunction&topology}
Let $t$ be an admissible topology.
Consider a morphism of abelian $\Pmor$-premotivic categories
$$
\varphi:\A \rightleftarrows \B:\psi
$$
such that:
\begin{itemize}
\item[(a)] For any scheme $S$, $\psi_S$ is exact. 
\item[(b)] The morphism $\varphi$ induces an isomorphism
 of the underlying set of twists of $\A$ and $\B$.
\end{itemize}
According to the last property, we identify the set of twists
 of $\A$ and $\B$ to a monoid $\tau$ in such a way that $\varphi$ acts on $\tau$
 by the identity.

Assume that $\A$ is compatible with $t$. Then the following conditions
 are equivalent:
\begin{enumerate}
\item[(i)] $\B$ is compatible with $t$.
\index{word}{compatible with (a topology) $t$}
\item[(ii)] $\B$ satisfies cohomological $t$-descent,
\end{enumerate}
\end{prop}
\begin{proof}
The fact $(i)$ implies $(ii)$ is clear from the definition,
 and we prove the converse using the following lemma:
\begin{lm} \label{abelian_adjunction&derived_Hom}
Consider a morphism of $\Pmor$-premotivic abelian categories
$$
\varphi:\A \rightleftarrows \B:\psi
$$
satisfying conditions (a) and (b) of the above proposition
 and a base scheme $S$.

Given a simplicial scheme $\cX$ which is degree-wise a sum of $\Pmor$-schemes
 over $S$,
 a twist $i \in \tau$ and a complex $C$ over $\B_S$,
we denote by
\begin{equation*} \label{eq:abelian_adjunction&derived_Hom1}
\epsilon_{\cX,i,C}:\Hom_{\Comp(\B_S)}\big(\Mab S \cX \B\{i\},C\big)
 \rightarrow
 \Hom_{\Comp(\A_S)}\big(\Mab S \cX \A\{i\},\psi_S(C)\big)
\end{equation*}
the adjunction isomorphism obtained for the adjoint pair $(\varphi_S,\psi_S)$. \\
Then there exists a unique isomorphism $\epsilon'_{\cX,i,C}$
 making the following diagram commutative:
\begin{equation*} \label{diag:abelian_adjunction&derived_Hom2}
\begin{split}
\xymatrix@C=40pt{
\Hom_{\Comp(\B_S)}\big(\Mab S \cX \B\{i\},C\big)\ar^-{\epsilon_{\cX,i,C}}[r]\ar[d]
 & \Hom_{\Comp(\A_S)}\big(\Mab S \cX \A\{i\},\psi_S(C)\big)\ar[d] \\
\Hom_{\K(\B_S)}\big(\Mab S \cX \B\{i\},C\big)\ar^-{\epsilon'_{\cX,i,C}}[r]
 & \Hom_{\K(\A_S)}\big(\Mab S \cX \A\{i\},\psi_S(C)\big).
}
\end{split}
\end{equation*}

Assume moreover that $\B$ satisfies cohomological $t$-descent. \\
Then there exists an isomorphism $\epsilon''_{\cX,i,C}$
 making the following diagram commutative:
\begin{equation} \label{diag:abelian_adjunction&derived_Hom3}
\begin{split}
\xymatrix@C=40pt{
\Hom_{\K(\B_S)}\big(\Mab S \cX \B\{i\},C\big)\ar^-{\epsilon'_{\cX,i,C}}[r]
 \ar_{\pi^\B_{\cX,i,C}}[d]
 & \Hom_{\K(\A_S)}\big(\Mab S \cX \A\{i\},\psi_S(C)\big)\ar^{\pi^\A_{\cX,i,C}}[d] \\
\Hom_{\Der(\B_S)}\big(\Mab S \cX \B\{i\},C\big)\ar^-{\epsilon''_{\cX,i,C}}[r]
 & \Hom_{\Der(\A_S)}\big(\Mab S \cX \A\{i\},\psi_S(C)\big),
}
\end{split}
\end{equation}
where $\pi^\A_{\cX,i,C}$ and $\pi^\B_{\cX,i,C}$ are induced by the
 obvious localization functors.
\end{lm}
The existence and unicity of isomorphism $\epsilon'_{\cX,i,C}$
follows from
the fact that the functors $\varphi_S$ and $\psi_S$ are additive.
Indeed, this implies that the isomorphism $\epsilon_{\cX,i,C}$
 is compatible with chain homotopies. 

Consider the injective model structure on $\Comp(\A_S)$ and 
$\Comp(\B_S)$ (see for example \cite[1.2]{CD1} for the definition).
We first treat the case when $C$ is fibrant for this model structure
 on $\Comp(\B_S)$.
Because the premotive $\Mab S \cX \B\{i\}$ is cofibrant for the injective model structure,
 we obtain that the canonical map $\pi^\B_{\cX,i,C}$
 is an isomorphism.
This implies there exists
 a unique map $\epsilon''_{\cX,i,C}$ making
 diagram \eqref{diag:abelian_adjunction&derived_Hom3} commutative.
On the other hand,
 the isomorphism $\epsilon'_{\cX,i,C}$ obtained previously is 
 obviously functorial in $\cX$.
Thus, because $\B$ satisfies $t$-descent,
 we obtain that $\psi_{S}(C)$ is $t$-flasque.
Because $\A$ is compatible with $t$, this implies $\psi_S(C)$ is $t$-local,
 and because $\Mab S \cX \B\{i\}$ is cofibrant for the $t$-descent model
 structure on $\Comp(\A_S)$, this implies $\pi^\B_{\cX,i,C}$ is an isomorphism.
 Thus finally, $\epsilon''_{\cX,i,C}$ is an isomorphism as required.

To treat the general case, we consider a fibrant resolution $C \rightarrow D$
 for the injective model structure on $\Comp(\B_S)$.
 Because $\psi_S$ is exact, it preserves isomorphisms. Using the previous case,
 We define $\epsilon''_{\cX,i,C}$ by the following commutative diagram:
$$
\xymatrix@C=40pt{
\Hom_{\Der(\B_S)}\big(\Mab S \cX \B\{i\},C\big)\ar^-{\epsilon''_{\cX,i,C}}[r]\ar_\sim[d]
 & \Hom_{\Der(\A_S)}\big(\Mab S \cX \A\{i\},\psi_S(C)\big)\ar^\sim[d] \\
\Hom_{\Der(\B_S)}\big(\Mab S \cX \B\{i\},D\big)\ar^-{\epsilon''_{\cX,i,D}}[r]
 & \Hom_{\Der(\A_S)}\big(\Mab S \cX \A\{i\},\psi_S(D)\big).
}
$$
The required property for $\epsilon''_{\cX,i,C}$ then follows easily and the lemma is proved.

\bigskip

To finish the proof that (ii) implies (i), we note the
 lemma immediately implies, under (ii),
 that the following two conditions are equivalent:
\begin{itemize}
\item $C$ is $t$-flasque (resp. local)
in $\Comp(\B_S)$;
\item $\psi_S(C)$ is $t$-flasque (resp. local)
in $\Comp(\A_S)$.
\end{itemize}
This concludes.
\end{proof}

\subsubsection{Constructible premotivic complexes}
\label{num:derived_constructible}

\begin{df}
\label{df:hyper-bounded}
Let $\A$ be an abelian $\Pmor$-premotivic category
compatible with an admissible topology $t$.
We will say that $t$ is \emph{bounded\index{word}{bounded (topology)} in $\A$}
 if for any scheme $S$,
there exists an essentially small family $\cN_S^t$ 
of bounded complexes which are direct factors of
finite sums of objects of type $\mab S X\{i\}$ in each degree,
such that, for any complex $C$ of $\A_{S}$, the following
conditions are equivalent.
\begin{itemize}
\item[(i)] $C$ is $t$-flasque.
\item[(ii)] For any $H$ in $\cN_S^t$, the abelian group
$\Hom_{\Htp(\A_{S})}(H,C)$ vanishes.
\end{itemize}
In this case, we say the family $\cN^t_S$ is a \emph{bounded
generating family\index{word}{bounded generating family} for
 $t$-hyper\-co\-verings in $\A_{S}$}.
\end{df}

\begin{ex} \label{ex:Zar&Nis=bounded}
\begin{enumerate}
\item Assume $\Pmor$ contains the open immersions
 so that the Zariski topology is admissible.
Let $MV_S$ to be the family of complexes of the form
$$
\Rc_S(U \cap V) \xrightarrow{l_*-k_*} \Rc_S(U) \oplus \Rc_S(V)
 \xrightarrow{i_*+j_*} \Rc_S(X)
$$
for any open cover $X=U \cup V$, where $i$,$j$,$k$,$l$ denotes the obvious
open immersions.
It follows then from \cite{BG} that
$MV_{S}$ is a bounded generating family of Zariski hypercovers in $\sh \zar {\Pmorx S}$.
\item
Assume $\Pmor$ contains the \'etale morphisms 
 so that the Nisnevich topology is admissible.
We let $BG_S$ be the family of complexes
of the form
$$
\Rc_S(W) \xrightarrow{g_*-l_*} \Rc_S(U) \oplus \Rc_S(V)
 \xrightarrow{j_*+f_*} \Rc_S(X)
$$
for a Nisnevich distinguished\index{word}{square!Nisnevich distinguished}
 square in $\sch$
 (\textit{cf.} \ref{ex:lower&upper_cd_structures})
$$
\xymatrix@=14pt{
W\ar^l[r]\ar_g[d] & V\ar^f[d] \\ U\ar^j[r] & X.
}
$$
Then, by applying \ref{BGNis}, we see that
$BG_{S}$ is a bounded
generating family for Nisnevich hypercovers in $\sh \nis {\Pmorx S}$.
\item Assume that $\Pmor=\sft$ is the class of morphisms of finite type in $\sch$.
We let $PCDH_S$ be the family of complexes of the form
$$
\Rc_S(T) \xrightarrow{g_*-k_*} \Rc_S(Z) \oplus \Rc_S(Y)
 \xrightarrow{i_*+f_*} \Rc_S(X)
$$
for a $\cdh$-distinguished\index{word}{square!cdhdistinguished@$\cdh$-distinguished} square in $\sch$
 (\textit{cf.} \ref{ex:lower&upper_cd_structures})
$$
\xymatrix@=14pt{
T\ar^k[r]\ar_g[d] & Y\ar^f[d] \\ Z\ar^i[r] & X.
}
$$
Then, by virtue of \ref{BGcdh}, $CDH_S=BG_{S} \cup PCDH_S$ is a bounded
generating family for $\cdh$-hypercovers in $\sh \cdh {\sft/S}$.
\item The \'etale topology is not bounded in $\sh {\et} {\sm}$
for an arbitrary ring $\Rc$.
However, if $\Rc=\QQ$, it is bounded:\renewcommand{\Rc}{\QQ}
by virtue of Theorem \ref{ratdescent},
a bounded generating family for \'etale hypercovers in $\sh {\et} {\sm}_S$
%%the category of \'etale sheaves of $\QQ$-vector spaces
is the union of the class $BG_S$ and that of complexes 
of the form $\Rc_S(Y)_G\To\Rc_S(X)$ for any Galois cover $Y\to X$ of group $G$.
%is just $BG_S$ (to be more precise, we should add the complexes of the form
%$\Rc_S(Y)_G\To\Rc_S(X)$
%for any Galois cover $Y\to X$ of group $G$, but
%these are contractible (because $Y/G\simeq X$),
%so that we can just drop them).
\item As in the case of the \'etale topology, the $\qfh$-topology
is not bounded in general, but it is so with rational coefficients.
 Let $PQFH_S$ be the family of complexes of the form
$$
\Rc_S(T)_G \xrightarrow{g_*-k_*} \Rc_S(Z) \oplus \Rc_S(Y)_G
 \xrightarrow{i_*+f_*} \Rc_S(X)
$$
for a $\qfh$-distinguished square
\index{word}{square!qfhdistinguished@$\qfh $-distinguished}
 of group $G$ in $\sch$ (\textit{cf.} \ref{settinghdesc})
$$
\xymatrix@=14pt{
T\ar^k[r]\ar_g[d] & Y\ar^f[d] \\ Z\ar^i[r] & X.
}
$$
Then, by virtue of Theorem \ref{carratqfhdescent},
$QFH_S=PQFH_S\cup BG_S$ is
a bounded generating family for $\qfh$-hypercovers in $\sh {\qfh} {\sft/S}$.
\item Similarly, by Theorem \ref{hdesceqqfhcdhdesc}, $H_S=CDH_S\cup QFH_S$ is
a bounded generating family for $\h$-hypercovers in $\sh {\h} {\sft/S}$.
\end{enumerate}
\end{ex}
\renewcommand{\Rc}{\Lambda}

\begin{prop}\label{boundedcompactgen}
Let $\A$ be an abelian $\Pmor$-premotivic category compatible
 with an admissible topology $t$.
We make the following assumptions:
\begin{itemize}
\item[(a)] $t$ is bounded in $\A$;
\item[(b)] for any $\Pmor$-morphism $X\To S$ and any $n\in\tau$,
the functor $\Hom_{\A_S}(M_S(X)\{n\},-)$ preserves filtered colimits.
\end{itemize}
Then $t$-local complexes are stable by filtering colimits.
\end{prop}

\begin{proof}
Let $\cN^t_S$ is a bounded
generating family for $t$-hypercovers in $\A_{S}$.
Then a complex $C$ of $\A_{S}$ is $t$-flasque if and only if
for any $H\in\cN^t_S$, the abelian group
$\Hom_{\Htp(\A_{S})}(H,C)$ is trivial.
Hence it is sufficient to prove that the functor
$$C\mapsto\Hom_{\Htp(\A_{S})}(H,C)$$
preserves filtering colimits of complexes.
This will follow from the fact that the functor
$$C\mapsto\Hom_{\Comp(\A_{S})}(H,C)$$
preserves filtering colimits. As $H$ a is bounded complex that is degreewise compact,
this latter property is obvious.
\end{proof}

\begin{num} \label{num:compact_objects_derived}
Consider an abelian $\Pmor$-premotivic category $\A$
 compatible with an admissible topology $t$,
 with generating set of twists $\tau$.
Assume that $t$ is bounded in $\A$ and consider a bounded
generating family $\N^t_S$ for $t$-hypercovers in $\A_S$.

Let $\Mab {} {\Pmorx S} \A$ be the full subcategory of $\A_S$
spanned by direct factors of finite sums of premotives of shape $\mab S X \{i\}$
for a $\Pmor$-scheme $X/S$ and a twist $i \in \tau$.
This category is additive and we can associate with
it its category of complexes up to chain homotopy.
We get an obvious triangulated functor
\begin{equation} \label{eq:derived_geometric_premotives_proto}
K^b\big(\Mab {} {\Pmorx S} \A\big) \rightarrow \Der(\A_S).
\end{equation}
Then the previous functor induces a triangulated functor
$$
K^b\big(\Mab {} {\Pmorx S} \A\big)/\cN^t_S
 \rightarrow \Der(\A_S)
$$
where the left hand side stands for the Verdier
quotient of $K^b\big(\Mab {} {\Pmorx S} \A\big)$
by the thick subcategory generated by $\cN^t_S$.

The category $K^b\big(\Mab {} {\Pmorx S} \A\big)/\cN^t_S$
may not be pseudo-abelian while the aim of the previous
functor is. Thus we can consider its pseudo-abelian 
envelope and the induced functor
\begin{equation} \label{eq:derived_geometric_premotives}
\Big(K^b\big(\Mab {} {\Pmorx S} \A\big)/\cN^t_S\Big)^\natural
 \rightarrow \Der(\A_S).
\end{equation}
According to Definition \ref{df:tau-geometric},
 the image of this functor is the subcategory
  of $\tau$-constructible premotives of
 the triangulated $\Pmor$-premotivic category $\Der(\A_{S})$.
Then the following proposition is a corollary of 
\cite[theorem 6.2]{CD1}:
\end{num}
\begin{prop}
\label{compact_objects_derived}
Consider the hypothesis and notations above.

If $\A$ is finitely $\tau$-presented
\index{word}{finitely presented!finitely $\tau$-presented}
 then $\Der(\A)$ is compactly $\tau$-generated.
\index{word}{generated!compactly $\tau$-generated!triangulated $\Pmor$-fibred}
Moreover,
 the functor \eqref{eq:derived_geometric_premotives} is fully faithful.
\end{prop}

Let us denote by $\Der_{c}(\A)$ the subcategory of $\Der(\A)$
 made of $\tau$-constructible premotives in the
 sense of Definition \ref{df:tau-geometric}.
 Taking into account Proposition \ref{constructequivcompact},
  the previous proposition admits the following corollary:
\begin{cor} \label{cor:Der_compact&constructible}
Consider the situation of \ref{num:compact_objects_derived},
 and assume that $\A$ is finitely $\tau$-presented.
For any premotive $\cM$ in $\Der(\A_S)$,
 the following conditions are equivalent:
\begin{enumerate}
\item[(i)] $\cM$ is compact.\index{word}{compact}
\item[(ii)] $\cM$ is $\tau$-constructible.
\index{word}{constructible!$\tau$-constructible}
\end{enumerate}
Moreover, the functor \eqref{eq:derived_geometric_premotives}
 induces an equivalence of categories:
$$
\Big(K^b\big(\Mab {} {\Pmorx S} \A\big)
      /\cN^t_S \Big)^\natural
 \rightarrow \Der_{c}(\A_S).
$$
\end{cor}

\begin{ex} \label{ex:Nis_Et_compacity}
According to example \ref{ex:Zar&Nis=bounded},
 we get the following examples:
\begin{enumerate}
%\item
%Let $\Rc(S_\zar)$ be the $\Rc$-linear envelope of the category of
%open subschemes of $S$.
%We obtain a fully faithful functor
%$$
%\left(\K^b\left(\Rc(S_\zar)\right)/MV_S\right)^\natural
% \rightarrow \Der(\sh {\zar} {S_\zar}).
%$$
%which is essentially surjective on compact objects.
\item
Let $\Rc(\sm/S)=\Mab{}{\sm/S}{\A}$ for
$\A=\sh {\nis} {\sm/S}$.
We obtain a fully faithful functor
$$
\left(\K^b\left(\Rc(\sm/S)\right)/BG_S\right)^\natural
 \rightarrow \Der\big(\sh {\nis} {\sm/S}\big)
$$
which is essentially surjective on compact objects.
\item Let $\Rc(\sft/S)=\Mab{}{\sm/S}{\A}$ for
$\A=\sh {\cdh} {\sft/S}$.
We obtain a fully faithful functor
$$
\left(\K^b\big(\Rc(\sft/S)\big)/BG_S \cup CDH_S\right)^\natural
 \rightarrow \Der\Big(\sh {\cdh} {\smash{\sft/S}}\Big)
$$
which is essentially surjective on compact objects.
\renewcommand{\Rc}{\QQ}
\item Let $\Rc_\et(\sm/S)=\Mab{}{\sm/S}{\A}$ for
$\A=\sh {\et} {\sm/S}$.
We obtain a fully faithful functor
$$
\left(\K^b\left(\Rc_\et(\sm/S)\right)/BG_S\right)^\natural
 \rightarrow \Der\big(\sh {\et} {\sm/S}\big).
$$
which is essentially surjective on compact objects.
\end{enumerate}
\end{ex}
\renewcommand{\Rc}{\Lambda}

\begin{num} \label{num:continuity_derived_premotivic}
Consider an abelian $\Pmor$-premotivic category $\A$.
 We introduce the following property of $\A$:
\begin{itemize}
\item[(C)] Consider a projective system $(S_\alpha)_{\alpha\in A}$
\index{word}{projective system, of schemes}
 of schemes in $\sch$ with affine transition maps such that
 $S=\varprojlim_{\alpha\in A}S_\alpha$ belongs to $\sch$.
 For any index $\alpha_0 \in A$, any object $A_{\alpha_0}$
 in $\A_{S_{\alpha_0}}$, and any twist $n \in \tau$,
 the canonical map
$$
\varinjlim_{\alpha \in A/\alpha_0} \Hom_{\A_{S_\alpha}}(\un_{S_\alpha}\{n\},A_\alpha)
 \rightarrow \Hom_{\A_S}(\un_{S}\{n\},A)
$$
is an isomorphism where $A_\alpha$ (resp. $A$)
 is the pullback of $A_{\alpha_0}$ along the canonical map
 $S_\alpha \rightarrow S_{\alpha_0}$
  (resp. $S \rightarrow S_{\alpha_0}$).
\end{itemize}
We will denote by (wC) the analogous property when one restricts
 pro-objects to thus with affine and dominant transition maps.
\end{num}
\begin{prop} \label{prop:continuity_derived_premotivic}
Consider an abelian $\Pmor$-premotivic category $\A$
 compatible with an admissible topology $t$
 and satisfying the assumption (C) (resp. (wC)) above.

Then the derived premotivic category $\Der(\A)$
 is $\tau$-continuous (resp. weakly $\tau$-continuous)
 --- see Definition \ref{df:continuous}.
\end{prop}
\begin{proof}
We use Proposition \ref{abstractcontinuity}
 applied to the $t$-descent model structure on $\Comp(\A_T)$
 for $T=S$ or $T=S_\alpha$.
 (see Paragraph \ref{num:derived&generators}).
 Recall from Paragraph \ref{t-descent_model_structure}
 that this model structure is associated with a descent structure.
 Thus according to \cite[2.3]{CD1},
 there exists an explicit generating set $I$ (resp. $J$)
 for cofibrations (resp. trivial cofibrations).
 Moreover, the source or target of any map in $I \cup J$ 
 is a complex $C$ satisfying the following assumption:
\begin{itemize}
\item[(rep)] for any integer $i \in \ZZ$,
$C^i$ is a sum of premotives of the form
 $\mab T X\{n\}$ where $X/T$ is a $\Pmor$-scheme and $n \in \tau$.
\end{itemize}
Thus,
 to check the assumption of \ref{abstractcontinuity} for $\Comp(\A)$, 
 we fix a projective system $(S_\alpha)_{\alpha \in A}$
  satisfying the assumptions of property (C) (resp. (wC)) above;
 we have to prove that for any index $\alpha_0 \in A$
 and any complexes $C_{\alpha_0}$ and $E_{\alpha_0}$ such that $C_{\alpha_0}$
 satisfies (rep), the natural map:
$$
\varinjlim_{\alpha\in A/\alpha_0}
 \Hom_{\Comp(\A_{S_\alpha})}(C_{\alpha},E_\alpha)
\To \Hom_{\Comp(\A_S)}(C,E)
$$
is bijective.

Given the definition of morphisms in a category of complexes,
 it is sufficient to check this when the Hom groups are computed
 as morphisms of $\ZZ$-graded objects.
 Thus it is sufficient to treat the case where $C_{\alpha_0}$
 and $E_{\alpha_0}$ are concentrated in degree $0$.
 Thus, as $C_{\alpha_0}$ satisfies property (rep), 
 we are exactly reduced to assumption (C) (resp. (wC)) on $\A$.
\end{proof}

\begin{ex} \label{ex:sheaves_continuity}
\begin{enumerate}
\item Assume $\Pmor$ is contained in the class of morphisms of finite type.

Then the abelian $\Pmor$-premotivic category $\psh{\Pmor}$ of example
 \ref{ex:abelian_premotivic_presheaves} satisfies assumption (C).
 Indeed, property (C) when $A$ is a representable presheaf 
 follows from the assumption on $\Pmor$: $\Pmor$-schemes over some base $S$
 always are of finite presentation over $S$ --
  $S$ is noetherian according to our general assumption
  \ref{num:assumption_sch_derived}.
 Then the case of a general presheaf $A$ follows because $A$
  is an inductive limit of representable presheaf and the global
  sections functor commutes with inductive limit of presheaves.
\item Let $\sft$ be the class of morphisms of finite type
 and let $t$ be one of the following topologies: 
 $\nis, \et, \cdh, \qfh, \h$. \\
Then the generalized abelian premotivic category $\sh t {\sft}$
 of example \ref{ex:abelian_premotivic_sheaves} satisfies assumption (C).

Indeed, according to the preceding example,
 we have only to prove that for any morphism $f:X \rightarrow S$,
 the functor 
$$
f^*:\psh{\sft_S} \rightarrow \psh{\sft_T}
$$
preserves the property of being a $t$-sheaf.

If $f$ is a morphism of finite type, the functor $f^*$
 admits as a left adjoint the functor $f_\sharp$, 
 which preserves $t$-covers. Thus the assertion is clear
 in that case.
 
In the general case, we use the fact that $X/S$ is a projective
 limit of a projective system $(X_\alpha)_{\alpha \in A}$
 where $X_\alpha$ is an $S$-scheme affine and of finite type over $S$.
 To check that for a $t$-sheaf $F$ over $S$,
  the presheaf $f^*(F)$ is a $t$-sheaf, we fix a $t$-cover
  $(W_i)_{i \in I}$ of $X$ in $\sft_X$.
 As $X$ is noetherian, we can assume $I$ is finite.
 Moreover, there exists an index $\alpha_0 \in A$ such that
  for the $t$-cover $(W_i)_{i \in I}$ can be lifted to $X_{\alpha_0}$.
 Then, using property (C) of $\psh{\sft}$ applied to $F$ and $(X_\alpha)$,
  we reduce to check that $f_\alpha^*(F)$ is a $t$-sheaf
  for $\alpha \geq \alpha_0$. This follows from the first case treated.
\item Let $\sm$ be the class of smooth morphisms and 
 $t$ be one of the topologies: $\nis, \et$. \\
 As we will see in Example \ref{ex:enlargementofsheaves},
 there exists a canonical enlargement of abelian premotivic categories
 (see \eqref{eq:enlargement_sheaves}):
$$
\rho_\sharp:\sh t {\sm} \rightleftarrows \sh t {\sft}:\rho^*.
$$
As the functor $\rho_\sharp$ is fully faithful and commutes
 with $f^*$ for any morphism of schemes $f$, we deduce
 from the preceding point that the abelian premotivic category
 $\sh t {\sm}$ satisfies the above condition (C).

As an application of the previous proposition,
 we thus obtain that
  the derived premotivic category $\Der(\sh t {\sm})$
  is $\tau$-continuous.
\end{enumerate}
\end{ex}
\subsection{The $\AA^1$-derived premotivic category} \label{sec:A^1-derived}

\subsubsection{Localization of triangulated premotivic categories}
\label{sec:loc_derived}

\renewcommand{\Rc}{\ZZ}

\begin{paragr}\label{par:t-model_comp_localized}
Let $\A$ be an abelian $\Pmor$-premotivic category compatible 
with an admissible topology $t$
and $\Der(\A)$ be the associated derived $\Pmor$-premotivic category.

Suppose given an essentially small family of morphisms $\W$
in $\Comp(\A)$ which is stable by the operations $f^*$, $f_\sharp$
(in other words, $\W$ is a sub-$\Pmor$-fibred category of $\Comp(\A)$).
Remark that the localizing subcategory $\cT$
of $\Der(\A)$ generated by the cones of arrows in $\W$ is again stable 
by these operations. Moreover, as for any $\Pmor$-morphism $f:X \rightarrow S$
 we have $f_\sharp f^*=\mab S X \otimes_S(-)$, the category $\cT$
is stable by tensor product with a geometric section.

%It follows from the Brown representability theorem proved
%by Neeman~\cite{Nee1} and from the fact $\Der(\A_S)$
%is well generated that the canonical functor
%$$
%\pi:\Der(\A_S) \rightarrow \Der(\A_S)[\W_S^{-1}]
%$$
%admits a right adjoint (and that $\Der(\A_S)[\W_S^{-1}]$
%is well generated, so that Brown representability
%holds in $\Der(\A_S)[\W_S^{-1}]$ as well).
%The following proposition is straightforward.
%\end{paragr}
%\begin{prop}
%\label{prop:localization_exists}
%Suppose we are under the previous hypothesis.
%
%Then there exists a unique structure of a $\Pmor$-premotivic
%triangulated category on $\Der(\A)[\W^{-1}]$ 
%such that the functor
%$$
%\pi:\Der(\A) \rightarrow \Der(\A)[\W^{-1}]
%$$
%is a morphism of triangulated $\Pmor$-premotivic categories.
%\end{prop}
%%The category $\Der(\A_S)[\W_S^{-1}]$ is identified
%%with the full subcategory of $\Der(\A_S)$ made of $\W$-local
%%complexes.
%
%\begin{paragr}
%\label{par:t-model_comp_localized}
%We use model category theory to describe
%the category obtained previously. Suppose that
%$\A$ is compatible with an admissible topology $t$.
%
We will say that a complex $K$ over $\A_S$ is \emph{$\W$-local}
\index{word}{local!$\W$-local}
 if for any object $T$ of $\cT$ and any integer $n \in \ZZ$,
 $\Hom_{\Der(\A_S)}(T,K[n])=0$.
A morphism of complexes $p:C \rightarrow D$ over $\A_S$ is a
 \emph{$\W$-equivalence}
\index{word}{equivalence!Wequivalence@$\W$-equivalence}
 if for any $\W$-local complex $K$ over $\A_S$,
 the induced map
$$
\Hom_{\Der(\A_S)}(D,K) \rightarrow \Hom_{\Der(\A_S)}(C,K)
$$
is bijective.

A morphism of complexes
over $\A_S$ is called a \emph{$\W$-fibration}
\index{word}{fibration!Wfibration@$\W$-fibration}
 if it is a $t$-fibration
with a $\W$-local kernel.
A complex over $\A_S$ will be called \emph{$\W$-fibrant}
 if it is $t$-local and $\W$-local.
\end{paragr}
As consequence of \cite[4.3, 4.11 and 5.6]{CD1}, we obtain:
\begin{prop}
\label{t-model_category_complexes_localized}
Let $\A$ be an abelian $\Pmor$-premotivic category compatible 
with an admissible topology $t$ and $\W$ be an essentially small
family of morphisms in $\Comp(\A)$ stable by $f^*$ and $f_\sharp$.

Then the category $\Comp(\A_S)$ is a proper closed symmetric monoidal 
category with the $\W$-fibrations as fibrations, 
the cofibrations as defined in \ref{t-descent_model_structure},
and the $\W$-equivalences as weak equivalences.
\end{prop}

The homotopy category associated with this model category 
will be denoted by $\Der(\A_S)[\W_S^{-1}]$.
It can be described as the Verdier quotient $\Der(\A_S)/\cT_S$.

In fact, the $\W$-local model category
\index{word}{model structure!Wlocal@$\W$-local}
 on $\Comp(\A_S)$
 is nothing else than the left Bousfield localization of the
 $t$-local model category structure.
As a consequence,
 we obtain an adjunction of triangulated categories:
\begin{equation} \label{eq:adjunction_bousfield_loc}
\pi_S:\Der(\A_S) \rightleftarrows \Der(\A_S)[\W_S^{-1}]:\cO_S
\end{equation}
such that $\cO_S$ is fully faithful with essential image
the $\W$-local complexes.
In fact, the model structure gives a functorial $\W$-fibrant resolution
$1 \rightarrow R_\W$
$$
R_\W:\Comp(\A_S) \rightarrow \Comp(\A_S) \, ,
$$
which induces $\cO_S$. \\
Note that the triangulated category $\Der(\A_S)[\W_S^{-1}]$ is
 generated by the complexes concentrated in degree $0$ of the form
 $\mab S X\{i\}$ --- or, equivalently, the $\W$-local complexes
 $R_\W(\mab S X \{i\})$ --- for a $\Pmor$-scheme $X$
 and a twist $i \in \tau$.
 
\begin{rem} Another very useful property is that $\W$-equivalences
 are stable by filtering colimits; see \cite[prop. 3.8]{CD1}.
 \end{rem}

\begin{num}
Recall from \ref{t-descent_model_structure_bis} that for any
morphism (resp. $\Pmor$-morphism) $f:T \rightarrow S$,
the functor $f^*$ (resp. $f_\sharp$) satisfies descent;
as it also preserves $\W$, it follows from \cite[4.9]{CD1} that the adjunction
\begin{align*}
f^*:\Comp(\A_S) &\rightarrow \Comp(\A_T):f_* \\
\text{(resp. } f_\sharp:\Comp(\A_S) &\rightarrow \Comp(\A_T):f^*)
\end{align*}
is a Quillen adjunction with respect
 to the $\W$-local model structures. 
This gives the following corollary.
\end{num}

\begin{cor} \label{cor:derived_premotivic_localization_exists}
The $\Pmor$-fibred category $\Comp(\A)$ with the $\W$-local model structure
 on its fibers defined above is a monoidal $\Pmor$-fibred model category,
 which is moreover stable, proper and combinatorial.
\end{cor}
We will denote by $\Der(\A)[\W^{-1}]$ the triangulated $\Pmor$-premotivic category
 whose fiber over a scheme $S$ is the homotopy category of the $\W_S$-local
 model category $\Comp(\A_S)$.
The adjunction \eqref{eq:adjunction_bousfield_loc} readily defines
 an adjuntion of triangulated $\Pmor$-premotivic categories
\begin{equation} \label{eq:localized_adjunction}
\pi:\Der(\A) \rightleftarrows \Der(\A)[\W^{-1}]:\cO.
\end{equation}
The $\Pmor$-fibred categories $\Der(\A)$ and $\Der(\A)[\W^{-1}]$
 are both $\tau$-generated (and this adjunction is compatible
 with $\tau$-twists in a strong sense).

\begin{rem}
For any scheme $S$,
 the category $\Der(\A_S)[\W_S^{-1}]$ is well generated and has
 a canonical dg-structure
\index{word}{dgstructure@dg-structure}
 (see also \ref{rem:DG-structure}).
\end{rem}
 
\begin{num} \label{num:loc_derived_premotivic_diagrams&descent}
With the notations above,
 let us put $\T=\Der(\A)[\W^{-1}]$  to clarify the following notations.
As in \ref{num:derived_premotivic_diagrams&descent},
 the fibred category $\T$ has a canonical extension
 to $\site$-diagrams
\index{word}{diagram!Sdiagram@$\site$-diagram}
 $\X:I \rightarrow \sch$.

If we define $\W_\X$ as the class of morphisms $(f_i)_{i\in I}$ in $\Comp(\A(\X,I))$
 such that for any object $i$, $f_i$ is a $\W$-equivalence,
 then $\T(X)$ is the triangulated category $\Der(\A(\X,I))[\W_\X^{-1}]$.

Again, this triangulated category is symmetric monoidal closed
 and for any morphism (resp. $\Pmor$-morphism) $\varphi:(\X,I) \rightarrow (\Y,J)$,
 we get (derived) adjunctions as in \ref{num:derived_premotivic_diagrams&descent}:
\begin{align} \label{eq:loc_derived_premotivic_diag1}
\derL \varphi^*:\T(\Y,J)
 & \rightleftarrows \T(\X,I):\derR \varphi_* \\
\label{eq:loc_derived_premotivic_diag2}
\text{(resp. } \derL \varphi_\sharp:\T(\X,I)
 & \rightleftarrows \T(\Y,J):\derL \varphi^*=\varphi^*)
\end{align}
In fact, $\T$ is then a complete monoidal $\Pmor_\cart$-fibred category
 over the category of diagrams of schemes and the adjunction \eqref{eq:localized_adjunction}
 extends to an adjunction of complete monoidal $\Pmor_\cart$-fibred categories.
\end{num}

\begin{ex} \label{ex:P-premotivic_shv&preshv}
\renewcommand{\Rc}{\Lambda}
Suppose we are under the hypothesis of Example
\ref{ex:associated_sheaf=adjunction}.

Let $\W_{t,S}$ denote the family of maps
which are of the form 
$\Rc_S(\cX) \rightarrow \Rc_S(X)$
for a $t$-hypercover $\cX \rightarrow X$ in $\Pmorx S$.
Then $\W_t$ is obviously stable by $f^*$ and $f_\sharp$.

Recall now that a complex of $t$-sheaves on $\Pmorx S$ is local
if and only if its $t$-hypercohomology and its hypercohomology
computed in the coarse topology agree 
(\textit{cf.} \ref{ex:t-local_sheaves}).

This readily implies that the adjunction considered in
Example \ref{ex:associated_sheaf=adjunction}
$$
a^*_t:\Der(\psh{\Pmor})\rightleftarrows\Der(\sh t {\Pmor}):\derR a_{t,*}
$$
induces an equivalence of triangulated $\Pmor$-premotivic categories
$$
\Der(\psh{\Pmor})[\W_t^{-1}]
 \rightleftarrows\Der(\sh t {\Pmor}).
$$
Recall $\derR a_{t,*}$ is fully faithful and identifies
 $\Der(\sh t S)$ with the full subcategory of
 $\Der(\psh S)$ made by $t$-local complexes.
\end{ex}

\begin{paragr} \label{descent_derived_premotives}
A triangulated $\Pmor$-premotivic category $(\T,M)$
 such that there exists:
\begin{enumerate}
\item an abelian $\Pmor$-premotivic category $\A$
 compatible with an admissible topology $t_0$ on $\sm$.
\item an essentially small family $\W$
 of morphisms in $\Comp(\A)$ stable by
  $f^*$ and $f_\sharp$
\item an adjunction of triangulated $\Pmor$-premotivic categories
$\Der(\A)[\W^{-1}] \simeq \T$
\end{enumerate}
will be called for short a \emph{derived $\Pmor$-premotivic category}.
\index{word}{derived!derived $\Pmor$-premotivic category}
According to convention \ref{num:assumption_sch_derived}(d)
 and from the above construction,
 $\T$ is $\tau$-generated for some set of twists $\tau$.
\footnote{
We will formulate in some remarks below universal properties
 of some derived $\Pmor$-premotivic categories. When doing so,
 we will restrict to morphisms of derived $\Pmor$-premotivic
 categories which can be written as 
$$
\derL \varphi:\Der(\A_1)[\W_1^{-1}] \rightarrow \Der(\A_2)[\W_2^{-1}]
$$
for a morphism $\varphi:\A_1 \rightarrow \A_2$ of abelian
$\Pmor$-premotivic categories compatible with suitable topologies.
More natural universal properties could be obtained if one considers
 the framework of dg-categories or triangulated derivator.} 

Let us denote simply by $M_S(X)$ the geometric sections of $\T$.
In this case, using the morphisms \eqref{eq:universal_derived_presheaves}
 and \eqref{eq:localized_adjunction}, we get
 a canonical morphism of triangulated $\Pmor$-premotivic categories:
\begin{equation} \label{eq:universal_derived_P-premotivic}
\varphi^*:\Der(\psh{\Pmor})
 \rightleftarrows \T:\varphi_*.
\end{equation}
By definition, for any premotive $\cM$, any scheme $X$ and any integer $n \in \ZZ$,
we get a canonical identification:
\begin{equation} \label{eq:compute_coh_by_psh}
\Hom_{\T(S)}(M_S(X),\cM[n])=H^n \Gamma(X,\varphi_*(\cM)).
\end{equation}
Given any simplicial scheme $\cX$,
 we put $M_S(\cX)=\varphi^*\big( \ZZ_S(\cX) \big)$, so
 that we also obtain:
\begin{equation} \label{eq:compute_scoh_by_psh}
\Hom_{\T(S)}(M_S(\cX),\cM[n])
 =H^n\big(\Totp \Gamma(\cX,\derR \gamma_*(\cM))\big).
\end{equation}
\end{paragr}
\begin{prop} \label{prop:descent&derived_P-premotivic}
Consider the above notations
 and $t$ an admissible topology.
The following conditions are equivalent.
\begin{enumerate}
\item[(i)] For any $t$-hypercover\index{word}{hypercover}
 $\cX \rightarrow X$ in $\Pmorx S$, the induced map 
 $M_S(\cX) \rightarrow M_S(X)$
 is an isomorphism in $\T(S)$.
 \item[(i$'$)]  For any $t$-hypercover $p:\cX \rightarrow X$ in $\Pmorx S$,
 the induced functor $\derL p^*:\T(X)\To\T(\cX)$ is fully faithful.
\item[(i$''$)] The triangulated $\Pmor$-premotivic category $\T$ satisfies $t$-descent \index{word}{descent!tdescent@$t$-descent}
(Definition \ref{defdescente}).
\item[(ii)] There exists an essentially unique map
 $\varphi_t^*:\Der(\sh t {\Pmorx S}) \rightarrow \T(S)$
  making the following diagram essentially commutative:
$$
\xymatrix@C=30pt@R=15pt{
\Der(\psh {\Pmorx S})\ar^-{\varphi^*}[r]\ar_{a_t}[d] & \T(S) \\
\Der(\sh t {\Pmorx S})\ar_-{\varphi^*_t}[ru]. &
}
$$
\item[(ii$^\prime$)] For any complex $C \in \Comp(\psh  {\Pmorx S})$ such that
 $a_t(C)=0$, $\varphi^*(C)=0$.
\item[(ii$''$)] For any map $f:C \rightarrow D$ in $\Comp(\psh {\Pmorx S})$ such that
 $a_t(f)$ is an isomorphism, $\varphi^*(f)$ is an isomorphism.
\item[(iii)] There exists an essentially unique map
 $\varphi_{t*}:\T(S) \rightarrow \Der(\sh t {\Pmorx S})$
  making the following diagram essentially commutative:
$$
\xymatrix@C=30pt@R=15pt{
\Der(\psh {\Pmorx S}) & \T(S)\ar_-{\varphi_*}[l]\ar^-{\varphi_{t*}}[ld] \\
\Der(\sh t {\Pmorx S})\ar^{\derR \cO_t}[u]. &
}
$$
\item[(iii$'$)] For any premotive $\cM$ in $\T(S)$,
 the complex $\varphi_*(\cM)$ is local.
 \item[(iii$''$)] For any premotive $\cM$ in $\T(S)$,
 any $\Pmor$-scheme $X/S$ and any integer $n \in \ZZ$,
$$
\Hom_{\T(S)}(M_S(X),\cM[n])=H^n_t(X,\varphi_*(\cM)).
$$
\end{enumerate}
When these conditions are fulfilled for any scheme $S$,
the functors appearing in (ii) and (iii) induce
a morphism of triangulated $\Pmor$-premotivic categories:
$$
\varphi_{t}^*:\Der(\sh t {\Pmor}) \rightleftarrows \T:\varphi_{t*}.
$$
\end{prop}
\begin{proof}
The equivalence between conditions $(i)$, $(i')$ and $(i'')$
is clear (we proceed as in the proof of \ref{equivcohdescentabstractdescent}).
The equivalences $(ii) \Leftrightarrow (ii') \Leftrightarrow (ii'')$
 and $(iii) \Leftrightarrow (iii')$ follows
 from example \ref{ex:P-premotivic_shv&preshv} and the definition
 of a localization.
The equivalence $(i) \Leftrightarrow (ii'')$ follows again
 from \emph{loc. cit.}
The equivalences $(i) \Leftrightarrow (iii') \Leftrightarrow (iii'')$
 follows finally from \eqref{eq:compute_coh_by_psh},
  \eqref{eq:compute_scoh_by_psh},
 and the characterization of a local complex of sheaves
  (\textit{cf.} \ref{ex:t-local_sheaves}).
\end{proof}

\begin{rem} \label{rem:descent&derived_P-premotivic}
The preceding proposition expresses the fact that
the category $\Der(\sh t {\Pmor})$ is the
universal derived $\Pmor$-premotivic
\index{word}{universal}
 category satisfying $t$-descent.
\end{rem}

\begin{num} \label{num:BG_derived-premotivic}
We end this section by making explicit two particular cases
 of the descent property for
 derived $\Pmor$-premotivic categories.

Consider a derived $\Pmor$-premotivic category $\T$
 with geometric sections $M$.
Considering any diagram $\cX:I \rightarrow \Pmorx S$ of $\Pmor$-schemes over $S$,
 with projection $p:\cX \rightarrow S$,
 we can associate a premotive in $\T$:
$$
M_S(\cX)=\derL p_\sharp(\un_S)=\derL \ilim_{i \in I} M_S(\cX_i).
$$
In particular, when $I$ is the category $\bullet \rightarrow \bullet$,
 we associate to every $S$-morphism  $f:Y \rightarrow X$ of $\Pmor$-schemes
 over $S$ a canonical\footnote{ 
 In fact, if $\T=\Der(\A)[\W^{-1}]$ for an abelian $\Pmor$-premotivic
 category $\A$, then we can define $M_S(X \rightarrow Y)$
 as the cone of the morphism of complexes (concentrated in degree $0$)
$\mab S X \xrightarrow{f_*} \mab S Y$.}
  \emph{bivariant premotive}
$$
M_S(X\xrightarrow f Y).
$$
When $f$ is an immersion, we will also write $M_S(Y/X)$ for this premotive.
Note that in any case, there is a canonical distinguished triangle in $\T(S)$:
$$
M_S(X) \xrightarrow{f_*} M_S(Y)
 \xrightarrow{\pi_f} M_S(X \xrightarrow f Y)
 \xrightarrow{\partial_f} M_S(X)[1].
$$
This triangle is functorial in the arrow $f$
 -- with respect to \emph{commutative} squares. \\
Given a commutative square of $\Pmor$-schemes over $S$
\begin{equation} \label{eq:BG_derived-premotivic}
\begin{split}
\xymatrix@=20pt{
B\ar^{e'}[r]\ar_g[d] & Y\ar^f[d] \\
A\ar^e[r] & X
}
\end{split}
\end{equation}
we will say that the image square in $\T(S)$
$$
\xymatrix@=20pt{
M_S(B)\ar^{e'_*}[r]\ar_{g_*}[d] & M_S(Y)\ar^{f_*}[d] \\
M_S(A)\ar^{e_*}[r] & M_S(X)
}
$$
is \emph{homotopy cartesian}\footnote{If $\T=\Der(\A)[\W^{-1}]$, this amount
 to say that the diagram obtained of complexes by applying the functor $\mab S -$
 is homotopy cartesian in the $\W$-local model category $\Comp(\A)$.}
\index{word}{homotopy cartesian}
 if the premotive associated with diagram \ref{eq:BG_derived-premotivic} is zero.
\end{num}
\begin{prop} \label{prop:BG_property_derived_premotivic}
Consider a derived $\Pmor$-premotivic category $\T$.
We assume that $\Pmor$ contains the \'etale morphisms
 (resp. $\Pmor=\sft$).
Then, with the above definitions, the following conditions are equivalent:
\begin{enumerate}
\item[(i)] $\T$ satisfies Nisnevich (resp. proper $\cdh$) descent.
\item[(ii)] For any scheme $S$ and
 any Nisnevich (resp. proper $\cdh$) distinguished square $Q$ of $S$-schemes,
 the square $M_S(Q)$ is homotopy cartesian in $\T(S)$.
\item[(iii)] For any Nisnevich (resp. proper $\cdh$) distinguished square
 of shape \eqref{eq:BG_derived-premotivic},
 the canonical map $M_S(Y/B) \xrightarrow{(f/g)_*} M_S(X/A)$ is an isomorphism.
\end{enumerate}
Moreover, under these conditions, to any Nisnevich (resp. proper $\cdh$)
 distinguished square $Q$
 of shape \eqref{eq:BG_derived-premotivic}, we associate a map
$$
\partial_Q:M_S(X) \xrightarrow{\pi_e} M_S(X/A) \xrightarrow{(f/g)_*^{-1}} M_S(Y/B)
 \xrightarrow{\partial_{e'}} M_S(Y)[1]
$$
which defines a distinguished triangle in $\T(S)$:
$$
M_S(B) \xrightarrow{\text{\tiny {$\begin{pmatrix} e'_* \\ -g_* \end{pmatrix}$}}}
 M_Z(Y) \oplus M_S(A) \xrightarrow{ (f_*,e_*) } M_S(X)
  \xrightarrow{\partial_Q} M_S(Y)[1].
$$
\end{prop}
\begin{proof}
The equivalence of (i) and (ii) follows from the theorem of Morel-Voevodsky
 \ref{BGNis} (resp. the theorem of Voevodsky \ref{BGcdh}).
To prove the equivalence of (ii) and (iii), we assume $\T=\Der(\A)[\W^{-1}]$.
Then,  the homotopy colimit of a square of shape \ref{eq:BG_derived-premotivic}
is given by the complex
$$
\mathrm{Cone}\big(\mathrm{Cone}(\mab S B \rightarrow \mab S Y )
 \rightarrow \mathrm{Cone}(\mab S A  \rightarrow \mab S X)\big).
$$
This readily proves the needed equivalence,
 together with the remaining assertion.  
\end{proof}

\begin{rem} In the first of the respective cases of the proposition,
 condition (ii) is what we usually called the \emph{Brown-Gersten} property \bg
 for $\T$, whereas condition (iii) can be called the \emph{excision property}.
 In the second respective case, condition (ii) will be called the
 \emph{proper $\cdh$} property for the \emph{generalized premotivic category} $\T$.
 We say also that $\T$ satisfies the \bgcdh property if it satisfies condition (ii)
 with respect to any $\cdh$ distinguished square $Q$. 
 \end{rem}

\subsubsection{The homotopy relation}

\begin{paragr}
Let $\A$ be an abelian $\Pmor$-premotivic category
compatible with an admissible topology $t$.

We consider $\W_{\AA^1}$ to be the family of morphisms
$\mab S {\AA^1_X}\{i\} \rightarrow \mab S X\{i\}$
for a $\Pmor$-scheme $X/S$ and a twist $i$ in $\tau$.
The family $\W_{\AA^1}$ is obviously stable by $f^*$ and $f_\sharp$.
\end{paragr}
\begin{df}
\label{df:effective_triangulated_premotives}
Let $\A$ be an abelian $\Pmor$-premotivic category
compatible with an admissible topology $t$.
With the notation above, 
we define $\DMue(\A)=\Der(\A)[\W_{\AA^1}^{-1}]$
\index{notat}{DMA1effA@$\DMue(\A)$}
and refer to it as the (effective)
 \emph{$\Pmor$-premotivic $\AA^1$-derived category}
\index{word}{premotivic!Ppremotivic@$\Pmor$-premotivic!$\AA^1$-derived category}
 with coefficients in $\A$.
%
%For a scheme $S$, a $\Pmor$-scheme $X/S$, and an
%element $i$ in $I$, we denote by $\Meff S X \A\{i\}$ the image
%of $\Mab {S} X \A\{i\}$ in $\DMue(\A)$ by the localization functor
%from $\Der(\A_{S})$ to $\DMue(\A)$.
\end{df}
By definition, the category $\DMue(\A)$ satisfies the homotopy property
 \htp (see \ref{df:ppty:htp}).
 According to the general facts about localization of derived premotivic
 categories, the triangulated premotivic category $\DMue(\A)$
 is $\tau$-generated.

\begin{ex}
\label{ex:AA^1-derived_categories} We can divide our examples
 into two types: \\
\renewcommand{\Rc}{\Lambda}
\noindent 1) Assume $\Pmor=\sm$:

Consider the admissible topology $t=\nis$.
Following F.~Morel, we define the 
\emph{(effective) $\AA^1$-derived category over $S$} to be
$\DMue\left(\sh {\nis} {\sm/S}\right)$.
Indeed we get a triangulated premotivic category 
(see also the construction of \cite{ayoub2}):
\begin{equation} \label{eq:DM_tilde_effectif}
\DMuex \Rc:=\DMue\left(\sh {\nis} {\sm}\right).
\end{equation}
\index{notat}{DMA1Lambdaeffa@$\DMuex \Lambda$}
We shall also write its fibers
\begin{equation} \label{eq:DM_tilde_effectifbis}
\DMue(S,\Rc):=\DMuex \Rc(S)=\DMue\left(\sh {\nis} {\sm/S}\right)
\end{equation}
%\index{notat}{DMA1effSLambda@$\DMue(S,\Lambda)$}
for a scheme $S$. For $\Rc=\ZZ$, we shall often write simply
\begin{equation} \label{eq:DM_tilde_effectifsimple}\renewcommand{\Rc}{\ZZ}
\DMue:=\DMue\left(\sh {\nis} {\sm}\right).
\end{equation}\renewcommand{\Rc}{\Lambda}

Another interesting case is when $t=\et$;
 we get a triangulated premotivic
 category of \emph{effective \'etale premotives}:
$$
\DMue\left(\sh {\et} {\sm}\right).
$$

In each of these cases, we denote by $\Rc_S^t(X)$
 the premotive associated with a smooth $S$-scheme $X$.

\noindent 2) Assume $\Pmor=\sft$:

Consider the admissible topology $t=\h$ (resp. $t=\qfh$).
In \cite{V1}, Voevodsky has introduced the category
 of $\h$-motives (resp. $\qfh$-motives).
In our formalism,
 one defines the category of \emph{effective $\h$-motives}
\index{word}{motive!effective $\h$-motives}
(resp. \emph{effective $\h$-motives})
\index{word}{motive!effective $\qfh$-motives}
 over $S$ with coefficients in $\Rc$ as:
\begin{align*}
& \underline{\mathit{DM}}^{\mathit{eff}}_\h(S,\Rc)
 =\DMue\!\left(\sh {\h} {\sft/S}\right), \\
\text{resp. } & \underline{\mathit{DM}}^{\mathit{eff}}_\qfh(S,\Rc)
 =\DMue\!\left(\sh {\qfh} {\sft/S}\right).
\end{align*}
\index{notat}{DMteffSL@$\underline{\mathit{DM}}^{\mathit{eff}}_{\h,\Rc}$}
\index{notat}{DMteffSL@$\underline{\mathit{DM}}^{\mathit{eff}}_{\qfh,\Rc}$}
In other words, this is the $\AA^1$-derived category of
 $\h$-sheaves (resp. $\qfh$-sheaves) of $\Rc$-modules.
 \index{word}{sheaf!hsheaf@$\h$-sheaf}
\index{word}{sheaf!qfhsheaf@$\qfh$-sheaf}
Moreover, these categories for various schemes $S$ are
 the fibers of a generalized premotivic triangulated category.
What we have added to the construction of Voevodsky
 is the functors of the generalized premotivic structure.

We will denote simply by $\uRc_S^t(X)$
\index{notat}{uLambdaStX@$\uRc_S^t(X)$}
 the corresponding premotive associated with $X$ in
 $\underline{\mathit{DM}}^{\mathit{eff}}_t(S,\Rc)$.

Another interesting case is obtained when $t=\cdh$.
 We get an $\AA^1$-derived generalized premotivic category
  $\DMue\!\left(\sh {\cdh} \sft\right)$
  whose premotives are simply denoted by $\uRc^\cdh_S(X)$
  for any finite type $S$-scheme $X$.
\end{ex}

\begin{num} \label{num:AA^1-equivalences_weak&strong}
Let $C$ be a complex with coefficients in $\A_S$.
According to the general case, we say that 
$C$ is \emph{$\AA^1$-local}
\index{word}{local!A1local@$\AA^1$-local}
if for any $\Pmor$-scheme $X/S$ and any $(i,n) \in \tau \times \ZZ$, 
the map induced by the canonical projection
$$
\Hom_{\Der(\A_S)}(\mab S X\{i\}[n],C)
 \rightarrow
  \Hom_{\Der(\A_S)}(\mab S {\AA^1_X}\{i\}[n],C)
$$
is an isomorphism.
The adjunction \eqref{eq:adjunction_bousfield_loc}
defines a morphism of triangulated $\Pmor$-premotivic categories
$$
\Der(\A) \rightleftarrows \DMue(\A)
$$
such that for any scheme $S$, $\DMue(\A_S)$
is identified with the full subcategory of $\Der(\A_S)$ made of 
$\AA^1$-local complexes.

Fibrant objects for the model category structure on $\Comp(\A_S)$
appearing in Proposition \ref{t-model_category_complexes_localized}
relatively to $\W_{\AA^1}$,
simply called \emph{$\AA^1$-fibrant}
\index{word}{fibrant!A1fibrant@$\AA^1$-fibrant}
 objects, are the $t$-flasque and $\AA^1$-local
complexes.

We say a morphism
$f:C \rightarrow D$ of complexes of $\A_S$
is an \emph{$\AA^1$-equivalence}
\index{word}{equivalence!A1equivalence@$\AA^1$-equivalence}
if it becomes an isomorphism in $\DMue(\A_S)$.
Considering moreover two morphisms
$f,g:C \rightarrow D$ of complexes of $\A_S$,
we say they are \emph{$\AA^1$-homotopic}
\index{word}{homotopic, $\AA^1$-homotopic}
 if there exists a morphism of complexes 
$$
H:\mab S {\AA^1_S} \otimes_S C \rightarrow D
$$
such that $H \circ (s_{0} \otimes 1_C)=f$
and $H \circ (s_{1} \otimes 1_C)=g$,
where $s_0$ and $s_1$ are respectively induced by the zero and 
the unit section of $\AA^1_S/S$. When $f$ and $g$
are $\AA^1$-homotopic, they are equal as morphisms of
$\DMue(\A_S)$. We say the morphism $p:C \rightarrow D$
is a \emph{strong $\AA^1$-equivalence}
\index{word}{equivalence!strong $\AA^1$-equivalence}
if there exists a morphism $q:D \rightarrow C$ such that
the morphisms $p \circ q$ and $q \circ p$ are 
$\AA^1$-homotopic to the identity. A complex $C$ is
\emph{$\AA^1$-contractible}
\index{word}{contractible, A1contractible@$\AA^1$-contractible}
 if the map $C\To 0$ is a strong $\AA^1$-equivalence.

As an example, for any integer $n \in \NN$,
and any $\Pmor$-scheme $X/S$,
the map
$$
p_{*}:\mab S {\AA^n_X} \rightarrow \mab {S} X
$$
induced by the canonical projection
is a strong $\AA^1$-equivalence
with inverse the zero section $s_{0,*}:\mab {S} X \rightarrow \mab {S} {\AA^n_X}$.
\end{num}

\begin{num} \label{num:functoriality_DMue}
The category $\DMue(\A)$ is functorial in $\A$.

Let $\varphi:\A \rightleftarrows \B:\psi$ be
 an adjunction of abelian $\Pmor$-premotivic categories.
Consider two topologies $t$ and $t$' such that $t'$ is finer than $t$.
Suppose $\A$ (resp. $\B$) is compatible with $t$ (resp. $t'$).

For any scheme $S$,
consider the evident extensions
$\varphi_S:\Comp(\A_S) \rightleftarrows \Comp(\B_S):\psi_S$ of
the above adjoint functors to complexes.
We easily check that
 the functor $\psi_S$ preserves $\AA^1$-local complexes.
Thus, applying \ref{num:derived_adjunction},
 the pair $(\varphi_S,\psi_S)$ is a Quillen adjunction 
for the respective $\AA^1$-localized model structure on 
$\Comp(\A_S)$ and $\Comp(\B_S)$; see \cite[3.11]{CD1}.
Considering the derived functors, 
it is now easy to check we have obtained an adjunction
$$
\derL \varphi:\DMue(\A)
 \rightleftarrows\DMue(\B):\derR\psi
$$
of triangulated $\Pmor$-premotivic categories.
\end{num}

\begin{ex} \label{ex:functoriality_basic_A^1_derived}
\renewcommand{\Rc}{\Lambda}
Consider the notations of \ref{ex:AA^1-derived_categories}.
In the case where $\Pmor=\sm$,
 we get from the adjunction of \eqref{ex:associated_sheaf=adjunction}
 the following adjunction of triangulated premotivic categories
$$
a^*_{\et}:\DMuex \Rc
 \rightleftarrows  \DMue\left(\sh {\et} {\sm}\right):\derR a_{\et,*}.
$$
%% In the case $\Pmor=\sft$, the adjunction \eqref{ex:associated_sheaf=adjunction}
%%  defines the following morphism of generalized triangulated premotivic
%%   categories:
%% $$
%% \DMue\left(\sh {\cdh} {\sft}\right) \xrightarrow{\underline a_{\qfh}}
%%  DM^{eff}_{\qfh,\Rc} \xrightarrow{a_{\h}}
%%  DM^{eff}_{\h,\Rc}.
%% $$
\end{ex}

\begin{ex} \label{ex:universality_(Htp),(BG)}
Let $\T$ be a derived $\Pmor$-premotivic category
 as in \ref{descent_derived_premotives}.
If $\T$ satisfies the property \htp, then the canonical morphism
 \eqref{eq:universal_derived_P-premotivic}
induces a morphism
$$
\DMue(\psh{\Pmor}) \rightleftarrows \T.
$$
If moreover $\T$ satisfies $t$-descent for an admissible topology $t$,
 we further obtain as in \ref{prop:descent&derived_P-premotivic} a morphism
$$
\DMue(\sh t {\Pmor}) \rightleftarrows \T.
$$
Particularly interesting cases are given by $\DMue$
(resp. $\DMue\!\left(\sh {\cdh} {\sft}\right)$)
which is the universal derived premotivic category
(resp. generalized premotivic category),
\emph{i.e.} initial premotivic category
satisfying Nisnevich descent (resp. $\cdh$ descent) and the homotopy property.
\end{ex}

\begin{num} \label{num:A^1-derived_chg_coef}
\renewcommand{\Rc}{\Lambda}
As in Example \ref{ex:derived_chg_coef},
 let $t$ be an admissible topology and 
 $\varphi:\Lambda \rightarrow \Lambda'$ be an extension of rings.
 Then, 
  from the $\Pmor$-premotivic adjunction \eqref{eq:abelian_chg_coef}
  and according to Paragraph \ref{num:functoriality_DMue},
  we get an adjunction of triangulated $\Pmor$-premotivic categories:
$$
\derL \varphi^*:\DMue\big(\sh t {\Pmor}\big)
 \leftrightarrows
 \renewcommand{\Rc}{\Lambda'} 
 \DMue\big(\sh t {\Pmor}\big):\derR \varphi_*.
$$
Consider also complexes $C$ and $D$ of $t$-sheaves
 of $\Rc$-modules over $\Pmor_S$.
 Then there exists a canonical morphism of $\Lambda'$-modules:
\begin{equation} \label{eq:A^1-derived_chg_coef}
\Hom_{\DMue(\sh t {\Pmor_S})}\big(C,D\big) \otimes_\Lambda \Lambda'
 \longrightarrow
 \renewcommand{\Rc}{\Lambda'} 
\Hom_{\DMue(\sh t {\Pmor_S})}\big(\derL \varphi^*(C),\derL \varphi^*(D)\big)
\end{equation}
There are two notable cases where this map is an isomorphism:
\end{num}
\begin{prop} \label{prop:A^1-derived_chg_coef}
Consider the above assumptions.
Then the map \eqref{eq:A^1-derived_chg_coef} is an isomorphism
 in the two following cases:
\begin{enumerate}
\item If $\Lambda'$ is a free $\Lambda$-module and $C$ is compact;
\item If $\Lambda'$ is a free $\Lambda$-module of finite rank.
\end{enumerate}
\end{prop}
\begin{proof} \renewcommand{\Rc}{\Lambda} 
Note that in any case, the functor $\varphi_*$ admits a right adjoint
 $\varphi^!$.\footnote{It is defined by the formula: 
$$\varphi^!(F)=\uHom_{\Rc}(\Rc',F)$$
equipped with its canonical structure of sheaf of $\Rc'$-modules.} \\
We can assume that $\Lambda'=I.\Lambda$ for a set $I$.
In this case, we get for any sheaf $F$ of $\Rc$-modules:
$$
\varphi_*\varphi^*(F)=F \otimes_\Lambda \Lambda'
 =I.F.
$$
Moreover, for any $\Pmor$-scheme $X/S$, we get:
$$
\renewcommand{\Rc}{\Lambda^{\prime t}}
\varphi_*(\rep S X)
=\rep S X=\renewcommand{\Rc}{\Lambda}
I.\repx t S X.
$$
In particular, the functor
$
\varphi_*: \renewcommand{\Rc}{\Lambda'} 
\Comp(\sh t {\Pmor_S}) \rightarrow
 \renewcommand{\Rc}{\Lambda} 
 \Comp(\sh t {\Pmor_S})
$
 satisfies descent in the sense of \cite[2.4]{CD1} and preserves
 the family $\W_{\AA^1}$. Thus it is a left Quillen functor with respect
 to the $\AA^1$-local model structures. In particular,
 because it is also a right Quillen functor, we get:
 $\derR \varphi_*=\varphi_*=\derL \varphi_*$.
 In particular, we get in $\DMue(\sh t {\Pmor_S})$:
$$
\derR \varphi_* \derL \varphi^*(D)=
\derL \varphi_* \derL \varphi^*(D)=
\derL(\varphi_* \varphi^*)(D)=I.D.
$$
Thus the Proposition follows as the functor $\Hom(C,-)$
 commutes with direct sums if $C$ is compact and with
 finite direct sums in any case.
\end{proof}

We remark the following useful property.
\begin{prop} \label{prop:A^1-derivation_exact_right_adj}
Consider a morphism
$$
\varphi^*:\A \rightleftarrows \B:\varphi_*
$$
of abelian $\Pmor$-premotivic categories such
that $\A$ (resp. $\B$) is compatible with an admissible topology $t$
 (resp. $t'$). Assume $t'$ is finer than $t$.

Let $S$ be a base scheme.
Assume that $\varphi_*:\A_S \rightarrow \B_S$ commutes
 with colimits\footnote{This amounts to ask
 that $\varphi_{*}$ is exact and commutes with direct sums.}.
Then $\varphi_*:\Comp(\A_S) \rightarrow \Comp(\B_S)$
 respects $\AA^1$-equivalences.
\end{prop}
In other words, the right derived functor
 $\derR \varphi_*:\DMue(\B_S) \rightarrow \DMue(\A_S)$
 satisfies the relation $\derR \varphi_*=\varphi_*$.
\begin{proof}
In this proof, we write $\varphi_*$ for $\varphi_{*,S}$.
We first prove that $\varphi_{*}$ preserves strong $\AA^1$-equivalences
 (see \ref{num:AA^1-equivalences_weak&strong}).

Consider two maps $u,v:K \rightarrow L$ in $\Comp(\B_S)$.
To give an $\AA^1$-homotopy
 $H:\Mab S {\AA^1_S} \B \otimes_S K \rightarrow L$
 between $u$ and $v$
 is equivalent by adjunction to give a map
 $H':K \rightarrow \uHom_{\B_S}(\Mab S {\AA^1_S} \B,L)$
 which fits into the following commutative diagram:
%$$
%\xymatrix@R=20pt@C=35pt{
%&&& L \\
%K\ar|/0pt/{H'}[rr]\ar@/^8pt/^u[rrru]\ar@/_8pt/_v[rrrd]
% && \uHom_{\B_S}(\Mab S {\AA^1_S} \B,L)\ar_/12pt/{s_0^*}[ru]\ar^/12pt/{s_1^*}[rd] \\
%&&& L
%}
%$$
$$
\xymatrix@R=30pt@C=35pt{
& K\ar^{H'}[d]\ar_u[ld]\ar^v[rd] & \\
L
 & \uHom_{\B_S}(\Mab S {\AA^1_S} \B,L)\ar^-{s_0^*}[l]\ar_-{s_1^*}[r]
 & L
}
$$
where $s_0$ and $s_1$ are the respective zero and unit section of $\AA^1_S/S$.

Because $\Mab S {\AA^1_S} \B=\varphi^*_S(\Mab S {\AA^1_S} \A)$,
 we get a canonical isomorphism (see paragraph \ref{num:mon_exchange4})
$$
\varphi_*(\uHom_{\B_S}(\Mab S {\AA^1_S} \B,L))
 \simeq \uHom_{\B_S}(\Mab S {\AA^1_S} \A,\varphi_*(L)).
$$
Thus, applying $\varphi_*$ to the previous commutative diagram and using
 this identification, we obtain that $\varphi_*(u)$ is $\AA^1$-homotopic to
 $\varphi_*(v)$.

As a consequence, for any $\Pmor$-scheme $X$ over $S$,
 and any $\B$-twist $i$, the map
$$
\varphi_*(\Mab S {\AA^1_X} \B\{i\}) \rightarrow \varphi_*(\Mab S X \B\{i\})
$$
induced by the canonical projection is a strong $\AA^1$-equivalence,
 thus an $\AA^1$-equivalence.
 
The functor $\varphi_*:\B_S \rightarrow \A_S$ commutes with colimits.
Thus it admits a right adjoint that we will denote by $\varphi^!$.
Consider the injective model structure on $\Comp(A_S)$ and $\Comp(\B_S)$
 (see \cite[2.1]{CD1}). Because $\varphi_*$ is exact, it is a left Quillen
 functor for these model structures. Thus,
  the right derived functor $\derR \varphi^!$ is well-defined. 
From the result we just get, we see that $\derR \varphi^!$
preserves $\AA^1$-local objects, and this readily implies $\derL \varphi_*=\varphi_*$
preserves $\AA^1$-equivalences.
\end{proof}

\begin{paragr}
\label{par:link_homotopy_homology}
%The category $\DMte(S)$ 
%is introduced by Morel in \cite[4.3.1]{dmtilde}
%in the case where $S$ is the spectrum of a field,
%though the author does not consider the same model structure
%as the one given by Proposition \ref{t-model_category_complexes}.
To relate the category $\DMue(S)$ with the homotopy category
\index{word}{homotopy category}
 of schemes of Morel and Voevodsky~\cite{MV}, we have to consider
the category of simplicial Nisnevich sheaves of sets
denoted by $\Delta^{op} \operatorname{Sh}(\sm/S)$.
Considering the free abelian sheaf functor, 
we obtain an adjunction of categories
$$
\Delta^{op} \operatorname{Sh}(\sm/S)
 \rightleftarrows \Comp(\sh {} {\sm/S}).
$$
If we consider Blander's projective $\AA^1$-model structure~\cite{BL} on
the category $\Delta^{op} \operatorname{Sh}(\sm/S)$, 
we can easily see that this is a Quillen pair,
so that we obtain a $\Pmor$-premotivic adjunction of simple $\Pmor$-premotivic categories
%According to \cite[part 2, p. 55-56]{MV}, 
%the normalized chain complex functor induces an adjunction
%$$
%\Delta^{op} \sh {} {\sm/S}
% \rightleftarrows \Comp(\sh {} {\sm/S}).
%$$
%By composition, we thus obtain for any scheme $S$ an adjunction
%$$
%N:\Delta^{op} \sh{}{\sm/S;\ens}
% \rightleftarrows \Comp(\sh {} {\sm/S}):K.
%$$
%There is an obvious structure of a ($\emptyset$-twisted) $\Pmor$-premotivic 
%category on the fibered category
% $\Delta^{op} \sh{}{\sm;\ens}$ and it is obvious
%the adjunction $(N,K)$ is an adjunction of $\Pmor$-premotivic categories.
%By definition of the normalized chain complex,
%the functor $N$ (resp. $K$)
%sends weak equivalences 
%(resp. quasi-isomorphisms) to quasi-isomorphisms 
%(resp. weak equivalences). In particular it induces
%a pair of adjoint functors
%$$
%N:Ho(\Delta^{op} \sh{}{\sm/S;\ens})
% \rightleftarrows \Der(\sh {} {\sm/S}):K.
%$$
%Moreover, this functor respects the $\AA^1$-localized 
%categories and induces a $\Pmor$-premotivic adjunction
$$
N:\H \rightleftarrows \DMue:K.
$$
Note that the functor $N$ sends cofiber sequences
 in $\H(S)$ to distinguished triangles in $\DMue(S)$.
\end{paragr}

%\begin{rem}
%Consider an abelian $\Pmor$-premotivic category $(\A,\ZZ^\A)$
%which satisfies Nisnevich descent.
%
%From the previous constructions, we thus obtain an adjunction
%of $\Pmor$-premotivic categories
%$$
%\iota N:\H \rightleftarrows \DMue(\A).
%$$
%This allows to deduce from the results of \cite{MV}
%nice properties in the abstract category $\DMue(\A_S)$ 
%for a base scheme $S$. 
%For example, this is the case of the purity theorem 2.23
%of \emph{loc. cit.}
%\end{rem}

\subsubsection{Explicit $\AA^1$-resolution}

\begin{num}
\label{num:Suslin_complex}
Consider an abelian $\Pmor$-premotivic category $\A$ compatible
 with an admissible topology $t$. \\
Consider the canonically split exact sequence
$$0\To \un_{S}\xrightarrow{s_{0}}\mab S {\AA^1_S} \To U\To 0$$
where the map $s_{0}:\un_{S}\To\mab S {\AA^1_S}$ is induced by the
zero section of $\AA^1$. The section corresponding to $1$
in $\AA^1$ defines another map
$$s_{1}:\un_{S}\To\mab S {\AA^1_S}$$
which does not factor through $s_{0}$, so that we get
canonically a non-trivial map
$u:\un_{S}\To U$. This defines for any complex $C$ of $\A_{S}$
a map, called the \emph{evaluation at $1$},
$$\sHom(U,C)=\un_{S}\otimes_{S}\sHom(U,C)\xrightarrow{u\otimes 1}
U\otimes_{S}\sHom(U,C)\xrightarrow{\mathit{ev}} C.$$
We define the complex $R^{(1)}_{\AA^1}(C)$ to be
$$R^{(1)}_{\AA^1}(C)=\mathrm{Cone}\big(\sHom(U,C)\To C\big).$$
We have by construction a map
$$r_{C}:C\To R^{(1)}_{\AA^1}(C).$$
This defines a morphism of functors from the identity functor
to $R^{(1)}_{\AA^1}$. For an integer $n\geq 1$, we define
by induction a complex
$$R^{(n+1)}_{\AA^1}(C)=R^{(1)}_{\AA^1}(R^{(n)}_{\AA^1}(C)),$$
and a map
$$r_{R^{(n)}_{\AA^1}(C)}:R^{(n)}_{\AA^1}(C)\To R^{(n+1)}_{\AA^1}.$$
We finaly define a complex $\sing(C)$ by the formula
$$\sing(C)=\ilim_{n}R^{(n)}_{\AA^1}(C).$$
\index{notat}{RA1@$\sing$}
We have a functorial map
$$C\To R_{\AA^1}(C).$$
\end{num}

\begin{lm}\label{lemmeA1eqsHom}
With the above hypothesis and notations,
 the map $C\To\sing(C)$ is an $\AA^1$-equivalence.
\end{lm}
\begin{proof}
For any closed symmetric monoidal category $\C$
and any objects $A$, $B$, $C$ and $I$ in $\C$, we have
$$\begin{aligned}
\Hom(I\otimes\sHom(B,C),\sHom(A,C))&=\Hom(\sHom(B,C),\sHom(I,\sHom(A,C)))\\
&=\Hom(\sHom(B,C),\sHom(I\otimes A,C)).
\end{aligned}$$
Hence any map $I\otimes A\To B$ induces a map
$I\otimes\sHom(B,C)\To\sHom(A,C)$ for any object $C$.
If we apply this to $\C=\Comp(\A_{S})$ and $I=\mab S{\AA^1}$, we see immediately that
the functor $\sHom(-,C)$ preserves strong $\AA^1$-homotopy equivalences.
In particular, for any complex $C$, the map
$C\To\sHom(\mab S {\AA^1_X},C)$ is a strong $\AA^1$-homotopy equivalence.
This implies that $\sHom(U,C)\To 0$ is an $\AA^1$-equivalence, so that
the map $r_{C}$ is an $\AA^1$-equivalence as well. As $\AA^1$-equivalences
are stable by filtering colimits, this implies our result.
\end{proof}
%\begin{lm}\label{lemmaA1locUzero}
%A complex $C$ of $\A_{S}$ is $\AA^1$-local if and only if
%for any smooth scheme $X$ over $S$ and any $i$ in $I$, the maps
%induced by the evaluation at $1$
%$$\Hom_{\Der(\A_{S})}(\ZZ_{S}(X;\A)\{i\}\otimes U,C[n])\To
%\Hom_{\Der(\A_{S})}(\ZZ_{S}(X;\A)\{i\},C[n])$$
%is zero.
%\end{lm}
%\begin{proof}
%It is obvious that $C$ is $\AA^1$-local if and only if the groups
%$\Hom_{\Der(\A_{S})}(\ZZ_{S}(X;\A)\{i\}\otimes U,C[n])$.
%Hence this is a necessary condition. Conversely, suppose
%that the maps induced by the evaluation at $1$ are zero.
%We want to prove that for any smooth \
%\end{proof}

\begin{prop} \label{prop:singular&AA^1_localization}
Consider the above notations and hypothesis, and
assume that $t$ is bounded in $\A$.

For any $t$-flasque complex $C$ of $\A_S$, the complex 
$\sing(C)$ is $t$-flasque and $\AA^1$-local.\index{word}{local!A1local@$\AA^1$-local}
Moreover, the morphism $C \rightarrow \sing(C)$
is an $\AA^1$-equivalence.
If furthermore $C$ is $t$-flasque, so is $\sing(C)$.
\end{prop}

\begin{proof}
The last assertion is a particular case of Lemma \ref{lemmeA1eqsHom}.
The functor $R^{(1)}_{\AA^1}$ preserves $t$-flasque complexes.
By virtue of \ref{boundedcompactgen}, the functor $\sing$
has the same gentle property. It thus remains to prove that
the functor $\sing$ sends $t$-flasque complexes on $\AA^1$-local ones.
We shall use that the derived category $\Der(\A_S)$ is
compactly generated; see \ref{boundedcompactgen}.

Let $C$ be a $t$-flasque complex of $\A_{S}$.
To prove $\sing(C)$ is $\AA^1$-local, we are reduced to prove that
the map
$$\sing(C)\To\sHom(\mab S {\AA^1_X},\sing(C))$$
is a quasi-isomorphism, or, equivalently, that the complex
$\sHom(U,\sing(C))$ is acyclic.
As $U$ is a direct factor of $\Mab S {\AA^1_X} \A$,
for any $\Pmor$-scheme $X$ over $S$ and any $i$ in $I$, the
object $\ZZ_{S}(X;\A)\{i\}\otimes_{S} U$ is compact.
This implies that the canonical map
$$\ilim_{n}\sHom(U,R^{(n)}_{\AA^1}(C))\To\sHom(U,\sing(C))$$
is an isomorphism of complexes.
As filtering colimits preserve quasi-isomorphisms, the complex
$\sHom(U,\sing(C))$ (resp. $\sing(C)$) can be considered as the homotopy colimit of the complexes
$\sHom(U,R^{(n)}_{\AA^1}(C))$ (resp. $R^{(n)}_{\AA^1}(C)$).
In particular, for any compact object $K$ of $\Der(\A_{S})$, the canonical morphisms
$$\ilim_{n}\Hom(K,\sHom(U,R^{(n)}_{\AA^1}(C)))
\To\Hom(K,\sHom(U,\sing(C)))$$
$$\ilim_{n}\Hom(K,R^{(n)}_{\AA^1}(C))
\To\Hom(K,\sing(C))$$
are bijective.

By construction, we have distinguished triangles
$$\sHom(U,R^{(n)}_{\AA^1}(C))\To R^{(n)}_{\AA^1}(C)\To R^{(n+1)}_{\AA^1}(C)
\To\sHom(U,R^{(n)}_{\AA^1}(C))[1].$$
This implies that the evaluation at $1$ morphism
$$\mathit{ev}_{1}:\sHom(U,\sing(C))\To\sing(C)$$
induces the zero map
$$\Hom_{\Der(\A_{S})}(K,\sHom(U,\sing(C)))\To\Hom_{\Der(\A_{S})}(K,\sing(C))$$
for any compact object $K$ of $\Der(\A_{S})$.
Hence the induced map
$$a=\sHom(U,\mathit{ev_{1}}):\sHom(U,\sHom(U,\sing(C)))
\To\sHom(U,\sing(C))$$
has the same property: for any compact object $K$, the map
$$\Hom_{\Der(\A_{S})}(K,\sHom(U,\sHom(U,\sing(C))))
\To\Hom_{\Der(\A_{S})}(K,\sHom(U,\sing(C)))$$
is zero.

The multiplication map $\AA^1\times\AA^1\To\AA^1$
induces a map
$$\mu:U\otimes_{S}U\To U$$
such that the composition of
$$\mu^*:\sHom(U,\sing(C))\To\sHom(U\otimes_{S}U,\sing(C))
=\sHom(U,\sHom(U,\sing(C)))$$
with $a$ is the identity of $\sHom(U,\sing(C))$.
As $\Der(\A_{S})$ is compactly generated,
this implies that $\sHom(U,\sing(C))=0$ in the derived category
$\Der(\A_{S})$.
\end{proof}

\begin{rem}
Consider a $t$-flasque resolution functor 
(\emph{i.e.} a fibrant resolution for the $t$-local model structure)
$R_t:\Comp(\A_S) \rightarrow \Comp(\A_S)$, $1 \rightarrow R_t$.
As a corollary of the proposition,
the composite functor $\sing \circ R_t$ is 
a resolution functor by $t$-local and $\AA^1$-local complexes.
\end{rem}

\begin{ex} \renewcommand{\Rc}{\Lambda}
Consider an admissible topology $t$
 and the $\Pmor$-premotivic $\AA^1$-derived category
  $D=\DMue\left(\sh t {\Pmor}\right)$. 
Suppose that $t$ is bounded for abelian $t$-sheaves (for example,
this is the case for the Zariski and the Nisnevich topologies,
see \ref{ex:Zar&Nis=bounded}).

Let $C$ be a complex of abelian $t$-sheaves on $\Pmorx S$.
If $C$ is $\AA^1$-local, then
$$\Hom_{D(S)}(\Rc^{t}_S(X),C)
=\HH^n_t(X;C)$$
(this is true without any condition on $t$). \\
Consider a $t$-local resolution $C_{t}$ of $C$
in $\Comp\big(\sh t {\Pmorx S}\big)$.
Then we get the following formula:
$$
\Hom_{D(S)}\left(\Rc^{t}_S(X),C[n]\right)
=\HH^n\big(\Gamma\big(X,\sing(C_t)\big)\big).
$$
\end{ex}

\begin{cor}
\label{rightadjpreserveA1eq}
Consider a morphism of abelian $\Pmor$-premotivic categories
$$\varphi : \A\rightleftarrows\B : \psi$$
Suppose there are admissible topologies $t$ and $t'$,
with $t'$ finer than $t$, such that the following conditions are verified.
\begin{itemize}
\item[(i)] $\A$ is compatible with $t$ and $\B$ is compatible with $t'$.
\item[(ii)] $\B$ and $\Der(\B)$ are compactly $\tau$-generated.
\item[(iii)] For any scheme $S$, the functor $\psi_{S}:\B_{S}\To\A_{S}$
preserves filtering colimits.
\end{itemize}
Then, $\psi_{S}:\Comp(\B_{S})\To\Comp(\A_{S})$
preserves $\AA^1$-equivalences between $t'$-flasque objects.
If moreover $\psi_{S}$ is exact, the functor $\psi_{S}$
preserves $\AA^1$-equivalences.
\end{cor}

\begin{proof}
We already know that $\psi_{S}$ is a right Quillen functor, so that
it preserves local objects and $\AA^1$-fibrant objects.
This implies also that $\psi_{S}$ preserves
$\AA^1$-equivalences between $\AA^1$-fibrant objects
(this is Ken Brown's lemma~\cite[1.1.12]{Hovey}).
Let $D$ be a $t'$-flasque complex of $\B_{S}$.
Then $\psi_{S}(D)$ is a $t$-flasque complex of $\A_{S}$.
It follows from Proposition \ref{prop:singular&AA^1_localization}
that $\sing(D)$ is $\AA^1$-local and that $D\To\sing(D)$
is an $\AA^1$-equivalence.
Lemma \ref{lemmeA1eqsHom} implies the map
$$\psi_{S}(D)\To\sing(\psi_{S}(D))=\psi_{S}(\sing(D))$$
is a an $\AA^1$-equivalence. This implies the first
assertion.

The last assertion is a direct consequence of the first one.
\end{proof}

\begin{paragr}\label{suslinhomology2}
Consider the usual cosimplicial scheme $\Delta^\bullet$
defined by
$$\Delta^n=\spec{\ZZ[t_{0},\dots,t_{n}]/(t_{1}+\dots+t_{n}-1)}\simeq\AA^n$$
(see \cite{MV}). For any scheme $S$, we get a cosimplicial object of $\A_{S}$,
namely $\mab S {\Delta^\bullet_{S}}$.
%%We still write $\Mab S {\Delta^\bullet_{S}} \A$
%%for the corresponding complex obtained by taking the alternated
%%sums of coface operators.
Given any complex $C$ of $\A_{S}$, we define its associated
\emph{Suslin singular complex}\index{word}{singular!Suslin singular complexe}
as
\begin{equation} \label{eq:Sulin_singular_complex}
\suslin{(C)}=\mathrm{Tot}^\oplus \sHom(\mab S {\Delta^\bullet_{S}},C), 
\end{equation}
\index{notat}{Cunderline@$\suslin$}
where $\sHom(\mab S {\Delta^\bullet_{S}},C)$ is considered as a bicomplex
by the Dold-Kan correspondence.
The canonical map $\mab S {\Delta^\bullet_{S}} \To \un_{S}$
induces a map
$$C\To\suslin(C).$$
\end{paragr}

\begin{lm}\label{suslinA1eqtriv1}
For any complex $C$ of $\A_{S}$, the map
$$\suslin(C)\To\sHom(\mab S {\AA^1_S},\suslin(C))
=\suslin(\sHom(\mab S {\AA^1_S},C))$$
is a chain homotopy equivalence.
\end{lm}

\begin{proof}
The composite morphism 
$$
(s_{0}p \times Id)_*:\mab S {\AA^1 \times \Delta^\bullet_{S}}
 \rightarrow \mab S {\AA^1 \times \Delta^\bullet_S},
$$
where $s_{0}$ is the map induced by the zero section,
and $p$ is the map induced by the obvious projection of $\AA^1$
on  its base, is chain homotopic to the identity. 
Indeed, the homotopy relation is given by the formula
$$
s_n=\sum_{i=0}^n (-1)^i.(1 \otimes_{S} \psi_i):
 \mab S {\AA^1 \times \Delta^{n+1}_{S}} 
  \rightarrow \mab S {\AA^1 \times \Delta^n_{S}} 
$$
where $1$ is the identity of $\mab S {\AA^1_S} $, and
$\psi_i$  is induced by the map
$\Delta^{n+1}_{S} \rightarrow \AA^1 \times \Delta^n_{S}$ which
sends the $j$-th vertex $v_{j,n+1}$ to either $0 \times v_{j,n}$,
if $j \leq i$, or to $1 \times v_{j-1,n}$ otherwise.
This implies the lemma.
\end{proof}

\begin{lm}\label{suslinA1eqtriv12}
%%If $t$ is bounded in $\A$, then,
For any $t$-flasque complex $C$ of $\A_S$, we have a canonical
isomorphism
$$\suslin(C)\simeq \derL\varinjlim_n \derR\uHom(\mab S {\Delta^{n}_S},C)$$
in $\Der(\A_S)$.
\end{lm}

This is a variation on the Dold-Kan correspondence.
As a direct consequence, we get: 

\begin{lm}\label{suslinA1eqtriv2}
For any complex $C$ of $\A_{S}$, the map $C\To\suslin(C)$
is an $\AA^1$-equivalence.
\end{lm}

\begin{prop}
If $t$ is bounded in $\A$, then, for any $t$-flasque complex $C$ of $\A_S$,
$\suslin (C)$ is $\AA^1$-local.
\end{prop}

\begin{proof}
Using the first premotivic adjunction of example \ref{ex:universality_(Htp),(BG)}
and the fact that $\Der(\A)$ is compactly generated (\ref{boundedcompactgen}),
we can reduce the proposition to the case where $\A_S$ is the category of presheaves of abelian groups
over $\Pmor/S$, in which case this is well-known.
\end{proof}

\subsubsection{Constructible $\AA^1$-local premotives}
\label{sec:constructible_A1derived}

\begin{num} \label{num:hyp_compact_DMue}
Consider an abelian $\Pmor$-premotivic category $\A$
 compatible with an admissible topology $t$.
Assume that $t$ is bounded in $\A$ (see Definition \ref{df:hyper-bounded})
 and consider a bounded
 generating family $\N^t_S$ for $t$-hypercovers in $\A_S$.

Let $\cT_{\AA^1_S}$ be
the family of complexes of $\Comp(\A_S)$ 
of shape
$$
\mab S {\AA^1_X}\{i\} \rightarrow \mab S X\{i\}
$$
for a $\Pmor$-scheme $X$ over $S$ and a twist $i \in I$.
Then the functor \eqref{eq:derived_geometric_premotives_proto}
 obviously induces the following functor
\begin{equation} \label{eq:A^1-derived_geometric_premotives}
\Big(K^b\big(\Mab {} {\Pmorx S} \A\big)
      /\cN^t_S \cup \cT_{\AA^1_S}\Big)^\natural
 \rightarrow \DMue(\A_S),
\end{equation}
where the category on the left is the pseudo-abelian
category associated to the Verdier
quotient of $K^b\big(\Mab {} {\Pmorx S} \A\big)$ by the
thick subcategory generated by $\cN^t_S \cup \cT_{\AA^1_S}$.
Applying Thomason's localization theorem~\cite{Nee1},
we get from Proposition \ref{compact_objects_derived} the following result:
\end{num}
\begin{prop} \label{prop:compact_DMue}
Consider the previous hypothesis and notations
 and assume that $\A$ is finitely $\tau$-presented.
 \index{word}{finitely presented!finitely $\tau$-presented}
 
Then $\DMue(\A)$ is compactly $\tau$-generated.
\index{word}{generated!compactly $\tau$-generated!triangulated $\Pmor$-fibred}
Moreover,
 the functor \eqref{eq:A^1-derived_geometric_premotives} is fully faithful.
\end{prop}

Let us denote by $\DMuex{c}(\A)$ the subcategory of $\DMue(\A)$
 made of $\tau$-constructible premotives in the
 sense of Definition \ref{df:tau-geometric}.
 Taking into account Proposition \ref{constructequivcompact},
  we deduce from the above proposition the following corollary:
\begin{cor} \label{cor:DMue_compact&constructible}
Under the assumptions of \ref{prop:compact_DMue},
for any premotive $\cM$ in $\DMue(\A_S)$, the following conditions are equivalent:
\begin{enumerate}
\item[(i)] $\cM$ is compact;\index{word}{compact}
\item[(ii)] $\cM$ is $\tau$-constructible.
\index{word}{constructible!$\tau$-constructible}
\end{enumerate}
Moreover, the functor \eqref{eq:A^1-derived_geometric_premotives}
 induces an equivalence of categories:
$$
\Big(K^b\big(\Mab {} {\Pmorx S} \A\big)
      /\cN^t_S \cup \cT_{\AA^1_S}\Big)^\natural
 \rightarrow \DMuex{c}(\A_S).
$$
\end{cor}

\begin{ex} \label{ex:A^1-Nis_Et_compacity}
\renewcommand{\Rc}{\Lambda}
With the notations of \ref{ex:Nis_Et_compacity},
 we get the following equivalences of categories:
\begin{align*}
\left(\K^b\left(\Rc(\sm/S)\right)/(BG_S \cup \cT_{\AA^1_S})\right)^\natural
 & \rightarrow \DMuex{c}(S,\Rc). \\
\left(\K^b\left(\Rc(\sft/S)\right)/CDH_S \cup \cT_{\AA^1_S}\right)^\natural
 & \rightarrow \DMuex{c}\Big(\sh {\cdh} {\smash{\sft/S}}\Big).
\end{align*}
This statement is the analog of the embedding theorem \cite[chap. 5, 3.2.6]{FSV}.
\end{ex}

\begin{prop} \label{prop:continuity_gen-A^1_premotivic}
Assume $\Pmor=\sft$ is the class of finite type
 (resp. separated and of finite type) morphisms.

Let $\A$ be an abelian generalized premotivic category
 compatible with an admissible topology $t$
 and
 satisfying the property (C) (resp. (wC)) of 
 Paragraph \ref{num:continuity_derived_premotivic}.

Then the triangulated generalized premotivic category
 $\DMue(\A)$ is $\tau$-continuous (resp. weakly $\tau$-continuous)
 --- see Definition \ref{df:continuous}.
\end{prop}
\begin{proof}
 The proof relies on the following lemma:
\begin{lm}
Under the assumptions of the preceding proposition,
 for any morphism of schemes $f:T \rightarrow S$,
 the functor 
$$\derL f^*:\Der(\A_S) \rightarrow \Der(\A_T)$$
preserves $\AA^1$-local complexes.
\end{lm}
When $f$ is a morphism of finite type
 (resp. separated of finite type),
 the functor $\derL f^*$ admits $\derL f_\sharp$ as a left adjoint
 and the lemma is clear.
 In the general case, one can write $f$ as a projective limit
 of a projective system of morphisms of scheme
 $(f_\alpha:T_\alpha \rightarrow S)_{\alpha \in A}$
 such that $f_\alpha$ is affine of finite type.
 Recall from Proposition \ref{prop:continuity_derived_premotivic},
 $\Der(\A)$ is $\tau$-continuous.
 Thus, to check that
 for an $\AA^1$-local complex $C$ in $\Der(\A_S)$,
 the complex $\derL f^*(C)$ is $\AA^1$-local,
 we thus are reduced to prove that 
 $\derL f_\alpha^*(C)$ is $\AA^1$-local
 which follows from the first treated case. The lemma is proven.

Given the full embedding $\DMue(\A) \rightarrow \Der(\A)$
 whose image is made of $\AA^1$-local complexes,
 the proposition now directly follows from
 the previous lemma and the fact $\Der(\A)$ is $\tau$-continuous.
\end{proof}

\begin{ex}
Taking into account the second point of Example \ref{ex:sheaves_continuity},
 the previous proposition can be applied to the category
 $\sh t {\sft}$ where $t=\nis, \et, \cdh, \qfh, \h$.
\end{ex}

\begin{rem}
The previous proposition will be extended to the (non
 generalized) premotivic case in Corollary
 \ref{cor:continuity_stable-A^1_premotivic}.
\end{rem}

\subsection{The stable $\AA^1$-derived premotivic category} \label{sec:P^1-stable-derived}

\subsubsection{Modules}

Let $\A$ be an abelian $\Pmor$-premotivic category
 with generating set of twists $\tau$.

A \emph{cartesian commutative monoid}
\index{word}{monoid!cartesian commutative monoid}
 $R$ of $\A$ is a cartesian section
\index{word}{section!cartesian}
 of the fibred category $\A$ over $\sch$ such that
 for any scheme $S$, $R_S$ has a commutative monoid structure
 in $\A_S$ and for any morphism of schemes $f:T \rightarrow S$,
 the structural transition maps $\phi_f:f^*(R_S) \rightarrow R_T$
 are isomorphisms of monoids.

Let us fix a cartesian commutative monoid $R$ of $\A$.

Consider a base scheme $S$.
We denote by $\Mod{R_S}$ the category of modules in the monoidal
category $\A_S$ over the monoid $R_S$. 
%% The forgetful functor $\cO_{A,S}:\Mod{R_S} \rightarrow \A_S$ has 
%% a left adjoint denoted by $L_{R,S}$, 
%% which sends an object $E$ of $\A_S$ to
%%  the object $R_S \otimes_S E$ with its canonical $R_S$-module structure.
For any $\Pmor$-scheme $X/S$ and any twist $i \in \tau$,
 we put 
$$
R_S(X)\{i\}=R_S \otimes_S \mab S X \{i\}
$$
 endowed with its canonical $R_S$-module structure.
The category $\Mod{R_S}$ is a Grothendieck abelian category 
such that the forgetful functor $U_{S}:\Mod{R_S} \rightarrow \A_S$
is exact and conservative.
A family of generators for $\Mod{R_S}$ is given by 
the modules $R_S(X)\{i\}$ for a $\Pmor$-scheme $X/S$ and a twist $i\in \tau$.
As $A_S$ is commutative, $\Mod{R_S}$ has a unique symmetric monoidal
structure such that the free $R_S$-module functor is symmetric monoidal.
We denote by $\otimes_R$ this tensor product.
Note that $R_S(X) \otimes_R R_S(Y)=R_S(X \times_S Y)$.
Finally, the categories of modules $\Mod{R_S}$
form a symmetric monoidal $\Pmor$-fibred category, such that
the following proposition holds (see \ref{PfibredRMod}).

\begin{prop} \label{module&premotivic}
Let $\A$ be a $\tau$-generated abelian $\Pmor$-premotivic category
and $R$ be a cartesian commutative monoid of $\A$.

Then the category $\Mod R$ equipped with the structures
introduced above is a $\tau$-generated abelian $\Pmor$-premotivic category.

Moreover, we have an adjunction of abelian $\Pmor$-premotivic categories:
\begin{equation} \label{eq:adj_premot_ab_mod}
R\otimes(-):\A \rightleftarrows \Mod R: U\, .
\end{equation}
\end{prop}

\begin{rem} \label{rem:modules_exchange&topology}
With the hypothesis of the preceding proposition,
 for any morphism of schemes $f:T \rightarrow S$,
 the exchange transformation
$f^* U_{S} \rightarrow U_{T} f^*$
is an isomorphism by construction of $\Mod R$ (\ref{PfibredRMod}).
\end{rem}

\begin{prop} \label{prop:modules&topology}
Let $\A$ be a $\tau$-generated abelian $\Pmor$-premotivic category compatible
with an admissible topology $t$.
Consider a cartesian commutative monoid $R$ of $\A$
such that for any scheme $S$, tensoring quasi-isomorphisms
between cofibrant complexes by $R_S$ gives quasi-isomorphisms
(e.g. $R_S$ might be cofibrant (as a complex concentrated in degree zero),
or flat). Then the abelian $\Pmor$-premotivic category $\Mod R$
is compatible with $t$.
\end{prop}

\begin{proof}
In view of Proposition \ref{adjunction&topology}, we have only
 to show that $\Mod R$ satisfies cohomological $t$-descent.
Consider a $t$-hypercover $p:\cX \rightarrow X$ in $\Pmorx S$.
We prove that the map
 $p_*:R_S(\cX) \rightarrow R_S(X)$ is a quasi-isomorphism in
 $\Comp(\Mod{R_S})$.
The functor $U_{S}$ is conservative, and $U_{S}(p_*)$
 is equal to the map:
$$
R_S \otimes_S \mab S \cX \rightarrow R_S \otimes_S \mab S X.
$$
But this is a quasi-isomorphism in $\Comp(\A_S)$ by assumption on $R_S$.
\end{proof}

\begin{rem} \label{rem:comput_derived_Hom_mod}
According to Lemma \ref{abelian_adjunction&derived_Hom},
 for any simplicial $\Pmor$-scheme $\cX$ over $S$, 
 any twist $i \in \tau$ and any $R_S$-module $C$, 
 we get canonical isomorphisms:
\begin{align}
\label{eq:comput_derived_Hom_mod1}
\Hom_{\K(\Mod{R_S})}\big(R_S(\cX)\{i\},C\big)
 & \simeq
  \Hom_{\K(\A_S)}(\mab S \cX \{i\},C) \\
\label{eq:comput_derived_Hom_mod2}
\Hom_{\Der(\Mod{R_S})}(R_S(\cX)\{i\},C)
 & \simeq
  \Hom_{\Der(\A_S)}(\mab S \cX \{i\},C).
\end{align}
\end{rem}

\subsubsection{Symmetric sequences}

Let $\A$ be an abelian category.

Let $G$ be a group.
An action of $G$ on an object $A \in \A_S$ is a morphism of groups
$G \rightarrow \mathrm{Aut}_\A(A), g \mapsto \gamma^A_g$.
We say that $A$ is a $G$-object of $\A$. A $G$-equivariant morphism 
$A \xrightarrow f B$ of $G$-objects of $\A$ is a morphism $f$ in $\A$
such that $\gamma_g^B \circ f=f \circ \gamma_g^A$.

If $E$ is any object of $\A$, we put
$G \times E=\bigoplus_{g \in G} E$ 
considered as a $G$-object 
via the permutation isomorphisms of the summands.

If $H$ is a subgroup of $G$, and $E$ is an $H$-object, $G \times E$
has two actions of $H$~: the first one, say $\gamma$, is obtained via
the inclusion $H \subset G$, and the second one denoted by $\gamma'$,
is obtained using the structural action of $H$ on $E$. 
We define $G \times_H E$ as the
coequalizer of the family of morphisms
$(\gamma_\sigma-\gamma'_\sigma)_{\sigma \in H}$, and consider it 
equipped with its induced action of $G$.
\begin{df}
\label{df:symmetric_sequences}
Let $\A$ be an abelian category.

A symmetric sequence
\index{word}{sequence!symmetric sequence}
 of $\A$ is a sequence $(A_n)_{n \in \NN}$ 
such that for each $n \in \NN$, $A_n$ is a $\fS_n$-object 
of $\A$. 
A morphism of symmetric sequences of $\A$ is 
a collection of $\fS_n$-equivariant morphism 
$(f_n:A_n \rightarrow B_n)_{n \in \NN}$.

We let $\A^\fS$
\index{notat}{ASigma@$\A^\fS$}
 be the category of symmetric sequences of $\A$.
\end{df}
It is straightforward to check $\A^\fS$ is abelian.
For any integer $n \in \NN$, we define the \emph{$n$-th evaluation functor}
\index{word}{functor!evaluation}
 as follows:
$$ev_n:\A^\fS \rightarrow \A, A_* \mapsto A_n.$$
Any object $A$ of $\A$ can be considered as the trivial symmetric
sequence $(A,0,\ldots)$. The functor $i_0:A \mapsto (A,0,\ldots)$ is 
obviously left adjoint to $ev_0$ and we obtain an adjunction
\begin{equation} \label{eq:sym_seq_iota&ev_0}
i_0:\A \rightleftarrows \A^\fS:ev_0.
\end{equation}
Remark $i_0$ is also right adjoint to $ev_0$. Thus, $i_0$
preserves every limit and colimit.

For any integer $n \in \NN$ and 
any symmetric sequence $A_*$ of $\A$, we put
\begin{equation} \label{eq:twist_sym_seq}
\begin{split}
(A_*\{-n\})_m=\left\{
\begin{array}{ll}
\fS_m \times_{\fS_{m-n}} A_{m-n}
 & \text{if } m \geq n \\
0 & \text{otherwise.} 
\end{array}\right.
\end{split}
\end{equation}
This define an endofunctor on $\A^\fS$, and we have 
$A_*\{-n\}\{-m\}=A_*\{-n-m\}$ (through a canonical isomorphism).
Remark finally that for any integer $n \in \NN$,
the functor 
$$i_n:\A \rightarrow \A^\fS, A \mapsto (i_0(A))\{-n\}$$
is left adjoint to $ev_n$.

\begin{rem}
Let $\fS$ be the category of finite sets 
with bijective maps as morphisms. 
Then the category of symmetric sequences is canonically equivalent 
to the category of functors $\fS \rightarrow \A$.
This presentation is useful to define a tensor product on $\A^\fS$.
\end{rem}

\begin{df}\label{df:symmtensorproduct}
Let $\A$ be a symmetric closed monoidal abelian category.

Given two functors $A_*,B_*:\fS \rightarrow \A$,
we put:
$$
\begin{array}{rcl}
E \otimes^\fS F:\fS & \mapsto & \A \\
N & \mapsto & \bigoplus_{N=P \sqcup Q} E(P) \otimes F(Q).
\end{array}
$$
\index{notat}{tensorSigma@$\otimes^\fS$}
\end{df}
If $\unit_\A$ is the unit object of the monoidal category
$\A$, the category $\A^\fS$ is then a symmetric closed monoidal category
with unit object $i_0(\unit_\A)$.

\begin{num} Let $A$ be an object of $\A$. 
Then the $n$-th tensor power
$A^{\otimes n}$ of $A$ is endowed with a canonical action
of the group $\fS_n$ through the structural permutation
isomorphism of the symmetric structure on $\A$.
Thus the sequence $\Sym(A)=(A^{\otimes n})_{n \in \NN}$ is a 
symmetric sequence.

Moreover, the isomorphism
$A^{\otimes n} \otimes A^{\otimes m}
 \rightarrow A^{\otimes n+m}$
is $\fS_n \times \fS_m$-equivariant. 
Thus it induces a morphism
$\mu:\Sym(A) \otimes^\fS \Sym(A) \rightarrow \Sym(A)$
 of symmetric sequences.
We also consider the obvious morphism 
$\eta:i_0(\unit_\A)=i_0(A^{\otimes 0}) \rightarrow \Sym(A)$.
One can check easily that $\Sym(A)$ equipped with the multiplication 
$\mu$ and the unit $\eta$ is a commutative monoid 
in the monoidal category $\A^\fS$.
\end{num}

\begin{df} \label{df:free_symmetric_monoid}
Let $\A$ be an abelian symmetric monoidal category.
The commutative monoid $\Sym(A)$
\index{notat}{SymA@$\Sym(A)$}
 of $\A^\fS$ defined above
will be called the symmetric monoid generated by $A$.
\end{df}

\begin{rem}\label{rem:caractcommringspectra}
One can describe $\Sym(A)$ by a universal property: given a
commutative monoid $R$ in $\A^\fS$, to give a morphism
of commutative monoids $\Sym(A)\To R$ is equivalent to give
a morphism $A\To R_1$ in $\A$.
\end{rem}

\begin{paragr}
Consider an abelian $\Pmor$-premotivic category $\A$.
 
Consider a base scheme $S$.
According to the previous paragraph,
the category $\A^\fS_S$ is an abelian category, 
endowed with a symmetric tensor product $\otimes_S^\fS$. 
For any $\Pmor$-scheme $X/S$ and any integer $n \in \NN$, 
using \eqref{eq:twist_sym_seq}, we put
$$
\Mab S X {\A^\fS}\{-n\}=i_0(\Mab S X \A)\{-n\}.
$$
It is immediate that the class of symmetric sequences 
of the form $\Mab S X {\A^\fS}\{-n\}$ for a smooth $S$-scheme $X$ and an integer
 $n \geq 0$ is a generating family for the abelian
category $\A^\fS_S$
 which is therefore a Grothendieck abelian category. 
 It is clear that for any $\Pmor$-scheme $X$ and $Y$ over $S$, 
$$
\Mab S X {\A^\fS}\{-n\} \otimes^\fS_S \Mab S Y {\A^\fS}\{-n\}
=\Mab S {X\times_S Y} {\A^\fS}\{-n\}.
$$

Given a morphism (resp. $\Pmor$-morphism) 
of schemes $f:T \rightarrow S$
and a symmetric sequence
$A_*$ of $\A_S$, we put $f^*_\fS(A_*)=(f^*A_n)_{n \in \NN}$
(resp. $f_\sharp^\fS(A_*)=(f_\sharp A_n)_{n \in \NN}$).
This defines a functor $f^*_\fS:\A_S^\fS \rightarrow \A_T^\fS$
(resp. $f_\sharp^\fS:A_T^\fS \rightarrow \A_S^\fS$)
which is obviously right exact. 
Thus, the functor $f^*_\fS$ admits a right adjoint 
which we denote by $f_*^\fS$. When $f$ is in $\Pmor$, we check easily
the functor $f_\sharp^\fS$ is left adjoint to $f^*_\fS$.

From criterion \ref{prop:carac_monoidal_P-fibred}
 and Lemma \ref{lm:carac_morph_monoidal_P-fibred},
 we check easily the following proposition:
\end{paragr}
\begin{prop}
\label{prop:premotivic_cat_of_symmetric_sequences}
Consider the previous hypothesis and notations.

The association $S \mapsto \A^\fS_S$ together with the structures
 introduced above defines an $\NN\times\tau$-generated abelian $\Pmor$-premotivic
 category.

Moreover, the different adjunctions of the form \eqref{eq:sym_seq_iota&ev_0}
 over each fibers over a scheme $S$ define an adjunction of $\Pmor$-premotivic
 categories:
\begin{equation} \label{eq:premotivic_sym_seq_iota&ev_0}
i_0:\A \rightleftarrows \A^\fS:ev_0
\end{equation}
\end{prop}
Indeed, $i_0$ is trivially compatible with twists.

\begin{prop} \label{prop:symmetric_seq&topology}
Let $\A$ be an abelian $\Pmor$-premotivic category,
and $t$ be an admissible topology.
If $\A$ is compatible with $t$ then $\A^\fS$ is compatible with $t$.
\end{prop}

\begin{proof} This is based on the following lemma (see \cite[7.5, 7.6]{CD1}):
\begin{lm} \label{lm:comput_derived_Hom_sseq}
%Assume that $\A$ is compatible with $t$.
For any complex $C$ of $\A_S$, any complex $E$ of $\A^\fS_S$
 and any integer $n \geq 0$, there are canonical isomorphisms:
\begin{align}
\label{eq:comput_derived_Hom_sseq1}
\Hom_{\K(\A_S^\fS)}(i_0(C)\{-n\},E) & \simeq \Hom_{\K(\A_S)}(C,E_n) \\
\label{eq:comput_derived_Hom_sseq2}
\Hom_{\Der(\A_S^\fS)}(i_0(C)\{-n\},E) & \simeq \Hom_{\Der(\A_S)}(C,E_n)
\end{align}
\end{lm}
If $\A$ is compatible with $t$, this implies 
 that $E$ is local (resp. $t$-flasque) if
 and only if for any $n \geq 0$, $E_n$ is local (resp. $t$-flasque).
 This concludes.
\end{proof}

\subsubsection{Symmetric Tate spectra}

\begin{num}\label{num:Tate_premot}
Consider an abelian $\Pmor$-premotivic category $\A$.

For any scheme $S$, the unit point of $\GGx S$ defines 
a split monomorphism of $\A$-premotives 
$\MabNX S \A \rightarrow \mab S {\GGx S}$.
We denote by $\MabNX S \A\{1\}$ the cokernel of this monomorphism
and call it the \emph{suspended Tate $S$-premotive} 
with coefficients in $\A$.
The collection of these objects for any scheme $S$
 is a cartesian section of $\A$ denoted by $\MabNX{} \A\{1\}$.
 For any integer $n \geq 0$, we denote by $\MabNX{}\A\{n\}$ its
 $n$-the tensor power.

With the notations of \ref{df:free_symmetric_monoid},
 we define the \emph{symmetric Tate spectrum} over $S$
 as the symmetric sequence $\MabNX S \A\{*\}=Sym(\MabNX S \A\{1\})$
 in $\A^\fS_S$.
The corresponding collection defines a cartesian commutative monoid
 of the fibred category $\A^\fS$, called the
 \emph{absolute Tate spectrum}.
\index{word}{spectrum!absolute Tate spectrum}
\end{num}
 
\begin{df} \label{df:Tate_sym_spectra}
Consider an abelian $\Pmor$-premotivic category $\A$.

We denote by $\Spt(\A)$ the abelian $\Pmor$-premotivic category 
 of modules over $\MabNX {}{\A}\{*\}$ in the category $\A^\fS$.
The objects of $\Spt(\A)$ are called
the abelian (symmetric) Tate spectra.\footnote{As we will almost never
consider non symmetric spectra, we will cancel the word "symmetric" 
in our terminology.}
\index{word}{spectrum!abelian Tate spectrum}
\index{word}{spectra|see{spectrum}}
\end{df}
The category $\Spt(\A)$ is $(\NN\times \tau)$-generated.
Composing the adjunctions \eqref{eq:adj_premot_ab_mod}
 and \eqref{eq:premotivic_sym_seq_iota&ev_0},
 we get an adjunction
\begin{equation} \label{eq:adj_premot_Tate_sp}
\sus:\A \rightleftarrows \Spt(\A):\lop
\end{equation}
of abelian $\Pmor$-premotivic categories.
%
%$$\sus(A)=\ZZ_S^\A\{*\}\otimes^\fS_{S}\iota(A)$$
%for any object $A$ of $\A_{S}$.

Let us explicit the definition.
An abelian Tate spectrum $(E,\sigma)$ is the data of~:
\begin{enumerate}
\item for any $n \in \NN$, an object $E_n$ of $\A_S$
endowed with an action of $\fS_n$
\item for any $n \in \NN$, a morphism
$\sigma_n:E_n\{1\} \rightarrow E_{n+1}$
in $\A_S$
\end{enumerate}
such that the composite map
$$
E_m\{n\}
 \xrightarrow{\sigma_m\{n-1\}} E_{m+1}\{n-1\}
   \rightarrow ... \xrightarrow{\sigma_{m+n-1}} E_{m+n}
$$
is $\fS_n \times \fS_m$-equivariant with respect to 
the canonical action of $\fS_n$ on $\MabNX {S} {\A} \{n\}$ 
and the structural action of $\fS_m$ on $E_m$.
By definition, $\lop(E)=E_0$. 
Thus, the functor $\lop$ is exact.

Given an object $A$ of $\A_S$, the abelian Tate spectrum
$\sus A$ is defined such that $(\sus A)_n=A\{n\}$ with the
action of $\fS_n$ given by its action on $\MabNX {S} {\A} \{n\}$ 
by permutations of the factors.

Be careful we consider the category $\Spt(\A_S)$ as $\NN$-twisted 
by negative twists. For any abelian Tate spectrum $E_*$,
$(E_*\{-n\})_m=\fS_m \times_{\fS_{m-n}} E_{m-n}$ for $n \geq m$.

\begin{num} \label{num:functoriality_abelian_spectra}
Consider a morphism
$$
\varphi:\A \rightarrow \B
$$
of abelian $\Pmor$-premotivic categories.
Then as $\varphi(\un^\A\{1\})=\un^\B\{1\}$,
 $\varphi$ can be extended to abelian Tate spectra
 in such a way that the following diagram commutes:
$$
\xymatrix@C=40pt{
\A\ar^\varphi[r]\ar_{\sus_\A}[d] & \B\ar^{\sus_\B}[d] \\
\Spt(\A)\ar^-{\Spt(\varphi)}[r] & \Spt(\B).
}
$$
(Of course the obvious diagram for the corresponding right adjoints also commutes.)
\end{num}

\begin{df}\label{df:complexTatespectra}
For any scheme $S$,
a complex of abelian Tate spectra over $S$
will be called simply a \emph{Tate spectrum}
\index{word}{spectrum!Tate spectrum}
 over $S$.
\end{df}
A Tate spectrum $E$ is a bigraded object. In the notation $E_n^m$,
the index $m$ corresponds to the (cochain) complex structure
and the index $n$ to the symmetric sequence structure.

From propositions \ref{prop:modules&topology}
 and \ref{prop:symmetric_seq&topology},
 we get the following:
\begin{prop}
Let $\A$ be an abelian $\Pmor$-premotivic category
compatible with an admissible topology $t$.
Then $\Spt(\A)$ is compatible with $t$.
\end{prop}
Note also that remark \ref{rem:comput_derived_Hom_mod}
 and Lemma \ref{lm:comput_derived_Hom_sseq} implies that 
 for any simplicial $\Pmor$-scheme $\cX$ over $S$, any integer $n \in \NN$,
 and any Tate spectrum $E$, we have canonical isomorphisms:
\begin{align}
\label{eq:comput_derived_Hom_Tate_spectrum1}
\Hom_{\K(\Spt(\A_S))}(\sus \Mab S \cX \A\{-n\},E)
 & \simeq \Hom_{\K(\A_S)}(\sus \Mab S \cX \A,E_n) \\
\label{eq:comput_derived_Hom_Tate_spectrum2}
\Hom_{\Der(\Spt(\A_S))}(\sus \Mab S \cX \A\{-n\},E)
 & \simeq \Hom_{\Der(\A_S)}(\sus \Mab S \cX \A,E_n)
\end{align}
According to the proposition,
 the category $\Comp(\Spt(\A_S))$ of Tate spectra over $S$
 has a $t$-descent model structure.
 The previous isomorphisms allow to describe this structure
 as follows:
\begin{enumerate}
\item For any simplicial $\Pmor$-scheme $\cX$ over $S$,
 and any integer $n \geq 0$,
 the Tate spectrum $\sus \Mab S \cX \A\{-n\}$ is cofibrant.
\item A Tate spectrum $E$ over $S$ is fibrant if and only if
 for any integer $n \geq 0$, the complex $E_n$ over $\A_S$ is local 
 (\emph{i.e.} $t$-flasque).
\item Let $f:E \rightarrow F$ be a morphism of Tate spectra over $S$.
Then $f$ is a fibration (resp. quasi-isomorphism)
 if and only if for any integer $n \geq 0$,
 the morphism $f_n:E_n \rightarrow F_n$ of complexes over $\A_S$
 is a fibration (resp. quasi-isomorphism).
\end{enumerate}
Note that properties (2) and (3) follows from \eqref{eq:comput_derived_Hom_mod1}
 and \eqref{eq:comput_derived_Hom_sseq1}.

\begin{num} \label{num:level_A^1-local_spectra}
We can also introduce the $\AA^1$-localization of this model structure.
The corresponding homotopy category is the $\AA^1$-derived $\Pmor$-premotivic category
 $\DMue(\Spt(\A))$ introduced in \ref{df:effective_triangulated_premotives}.
The isomorphism \eqref{eq:comput_derived_Hom_Tate_spectrum2} gives
 the following assertion:
From the above, a Tate spectrum $E$ is $\AA^1$-local if and only
 if for any integer $n \geq 0$, $E_{n}$ is $\AA^1$-local.
\begin{enumerate}
\item A Tate spectrum $E$ over $S$ is $\AA^1$-local if and only if
 for any integer $n \geq 0$, the complex $E_n$ over $\A_S$ is $\AA^1$-local.
\item Let $f:E \rightarrow F$ be a morphism of Tate spectra over $S$.
Then $f$ is a $\AA^1$-local fibration (resp. weak $\AA^1$-equivalence)
 if and only if for any integer $n \geq 0$,
 the morphism $f_n:E_n \rightarrow F_n$ of complexes over $\A_S$
 is a $\AA^1$-local fibration (resp. weak $\AA^1$-equivalence).
\end{enumerate}
As a consequence, the isomorphism \eqref{eq:comput_derived_Hom_Tate_spectrum2}
 induces an isomorphism
\begin{align}
\label{eq:comput_A^1-derived_Hom_Tate_spectrum2}
\Hom_{\DMue(\Spt(\A_S))}(\sus \Mab S \cX \A\{-n\},E)
 & \simeq \Hom_{\DMue(\A_S)}(\sus \Mab S \cX \A,E_n).
\end{align}
Similarly, the adjunction \eqref{eq:adj_premot_Tate_sp} induces an adjunction
 of triangulated $\Pmor$-premotivic categories
\begin{equation} \label{eq:adjunction_derived_spectra}
\derL \sus:\DMue(\A)
 \rightleftarrows
  \DMue(\Spt(\A)):\derR \lop.
\end{equation}
\end{num}

\subsubsection{Symmetric Tate $\Omega$-spectra}

\begin{paragr}
The final step is to localize further the category $\DMue(\Spt(\A))$.
The aim is to relate 
the positive twists on $\DMue(\A)$ 
obtained by tensoring with $\MabNX {S} {\A} \{1\}$
and the negative twists on $\DMue(\Spt(\A))$
induced by the consideration of symmetric sequences.

Let $X$ be a $\Pmor$-scheme over $S$.
From the definition of $\sus$,
there is a canonical morphism of abelian Tate spectra:
$$
\big\lbrack\sus\big(\MabNX {S} {\A} \{1\}\big)\big\rbrack\{-1\}
 \rightarrow \sus \MabNX {S} {\A}.
$$
Tensoring this map by $\sus\Mab S X {\A} \{-n\}$
 for any $\Pmor$-scheme $X$ over $S$
 and any integer $n \in \NN$,
 we obtain a family of morphisms of Tate spectra
 concentrated in cohomological degree $0$:
$$
\big\lbrack\sus\big(\Mab S X {\A} \{1\}\big)\big\rbrack\{-n-1\}
 \rightarrow \sus \Mab S X {\A}\{-n\}.
$$
We denote by $\W_\Omega$ this family
 and put $\W_{\Omega,\AA^1}=\W_{\Omega}\cup\W_{\AA^1}$.
Obviously, $\W_{\Omega,\AA^1}$ is stable by the operations
$f^*$ and $f_\sharp$.
\end{paragr}
\begin{df}
\label{df:triangulated_premotives}
Let $\A$ be an abelian $\Pmor$-premotivic category compatible
with an admissible topology $t$.
With the notations introduced above, we define the
 \emph{stable $\AA^1$-derived $\Pmor$-premotivic category
  with coefficients in $\A$}
\index{word}{premotivic!Ppremotivic@$\Pmor$-premotivic!stable $\AA^1$-derived category} 
  as the derived $\Pmor$-premotivic category 
$$
\DMu(\A):=\Der(\Spt(\A))[\W_{\Omega,\AA^1}^{-1}]
$$ 
\index{notat}{DA1A@$\DMu(\A)$}
defined in Corollary \ref{cor:derived_premotivic_localization_exists}.
\end{df}

\begin{num} \label{num:basic_DMu}
According to this definition, we get the following
 identification:
$$\DMu(\A)=\DMue(\Spt(\A))[\W_\Omega^{-1}].$$
Using the left Bousfield localization of the $\AA^1$-local
model structure on $\Comp(\Spt(\A))$, we thus obtain
a canonical adjunction of triangulated $\Pmor$-fibred
premotivic categories
$$
\DMue(\Spt(\A)) \rightleftarrows \DMue(\Spt(\A))[\W_\Omega^{-1}]
$$
which allows us to describe $\DMu(\A_S)$ as the full subcategory of $\DMue(\Spt(\A_S))$
made of Tate spectra which are $\W_\Omega$-local in $\DMue(\Spt(\A_S))$.
Recall a Tate spectrum $E$ is a sequence of complexes $(E_n)_{n \in \NN}$
 over $\A_S$ together with suspension maps in $\Comp(\A_S)$
$$
\sigma_n:\MabNX S \A\{1\} \otimes E_n \rightarrow E_{n+1}.
$$
From this, we deduce a canonical morphism
 $\MabNX S \A\{1\} \otimes^\derL E_n \rightarrow E_{n+1}$ in $\DMue(\A)$
 whose adjoint morphism we denote by
\begin{equation} \label{eq:A^1-local_adjoint_of_sus}
u_n:E_n \rightarrow \derR \uHom_{\DMue(\A_S)}(\MabNX S \A\{1\},E_{n+1})
\end{equation}
According to \eqref{eq:comput_A^1-derived_Hom_Tate_spectrum2}, 
 the condition that $E$ is $\W_\Omega$-local in $\DMue(\Spt(\A))$ is
 equivalent to ask that for any integer $n \geq 0$, the map \eqref{eq:A^1-local_adjoint_of_sus}
 is an isomorphism in $\DMue(\Spt(\A))$.

Considering the adjunction \eqref{eq:adjunction_derived_spectra},
 we obtain finally an adjunction of triangulated $\Pmor$-fibred
 categories:
\begin{equation} \label{eq:adj_DMue-DMu}
\sus:\DMue(\A) \rightleftarrows \DMue(\Spt(\A))
 \rightleftarrows \DMu(\A):\lop.
\end{equation}
Note that tautologically, the Tate spectrum $\sus(\MabNX S \A\{1\})$
 has a tensor inverse given by the spectrum
  $(\sus\MabNX S \A)\{-1\}$ in $\DMu(\A_S)$.
Thus, we have obtained from the abelian premotivic category $\A$
 a triangulated premotivic category $\DMu(\A_S)$ which satisfies the
properties:
\begin{itemize}
\item the homotopy property \htp;
\item the stability property \stab;
\item the $t$-descent property.
\end{itemize}
As we will see in the followings, the construction satisfies
 a universality property that the reader can already guess.
\end{num}

\begin{df}
Consider the assumptions of definition \ref{df:triangulated_premotives}.

For any scheme $S$, we say that a Tate spectrum $E$ over $S$
 is a \emph{Tate $\Omega$-spectrum}
\index{word}{spectrum!Tateomega@Tate $\Omega$-spectrum}
 if the following conditions are fulfilled:
\begin{enumerate}
\item[(a)] For any integer $n\geq 0$, $E_n$ is $t$-flasque and $\AA^1$-local.
\item[(b)] For any integer $n \geq 0$,
the adjoint of the structural suspension map
$$
E_n \rightarrow \uHom_{\Comp(\A_S)}(\MabNX S \A\{1\},E_{n+1})
$$
is a quasi-isomorphism.
\end{enumerate}
\end{df}
In particular, a Tate $\Omega$-spectrum is $\W_\Omega$-local in
 $\DMue(\Spt(\A_S))$. In fact, it is also $\W_{\Omega,\AA^1}$-local
 in the category $\Der(\Spt(\A_S))$ so that the category $\DMu(\A)$
 is also equivalent to the full subcategory of $\Der(\Spt(\A_S))$
 spanned by Tate $\Omega$-spectra.

Fibrant objects of the $\W_{\Omega,\AA^1}$-local model category 
 on $\Comp(\Spt(\A))$ obtained in definition \ref{df:triangulated_premotives}
 are exactly the Tate $\Omega$-spectra. 

%Essentially by definition, we have: 
% 
%\begin{prop} \label{prop:fomal_ppty_A1-local_spectra}
%Let $\A$ be as definition \ref{df:triangulated_premotives}.
%Then $\DMu(\A)$ is a $\Pmor$-premotivic triangulated category which satisfies
%$t$-descent, the homotopy property, and the stability property.
%\end{prop} 
 
\begin{prop} \label{prop:universality_stabilisation}
Consider the above notations. Let $S$ be a base scheme.
\begin{enumerate}
\item If the endofunctor 
$$
\DMue(\A_S) \rightarrow \DMue(\A_S),
C \mapsto \derR \uHom_{\DMue(\A_S)}(\MabNX S \A\{1\},C)
$$
is conservative, then the functor $\lop_S$ is conservative.
\item If the Tate twist
\index{word}{twist!Tate}
 $E \mapsto E(1)$ is fully faithful in $\DMue(\A_S)$,
 then $\sus_S$ is fully faithful.
\item If the Tate twist $E \mapsto E(1)$ induces
 an auto-equivalence of $\DMue(\A_S)$,
 then $(\sus_S,\lop_S)$ are adjoint equivalences of categories.
\end{enumerate}
\end{prop}

\begin{rem} Similar statements can be obtained for the derived categories
 rather than the $\AA^1$-derived categories. We left their formulation
 to the reader.
 \end{rem}
 
\begin{proof}
Consider point (1). We have to prove that for any $\W_\Omega$-local
 Tate spectrum $E$ in $\DMue(\Spt(\A_S))$, if $\derR \lop(E)=0$,
 then $E=0$. But $\derR \lop(E)=\lop (E)=E_0$
  (see \ref{num:level_A^1-local_spectra}).
 Because for any integer $n \geq 0$,
 the map \eqref{eq:A^1-local_adjoint_of_sus} is an $\AA^1$-equivalence,
 we deduce that for any integer $n \in \ZZ$,
 the complex $E_n$ is (weakly) $\AA^1$-acyclic.
 According to \eqref{eq:comput_A^1-derived_Hom_Tate_spectrum2},
  this implies $E=0$ --- because $\DMu(\A_S)$ is $\NN$-generated.

Consider point (2). We want to prove that for any complex $C$
 over $\A_S$,
 the counit map $C \rightarrow \derR \lop \derL \sus(C)$ is an isomorphism.
It is enough to treat the case where $C$ is cofibrant.

Considering the left adjoint $\derL \sus$ 
 of \eqref{eq:adjunction_derived_spectra},
 we first prove that $\derL \sus(C)$ is $\W_\Omega$-local.
Because $C$ is cofibrant, this Tate spectrum is equal in degree $n$
 to the complex $C\{n\}$ (with its natural action of $\fS_n$).
 Moreover, the suspension map is given by the isomorphism
 (in the monoidal category $\Comp(\A_S)$)
$$
\sigma_n:\MabNX S \A\{1\} \otimes_S C\{n\} \rightarrow C\{n+1\}.
$$
In particular, the corresponding map in $\DMue(\A_S)$
$$
\sigma'_n:\MabNX S \A\{1\} \otimes^L_S C\{n\} \rightarrow C\{n+1\}.
$$
is canonically isomorphic to
$$
\MabNX S \A\{1\} \otimes^L_S C\{n\}
 \xrightarrow{1 \otimes 1} \MabNX S \A\{1\} \otimes^L_S C\{n\}.
$$
Thus, because the Tate twist is fully faithful in $\DMue(\A_S)$,
 the adjoint map to $\sigma'_n$ is an $\AA^1$-equivalence.
 In other words, $\derL \sus(C)$ is $\W_\Omega$-local.
 But then, as $C$ is cofibrant, $C=\lop \sus(C)=\derR \lop \derL \sus (C)$, and this concludes.

Point (3) is then a consequence of (1) and (2).
\end{proof}

\begin{rem}
\begin{enumerate}
%% \item In the situation of definition \ref{df:triangulated_premotives},
%%  $\DMu(\A)$, as a derived $\Pmor$-premotivic category,
%%  has all the enrichments defined
%%  in \ref{num:P-premotivic_derived&a_lot_of_structures}.
\item The construction of the triangulated category $\DMu(\A)$
 can also be obtained using the more general construction
 of \cite[\textsection 7]{CD1} --- see also \cite[7.11]{Hov}
 and \cite[chap. 4]{ayoub2} for even more general accounts.
Here, we exploit the simplification arising from the fact that
 we invert a complex concentrated in degree $0$: this allowed us to
describe $\DMu(\A)$ simply as a Verdier quotient of the derived category
of an abelian category. However, we can also consider the category of symmetric
 spectra in $\Comp(\A_S)$ with respect to one of the complexes
 $\MabNX S \A(1)[2]$ or $\MabNX S \A(1)$ and this leads to
 the equivalent categories; see \cite[8.3]{Hov}.
\item Point (3) of Proposition \ref{prop:universality_stabilisation} 
is a particular case of \cite[8.1]{Hov}.
\end{enumerate}
\end{rem}

\begin{num} \label{num:functoriality_DMu}
Consider a morphism of abelian $\Pmor$-premotivic categories
$$
\varphi:\A \rightleftarrows \B:\psi
$$ 
such that $\A$ (resp. $\B$) is compatible with a system of topology $t$
(resp. $t'$). Suppose $t'$ is finer than $t$.
According to \ref{num:functoriality_abelian_spectra}, we obtain
an adjunction of abelian $\Pmor$-premotivic categories
$$\varphi:\Comp(\Spt(\A)) \rightleftarrows \Comp(\Spt(\B)):\psi.$$
The pair $(\varphi_S,\psi_S)$ is a Quillen adjunction 
for the stable model structures (apply again \cite[prop. 3.11]{CD1}).
Thus we obtain a morphism of triangulated $\Pmor$-premotivic categories:
$$
\derL\varphi:\DMu(\A)
 \rightleftarrows \DMu(\B):\derR\psi.
$$
\end{num}

\begin{rem} \label{rem:stabilisation&universality}
Under the light of Proposition \ref{prop:universality_stabilisation},
the category $\DMu(\A)$ might be considered as the universal derived
$\Pmor$-premotivic category $\T$ with a morphism $\Der(\A) \rightarrow \T$,
and such that $\T$ satisfies the homotopy and the stability property.
This can be made precise in the setting of algebraic derivators
or of dg-categories (or any other kind of stable $\infty$-categories).
\end{rem}
%
%\begin{prop}\label{spanierwhitehead1}
%Assume that $\Pmor$ contains the class of smooth morphisms
%of finite type. Let $\A$ be an abelian $\Pmor$-premotivic category.
%Assume that $\DMue(\A)$ is compactly generated
%by its geometric sections $\mab S X$,  and satisfies Nisnevich
%descent. Then, for any scheme $S$,
%for any compact object $C$ of $\DMue(\A_S)$
%and for any Tate spectrum $E$ in $\A_S$, we have a canonical isomorphisms
%$$\varinjlim_{r >>0}\Hom_{\DMue(\A_S)}(C\{a+r\},E_{r}[i])
%\simeq \Hom_{\DMu(\A_S)}(\derL \sus (C)\{a\},E[i])$$
%for any integers $a$ and $i$.
%\end{prop}
%
%\begin{proof}
%The permutation $\sigma=(123)$ acts as the identity on
%$\MabNX S \A\{1\}^{\otimes 3}=\MabNX S \A\{3\}$
%in $\DMue(\A_S)$: using example \ref{ex:universality_(Htp),(BG)}, it
%is sufficient to prove this in $\DMt(S)$, which is well known; see \cite[4.5.65]{ayoub2}.
%This proposition is then a direct consequence of \cite[theorems 4.3.61 and 4.3.79]{ayoub2}.
%\end{proof}
%
%\begin{cor}\label{spanierwhitehead2}\label{nonsymtosymequ2}
%Under the assumptions of the preceding proposition, the triangulated
%category $\DMu(\A_S)$ is compactly generated.
%More precisely, if $\DMue(\A_S)_c$ denotes the category
%of compact objects in $\DMue(\A_S)$, then the category of
%compact objects in $\DMu(\A_S)$ is canonically equivalent to the
%pseudo-abelian completion of the category obtained as the $2$-colimit
%of the following diagram.
%$$\DMue(\A_S)_c\xrightarrow{\otimes \MabNX S \A\{1\}}
%\DMue(\A_S)_c\xrightarrow{\otimes \MabNX S \A\{1\}}\cdots
%\DMue(\A_S)_c\xrightarrow{\otimes \MabNX S \A\{1\}} \cdots $$
%\end{cor}

\renewcommand{\Rc}{\Lambda}
\begin{prop}\label{caracdescentbyhmtppties}
Let $t$ and $t'$ be two admissible topologies, with $t'$
finer than $t$. Then $\DMu\left(\sh {t'} {\Pmor}\right)$
is canonically equivalent to the the full subcategory of
$\DMu\left(\sh {t} {\Pmor}\right)$ spanned by the objects
which satisfy $t'$-descent.\index{word}{descent!tdescent@$t$-descent}
\end{prop}

\begin{proof}
It is sufficient to prove this proposition in the case where $t$
is the coarse topology. We deduce from \cite[4.4.42]{ayoub2} that,
for any scheme $S$ in $\sch$, we have
$$\DMu\left(\sh {t'} {\Pmor/S}\right)=\Der\left(\psh {\Pmor/S}\right)[\W^{-1}]\, ,$$
with $\W=\W_{t'}\cup\W_{\AA^1}\cup\W_\Omega$, where
$\W_{t'}$ is the set of maps of shape
$$\sus \mab{S}{\cX}\{n\}[i] \To \sus \mab{S}{X}\{n\}[i]\, ,$$
for any $t'$-hypercover $\cX \To X$ and any integers $n\leq 0$
and $i$. The assertion is then a particular case of the description of
the homotopy category of a left Bousfield localization.
\end{proof}

%\begin{prop}\label{sheavesTatespectranicederivedtensorproduct}
%Let $t$ be an admissible topology. Then, for any scheme $S$ in $\sch$,
%the symmetric monoidal model structure on
%$\Comp(\Spt(\sh {t} {\Pmor/S}))$
%underlying $\DMu\left(\sh {t} {\Pmor/S}\right)$
%is perfect (see \ref{defperfectsymmonidcmf}).
%\end{prop}
%
%\begin{proof}
%The generating family of $\sh {t} {\Pmor/S}$ is flat
%in the sense of \cite[3.1]{CD1}, so that, by virtue of
%\cite[prop. 7.22 and cor. 7.24]{CD1}, the assumptions
%of Proposition \ref{modulesbasechangestable} are fulfilled.
%\end{proof}

\begin{ex}
\label{ex:stable_AA^1-derived_categories}
We have the stable versions of the $\Pmor$-premotivic categories
 introduced in example \ref{ex:AA^1-derived_categories}: \\
\noindent 1) Consider the admissible topology $t=\nis$.
Following F.~Morel, we define the 
\emph{stable $\AA^1$-derived premotivic category}
\index{word}{premotivic!stableA1@stable $\AA^1$-derived premotivic category}
 as (see also the construction of \cite{ayoub2}):
\begin{equation*}
\DMux \Rc:=\DMu\left(\sh {\nis} {\sm}\right)
\quad\text{and}\quad\uDMux \Rc :=\DMu\left(\sh {\nis} {\sft}\right),
\end{equation*}
\index{notat}{DA1Lambda@$\DMux \Lambda$}
as well as the \emph{generalized stable $\AA^1$-derived premotivic
 category}\footnote{We will see in Example \ref{ex:enlargement_of_DMtilde}
 that the generalized version contains the usual one as a full subcategory.}
\begin{equation}
\uDMux \Rc :=\DMu\left(\sh {\nis} {\sft}\right).
\end{equation}
\index{notat}{DA1Lambda@$\uDMux \Lambda$}

Given a scheme $S$, we shall also write:
\begin{equation} \label{eq:DM_tildebis}
\DMu(S,\Rc):=\DMux \Rc(S)
\quad\text{and}\quad\uDMu(S,\Rc):=\uDMux \Rc(S).
\end{equation}

In the case when $t=\et$,
 we get the triangulated premotivic
 categories of \emph{\'etale premotives}:
$$
\DMu\left(\sh {\et} {\sm}\right)\quad\text{and}\quad\DMu\left(\sh {\et} {\sft}\right).
$$
In each of these cases, we denote by $\sus \Rc_S^t(X)$
 the premotive associated with a smooth $S$-scheme $X$.

From the adjunction \eqref{ex:associated_sheaf=adjunction},
 we get an adjunction of triangulated premotivic categories:
$$
a_{\et}:\DMux \Rc
 \rightleftarrows  \DMu\left(\sh {\et} {\sm}\right):\derR O_{\et}.
$$

\noindent 2) Assume $\Pmor=\sft$:
 
Consider the $\sft$-admissible topology $t=\h$ (resp. $t=\qfh$).
In \cite{V1}, Voevodsky has introduced the category of effective
 $\h$-motives (resp. $\qfh$-motives).
According to the theory presented above,
 one can extend this definition to the stable setting:
 one defines the category of stable \emph{$\h$-motives}
\index{word}{motive!hmotive@$\h$-motive}
(resp. \emph{$\qfh$-motives})
\index{word}{motive!hmotive@$\qfh$-motive}
over $S$ with coefficients in $\Rc$ as:
\begin{align*}
&\uDMV_{\h}(S,\Rc):=\DMu\left(\sh {\h} {\sft/S}\right). \\
\text{resp. }
 & \uDMV_{\qfh}(S,\Rc):=\DMu\left(\sh {\qfh} {\sft/S}\right).
\end{align*}
In other words, this is the stable $\AA^1$-derived category of
 $\h$-sheaves (resp. $\qfh$-sheaves) of $\Rc$-modules.
\index{word}{sheaf!hsheaf@$\h$-sheaf}
\index{word}{sheaf!qfhsheaf@$\qfh$-sheaf}
Moreover, we get the
\emph{generalized triangulated premotivic category of $\h$-motives
 (resp. $\qfh$-motives)}
 with coefficients in $\Rc$ over $\site$:
\begin{align*}
&\uDMV_{\h,\Rc}:=\DMu\left(\sh {\h} {\sft}\right). \\
\text{resp. }
 & \uDMV_{\qfh,\Rc}:=\DMu\left(\sh {\qfh} {\sft}\right).
\end{align*}
For an $S$-scheme of finite type $X$,
we will denote by $\sus \uRc_S^\h(X)$ (resp $\sus\uRc_S^\qfh(X)$)
the corresponding premotive associated with $X$ in
$\uDMV_t(S,\Rc)$.
Note that the $h$-sheafification functor induces
 a premotivic adjunction (see Paragraph \ref{num:functoriality_DMu}):
\begin{equation} \label{eq:DM_qfh/h1etdemi}
\uDMV_{\qfh,\Rc}\rightleftarrows\uDMV_{\h,\Rc}\, .
\end{equation}

These generalized premotivic categories are too big to be reasonable
 (in particular for the localization property
 --- see Remark \ref{rem:i^!j_*=0}).
Therefore, we introduce the triangulated category $\DMV_t(S,\Rc)$
as the localizing subcategory of
$\uDMV_t(S,\Rc)$ generated by objects of shape
$\sus \Rc_S^t(X)(p)[q]$ for any smooth $S$-scheme of finite type $X$ and any
integers $p$ and $q$. 
The fibred category $\DMV_{\h,\Rc}$ (resp. $\DMV_{\qfh,\Rc}$)
\index{notat}{DMhLambda@$\DMV_{\h,\Rc}$}
\index{notat}{DMqfhLambda@$\DMV_{\qfh,\Rc}$}
 defined above is premotivic.
We call it the \emph{premotivic category of $\h$-motives (resp. $\qfh$-motives)}.
\index{word}{premotivic!category of $\h$-motives}
\index{word}{premotivic!category of $\qfh$-motives}
The family of inclusions
\begin{equation} \label{eq:DM_qfh/h2}
\DMV_t(S,\Rc)\To \uDMV_t(S,\Rc)
\end{equation}
indexed by a scheme $S$ defines a premotivic morphism
(the existence of right adjoints is ensured by
 the Brown representability theorem).
\end{ex}

\begin{rem}
When $\Rc=\QQ$, we will show that the categories
 $\DMV_{\h,\QQ}$ and $\DMV_{\qfh,\QQ}$ are equivalent
 and satisfy the axioms of a motivic category.
In fact, they are equivalent to the category
 of Beilinson motives.
 See Theorem \ref{plongement01} for all these results.
\end{rem}

\begin{prop}\label{tautologicalqfhdescent}
Consider the notations of the second point in the above example.
Then the premotivic category $\DMV_{t,\Rc}$ satisfies $t$-descent.
\end{prop}

\begin{proof}
This is true for $\uDMV_{t,\Rc}$ by construction, which implies formally
the assertion for $\DMV_{t,\Rc}$.
\end{proof}

%% 
%% The adjunction \eqref{ex:associated_sheaf=adjunction}
%%  induces morphisms of generalized triangulated premotivic
%%   categories:
%% $$
%% \DMue\left(\sh {\cdh} {\sft}\right) \xrightarrow{\underline a_{\qfh}}
%%  DM^{eff}_{\qfh,\Rc} \xrightarrow{a_{\h}}
%%  DM^{eff}_{\h,\Rc}.
%% $$

\begin{rem} \label{ex:universality_htp+stab+bg}
\renewcommand{\Rc}{\ZZ}
According to Proposition \ref{prop:descent&derived_P-premotivic}
 and Remark \ref{rem:stabilisation&universality},
 for any admissible topology $t$, $\DMu(\sh t \Pmor)$ is the
 universal derived $\Pmor$-premotivic category satisfying
 $t$-descent as well as the homotopy and stability properties.

A crucial example for us: the stable $\AA^1$-derived premotivic category $\DMu$ 
 is the universal derived premotivic 
 category satisfying the properties of homotopy, of stability and of Nisnevich descent.
\end{rem}

%Thus, as in the effective case (\textit{cf.} \ref{ex:functoriality_DMt}), 
%given any abelian $\Pmor$-premotivic category $\A$ satisfying Nisnevich descent,
%we obtain a canonical functor
%$$
%\derL\iota:\DMt
% \rightleftarrows \DMu(\A):\derR\delta.
%$$

\begin{paragr}
\label{par:link_stablehomotopy_homology}
\renewcommand{\Rc}{\ZZ}
We assume $\Pmor=\sm$. \\
Let $\operatorname{Sh}_\bullet(\sm)$
 be the category of pointed Nisnevich sheaves of sets.
Consider the pointed version of the adjunction of $\Pmor$-premotivic
categories
$$
N:\Delta^{op} \operatorname{Sh}_\bullet(\sm)
 \rightleftarrows \Comp(\sh{\nis}\sm):K
$$
constructed in \ref{par:link_homotopy_homology}.

If we consider on the left-hand side the $\AA^1$-model category
defined by Blander~\cite{BL}, $(N_S,K_S)$ is a Quillen
adjunction for any scheme $S$.

We consider $(\GG,1)$ as a constant pointed simplicial sheaf.
The construction of symmetric $\GG$-spectra respectively to
the model category
$\Delta^{op} \operatorname{Sh}_\bullet(\sm)$ 
can now be carried out following \cite{Jar} or \cite{ayoub2}
 and yields a symmetric monoidal model category
 whose homotopy category
 is the stable homotopy category of Morel and Voevodsky $\SH(S)$.

Using the functoriality statements \cite[th. 8.3 and 8.4]{Hov}, 
we finally obtain a $\Pmor$-premotivic adjunction
\begin{equation} \label{eq:adj_SH_DMt}
N:\SH \rightleftarrows \DMu:K.
\end{equation}
%Moreover, it can be checked this is an adjunction of 
%triangulated $\ZZ \times \ZZ$-twisted categories.

The functor $K$ is the analog of the Eilenberg-Mac Lane 
functor in algebraic topology; in fact, this adjunction is actually induced
by the Eilenberg-MacLane functor (see \cite[chap. 4]{ayoub2}).
In particular, as the rational model category of topological
(symmetric) $S^1$-spectra is Quillen equivalent to the model category
of complexes of $\QQ$-vector spaces, we have a natural equivalence
of premotivic categories
\begin{equation} \label{eq:adj_SH_DMt_rational}
\renewcommand{\Rc}{\QQ}
\SH_\Rc \rightleftarrows \DMt\, ,
\end{equation}
(where $\SH_\QQ(S)$ denotes the Verdier quotient of $\SH(S)$
by the localizing subcategory generated by compact torsion objects).
\end{paragr}

\begin{num} \label{num:stable-A^1-derived_chg_coef}
\renewcommand{\Rc}{\Lambda}
We can extend the considerations of Example \ref{ex:derived_chg_coef}
 and Paragraph \ref{num:A^1-derived_chg_coef} on changing
 coefficients in categories of sheaves.

Let $t$ be an admissible topology and 
 $\varphi:\Lambda \rightarrow \Lambda'$ be an extension of rings.
 Using the $\Pmor$-premotivic adjunction \eqref{eq:abelian_chg_coef}
  and according to Paragraph \ref{num:functoriality_DMu},
  we get an adjunction of triangulated $\Pmor$-premotivic categories:
$$
\derL \varphi^*:\DMu\big(\sh t {\Pmor}\big)
 \rightleftarrows
 \renewcommand{\Rc}{\Lambda'} 
 \DMu\big(\sh t {\Pmor}\big):\derR \varphi_*.
$$
Given two Tate spectra $C$ and $D$ of $t$-sheaves
 of $\Rc$-modules over $\Pmor_S$,
 we get a canonical morphism of $\Lambda'$-modules:
\begin{equation} \label{eq:stable-A^1-derived_chg_coef}
\Hom_{\DMu(\sh t {\Pmor_S})}\big(C,D\big) \otimes_\Lambda \Lambda'
 \longrightarrow
 \renewcommand{\Rc}{\Lambda'} 
\Hom_{\DMu(\sh t {\Pmor_S})}\big(\derL \varphi^*(C),\derL \varphi^*(D)\big)
\end{equation}
Then the stable version of Proposition \ref{prop:A^1-derived_chg_coef}
 holds (the proof is the same):
\end{num}
\begin{prop} \label{prop:stable-A^1-derived_chg_coef}
Consider the above assumptions.
Then the map \eqref{eq:stable-A^1-derived_chg_coef} is an isomorphism
 in the two following cases:
\begin{enumerate}
\item If $\Lambda'$ is a free $\Lambda$-module and $C$ is compact;
\item If $\Lambda'$ is a free $\Lambda$-module of finite rank.
\end{enumerate}
\end{prop}

\subsubsection{Constructible premotivic spectra}
\label{sec:constructible_stableA1derived}

\begin{lm}
Let $\A$ be an abelian $\Pmor$-premotivic category compatible with a topology $t$
 and such that the category $\AA^1$-derived category $\DMue(\A)$
 satisfies Nisnevich descent.

Then, for any scheme $S$,
 the non-trivial cyclic permutation $(123)$ of order $3$ acts as the identity on
 the premotive
$\MabNX S \A\{1\}^{\otimes 3}$
in $\DMue(\A_S)$.
\end{lm}
\begin{proof}
Using example \ref{ex:universality_(Htp),(BG)},
 it is sufficient to prove this in $\DMte(S)$, which is well-known;
  see for example \cite[4.5.65]{ayoub2}.
\end{proof}

\begin{prop}\label{spanierwhitehead1}
Consider the hypothesis of the previous lemma
 and assume that the triangulated premotivic category $\DMue(\A)$
 is compactly $\tau$-generated.

Then, for any scheme $S$,
 any couple of integers $(i,a)$,
 any compact object $C$ of $\DMue(\A_S)$
  and any Tate spectrum $E$ in $\A_S$, we have a canonical isomorphism
$$\Hom_{\DMu(\A_S)}(\derL \sus (C)\{a\},E[i])
 \simeq \varinjlim_{r >>0}\Hom_{\DMue(\A_S)}(C\{a+r\},E_{r}[i]).$$
\end{prop}
\begin{proof}
Given the previous lemma,
 this is a direct consequence of \cite[theorems 4.3.61 and 4.3.79]{ayoub2}.
\end{proof}

\begin{cor}\label{spanierwhitehead2}\label{nonsymtosymequ2}
Under the assumptions of the preceding proposition, the triangulated
 category $\DMu(\A_S)$ is compactly $(\ZZ \times \tau)$-generated
\index{word}{generated!compactly $(\ZZ \times \tau)$-generated}
 where the factor $\ZZ$ corresponds to the Tate twist.

More precisely, if $\DMuex c(\A_S)$ denotes the category
of compact objects in $\DMue(\A_S)$, then the category of
compact objects in $\DMu(\A_S)$ is canonically equivalent to the
pseudo-abelian completion of the category obtained as the $2$-colimit
of the following diagram:
$$\DMuex{c}(\A_S)\xrightarrow{\otimes \MabNX S \A\{1\}}
\DMuex{c}(\A_S)
\longrightarrow \cdots \longrightarrow
\DMuex{c}(\A_S)\xrightarrow{\otimes \MabNX S \A\{1\}} \DMuex{c}(\A_S) \longrightarrow \cdots $$
\end{cor}

\begin{num}
Let $\A$ be an abelian $\Pmor$-premotivic category
 compatible with an admissible topology $t$.
Assume that:
\begin{itemize}
\item The topology $t$ is bounded in $\A$ (Definition \ref{df:hyper-bounded}).
\item The abelian $\Pmor$-premotivic category $\A$ is
 finitely $\tau$-presented.
\index{word}{finitely presented!finitely $\tau$-presented}
\end{itemize}
We will denote by $\N^t_S$ a bounded generating family for $t$-hypercovers in $\A_S$.

Recall from Proposition \ref{prop:compact_DMue} that
 the category of compact objects of 
 the triangulated category $\DMue(\A_S)$
 is canonically equivalent to the triangulated monoidal category:
$$
\Big(K^b\big(\ZZ_S(\sm/S;\A)\big)
 /(\cN^t_S \cup \T_{\AA^1_S})\Big)^\natural
$$
Let us denote by $\DMux{gm}(\A_S)$
\index{notat}{DA1gmAS@$\DMux{gm}(\A_S)$}
 the category obtained
 from the monoidal category on the left-hand side of the above functor
 by formally inverting the Tate twist $\ZZ_S^\A(1)$.
Because $\DMu(\A)$ satisfies the stability property by construction,
 we readily obtains a canonical monoidal functor
\begin{equation} \label{eq:A^1-spectra_geometric_premotives}
\DMux{gm}(\A_S) \rightarrow \DMu(\A_S).
\end{equation}
Then applying Proposition \ref{prop:compact_DMue},
 the above corollary and Proposition \ref{constructequivcompact},
 we deduce:
\end{num}
\begin{cor} \label{cor:geometric_premotives}
Consider the above hypothesis and notations.

Then the triangulated premotivic category $\DMu(\A)$ is compactly
 $(\ZZ \times \tau)$-generated.
For any premotive $\cM$ in $\DMu(\A_S)$
 the following conditions are equivalent:
\begin{enumerate}
\item[(i)] $\cM$ is compact;
\item[(ii)] $\cM$ is $(\ZZ \times \tau)$-constructible.
\index{word}{constructible!$(\ZZ \times \tau)$-constructible}
\end{enumerate}
Moreover, the functor \eqref{eq:A^1-spectra_geometric_premotives}
 is fully faithful and has for essential image the compact
 (\emph{i.e.} $\tau$-constructible)
 objects of $\DMu(\A_S)$.
\end{cor}

\begin{ex} \label{ex:A^1-Nis_Et_compacity_stable}
\renewcommand{\Rc}{\Lambda}
From the considerations of Example \ref{ex:A^1-Nis_Et_compacity},
 we get that, for any scheme $S$,
  the category of compact objects of $\DMu(S,\Rc)$
   (resp., of its $\cdh$-local counterpart $\DMu\Big(\sh {\cdh} {\smash{\sft/S}}\Big)$)
 is obtained from the monoidal triangulated category
$$
\K^b\left(\Rc(\sm/S)\right) \text{ (resp. } \K^b\left(\Rc(\sft/S)\right)
 \text{)}
$$
by the following steps:
\begin{itemize}
\item one mods out by the triangulated subcategories $\cT_{\AA^1_S}$
 and $BG_S$ (resp. $CDH_S$)
 corresponding to the $\AA^1$-homotopy property and 
 the Brown-Gersten triangles (resp. cdh-triangles),
\item one takes the pseudo-abelian envelope,
\item one formally inverts the Tate twist.
\end{itemize}
\end{ex}

\begin{prop} \label{prop:continuity_gen-stable-A^1_premotivic}
Assume $\Pmor=\sft$ is the class of finite type
 (resp. separated and of finite type) morphisms.

Let $\A$ be an abelian generalized premotivic category
 compatible with an admissible topology $t$ such that:
\begin{itemize}
\item $\A$ satisfies property (C) (resp. (wC)) of 
 Paragraph \ref{num:continuity_derived_premotivic}.
\item The $\AA^1$-derived category $\DMue(\A)$ is
 compactly $\tau$-generated and satisfies Nisnevich descent.
\end{itemize}
Then the stable $\AA^1$-derived premotivic category
 $\DMu(\A)$ is $(\ZZ \times \tau)$-continuous
 (resp. weakly $(\ZZ \times \tau)$-continuous)
 --- see Definition \ref{df:continuous}.
\end{prop}
\begin{proof}
This is an immediate corollary of Proposition
 \ref{prop:continuity_gen-A^1_premotivic}
 combined with Proposition \ref{spanierwhitehead1}.
\end{proof}

\begin{ex}
According to the previous proposition
 and the second point of Example \ref{ex:sheaves_continuity},
 the generalized triangulated premotivic category
 $\uDMux \Rc$ is continuous.
We also refer the reader to
 Corollary \ref{cor:continuity_stable-A^1_premotivic}
 for an extension of this result to the non generalized case.
\end{ex}
\section{Localization and the universal derived example}
\label{sec:dmtilde}

\begin{assumption}
In this section, $\sch$ is an adequate category
 of $\base$-schemes as in \ref{num:assumption1_sch}.
In sections \ref{sec:Nisnevich&loc} and \ref{sec:nearly_Nisnevich&loc},
we assume in addition that the schemes in $\sch$ are finite dimensional.
 
We will apply the definitions of the preceding section
 to the admissible class made of morphisms of finite type
  (resp. smooth morphisms of finite type) in $\sch$,
 denoted by $\sft$ (resp. $\sm$).

Recall the general convention of section \ref{sec:premotivic_convention}:
\begin{itemize}
\item \emph{premotivic} means $\sm$-premotivic;
\item \emph{generalized premotivic} means $\sft$-premotivic.
\end{itemize}
\end{assumption}

\subsection{Generalized derived premotivic categories}
\label{sec:generalized_derived}

\renewcommand{\Rc}{\Lambda}

%
%\begin{ex}
%The notion of generalized premotivic categories is particularly relevant
% for the following topologies on $\sft$, first introduced by Voevodsky~: 
% $t=\cdh, \qfh, \h$. 
%In each of this case, the category of abelian sheaves
%$\sh t {\sft/S}$ is a Grothendieck abelian category generated
%by the representable sheaves $\ZZ_S^t(X)$ for a finite type  $S$-scheme $X$. 
%As in the case of example \ref{ex:abelian_premotivic_sheaves},
%we finally get a generalized abelian premotivic category $\sh t {\sft}$.
%It is obviously compatible with the topology $t$, in the generalized
%sense analogous to that of definition \ref{df:basic_complexes&topology}.
%\end{ex}
%
\begin{ex} \label{ex:enlargementofsheaves}
Let $t$ be a $\sft$-admissible topology. For a scheme $S$,
 we denote by $\sh t {\sft/S}$ the category of sheaves of abelian groups 
 on $\sft/S$ for the topology $t_S$.
For an $S$-scheme of finite type $X$,
we let $\urepx  t S X$ be the free t-sheaf of $\Rc$-modules represented by $X$.
Recall $\sh t {\sft}$ is a generalized abelian premotivic category
 (see \ref{ex:abelian_premotivic_sheaves}). \\
Let $\rho:\sm/S \rightarrow \sft/S$ be the obvious inclusion functor
 and let us denote by $t_S$ the initial topology on $\sm/S$ such that $\rho$
 is continuous.
Then it induces (\textit{cf.} \cite[IV, 4.10]{SGA4}) a sequence of adjoint functors
$$
\xymatrix@C=48pt@R=10pt{
\sh t {\sm/S}\ar@/^10pt/^{\rho_\sharp}[r]\ar@/_10pt/_{\rho_*}[r] 
 & \sh t {\sft/S}\ar|/-2pt/{\rho^*}[l]  \\
&& 
}
$$
and we checked easily that this induces an enlargement
 of abelian premotivic categories:
\begin{equation} \label{eq:enlargement_sheaves}
\rho_\sharp:\sh t {\sm} \rightleftarrows \sh t {\sft}:\rho^*.
\end{equation}
\end{ex}

\begin{rem} \label{rem:Vologodsky}
Note that for any scheme $S$, the abelian category $\sh t {\sm/S}$
 can be described as the Gabriel quotient\index{word}{quotient!Gabriel}
 of the abelian category $\sh t {\sft/S}$
 with respect to the sheaves $\underline{F\!\!}\;$ over $\sft/S$ such that
  $\rho^*(\underline{F\!\!}\;)=0$.

An example of such a sheaf in the case where $t=\nis$ and $\dim(S)>0$
is the Nisnevich sheaf $\urep S Z$ on $\sft/S$ represented
 by a nowhere dense closed subscheme $Z$
 of $S$ is zero when restricted to $\sm/S$.
 \end{rem}

\begin{num} \label{abelian_enlargement}
Consider an abelian premotivic category $\A$ compatible 
 with an admissible topology $t$ on $\sm$ 
 and a generalized abelian premotivic category $\uA$
 compatible with an admissible topology $t'$ on $\sch$.
We denote by $M$ (resp. $\underline{M}$)
 the geometric sections of $\A$ (resp. $\uA$).
We assume that $t'$ restricted to $\sm$ is finer that $t$,
 and consider an adjunction
%% and consider an enlargement (see definition \ref{df:enlargement})
 of abelian premotivic categories:
$$
\rho_\sharp:\A \rightleftarrows \uA:\rho^*.
$$

Let $S$ be a scheme in $\sch$. The functors $\rho_\sharp$ and $\rho^*$
induce a derived adjunction (see \ref{num:functoriality_DMue}):
%% $$
%% \derL \rho_\sharp:\Der(\A_S) \rightleftarrows \Der(\uA_S):\derR \rho^*\, .
%% $$
%% Furthermore, as $\rho$ is exact, for any object $M$ of $\Der(\uA_S)$,
%% we have $\rho^*(M)=\derR\rho^*(M)$.
%% Moreover, we deduce from the properties of an enlargement that 
%% $\rho^*$ (resp. $\rho_\sharp$) 
%% preserves $\AA^1$-local objects (resp. $\AA^1$-equivalences);
%% see Proposition \ref{prop:A^1-derivation_exact_right_adj}.
%% In particular, they induces an adjunction of the $\AA^1$-localized
%% categories~:
$$
\derL \rho_\sharp:\DMue(\A_S) \rightleftarrows \DMue(\uA_S):\derR \rho^*
$$
(where $\uA$ is considered as an $\sm$-fibred category).
\end{num}

\begin{prop} \label{prop:enlargement_eff}
Consider the previous hypothesis, and fix a scheme $S$.
Assume furthermore that we have the following properties.
\begin{itemize}
\item[(i)] The functor $\rho_\sharp:\A_S\To \uA_S$
is fully faithful.
\item[(ii)] The functor $\rho^*:\uA_S\To \A_S$
commutes with small colimits.
\end{itemize}
Then, the following conditions hold~:
\begin{enumerate}
\item[(a)]
The induced functor
$$
\rho^*:\Comp(\uA_S) \rightarrow \Comp(\A_S)
$$
preserves $\AA^1$-equivalences.
\item[(b)] The $\AA^1$-derived functor 
$\derL \rho_\sharp:\DMue(\A_S) \rightarrow \DMue(\uA_S)$
is fully faithful.
\end{enumerate}
\end{prop}
\begin{proof}
Point (a) follows from Proposition \ref{prop:A^1-derivation_exact_right_adj}.
To prove (b), we have to prove that the unit map
$$M\To \rho^*\derL\rho_\sharp(M)$$
is an isomorphism for any object $M$ of $\DMue(\A_S)$.
For this purpose, we may assume that $M$ is cofibrant, so that
we have
$$M\simeq  \rho^* \rho_\sharp(M)\simeq  \rho^*\derL\rho_\sharp(M)$$
(where the first isomorphism holds already in $\Comp(\A_S)$).
\end{proof}

\begin{cor} \label{cor:prop:enlargement_eff}
Consider the hypothesis of the previous proposition.
Then the family of adjunctions
$
\derL \rho_\sharp:\DMue(\A_S)
 \rightarrow \DMue(\uA_S):\derR \rho^*
$
 indexed by a scheme $S$
 induces an enlargement of triangulated premotivic categories
$$
\derL \rho_\sharp:\DMue(\A)
 \rightleftarrows \DMue(\uA):\derR \rho^*.
$$
\end{cor}

\begin{ex} Considering the situation of \ref{ex:enlargementofsheaves},
we will be particularly interested in the case of the Nisnevich topology.
We denote by $\uDMte$ the generalized $\AA^1$-derived premotivic category
associated with $\sh {} {\sft}$
 (see also Example \ref{ex:stable_AA^1-derived_categories}). 
The preceding corollary gives a canonical enlargement:
\begin{equation} \label{eq:enalrgement_Nis_eff}
\DMte \rightleftarrows \uDMte
\end{equation}
\end{ex}

\begin{num} Consider again the hypothesis of \ref{abelian_enlargement}.
%% We assume moreover(mainly for simplicity):
%% \begin{enumerate}
%% \item[(C)] $\DMu(\A)$ and $\uA$ are geometrically
%%  generated and that their geometric sections are compact.
%% \end{enumerate}
We denote simply by $\motNP$ (resp. $\umotNP$)
 the geometric sections of the premotivic triangulated
 category $\DMu(\A)$ (resp. $\DMu(\uA)$).

Recall from \ref{num:Tate_premot} that we have defined
$\MabNX S {\A}\{1\}$ (resp. $\uMabNX S {\A}\{1\}$)
as the cokernel of the canonical map
$\MabNX S \A \rightarrow \mab S {\GGx S}$
(resp. $\uMabNX S {\A} \rightarrow \umab S {\GGx S}$).
Thus, it is obvious that we get a canonical identification
$\rho_\sharp(\MabNX S {\A}\{1\})=\uMabNX S {\A}\{1\}$.
Therefore, the enlargement $\rho_\sharp$ can be extended canonically
 to an enlargement 
$$
 \rho_\sharp:\Spt(\A) \rightleftarrows \Spt(\uA):  \rho^*
$$
 of abelian premotivic categories
 in such a way that for any scheme $S$,
 the following diagram commutes:
$$
\xymatrix@C=40pt{
\A_S\ar^{\rho_\sharp}[r]\ar_{\sus_\A}[d]
 & \uA_S\ar^{\sus_{\underline{\mathscr A\!}}}[d] \\
\Spt(\A_S)\ar^-{ \rho_\sharp}[r] & \Spt(\uA_S).
}
$$
According to Proposition \ref{prop:symmetric_seq&topology},
 $\Spt(\A)$ (resp. $\Spt(\uA)$) is compatible with $t$
 (resp. $t'$), and we obtain an adjoint pair of functors (\ref{num:functoriality_DMu}):
$$
\derL \rho_\sharp:\DMu(\A_S)
 \rightleftarrows \DMu(\uA_S):\derR \rho^*.
$$

From the preceding commutative square, 
 we get the identification:
\begin{equation} \label{eq:enlargement&infinite_sus}
\derL \rho_\sharp \circ \sus_\A
 =\sus_{\underline{\mathscr A\!}} \circ \derL \rho_\sharp
\end{equation}
As in the non-effective case, we get the following result:
\end{num}

\begin{prop} \label{prop:enlargement}
Keep the assumptions of Proposition \ref{prop:enlargement_eff},
and suppose furthermore that
both $\DMue(\A)$ and $\DMue(\uA)$ are
compactly $\tau$-generated.
Then the derived functor
$\derL \rho_\sharp:\DMu(\A_S) \rightarrow \DMu(\uA_S)$
is fully faithful.
\end{prop}
\begin{proof}
We have to prove that for any Tate spectrum $E$ of $\DMu(\A_S)$,
 the adjunction morphism
$$
E \rightarrow \derL \rho^* \derR \rho_\sharp(E)
$$
is an isomorphism.
According to Proposition \ref{prop:exist_right_adjoint},
 the functor $\derL \rho^*$ admits a right adjoint.
 Thus, applying Lemma \ref{lm:generators&adjoints}, 
 it is sufficient to consider the case where 
 $E=\mot S X\{i\}[n]$ for a smooth $S$-scheme $X$,
 and a couple $(n,i) \in \ZZ \times \tau$.

Moreover, it is sufficient to prove that for another
 smooth $S$-scheme $Y$ and an integer $j \in \ZZ$,
 the induced morphism
$$
\Hom(\sus \mab S Y \{j\},\sus \mab S X \{i\}[n])
 \rightarrow \Hom(\sus \umab S Y \{j\},\sus \umab S X \{i\}[n])
$$
is an isomorphism. Using the identification \eqref{eq:enlargement&infinite_sus},
propositions \ref{spanierwhitehead1} and \ref{prop:enlargement_eff} allows us to conclude.
\end{proof}

\begin{cor}\label{cor:derived_enlargement}
If the assumptions of Proposition \ref{prop:enlargement}
hold for any scheme $S$ in $\sch$, then
we obtain an enlargement of triangulated premotivic categories
\index{word}{premotivic!enlargement of ---- category}
$$
\derL \rho_\sharp:\DMu(\A)
 \rightleftarrows \DMu(\uA):\derR \rho^*.
$$
\end{cor}

\begin{ex}\label{ex:enlargement_of_DMtilde}
 Considering again the situation of \ref{ex:enlargementofsheaves},
in the case of the Nisnevich topology.
We denote by $\uDMt$ the generalized stable $\AA^1$-derived premotivic category
associated with $\sh {} {\sft}$. The preceding corollary gives
a canonical enlargement:
\begin{equation} \label{eq:enalrgement_Nis}
\derL \rho_\sharp : \DMt \rightleftarrows \uDMt : \derR \rho^*
\end{equation}
which is compatible with the enlargement \eqref{eq:enalrgement_Nis_eff}
 in the sense that the following diagram is essentially commutative:
$$
\xymatrix@=26pt{
\DMte\ar[r]\ar_{\sus}[d] & \uDMte\ar^{\underline \Sigma^\infty}[d] \\
\DMt\ar[r] & \uDMt
}
$$
\end{ex}

\begin{cor}\label{enlargeddescent}
Consider a Grothendieck topology $t$ on our category of schemes $\sch$.
Let $S$ be a scheme in $\sch$, and $M$ an object of $\DMt(S)$.
Then $M$ satisfies $t$-descent in $\DMt(S)$ if and only if
$\derL \rho_\sharp(M)$ satisfies $t$-descent in $\uDMt(S)$.
\end{cor}
\begin{proof}
Let $f:\X\To S$ be a diagram of $S$-schemes of finite type.
Define
$$H^q(\X,M(p))=\Hom_{\DMt(S)}(\Rc_\X,\derL f^*(M)(p)[q])$$
$$\underline{H}^q(\X,M(p))=\Hom_{\uDMt(S)}(\underline{\Rc}_\X,\derL f^*\, \derL\rho_\sharp (M)(p)[q])$$
for any integers $p$ and $q$. The full faithfulness of $\derL\rho_\sharp$
ensures that the comparison map
$$H^q(\X,M(p))\To \underline{H}^q(\X,M(p))$$
is always bijective. This proposition follows then from
the fact that $M$ (resp. $\derL\rho_\sharp(M)$) satisfies $t$-descent if and only
if, for any integers $p$ and $q$, for any $S$-scheme of finite type $X$, and
any $t$-hypercover $\X\To X$, the induced map
$$H^q(X,M(p))\To H^q(\X,M(p)) \ \text{(resp.}\
\underline{H}^q(X,M(p))\To \underline{H}^q(\X,M(p)) \, )$$
is bijective.
\end{proof}

We end-up this section with another interesting application
 of the preceding results.
\begin{cor} \label{cor:continuity_stable-A^1_premotivic}
Consider the hypothesis and assumptions
 of Proposition \ref{prop:enlargement_eff}.
 We suppose furthermore that the generalized abelian premotivic
 category $\uA$ satisfies condition (C)
 of Paragraph \ref{num:continuity_derived_premotivic}.
\begin{enumerate}
\item Then the triangulated premotivic category
 $\DMue(\A)$ is $\tau$-continuous.
\item Assume furthermore that $\DMue(\A)$ and $\DMue(\uA)$ 
 are compactly $\tau$-generated.
Then the triangulated premotivic category $\DMu(\A)$
 is $\tau$-continuous.
\end{enumerate}
\end{cor}
\begin{proof}
According to Proposition \ref{prop:continuity_gen-A^1_premotivic},
 the category $\DMue(\uA)$ is $\tau$-continuous.
 According to Corollary \ref{cor:prop:enlargement_eff},
 the functor $\derL \rho_\sharp:\DMue(\A)
 \rightarrow \DMue(\uA):\derR \rho^*$ is fully faithful and commutes
 with $\derL f^*$. Thus Point (1) follows.

In the assumption of Point (2),
 we deduce from 
 Proposition \ref{prop:continuity_gen-stable-A^1_premotivic}
 that $\DMu(\uA)$ is $(\ZZ \times \tau)$-continuous.
 Thus it is sufficient to apply Corollary \ref{cor:derived_enlargement}
 as in the effective case to get the assertion of Point (2).
\end{proof}

\begin{ex} \label{DMtcontinuous}
According to the second point of Example \ref{ex:sheaves_continuity},
 we can apply this corollary to the enlargement
$$
\sh {\nis} {\sm} \rightarrow \sh {\nis} {\sft}.
$$
Thus, we deduce that the triangulated premotivic categories
 $\DMte$ and $\DMt$ both are continuous.
\end{ex}

\subsection{The fundamental example} \label{sec:Nisnevich&loc}

Recall the following theorem of Ayoub \cite{ayoub2}:
\begin{thm} \label{thm:DMt_localization}
The triangulated premotivic categories $\DMte$ and $\DMt$ satisfy the
 localization property.
\end{thm}

\begin{cor} \label{cor:DM_tilde&univ_derived_motivic}
\begin{enumerate}
\item The premotivic category $\DMt$ is a motivic category.
\item It is compactly generated by the Tate twist.
\item  \renewcommand{\Rc}{\ZZ}
Suppose that $\T$ is a derived premotivic category
 (see \ref{descent_derived_premotives}) which is a motivic category.
Then there exists a canonical morphism
 of derived premotivic categories:
$$
\DMt \rightarrow \T.
$$
\end{enumerate}
\end{cor}
\begin{proof}
The first assertion follows from the previous theorem
 and Remark \ref{rem:compatly&motivic}.
The second one follows from Corollary \ref{cor:geometric_premotives}.
The last one follows from Proposition \ref{localizationNisnevichdescent}
 and Example \ref{ex:universality_htp+stab+bg}.
\end{proof}

\begin{rem} Thus, Theorem \ref{thm:cor3_Ayoub} can be applied to $\DMt$.
In particular, for any separated morphism of finite type $f:T \rightarrow S$,
 there exists a pair of adjoint functors
$$
f_!:\DMt(T) \rightleftarrows \DMt(S):f^!
$$
as in the theorem \emph{loc. cit.} so that we have removed the 
quasi-projective assumption in \cite{ayoub}.
\end{rem}

\begin{num} Because the cdh topology is finer than the Nisnevich topology,
 we get an adjunction of generalized premotivic categories:
$$
a^*_\cdh:\uDMt \rightleftarrows \DMu\left(\sh {\cdh} {\sft}\right):\derR a_{\cdh,*}.
$$
\end{num}

\begin{cor}
For any scheme $S$, the composite functor
$$
\DMtx S \rightarrow \uDMtx S \xrightarrow{a_\cdh} \DMu\left(\sh {\cdh} {\sft/S}\right)
$$
is fully faithful.

Moreover, it induces an enlargement of premotivic categories:
\begin{equation} \label{eq:enlargement_Nis_cdh}
\DMt \rightleftarrows \DMu\left(\sh {\cdh} {\sft}\right)
\end{equation}
\end{cor}

\begin{rem} This corollary is a generalization in our derived setting 
 of the main theorem of \cite{voecd2}. Note that if $\dim(S)>0$,
 there is no hope that the above composite functor 
 is essentially surjective because as soon as $Z$ is a nowhere dense
 closed subscheme of $S$, 
 the premotive $\underline M^\cdh_S(Z,\Rc)$ does not belong to its image
 (\textit{cf.} remark \ref{rem:Vologodsky}).
 \end{rem}
 
\begin{proof}
According to Corollary \ref{cor:DM_tilde&univ_derived_motivic}
 and Proposition \ref{cdhdescent},
 any Tate spectrum $E$ of $\DMtx S$ satisfies $\cdh$-descent in
 the derived premotivic category $\DMt$, and this implies the first
 assertion by \ref{caracdescentbyhmtppties} and \ref{enlargeddescent}.
 The second one then follows from the fact the forgetful functor
$$
\DMu\left(\sh {\cdh} {\sft/S}\right) \rightarrow \uDMtx S.
$$
commutes with direct sums (its left adjoint preserves compact objects).
\end{proof}

\subsection{Nearly Nisnevich sheaves} \label{sec:nearly_Nisnevich&loc}

\renewcommand{\Rc}{\ZZ}

\begin{num}\label{assumption:nearly_nisnevich}
In all this section, we fix an abelian premotivic category $\A$
 and we consider the canonical premotivic adjunction \eqref{eq:universality_presheaves}
 associated with $\A$.

We assume $\A$ satisfies the following properties.
\begin{itemize}
\item[(i)] $\A$ is compatible with Nisnevich topology, so that
we have from  \eqref{eq:universality_presheaves} a premotivic adjunction:
\begin{equation} \label{eq:Nearly_gamma_eff}
\gamma^*:\sh{\nis}{\sm} \rightleftarrows \A:\gamma_*.
\end{equation}
\item[(ii)] $\A$ is finitely presented (\emph{i.e.}
 the functors $\Hom_{\A_S}(\mab S X,-)$ preserve
 filtered colimits and form a conservative family,
 Def. \ref{df:abelian_P-premotivic_well&compactly_generated}).
\item[(iii)] For any scheme $S$, and for any
open immersion $U\To X$ of smooth $S$-schemes, the map
$\mab S U \To\mab S X$ is a monomorphism.
\item[(iv)] For any scheme $S$, the functor
$\gamma_{*}:\A_{S} \To \sh{\nis}{\sm/S}$ is exact.
\end{itemize}
Note that the functor $\gamma_{*}:\A_{S} \To \sh{\nis}{\sm/S}$ is exact
and conservative. As it also preserves filtered colimits, this functor
preserves in fact small colimits.

Observe also that, according to assumptions (i)-(iv),
 the abelian premotivic category of Tate spectra $\Spt(\A)$
 is compatible with Nisnevich topology and $\NN$-generated.
 Moreover, we get a canonical premotivic adjunction
\begin{equation} \label{eq:Nearly_gamma}
\gamma^*:\Spt(\sh{\nis}{\sm}) \rightleftarrows \Spt(\A):\gamma_*
\end{equation}
such that $\gamma_*$ is conservative and preserves small colimits.
\end{num}

In the following, we show how one can deduce properties
 of the premotivic triangulated categories $\DMue(\A)$ and $\DMu(\A)$
 from the good properties of $\DMte$ and $\DMt$.

\subsubsection{Support property (effective case)}

\begin{prop}\label{oubligentil}
For any scheme $S$, the functor
$\gamma_{*}:\Comp(\A_{S}) \To\Comp(\sh{\nis}{\sm/S})$
preserves and detects $\AA^1$-equivalences.
\end{prop}
\begin{proof}
It follows immediately from Corollary \ref{rightadjpreserveA1eq}
that $\gamma_{*}$ preserves $\AA^1$-equivalences.
The fact it detects them can be rephrased by saying that the
induced functor
$$
\gamma_{*}:\DMue(\A_{S})\To\DMte(S)
$$
is conservative. This is obviously true once we noticed
that its left adjoint is essentially surjective on generators.
\end{proof}

\begin{cor}\label{der_oubligentil}
The right derived functor
$$
\derR \gamma_{*}=\gamma_*:\DMue(\A_{S})\To\DMte(S)
$$
is conservative.
\end{cor}

\begin{prop}\label{morphfiniexact}
Let $f:S'\To S$ be a finite morphism of schemes.
Then the induced functor
$$f_{*}:\Comp(\A_{S'})\To\Comp(\A_{S})$$
preserves colimits and $\AA^1$-equivalences.
\end{prop}

\begin{proof}
We first prove $f_{*}$ preserves colimits.
We know the functors
$\gamma_{*}$ preserve colimits and are conservative.
As we have the identification
$\gamma_{*}f_{*}=f_{*}\gamma_{*}$,
it is sufficient to prove the property for $\A=\sh {\nis}{\sm}$.
Let $X$ be a smooth $S$-scheme. It is sufficient to prove
that, for any point $x$ of $X$, if $X^h_{x}$ denotes the
henselization of $X$ at $x$, the functor
$$\sh {\nis}{\sm/{S'}}\To\ab\quad ,
 \qquad F\mapsto f_{*}(F)(X^h_{x})=F(S'\times_{S}X^h_{x})$$
commutes to colimits. Moreover, the scheme $S'\times_{S}X^h_{x}$
is finite over $X^h_{x}$, so that we have
$S'\times_{S}X^h_{x}=\amalg_{i} Y_{i},$
where the $Y_{i}$'s are a finite family of
henselian local schemes over $S'\times_{S}X^h_{x}$.
Hence, we have to check that the functor
$F\mapsto \bigoplus_{i}F(Y_{i})$
preserves colimits. As colimits commute to sums,
it is thus sufficient to prove that the functors
$F\mapsto F(Y_{i})$
commute to colimits. This follows from the fact
that the local henselian schemes $Y_{i}$ are
points of the topos of sheaves over
the small Nisnevich site of $X$.

We are left to prove that the functor
$f_{*}:\Comp(\A_{S'})\To\Comp(\A_{S})$
respects $\AA^1$-equivalences.
For this, we shall study the behavior of $f_*$
with respect to the $\AA^1$-resolution
functor constructed in \ref{num:Suslin_complex}.
Note that $f_{*}$ commutes to limits because it has a left adjoint.
In particular, we know that $f_{*}$ is exact.
Moreover, one checks easily that
$f_{*}R^{(n)}_{\AA^1}=f_{*}R^{(n)}_{\AA^1}$.
As $f_{*}$ commutes to colimits, this gives
the formula $f_{*}\sing=\sing f_{*}$.
Let $C$ be a complex of Nisnevich sheaves
of abelian groups on $\sm/{S'}$.
Choose a quasi-isomorphism $C\To C'$
with $C'$ a $\nis$-flasque complex.
Applying Proposition \ref{prop:singular&AA^1_localization},
we know that $\sing(C')$ is $\AA^1$-fibrant
and that we get a canonical
$\AA^1$-equivalence
$$f_{*}(C)\To f_{*}(C')\To f_{*}(\sing(C'))=\sing(f_{*}(C')).$$
Hence, we are reduced to prove that
$f_{*}$ preserves $\AA^1$-equivalences
between $\AA^1$-fibrant objects.
But such $\AA^1$-equivalences are quasi-isomorphisms,
so that we can conclude using the exactness of $f_{*}$.
\end{proof}

\begin{prop}\label{immouvechangeAS}
For any open immersion of schemes $j:U\To S$,
 the exchange transformation 
$j_{\sharp}\gamma_{*}\To\gamma_{*}j_{\sharp}$
is an isomorphism of functors.
\index{word}{exchange!isomorphism}
\end{prop}

\begin{proof}
Let $X$ be a scheme, and $F$ a Nisnevich sheaf of abelian groups
on $\sm/{X}$. Define the category $\C_{F}$ as follows. The objects
are the couples $(Y,s)$, where $Y$ is a smooth scheme over $X$,
and $s$ is a section of $F$ over $Y$. The arrows $(Y,s)\To(Y',s')$
are the morphisms $f\in\Hom_{\sh {\nis} {\sm/X}}(\rep X Y,\rep X {Y'})$
such that $f^*(s')=s$. We have a canonical functor
$$\varphi_{F}:\C_{F}\To\sh {\nis} {\sm/X}$$
defined by $\varphi_{F}(Y,s)=\rep X Y$, and one easily checks
that the canonical map
$$\ilim_{\C_{F}}\varphi_{F}=\ilim_{(Y,s)\in\C_{F}}\rep X Y
\To F$$
is an isomorphism in $\sh {\nis} {\sm/X}$ (this is essentially a
reformulation of the Yoneda lemma).

Consider now an object $F$ in the category $\A_{U}$.
We get two categories $\C_{\gamma_{*}(F)}$
and $\C_{\gamma_{*}(j_{\sharp}(F))}$. There is a functor
$$i:\C_{\gamma_{*}(F)}\To\C_{\gamma_{*}(j_{\sharp}(F))}$$
which is defined by the formula $i(Y,s)=(Y,j_{\sharp}(s))$.
To explain our notations, let us say that we see $s$ as
a morphism from $\Mab S U \A$ to $F$, so that $j_{\sharp}(s)$
is a morphism from $\Mab S Y \A=j_{\sharp}\Mab S U \A$ to $j_{\sharp}(F)$.
This functor $i$ has right adjoint
$$i':\C_{\gamma_{*}(j_{\sharp}(F))}\To\C_{\gamma_{*}(F)}$$
defined by $i'(Y,s)=(Y_{U},s_{U})$, where
$Y_{U}=Y\times_{S}U$, and $s_{U}$ is the section
of $\gamma_{*}(F)$ over $Y_{U}$ that corresponds
to the section $j^*(s)$ of $j^*j_{\sharp}\gamma_{*}(F)$
over $Y_{U}$ under the canonical isomorphism
$\gamma_{*}(F)\simeq j^*j_{\sharp}\gamma_{*}(F)$
(here, we use strongly the fact the functor $j_{\sharp}$
is fully faithful). The existence of a right adjoint
implies $i$ is cofinal.
This latter property is sufficient for the canonical morphism
$$\ilim_{C_{\gamma_{*}(F)}} \varphi_{\gamma_{*}(j_{\sharp}(F))}\circ i
\To\ilim_{\C_{\gamma_{*}(j_{\sharp}(F))}}\varphi_{\gamma_{*}(j_{\sharp}(F))}
=\gamma_{*}(j_{\sharp}(F))$$
to be an isomorphism. But the functor
$\varphi_{\gamma_{*}(j_{\sharp}(F))}\circ i$
is exactly the composition
of the functor $j_{\sharp}$ with $\varphi_{\gamma_{*}(F)}$.
As the functor $j_{\sharp}$ commutes with colimits, we have
$$%\begin{aligned}
\ilim_{C_{\gamma_{*}(F)}} \varphi_{\gamma_{*}(j_{\sharp}(F))}\circ i
= \ilim_{C_{\gamma_{*}(F)}}j_{\sharp}\, \varphi_{\gamma_{*}(F)}
%&
\simeq j_{\sharp}\ilim_{C_{\gamma_{*}(F)}}\varphi_{\gamma_{*}(F)}%\\
%&
\simeq j_{\sharp}(\gamma_{*}(F)).
%\end{aligned}
$$
Hence we obtain a canonical isomorphism
$j_{\sharp}(\gamma_{*}(F))\simeq\gamma_{*}(j_{\sharp}(F))$.
It is easily seen that the corresponding map
$\gamma_{*}(F)\To j^{*}(\gamma_{*}(j_{\sharp}(F)))
=\gamma_{*}(j^{*}j_{\sharp}(F))$
is the image by $\gamma_{*}$ of the unit
map $F\To j^{*}j_{\sharp}(F)$. This shows the isomorphism
we have constructed is the exchange morphism.
\end{proof}

\begin{cor}\label{jsharpexactAS}
For any open immersion of schemes $j:U\To S$,
the functor $j_{\sharp}:\A_{U}\To\A_{S}$ is exact.
Moreover, the induced functor
$$j_{\sharp}:\Comp(\A_{U})\To\Comp(\A_{S})$$
preserves $\AA^1$-equivalences.
\end{cor}

\begin{proof}
Using the fact $\gamma_{*}$ is exact and conservative,
and propositions \ref{oubligentil} and \ref{immouvechangeAS},
it is sufficient to prove this corollary when $\A=\sh {\nis}{\sm}$.
It is straightforward to prove exactness using
Nisnevich points. The fact $j_{\sharp}$ preserves
$\AA^1$-equivalences follows from the exactness
property and from the obvious fact it preserves
strong $\AA^1$-equivalences.
\end{proof}

\begin{cor}\label{commoublitransimmouv}
Let $j:U\To S$ be an open immersion
of schemes. For any object $M$ of $\DMue(\A_{U})$
the exchange morphism
\begin{equation} \label{eq:j_gamma_commutes}
\derL j_{\sharp}(\derR\gamma_{*}(M))\To\derR\gamma_{*}(\derL j_{\sharp}(M))
\end{equation}
is an isomorphism in $\DMtex S$.
\index{word}{exchange!isomorphism}
\end{cor}

%\begin{cor} \label{cor:supp_nearly_nis_eff}
%The triangulated premotivic category $\DMue(\A)$
% satisfies the support property.
%\end{cor}
%\begin{proof}
%According to corollary  \ref{der_oubligentil},
% the functor $\derR \gamma_*$ is conservative.
% Thus, by virtue of  the preceding corollary,
% to prove the support property in the case of $\DMue(\A)$
% it is sufficient to prove it in the case where $\A=\sh{\nis}{\sm}$.
%This follows from theorems \ref{thm:DMt_localization}
% and \ref{thm:cor1_Ayoub}.
%\end{proof}

\subsubsection{Support property (stable case)}

\begin{num} Recall from \ref{num:functoriality_abelian_spectra}
 that the premotivic adjunction $(\gamma^*,\gamma_*)$
 induces a canonical adjunction of abelian premotivic categories
 that we denote by:
$$
\tilde \gamma^*:\Spt(\sh{\nis}{\sm}) \rightleftarrows \Spt(\A_S):\tilde \gamma_*
$$
\end{num}

\begin{prop}\label{oubligentil_stable}
For any scheme $S$,
 the functor induced functor
$$
\tilde \gamma_*:\Comp\big(\Spt(\A_S)\big)
   \rightleftarrows \Comp\big(\Spt(\sh{\nis}{\sm/S})\big)
$$
preserves and detects stable $\AA^1$-equivalences.
\end{prop}
\begin{proof}
Using the equivalence between symmetric Tate spectra and non symmetric Tate
spectra, we are reduced to prove this for complexes of non symmetric
Tate spectra. Consider a non symmetric Tate spectrum $(E_n)_{n \in \NN}$
with suspension maps
$\sigma_n:E_n\{1\} \rightarrow E_{n+1}$.
The non symmetric Tate spectrum $\tilde \gamma_*(E)$ is
 equal to $\gamma_*(E_n)$ in degree $n \in \ZZ$, 
 and the suspension map is given by the composite:
$$
\MabNX S {}\{1\} \otimes_S \gamma_*(E_n)
 \rightarrow \gamma_*(\gamma^*(\MabNX S {}\{1\}) \otimes_S E_n)=\gamma_*(E_n\{1\})
 \xrightarrow{\gamma_*(\sigma_n)} E_{n+1}.
$$
Thus, propositions \ref{oubligentil} and \ref{nonsymtosymequ2}
 allows us to conclude.
\end{proof}

\begin{cor}\label{der_oubligentil_stable}
The right derived functor
$$
\derR \gamma_{*}=\gamma_*:\DMu(\A_{S})\To\DMt(S)
$$
is conservative.
\end{cor}

\begin{prop}\label{comm_oubli_ouv_stable}
Let $j:U\To X$ be an open immersion of schemes.
For any object $M$ of $\DMu(\A_{U})$, the exchange morphism
$$\derL j_{\sharp}(\derR\gamma_{*}(M))\To\derR\gamma_{*}(\derL j_{\sharp}(M))$$
is an isomorphism in \smash{$\DMt(X)$}.
\index{word}{exchange!isomorphism}
\end{prop}
\begin{proof}
From Corollary \ref{jsharpexactAS} and
 the $\Pmor$-base change formula for the open immersion $j$,
 one deduces easily that $j_\sharp$ preserves stable $\AA^1$-equivalences
 of (non symmetric) Tate spectra.
 Moreover, Proposition \ref{immouvechangeAS} shows that
  $j_\sharp \gamma_*=\gamma_*j_\sharp$ at the level of Tate spectra.
 This concludes.
\end{proof}

\begin{cor} \label{cor:supp_nearly_nis}
The triangulated premotivic category $\DMu(\A)$
 satisfies the support property.
\end{cor}
\begin{proof}
According to corollary  \ref{der_oubligentil_stable},
 the functor $\derR \gamma_*$ is conservative.
 Thus, by virtue of  the preceding proposition,
 to prove the support property in the case of $\DMu(\A)$
 it is sufficient to prove it in the case where $\A=\sh{\nis}{\sm}$.
This follows from theorems \ref{thm:DMt_localization}
 and \ref{thm:cor3_Ayoub}.
\end{proof}

\subsubsection{Localization for smooth schemes}

\begin{lm}\label{lmsmoothlocabstract}
Let $i:Z\To S$ be a closed immersion
which admits a smooth retraction $p:S\To Z$.
Then the exchange transformation
$$\derL\gamma^*\derR i_{*}\To\derR i_{*}\derL\gamma^*$$
is an isomorphism in $\DMue(\A_{S})$ (resp. $\DMu(\A_{S})$).
\index{word}{exchange!isomorphism}
\end{lm}

\begin{proof}
We first remark that for any object $C$ of $\Comp(\A_{Z})$
 (resp. $\Comp(\Spt(\A_Z))$) the canonical sequence
$$
j_\sharp (pj)^*(C) \rightarrow p^*(C)
 \rightarrow i_*(C)
$$
is a cofiber sequence in $\DMue(\A_S)$
 (resp. $\DMu(\A)_S)$).
Indeed, we can check this after applying 
the exact conservative functor $\gamma_{*}$. 
The sequence we obtain is canonically isomorphic 
through exchange transformations to
$$
j_\sharp j^* p^*(\gamma_{*} C)
 \rightarrow p^*(\gamma_{*}C)
 \rightarrow i_*i^*p^*(\gamma_{*}C)
$$
using Corollary \ref{commoublitransimmouv},
 the commutation of $\gamma_*$ with $j^*$, $p^*$ and $i_*$
 (recall it is the right adjoint of a premotivic adjunction)
 and the relation $pi=1$.
But this last sequence is a cofiber sequence in $\DMte(S)$
 (resp. $\DMt(S)$) because it satisfies the localization property
 (see \ref{thm:DMt_localization}).

Using exchange transformations,
 we obtain a morphism of distinguished triangles in $\DMe(S)$
$$
\xymatrix@R=10pt@C=14pt{
\gamma^* j_\sharp j^* p^*(C)\ar[r]\ar@{=}[d]
 & \gamma^* p^*(C) \ar[r]\ar@{=}[d]
  & \gamma^* i_*(C)\ar^/15pt/{}[r]\ar^{Ex(\gamma^*,i_*)}[d]
 &\gamma^* j_\sharp j^* p^*(C)[1] \ar@{=}[d]\\
j_\sharp j^* p^*(\gamma^* C)\ar[r]
 & p^*(\gamma^* C) \ar[r]
  & i_*(\gamma^* C)\ar^/15pt/{}[r] &j_\sharp j^* p^*(\gamma^* C)[1]
}
$$
The first two vertical arrows are isomorphisms
 as $\gamma^*$ is the left adjoint of a premotivic adjunction;
thus the morphism $Ex(\gamma^*,i_*)$ is also an isomorphism.
\end{proof}

\begin{prop}\label{smoothlocabstract}
Let $i:Z\To S$ be a closed immersion.
If $i$ admits a smooth retraction,
then \smash{$\DMue(\A)$}
satisfies \locx i.
\end{prop}

\begin{proof}
This follows from Proposition \ref{prop:localization&exchange}
and the preceding lemma.
\end{proof}

\begin{cor}\label{smoothloc}
Let $S$ be a scheme.
Then the premotivic category \smash{$\DMue(\A)$} (resp. $\DMu(\A)$)
satisfies localization with respect to any closed immersion
between smooth $S$-schemes.
\end{cor}
\begin{proof}
Let $i:Z\To X$ be closed immersion between
smooth $S$-schemes.
We want to prove that $\DMue(\A)$ (resp. $\DMu(\A)$)
satisfies localization with respect to $i$.
According to \ref{cor:premotivic&i_*_bis},
 it is sufficient to prove that for any smooth $S$-scheme $S$,
 the canonical map
$$
M_S(X/X-X_Z) \rightarrow i_*M_Z(X_Z)
$$
is an isomorphism where we use the notation of \emph{loc. cit.}
 and $M(.,\A)$ denotes the geometric sections of $\DMue(\A)$
  (resp. $\DMu(\A)$).
But the premotivic triangulated category $\DMu(\A)$
 (resp. $\DMue(\A)$) satisfies the Nisnevich separation property
 and the $\sm$-base change property.
 Thus, we can argue locally in $S$ for the Nisnevich topology.
 Thus, the statement is reduced to the preceding proposition
 as $i$ admits locally for the Nisnevich topology a smooth retraction
 (see for example \cite[4.5.11]{Deg7}).
\end{proof}

\section{Basic homotopy commutative algebra} \label{sec:modules_et_anneaux}

\subsection{Rings}

\begin{df}
A symmetric monoidal model category $\V$ satisfies the
\emph{monoid axiom}
\index{word}{monoid axiom}
 if, for any trivial cofibration $A\To B$
and any object $X$, the smallest class of maps of $\V$ which contains
the map $X\otimes A\to X\otimes B$ and is stable by pushouts and
transfinite compositions is contained in the class of weak equivalences.
\end{df}

\begin{paragr}
Let $\V$ be a symmetric monoidal category.
We denote by $\Alg(\V)$ the category of monoids
\index{word}{monoid}
 in $\V$.
If $\V$ has small colimits, the forgetful functor
$$U:\Alg(\V)\To \V$$ %%\quad , \qquad (R,\mu,\eta)\mapsto R$$
has a left adjoint
$$F:\V\To\Alg(\V)\, .$$
\end{paragr}

\begin{thm}\label{cmfmonoids}
Let $\V$ a symmetric monoidal combinatorial model category
which satisfies the monoid axiom. The category of monoids $\Alg(\V)$
is endowed with the structure of a combinatorial model category
whose weak equivalences 
\index{word}{equivalence!weak equivalence of monoids}
(resp. fibrations)
\index{word}{fibration!of monoids}
 are the morphisms of commutative monoids
which are weak equivalences (resp. fibrations) in $\V$.
In particular, the forgetful functor $U:\Alg(\V)\To\V$ is a right Quillen functor.
Moreover, if the unit object of $\V$ is cofibrant, then
any cofibrant object of $\Alg(\V)$ is cofibrant as an object of $\V$.
\end{thm}

\begin{proof}
This is very a particular case of the third assertion of \cite[Theorem 4.1]{SS}
(the fact that $\Alg(\V)$ is combinatorial whenever $\V$ is so
comes for instance from \cite[Proposition 2.3]{beke1}).
\end{proof}

\begin{df}\label{def:stronglyQQlinear}
A symmetric monoidal model category $\V$ is \emph{strongly $\QQ$-linear}
\index{word}{linear!strongly $\QQ$-linear}
if the underlying category of $\V$ is additive and $\QQ$-linear
(i.e. all the objects of $\V$ are uniquely divisible).
\end{df}

\begin{rem}
If $\V$ is a strongly $\QQ$-linear stable model category, then it is
$\QQ$-linear in the sense of \ref{defQlinear}.
\end{rem}

\begin{lm}\label{trivquotients}
Let $\V$ be a strongly $\QQ$-linear model category, $G$ a finite group,
and $u:E\To F$ an equivariant morphism of representations of $G$ in $\V$.
Then, if $u$ is a cofibration in $\V$, so is the induced map
$E_G\To F_G$ (where the subscript $G$ denotes the coinvariants under the action
of the group $G$).
\end{lm}

\begin{proof}
The map $E_G\To F_G$ is easily seen to be a direct factor (retract) of the
cofibration $E\to F$.
\end{proof}

\begin{paragr}
If $\V$ is a symmetric monoidal category,
we denote by $\Calg(\V)$ the category of commutative monoids
\index{word}{monoid!commutative monoid}
 in $\V$.
If $\V$ has small colimits, the forgetful functor
$$U:\Calg(\V)\To \V$$ %%\quad , \qquad (R,\mu,\eta)\mapsto R$$
has a left adjoint
$$F:\V\To\Calg(\V)\, .$$
\end{paragr}

\begin{thm}\label{cmfQlinearcomm}
Let $\V$ a symmetric monoidal combinatorial model category.
Assume that $\V$ is left proper and tractable,\index{word}{tractable}
 satisfies the monoid axiom, and is strongly $\QQ$-linear. 
Then the category of commutative
monoids $\Calg(\V)$ is endowed with the structure of a combinatorial model category
whose weak equivalences\index{word}{equivalence!weak equivalence of commutative monoids}
 (resp. fibrations)
\index{word}{fibration!of commutative monoids}
 are the morphisms of commutative monoids
which are weak equivalences (resp. fibrations) in $\V$.
In particular, the forgetful functor $U:\Calg(\V)\To\V$ is a right Quillen functor.

If moreover the unit object of $\V$ is cofibrant, then
any cofibrant object of $\Calg(\V)$ is cofibrant as an object of $\V$.
\end{thm}

\begin{proof}
We will observe first that $\V$ is freely powered
in the sense of \cite[Definition 4.5.4.2]{DAG3}. Therefore, the existence of this model category structure will follow from
a general result of Lurie \cite[Proposition 4.5.4.6]{DAG3}.
For this, it is sufficient to check that a $G$-equivariant
map $f:A\To B$ in $\V$ which is a trivial cofibration when we forget
the $G$-action has the left lifting property with respect to
any $G$-equivariant map $p:X\To Y$ which is a fibration in $\V$
(after forgetting the $G$-action). In other words, we have to check
that the map induced by $f$ and $p$ in $\V$
$$\Hom_V(B,X)\To\Hom_\V(A,X)\times_{\Hom_\V(A,X)}\Hom_\V(Y,B)$$
will induce a surjective map after we apply the $G$-invariants
functor (we let the reader construct a natural $G$-action on
$\Hom_V(B,X)$, the $G$-invariants of which gives the $\QQ$-vector
space of $G$-equivariant maps from $B$ to $X$).
Since $G$ is a finite group, the $G$-invariant subspace functor
is exact, hence this is obvious. This proves the
first assertion. The second assertion of the theorem is true by definition.

The last assertion is proved by a careful analysis of pushouts by free maps
in $\Calg(\V)$ as follows.
For two cofibrations $u:A\To B$ and $v:C\To D$ in $\V$, write $u\wedge v$
for the map
$$u\wedge v: A\otimes D\amalg_{A\otimes C} B\otimes C\To B\otimes D$$
(which is a cofibration by definition of monoidal model categories).
By iterating this construction, we get, for a cofibration $u:A\to B$ in $\V$, a
cofibration
$$\wedge^n(u)=
\underset{\text{$n$ times}}{\underbrace{u\wedge \dots \wedge u}}
:\square^n(u)\To B^{\otimes n}\, .$$
Note that the symmetric group $\mathfrak{S}_n$ acts naturally on $B^{\otimes n}$
and $\square^n(u)$. We define
$$\mathit{Sym}^n(B)=(B^{\otimes n})_{\mathfrak{S}_n}
\quad\text{and}\quad
\mathit{Sym}^n(B,A)=\square^n(u)_{\mathfrak{S}_n}\, .$$
By virtue of Lemma \ref{trivquotients}, we get a cofibration of $\V$:
$$\sigma^n(u):\mathit{Sym}^n(B,A)\To \mathit{Sym}^n(B)\, .$$
Consider now the free map
$F(u):F(A)\to F(B)$ can be filtered by $F(A)$-modules as follows.
Define $D_0=F(A)$. As $A=\mathit{Sym}^1(B,A)$, we have
a natural morphism $F(A)\otimes \mathit{Sym}^1(B,A)\To F(A)$.
The objects $D_n$ are then defined by induction with the pushouts
below.
$$\xymatrix{
F(A)\otimes \mathit{Sym}^n(B,A)\ar[rr]^{1_{F(A)\otimes \sigma^n(u)}}\ar[d]&&
F(A)\otimes \mathit{Sym}^n(B)\ar[d]\\
D_{n-1}\ar[rr]&&D_n
}$$
We get natural maps $D_n\To F(B)$ which induce an isomorphism
$$\varinjlim_{n\geq 0} D_n\simeq F(B)$$
in such a way that the morphism $F(u)$ correspond to the canonical map
$$F(A)=D_0\To \varinjlim_{n\geq 0} D_n\, .$$
Hence, if $F(A)$ is cofibrant, all the maps $D_{n-1}\to D_n$
are cofibrations, so that the map $F(A)\To F(B)$ is
a cofibration in $\V$. In the particular case where $A$ is the initial
object of $\V$, we see that for any cofibrant object $B$ of $\V$, the
free commutative monoid $F(B)$ is cofibrant as an object of $\V$
(because the initial object of $\Calg(\V)$ is the unit object of $\V$).
This also implies that, if $u$ is a cofibration between cofibrant
objects, the map $F(u)$ is a cofibration in $\V$.

This description of $F(u)$ also allows to compute the pushouts of $F(u)$ in
$\Calg(\V)$ in $\V$ as follows. Consider a pushout
$$\xymatrix{
F(A)\ar[r]^{F(u)}\ar[d]&F(B)\ar[d]\\
R\ar[r]_v& S
}$$
in $\Calg(\V)$.
For $n\geq 0$, define $R_n$ by the pushouts of $\V$:
$$\xymatrix{
F(A)\ar[r]\ar[d]& D_n\ar[d]\\
R\ar[r]& R_n
}$$
We then have an isomorphism
$$\varinjlim_{n\geq 0}R_n\simeq S\, .$$
In particular, if $u$ is a cofibration between cofibrant
objects, the morphism of commutative monoids $v:R\to S$
is then a cofibration in $\V$. As the forgetful functor $U$
preserves filtered colimits, we conclude easily from there
(with the small object argument \cite[Theorem 2.1.14]{Hovey})
that any cofibration of $\Calg(\V)$ is a cofibration of $\V$.
Using again that the unit object of $\V$ is cofibrant in $\V$
(i.e. that the initial object of $\Calg(\V)$ is cofibrant in $\V$)
this proves the last assertion of the theorem.
\end{proof}

\begin{cor}\label{strictifyloccommalg}
Let $\V$ a symmetric monoidal combinatorial model category.
Assume that $\V$ is left proper and tractable, satisfies the monoid axiom, and
is strongly $\QQ$-linear. Consider a small set $H$ of maps of $\V$, and
denote by $L_H\V$ the left Bousfield localization of $\V$ by $H$; see
\cite[Theorem 4.7]{Bar}. Define the class of $H$-equivalences in $\ho(\V)$
to be the class of maps which become invertible in $\ho(L_H\V)$.
If $H$-equivalences are stable by (derived) tensor product in $\ho(\V)$, then
$L_H\V$ is a symmetric monoidal combinatorial model category
(which is again left proper and tractable, satisfies the monoid axiom, and
is strongly $\QQ$-linear).

In particular, under these assumptions,
there exists a morphism of commutative monoids $\unit \To R$
in $\V$ which is a weak equivalence of $L_H\V$, with $R$
a cofibrant and fibrant object of $L_H\V$.
\end{cor}

\begin{proof}
The first assertion is a triviality. The last assertion follows
immediately: the map $\unit\To R$ is simply obtained as
a fibrant replacement of $\unit$ in the model category $\Calg(L_H\V)$
obtained from Theorem \ref{cmfQlinearcomm} applied to $L_H\V$.
\end{proof}

\begin{paragr}\label{fibredmodcatcomm0}
Consider now a category $\site$,
%%endowed with the admissible class $\mathit{Iso}$
%%of isomorphisms of $\C$,
as well as a closed symmetric monoidal
bifibred category $\M$ over $\site$.
We shall also assume that the fibers of $\M$ admit limits and colimits.

Then the categories $\Alg(\M(X))$ (resp. $\Calg(\M(X))$)
define a bifibred category
over $\site$ as follows. Given a morphism $f:X\To Y$, the functor
$$f^*:\M(Y)\To\M(X)$$
is symmetric monoidal, so that it preserves monoids
(resp. commutative monoids) as well as morphisms between them.
It thus induces a functor
\begin{equation}\label{CommPfibred1}
\begin{aligned}
& f^*:\Alg(\M(Y))\To\Alg(\M(X))\\
(\text{resp. } & f^*:\Calg(\M(Y))\To\Calg(\M(X))\, ).
\end{aligned}
\end{equation}
As $f^*:\M(Y)\To\M(X)$ is symmetric monoidal, its right adjoint $f_*$ is lax monoidal:
there is a natural morphism
\begin{equation}\label{CommPfibred3}
\unit_Y\To f_*(\unit_X)=f_*\, f^*(\unit_Y)\, ,
\end{equation}
and, for any objects $A$ and $B$ of $\M(X)$, there is
a natural morphism
\begin{equation}\label{CommPfibred4}
f_*(A)\otimes_Y f_*(B)\To f_*(A\otimes_X B)
\end{equation}
which corresponds by adjunction to the map
$$f^*(f_*(A)\otimes_Y f_*(B))\simeq f^*\, f_*(A)\otimes f^*\, f_*(B)\To A\otimes B\, .$$
Hence the functor $f_*$ preserves also monoids (resp. commutative monoids)
as well as morphisms between them, so that we get a functor
\begin{equation}\label{CommPfibred5}
\begin{aligned}
& f_*:\Alg(\M(X))\To\Alg(\M(Y))\\
(\text{resp. } & f_*:\Calg(\M(X))\To\Calg(\M(Y))\, ).
\end{aligned}
\end{equation}
By construction, the functor $f^*$ of \eqref{CommPfibred1}
is a left adjoint ot the functor $f_*$ of \eqref{CommPfibred5}.
These constructions extend to morphisms of $\site$-diagrams in a similar
way.
\end{paragr}

\begin{prop}\label{Qlinearhomotopyalgebras}
Let $\M$ be a symmetric monoidal combinatorial fibred model
category over $\site$. Assume that, for any object $X$ of $\site$, the
model category $\M(X)$ satisfies the monoid axiom
(resp. is left proper and tractable, satisfies the monoid axiom, and
is strongly $\QQ$-linear).
\begin{itemize}
\item[(a)] For any object $X$ of $\site$, the category
$\Alg(\M)(X)$ (resp. $\Calg(\M)(X)$) of monoids (resp. of
commutative monoids) in $\M(X)$ is a combinatorial
model category structure whose weak equivalences
(resp. fibrations) are the morphisms of commutative monoids
which are weak equivalences (resp. fibrations) in $\M(X)$.
This turns $\Alg(\M)$ (resp. $\Calg(\M)$) into a combinatorial fibred
model category over $\site$.
\item[(b)] For any morphism of $\site$-diagrams
$\varphi:(\X,I) \To (Y,J)$, the adjunction
$$\varphi^*:\Alg(\M)(\Y,J)\rightleftarrows\Alg(\M)(\X,I):\varphi_*$$
$$(\text{resp. } \varphi^*:\Calg(\M)(\Y,J)\rightleftarrows\Calg(\M)(\X,I):\varphi_*)$$
is a Quillen adjunction (where the categories of monoids
$\Alg(\M)(\X,I)$ (resp. of commutative monoids $\Calg(\M)(\X,I)$)
are endowed with the injective model
category structure obtained from Proposition \ref{injdiamodcat}
applied to $\Alg(\M)$ (resp. to $\Calg(\M)$).
\item[(d)] If moreover, for any object $X$ of $\site$, the unit $\unit_X$
is cofibrant in $\M(X)$, then, for morphism of $\site$-diagrams
$\varphi:(\X,I) \To (Y,J)$, the square
\begin{equation}\label{derivedCommPfibred21}\begin{split}
\xymatrix{
\ho(\Alg(\M))(\Y,J)\ar[r]^{\derL \varphi^*}\ar[d]_U&\ho(\Alg(\M))(\X,I)\ar[d]^U\\
\ho(\M)(\Y,J)\ar[r]^{\derL \varphi^*}&\ho(\M)(\X,I)
}\end{split}
\end{equation}
is essentially commutative. Similarly, in the respective case, the square
\begin{equation}\label{derivedCommPfibred22}\begin{split}
\xymatrix{
\ho(\Calg(\M))(\Y,J)\ar[r]^{\derL \varphi^*}\ar[d]_U&\ho(\Calg(\M))(\X,I)\ar[d]^U\\
\ho(\M)(\Y,J)\ar[r]^{\derL \varphi^*}&\ho(\M)(\X,I)
}\end{split}
\end{equation}
is essentially commutative.
\end{itemize}
\end{prop}

\begin{proof}
Assertion (a) is an immediate consequence of
Theorem \ref{cmfmonoids} (resp. of Theorem \ref{cmfQlinearcomm}),
and assertion (b) is a particular case of Proposition \ref{basicfunctdiag}
(beware that the injective model category structure on $\Calg(\M)(\X,I)$
does not necessarily coincide with the model category structure
given by Theorem \ref{cmfmonoids} (resp. of Theorem \ref{cmfQlinearcomm})
applied to the injective model structure on $\M(\X,I)$).
For assertion (d), we see by the second assertion of Proposition \ref{projdiamodcat}
that it is sufficient to prove it when $\varphi:X\to Y$ is simply a morphism of $\site$.
In this case, by construction of the total left derived functor
of a left Quillen functor, this follows from the fact that $\varphi^*$
commutes with the forgetful functor and
from the fact that, by virtue of the last assertion of
Theorem \ref{cmfmonoids} (resp. of Theorem \ref{cmfQlinearcomm}),
the forgetful functor $U$ preserves weak equivalences and cofibrant objects.
\end{proof}

\begin{rem}\label{adjuntionsmonoids}
The main application of the preceding corollary will come from
assertion (d): it says that, given a monoid (resp. a commutative monoid)
$R$ in $\M(Y)$ and a morphism $f:X\To Y$, the image of $R$ by the functor
$$\derL f^*:\ho(\M)(Y)\To\ho(\M)(X)$$
is canonically endowed with a structure of monoid (resp. of commutative monoid)
in the strongest sense possible. Under the assumptions of
assertion (c) of Proposition \ref{Qlinearhomotopyalgebras},
we shall often make the abuse of
saying that $\derL f^*(R)$ is a monoid (resp. a commutative monoid)
in $\M(X)$ without refereeing explicitly to the model category structure on
$\Alg(\M)(X)$ (resp. on $\Calg(\M)(X)$).
Similarly, for any monoid (resp. commutative monoid) $R$ in $\M(X)$,
$\derR f_*(R)$ will be canonically endowed with a structure
of a monoid (resp. a commutative monoid) in $\M(Y)$.
In particular, for any monoid (resp. commutative monoid) $R$ in $\M(Y)$,
the adjunction map
$$R\To \derR f_*\, \derL f^*(R)$$
is a morphism of monoids (i.e. is a map in the homotopy category
$\ho(\Alg(\M))(X)$ (resp. $\ho(\Calg(\M))(X)$)),
and, for any monoid (resp. commutative monoid) $R$ in $\M(X)$,
the adjunction map
$$\derL f^*\, \derR f_*(R)\To R$$
is a morphism of monoids (i.e. is a map in the homotopy category $\ho(\Alg(\M))(Y)$
(resp. $\ho(\Calg(\M))(Y)$)).
\end{rem}

\begin{rem}
In order to get a good homotopy theory of
commutative monoids wihout the strongly $\QQ$-linear assumption,
we should replace commutative monoids by $E_\infty$-algebras
\index{word}{algebra!einfinityalgebra@$E_\infty$-algebra}
(i.e. objects endowed with a structure of commutative
monoid up to a bunch of coherent homotopies).
More generally, we should prove the analog of Theorem \ref{cmfmonoids}
and of Theorem \ref{cmfQlinearcomm}
by replacing $\Alg(\V)$ by the category of algebras of some 
`well-behaved' operad, and then get as a consequence the analog of Proposition \ref{Qlinearhomotopyalgebras}.
All this is a consequence of the general constructions and results
of \cite{spit,BerMoe1,BerMoe2}.

However, in the case we are interested in the homotopy theory of
commutative monoids in some category of spectra $\V$, it seems
that some version of Shipley's \emph{positive stable model structure}
\index{word}{model structure!positive stable model structure}
(\textit{cf.} \cite[Proposition 3.1]{shipley}) would provide a good
model category for commutative monoids, which,
 by Lurie's strictification theorem
\index{word}{strictification theorem}
 \cite[Theorem 4.5.4.7]{DAG3}, would be equivalent
to the homotopy theory of $E_\infty$-algebras in $\V$.
This kind of technics is now available in the context
of stable homotopy theory of schemes, which provides a good
setting to speak of motivic commutative ring spectra; see \cite{hornbostel,gulet1,gulet2,positive1}.
Therefore, Theorem \ref{cmfQlinearcomm} and Proposition \ref{Qlinearhomotopyalgebras}
are in fact true in $\SH$ for genuine commutative monoids without any $\QQ$-linearity
assumption.
\end{rem}

\subsection{Modules}

\begin{paragr}
Given a monoid $R$ in a symmetric monoidal
category $\V$, we shall write $\Mod R(\V)$ for the category of (left) $R$-modules.
\index{word}{modules!over a monoid}
The forgetful functor
$$U:\Mod R(\V) \To \V$$
is a left adjoint to the free $R$-module functor
$$R\otimes (-):\V\To\Mod R(\V) \, .$$
If $\V$ has enough small colimits, and if $R$ is a commutative monoid,
the category $\Mod R(\V)$ is endowed with a unique symmetric monoidal structure
such that the functor $R\otimes (-)$ is naturally symmetric monoidal. We shall
denote by $\otimes_R$ the tensor product of $\Mod R(\V)$.
\end{paragr}

\begin{thm}\label{abstractcmfmodules}
Let $\V$ be a combinatorial symmetric model category which
satisfies the monoid axiom.\index{word}{monoid axiom}
\begin{itemize}
\item[(i)] For any monoid $R$ in $\V$, the category of right (resp. left) $R$-modules
is a combinatorial model category with weak equivalences
\index{word}{equivalence!weak equivalence of modules}
 (resp. fibrations)
\index{word}{fibration!of modules}
the morphisms of $R$-modules which are weak equivalences (resp. fibrations)
in $\V$.
\item[(ii)] For any commutative monoid $R$ in $\V$, the model category of $R$-modules
given by (i) is a combinatorial symmetric monoidal model category which satisfies the
monoid axiom.
\end{itemize}
\end{thm}

\begin{proof}
Assertions (i) and (ii) are particular cases of the first two assertions of \cite[Theorem 4.1]{SS}.
\end{proof}

%% \begin{thm}\label{modulesbasechange}
%% Let $\V$ a combinatorial symmetric model category which
%% satisfies the monoid axiom, and $R\To S$ a weak equivalence of monoids $R\To S$.
%% If, for any cofibrant right (resp. left) $R$-module $M$, the functor $M \otimes_R (-)$
%% (resp. $(-)\otimes_R M$) sends weak equivalences of left (resp. right)
%% $R$-modules to weak equivalences of $\V$, then the functor
%% $M\mapsto S\otimes_R M$ (resp. $M\mapsto M\otimes_R S$)
%% is a left Quillen equivalence from the category of left (resp. right)
%% $R$-modules to the category of left (resp. right) $S$-modules.
%% \end{thm}
%% 
%% \begin{proof}
%% This is a special case of \cite[Theorem 4.3]{SS}.
%% \end{proof}

\begin{df}\label{defperfectsymmonidcmf}
A symmetric monoidal model category $\V$ is \emph{perfect}
\index{word}{perfect}
if it has the following properties.
\begin{itemize}
\item[(a)] $\V$ is combinatorial and tractable (\ref{deftractable});
\item[(b)] $\V$ satisfies the monoid axiom;
\item[(c)] For any weak equivalence of monoids $R\To S$,
the functor $M\mapsto S\otimes_R M$
is a left Quillen equivalence from the category of left
$R$-modules to the category of left $S$-modules.
\item[(d)] weak equivalences are stable by small sums in $\V$.
\end{itemize}
\end{df}

\begin{rem}\label{perfectarenice1}
If $\V$ is a perfect symmetric monoidal model category, then, for any
commutative monoid $R$, the symmetric monoidal model category
of $R$-modules in $\V$ given by Theorem~\ref{abstractcmfmodules}~(ii)
is also perfect: condition (c) is quite obvious, and condition (d) comes
from the fact that the forgetful functor $U:\Mod R\To \V$
commutes with small sums, while it preserves and detects weak equivalences.
Note that condition (d) implies that the functor
$U:\ho(\Mod R)\To \ho(\V)$ preserves small sums.
\end{rem}

\begin{rem}\label{perfectarenice2}
If $\V$ is a stable symmetric monoidal model category which satisfies the monoid axiom, then
for any monoid $R$ of $\V$, the model category of (left) $R$-modules
given by Theorem \ref{abstractcmfmodules} is stable as well: the suspension
functor of $\ho(\Mod R)$ is given by the derived tensor product by the
$R$-bimodule $R[1]$, which is clearly invertible with inverse $R[-1]$.
\end{rem}

In this work,
 a basic example of perfect model categories are those coming
 from stable $\AA^1$-derived premotivic categories
  (cf Def. \ref{df:triangulated_premotives}):
\begin{prop}\label{sheavesTatespectranicederivedtensorproduct}
Let $t$ be an admissible topology. Then, for any scheme $S$ in $\sch$,
the symmetric monoidal model structure on
$\Comp(\Spt(\sh {t} {\Pmor/S}))$
underlying the triangulated category $\DMu\left(\sh {t} {\Pmor/S}\right)$
  is perfect.
\end{prop}

\begin{proof}
The generating family of $\sh {t} {\Pmor/S}$ is flat
in the sense of \cite[3.1]{CD1}, so that, by virtue of
\cite[prop. 7.22 and cor. 7.24]{CD1}, the assumptions
of Proposition \ref{modulesbasechangestable} are fulfilled.
\end{proof}

\begin{prop}\label{Modcompactgenerators}
Let $\V$ be a stable perfect symmetric monoidal model category.
Assume furthermore that $\ho(\V)$ admits a small family $\mathcal{G}$ of
compact generators (as a triangulated category).
For any monoid $R$ in $\V$, the triangulated category $\ho(\Mod R(\V))$
admits the set $\{R\otimes^\derL E\ | \ E\in \mathcal{G}\}$
as a family of compact generators.
\end{prop}

\begin{proof}
We have a derived adjunction
$$R\otimes^\derL(-):\ho(\V)\rightleftarrows\ho(\Mod R(\V)):U\, .$$
As the functor $U$ preserves small sums the functor $R\otimes^\derL(-)$
preserves compact objects. But $U$ is also
conservative, so that $\{R\otimes^\derL E\ | \ E\in \mathcal{G}\}$
is a family of compact generators of $\ho(\Mod R(\V))$.
\end{proof}

\begin{rem}\label{problemderivedtensor}
If $\V$ is a combinatorial symmetric model category which
satisfies the monoid axiom, then there are two ways to derive the tensor product.
The first one consists in deriving the left Quillen bifunctor $(-)\otimes(-)$, which
gives the usual derived tensor product
$$(-)\otimes^\derL (-):\ho(\V)\times\ho(\V)\To\ho(\V)\, .$$
Remember that, by construction, $A\otimes^\derL B=A'\otimes B'$,
where $A'$ and $B'$ are cofibrant replacements of $A$ and $B$
respectively. On the other hand, the monoid axiom gives that, for
any object $A$ of $\V$, the functor $A\otimes (-)$ preserves
weak equivalences between cofibrant objects, which implies that
it has also a total left derived functor
$$A\otimes^\derL(-):\ho(\V)\To\ho(\V)\, .$$
Despite the fact we have adopted very similar (not to say identical) notations
for these two derived functor, there is no reason they would coincide
in general: by construction, the second one is defined by
$A\otimes^\derL B=A\otimes B'$, where $B'$ is some cofibrant
replacement of $B$. However, they coincide quite often in practice
(e.g. for simplicial sets, for the good reason that all of them are cofibrant,
or for symmetric $S^1$-spectra, or for complexes of
quasi-coherent $\mathcal{O}_X$-modules over a quasi-compact and quasi-separated
scheme $X$).
\end{rem}

\begin{prop}\label{modulesbasechangestable}
Let $\V$ be a stable combinatorial symmetric monoidal model category
which satisfies the monoid axiom.\index{word}{monoid axiom}
Assume furthermore that, for any cofibrant object $A$ of $\V$,
the functor $A\otimes(-)$ preserve weak equivalences (in other words,
that the two ways to derive the tensor product explained in Remark
\ref{problemderivedtensor} coincide), and that weak
equivalences are stable by small sums in $\V$.
Then the symmetric monoidal model category $\V$ is perfect.
\end{prop}

\begin{proof}
We just have to check condition (c) of Definition~\ref{defperfectsymmonidcmf}.
Consider a weak equivalence of monoids $R\To S$.
We then get a derived adjunction
$$S\otimes^\derL_R(-):\ho(\Mod R(\V))\rightleftarrows
\ho(\Mod S(\V)):U\, ,$$
where $S\otimes^\derL_R(-)$ is the left derived functor
of the functor $M\mapsto S\otimes_RM$.
We have to prove that, for any left $R$-module $M$, the map
$$M\To S\otimes^\derL_R M$$
is an isomorphism in $\ho(\V)$.
As this is a morphism of triangulated functors which commutes with sums,
and as $\ho(\Mod R(\V))$ is well generated in the sense of Neeman~\cite{Nee1}
(as the localization of a stable combinatorial model category),
it is sufficient to check this when $M$ runs over a small
family of generators of $\ho(\Mod R(\V))$. Let us chose
is a small family of generators  $\mathcal{G}$ of $\ho(\V)$.
As the forgetful functor from $\ho(\Mod R(\V))$ to $\ho(\V)$ is conservative,
we see that $\{R\otimes^\derL E\ | \ E\in \mathcal{G}\}$ is
a small generating family of $\ho(\Mod R(\V))$.
We are thus reduced to prove that the map
$$R\otimes^\derL E\To S\otimes^\derL_R(R\otimes^\derL E)
\simeq S\otimes^\derL E$$
is an isomorphism for any object $E$ in $\mathcal{G}$.
For this, we can assume that $E$ is cofibrant, and this follows then from the fact
that the functor $(-)\otimes E$ preserves weak equivalences by
assumption.
\end{proof}

\begin{paragr}\label{PfibredRMod}
Let $\site$ be a category endowed with an admissible class of morphisms $\Pmor$,
and $\M$ a cocomplete symmetric monoidal $\Pmor$-fibred category.
%%which admits limits and colimits.
Consider a monoid $R$ in the symmetric
monoidal category $\M(1_\site,\site)$ (i.e. a section of the fibred category $\Alg(\M)$
over $\site$). In other words, $R$ consists of the data
of a monoid $R_X$ for each object $X$ of $\site$, and of a morphism of monoids
$a_f:f^*(R_Y)\To R_X$ for each map $f:X\To Y$ in $\site$, subject to
coherence relations; see \ref{defevaluation}.

For an object $X$ of $\site$, we shall write $\Mod R(X)$ for the
category of (left) $R_X$-modules in $\M(X)$, i.e.
$$\Mod R(X)=\Mod{R_X}(\M(X))\, .$$
This defines a fibred category $\Mod R$ over $\site$ as follows.
%(to avoid possible confusions, we shall write somtimes
%with an index $R$ the operations in $\Mod R$).

For a morphism $f:X\To Y$, the inverse image functor
\begin{equation}\label{PfibredRMod1}
f^*:\Mod R(Y)\To\Mod R(X)
\end{equation}
is defined by
\begin{equation}\label{PfibredRMod2}
%%f^*(M)=
M\mapsto R_X\otimes_{f^*(R_Y)}f^*(M)
\end{equation}
(where, on the right-hand side, $f^*$ stands for the
inverse image functor in $\M$). The functor \eqref{PfibredRMod1}
has a right adjoint
\begin{equation}\label{PfibredRMod3}
f_*:\Mod R(X)\To\Mod R(Y)
\end{equation}
which is simply the functor induced by $f_*:\M(X)\To\M(Y)$ (as the latter
sends $R_X$-modules to $f_*(R_X)$-modules, which are themselves
$R_Y$-modules via the map $a_f$).

If the map $f$ is a $\Pmor$-morphism, then, for any
$R_X$-module $M$, the object $f_\sharp(M)$ has a natural structure of
$R_Y$-module: using the map $a_f$, $M$ has a natural
structure of $f^*(R_Y)$-module
$$f^*(R_Y)\otimes_X M\To M\, ,$$
and applying $f_\sharp$, we get by the $\Pmor$-projection formula \eqref{basicprojformula}
a morphism
$$R_Y\otimes f_\sharp(M)\simeq f_\sharp(f^*(R_Y)\otimes M)\To f_\sharp(M)$$
which defines a natural $R_Y$-module structure on $f_\sharp(M)$.
For a $\Pmor$-morphism $f:X\To Y$, we define a functor
\begin{equation}\label{PfibredRMod4}
f_\sharp:\Mod R(X)\To\Mod R(Y)
\end{equation}
as the functor induced by $f_\sharp:\M(X)\To\M(Y)$.
Note that the functor \eqref{PfibredRMod4} is a left adjoint to the functor
\eqref{PfibredRMod1} whenever the map $a_f:f^*(R_Y)\To R_X$ is an
isomorphism in $\M(X)$.

We shall say that $R$ is a \emph{cartesian monoid in $\M$ over $\site$}
\index{word}{monoid!cartesian}
 if
$R$ is a monoid of $\M(1_\C,\C)$ such that
all the structural maps $f^*(R_Y)\To R_X$ are isomorphisms
(i.e. if $R$ is a cartesian section of the fibred category $\Alg(\M)$ over $\site$)

If $R$ is a cartesian monoid in $\M$ over $\site$, then $\Mod R$ is
a $\Pmor$-fibred category over $\site$: to see this,
it remains to prove that, for any pullback square of $\site$
$$\xymatrix{
X'\ar[r]^g\ar[d]_{f'}&X\ar[d]^f\\
Y'\ar[r]_{h}&Y
}$$
in which $f$ is a $\Pmor$-morphism, and for any $R_X$-module $M$,
the base change map
$$f'_\sharp \, g^*(M)\To h^*\, f_\sharp(M)$$
is an isomorphism, which follows immediately from the analogous formula for $\M$.
%%and from the projection formula \eqref{basicprojformula}.

Similarly, we see that whenever $R$ is a commutative monoid of $\M(1_\site,\site)$
(i.e. $R_X$ is a commutative monoid in $\M(X)$ for all $X$ in $\site$), then
$\Mod R$ is a symmetric monoidal $\Pmor$-fibred category.
\index{word}{modules!over a homotopy cartesian commutative monoid}
\end{paragr}

\begin{prop}\label{cmfPfibredRMod}
Let $\M$ be a combinatorial symmetric monoidal $\Pmor$-fibred model category
over $\site$ which satisfies the monoid axiom,
\index{word}{monoid}
 and $R$ a monoid in $\M(1_\site,\site)$ (resp. a cartesian monoid
in $\M$ over $\site$). Then \ref{abstractcmfmodules}~(i) applied termwise turns
$\Mod R$ into a combinatorial fibred model category
(resp. a combinatorial $\Pmor$-fibred model category).

If moreover $R$ is commutative, then
$\Mod R$ is a combinatorial symmetric monoidal fibred model category
(resp. a combinatorial symmetric monoidal $\Pmor$-fibred model category).
\end{prop}

\begin{proof}
Choose, for each object $X$ of $\site$, two small sets of maps $I_X$ and $J_X$
which generate the class of cofibrations and the class of trivial cofibrations
in $\M(X)$ respectively. Then $R_X\otimes_X I_X$ and
$R_X\otimes_X J_X$ generate the class of cofibrations and the class of trivial cofibrations
in $\Mod R(X)$ respectively. For a map $f:X\To Y$ in $\site$, we see from formula
\eqref{PfibredRMod2} that the functor \eqref{PfibredRMod1} sends these
generating cofibrations and trivial cofibrations to cofibrations and trivial cofibrations respectively,
from which we deduce that the functor \eqref{PfibredRMod1} is a left Quillen
functor. In the respective case,
if $f$ is a $\Pmor$-morphism, then we deduce similarly from the projection formula
\eqref{basicprojformula} in $\M$ that the functor \eqref{PfibredRMod4}
sends generating cofibrations and trivial cofibrations to cofibrations and trivial cofibrations
respectively. The last assertion follows easily by applying \ref{abstractcmfmodules}~(ii) termwise.
\end{proof}

\begin{df}
Let $\M$ be a symmetric monoidal $\Pmor$-fibred model category over $\site$.
A \emph{homotopy cartesian monoid}
\index{word}{monoid!homotopy cartesian}
$R$ in $\M$ will be a homotopy cartesian section
 of $\Alg(\M)$.
\end{df}

\begin{prop}\label{abstractmotivicmodules}
Let $\M$ be a perfect\index{word}{perfect}
 symmetric monoidal $\Pmor$-fibred model category over $\site$, and
consider a homotopy cartesian monoid $R$ in $\M$ over $\site$.

Then $\ho(\Mod R)$ is a $\Pmor$-fibred category over $\site$, and
$$R\otimes^\derL(-):\ho(\M)\To\ho(\Mod R)$$
is a morphism of $\Pmor$-fibred categories.
In the case where $R$ is commutative, $\ho(\Mod R)$ is even a
symmetric monoidal $\Pmor$-fibred category.

Moreover, for any weak equivalence between homotopy cartesian monoids
$R\To S$ over $\site$, the Quillen morphism
$$S\otimes_R(-): \Mod R\To \Mod S$$
induces an equivalence of $\Pmor$-fibred categories over $\site$
$$S\otimes^\derL_R(-):\ho(\Mod R)\To\ho(\Mod S)\, .$$
\end{prop}

\begin{proof}
It is sufficient to prove these assertions by restricting
everything over $\site/S$, where $S$ runs over all the objects of $\site$.
In particular, we may (and shall) assume that $\site$ has a terminal
object $S$. As $\M$ is perfect, it follows from
condition (c) of Definition \ref{defperfectsymmonidcmf} that we can replace
$R$ by any of its cofibrant resolution.
In particular, we may assume that $R_S$ is a cofibrant object of
$\Alg(\M)(S)$. We can thus define a termwise cofibrant cartesian monoid $R'$
as the family of monoids $f^*(R_S)$, where $f:X\To S$ runs over all the
objects of $\site\simeq \site/S$. There is a canonical morphism of homotopy
cartesian monoids $R'\To R$ which is a termwise weak equivalence.
We thus get, by condition (c) of Definition \ref{defperfectsymmonidcmf},
an equivalence of fibred categories
$$R\otimes^\derL_{R'}(-):\ho(\Mod{R'})\To\ho(\Mod R)\, .$$
We can thus replace $R$ by $R'$, which just means that we can assume
that $R$ is cartesian and termwise cofibrant. The first assertion follows then
easily from Proposition \ref{cmfPfibredRMod}.
In the case where $R$ is commutative, we prove that $\ho(\Mod R)$
is a $\Pmor$-fibred symmetric monoidal category as follows.
Let $f:X \To Y$ a morphism of $\site$. We would like to
prove that, for any object $M$ in $\ho(\Mod R)(X)$ and any object
$N$ in $\ho(\Mod R)(Y)$, the
canonical map
\begin{equation}\label{projformulaRmod1}
\derL f_\sharp(M\otimes^\derL_R f^*(N))\To \derL f_\sharp(M)\otimes^\derL_R N
\end{equation}
is an isomorphism. By adjunction, this is equivalent to prove that, for
any objects $N$ and $E$ in $\ho(\Mod R)(Y)$, the map
\begin{equation}\label{projformulaRmod2}
f^*\derR\sHom_R(N,E)\To \derR\sHom_R(f^*(N),f^*(E))
\end{equation}
is an isomorphism in $\ho(\Mod R)(X)$ (where $\derR \sHom_R$ stands for the
internal Hom of $\ho(\Mod R)$). But the forgetful functors
$$U:\ho(\Mod R)(X)\To \ho(\M)(X)$$
are conservative, commute with $f^*$ for any $\Pmor$-morphism $f$, and
commute with internal Hom: by adjunction, this follows immediately
from the fact that the functors
$$R\otimes^\derL (-):\ho(\M)(X)\To\ho(\Mod R)(X)\simeq \ho(\Mod{R'})(X)$$
are symmetric monoidal and define a morphism of $\Pmor$-fibred categories
(and thus, in particular, commute with $f_\sharp$ for any $\Pmor$-morphism $f$).
Hence, to prove that \eqref{projformulaRmod2} is an isomorphism, it is
sufficient to prove that its analog in $\ho(\M)$ is so, which
follows immediately from the fact that the analog of \eqref{projformulaRmod1}
is an isomorphism in $\ho(\M)$ by assumption.

For the last assertion, we are also
reduced to the case where $R$ and $S$ are cartesian and termwise
cofibrant, in which case this follows easily again from condition (c)
of Definition \ref{defperfectsymmonidcmf}.
\end{proof}

\begin{prop}\label{functorialityforgetful}
Let $\M$ be a combinatorial symmetric monoidal model category over $\site$
which satisfies the monoid axiom.\index{word}{monoid axiom}
Then, for any cartesian monoid\index{word}{monoid!cartesian}
 $R$ in $\M$ over $\site$ we have a Quillen morphism
$$R\otimes(-):\M\To \Mod R\, .$$
If, for any object $X$ of $\site$, the unit object $\unit_X$ is cofibrant in $\M(X)$
and the monoid $R_X$ is cofibrant in $\Alg(\M)(X)$, then
the forgetful functors also define a Quillen morphism
$$U:\Mod R\To \M\, .$$
\end{prop}

\begin{proof}
The first assertion is obvious. For the second one, note that, for any object $X$ of
$\site$, the monoid $R_X$ is also cofibrant as an object of $\M(X)$;
see Theorem \ref{cmfmonoids}. This implies that the forgetful functor
$$U:\Mod{R_X}\To \M(X)$$
is a left Quillen functor: by the small object argument and
by definition of the model category structure of Theorem
 \ref{abstractcmfmodules}~(i),
this follows from the trivial fact that the endofunctor
$$R_X\otimes(-):\M(X)\To\M(X)$$
is a left Quillen functor itself whenever $R_X$ is cofibrant in $\M(X)$.
\end{proof}

\begin{rem}
The results of the preceding proposition (as well as their proofs) are also true
in terms of $\Pmor_\cart$-fibred categories (\ref{remPacrtfibredcat})
over the category of $\site/S$-diagrams
for any object $S$ of $\site$ (whence over all $\site$-diagrams whenever $\site$
has a terminal object).
\end{rem}

\begin{paragr}
Consider now a noetherian scheme $S$ of finite dimension.
We choose a full subcategory of
the category of separated noetherian $S$-schemes of finite dimension
which is stable by finite limits, contains separated $S$-schemes of finite type,
and such that, for any \'etale $S$-morphism $Y\To X$, if $X$ is in $\sch/S$, so is $Y$.
We denote by  $\sch/S$ this chosen category of $S$-schemes.

We also fix an admissible class $\Pmor$
of morphisms of $\sch/S$ which contains the class
of \'etale morphisms.
\end{paragr}

\begin{df}
A property $\mathsf{P}$ of $\ho(\M)$, for $\M$
a stable combinatorial $\Pmor$-fibred model category over $\sch/S$,
is \emph{homotopy linear}\index{word}{homotopy linear}
 if the following implications are true.
\begin{itemize}
\item[(a)] If $\gamma:\M\To \M'$ is a Quillen equivalence
(i.e. a Quillen morphism which is termwise a Quillen equivalence)
between stable combinatorial $\Pmor$-fibred model category over $\sch/S$, then
$\M$ has property $\mathsf{P}$ is and only if $\M'$ has property $\mathsf{P}$.
\item[(b)] If $\M$ is a stable combinatorial symmetric monoidal $\Pmor$-model category
which satisfies the monoid axiom, and such that the unit $\unit_X$ of $\M(X)$
is cofibrant, then, for any cartesian and termwise cofibrant monoid $R$ in $\M$
over $\sch/S$, $\Mod R$ has property $\mathsf{P}$.
\end{itemize}
\end{df}

\begin{prop}\label{hmtlinearproperties}
The following properties are homotopy linear: $\AA^1$-homotopy invariance,
$\PP^1$-stability, the localization property, the property of proper transversality,
separability, semi-separability, $t$-descent (for a given Grothendieck topology $t$ on $\sch/S$).
\end{prop}

\begin{proof}
Property (a) of the definition above is obvious.
Property (b) comes from the fact that the forgetful functors
$$U:\ho(\Mod R)\To \ho(\M)$$
are conservative and commute with all the operations: $\derL f^*$
and $\derR f_*$ for any morphism $f$, as well as $\derL f_\sharp$ for any
$\Pmor$-morphism (by Proposition \ref{functorialityforgetful}).
Hence any property formulated in terms of equations involving
only these operations is homotopy linear.
\end{proof}

\part{Motivic complexes and relative cycles}
\markboth{Motivic complexes and relative cycles}{}
In this entire part,
 we adopt the special convention that smooth
 means smooth separated\index{word}{morphism!separated} of finite type.
 This concerns also the framework of
 premotivic categories: we assume the admissible class $\sm$
 is made of smooth separated morphisms of finite type.

This assumption is required by the use of the theory
 of finite correspondences
 (see more precisely Example \ref{ex:cor_graph&transpose}). 

\section{Relative cycles} \label{sec:cycles}

\begin{assumption}\label{ass:cycles}
In this entire section, $\sch$ is the category of noetherian schemes;
 any scheme is assumed to be noetherian.
\renewcommand{\Rc}{\Lambda}
We fix a subring $\Rc \subset \QQ$ which will be the ring of coefficients
of the algebraic cycles considered in the following section.
 When we want to be precise, we say \emph{$\Rc$-cycle} for
 "algebraic cycle with coefficients in $\Rc$". Otherwise,
 we simply say \emph{cycle}
 and the reader must assume that all algebraic cycles have their
 coefficients in the ring $\Rc$.
\end{assumption}

\subsection{Definitions}

\subsubsection{Category of cycles}

\begin{num} Let $X$ be a scheme. As usual, an element of the underlying set
 of $X$ will be called a \emph{point}\index{word}{point}
  and a morphism $\spec k \rightarrow X$ where $k$ is a field
   will be called a \emph{geometric point}.\index{word}{point!geometric}
We often identify a point $x \in X$ with the corresponding geometric point
  $\spec{\kappa_x} \rightarrow X$.
However, the explicit expression 
 "the point $\spec k \rightarrow X$" always refers
  to a geometric point. \\
As our schemes are assumed to be noetherian,
 any immersion $f:X \rightarrow Y$ is quasi-compact.
Thus, according to \cite[9.5.10]{EGA1},
 the \emph{schematic closure}\index{word}{schematic closure}
 $\bar X$ of $X$ in $Y$ exists
  which gives a unique factorization of $f$
$$
X \xrightarrow j \bar X \xrightarrow i Y
$$
such that $i$ is a closed immersion and $j$ is an open
immersion with dense image\footnote{Recall the scheme $\bar X$ 
 is characterized by the property of being the smallest sub-scheme of $Y$ 
 with the existence of such a factorization.}.
Note that when $Y$ is reduced, $\bar X$ coincide with the topological closure of
$X$ in $Y$ with its induced reduced subscheme structure.
In this case, we simply call $\bar Y$ the closure of $Y$ in $X$.
\end{num}
\renewcommand{\Rc}{\Lambda}
\begin{df} \label{df:cat_cycles}
A \emph{$\Rc$-cycle}\index{word}{cycle!$\Rc$-cycle}
 is a couple $(X,\alpha)$ such that
$X$ is a scheme and $\alpha$ is a $\Rc$-linear combination
of points of $X$. 
A generic point\index{word}{point!generic (of a cycle)}
 of $(X,\alpha)$ is a point
 which appears in the $\Rc$-linear combination $\alpha$
  with a non zero coefficient.
The support $\suppc \alpha$ of $\alpha$ is the closure
 of the generic points of $\alpha$, seen as a reduced closed subscheme
 of $X$.

A morphism of $\Rc$-cycles\index{word}{morphism!of $\Rc$-cycles}
 $(Y,\beta) \rightarrow (X,\alpha)$
is a morphism of scheme $f:Y \rightarrow X$ such that 
$f(\suppc \beta) \subset \suppc \alpha$.
We say this morphism is \emph{pseudo-dominant}
\index{word}{morphism!pseudo-dominant}
 if for any generic point $y$ of $(Y,\beta)$, $f(y)$
  is a generic point of $(X,\alpha)$.
\end{df}
When considering such a pair $(X,\alpha)$, we will denote it simply 
by $\alpha$ and refer to $X$ as the \emph{domain}
\index{word}{domain (of a $\Rc$-cycle)}
 of $\alpha$.
We also use the notation $\alpha \subset X$ to mean the domain
of the cycle $\alpha$ is the scheme $X$.
 
The category of $\Rc$-cycle is functorial in $\Rc$ with respect
 to morphisms of integral rings.
In what follows,
 cycles are assumed to have coefficients in $\Rc$ unless explicitly stated (following our conventions for this section,
 see Paragraph \ref{ass:cycles}).
\index{word}{cycle!$\Rc$-cycle|seealso{cycle}}
\begin{num} \label{num:properties_morphism_cycles}
Given a property $(\mathcal P)$ of morphisms of schemes, 
we will say that a morphism $f:\beta \rightarrow \alpha$ of cycles 
satisfies property $(\mathcal P)$ 
if the induced morphism $f|_{\suppc \beta}^{\suppc \alpha}$ satisfies
property $(\mathcal P)$. 
%A morphism $f:Y \rightarrow X$ will be said to be \emph{pseudo-dominant}
%if any irreducible component of $Y$ dominates an irreducible component of $X$. 
%Thus a morphism of cycles 
%$\sum_{j \in J} m_j.y_j \xrightarrow f \sum_{i \in I} n_i.x_i$
%is pseudo-dominant if and only if for any $j \in J$ there exist
%$i \in I$ such that $f(y_j)=x_i$.
\end{num}

\begin{df}
Let $X$ be a scheme. We denote by $X^{(0)}$ the set of generic
points of $X$. 
We define as usual the \emph{cycle associated with $X$}
\index{word}{cycle!associated}
 as the cycle 
with domain $X$:
$$
\acycl X=\sum_{x \in X^{(0)}} \lg(\mathcal O_{X,x}).x.
$$
The integer $\lg(\mathcal O_{X,x})$, length of an artinian local ring,
is called the \emph{geometric multiplicity}
\index{word}{multiplicity!geometric}
 of $x$ in $X$. 
\end{df}
%A convention which will come often is when we consider the case
%$X=\spec{L \otimes_k K}$ where $L/k$ and $K/k$ are extensions field:
%we simply write $\acycl{L \otimes_k K}=\acycl X$.
When no confusion is possible, we usually omit the
delimiters in the  notation $\acycl X$. 
As an example, we say that \emph{$\alpha$ is a cycle over $X$}
 to mean the existence of a structural morphism of cycles 
  $\alpha \rightarrow \acycl X$.

\begin{num} \label{num:std_form}
When $Z$ is a closed subscheme of a scheme $X$, 
we denote by $\acycl Z_X$
\index{notat}{bracketZX@$\acycl Z_X$}
 the cycle $\acycl Z$ considered as a cycle with domain $X$. \\
Consider a cycle $\alpha$ with domain $X$. Let $(Z_i)_{i \in I}$ be
the family of the reduced closure of generic points of $\alpha$.
Then we can write $\alpha$ uniquely as
 $\alpha=\sum_{i \in I} n_i.\acycl{Z_i}_X$. 
We call this writing the \emph{standard form}
\index{word}{cycle!standard form}
 of $\alpha$ for short.
\end{num}

\begin{df}
Let $\alpha=\sum_{i \in I} n_i.x_i$ be a cycle with domain $X$ and
$f:X \rightarrow Y$ be any morphism.

For any $i \in I$, put $y_i=f(x_i)$. Then $f$ induces an extension
field $\kappa(x_i)/\kappa(y_i)$ between the residue fields.
We let $d_i$ be the degree of this extension field in case it is
finite and $0$ otherwise.

We define the \emph{pushforward}
\index{word}{cycle!pushforward}
 of $\alpha$ by $f$ as the cycle with domain $Y$
$$
f_*(\alpha)=\sum_{i \in I} n_id_i.f(x_i).
$$
\end{df}
Thus, when $f$ is an immersion, $f_*(\alpha)$ is the same cycle as
$\alpha$ but seen as a cycle with domain $X$.
Remark also that we obtain the following equality
\begin{equation} \label{pushout&schematic_closure}
f_*\big(\langle X \rangle\big)=\left\langle \bar X \right\rangle_Y
\end{equation}
where $\bar X$ is the schematic closure of $X$ in $Y$
(indeed $X$ is a dense open subscheme in $\bar X$).
When $f$ is clear,
 we sometimes abusively put: $\acycl X_Y:=f_*(\acycl X)$.

%We always have a canonical morphism $\alpha \rightarrow f_*(\alpha)$.
%In case we have a morphism of cycles $\alpha \rightarrow \gamma$ 
%with domain the morphism $p:X \rightarrow S$, any commutative
%diagram of schemes
%$$
%\xymatrix@=10pt{
%X\ar^f[rr]\ar_p[rd] & & Y\ar^q[ld] \\
%& S &
%}
%$$
%induces a unique commutative diagram of cycles
%$$
%\xymatrix@=10pt{
%\alpha\ar[rr]\ar[rd] & & f_*(\alpha).\ar[ld] \\
%& \gamma &
%}
%$$
By transitivity of degrees, we obviously have $f_*g_*=(fg)_*$ for a
composable pair of morphisms $(f,g)$.

%\rem The obvious functor $\sch \rightarrow \cycl_A, 
%X \mapsto \langle X \rangle$ is fully faithful. But the notion of 
%a morphism of cycles does not truely reflect the structure of cycles.
%For example, it does not take into account multiplicities : 
%given any cycle $\alpha$, and any element $\lambda \in A-\{0\}$,
%$\alpha$ is isomorphic to $\lambda.\alpha$ according to our definition.
%A notion of morphism which will be more accurate to the study of cycles
%will be introduced later in definition \ref{df:pre-special}.

\begin{df} \label{df:open_restriction}
Let $\alpha=\sum_{i \in I} n_i.x_i$ be a cycle over a scheme $S$
 with domain $f:X\rightarrow S$
  and $U \subset S$ be an open subscheme.
Let $I'=\{i \in I \mid f(x_i) \in U\}$.
We define the \emph{restriction}
\index{word}{cycle!restriction}
 of $\alpha$ over $U$ as the cycle 
 $\alpha|_U=\sum_{i \in I'} n_i.x_i$ with domain $X \times_S U$
  considered as a cycle over $U$.
\end{df}
If $\alpha=\sum_{i \in I} n_i.\langle Z_i \rangle_X$,
  then obviously $\alpha|_U=\sum_{i \in I} n_i.\langle Z_i \times_S U \rangle_{X_U}$.
We state the following obvious lemma for convenience:
\begin{lm} \label{lm:open_restriction}
Let $S$ be a scheme, $U \subset S$ an open subscheme
 and $X$ be an $S$-scheme. Let $j:X_U \rightarrow X$ be
  the obvious open immersion.
\begin{enumerate}
\item[(i)] For any cycle $(X_U,\alpha')$,
 $\big(j_*(\alpha')\big)|_U=\alpha'$.
\item[(ii)] Assume $\bar U=S$.
For any cycle $(X,\alpha)$ pseudo-dominant over $S$,
 $j_*(\alpha\mid_U)=\alpha$.
\end{enumerate}
\end{lm}

%%%%%%%%%%%%%%%%%%%%%%%%%%%%%%%%%%%%%%%%%%%%%%%%%%%%%%%%%%%%%%%%%%%%%%%
%
% Hilbert
%

\subsubsection{Hilbert cycles}

\begin{num} \label{num:convention_flat_equidim}
Recall that a finite dimensional scheme $X$ is equidimensional
 -- we will say  \emph{absolutely equidimensional} -- 
\index{word}{equidimensional!absolutely}
 if its irreducible components have all the same dimension.

We will say that a flat morphism $f:X \rightarrow S$ is equidimensional
\index{word}{equidimensional!flat morphism}
 if it is of finite type
 and for any connected component $X'$ of $X$,
 there exists an integer $e \in \NN$ such that for any generic point $\eta$
 in $X'$, the fiber $f^{-1}\lbrack f(\eta) \rbrack$ is absolutely 
 equidimensional of dimension $e$.
\end{num}

\begin{df}
Let $S$ be a scheme.

Let $\alpha$ be a cycle over $S$ with domain $X$.
We say that $\alpha$ is a \emph{Hilbert cycle}
\index{word}{cycle!Hilbert}
 over $S$ if there exists a finite family $(Z_i)_{i \in I}$ of 
closed subschemes of $X$ which are flat equidimensional over $S$ 
and a finite family $(n_i)_{i \in I} \in \Rc^I$ such that
$$
\alpha=\sum_{i \in I} n_i.\acycl{Z_i}_X.
$$
\end{df}

\begin{ex} \label{ex_Hilbert_cycl}
Any cycle over a field $k$ is a Hilbert cycle over $\spec k$. 
Let $S$ be the spectrum of a	discrete valuation ring.
A cycle $\alpha=\sum_{i \in I} n_i.x_i$ over $S$ is a Hilbert cycle
if and only if each point $x_i$ lies over the generic points of $S$.
Indeed, an integral $S$-scheme is flat if and only if it is dominant.
\end{ex}

The following lemma follows almost directly
 from a result of \cite{SV1}:
\begin{lm}
Let $f:S ' \rightarrow S$ be a morphism of schemes and
 $X$ be an $S$-scheme of finite type.
 Put $X'=X \times_S S'$.
 
Let $(Z_i)_{i \in I}$ be a finite family
 of closed subschemes of $X$ such that each $Z_i$
 is flat equidimensional over $S$.
 We assume the following relation:
\begin{equation} \label{eq:lm:Hilbert_pullback}
\sum_{i \in I} n_i.\acycl{Z_i}_X=0
\end{equation}
Then we the following equality holds:
$$
\sum_{i \in I} n_i.\acycl{Z_i \times_S S'}_{X'}=0.
$$
\end{lm}
\begin{proof}
When we assume that for any index $i \in I$,
 $Z_i/S$ is equidimensional of dimension $e$,
 this lemma is exactly \cite[Prop. 3.2.2]{SV1}.
We show how to reduce to that case.

Up to adding more members to the family $(Z_i)$,
 we can always assume that $Z_i$ is connected.
 Then, because $Z_i/S$ is equidimensional by assumption,
 there exists an integer $e_i$ such that for any point $x \in Z_i^{(0)}$,
 the fiber $f^{-1}\lbrack f(x) \rbrack$ is absolutely equidimensional
 of dimension $e_i$. In particular the transcendence degree $d_x$ of the
 residual extension $\kappa_x/\kappa_{f(x)}$ satisfies the relation:
 $d_x=e_i$.

For any integer $e \in \NN$, we define the following subset of $I$:
$$
I_e=\{ i \in I \mid \forall x \in Z_i^{(0)}, d_x=e \}.
$$
Thus $(I_e)_{e \in \NN}$ is a partition of $I$.

One can rewrite the assumption \eqref{eq:lm:Hilbert_pullback} as follows: 
 for any point $x \in X$,
$$
\sum_{i \in I \mid x \in Z_i^{(0)}} n_i.\mathrm{lg}(\mathcal O_{Z_i,x})=0.
$$
In particular, given any integer $e \in \NN$,
 we deduce that the family $(Z_i)_{i \in I_e}$ still
 satisfies the relation \eqref{eq:lm:Hilbert_pullback}.
 As any member of this family is equidimensional of dimension $e$,
 we can apply \cite[Prop. 3.2.2]{SV1} to $(Z_i)_{i \in I_e}$.
 This concludes.
\end{proof}

\begin{num}
Consider a Hilbert $S$-cycle $\alpha \subset X$ and 
a morphism of schemes $f:S' \rightarrow S$. Put $X'=X \times_S S'$.
We choose a finite family $(Z_i)_{i \in I}$ of flat equidimensional $S$-schemes
and a finite family $(n_i)_{i \in I} \in \Rc^I$
such that $\alpha=\sum_{i \in I} n_i.\acycl{Z_i}_X$.
The previous lemma says exactly that the cycle
$$
\sum_{i \in I} n_i.\acycl{Z_i \times_S S'}_{X'}
$$
depends only on $\alpha$ and not on the chosen families.
\end{num}
\begin{df}
\label{df:pullback_Hilbert}
Adopting the preceding notations and hypothesis, we define the
 \emph{pullback cycle}
\index{word}{cycle!pullback!of Hilbert cycles}
 of $\alpha$ along the morphism $f:S' \rightarrow S$ as the cycle
with domain $X'$
$$
\alpha \otimesf_S S'=\sum_{i \in I} n_i.\acycl{Z_i \times_S S'}_{X'}.
$$
\index{notat}{alphatensorflat@$\alpha \otimesf_S S'$}
\end{df}

In this setting the following lemma is obvious:
\begin{lm}
\label{transitivity_Hilbert}
Let $\alpha$ be a Hilbert cycle over $S$, 
and $S'' \rightarrow S' \rightarrow S$ be morphisms of schemes.

Then $(\alpha \otimesf_S S') \otimesf_{S'} S''=\alpha \otimesf_S S''$.
\end{lm}

We will use another important computation from \cite{SV1} 
(it is a particular case of \emph{loc. cit.}, 3.6.1).
\begin{prop}
\label{Hilbert_proj_formula}
Let $R$ be a discrete valuation ring with residue field $k$. \\
Let $\alpha \subset X$ be a Hilbert cycle over $\spec R$
and $f:X \rightarrow Y$ a morphism over $\spec R$. 
We denote by $f':X' \rightarrow Y'$ the pullback of $f$
 over $\spec k$.

Suppose that the support of $\alpha$
 is proper with respect to $f$.

Then $f_*(\alpha)$ is a Hilbert cycle over $R$
and 
the following equality of cycles holds in $X'$:
$$
f'_*(\alpha \otimesf_S k)=f_*(\alpha) \otimesf_S k.
$$
\end{prop}

\begin{df}
Let $p:\tilde S \rightarrow S$ be a birational morphism.
Let $C$ be the minimal closed subset of $S$ such that
$p$ induces an isomorphism
$(\tilde S-\tilde S \times_S C) \rightarrow (S-C)$.

Consider $\alpha=\sum_{i \in I} n_i.\acycl{Z_i}_X$ a cycle over $S$ written
in standard form.

We define the strict transform\index{word}{strict transform}
 $\tilde Z_i$ of 
the closed subscheme $Z_i$ in $X$ along $p$ as the schematic closure of 
$(Z_i-Z_i \times_S C) \times_S \tilde S$ in $X \times_S \tilde S$.
We define the strict transform of $\alpha$ along $p$ as the cycle over
$\tilde S$
$$
\tilde \alpha=\sum_{i \in I} n_i.\acycl{\tilde Z_i}_{X \times_S \tilde S}.
$$
\index{notat}{alphatilde@$\tilde \alpha$}
\end{df}

As in \cite{SV1}, we remark that a corollary of the platification
theorem of Gruson-Raynaud is the following:
\begin{lm}
\label{lm:Hilbertification}
Let $S$ be a reduced scheme and $\alpha$ be 
a pseudo-dominant cycle over $S$.

Then there exists a dominant blow-up $p:\tilde S \rightarrow S$ such that the
strict transform $\tilde \alpha$ of $\alpha$ along $p$ is a Hilbert
cycle over $\tilde S$.
\end{lm}

We conclude this part by recalling an elementary lemma
 about cycles and Galois descent\index{word}{descent!Galois}
 which will be used
  extensively in the next sections:
\begin{lm}
\label{lm:basic_cycles}
Let $L/K$ be an extension of fields and $X$ be a $K$-scheme.
We put $X_L=X \times_K \spec L$ and consider the faithfully
flat morphism $f:X_L \rightarrow X$.

Denote by $\cycl(X)$ (resp. $\cycl(X_L)$) the cycles with domain $X$
(resp. $X_L$).
\begin{enumerate}
\item The morphism $f^*:\cycl(X) \rightarrow \cycl(X_L),
\beta \mapsto \beta \otimesf_K L$
is a monomorphism.
\item Suppose $L/K$ is finite.
For any $K$-cycle $\beta \in \cycl(X)$, \\
$f_*(\beta \otimesf_K L)=[L:K].\beta$.
\item Suppose $L/K$ is finite normal with Galois group $G$.

The cycles in the image of $f^*$ are invariant under the
action of $G$. For any cycle $\beta \in \cycl(X_L)^G$, 
there exists a unique cycle $\beta_K \in \cycl(X)$ such that
$$\beta_K \otimesf_K L=[L:K]_i.\beta$$
where $[L:K]_i$ is the
\emph{inseparable degree} of $L/K$.
\end{enumerate}
\end{lm}

%%%%%%%%%%%%%%%%%%%%%%%%%%%%%%%%%%%%%%%%%%%%%%%%%%%%%%%%%%%%%%%%%%%%%%
%
% Special
%
%

\subsubsection{Specialization}

The aim of this section is to give conditions on cycles so that
one can define a \emph{relative tensor product} on them.
\begin{df} \label{df:pre-special}
Consider two cycles $\alpha=\sum_{i \in I} n_i.s_i$ 
and $\beta=\sum_{j \in J} m_i.x_j$. Let $S$ be the support of
$\alpha$.

A morphism $\beta \xrightarrow f \alpha$ of cycles is said to be
\emph{pre-special}
\index{word}{cycle!pre-special (morphism of)}
 if it is of finite type and for any $j \in J$, 
there exists $i \in I$ such that $f(x_j)=s_i$ and $n_i|m_j$ in $\Rc$.
We define the reduction of $\beta/\alpha$ as the cycle over $S$
$$
\beta_0=\sum_{j \in J, f(x_j)=s_i} \frac{m_j}{n_i}.x_j.
$$
\end{df}

\begin{ex}
Let $S$ be a scheme and $\alpha$ a Hilbert $S$-cycle.
Then the canonical morphism of cycles
$\alpha \rightarrow \acycl S$ is pre-special.
If $S$ is the spectrum of a discrete valuation ring,
an $S$-cycle $\alpha$ is pre-special if and only if
it is a Hilbert $S$-cycle.
\end{ex}

\begin{df}\label{df:fat_points}
Let $\alpha$ be a cycle.

A \emph{point}
\index{word}{point!of a cycle}
(resp. \emph{trait})
\index{word}{trait!of a cycle}
 of $\alpha$ will be a morphism of the form
$\spec k \xrightarrow x \alpha$ 
(resp. $\spec R \xrightarrow \tau \alpha$)
such that $k$ is a field 
(resp. $R$ is a discrete valuation ring).
We simply say that $x$ (resp. $\tau$) is \emph{dominant}
 if the image of the generic point in the domain of $\alpha$
  is a generic point of $\alpha$. \\
Let $x:\spec{k_0} \rightarrow \alpha$ be a point. An extension of $x$ 
 will be a point $y$ on $\alpha$
  of the form $\spec k \rightarrow \spec{k_0} \xrightarrow x \alpha$.

A \emph{fat point}
\index{word}{point!fat point (of a cycle)}
 of $\alpha$ will be a couple of morphisms
$$
\spec k \xrightarrow s \spec R \xrightarrow \tau \alpha
$$
such that $\tau$ is a dominant trait 
 and the image of $s$ is the closed point of $\spec R$. \\
Given a point $x:\spec k \rightarrow \alpha$,
 a fat point over $x$ is a factorization of $x$ through a dominant
 trait as above.
\end{df}
In the situation of the last definition,
 we denote simply by $(R,k)$ a fat point over $x$, 
 without indicating in the notation the morphisms $s$ and $\tau$.

\begin{rem} With our choice of terminology, 
 a point of $\alpha$ is in general an extension of a
 specialization of a generic point of $\alpha$.
As a further example, a dominant point of $\alpha$ is an extension of
 a generic point of $\alpha$.
\end{rem}

\begin{lm} \label{lm:exists_fat_pts}
For any cycle $\alpha$ and any non dominant point $x:\spec{k_0} \rightarrow \alpha$,
 there exists an extension $y:\spec k \rightarrow \alpha$ of $x$
  and a fat point $(R,k)$ over $y$.
\end{lm}
\begin{proof}
Replacing $\alpha$ by its support $S$, we can assume $\alpha=\langle S \rangle$.
Let $s$ be the image of $x$ in $S$, $\kappa$ its residue field.
We can assume $S$ is reduced, irreducible by taking one irreducible component
containing $s$, and local with closed point $s$.
Let $S=\spec A$, $K=\mathrm{Frac}(A)$. According to \cite[7.1.7]{EGA2}, there
exists a discrete valuation ring $R$ such that $A \subset R \subset K$, and $R/A$
is an extension of local rings. Then any composite extension $k/\kappa$ of
$k_0$ and the residue field of $R$ over $\kappa$ gives the desired fat point $(R,k)$.
\end{proof}

\begin{df} \label{df:foncteur_sp}
Let $\beta \rightarrow \alpha$ be a pre-special morphism of cycles.
Consider $S$ the support of $\alpha$ and $X$ the domain of $\beta$.
Let $\beta_0=\sum_{j \in J} m_j.\acycl{Z_j}_X$ be the reduction 
of $\beta/\alpha$ written in standard form.
\begin{enumerate}
\item
Let $\spec K \rightarrow \alpha$ be a dominant point.
We define the following cycle over $\spec K$ with domain 
$X_K=X \times_S \spec K$:
$$
\beta_K=\sum_{j \in J} m_j.\acycl{Z_j \times_S \spec K}_{X_K}.
$$
\item
Let $\spec R \xrightarrow \tau S$ be a dominant trait,
$K$ be the fraction field of $R$ and $j:X_K \rightarrow X_R$
be the canonical open immersion.
We define the following cycle over $R$ with domain $X_R$:
$$
\beta_R=j_*(\beta_K).
$$
According to example \ref{ex_Hilbert_cycl}, $\beta_R$ is a Hilbert
cycle over $R$.
\item
Let $x:\spec k \rightarrow \alpha$
be a point on $\alpha$
and $(R,k)$ be a fat point over $x$. \\
We define the \emph{specialization of $\beta$ along the fat point}
$(R,k)$
\index{word}{cycle!specialization}
 as the cycle 
$$\beta_{R,k}:=\beta_R \otimesf_R k$$
\index{notat}{betaindexRk@$\beta_{R,k}$}
using the above notation and definition \ref{df:pullback_Hilbert}. 
It is a cycle over $\spec k$ with domain $X_k=X \times_S \spec k$.
\end{enumerate}
\end{df}

\begin{rem} \label{sp&open_ngb_gen_pt}
Let $\beta$ be an $S$-cycle,
 $x:\spec K \rightarrow S$ be a dominant point
 and $U$ be an open neighborhood of $x$ in $S$. \\
Then if $\beta$ is pre-special over $S$, $\beta|_U$ is pre-special over $U$
  and $\beta_K=(\beta|_U)_K$. \\
If $\tau:\spec R \rightarrow S$ 
 (resp. $(R,k)$) is a trait 
  (resp. fat point) with generic point $x$,
 we also get $\beta_R=(\beta|_U)_R$
  (resp. $\beta_{R,k}=(\beta|_U)_{R,k}$).
\end{rem}

%\rem Following a universal abuse of notations in algebraic geometry,
% we have neglected the morphism $x:\spec K \rightarrow \alpha$
%  in the notation $\beta_K$ of the first point. Sometimes, it will be
%   necessary on the countrary to be precised about this morphism ; 
%   then we will replace the notation $\beta_K$ by $\beta_x$.
%   (A similar remark holds for the two other cases of the above notation
%    but we will not be placed in a situation where a confusion can appear
%     for these latter cases).  

\begin{num}
%\label{comput_special1}
Let $S$ be a reduced scheme,
 and $\beta=\sum_{i \in I} n_i.x_i$ be an $S$-cycle with domain $X$.
For any index $i \in I$, let $\kappa_i$ be the residue field of $x_i$.

Consider a dominant point $x:\spec K \rightarrow S$.
Let $\eta$ be its image in $S$ and $F$ be the residue field of $\eta$.
We put $I'=\{i \in I \mid f(x_i)=\eta \}$
 where $f:X \rightarrow S$ is the structural morphism.
With these notations, we get
$$
\beta_K=\sum_{i \in I'} n_i.\langle \spec{\kappa_i \otimes_F K} \rangle_{X_K},
$$
and for a dominant trait $\spec R \rightarrow S$ with generic point $x$,
\begin{equation}\label{eq:comput_specialR}
\beta_R=\sum_{i \in I'} n_i.\langle \spec{\kappa_i \otimes_F K} \rangle_{X_R},
\end{equation}
where $\spec{\kappa_i \otimes_F K}$ is seen as a subscheme of $X_K$
(resp. $X_R$).

Consider a fat point $(R,k)$ with generic point $x$
 and write $\beta=\sum_{i \in I} n_i.\langle Z_i \rangle_X$
  in standard form (\emph{i.e.} $Z_i$ is the closure of $\{x_i\}$ in $X$).
Then according to \eqref{pushout&schematic_closure}, we obtain\footnote{
This shows that
our definition coincide with the one given in \cite{SV1}
 (p. 23, paragraph preceding 3.1.3)
  in the case where $\alpha=\langle S \rangle$, $S$ reduced.
}
$$
\beta_{R,k}
 =\sum_{i \in I'} n_i.
  \left\langle \overline{Z_{i,K}} \times_R \spec k \right\rangle_{X_k}
$$
where $Z_{i,K}=Z_i \times_S \spec K$ is considered as a subscheme of $X_K$
 and the schematic closure is taken in $X_R$.

Considering the description of the schematic closure for the generic
fiber of an $R$-scheme (\textit{cf.} \cite[2.8.5]{EGA4}), we obtain the following
way to compute $\beta_{R,k}$. By definition, $R$ is an $F$-algebra.
For $i \in I'$,
 let $A_i$ be the image of the canonical morphism
$$
\kappa_i \otimes_F R \rightarrow \kappa_i \otimes_F K.
$$
It is an $R$-algebra without $R$-torsion. Moreover,
the factorization
$$
\spec{\kappa_i \otimes_F K} \rightarrow \spec{A_i}
 \rightarrow \spec{\kappa_i \otimes_F R}
$$
defines $\spec{A_i}$ as the schematic closure
 of the left hand side in the right hand
side (\textit{cf.} \cite[2.8.5]{EGA4}). In particular,
 we get an immersion $\spec{A_i \otimes_R k} \rightarrow X_k$
  and the nice formula:
$$
\beta_{R,k}
 =\sum_{i \in {I'}} n_i.
  \left\langle \spec{A_i \otimes_R k} \right\rangle_{X_k}.
$$
\end{num}
%
%\begin{rem}
%\label{comput_special2}
%Moreover,
%we can give a purely algebraic description of this specialization.
%Suppose in the preceding notations that $S$ (resp. $X$) is integral 
%and $\alpha=\acycl S$ (resp. $\beta=\acycl X$). Let $E$ (resp. $L$)
%be the function field of $S$ (resp. $X$). \\
%First consider the immersion 
%$i:\spec{L \otimes_E K} \rightarrow X \times_S \spec R$. Then according
%to the previous computation,
%$$
%\beta_R=i_*\acycl{L \otimes_E K}.
%$$
%Moreover, consider the $R$-linear morphism 
%$\phi:L \otimes_E R \rightarrow L \otimes_E K$. Its image $\Im(\phi)$
%if a torsion free $R$-module which spectrum is the subscheme of generic 
%points of the schematic closure of $X_K$ in $X_R$. \\
%If now we consider the immersion
%$j:\spec{\Im(\phi) \otimes_R k} \rightarrow X \times_S \spec k$, we obtain
%the final formula
%$$
%\beta_{R,k}
%=j_*\acycl{\Im(L \otimes_E R \rightarrow L \otimes_E K) \otimes_R k}.
%$$
%
%\end{rem}

\begin{df} \label{df:be_special}
Consider a morphism of cycles $f:\beta \rightarrow \alpha$
 and a point $x:\spec{k_0} \rightarrow \alpha$. \\
We say that $f$ is special
\index{word}{cycle!special (morphism of)}
 at $x$
if it is pre-special
 and for any extension $y:\spec k \rightarrow \alpha$ of $x$,
  for any fat points $(R,k)$ and $(R',k)$ over $y$,  
   the equality $\beta_{R,k}=\beta_{R',k}$ holds in $X_k$. 
Equivalently, we say that $\beta/\alpha$ is special at $x$. \\
We say that $f$ is special (or that $\beta$ is special over $\alpha$)
 if it is special at every point of $\alpha$.
\end{df}

\begin{num}\label{num:special&SV}
Here is a dictionary to compare the above definition with that
 of Suslin and Voevodsky in \cite[3.1.3]{SV1}.

Consider a pre-special morphism $\beta/\alpha$.
 Let $X$ be the domain of $\beta$,
 $S$ be the support of $\alpha$ and $\beta_0$ be the reduction
 (see Definition \ref{df:pre-special})
 of $\beta/\alpha$, seen as a pre-special $S$-cycle.

Then the following conditions are equivalent:
\begin{enumerate}
\item[(i)] $\beta/\alpha$ is special;
\item[(ii)] $\beta_0/S$ is special.
\end{enumerate}
This follows from the very definition of the specialization of
 $\beta/\alpha$ along fat points (Definition \ref{df:foncteur_sp}).
 
Moreover, condition (ii) says exactly that $\beta_0$
 is a relative cycle on $X$ over $S$ in the sense
 of Definition 3.1.3 of \cite{SV1}.
\end{num}

\begin{rem} \label{rm:trivial_special}
\begin{enumerate}
\item Trivially, $f$ is special at every dominant point of $\alpha$.
\item Given an extension $y$ of $x$,
 it is equivalent for $f$ to be special at $x$ or at $y$ 
 (use Lemma \ref{lm:basic_cycles}(1)).
Thus, in the case where $\alpha=\langle S \rangle$,
 we can restrict our attention to the points $s \in S$.
\item According to \ref{sp&open_ngb_gen_pt},
 the property that $\beta/S$ is special at $s \in S$ depends only
  on an open neighborhood $U$ of $s$ in $S$.
More precisely, the following conditions are equivalent:
\begin{enumerate}
\item[(i)] $\beta$ is special at $s$ over $S$.
\item[(ii)] $\beta|_U$ is special at $s$ over $U$.
\end{enumerate}
\end{enumerate}
\end{rem}

\begin{ex}
Let $S$ be a scheme and $\beta$ be a Hilbert cycle over $S$.
We have already seen that $\beta \rightarrow \acycl S$ is pre-special.
The next lemma shows this morphism is in fact special.
\end{ex}

\begin{lm}
\label{Hilbert&special}
Let $S$ be a scheme and $\beta$ be a Hilbert cycle over $S$.
Consider a point $x:\spec k \rightarrow S$ and a fat point $(R,k)$
over $x$.

Then $\beta_{R,k}=\beta \otimesf_S k$.
\end{lm}
\begin{proof}
According to the preceding definition and 
Lemma \ref{transitivity_Hilbert} it is sufficient to prove
$\beta_R=\beta \otimesf_S R$.
As the two sides of this equation are unchanged when replacing 
$\beta$ by the reduction $\beta_0$ of $\beta/S$, we can assume that
$S$ is reduced. By additivity, we are reduced to the case where 
$\beta=\acycl X$ is the fundamental cycle associated with a flat 
$S$-scheme $X$.
According to \ref{pushout&schematic_closure},
$\beta_R=\big\langle \overline{X_K} \big\rangle_{X_R}$.
Applying now \cite[2.8.5]{EGA4}, $\overline{X_K}$
is the unique closed subscheme $Z$ of $X_R$ such that
 $Z$ is flat over $\spec R$ and $Z \times_R \spec K=X_K$.
Thus, as $X_R$ is flat over $\spec R$,
 we get $\overline{X_K}=X_R$ and this concludes.
\end{proof}

\begin{lm}
\label{proper_transform&special}
Let $p:\tilde S \rightarrow S$ be a birational morphism 
 and consider a commutative diagram
$$
\xymatrix@R=-2pt@C=16pt{
& & \tilde S\ar^p[dd] \\
\spec k\ar[r] & \spec R\ar[ru]\ar[rd] & \\
& & S
}
$$
such that $(R,k)$ is a fat point of $\tilde S$ and $S$.

Consider a pre-special cycle $\beta$ over $S$ and 
$\tilde \beta$ its strict transform\index{word}{strict transform}
 along $p$.
Then, $\tilde \beta$ is pre-special and 
$\tilde \beta_{R,k}=\beta_{R,k}$.
\end{lm}
\begin{proof}
Using \ref{sp&open_ngb_gen_pt}, we reduce to the case where
 $p$ is an isomorphism which is trivial.
%Write $\beta=\sum_{j \in J} m_j.\acycl{Z_j}_X$ in standard form.
%Let $\tilde Z_j$ be the strict transform of $Z_j$ along $p$.
%The morphism $\tilde Z_j \rightarrow Z_j$ induced by $p$
%is birational. Let $K$ be the fraction field of $R$.
%As $\spec K$ dominates an irreducible component of $\tilde S$, 
%which is birational to an irreducible component of $S$,
%the induced morphism 
%$\tilde Z_j \otimes_{\tilde S} K \rightarrow Z_j \otimes_S K$ 
%is an isomorphism. This concludes by the very definition.
\end{proof}

\begin{lm} \label{lm:special&Hilbert_blow-up}
Let $S$ be a reduced scheme,
 $x:\spec{k_0} \rightarrow S$ be a point
  and $\alpha$ be a pre-special cycle over $S$.
Let $p:\tilde S \rightarrow S$ be a dominant blow-up
 such that the strict transform $\tilde \alpha$ of $\alpha$ along $p$
  is a Hilbert cycle over $\tilde S$.
Then the following conditions are equivalent:
\begin{enumerate}
\item[(i)] $\alpha$ is special at $x$.
\item[(ii)] for every couple of points $x_1,x_2:\spec k \rightarrow \tilde S$
 such that $p \circ x_1=p \circ x_2$ and $p \circ x_1$ is an extension of $x$,
  $\tilde \alpha \otimesf_{\tilde S} x_1=\tilde \alpha \otimesf_{\tilde S} x_2$.
\end{enumerate}
\end{lm}
\begin{proof}
The case where $x$ is a dominant point follows from the definitions
 and the fact $p$ is an isomorphism at the generic point. 
We thus assume $x$ is non dominant. \\
$(i) \Rightarrow (ii)$: 
Applying Lemma \ref{lm:exists_fat_pts} to $x_i$, $i=1,2$,
 we can find an extension $x'_i:\spec{k_i} \rightarrow \tilde S$
  of $x_i$ and a fat point $(R_i,k_i)$ over $x'_i$.
Taking a composite extension $L$ of $k_1$ and $k_2$ over $k$,
 we can further assume $L=k_1=k_2$ and $p \circ x'_1=p\circ x'_2$.
Then for $i=1,2$, we get
$$
\xymatrix@=18pt{
\big(\tilde \alpha \otimesf_{\tilde S} x_i\big)
  \otimesf_{k} L\ar@{=}^-{\ref{transitivity_Hilbert}}[r]
 & \tilde \alpha \otimesf_{\tilde S} x'_i\ar@{=}^-{\ref{Hilbert&special}}[r]
  & \tilde \alpha_{R_i,L}\ar@{=}^-{\ref{proper_transform&special}}[r]
   & \alpha_{R_i,L},
}
$$
and this concludes according to \ref{lm:basic_cycles}(1). \\
$(ii) \Rightarrow (i)$:
Consider an extension $y:\spec k \rightarrow \alpha$
 over $x$ and two fat point $(R_1,k)$, $(R_2,k)$ over $y$.
Fix $i \in \{1,2\}$.
As $p$ is proper birational, the trait $\spec{R_i}$ on $S$
 can be extended (uniquely) to $\tilde S$. 
Let $x_i:\spec k \rightarrow \spec{R_i} \rightarrow \tilde S$ be the induced point.
Then the following computation allows concluding: \\
$
\xymatrix@=18pt{
&&&&& \alpha_{R_i,k}\ar@{=}^-{\ref{proper_transform&special}}[r]
 & \tilde \alpha_{R_i,k}\ar@{=}^-{\ref{Hilbert&special}}[r]
  & \tilde \alpha \otimesf x_i
}
$
\end{proof}

%Note that specializations defines a conservative family for cycles 
%in the following sense:
%\begin{lm}
%Let $\alpha$ be a cycle and
%suppose given a surjective family of points 
%$(x_i:\spec{k_i} \rightarrow \alpha)_{i \in I}$
%and for any $i \in I$ a fat point $(R_i,k_i)$ over $x_i$.
%
%Let $S$ be the domain of $\alpha$, $X$ an $S$-scheme
%and $\beta \rightarrow \alpha$, $\beta' \rightarrow \alpha$
%two pre-special morphisms with underlying domain $X \rightarrow S$.
%
%Then we have the implication
%$$
%\forall i \in I, \beta_{R_i,k_i}=\beta'_{R_i,k_i}
%\Rightarrow \beta=\beta'.
%$$
%\end{lm}
%\begin{proof}
%By considering the reduction of $\beta/\alpha$ and $\beta'/\alpha$, 
%We are reduced to the case where $\alpha=[S]$ and $S$ is a reduced scheme.
%
%Let $s$ be any set theoretic point of $S$. 
%Choose an irreducible component of $S$
%which contains $s$ and let $K$ be its residue field.
%Then according to \cite{EGA4}, there exists a dominant trait 
%$\spec R \xrightarrow \tau S$ such that the generic point lies 
%isomorphically over $\spec K$, the closed point $\spec k$ lies over $s$ 
%and $k/\kappa(s)$ is finite.
%\end{proof}

%%%%%%%%%%%%%%%%%%%%%%%%%%%%%%%%%%%%%%%%%%%%%%%%%%%%%%%%%%%%%%%%%%%%%%%
%
% the pullback
%

\subsubsection{Pullback}

\begin{num}  \label{fund_inter_hyp}
In this part, we construct a \emph{pullback}
\index{word}{cycle!pullback}
 which extends the pullback defined by Suslin et Voevodsky
 in \cite[3.3.1]{SV1} to the case of morphism of cycles.
Consider the situation of a diagram of cycles
$$
\xymatrix@C=18pt@R=4pt{
\qquad \quad & \beta\ar^f[dd] & & \qquad \qquad  & X\ar[dd] \\
&& \subset && \\
\alpha'\ar[r] & \alpha & & S'\ar[r] & S
}
$$
where the diagram on the right is the domain of the one on the left.
Let $n$ be exponential characteristic of $\suppc{\alpha'}$.

The pullback of $\beta$, considered as an $\alpha$-cycle, over $\alpha'$
 will be a $\Rc[1/n]$-cycle denoted by $\beta \otimes_\alpha \alpha'$.
\index{notat}{betatensoralphaalphaprime@$\beta \otimes_\alpha \alpha'$}
It will fits into the following commutative diagram of cycles
$$
\xymatrix@C=18pt@R=4pt{
\beta \otimes_\alpha \alpha'\ar[r]\ar[dd] & \beta\ar[dd]
 & & X \times_S S'\ar[r]\ar[dd] & X\ar[dd] \\
&& \subset && \\
\alpha'\ar[r] & \alpha & & S'\ar[r] & S
}
$$
where the right commutative square is again the support of the left one.

It will be defined under an assumption on $\beta/\alpha$
 and is therefore non symmetric\footnote{See further \ref{commutativity}
  for this question.}. This assumption will imply that $\beta/\alpha$
is pre-special, and the first property of $\beta\otimes_\alpha \alpha'$
is that it is pre-special over $\alpha'$. \\
We define this product in three steps
  in which the following properties\footnote{
All these properties except (P3) will be particular cases of the associativity
of the pullback. }
 will be a guideline:
\begin{enumerate}
\item[(P1)] Let $S_0$ be the support of $\alpha$
 and $\beta_0$ be the reduction (see Definition \ref{df:pre-special})
 of $\beta/\alpha$, as an $S_0$-cycle. 
Consider the canonical factorization
  $\alpha' \rightarrow S_0 \rightarrow \alpha$. \\
Then, $\beta \otimes_\alpha \alpha'=\beta_0 \otimes_{S_0} \alpha'$.
\item[(P2)] Consider a commutative diagram
$$
\xymatrix@C=16pt@R=12pt{
\spec E\ar[r] & \spec{R'}\ar[r]\ar[d]\ar@{}|{(*)}[rd] & \spec R\ar[d] \\
& \alpha'\ar[r] & \alpha
}
$$
such that $(R,E)$ (resp. $(R',E)$) is a fat point on $\alpha$
 (resp. $\alpha'$). \\
Then, $(\beta \otimes_\alpha \alpha')_{R',E}=\beta_{R,E}$.
\end{enumerate}
\noindent Assume $\alpha' \rightarrow \alpha=\langle S' \rightarrow S \rangle$.
\begin{enumerate}
\item[(P3)] If $\beta$ is a Hilbert cycle over $S$, $\beta \otimes_S S'=\beta \otimesf_S S'$.
\item[(P4)] Consider a factorization $S' \rightarrow U \xrightarrow j S$
 such that $j$ is an open immersion. 
  Then $\beta \otimes_S S'=\beta|_U \otimes_U S'$.
\item[(P5)] Consider a factorization $S' \rightarrow \tilde S \xrightarrow p S$
 such that $p$ is a birational morphism. 
 Then $\beta \otimes_S S'=\tilde \beta \otimes_{\tilde S} S'$.
\end{enumerate}
\end{num}

\begin{lm}
\label{lm:existence_pullback_field}
Consider the hypothesis of \ref{fund_inter_hyp}
 in the case where $\alpha'=\spec k$ is a point $x$ of $\alpha$.

We suppose that $f$ is special at $x$.

Then the pre-special $\Rc[1/n]$-cycle $\beta \otimes_\alpha k$
 exists and is uniquely determined by property (P2) above. 
We also put $\beta_k:=\beta \otimes_\alpha k$.

The properties (P1) to (P5) are fulfilled and in addition : \\
(P6) For any extension fields $L/k$, $\beta_L=\beta_k \otimesf_k L$.
\end{lm}
\begin{proof}
According to Lemma \ref{lm:exists_fat_pts} there always exists
 a fat point $(R,E)$ over an extension of $x$. Thus the unicity
  statement follows from \ref{lm:basic_cycles}(1).

For the existence, 
 we first consider the case where $\alpha=\langle S \rangle$
  is a reduced scheme.
Applying Lemma \ref{lm:Hilbertification}, there exists
a blow-up $p:\tilde S \rightarrow S$ such that the strict transform 
$\tilde \beta$ of $\beta$ along $p$ is a Hilbert cycle over 
$\tilde S$.

As $p$ is surjective, the fiber $\tilde S_k$ is a non-empty
 algebraic $k$-scheme.
Thus, it admits a closed point given by a finite extension $k'_0$ of $k$.
Let $k'/k$ be a normal closure of $k'_0/k$
 and $G$ be its Galois group.
As $\beta/S$ is special at $x$ by hypothesis,
 Lemma \ref{lm:special&Hilbert_blow-up} implies that
  $\tilde \beta \otimesf_{\tilde S} k'$ is $G$-invariant.
Thus, applying Lemma \ref{lm:basic_cycles}, there exists a unique
cycle $\beta_k \subset X_k$ with coefficients in $\Rc[1/n]$ such that
$\beta_k \otimesf_k k'=\tilde \beta \otimesf_{\tilde S} k'$.

We prove (P2). Given a diagram $(*)$ with $\alpha'=\spec k$,
 we first remark that $(\beta_k)_{R',E}=\beta_k \otimesf_k E$.
As $p$ is proper birational,
 the dominant trait $\spec R \rightarrow S$ lifts 
  to a dominant trait $\spec R \rightarrow \tilde S$.
Let $E'/k$ be a composite extension of $k'/k$ and $E/k$. 
With these notations, we get the following computation:
$$
\xymatrix@=16pt{
\beta_{R,E} \otimesf_E E'\ar@{=}^-{\ref{proper_transform&special}}[r]
 & \tilde \beta_{R,E} \otimesf_E E'\ar@{=}^-{\ref{Hilbert&special}}[r]
 & \tilde \beta \otimesf_{\tilde S} E'\ar@{=}^-{\ref{transitivity_Hilbert}}[r]
 & (\tilde \beta \otimesf_{\tilde S} k') \otimesf_E E'\ar@{=}[r]
 & \beta_k \otimesf_k E',
}
$$
so that we can conclude by applying \ref{lm:basic_cycles}(1).

In the general case, we consider he support $S$ of $\alpha$
 abd $\beta_0/S$ the reduction of $\beta/\alpha$. 
According to (P1), we are led to put $\beta_k:=(\beta_0)_k$
 with the help of the preceding case. 
Considering the definition of specialization along fat points,
 we easily check this cycle satisfies property (P2).

Finally, property (P6) (resp. (P3), (P5))
 follows from the unicity statement applying
  lemmas \ref{lm:exists_fat_pts}, \ref{lm:basic_cycles}(1)
   (resp. and moreover Lemma \ref{Hilbert&special},
    \ref{proper_transform&special}).
\end{proof}

\begin{rem} \label{rm:pullback_dominant_pt}
In the case where $x$ is a dominant point, the cycle $\beta_k$ defined
 in the previous proposition agrees with the one defined in \ref{df:foncteur_sp}(1).
 \end{rem}

\begin{lm}
\label{lm:existence_pullback_DVR}
Consider the hypothesis of \ref{fund_inter_hyp}
 in the case where $\alpha'=\spec O$ is a trait of $\alpha$.
Let $K$ be the fraction field of $O$ and $x$ the corresponding
 point on $\alpha$.

We suppose that $f$ is special at $x$.

Then the pre-special $\Rc[1/n]$-cycle $\beta \otimes_\alpha O$
 exists and is uniquely defined by the property
   $(\beta \otimes_\alpha O) \otimesf_O K=\beta_K$
    with the notations of the preceding lemma.
We also put $\beta_O:=\beta \otimes_\alpha O$.

The properties (P1) to (P5) are fulfilled and in addition : \\
(P6') For any extension $O'/O$ of discrete valuation rings,
 $\beta_{O'}=\beta_O \otimesf_O O'$.
\end{lm}
\begin{proof}
Remark that, with the notation of definition \ref{df:open_restriction},
 $\beta_O \otimesf_O K=\beta_O|_{\spec K}$.
For the first statement, we simply apply Lemma \ref{lm:open_restriction} 
 and put  $\beta_O=j_*(\beta_K)$
  where $j:X_K \rightarrow X_O$ is the canonical open immersion.

Then properties (P1), (P3), (P4), (P5) and (P6') of the case
 considered in this lemma follows easily from the uniqueness statement
  and the corresponding properties in the preceding lemma 
   (applying again \ref{lm:open_restriction}).

It remains to prove (P2). According to (P1), we reduce to the case
$\alpha=\langle S\rangle$ for a reduced scheme $S$.
We choose a birational morphism $p:\tilde S \rightarrow S$ such
that the proper transform $\tilde \beta$ is a Hilbert $\tilde S$-cycles.
Consider a diagram of the form $(*)$ in this case. 
According to property (P3), we can assume $R'=O$. \\
Remark the trait $\spec R \rightarrow S$ admits an extension
 $\spec R \rightarrow \tilde S$ as $p$ is proper.
The point $x$ admits an extension $K'/K$
 which lifts to a point $x':\spec{K'} \rightarrow \tilde S$
 -- again $\tilde S_K$ is a non empty algebraic scheme.
The discrete valuation corresponding to $O \subset K$ 
extends to a discrete valuation on $K'$ as $K'/K$ is finite. Let
$O' \subset K'$ be the corresponding valuation ring. 
The corresponding trait $\spec{O'} \rightarrow S$
 thus admits a lifting to $\tilde S$
  corresponding to the point $x'$ as $p$ is proper.
Considering a composite extension $E'/K$ of $K'/K$ and $E/K$,
 we have obtained a commutative diagram
$$
\xymatrix@R=6pt@C=18pt{
\spec{E'}\ar[r] & \spec{O'}\ar@{=}[d]\ar[r] & \spec R\ar[d] \\
& \spec{O'}\ar[r] & \tilde S \\
}
$$
which lifts our original diagram $(*)$. Let $x_1$ (resp. $x_2$)
 be the point $\spec E' \rightarrow \tilde S$ corresponding
  to the the composite through the upper way (resp. lower way)
   in the preceding diagram.
  
Then, $\beta_{R,E} \otimesf_E E'=\tilde \beta_{x_1}$.
Moreover, we get
$$
\xymatrix@=16pt{
(\beta \otimes_S O)_{O,E} \otimesf_E E'\ar@{=}^-{\ref{Hilbert&special}}[r]
 & (\beta \otimes_S O) \otimesf_O E'\ar@{=}^-{(P5)+(P6')}[rr]
 && (\tilde \beta \otimes_{\tilde S} {O'}) \otimesf_{O'} E'\ar@{=}^-{(P3)}[r]
 & \tilde \beta_{x_2}.
}
$$
By hypothesis, $\beta/\alpha$ is special at $\spec{K'} \rightarrow S$.
Thus Lemma \ref{lm:special&Hilbert_blow-up} concludes.
\end{proof}

\begin{thm}
\label{existence_pullback}
Consider the hypothesis of \ref{fund_inter_hyp}.

Assume $f$ is special at the generic points of $\alpha'$.

Then the pre-special $\Rc[1/n]$-cycle $\beta \otimes_\alpha \alpha'$
 exists and is uniquely determined by property (P2).

It satisfies all the properties (P1) to (P5). 
\end{thm}
\begin{proof}
According to Lemma \ref{lm:exists_fat_pts},
for any point $s$ of $S'$ with residue field $\kappa$, 
there exists an extension $E/\kappa$ and a fat point
$(R,E)$ (resp. $(R',E)$) of $\alpha$ (resp. $\alpha'$)
over $\spec E \rightarrow \alpha$ 
(resp. $\spec E \rightarrow \alpha'$).
The uniqueness statement follows
 by applying Lemma \ref{lm:basic_cycles}(1).

For the existence, we write 
$\alpha'=\sum_{i \in I} n_i.\acycl{Z_i}_{S'}$ in standard form.

For any $i \in I$, let $K_i$ be the function field of $Z_i$
and consider the canonical morphism $\spec{K_i} \rightarrow \alpha$.
Let $\beta_{K_i} \subset X_{K_i}$ be the $\Rc[1/n]$-cycle defined in lemma
\ref{lm:existence_pullback_field}.
Let $j_i:X_{K_i} \rightarrow X'$ be the canonical 
immersion and put:
\begin{equation}
\label{formula:external_product}
\beta \otimes_\alpha \alpha'=\sum_{i \in I} n_i.j_{i*}(\beta_{K_i}).
\end{equation}

Then properties (P1), (P3), (P4) and (P5) are direct
 consequences of this definition and of the corresponding
  properties of Lemma \ref{lm:existence_pullback_field}.

We check property (P2).
Given a diagram of the form $(*)$, there exists
a unique $i \in I$ such that $\spec{R'}$ dominates $Z_i$.
Thus we get for this choice of $i \in I$ that
 $(\beta \otimes_\alpha \alpha')_{R',E}
  =\big(j_{i*}(\beta_{K_i})\big)_{R',E}$.
Let $K'$ be the fraction field of $R'$
 and consider the open immersion $j':X_{K'} \rightarrow X_{R'}$.
The following computation then concludes:
$$
\xymatrix@=3pt{
\big(j_{i*}(\beta_{K_i})\big)_{R',E}\ar@{=}[r]
 & j'_*\big(j_{i*}(\beta_{K_i})_{K'}) \otimesf_{R'} E
        \ar@{=}^-{\ref{sp&open_ngb_gen_pt}}[rr]
 && j'_*(\beta_{K'}) \otimesf_{R'} E\ar@{=}^-{\ref{lm:existence_pullback_DVR}}[rr]
 && \beta_{R'} \otimesf_{R'} E\ar@{=}^-{\ref{lm:existence_pullback_DVR}(P2)}[rrr]
 &&& \beta_{R,E}.
}
$$
\end{proof}

\begin{df} \label{df:tensor_product_cycles}
In the situation of the previous theorem,
 we call the $\Rc[1/n]$-cycle $\beta \otimes_\alpha \alpha'$
  the pullback of $\beta/\alpha$ by $\alpha'$.
\end{df}

\begin{num} \label{additivity_relative_product}
By construction, the cycle $\beta \otimes_\alpha \alpha'$
 is bilinear with respect to addition of cycles in the following sense:
\begin{enumerate}
\item[(P7)] Consider the hypothesis of \ref{fund_inter_hyp}.
Let $\alpha'_1$, $\alpha'_2$ be cycles with domain $S'$ such
that $\alpha=\alpha'_1+\alpha'_2$. If $\beta/\alpha$ is special
at the generic points of $\alpha_1$ and $\alpha_2$, then
the following cycles are equal in $X \times_S S'$:
$$
\beta \otimes_\alpha (\alpha'_1+\alpha'_2)
 =\beta \otimes_\alpha \alpha'_1+\beta \otimes_\alpha \alpha'_2.
$$
\item[(P7')] Consider the hypothesis of \ref{fund_inter_hyp}.
Let $\beta_1$, $\beta_2$ be cycles with domain $X$ such
that $\beta=\beta_1+\beta_2$. If $\beta_1$ and $\beta_2$
are special over $\alpha$ at the generic points of $\alpha'$,
then $\beta/\alpha$ is special at the generic points of $\alpha'$
and the following cycles are equal in $X \times_S S'$:
$$
(\beta_1+\beta_2) \otimes_\alpha \alpha'
 =\beta_1 \otimes_\alpha \alpha'+\beta_2 \otimes_\alpha \alpha'.
$$
\end{enumerate}

In the theorem above, we can assume that $X$ (resp. $S$, $S'$)
 is the support of $\beta$ (resp. $\alpha$, $\alpha'$).
Thus the support of $\beta \otimes_\alpha \alpha'$ is included
 in $X \times_S S'$. More precisely:
 \end{num}
 
\begin{lm}
Consider the hypothesis of \ref{fund_inter_hyp} and
 assume that $X$ (resp. $S$, $S'$)
 is the support of $\beta$ (resp. $\alpha$, $\alpha'$).
Then, if $\beta/\alpha$ is special at the generic points
 of $\alpha'$, we obtain:
\begin{enumerate}
\item[(i)] Let $(X \times_S S')^{(0)}$ be the generic points of
$X \times_S S'$. Then, we can write
$$
\beta \otimes_\alpha \alpha'=\sum_{x \in (X \times_S S')^{(0)}} m_x.x
$$
\item[(ii)] For any generic point $x$ of $X \times_S S'$,
 if $m_x \neq 0$, the image of $x$ in $S'$ is a generic point $s'$ and
 the multiplicity of $s'$ in $\alpha'$ divides $m_x$ in $\Rc[1/n]$.
\end{enumerate}
\end{lm}
\begin{proof}
Point (ii) is just a traduction that $\beta \otimes_\alpha \alpha'$
is pre-special over $\alpha'$. For point (i), we reduce easily 
to the case where $\alpha$ is the scheme $S$ and $S$ is reduced.
We can also assume that $\alpha'$ is the spectrum of a field $k$.
It is sufficient to check point (i) after an extension of $k$.
Thus we can apply Lemma \ref{lm:Hilbertification} to reduce to
that case where $\beta$ is a Hilbert cycle over $S$.
This case is obvious.
\end{proof}

\begin{df} \label{SV_multiplicities}
In the situation of the previous lemma,
 we put $$m^{SV}(x;\beta \otimes_\alpha \alpha'):=m_x \in \Rc[1/n]$$
\index{notat}{mSVx@$m^{SV}(x;\beta \otimes_\alpha \alpha')$}
 and we call them the Suslin-Voevodsky multiplicities
\index{word}{multiplicity!Suslin-Voevodsky}
 (in the operation of pullback).
\end{df}

\begin{rem} \label{rem:SV_mult}
Consider the notations of the previous lemma:
\begin{enumerate}
\item Assume that $\alpha$ is the spectrum of a field $k$. 
Then the product $\beta \otimes_k \alpha'$ is always defined
 and agrees with the classical \emph{exterior product} (according to (P3)).
\item According to the previous lemma,
the irreducible components of $X \times_S S'$ which does not dominate
an irreducible component of $S'$ have multiplicity $0$: they correspond
to the "non proper components" with respect to the operation 
$\beta \otimes_\alpha \alpha'$.
\item Assume $\alpha' \rightarrow \alpha=\acycl{S' \xrightarrow p S}$,
$\beta=\sum_{i \in I} n_i.x_i$.
Let $y$ be a generic point of $X \times_S S'$ lying over a generic
point $s'$ of $S'$. Let $S'_0$ be the irreductible component of $S'$
corresponding to $s'$. Consider \emph{any} irreducible component $S_0$
of $S$ which contains $p(s')$ and let $\beta_0=\sum_i n_i.x_i$ 
where the sum runs over the indexes $i$ such that $x_i$ lies
over $S_0$. Then, according to \eqref{formula:external_product},
$$
m^{SV}(y;\beta \otimes_S \acycl{S'})
 =m^{SV}(y;\beta_0 \otimes_{S_0} \acycl{S'_0}).
$$
This is a key property of the Suslin-Voevodsky multiplicities
which explains why we have to consider the property that $\beta/\alpha$ 
is special at $s'$ (see \ref{special&Samuel} for a refined statement).
\end{enumerate}
\end{rem}

\begin{lm} \label{lm:special&specialisation}
Consider a morphism of cycles $\alpha' \rightarrow \alpha$
 and a pre-special morphism $f:\beta \rightarrow \alpha$ which is special at the
  generic points of $\alpha$.
Consider a commutative square
$$
\xymatrix@R=8pt@C=18pt{
\spec{k'}\ar^-{x'}[r]\ar[d] & \alpha'\ar[d] \\
\spec k\ar^-x[r] & \alpha
}
$$
such that $k$ and $k'$ are fields.
Then the following conditions are equivalent:
\begin{enumerate}
\item[(i)] $f$ is special at $x$.
\item[(ii)] $\beta \otimes_\alpha \alpha' \rightarrow \alpha'$ is special at $x'$.
\end{enumerate}
\end{lm}
\begin{proof}
This follows easily from Lemma \ref{lm:exists_fat_pts}
 and property (P2).
\end{proof}

%\begin{cor}
%Let $\alpha$ be an $S$-cycle which is special at a point $s \in S$.
%Let $Z$ be the closure of $s$ in $S$ and $t \in Z$ be a specialisation of $s$.
%Then the following conditions are equivalent:
%\begin{enumerate}
%\item[(i)] $\beta/S$ is special at $t$.
%\item[(ii)] $\beta \otimes_S Z/Z$ is special at $t$.
%\end{enumerate}
%\end{cor}

\begin{cor}\label{cor:special_base_change}
Let $f:\beta \rightarrow \alpha$ be a special morphism. \\
Then for any morphism $\alpha' \rightarrow \alpha$,
 $\beta \otimes_\alpha \alpha' \rightarrow \alpha'$ is special.
\end{cor}

%The next lemma follows easily from the construction:
%\begin{lm}
%Let $\beta \rightarrow \alpha$ be a morphism of cycles,
% and $\spec k \rightarrow \alpha$ be a point with image $s \in S$.
%
%The following conditions are equivalent:
%\begin{enumerate}
%\item For any point $s \in S$ with residue field $\kappa_s$,
%considering the canonical point 
%$\spec{\kappa_s} \rightarrow \alpha$, the cycle 
%$\beta_{\kappa_s}$ has coefficients in $\Rc$.
%\item For any morphism of cycles $\alpha' \rightarrow \alpha$, 
%$\beta \otimes_\alpha \alpha'$ has coefficients in $\Rc$.
%\end{enumerate}
%\end{lm}
\begin{df} \label{df:universal}
Let $f:\beta \rightarrow \alpha$ be a morphism of cycles
 and $x:\spec k \rightarrow \alpha$ be a point. \\
We say that $f$ is $\Rc$-universal
\index{word}{cycle!$\Rc$-universal (morphism of)}
 at $x$
 if it is special at $x$ and the cycle $\beta \otimes_\alpha k$
  has coefficients in $\Rc$.
\end{df}
In the situation of this definition,
 let $s$ be the image of $x$ in the support of $\alpha$,
  and $\kappa_s$ be its residue field.
Then according to (P6),
 $\beta_k=\beta_{\kappa_s} \otimesf_{\kappa_s} k$.
Thus $f$ is $\Rc$-universal at $x$
 if and only if it is $\Rc$-universal at $s$.
Furthermore, the following lemma follows easily:
\begin{lm}
Let $f:\beta \rightarrow \alpha$ be a morphism of cycles.
The following conditions are equivalent:
\begin{enumerate}
\item[(i)] For any point $s \in \overline \alpha$,
 $f$ is $\Rc$-universal at $s$.
\item[(ii)] For any point $x:\spec k \rightarrow \alpha$,
 $f$ is $\Rc$-universal at $x$.
\item[(iii)] For any morphism of cycles $\alpha' \rightarrow \alpha$,
 $\beta \otimes_\alpha \alpha'$ has coefficients in $\Rc$.
\end{enumerate}
\end{lm}

\begin{df}\label{df:universal_cycles}
We say that a morphism of cycles $f$ is $\Rc$-universal if it 
 satisfies the equivalent properties of the preceding lemma.
\end{df}
Of course, $\Rc$-universal morphisms are stable by base change.
These definitions will be applied similarly to morphisms of schemes
 by considering the associated morphism of cycles.

\begin{ex} \label{ex:universal&flat}
According to property (P3) of the pullback,
 a flat equidimensional morphism of schemes
 is $\Rc$-universal.
\end{ex}

\begin{num}\label{num:universal_reduce}
Let $\beta/\alpha$ be a morphism of $\Rc$-cycles.

Let $S$ be the support of $\alpha$ and consider the obvious
 morphism of cycles $S \rightarrow \alpha$.
 Recall from property (P1) of Paragraph \ref{fund_inter_hyp}
 that the cycle
$$
\beta_0:=\beta \otimes_\alpha S
$$
is the reduction of $\beta/\alpha$ (Definition \ref{df:pre-special}).
 This is a special $\Rc$-cycle over $S$
 (see Paragraph \ref{num:special&SV})

Moreover, it follows from the definition of the product that the
 following conditions are equivalent:
\begin{enumerate}
\item[(i)] $\beta/\alpha$ is $\Rc$-universal;
\item[(ii)] $\beta_0/S$ is $\Rc$-universal.
\end{enumerate}
In particular, condition (ii) appear in Lemma 3.3.9 of \cite{SV1}
 (with a restriction on the relative dimension that is not needed in
 fact).
\end{num}

\begin{rem}
Though Lemma 3.3.9 of \cite{SV1} does not give rise to any definition
 in \emph{loc. cit.}, it is central in the theory of Suslin and Voevodsky.
 In particular, it appears in the definition of the groups
 $z(X/S,r)$, $c(X/S,r)$,...
 that takes place right after Lemma 3.3.9.

Our definition has the advantage to:
\begin{itemize}
\item work properly over non reduced schemes;
\item have a local formulation (this is essential for the theorems
 of constructibility in subsection \ref{sec:constructibility});
\item being free of unnecessary assumptions such has the relative
 dimension of fibers (the integer $r$ that appear in $z(X/S,r)$).
\end{itemize} 
Besides, the categorical language introduced,
 obviously inspired by E.G.A., is very natural and will 
 prove to be useful in the treatment of finite correspondences
 (see for example the definition of the composition
 product, \ref{df:composition_corr}, and the short proof of
 the properties of this composition product, \ref{basic_propert_corr}).
\end{rem}

The following proposition shows that one can bound the denominators that can happen after 
 an arbitrary number of base changes.
\begin{prop}\label{prop:multiplicity_denominators_bounded}
Let $\beta/\alpha$ be a special morphism of $\Rc$-cycles.
 Then there exists an integer $N>0$ such that $N.\beta$ is $\Rc$-universal.
\end{prop}
\begin{proof}
According to Paragraph \ref{num:universal_reduce},
 one can reduce to the case where $\alpha$ is a reduced scheme $S$.
 We then prove by noetherian induction on $S$ the following assertion:
 for any closed subscheme $Z \subset S$, and any special $\Rc$-cycle $\alpha$ on $S$,
 there exists an integer $N>0$ such that $N.\alpha$ is $\Rc$-special

Take a special $\Rc$-cycle $\alpha$ on $S$.
 According to Lemma \ref{lm:Hilbertification}, there exists
 a birational morphism $p:\tilde S \rightarrow S$ such that
 the strict transform $\tilde \alpha$ of $\alpha$ along $p$
 is a Hilbert cycle, thus $\Rc$-universal.
 Let $U$ be a dense open subscheme of $S$ above which $p$ is an isomorphism.
 Thus for any point $s \in U$, with inverse image $t$ in $p^{-1}(U)$,
 we obtain that the cycle $\alpha_s=\tilde \alpha_t$ has $\Rc$-coefficients.

Let $Z$ be the complement of $U$ in $S$, with its reduced schematic structure.
 Then, by construction, the pullback $\alpha \otimes_S Z$ is an $\Rc[1/N]$-cycle.
 In particular $\alpha_0=N.\alpha \otimes_S Z$ is a special $\Rc$-cycle over $Z$.
 As $Z$ is a proper closed subscheme of $S$, we can apply to the Noetherian induction
 hypothesis to $Z$ and $\alpha_0$. We find an integer $N'>0$ such that $N'.\alpha_0$
 is $\Rc$-universal. By transitivity of pullbacks (which follows easily
 from the uniqueness statement of Theorem \ref{existence_pullback};
 see Proposition \ref{prop:transitivity1}),
 we thus obtain that $(NN').\alpha$ is $\Rc$-universal over $S$.
\end{proof}

%\begin{rem}
%This statement is obviously false without any finiteness assumption.
% Indeed, take any $\ZZ$-cycle $\alpha$ over a scheme $\alpha$ such that
% $\alpha$ is not universal
%\end{rem}

Recall that $\Rc$ is a sub-ring of the ring of rationals. One
 easily deduce from the preceding proposition the following result.
\begin{cor}\label{cor:multiplicity_denominators_bounded}
For any $\Rc$-cycle $\alpha$ special over a (noetherian) scheme $S$,
 there exists an integer $N>0$ such that $N.\alpha$
 is $\ZZ$-universal over $S$.
\end{cor}

%\begin{rem} \label{rem:additivity_tensor_cycles}
%Let $X$ be a scheme and $\beta$ be a cycle with domain $X$.
%Let $(X_\lambda)_{\lambda \in \Rc}$ be the connected components of $X$
%and $\beta_\lambda$ be the part of the cycle $\beta$ which is contained
%in $X_\lambda$.
%Then, a morphism of cycles $\beta \rightarrow \alpha$
%is special (resp. universal) if and only if 
%for any $\lambda$, the restriction $\beta_\lambda \rightarrow \alpha$
%is special (resp. universal).
%Moreover, considering a special morphism $\beta \rightarrow \alpha$
%and any morphism $\alpha' \rightarrow \alpha$, we check easily 
%the equality (as cycles of domain $X \otimes_{|\alpha|} |\alpha'|$):
%$$
%\beta \otimes_\alpha \alpha'
% =\sum_\lambda \beta_\lambda \otimes_\alpha \alpha'.
%$$
%\end{rem}

%%%%%%%%%%%%%%%%%%%%%%%%%%%%%%%%%%%%%%%%%%%%%%%%%%%%%%%%%%%%%%%%%%%
%
% The properties
%

\subsection{Intersection theoretic properties} \label{sec:intersection_theory}

\subsubsection{Commutativity}
\index{word}{cycle!pullback!commutativity}

\begin{lm}
\label{lm:pullback_special&dominant}
Consider morphisms of cycles with support in the left diagram
$$
\xymatrix@R=0pt@C=18pt{
 & \beta\ar[dd] 
 & & & X\ar^f[dd] \\
& & \subset & & \\
\gamma\ar[r]
 & \alpha & & T\ar^g[r] & S
}
$$
such that $\beta/\alpha$ is pre-special
 and $\gamma/\alpha$ is pseudo-dominant.

Assume
$$
\alpha=\sum_{i \in I} n_i.s_i, \ \
\beta=\sum_{j \in J} m_j.x_j, \ \
\gamma=\sum_{l \in H} p_l.t_l
$$
and denote by 
$\kappa_{s_i}$ (resp. $\kappa_{x_j}$, $\kappa_{t_l}$) the
residue field of $s_i$ (resp. $x_j$, $t_l$) in $S$
(resp. $X$, $T$).
Considering $(i,j,l) \in I \times J \times H$
 such that $f(x_j)=g(t_l)=s_i$,
we denote by $\nu_{j,l}:\spec{\kappa_{x_j} \otimes_{\kappa_{s_i}} \kappa_{t_l}} 
 \rightarrow X \times_S T$ the canonical immersion.

Then the following assertions hold:
\begin{enumerate}
\item[(i)] $\beta$ is special at the generic points of $\gamma$.
\item[(ii)] The cycle $\beta \otimes_\alpha \gamma$ has coefficients in $\Rc$.
\item[(iii)] The following equality of cycles holds
$$
\beta \otimes_\alpha \gamma=\sum_{i,j,l}
 \frac{m_j}{n_i}p_l.\nu_{j,l*}
  \big(\acycl{\spec{\kappa_{y_j} \otimes_{\kappa_{x_i}} \kappa_{z_l}}}\big)
$$
where the sum runs over $(i,j,l) \in I \times J \times H$
such that $f(x_j)=g(t_j)=s_i$.
\end{enumerate}
\end{lm}
\begin{proof}
Assertion (i) is in fact the first point of \ref{rm:trivial_special}.
Assertion (ii) follows from assertion (iii), which is a consequence
of the defining formula \eqref{formula:external_product}
and remark \ref{rm:pullback_dominant_pt}.
\end{proof}

\begin{cor}\label{cor:flat_pullback}
Let $g:T \rightarrow S$ be a flat morphism
 and $\beta=\sum_{j \in J} m_j.\acycl{Z_j}_X$ be a pre-special $S$-cycle
 written in standard form.

Then $\beta/S$ is pre-special at the generic points of $T$
 and 
$$
\beta \otimes_S \acycl T=\sum_{j \in J} m_j.\acycl{Z_j \times_S T}.
$$
\end{cor}

The pullback $\beta \otimes_\alpha \gamma$,
 at it is defined only when $\beta/\alpha$ is special,
 is in general non symmetric in $\beta$ and $\gamma$.
 However the previous lemma implies it is symmetric whenever it makes sense:
\begin{cor}
\label{commutativity}
Consider pre-special morphisms of cycles
$\beta \rightarrow \alpha$ and $\gamma \rightarrow \alpha$.

Then $\beta$ (resp. $\gamma$) 
 is special at the generic points of $\gamma$ (resp. $\beta$)
  and the following equality holds: 
$\beta \otimes_\alpha \gamma=\gamma \otimes_\alpha \beta$.
\end{cor}

\subsubsection{Associativity}
\index{word}{cycle!pullback!associativity}

\begin{prop}
\label{prop:transitivity1}
Consider morphism of cycles 
$\beta \xrightarrow f \alpha$, 
$\alpha'' \rightarrow \alpha' \rightarrow \alpha$ 
such that $f$ is special at the generic points of $\alpha'$
 and of $\alpha''$. 
Let $n$ be the exponential characteristic of $\alpha''$.

Then the following assertions hold:
\begin{enumerate}
\item[(i)] The relative cycle $(\beta \otimes_\alpha \alpha')/\alpha'$ is special 
at the generic points of $\alpha''$.
\item[(ii)] The cycle $(\beta \otimes_\alpha \alpha') \otimes_{\alpha'} \alpha''$
has coefficients in $\Rc[1/n]$.
\item[(iii)] 
$(\beta \otimes_\alpha \alpha') \otimes_{\alpha'} \alpha''
=\beta \otimes_\alpha \alpha''$.
\end{enumerate}
\end{prop}
\begin{proof}
Assertion (i) is a corollary of Lemma \ref{lm:special&specialisation}.
Assertion (ii) is in fact a corollary of assertion (iii),
 which in turn follows easily from the uniqueness statement in theorem
 \ref{existence_pullback}.
\end{proof}

\begin{lm}
Let $\gamma \xrightarrow g \beta \xrightarrow f \alpha$
be two pre-special morphisms of cycles with domains
$Y \rightarrow X \rightarrow S$.
Consider a fat point $(R,k)$ over $\alpha$ such that
 $\gamma/\beta$ is special at the generic points of
 $\beta_{R,k}$.
 
Then $\gamma/\alpha$ is pre-special and the following equality 
of cycles holds in $Y_k$:
$$\gamma_{R,k}=\gamma \otimes_\beta (\beta_{R,k}).$$
\end{lm}
\begin{proof}
The first statement is obvious.

We first prove: $\gamma_R=\gamma \otimes_\beta \beta_R$. \\
Remark that $\beta_R \rightarrow \beta$ is pseudo-dominant.
Thus $\gamma/\beta$ is special at the generic points of $\beta_R$ and the
right hand side of the preceding equality is well defined.
Moreover, according to Lemma \ref{lm:pullback_special&dominant},
we can restrict to the case where $\alpha=s$, $\beta=x$ and $\gamma=y$,
 with multiplicity $1$.
Let $\kappa_s$, $\kappa_x$, $\kappa_y$ be the corresponding respective
residue fields, and $K$ be the fraction field of $R$. \\
Then, according to \eqref{eq:comput_specialR},
 $\gamma_R=\acycl{\kappa_y \otimes_{\kappa_s} K}_{Y_R}$
 and $\beta_R=\acycl{\kappa_x \otimes_{\kappa_s} K}_{X_R}$.
But Lemma \ref{lm:pullback_special&dominant} implies
that $\gamma \otimes_\beta \beta_R=\acycl{\kappa_y \otimes_{\kappa_x} (\kappa_x \otimes_{\kappa_s} K)}_{X_R}$. Thus the associativity of the tensor product of fields
allows concluding.

From this equality and Proposition \ref{prop:transitivity1}, we deduce that:
$$
\gamma_R \otimes_{\beta_R} \beta_{R,k}
 =(\gamma \otimes_\beta \beta_R) \otimes_{\beta_R} \beta_{R,k}
 =\gamma \otimes_{\beta} \beta_{R,k}.
$$
Thus, the equality we have to prove can be written
$\gamma_R \otimesf_R k=\gamma_R \otimes_{\beta_R} (\beta_R \otimesf_R k)$
and we are reduced to the case $\alpha=\spec R$. \\
In this case, we can assume $\beta=\acycl X$ with $X$ integral.
Let us consider a blow-up $\tilde X \xrightarrow p X$ 
such that the proper transform $\tilde \gamma$ of $\gamma$ along $p$
is a Hilbert cycle over $\tilde X$ (\ref{lm:Hilbertification}).
We easily get (from (P3) and \ref{transitivity_Hilbert}) that
$$
\tilde \gamma_k
 =\tilde \gamma \otimes_{\tilde X} \acycl{\tilde X_k}.
$$
Let $Y$ (resp. $\tilde Y$) be the support of $\gamma$ 
 (resp. $\tilde \gamma$), 
$q:\tilde Y \rightarrow Y$ the canonical projection.
We consider the cartesian square obtained by pullback along 
 $\spec k \rightarrow \spec R$:
$$
\xymatrix{
\tilde Y_k\ar^{q_k}[r]\ar[d] & Y_k\ar[d] \\
\tilde X_k\ar^{p_k}[r] & X_k.
}
$$
As $X_k \subset X$ (resp. $Y_k \subset Y$) is purely of codimension $1$,
the proper morphism $p_k$ (resp. $q_k$) is still birational.
As a consequence, $q_{k*}(\tilde \gamma)=\gamma$.
Let $y$ be a point in $\tilde Y_k^{(0)} \simeq Y_k^{(0)}$
 which lies above a point $x$ in $\tilde X_k^{(0)} \simeq X_k^{(0)}$
Then, according to (P5) and using the notations of \ref{SV_multiplicities},
we get
$$
m^{SV}(y;\tilde \gamma \otimes_{\tilde X} \acycl{\tilde X_k})
=m^{SV}(y;\gamma \otimes_{X} \acycl{X_k}).
$$
This readily implies
$q_{k*}(\tilde \gamma \otimes_{\tilde X} \acycl{\tilde X_k})
=\gamma \otimes_{X} \acycl{X_k}$ and allows us to conclude.
\end{proof}

As a corollary of this lemma using the uniqueness statement in
Theorem \ref{existence_pullback}, we obtained:
\begin{cor}
\label{transitivity2}
Let $\gamma \xrightarrow g \beta \xrightarrow f \alpha$ be
pre-special morphisms of cycles.

Let $x:\spec k \rightarrow \alpha$ be a point.
If $\beta/\alpha$ is special (resp. $\Rc$-universal) at $x$ 
 and 
 $\gamma/\beta$ is special (resp. $\Rc$-universal)
 at the generic points of $\beta_k$, 
 then $\gamma/\alpha$ is special at $x$.

Let $\alpha' \rightarrow \alpha$ be any morphism of cycles 
 with domain $S' \rightarrow S$
 and $n$ be the exponential characteristic of $\alpha'$. 
Then, whenever it is well defined, 
the following equality of $\Rc[1/n]$-cycles holds:
$$
\gamma \otimes_\beta (\beta \otimes_\alpha \alpha')
=\gamma \otimes_\alpha \alpha'.
$$
\end{cor}

A consequence of the transitivity formulas is the associativity
of the pullback:
\begin{cor}
\label{associativity}
Suppose given the following morphisms of cycles
$$
\xymatrix@=10pt{
\alpha\ar[rd] & & \beta\ar^f[ld]\ar[rd] & & \gamma\ar^g[ld] \\
& \delta & & \sigma &
}
$$
such that $f$ and $g$ are pre-specials.

Then, whenever it is well defined,
the following equality of cycles hold:
$$
\gamma \otimes_\sigma (\beta \otimes_\delta \alpha)
=(\gamma \otimes_\sigma \beta) \otimes_\delta \alpha
$$
\end{cor}
\begin{proof}
Indeed, by the transitivity formulas \ref{prop:transitivity1}
and \ref{transitivity2}, both members of the equation are
equal to
$
(\gamma \otimes_\sigma \beta) \otimes_\beta
(\beta \otimes_\delta \alpha)$.
\end{proof}

\subsubsection{Projection formulas}
\index{word}{cycle!pullback!projection formulas}

\begin{prop}
\label{projection1}
Consider morphisms of cycles with support in the left diagram
$$
\xymatrix@R=0pt@C=18pt{
 & \beta\ar[dd] 
 & & & X\ar[dd] \\
& & \subset & & \\
\alpha'\ar[r]
 & \alpha & & S'\ar^q[r] & S
}
$$
such that $\beta/\alpha$ is special at the generic points of $\alpha'$. \\
Consider a factorization $S' \xrightarrow g T \rightarrow S$.

Then $\beta/\alpha$ is special at the generic points of
$g_*(\alpha)$ and the following equality of cycles holds in 
$X \times_S T$:
$$
\beta \otimes_\alpha g_*(\alpha')=(1_X \times_S g)_*(\beta \otimes_\alpha \alpha').
$$
\end{prop}
\begin{proof}
The first assuption is obvious.
By linearity, we can assume $S'$ is integral and $\alpha'$
is the generic point $s$ of $S'$ with multiplicity $1$.
Let $L$ (resp. $E$) be the residue field of $s$ (resp. $g(s)$).

Consider the pullback square
$
\xymatrix@C=10pt@R=10pt{
X_L\ar^{g_0}[r]\ar_j[d] & X_E\ar^i[d] \\
X \times_S S'\ar^/1pt/{g_X}[r] & X \times_S T
}
$
where $i$ and $j$ are the natural immersions.

Let $d$ be the degree of $L/E$ if it is finite and $0$ otherwise.
We are reduced to prove the equality
$g_{X*}(j_*(\beta_L))=d.i_*(\beta_E)$.
Using the functoriality of pushforward and property (P6),
it is sufficient to prove 
the equality $g_{0*}(\beta_E \otimesf_E L)=d.\beta_E$. 
If $d=0$, the morphism $g_0$ induces an infinite extension
of fields on any point of $X_L$ which concludes.
If $L/E$ is finite,
$g_0$ is finite flat and $\beta_E \otimesf_E L$ is the usual
pullback by $g_0$.
Then the needed equality follows easily (see \cite[1.7.4]{Ful}).
\end{proof}

\begin{lm} \label{lm:projection2}
Let $\beta \rightarrow \alpha$ be a pre-special morphism 
of cycles with domain $X \xrightarrow p S$.
Let $(R,k)$ a fat point over $\alpha$
 and $X \xrightarrow f Y \rightarrow S$ be a factorization of $p$.
Let $f_k$ be the pullback of $f$ over $\spec k$.

Suppose that the support of $\beta$ is proper with
respect to $f$.
Then $f_*(\beta)$ is pre-special over $\alpha$
and the equality of cycles
$\big(f_*(\beta)\big)_{R,k}=f_{k*}(\beta_{R,k})$ holds
in $Y_k$.
\end{lm}
\begin{proof}
As usual, considering the support $S$ of $\alpha$,
we reduce to the case where $\alpha=\acycl S$. Let $K$ be the fraction
field of $R$. As $\spec K$ maps to a generic point of $S$, 
we can assume $S$ is integral. Let $F$ be its function field.
We can assume by linearity that $\beta$ is a point $x$ in $X$
with multiplicity $1$. 

Let $L$ (resp. $E$) be the residue field of $x$ (resp. $y=f(x)$).
Let $d$ be the degree of $L/E$ if it is finite and $0$ otherwise.
Consider the following pullback square
$$
\xymatrix@R=10pt@C=20pt{
\spec{L \otimes_F K}\ar^/-6pt/j[r]\ar_{f_0}[d]
 & X \times_S \spec R=X_R\ar^{f_R}[d] \\
\spec{E \otimes_F K}\ar^/-6pt/i[r] & Y \times_S \spec R=Y_R.
}
$$
According to the formula \eqref{eq:comput_specialR}, we obtain:
\begin{align*}
f_{R*}(\beta_R)
 &=f_{R*}j_*(\acycl{L \otimes_F K})
 =i_*f_{0*}(\acycl{L \otimes_F K}) \\
& =i_*f_{0*}(f_0^*(\acycl{E \otimes_F K})=i_*(d.\acycl{E \otimes_F K})
 =\acycl{f_*(\beta)}_R.
\end{align*}

We are finally reduced to the case $S=\spec R$ and $\beta$ is
a Hilbert cycle over $\spec R$. Note that $f_*(\beta)$ is still
a Hilbert cycle over $\spec R$. 
As $\beta_{R,k}=\beta \otimesf_R k$, the result follows
 now from Proposition \ref{Hilbert_proj_formula}.
\end{proof}

\begin{cor}
\label{projection2}
Consider morphisms of cycles with support in the left diagram
$$
\xymatrix@R=0pt@C=18pt{
 & \beta\ar[dd] 
 & & & X\ar^p[dd] \\
& & \subset & & \\
\alpha'\ar[r]
 & \alpha & & S'\ar[r] & S
}
$$
such that $\beta/\alpha$ is special at the generic points of $\alpha'$
 (resp. $\Rc$-universal).

Consider a factorization $X \xrightarrow f Y \rightarrow S$
of $p$.

Suppose that the support of $\beta$ is proper with respect to $f$.
Then $f_*(\beta)/\alpha$ is special at the generic points of $\alpha'$
(resp. $\Rc$-universal)
and the following equality of cycles holds in $X \times_S S'$:
$$
(f \times_S 1_{S'})_*(\beta \otimes_\alpha \alpha')
=\big(f_*(\beta)\big) \otimes_\alpha \alpha'.
$$
\end{cor}

\subsection{Geometric properties}

\begin{num} \label{num:localization_cycles}
We introduce a notation which will often come
in the next section.
Let $S$ be a scheme
 and $\alpha=\sum_{i \in I} n_i.\acycl{Z_i}_X$ an $S$-cycle
 written in standard form.
 
Let $s$ be a point of $S$
 and $\spec k \xrightarrow{\bar s} S$ be a geometric point of $S$
 with $k$ separably closed.
Let $S'$ be one of the following local schemes:
the localization of $S$ at $s$,
the Hensel localization of $S$ at $s$,
the strict localization of $S$ at $\bar s$.

We then define the cycle with coefficients in $\Rc$
 and domain $X \times_S S'$ as:
$$
\alpha|_{S'}=\sum_{i \in I} n_i\acycl{Z_i \times_S S'}_{X \times_S S'}.
$$
\end{num}

\begin{rem} The canonical morphism $S' \rightarrow S$ is flat.
In particular, $\alpha/S$ is special at the generic points of $S'$
 and we easily get: $\alpha|_{S'}=\alpha \otimes_S {S'}$.
 \end{rem}

\subsubsection{Constructibility}\label{sec:constructibility}

\begin{df}
Let $S$ be a scheme and $s \in S$ a point.
We say that a pre-special $S$-cycle $\alpha$ is emph{trivial at $s$}
\index{word}{cycle!trivial}
 if it is special at $s$ and $\alpha \otimes_S s=0$.
\end{df}
Naturally, we say that $\alpha$ is trivial if it is zero.
Thus $\alpha$ is trivial if and only if it is trivial at
the generic points of $S$.

Recall from \cite[1.9.6]{EGA4} that an ind-constructible
\index{word}{indconstructible@ind-constructible}
subset of a noetherian scheme $X$ is a union of locally
closed subset of $X$.
\begin{lm} \label{ind-constructible}
Let $S$ be a noetherian scheme,
 and $\alpha/S$ be a pre-special cycle.
Then the set
$$
T=\big\{s \in S \mid
 \alpha/S \text{ is special (resp. trivial, $\Rc$-universal) at } s\big\}
$$
is ind-constructible in $S$.
\end{lm}
\begin{proof}
Let $s$ be a point of $T$, and $Z$ be its closure in $S$ with
its reduced subscheme structure. Put $\alpha_Z=\alpha \otimes_S Z$,
defined because $\alpha$ is special at the generic point of $Z$.
Given any point $t$ of $Z$, we know that $\alpha/S$ is special at
$t$ if and only if $\alpha_Z/Z$ is special at $t$ 
(\textit{cf.} \ref{lm:special&specialisation}). 
But there exists a dense open subset $U_s$ of $Z$ such that 
 $\alpha_Z|_{U_Z}$ is a Hilbert cycle over $U_Z$. Thus,
 $\alpha/S$ is special at each point of $U_s$ and $U_s \subset T$.
This concludes and the same argument proves the respective 
 statements.
\end{proof}

\begin{num} \label{projective_system}
Let $I$ be a left filtering category and
$(S_i)_{i \in I}$ be a projective system of noetherian schemes
\index{word}{projective system, of schemes}
 with affine transition morphisms.
We let $S$ be the projective limit of $(S_i)$ and we assume the followings:
\begin{enumerate}
\item $S$ is noetherian.
\item There exists an index $i \in I$ such that for all $j \geq i$,
 the canonical projection
 $S \xrightarrow{p_i} S_j$ is dominant.
\end{enumerate}
In this case, there exists an index $j/i$ such that 
for any $k/j$,
the map $p_k$ induces an isomorphism $S^{(0)} \rightarrow S_k^{(0)}$ 
on the generic points (\textit{cf.} \cite[8.4.2, 8.4.2.1]{EGA4}).
Thus, replacing $I$ by $I/j$, we can assume that this property is
satisfied for all index $i \in I$. As a consequence,
the following properties are consequences of the previous ones:
\begin{enumerate}
\item[(3)] For any $i \in I$, $p_i:S \rightarrow S_i$ is pseudo-dominant
 and $p_i$ induces an isomorphism $S^{(0)} \rightarrow S_i^{(0)}$.
\item[(4)] For any arrow $j \rightarrow i$ of $I$,
 $p_{ji}:S_j \rightarrow S_i$ is pseudo-dominant
 and $p_{ji}$ induces an isomorphism $S_j^{(0)} \rightarrow S_i^{(0)}$.
\end{enumerate}
\end{num}

\begin{prop} \label{lifting_ppty_special&universal}
Consider the notations and hypothesis above.
Assume we are given a projective system of cycles $(\alpha_i)_{i \in I}$
 such that $\alpha_i$ is a pre-special cycle over $S_i$
  and for any $j \rightarrow i$, $\alpha_j=\alpha_i \otimes_{S_i} S_j$.
Put $\alpha=\alpha_i \otimes_{S_i} S$ for an index $i \in I$.\footnote{The
pullback is well defined because of point (3) and (4) of the hypothesis
above.}

The following conditions are equivalent:
\begin{enumerate}
\item[(i)] $\alpha/S$ is special (resp. $\Rc$-universal).
\item[(ii)] There exists $i \in I$ such that 
 $\alpha_i/S_i$ is special (resp. $\Rc$-universal).
\item[(iii)] There exists $i \in I$ such that for all $j/i$,
 $\alpha_j/S_j$ is special (resp. $\Rc$-universal).
\end{enumerate}

Let $s$ be point of $S$ and $s_i$ its image in $S_i$.
Then the following conditions are equivalent:
\begin{enumerate}
\item[(i)] $\alpha/S$ is special (resp. $\Rc$-universal) at $s$.
\item[(ii)] There exists $i \in I$ such that 
 $\alpha_i/S_i$ is special (resp. $\Rc$-universal) at $s_i$.
\item[(iii)] There exists $i \in I$ such that for all $j/i$,
 $\alpha_j/S_j$ is special (resp. $\Rc$-universal) at $s_j$.
\end{enumerate}
\end{prop}
\begin{proof}
Let $\mathsf P$ be one of the respective properties:
``special'', ``trivial'', ``$\Rc$-universal''.
Using the fact that being $\mathsf P$ at $s$
is an ind-constructible property (from Lemma \ref{ind-constructible}),
it is sufficient to apply \cite[th. 8.3.2]{EGA4} to
 the following family of sets:
$$
F_i=\{s_i \in S_i \mid \text{ $\alpha_i$ satisfies $\mathsf P$ at $s_i$}\},
 \quad F=\{s \in S \mid \text{ $\alpha$ satisfies $\mathsf P$ at $s$}\}.
$$
To get the two sets of equivalent conditions of the statement from \emph{op. cit.}
we have to prove the following relations:
\begin{align*}
(1): \ & \forall (j \rightarrow i) \in \mathrm{Fl}(I),
  p_{ji}^{-1}(F_i) \subset F_j, \\
(2): \ & F=\cup_{i \in I} p_i^{-1}(F_i).
\end{align*}
We consider the case where $\mathsf P$ is the property ``special''.
For relation (1), we apply \ref{lm:special&specialisation} 
 which implies the stronger relation $p_{ji}^{-1}(F_i) = F_j$.
For relation (2), 
 another application of \ref{lm:special&specialisation} gives 
 in fact the stronger relation $F=p_i^{-1}(F_i)$ for any $i \in I$.

Consider a point $s_j \in S$ and put $s_i=p_{ji}(s_j)$.
Assume $\alpha_i$ is special at $s_i$. Then,
 applying \ref{prop:transitivity1} and (P3), we get:
\begin{equation} \label{eq1:proof_lifting}
\alpha_j \otimes_{S_j} s_j
 =(\alpha_i \otimes_{S_i} s_i) \otimesf_{\kappa(s_i)} \kappa(s_j).
\end{equation}
Similarly, given $s \in S_j$, $s_i=p_{i}(s)$,
 and assuming $\alpha_i$ is special at $s_i$, we get:
\begin{equation} \label{eq2:proof_lifting}
\alpha \otimes_S s=(\alpha_i \otimes_{S_i} s_i)
 \otimesf_{\kappa(s_i)} \kappa(s).
\end{equation}
 
We consider now the case where $\mathsf P$ is the property ``trivial''.
Then relation (1) follows from \eqref{eq1:proof_lifting}.
Relation (2) follows from \eqref{eq1:proof_lifting}
 and \ref{lm:basic_cycles}(1).

We finally consider the case $\mathsf P$ is the property ``$\Rc$-universal''.
Relation (1) in this case is again a consequence of \eqref{eq1:proof_lifting}.
According to \eqref{eq2:proof_lifting},
we get the inclusion $\cup_{i \in I} f_i^{-1}(F_i) \subset F$.
We have to prove the reciprocal inclusion. \\
Consider a point $s \in S$ with residue field $k$ such that $\alpha/S$
is $\Rc$-universal at $s$. For any $i \in I$, we put
$s_i=p_i(s)$ and denote by $k_i$ its residue field. 
It is sufficient to find an index $i \in I$ such
that $\alpha_i \otimes_{S_i} s_i$ has coefficients in $\Rc$.
Thus we are reduced to the following lemma:
\begin{lm}
Let $(k_i)_{i \in I^{op}}$ be an ind-field
 and put: $k=\ilim_{i \in I^{op}} k_i$. \\
Consider a family $(\beta_i)_{i \in I}$ such that
$\beta_i$ is a $k_i$-cycle of finite type with coefficients in $\QQ$
 and for any $j/i$, $\beta_j=\beta_i \otimesf_{k_i} k_j$.
 We put $\beta= \beta_i \otimesf_{k_i} k$.
 
If for an index $i \in I$, $\beta_i \otimesf_{k_i} k$
 has coefficients in $\Rc$,
 then there exists $j/i$ such that $\beta_j$ has coefficients in $\Rc$.
\end{lm}
We can assume that for any $j/i$, $\beta_j$ has positive coefficients.
Let $X_j$ (resp. $X$) be the support of $\beta_j$ (resp. $\beta$).
We obtain a pro-scheme $(X_j)_{j/i}$ such that
 $X=\plim_{i \in I} X_i$. The transition maps of $(X_j)_{j/i}$ are
 dominant. 
Thus, by enlarging $i$, we can assume that for any $j/i$, 
 the induced map $\pi_0(X_i) \rightarrow \pi_0(X_j)$ is a
 bijection. 
Thus we can consider each element of $\pi_0(X)$ separately
 and assume that all the $X_i$ are integrals:
 for any $j/i$, $\beta_j=n_j.\acycl{X_j}$ for a positive element $n_j \in \QQ$.
Arguing generically, we can further assume $X_j=\spec{L_j}$ 
for a field extension of finite type $L_j$ of $k_j$.
By assumption now, for any $j/i$, $L_i \otimes_{k_i} k_j$
 is an Artinian ring whose reduction is the field $L_j$.
Moreover, $n_j=n_i.\mathrm{lg}(L_i \otimes_{k_i} k_j)$
 and we know that $n:=n_i.\mathrm{lg}(L_i \otimes_{k_i} k)$ belongs to $\Rc$.

Let $p$ be a prime not invertible in $\Rc$ such that $v_p(n_i)<0$
where $v_p$ denotes the $p$-adic valuation on $\QQ$.
It is sufficient to find an index $j/i$ such that $v_p(n_j) \geq 0$.
Let $L=(L_i \otimes_{k_i} k)_{red}$. Remark that $L=\ilim_{i \in I^{op}} L_i$.
It is a field extension of finite type of $k$.
Consider elements $a_1,...,a_n$ algebraically independent over $k$
 such that $L$ is a finite extension of $k(a_1,...,a_n)$.
By enlarging $i$, we can assume that $a_1,...,a_n$ belongs to $L_i$.
Thus $L_i$ is a finite extension of $k_i(a_1,...,a_n)$:
replacing $k_i$ by $k_i(a_1,...,a_n)$, we can assume that $L_i/k_i$
is finite. \\
Let $L'$ be the subextension of $L$ over $k$
generated by the $p$-th roots of elements of $k$.
As $L/k$ is finite, $L'/k$ is finite,
 generated by elements $b_1,...,b_r \in L$.
consider an index $j/i$ such that $b_1,...,b_r$ belongs to $L_j$.
It follows that 
 $v_p(\mathrm{lg}(L_i \otimes_{k_i} k_j))=v_p(\mathrm{lg}(L_i \otimes_{k_i} k))$.
Thus $v_p(n_j)=v_p(n) \geq 0$ and we are done.
\end{proof}

\begin{cor} \label{cor:special_etale_local}
Let $S$ be a scheme and $\alpha$ be a pre-special $S$-cycle. \\
Let $\bar s$ be a geometric point of $S$, with image $s$ in $S$,
and $S'$ be the strict localization of $S$ at $\bar s$.

Then the following conditions are equivalent:
\begin{enumerate}
\item[(i)] $\alpha/S$ is special at $s$.
\item[(i')] $\alpha/S$ is special at $\bar s$.
\item[(ii)] $\big(\alpha|_{S'}\big)/S'$ is special at $\bar s$ 
(notation of \ref{num:localization_cycles}).
\item[(iii)] There exists an \'etale neighborhood $V$ of $\bar s$ in $S$
 such that $(\alpha \otimes_S V)/V$ is special at $\bar s$.
\end{enumerate}
\end {cor}
\begin{proof}
The equivalence of (i) and (i') follows trivially from definition (\textit{cf.} \ref{rm:trivial_special}). 
Recall from \ref{num:localization_cycles} that $\alpha|_{S'}=\alpha \otimes_S S'$.
Thus $(i') \Rightarrow (ii)$ is easy (see \ref{lm:special&specialisation}).
Moreover, $(ii) \Rightarrow (iii)$ is a consequence of the previous proposition
applied to the pro-scheme of \'etale neighborhood of $\bar s$.
Finally, $(iii) \Rightarrow (i)$ follows from Lemma \ref{lm:special&specialisation}.
\end{proof}

\begin{prop} \label{prop:lifting_cycles}
Consider the notations and hypothesis of \ref{projective_system}.
Assume that $S$ and $S_i$ are reduced for any $i \in I$. 

Suppose given a projective system $(X_i)_{i \in I}$
of $S_i$-schemes of finite type such that 
 for any $j/i$, $X_j=X_i \times_{S_i} S_j$.
We let $X$ be the projective limit of $(X_i)$.

Then for any pre-special (resp. special, $\Rc$-universal) 
$S$-cycle $\alpha \subset X$,
there exists $i \in I$ and a pre-special 
 (resp. special, $\Rc$-universal) $S_i$-cycle $\alpha_i \subset X_i$
 such that $\alpha=\alpha_i \otimes_{S_i} S$.\footnote{This
pullback is defined in any case because of point (3) of the hypothesis
above.}
\end{prop}
\begin{proof}
Using Proposition \ref{lifting_ppty_special&universal},
 we are reduced to consider the first of the respective cases of the proposition.
Write $\alpha=\sum_{r \in \Theta} n_r.\acycl{Z_r}_{X}$
 in standard form. \\
Consider $r \in \Theta$. As $X$ is noetherian, there exists an index $i \in I$
 and a closed subscheme $Z_{r,i} \subset X_i$ such that
 $Z_r=Z_{r,i} \times_{S_i} S$.
Moreover, replacing $Z_{r,i}$ by the reduced closure 
of the image of the canonical map $Z_r \xrightarrow{(*)} Z_{r,i}$,
we can assume that the map $(*)$ is dominant.
For any $j\in I/i$, we put $Z_{r,j}=Z_{r,i} \times_{S_i} S_j$.
The limit of the pro-scheme $(Z_{r,j})_{j \in I/i^{op}}$ is the
integral scheme $Z_r$.
Thus, applying \cite[8.2.2]{EGA4},
we see that by enlarging $i$, we can assume that for any $j \in I/i$,
$Z_{r,j}$ is irreducible (but not necessarily reduced). \\
We repeat this construction for every $r \in \Theta$,
 enlarging $i$ at each step.
Fix now an element $j \in I/i$.
The scheme $Z_{r,j}$ may not be reduced. However, its reduction $Z'_{r,j}$
is an integral scheme such that $Z'_{r,j} \times_{S_j} S=Z_r$.
We put
$$
\alpha_j=\sum_{r \in \Theta} n_r \acycl{Z'_{r,j}}_{X_j}.
$$
Let $z_{r,j}$ be the generic point of $Z'_{r,j}$,
 and $s_{r,j}$ be its image in $S_j$. It is a generic point
 and corresponds uniquely to a generic point $s_r$ of $S$
 according to the point (3) of the hypothesis \ref{projective_system}.
Thus $\alpha_j/S_j$ is pre-special.
Moreover, we get from the above that 
 $\kappa(z_{r,j}) \otimes_{\kappa(s_{r,j})} \kappa(s_r)=\kappa(z_r)$
 where $z_r$ is the generic point of $Z_r$.
Thus the relation $\alpha_j \otimes_{S_j} S=\alpha$ follows from lemma
 \ref{lm:pullback_special&dominant}.
\end{proof}

\subsubsection{Samuel multiplicities}

\begin{num} We give some recall on Samuel multiplicities,
\index{word}{multiplicity!Samuel (of a module)}
following as a general reference 
\cite[VIII.\textsection 7]{BourbakiAC8&9}. \\
Let $A$ be a noetherian local ring with maximal ideal $\fm$.
Let $M\neq 0$ be a $A$-module of finite type
 and $\fq \subset \fm$ an ideal of $A$ such that 
 $M/\fq M$ has finite length.
Let $d$ be the dimension of the support of $M$.
Recall from \emph{loc. cit.} that \emph{Samuel multiplicity} of $M$ 
at $\fq$ is defined as the integer:
$$
e_\fq^A(M):=\lim_{n \to \infty} \left(\frac{d!}{n^d}\lg_A(M/\fq^n M)\right)
$$
\index{notat}{eqAM@$e_\fq^A(M)$}
In the case $M=A$, 
 we simply put $e_\fq(A):=e^A_\fq(A)$ and $e(A):=e_\fm^A(A)$.

We will use the following properties of these multiplicities
that we recall for the convenience of the reader; let $A$ be
a local noetherian ring with maximal ideal $\fm$:

Let $\Phi$ be the generic points $\fp$ of $\spec A$ such that
 $\dim(A/\fp A)=\dim A$. Then 
 according to proposition 3 of \emph{loc. cit.}:
\begin{equation} \tag{$\cS1$} \label{Samuel_formula_irr_comp}
e_\fq(A)=\sum_{\fp \in \Phi} \lg(A_\fp).e_\fq(A/\fp).
\end{equation}

Let $B$ be a local flat $A$-algebra such that
$B/\fm B$ has finite length over $B$. 
Then according to proposition 4 of \emph{loc. cit.}:
\begin{equation} \tag{$\cS2$} \label{Samuel_formula_flat}
\frac{e_{\fm \! B}(B)}{e(A)}=\lg_B(B/\fm B).
\end{equation}

Let $B$ be a local flat $A$-algebra such that $\fm B$ is the maximal ideal of $B$. 
Let $\fq \subset A$ be an ideal such that $A/\fq A$ has finite length.
Then according to the corollary of proposition 4 in \emph{loc. cit.}:
\begin{equation} \tag{$\cS3$} \label{Samuel_inv_gonflement}
e_{\fq \! B}(B)=e_{\fq \!}(A).
\end{equation}

Assume $A$ is integral with fraction field $K$.
Let $B$ be a finite local $A$-algebra such that $B \supset A$.
Let $k_B/k_A$ be the extension of the residue fields of $B/A$.
Then, according to proposition 5 
and point b) of the corollary of proposition 4 in \emph{loc. cit.},
\begin{equation} \tag{$\cS4$} \label{Samuel_formula_finite_dom_ext}
\frac{e_{\fm \! B}(B)}{e(A)}
=\frac{\dim_K(B \otimes_A K)}{\lbrack k_B:k_A \rbrack}.
\end{equation}
\end{num}

\begin{df}
\begin{asparaenum}[(i)]
\item Let $S=\spec A$ be a local scheme,
 $s=\fm$ the closed point of $S$.

Let $Z$ be an $S$-scheme of finite type with special fiber $Z_s$.
For any generic point $z$ of $Z_s$,
denoting by $B$ the local ring of $Z$ at $z$,
we define the \emph{Samuel multiplicity of $Z$ at $z$ over $S$} as
the rational number:
$$
m^\cS(z,Z/S)=\frac{e_{\fm \! B}(B)}{e(A)}.
$$

In the case where $Z$ is integral,
 we define the \emph{Samuel specialization of the $S$-cycle $\acycl Z$ 
 at $s$} as the cycle with rational coefficients and domain $Z_s$:
$$
\acycl Z \otimes_S^\cS s
 =\sum_{z \in Z_s^{(0)}} m^\cS(z,Z/S).z.
$$

Consider an $S$-cycle of finite type
 $\alpha=\sum_{i \in I} n_i.\acycl{Z_i}_X$ written in standard form.
We define the \emph{Samuel specialization of the $S$-cycle $\alpha$ at $s$}
as the cycle with domain $X_s$:
$$
\alpha \otimes^\cS_S s=\sum_{i \in I} n_i.\acycl{Z_i} \otimes_S^\cS s.
$$
\index{notat}{alphatensorS@$\alpha \otimes^\cS_S s$}
\item Let $S$ be a scheme. For any point $s$ of $S$,
 we let $S_{(s)}$ be the localized scheme of $S$ at $s$.

Let $f:Z \rightarrow S$ be an $S$-scheme of finite type, 
 and $z$ a point of $Z$ which is generic in its fiber.
Put $s=f(z)$.
We define the \emph{Samuel multiplicity of $Z/S$ at $z$}
\index{word}{multiplicity!Samuel (of a cycle)}
 as the integer
$$
m^\cS(z,Z/S):=m^\cS(z,Z \times_S S_{(s)}/S_{(s)}).
$$

Consider an $S$-cycle of finite type
 $\alpha$ with domain $X$ and a point $s$ of $S$.
We define the
 \emph{Samuel specialization of the $S$-cycle $\alpha$ at $s$}
\index{word}{cycle!Samuel specialization}
 as the cycle with rational coefficients:
$$
\alpha \otimes^\cS_S s
 =\big(\alpha|_{S_{(s)}}\big) \otimes^\cS_{S_{(s)}} s.
$$
\end{asparaenum}
\end{df}

\begin{lm} \label{lm:Samuel_birat_nilp}
Let $S$ be a scheme, and $p:Z' \rightarrow Z$ an $S$-morphism
 which is a birational universal homeomorphism.
\index{word}{homeomorphism, universal}
Then for any point $s \in S$,
$$
\acycl{Z'} \otimes^\cS_S s=\acycl Z \otimes^\cS_S s
$$
in $(Z'_s)_{red}=(Z_s)_{red}$.
\end{lm}
\begin{proof}
By hypothesis, $p$ induces an isomorphism $Z'^{(0)} \simeq Z^{(0)}$
 between the generic points.
Given any irreducible component $T'$ of $Z'$ corresponding
 to the irreducible component $T$ of $Z$,
 we get by hypothesis: 
$$
T'_{red} \simeq T_{red} \text{ (as schemes), }
\lg\big(\mathcal O_{Z',T'}\big)=\lg\big(\mathcal O_{Z,T}\big).
$$
Thus, we easily conclude from the definition.
\end{proof}

\begin{num} Let $Z \xrightarrow f S$ be a morphism of finite type
 and a $z$ a point of $Z$, $s=f(z)$.
Assume $z$ is a generic point of $Z_s$.
We introduce the following condition:
$$
\cD(z,Z/S):\left\{
\begin{array}{ll}
\text{For any irreducible component $T$ of $Z_{(z)}$,} \\
\qquad T_s=\emptyset \text{ or } \dim(T)=\dim(Z_{(z)}).
\end{array}\right.
$$
\end{num}

\begin{rem} This condition is in particular satisfied if $Z_{(z)}$
 is absolutely equidimensional 
(and {\it a fortiori} if $Z$ is absolutely equidimensional).

An immediate translation of \eqref{Samuel_formula_irr_comp} gives:
\end{rem}

\begin{lm} \label{lm:Samuel&quasi-equidimensional}
Let $S$ be a local scheme with closed point $s$
 and $Z$ be an $S$-scheme of finite type such that $Z_s$
 is irreducible with generic point $z$.

If the condition $\cD(z,Z/S)$ is satisfied,
 then $\acycl Z \otimes^\cS_S s=m^\cS(z,Z/S).z$.
\end{lm}

We get directly from \eqref{Samuel_formula_flat}
 the following lemma:
\begin{lm} \label{lm:Samuel_formula_flat}
Let $S$ be a scheme, $s$ be a point of $S$,
 and $\alpha=\sum_{i \in I} n_i.\acycl{Z_i}_X$ be an $S$-cycle
 in standard form such that $Z_i$ is a flat $S$-scheme of finite type.

Then $\alpha$ is a Hilbert $S$-cycle and 
 $\alpha \otimes^{\cS}_S s=\alpha \otimesf_S s$.
\end{lm}

With the notations of \ref{num:localization_cycles},
we get from \eqref{Samuel_inv_gonflement}:
\begin{lm} \label{lm:Samuel_inv_gonflement}
Let $S$ be a scheme, $s$ a point of $S$ with residue field $k$
 and $\alpha$ an $S$-cycle of finite type.
\begin{asparaenum}[(i)]
\item Let $S'$ be the Hensel localization of $S$ at $s$.
Then,
$\alpha \otimes^{\cS}_S s=\big(\alpha|_{S'}\big) \otimes^{\cS}_{S'} s$.
\item Let $\bar k$ a separable closure corresponding and $\bar s$
  the corresponding geometric point of $S$.
Let $S_{(\bar s)}$ be the strict localization of $S$ at $\bar s$.
Then,
$$
\left(\alpha \otimes^{\cS}_S s\right) \otimesf_k \bar k
 =\big(\alpha|_{S_{(\bar s)}}\big) \otimes^{\cS}_{S_{(\bar s)}} \bar s.
$$
\end{asparaenum}
\end{lm}

%\num Given a scheme $X$ and a point $x \in X$,
% we put according to \cite[chap. 0, 14.1.1]{EGA4}:
%$$
%\dim_x(X)=\inf_{x \in U \subset X} \big(\dim(U)\big)
%$$
%where $U$ runs over the open neighbourhood of $x$ in $X$,
% and $\dim(U)$ is the dimension of $U$.
Let us recall from \cite[13.3.2]{EGA4} the following definition:
\begin{df} \label{df:relative_equidim}
Let $f:X \rightarrow S$ be a morphism of finite type between
noetherian schemes, and $x$ a point of $X$.
 
We say $f$ is equidimensional at $x$
 if there exists an open neighborhood $U$ of $x$ in $X$
  and a quasi-finite pseudo-dominant $S$-morphism
   $U \rightarrow \AA^d_S$ for $d \in \NN$.
 The integer $d$ is independent of the choice of $U$:
 it is called the relative dimension of $f$ at $x$.
 
We say $f$ is equidimensional if it is equidimensional
at every point of $X$.
\end{df}

\begin{rem} A quasi-finite morphism is equidimensional
 if and only if it is pseudo-dominant.
According to \cite[12.1.1.5]{EGA4}, this definition
 agrees with the convention stated in paragraph
 \ref{num:convention_flat_equidim}
 in the case of flat morphisms.
 \end{rem}

Note that a direct translation
 of \eqref{Samuel_formula_finite_dom_ext} gives:
\begin{lm} \label{lm:Samuel_formula_finite_dom_ext}
Let $S=\spec A$ be an integral local scheme
 with closed point $s$ and fraction field $K$.
Let $Z$ be a finite equidimensional $S$-scheme
 and $z$ a generic point of $Z_s$. Let $B$ be the local
 ring of $Z$ at $z$.

Then, 
$$
m^\cS(z,Z/S)=\frac{\dim_K(B \otimes_A K)}{\lbrack \kappa(z):\kappa(s)\rbrack}.
$$
\end{lm}

\begin{num}
Recall that a scheme $S$ is said to be
 \emph{unibranch (\emph{resp.} geometrically unibranch)
 at a point $s \in S$}
\index{word}{scheme!unibranch}
\index{word}{scheme!geometrically unibranch}
 if the henselisation (resp. strict henselisation)
 of the local ring $\mathcal O_{S,s}$ is irreducible
 (see \cite[6.15.1, 18.8.16]{EGA4}).
 The scheme $S$ is said to be \emph{unibranch}
 (resp. \emph{geometrically unibranch})
 if it is so at any point $s \in S$.

The following result is the key point
 of this subsection.
\end{num}
\begin{prop} \label{prop:computation_Samuel_mult}
Consider a cartesian square
$$
\xymatrix@=10pt{
Z'\ar^{g'}[r]\ar_{f'}[d] & Z\ar^f[d] \\
S'\ar^g[r] & S
}
$$
and a point $s'$ of $S'$, $s=g(s')$. 
Let $k$ (resp. $k'$) be the residue field of $s$ (resp. $s'$).
We assume the following conditions:
\begin{enumerate}
\item $S$ (resp. $S'$) is geometrically unibranch at $s$ (resp. $s'$).
\item $f$ and $f'$ are equidimensional of dimension $n$.
\item For any generic point $z$ of $Z_s$ (resp. $z'$ of $Z_{s'}$)
 the condition $\cD(z,Z/S)$ (resp. $\cD(z',Z'/S')$) is satisfied.
\end{enumerate}
Then, the following equality holds in $Z_{s'}$:
$$
\acycl{Z'} \otimes_{S'}^\cS s'=
(\acycl Z \otimes_S^\cS s) \otimesf_k k'.
$$
\end{prop}
\begin{proof}
According to Lemma \ref{lm:Samuel&quasi-equidimensional},
 we have to prove the equality:
\begin{equation} \label{proof:hard_Samuel_comput:rel}
\sum_{z' \in Z_{s'}^{(0)}} m^\cS(z',Z'/S').z'=
\sum_{z \in Z_s^{(0)}} m^\cS(z,Z/S).\acycl{\spec{\kappa(z) \otimes_k k'}}_{Z_{s'}}.
\end{equation}

As $f$ is equidimensional of dimension $n$, 
 we can assume according to \ref{df:relative_equidim} that there exists
 a quasi-finite pseudo-dominant $S$-morphism $p:Z \rightarrow \AA^n_S$.
For any generic point $z$ of $Z_s$, $t=p(z)$ is the generic point of
$\AA^n_s$. Thus applying \eqref{Samuel_inv_gonflement},
we get:
$$
m^\cS(z,Z/S)=m^\cS(z,Z/\AA^n_S).
$$

Consider the $S'$ morphism $p':Z' \rightarrow \AA^n_{Z'}$ obtained
by base change. It is quasi-finite. As $Z'/S'$ is equidimensional of
dimension $n$, $p'$ must be pseudo-dominant.
For any generic point $z'$ of $Z_{s'}$, $t'=p'(z')$ is the generic
point of $\AA^n_{s'}$ and as in the preceding paragraph, we get
$$
m^\cS(z',Z'/S')=m^\cS(z',Z'/\AA^n_{S'}).
$$

Moreover, the residue field $\kappa_t$ of $t$ (resp. $\kappa_{t'}$ of $t'$) 
is $k(t_1,...,t_n)$ (resp. $k'(t_1,...,t_n)$) and this implies
$\spec{\kappa(z) \otimes_{\kappa_t} \kappa_{t'}}$
is homeomorphic to $\spec{\kappa(z) \otimes_k k'}$
and has the same geometric multiplicities.
Putting this and the two preceding relations
in \eqref{proof:hard_Samuel_comput:rel}, we get reduced
to the case $n=0$ --
indeed, according to \cite[14.4.1.1] {EGA4},
 $\AA^n_S$ (resp. $\AA^n_{S'}$) is geometrically unibranch at $t$
 (resp. $t'$).

\bigskip

Assume now $n=0$, so that $f$ and $f'$ are quasi-finite pseudo-dominant.

Let $\bar k$ be a separable closure of $k$ and $\bar k'$ a separable
closure of a composite of $\bar k$ and $k'$.
It is sufficient to prove relation \eqref{proof:hard_Samuel_comput:rel}
after extension to $\bar k'$ (Lemma \ref{lm:basic_cycles}).
Thus according to \ref{lm:Samuel_inv_gonflement} and hypothesis (3),
we can assume $S$ and $S'$ are integral strictly local schemes.

For any $z \in Z_s^{(0)}$, 
 the extension $\kappa(z)/k$ is totally inseparable. Moreover, $z$ corresponds
 to a unique point $z' \in Z_{s'}^{(0)}$ and we have to prove for any
  $z \in Z^{(0)}_s$:
$$
m^\cS(z',Z'/S')=
m^\cS(z,Z/S).\lg(\kappa(z) \otimes_k k').
$$
Let $S=\spec A$, $K=\mathrm{Frac}(A)$ and $B=\mathcal O_{Z,z}$ 
 (resp. $S'=\spec{A'}$, $K'=\mathrm{Frac}(A')$ and $B'=\mathcal O_{Z',z'}$).
As $B$ is quasi-finite dominant over $A$ and $A$ is henselian, 
 $B/A$ is necessarily finite dominant. The same is true
 for $B'/A'$ and \eqref{Samuel_formula_finite_dom_ext} gives the formulas:
$$
m^\cS(z,Z/S)
=\frac{\dim_K(B \otimes_A K)}{\lbrack \kappa(z):k \rbrack},
 \qquad 
m^\cS(z',Z'/S')
=\frac{\dim_{K'}(B' \otimes_{A'} K')}{\lbrack \kappa(z'):k' \rbrack}.
$$
As $B' \otimes_{A'} K'=(B \otimes_{A} K) \otimes_K K'$, 
the numerator of these two rationals are the same.
To conclude, we are reduced to the easy relation
$$
\lbrack \kappa(z'):k' \rbrack.\lg(\kappa(z) \otimes_k k')
 =\lbrack \kappa(z):k \rbrack.
$$
\end{proof}

\begin{df}
Let $S$ be a scheme
 and $\alpha=\sum_{i \in I} n_i.\acycl{Z_i}_X$ be an $S$-cycle 
 in standard form.

We say $\alpha/S$ is \emph{pseudo-equidimensional over $s$}
\index{word}{cycle!pseudo-equidimensional}
 if it is pre-special and for any $i \in I$, 
 the structural map $Z_i \rightarrow S$
 is equidimensional at the generic points of the fiber $Z_{i,s}$.
\end{df}

\begin{prop}
Let $S$ be a strictly local integral scheme
\index{word}{scheme!strictly local}
 with closed point $s$ and residue field $k$
  and $\alpha$ be an $S$-cycle pseudo-equidimensional at $s$.

Then
 for any extension $\spec{k'} \xrightarrow{s'} S$ of $s$
 and any fat point $(R,k')$ of $S$ over $s'$,
the following relation holds:
$$
\alpha_{R,k'}=\left(\alpha \otimes^\cS_S s\right) \otimesf_k k'.
$$
\end{prop}
\begin{proof}
We put $S'=\spec R$ and denote by $s'$ its closed point. \bigskip \\ 
\textit{Reductions}.-- 
By additivity, we reduce to the case $\alpha=\acycl Z$,
 $Z$ is integral and the structural morphism $f:Z \rightarrow S$
 is equidimensional at the generic points of $Z_s$.
 Any generic point of $S'_{s'}$ dominate a generic
point of $Z_s$ so that we can argue locally
 at each generic point $x$ of $Z_s$.
Thus we can assume $Z_s$ is irreducible with generic point $x$.
Moreover, as $Z$ is equidimensional at $x$, we can assume
 according to \ref{df:relative_equidim} there exists
 a quasi-finite pseudo-dominant $S$-morphism 
\begin{equation} \label{proof:pseudo-finite}
Z \xrightarrow p \AA^n_S.
\end{equation}

Note that $S$ is geometrically unibranch at $s$.
Thus, applying \cite[14.4.1]{EGA4} ("crit\`ere de Chevalley"),
 $f$ is universally open at $x$.
As $S'$ is a trait whose close point goes to $s$ in $S$, 
 it follows from \cite[14.3.7]{EGA4}
 that the base change $f':Z' \rightarrow S'$ of $f$
  along $S'/S$ is pseudo-dominant.

Let $T$ be an irreducible component of $Z'$,
 with special fiber $T_{s'}$ and generic fiber $T_{K'}$
 over $S'$.
Then $T \rightarrow S'$ is a dominant 
morphism of finite type. Thus, according to \cite[14.3.10]{EGA4},
either $T_{s'}=\emptyset$ or $\dim(T_{s'})=\dim(T_{K'})$.
Moreover, the dimension of $T_\eta$ is equal
to the transcendental degree of the function field of $T$
over $K'$, which is equal to the transcendental degree of
$Z$ over $K$. This is $n$ according
 to \eqref{proof:pseudo-finite}.
Thus, in any case, $T$ is equidimensional of dimension $n$ over $S'$ 
and this implies $Z'$ is equidimensional
 of dimension $n$ over $S'$.
Moreover, either $T_{s'}=\emptyset$ or $\dim(T)=n+1=\dim(Z')$.
Note this implies that for any generic point $z'$ of $Z_{s'}$,
 the condition $\cD(z',Z'/S')$ is satisfied. \bigskip \\ 
\textit{Middle step}.-- 
We prove: $\alpha_{R,k}=\acycl{Z'} \otimes^\cS_{S'} s'$. \\
According to Lemma \ref{lm:Samuel_formula_flat},
$$
\alpha_{R,k}=\acycl{\overline{Z'_K}} \otimesf_R k'=
 \acycl{\overline{Z'_K}} \otimes^{\cS}_{S'} s'.
$$
But the canonical map $\overline{Z'_K} \rightarrow Z'$ is a birational
 universal homeomorphism so that we conclude this step by
  Lemma \ref{lm:Samuel_birat_nilp}.\bigskip \\ 
\textit{Final step}.-- We have only to point out that
the conditions of Proposition \ref{prop:computation_Samuel_mult}
 are fulfilled for the obvious square; this is precisely what we need.
\end{proof}

\begin{cor} \label{special&Samuel}
Let $S$ be a reduced scheme, $s$ a point of $S$
 and $\alpha$ an $S$-cycle which is pseudo-equidimensional over $s$.

Let $\bar s$ be a geometric point of $S$ with image $s$ in $S$
 and $S'$ be the strict localization of $S$ at $\bar s$.
We let $S'=\cup_{i \in I} S'_i$ be the irreducible
 components of $S'$ and $\alpha_i$ 
 be the cycle made by the part of the cycle $\alpha \otimes^\flat_S S'$
 whose points dominate $S'_i$.

Then the following conditions are equivalent:
\begin{asparaenum}[(i)]
\item $\alpha/S$ is special at $s$.
\item the cycle $\alpha_\lambda \otimes^\cS_{S'_i} \bar s$ does not depend
 on $i \in I$.
\end{asparaenum}
Moreover, when these conditions are fulfilled,
$\alpha \otimes_S \bar s=\alpha_\lambda \otimes^\cS_{S'_i} \bar s$.
\end{cor}
\begin{proof}
According to Corollary \ref{cor:special_etale_local},
 we reduce to the case $S=S'$. Then this follows directly from
 the preceding proposition.
\end{proof}

\begin{cor} \label{cor:geo_unibranch&special}
Let $S$ be a reduced scheme,
 geometrically unibranch at a point $s \in S$,
 and $\alpha$ an $S$-cycle. The following conditions
 are equivalent:
\begin{asparaenum}[(i)]
\item $\alpha/S$ is pseudo-equidimensional over $s$.
\item $\alpha/S$ is special at $s$.
\end{asparaenum}
Under these conditions, 
 $\alpha \otimes_S s=\alpha \otimes^\cS_S s$.
\end{cor}

\begin{rem} \label{rem:geo_unibranch&special}
In particular, over a reduced geometrically unibranch scheme $S$,
 every cycle whose support is equidimensional over $S$
 is special.
 \end{rem}
 
\begin{cor} \label{cor:universal&regular}
Let $S$ be a reduced scheme and $s \in S$ a point 
 such that $S$ is geometrically unibranch at $s$ 
 and $e(\mathcal O_{S,s})=1$.
Then for any $S$-cycle $\alpha$, the following conditions
are equivalent:
\begin{enumerate}
\item[(i)] $\alpha/S$ is pseudo-equidimensional over $s$.
\item[(ii)] $\alpha/S$ is $\Rc$-universal at $s$.
\end{enumerate}
\end{cor}

\begin{rem} \label{rem:regular&universal}
In particular, over a regular scheme $S$,
 every cycle whose support is equidimensional over $S$
 is $\Rc$-universal.
Remark also the following theorem:
\end{rem}

\begin{thm}
Let $S$ be an excellent scheme, $s \in S$ a point.
The following conditions are equivalent:
\begin{enumerate}
\item[(i)] $S$ is regular at $s$.
\item[(ii)] $S$ is geometrically unibranch at $s$ and $e(\mathcal O_{S,s})=1$.
\item[(iii)] $S$ is unibranch at $s$ and $e(\mathcal O_{S,s})=1$.
\end{enumerate}
\end{thm}
\noindent \textit{Bibliographical references for the proof}.
We can assume $S$ is the spectrum of an excellent local ring $A$
 with closed point $s$.
The implication $(i) \Rightarrow (ii)$ follows from the fact that 
a normal local ring is geometrically unibranch (at its closed point)
and from \cite[AC.VIII.\textsection 7, prop. 2]{BourbakiAC8&9}.
$(ii) \Rightarrow (iii)$ is trivial.
Concerning the implication $(iii) \Rightarrow (i)$,
 let $\hat A$ be the completion of the local ring $A$.
We know from \cite[AC.VIII.108, ex. 24]{BourbakiAC8&9} that 
when $e(A)=1$ and $\hat A$ is integral, $A$ is regular.
 Note $e(A)=1$ implies $A$ is reduced.
 To conclude, we refer to \cite[7.8.3, (vii)]{EGA4} which
 established that if $A$ is local excellent reduced,
  $\hat A$ is integral if and only if $A$ is unibranch.

\bigskip

Finally, we get the following theorem
 already proved by Suslin and Voevodsky
 (\cite[3.5.9]{SV1}):
\begin{thm} \label{th:mSV&tor}
Let $S$ be a scheme and $s$ a point with residue field $\kappa_s$
such that the local ring $A$ of $S$ at $s$ is regular.
Then for any equidimensional $S$-scheme $Z$ 
 and any generic point $z$ of $Z_s$,
$$
m^{SV}(z,\acycl Z \otimes_S s)
=\sum_i (-1)^i \mathrm{lg}_A \mathrm{Tor}_i^A(\mathcal O_{Z,z},\kappa_s).
$$
\end{thm}
\begin{proof}
We reduce to the case $S=\spec A$. Then $Z$ is absolutely equidimensional,
 and we can apply Lemma \ref{lm:Samuel&quasi-equidimensional}
 together with Corollary \ref{cor:geo_unibranch&special} to get that
 $m^{SV}(z,\acycl Z \otimes_S s)=m^\cS(z,Z/S)$.
Then the result follows from a theorem of Serre
  \cite[IV.12, th. 1]{Ser}.
\end{proof}

\begin{rem} \label{rem:intersection&Serre}
Let $S$ be a regular scheme,
\index{word}{scheme!regular}
 $X$ a smooth $S$-scheme
 and $\alpha \subset X$ an $S$-cycle whose support is
 equidimensional over $S$.
Let $s$ be a point of $S$ and $i:X_s \rightarrow X$ the
closed immersion of the fiber of $X$ at $s$.
Then the cycle $i^*(\alpha)$ of \cite[V-28, par. 7]{Ser}
is well-defined and we get:
$$
\alpha \otimes_S s=i^*(\alpha).
$$
\end{rem}

\section{Finite correspondences} \label{sec:finite_cor}

\begin{assumption} \label{num:cor_convention_class_P}
In this section, $\sch$ is the category of all noetherian schemes.
We fix an admissible class $\Pmor$ of morphisms
 in $\S$ and assume in addition that $\Pmor$ is contained
 in the class of separated morphisms of finite type.
%In this section,
% any relative scheme is assumed to have its structural
% morphism in $\Pmor$.

Consider two $S$-schemes $X$ and $Y$.
To clarify certain formulas,
 we will denote $X \times_S Y$ simply by $XY$
 and let $p_{XY}^X:XY \rightarrow X$ 
 be the canonical projection morphism.

\renewcommand{\Rc}{\Lambda}
We fix a ring of coefficients $\Rc \subset \QQ$.
\end{assumption}

\subsection{Definition and composition}

\begin{num} \label{num:c_0&pushout}
Let $S$ be a base scheme.
For any $\Pmor$-scheme $X/S$, we let $c_0(X/S,\Rc)$
\index{notat}{czeroXS@$c_0(X/S,\Rc)$}
 be the $\Rc$-module
 made of the finite and $\Rc$-universal $S$-cycles with domain $X$.\footnote{
  With the notations of \cite{SV1},
  $c_0(X/S,\ZZ)=c_{equi}(X/S,0)$ when $S$ is reduced.}
Consider a morphism $f:Y \rightarrow X$ of $\Pmor$-schemes over $S$.
Then the pushforward of cycles induces a well-defined morphism:
$$
f_*:c_0(Y/S,\Rc) \rightarrow c_0(X/S,\Rc).
$$
Indeed, consider a cycle $\alpha \in c_0(Y/S)$.
Let us denote by $Z$ its support in $Y$
 and by $f(Z) \subset X$ image of the latter
 by $f$. We consider these subsets as reduced subschemes.
Note that $f(Z)$ is separated and of finite type over $S$
 because $X/S$ is noetherian, separated, and of finite type, by assumption
 \ref{num:cor_convention_class_P}.
Because $Z/S$ is proper, \cite[5.4.3(ii)]{EGA2} shows that
 $f(Z)$ is indeed proper over $S$.
Thus, the cycle $f_*(\alpha)$ is $\Rc$-universal according to
 Corollary \ref{projection2}.
Finally, $Z/S$ is finite,
 we deduce that $f(Z)$ is quasi-finite, thus finite, over $S$.
 This implies the result.
\end{num}
\begin{df} \label{df:finite_cor}
Let $X$ and $Y$ be two $\Pmor$-schemes over $S$.

A finite $S$-correspondence
\index{word}{finite correspondence}
\index{word}{finite correspondence!finite $S$-correspondence|see{finite correspondence}}
 from $X$ to $Y$ with coefficients in $\Rc$ is an element of
$$
\corrg S X Y:=c_0(X \times_S Y/X).
$$
\index{notat}{cSXY@$\corrg S X Y$}
We denote such a correspondence by the symbol
$X \xdoto \alpha Y$.
\index{notat}{dotarrow@$\doto$}
\end{df}
In the case $\Rc=\ZZ$, we simply put
\renewcommand{\Rc}{\ZZ}
 $\corr S X Y:=\corrg S X Y$.
\renewcommand{\Rc}{\Lambda}
Through the rest of this section,
 unless explicitly stated,
 any cycle and any finite $S$-correspondence are assumed
  to have coefficients in $\Rc$.

\begin{rem} \label{rem:basic_ppties_corr}
\begin{enumerate}
\item According to properties (P7) and (P7') 
(\textit{cf.} \ref{additivity_relative_product})
of the pullback, $\corrg S X Y$ commutes with
finite sums in $X$ and $Y$.
\item Consider $\alpha \in \corrg S X Y$. 
Let $Z$ be the support of $\alpha$.
Then, 
 $Z$ is finite pseudo-dominant over $X$ (by definition \ref{df:pre-special}).
This means that $Z$ is finite equidimensional over $X$.

When $X$ is regular 
\index{word}{scheme!regular}
 (resp. $X$ is reduced geometrically unibranch
 \index{word}{scheme!geometrically unibranch}
  and $\car X\subset \Rc^\times$),
 a cycle $\alpha \subset X \times_S Y$
 written in standard form:
$$
\alpha=\sum_i n_i\acycl{Z_i}_{X \times_S Y}
$$
defines a finite $S$-correspondence from $X$ to $Y$ if and only if
 for any index $i \in I$, the scheme $Z_i$ is finite equidimensional over $X$
 (\emph{i.e.} finite and dominant over an irreducible component of $X$)
 -- \textit{cf.} \ref{rem:regular&universal} (resp. \ref{rem:geo_unibranch&special}).

Moreover, in each respective case, $c_S(X,Y)_\Rc$ is the
free $\Rc$-module generated by the closed integral subschemes $Z$
of $X \times_S Y$ which are finite equidimensional over $X$.
\item By definition, we get an inclusion:
$$
\corr S X Y \subset \corr S X Y_\Rc
$$
which induces an injective map:
$$
\corr S X Y \otimes_\ZZ\Rc \rightarrow \corr S X Y_\Rc.
$$
According to Corollary \ref{cor:multiplicity_denominators_bounded},
 this map is a bijection. Indeed, given any finite $\Rc$-linear $S$-correspondence
 $\alpha:X \doto Y$, applying the mentioned corollary,
 there exists an integer $N>0$ such that $N.\alpha$ is $\ZZ$-universal,
 so in particular an element of $\corr S X Y$. If we assume that $N$ is minimal,
 as $\alpha$ is $\Rc$-universal by assumption, $N$ must be invertible in $\Rc$.
 Therefore, $(N.\alpha) \otimes \frac 1 N$ belongs to $\corr S X Y \otimes_\ZZ\Rc$
 and is sent to $\alpha$ by the preceding map, which concludes.

Given more generally inclusions of rings
 $\Rc \subset \Rc' \subset \QQ$,
 we get an inclusion of groups
$$
\corr S X Y_\Lambda \subset \corr S X Y_{\Rc'}
$$
which induces an injection:
\begin{equation} \label{eq:cor_change_of_rings}
\corr S X Y_\Lambda \otimes_\Rc \Rc' \rightarrow \corr S X Y_{\Rc'}.
\end{equation}
Applying Proposition \ref{prop:multiplicity_denominators_bounded}
 and the same argument as above, we get that this map is in fact surjective,
 and therefore a bijection.
\end{enumerate}
\end{rem}

\begin{ex} \label{ex:cor_graph&transpose}
\begin{enumerate}
\item Let $f:X \rightarrow Y$ be a morphism in $\Pmorx S$.

Because $X/S$ is separated (assumption \ref{num:cor_convention_class_P}),
 the graph $\Gamma_f$ of $f$ is a closed subscheme of $X \times_S Y$. 
The canonical projection $\Gamma_f \rightarrow X$ is an isomorphism.
Thus $\acycl{\Gamma_f}_{XY}$ is a Hilbert cycle over $X$.
In particular, 
it is $\Rc$-universal and also finite over $X$,
thus it defines a finite $S$-correspondence from $X$ to $Y$.
\item Let $f:Y \rightarrow X$ be a finite $S$-morphism
which is $\Rc$-universal (as a morphism of the 
associated cycles).
Then the graph $\Gamma_f$ of $f$ is closed in $X \times_S Y$
and the projection $\Gamma_f \rightarrow X$ is isomorphic
to $f$.
Thus the cycle $\acycl{\Gamma_f}_{XY}$ is a finite 
$\Rc$-universal cycle over $X$ which therefore define
a finite $S$-correspondence $\tra f:X \doto Y$
\index{notat}{trasposef@$\tra f$}
called the \emph{transpose}
\index{word}{morphism!transpose}
\index{word}{finite correspondence!transpose|see{morphism}}
 of the finite $\Rc$-universal morphism $f$.
\end{enumerate}
\end{ex}

Suppose we are given finite $S$-correspondences
$X \xdoto \alpha Y \xdoto \beta Z$. Consider the
following diagram of cycles~:
\begin{equation}\label{eq:compose_cor}
\begin{split}
\xymatrix@=8pt{
\beta \otimes_Y \alpha\ar[r]\ar[d] & \beta\ar[r]\ar[d] & Z. \\
\alpha\ar[r]\ar[d] & Y & \\
X
}
\end{split}
\end{equation}
The pullback cycle is well-defined and has coefficients in $\Rc$
as $\beta$ is $\Rc$-universal over $Y$.
Moreover, according to the definition of pullback 
(\textit{cf.} \ref{existence_pullback})
and Corollary \ref{transitivity2},
$\beta \otimes_Y \alpha$ is a finite $\Rc$-universal
cycle over $X$ with domain $XYZ$.
Note finally that according to \ref{num:c_0&pushout},
 the pushforward of this latter cycle by $p_{XYZ}^{XZ}$
 is an element of $\corrg S X Z$.
\begin{df} \label{df:composition_corr}
Using the preceding notations,
we define the \emph{composition product}
\index{word}{finite correspondence!composition}
 of $\beta$ and $\alpha$ as the finite $S$-correspondence
$$
\beta \circ \alpha=p_{XYZ*}^{XZ}(\beta \otimes_Y \alpha)
:X \ldoto Z.
$$
\index{notat}{betacomposealpha@$\beta \circ \alpha$}
\end{df}

\begin{rem} \label{composition_reg=Tor}
In the case where $S$ is regular
\index{word}{scheme!regular}
 and $X$, $Y$, $Z$ are smooth over $S$,
 the composition product defined above agree with the one defined
  in \cite[4.1.16]{Deg7} in terms of the Tor-formula of Serre.
In fact, this is a direct consequence of \ref{th:mSV&tor}
 after reduction to the case where $\alpha$ and $\beta$
 are represented by closed integral subschemes 
 (see also point (2) of remark \ref{rem:basic_ppties_corr}).
\end{rem}
 
We sum up the main properties of the composition for finite 
correspondences in the following proposition~:
\begin{prop}\label{basic_propert_corr}
Let $X$, $Y$, $Z$ be $\Pmor$-schemes over $S$.
\begin{enumerate}
\item For any finite $S$-correspondences 
$X \xdoto \alpha Y \xdoto \beta Z \xdoto \gamma T$,
we have \\
\indent
$(\gamma \circ \beta) \circ \alpha=\gamma \circ (\beta \circ \alpha)$.
\item For any $X \xdoto \alpha Y \xrightarrow g Z$,
$\acycl{\Gamma_g}_{YZ} \circ \alpha=(1_X \times_S g)_*(\alpha)$.
\item For any $X \xrightarrow f Y \xdoto \beta Z$, 
$\beta \circ \acycl{\Gamma_f}_{XY}=\beta \otimes_Y \acycl X$.

\noindent Moreover, if $f$ is flat,
$\beta \circ \acycl{\Gamma_f}_{XY}=(f \times_S 1_Z)^*(\beta)$
considering the flat pullback of cycles in the classical sense.
\item For any $X \xleftarrow f Y \xdoto \beta Z$ such that
$f$ is finite $\Rc$-universal, \\
\indent $\beta \circ \tra f=(f \times_S 1_Z)_*(\beta)$.
\item For any $X \xdoto \alpha Y \xleftarrow g Z$ such that
$g$ is finite $\Rc$-universal, \\
\indent $\tra g \circ \alpha=\acycl Z \otimes_Y \alpha$.

\noindent If we suppose that $g$ is finite flat, then
$\tra g \circ \alpha=(1_X \times_S g)^*(\alpha)$.
\end{enumerate}
\end{prop}
\begin{proof}
(1) Using respectively the projection formulas \ref{projection2} 
and \ref{projection1}, we obtain
\begin{align*}
(\gamma \circ \beta) \circ \alpha
&=p_{XYZT*}^{XT}\big((\gamma \otimes_Z \beta) \otimes_Y \alpha\big) \\
\gamma \circ (\beta \circ \alpha)
&=p_{XYZT*}^{XT}\big(\gamma \otimes_Z (\beta \otimes_Y \alpha)\big).
\end{align*}
Thus this formula is a direct consequence of the associativity 
\ref{associativity}.

\bigskip

(2) Let $\epsilon:\Gamma_g \rightarrow Y$ 
and $p_{X\Gamma_g}^{XZ}:X\Gamma_g \rightarrow XZ$ be the canonical
projections. As $\epsilon$ is an isomorphism, we have tautologically
$\acycl Y=\epsilon_*(\acycl{\Gamma_g})$.
We conclude by the following
computation~:
\begin{align*}
(1_X \times_S g)_*(\alpha)
&=(1_X \times_S g)_*(\acycl Y \otimes_Y \alpha)
=(1_X \times_S g)_*(\epsilon_*\acycl{\Gamma_g} \otimes_Y \alpha) \\
&\stackrel{(*)}=(1_X \times_S g)_*(1_X \times_S \epsilon)_*
(\acycl{\Gamma_g} \otimes_Y \alpha)
=p_{X\Gamma_g*}^{XZ}(\acycl{\Gamma_g} \otimes_Y \alpha) \\
&\stackrel{(*)}=p_{XYZ*}^{XZ}(\acycl{\Gamma_g}_{YZ} \otimes_Y \alpha)
\end{align*}
The equalities labeled $(*)$ follow from the projection formula 
of \ref{projection2}.

(3) The first assertion follows from projection formula of 
\ref{projection1} and the fact that $\Gamma_f$ is isomorphic
to $X$~:
$$
\beta \circ \acycl{\Gamma_f}_{XY}
=p_{XYZ*}^{XZ}(\beta \otimes_Y \acycl{\Gamma_f}_{XY})
=\beta \otimes_Y p_{XY*}^X(\acycl{\Gamma_f}_{XY})
=\beta \otimes_Y \acycl{X}
$$
The second assertion follows from Corollary \ref{cor:flat_pullback}.

(4) and (5):
The proof of these assertions is strictly similar to that
of (2) and (3) instead that we use the projection formula
of \ref{projection1} (and do not need the commutativity
\ref{commutativity}).
\end{proof}

As a corollary, we obtain that the composition of $S$-morphisms
coincide with the composition of the associated graph considered
as finite $S$-correspondences. 
For any $S$-morphism $f:X \rightarrow Y$, we will still denote by
$f:X \doto Y$ the finite $S$-correspondence equal to 
$\acycl{\Gamma_f}_{XY}$.
Note moreover that for any $\Pmor$-scheme $X/S$, the identity
morphism of $X$ is the neutral element for the composition
of finite $S$-correspondences.
\begin{df} \label{df:Pmorc}
We let $\Pmorcx S$
\index{notat}{PcorlambdaS@$\Pmorcx S$}
 be the category of $\Pmor$-schemes over $S$
with morphisms the finite $S$-correspondences and the composition
product of definition \ref{df:composition_corr}.
\end{df}
An object of $\Pmorcx S$ will be denoted by $[X]$.
The category $\Pmorcx S$ is additive, and the direct sum is given
by the disjoint union of $\Pmor$-schemes over $S$.
We have a canonical faithful functor
\begin{equation} \label{eq:graph_functor}
\gamma:\Pmorx S \rightarrow \Pmorcx S
\end{equation}
which is the identity on objects and the graph on morphisms.
We call it the \emph{graph functor}.
\index{word}{functor!graph}
\index{word}{finite correspondence!graph functor|see{functor}}

\begin{num}
Given extension of rings $\Rc \subset \Rc' \subset \QQ$,
 we get according to Remark \ref{rem:basic_ppties_corr}(3)
 and the definition of composition of finite correspondences
 a functor of $\Rc'$-linear categories:
\begin{equation} \label{eq:Pcor_change_of_rings}
\Pmorcx S \otimes_\Rc \Rc'
 \rightarrow 
 \renewcommand{\Rc}{\Lambda'} \Pmorcx S
\end{equation}
which is the identity on objects and the maps
 of the form \eqref{eq:cor_change_of_rings} on morphisms.
 According to Remark \ref{rem:basic_ppties_corr}(3),
 the later maps are bijections and we get the following 
 result about changing coefficients.
\end{num}
\begin{prop} \label{prop:Pcor_change_of_rings}
With the notations above,
 the functor \eqref{eq:Pcor_change_of_rings} is an equivalence
 of categories.
\end{prop}

\begin{num} \label{num:comput_tra_circ_f}
Given two $S$-morphisms $f:Y \rightarrow X$ and $g:X' \rightarrow X$
such that $g$ is finite $\Rc$-universal, we get from the previous proposition
the equality of cycles in $YX'$:
$$
\tra g \circ f=\acycl{X'} \otimes_X \acycl{Y}_{YX}
$$
where $Y$ is seen as a closed subscheme of $YX$ through the
graph of $f$.

In particular, when either $f$ or $g$ is flat,
we get (use property (P3) of  \ref{fund_inter_hyp}
 or Corollary \ref{cor:flat_pullback}):
$$
\tra g \circ f=\acycl{X' \times_X Y}_{YX'}.
$$

To state the next formulas (the generalized degree formulas),
 we introduce the following notion:
 \end{num}
 
\begin{df} \label{df:degree_corr}
Let $f:X' \rightarrow X$ be a finite pseudo-dominant morphism
 (recall Definition \ref{df:cat_cycles}).
For any generic point $x$ of $X$,
 we define the degree of $f$ at $x$ as the integer:
$$
\deg_x(f)=\sum_{x'/x} [\kappa_{x'}:\kappa_x]
$$
\index{notat}{degreexf@$\deg_x(f)$}
where the sum runs over the generic points of $X'$ lying above $x$.
\end{df}

\begin{prop} \label{prop:degree_corr}
Let $X$ be a connected $S$-scheme
 and $f:X' \rightarrow X$ be a finite $S$-morphism.

If $f$ is special then there exists an integer $d \in \NN^*$
 such that for any generic point $x$ of $X$, $\deg_x(f)=d$.

Moreover, $f \circ \tra f=d.1_X$.
\end{prop}
We simply call $d$ the \emph{degree}
\index{word}{morphism!degree}
 of the finite special morphism $f$.
\begin{proof}
Let $\Delta'$ be the diagonal of $X'/S$.
For any generic point $x$ of $X$, we let $\Delta_x$ be the
 diagonal of the corresponding irreducible component of $X$,
 seen as a closed subscheme of $X$.
According to Proposition \ref{basic_propert_corr},
 and the definition of pushforwards, we get
$$
\alpha:=f \circ \tra f=(f \times_S f)_*(\acycl{\Delta'}_{X'X'})
 =\sum_{x \in X^{(0)}} \deg_x(f).\acycl{\Delta_x}_{XX}.
$$
Considering generic points $x$, $y$ of $X$,
 we prove $\deg_x(f)=\deg_y(f)$.
By induction,
 we can reduce to the case where $x$ and $y$ have a common
  specialization $s$ in $X$ because $X$ is connected and noetherian.
Then, as $\alpha/X$ is special,
  we get by definition of the pullback (see 
  more precisely \ref{rem:SV_mult})
$$
\alpha \otimes_S s=\deg_x(f).s=\deg_y(f).s
$$
as required. The remaining assertion then follows.
\end{proof}

\begin{prop}
Let $f:X' \rightarrow X$ be an $S$-morphism
 which is finite, radicial and $\Rc$-universal.

Assume $X$ is connected, and let $d$ be the degree of $f$.

Then $\tra f \circ f=d.1_{X'}$.
In particular, if $d$ is invertible in $\Rc$,
 $f$ is an isomorphism in $\Pmorcx S$.
\end{prop}
\begin{proof}
According to \ref{num:comput_tra_circ_f},
 $\tra f \circ f=\acycl{X'} \otimes_X \acycl{X'}$
 as cycles in $X'X'$.
Let $x$ be the generic point of $X$ and $k$ be its residue
field. Let $\{x'_i, i \in I\}$ be the set of generic
points of $X$, and for any $i\in I$, $k'_i$ be the residue
field of $x'_i$.
According to \ref{lm:pullback_special&dominant},
we thus obtain:
$$
\tra f \circ f
=\sum_{(i,j) \in I^2} \acycl{\spec{k'_i \otimes_k k'_j}}_{X'X'}.
$$
The result now follows by the definition of the degree
 and the fact that for any $i \in I$, $k'_i/k$ is radicial.
\end{proof}

\subsection{Monoidal structure}

Fix a base scheme $S$.
Let $X$, $X'$, $Y$, $Y$' be $\Pmor$-schemes over $S$.

Consider finite $S$-correspondences 
$\alpha:X \doto Y$ and $\alpha':X' \doto Y'$.
Then $\alpha X':=\alpha \otimes_X \acycl{XX'}$
and $\alpha' X:=\alpha' \otimes_{X'} \acycl{XX'}$ 
are both finite $\Rc$-universal cycles over $XX'$. 
Using stability by composition of finite $\Rc$-universal morphisms 
(\textit{cf.} Corollary \ref{transitivity2}), 
the cycle $(\alpha X') \otimes_{XX'} (\alpha' X)$ is
finite $\Rc$-universal over $XX'$.
\begin{df}
\label{df:tensor_corr}
Using the above notation, we define the \emph{tensor product}
\index{word}{tensor product!of finite correspondences}
\index{word}{finite correspondence!tensor porduct|see{tensor product}}
of $\alpha$ and $\alpha'$ over $S$ as the finite $S$-correspondence
$$
\alpha \otr_S \alpha'
=(\alpha X') \otimes_{XX'} (\alpha' X):
XX' \doto YY'.
$$
\index{notat}{alphatensortr@$\alpha \otr_S \alpha'$}
\end{df}
Let us first remark that this tensor product is commutative
(use commutativity of the pullback \ref{commutativity})
and associative (use associativity of the pullback
\ref{associativity}). 
Moreover, it is compatible with composition~:
\begin{lm}
Suppose given finite $S$-correspondences~: \\
$
\alpha:X \rightarrow Y, \
\beta:Y \rightarrow Z, \
\alpha':X' \rightarrow Y', \
\beta':Y' \rightarrow Z'
$.
Then
$$
(\beta \circ \alpha) \otr_S (\beta' \circ \alpha')
=(\beta \otr_S \beta') \circ (\alpha \otr_S \alpha').
$$
\end{lm}
\begin{proof}
We put
$\alpha X'=\alpha \otimes_X \acycl{XX'}$,
$\alpha' X=\alpha' \otimes_X \acycl{XX'}$
and 
$\beta Y'=\beta \otimes_Y \acycl{YY'}$,
$\beta' Y=\beta' \otimes_Y \acycl{YY'}$.
We can compute the right hand side of the above equation
as follows~:
\begin{align*}
&p_{XX'YY'ZZ'*}^{XX'ZZ'}\Big(
(\beta Y' \otimes_{YY'} \beta' Y)
\otimes_{YY'}
(\alpha X' \otimes_{XX'} \alpha' X) \Big) \\
&\stackrel{(1)}=p_{XX'YY'ZZ'*}^{XX'ZZ'}\Big(
(\beta Y' \otimes_{YY'} \beta' Y)
\otimes_{YY'}
(\alpha' X \otimes_{XX'} \alpha X') \Big) \\
&\stackrel{(2)}=p_{XX'YY'ZZ'*}^{XX'ZZ'}\Big(
\beta Y' \otimes_{YY'} ((\beta' Y
\otimes_{YY'}
\alpha' X) \otimes_{XX'} \alpha X') \Big) \\
&\stackrel{(3)}=p_{XX'YY'ZZ'*}^{XX'ZZ'}\Big(
(\beta Y' \otimes_{YY'} \alpha X') \otimes_{XX'}
(\beta' Y \otimes_{YY'} \alpha' X) ) \Big). \\
\end{align*}
Equality (1) follows from commutativity \ref{commutativity},
equality (2) from associativity \ref{associativity}
and equality (3) by both commutativity and associativity.

For the left hand side, we note that using the projection 
formula \ref{projection2},
the left hand side is equal to
$$
p_{XX'YY'ZZ'*}^{XX'ZZ'}\Big(
\big( (\beta \otimes_Y \alpha) \otimes_X \acycl{XX'}\big)
\otimes_{XX'}
\big( (\beta' \otimes_{Y'} \alpha') \otimes_{X'} \acycl{XX'}\big)\Big).
$$

We are left to remark that
$$
(\beta \otimes_Y \alpha) \otimes_X \acycl{XX'}
=\big((\beta Y') \otimes_{YY'} \alpha\big) \otimes_X \acycl{XX'}
=\beta Y' \otimes_{YY'} \alpha X',
$$
using transitivity \ref{prop:transitivity1} and associativity
\ref{associativity}. We thus conclude by symmetry of the 
other part in the left hand side.
\end{proof}

\begin{df} \label{df:pcor_tensor}
We define a symmetric monoidal structure on the category
$\Pmorcx S$ by putting
$[X] \otr_S [Y]=[X \times_S Y]$ on objects
and using the tensor product of the previous definition
for morphisms.
\end{df}

\begin{num} \label{num:gamma_monoidal}
Note that the functor $\gamma:\Pmorx S \rightarrow \Pmorcx S$
is monoidal for the cartesian structure
on the source category.
Indeed, this is a consequence of property (P3) of the relative
product (see \ref{fund_inter_hyp}) and the remark that
for any morphisms $f:X \rightarrow Y$ and $f':X' \rightarrow Y'$,
$(\Gamma_f \times_S X') \times_{XX'} (\Gamma_f' \times_S X)
=\Gamma_{f \times_S f'}$.
\end{num}

\subsection{Functoriality}

Fix a morphism of schemes $f:T \rightarrow S$. 
For any $\Pmor$-scheme $X/S$, we put $X_T=X \times_S T$.
For a pair of $\Pmor$-schemes over $S$ (resp. $T$-schemes) $(X,Y)$, 
we put $XY=X \times_S Y$ (resp. $XY_T=X \times_T Y$).

\subsubsection{Base change}

Consider a finite $S$-correspondence $\alpha:X \doto Y$.
The cycle $\alpha \otimes_X \acycl{X_T}$ defines
a finite $T$-correspondence from $X_T$ to $Y_T$
denoted by $\alpha_T$.
\begin{lm}
Consider finite $S$-correspondences
$X \xdoto \alpha Y \xdoto \beta Y$.

Then
$(\beta \circ \alpha)_T=\beta_T \circ \alpha_T$.
\end{lm}
\begin{proof}
This follows easily using 
the projection formula \ref{projection2},
the associativity formula \ref{associativity}
and the transitivity formula \ref{prop:transitivity1}~:
\begin{align*}
&p_{XYZ*}^{XZ}(\beta \otimes_Y \alpha) \otimes_X \acycl{X_T}
=p_{XYZ_T*}^{XZ_T}\big((\beta \otimes_Y \alpha) \otimes_X \acycl{X_T}\big) \\
&=p_{XYZ_T*}^{XZ_T}\big(\beta \otimes_Y (\alpha \otimes_X \acycl{X_T})\big)
=p_{XYZ_T*}^{XZ_T}\big((\beta \otimes_Y \acycl{Y_T}) \otimes_{Y_T}
(\alpha \otimes_X \acycl{X_T})\big).
\end{align*}
\end{proof}

\begin{df}
\label{df:pullback_cor}
Let $f:T \rightarrow S$ be a morphism of schemes.
Using the preceding lemma, we define the base change functor
$$
\begin{array}{rcl}
f^*:\Pmorcx S & \rightarrow & \Pmorcx T \\
{}[X/S] & \mapsto & [X_T/T] \\
\corrg S X Y \ni \alpha & \mapsto & \alpha_T.
\end{array}
$$
\end{df}

We sum up the basic properties of the base change for correspondences
in the following lemma.
\begin{lm} \label{lm:coherences_pullback_corr}
Take the notation and hypothesis of the previous definition.
\begin{enumerate}
\item The functor $f^*$ is symmetric monoidal.
\item Let $f^*_0:\Pmorx S \rightarrow \Pmorx T$ be the classical
base change functor on $\Pmor$-schemes over $S$. 
Then the following diagram is commutative:
$$
\xymatrix@R=10pt@C=20pt{
\Pmorx S\ar^/-2pt/{\gamma_S}[r]\ar_{f_0^*}[d] & \Pmorcx S\ar^{f^*}[d] \\
\Pmorx T\ar^/-2pt/{\gamma_T}[r] & \Pmorcx T.
}
$$
\item Let $\sigma:T' \rightarrow T$ be a morphism of schemes.
Through the canonical isomorphisms $(X_T)_{T'} \simeq X_{T'}$,
equality 
$(f \circ \sigma)^*=\sigma^* \circ f^*$ holds.
\end{enumerate}
\end{lm}
\begin{proof}
(1) This point follows easily using the associativity formula
\ref{associativity} and the transitivity formulas \ref{prop:transitivity1},
\ref{transitivity2}. \\
(2) This point follows from the fact that for any $S$-morphism 
$f:X \rightarrow Y$, there is a canonical isomorphism
$\Gamma_{f_T} \rightarrow \Gamma_f \times_S T$. \\
(3) This point is a direct application of the transitivity \ref{prop:transitivity1}.
\end{proof}

\begin{lm}
Let $f:T \rightarrow S$ be a universal homeomorphism. \\
Then $f^*:\Pmorcx S \rightarrow \Pmorcx T$ is fully faithful.
\end{lm}
\begin{proof}
Let $X$ and $Y$ be $\Pmor$-schemes over $S$. Then $X_T \rightarrow X$ is a universal
homeomorphism. Any generic point $x$ of $X$ corresponds uniquely to
a generic point of $X_T$. Let $m_x$ (resp. $m'_x$) be the geometric
multiplicity of $x$ in $X$ (resp. $X_T$).
Consider a finite $S$-correspondence $\alpha=\sum_{i \in I} n_i.z_i$.
For each $i\in I$, let $x_i$ be the generic point of $X$ dominated
 by $z_i$. Then we get by definition:
$$
f^*(\alpha)=\sum_{i \in I} m'_{x_i}\frac{n_i}{m_{x_i}}.z_i
$$
and the lemma is clear.
\end{proof}

\subsubsection{Restriction}

Consider a $\Pmor$-morphism $p:T \rightarrow S$.
For any pair of $T$-schemes $(X,Y)$, we denote by
$\delta_{XY}:X \times_T Y \rightarrow X \times_S Y$ the canonical
closed immersion deduced by base change from the
diagonal immersion of $T/S$.

Consider a finite $T$-correspondence $\alpha:X \doto Y$.
The cycle 
${\delta_{XY}}_*(\alpha)$ is the cycle $\alpha$ considered 
as a cycle in $X \times_S Y$. It defines a finite
$S$-correspondence from $X$ to $Y$.
\begin{lm}
\label{lm:res_cor_morph_lis}
Let $X$, $Y$ and $Z$ be $T$-schemes.
The following relations are true~:
\begin{enumerate}
\item For any $T$-morphism $f:X \rightarrow Y$,
  ${\delta_{XY}}_*\big(\acycl{\Gamma_f}_{XY_T}\big)=\acycl{\Gamma_f}_{XY}$.
\item For all $\alpha \in \corrg T X Y$ and $\beta \in \corrg T Y Z$,
$$
{\delta_{XZ}}_*(\beta \circ \alpha)=
({\delta_{YZ}}_*(\beta)) \circ
 ({\delta_{XY}}_*(\alpha)).
$$
\end{enumerate}
\end{lm}
\begin{proof}
The first assertion is obvious.

The second assertion is a consequence of the projection formulas 
\ref{projection1} and \ref{projection2}, and the functoriality
 of pushforwards~:
\begin{align*}
({\delta_{YZ}}_*(\beta)) \circ
 ({\delta_{XY}}_*(\alpha))
 &=p_{XYZ*}^{XZ}
\big({\delta_{YZ}}_*(\beta) \otimes_Y {\delta_{XY}}_*(\alpha)\big) \\
 &=p_{XYZ*}^{XZ} \delta_{XYZ*}(\beta \otimes_Y \alpha)
 =\delta_{XZ*} p_{XYZ_T*}^{XZ_T}(\beta \otimes_Y \alpha).
\end{align*}
\end{proof}

\begin{df}
\label{df:forget_base_corr}
Let $p:T \rightarrow S$ be a $\Pmor$-morphism.

Using the preceding lemma, we define a functor
$$
\begin{array}{rcl}
p_\sharp:\Pmorcx T & \rightarrow & \Pmorcx S \\
{}[X \rightarrow T] & \mapsto & [X \rightarrow T \xrightarrow p S] \\
\corrg T X Y \ni \alpha & \mapsto & \delta_{XY*}(\alpha).
\end{array}
$$
\end{df}

%\begin{num}
%\label{restriction_coincide&functoriality}
%Note that from the first point of the preceding lemma, the restriction
%of $p_\sharp$ to $\Pmorx T$ is the classical functor forgetting the base.
%Moreover, for a sequence
%$R \xrightarrow \sigma T \xrightarrow p S$ of morphisms in $\S$,
%we have $(p \circ \sigma)_\sharp
% = p_\sharp \circ \sigma_\sharp$.
%\end{num}
This functor enjoys the following properties:
\begin{lm}
\label{lm:forget_base_corr}
Let $p:T \rightarrow S$ be a $\Pmor$-morphism.
\begin{enumerate}
\item The functor $p_\sharp$ is left adjoint to the functor
$ p^*$.
\item For any composable $\Pmor$-morphisms $Z \xrightarrow q T \xrightarrow p S$,
 $(pq)_\sharp=p_\sharp q_\sharp$.
\item  Let $p_\sharp^0:\Pmorx T \rightarrow \Pmorx S$ be the functor
 induced by composition with $p$. 
Then the following diagram is commutative:
$$
\xymatrix@R=10pt@C=20pt{
\Pmorx T\ar^/-2pt/{\gamma_T}[r]\ar_{p_\sharp^0}[d] & \Pmorcx T\ar^{p_\sharp}[d] \\
\Pmorx S\ar^/-2pt/{\gamma_S}[r] & \Pmorcx S.
}
$$
% For any $T$-scheme $X$ and $S$-scheme $Y$, the
%obvious morphism obtained by adjunction
%$$
% p_\sharp( p^* [X] \otr_T [Y])
% \rightarrow [X] \otr_S  p_\sharp([Y])
%$$
%is an isomorphism.
\end{enumerate}
\end{lm}
\begin{proof}
For point (1), we have to construct for schemes $X/T$ and $Y/S$ a natural isomorphism
$\corrg S {p_\sharp X} Y \simeq \corrg T X {p^* Y}$. It is
induced by the canonical isomorphism of schemes
$( p_\sharp X) \times_S Y \simeq X \times_T ( p^* Y)$. \\
Point (2) follows from the associativity of the pushforward functor on cycles.
Note also that this identification is compatible with the transposition of
the identification of \ref{lm:coherences_pullback_corr}(3) according to
the adjunction property just obtained. \\
Point (3) is a reformulation of \ref{lm:res_cor_morph_lis}(2).
\end{proof}

\subsubsection{A finiteness property}

\begin{num} We assume here that $\Pmor$ is the class of
 all separated morphisms of finite type in $\sch$.

Let $I$ be a left filtering category and
$(X_i)_{i \in I}$ be a projective system
 of separated $S$-schemes of finite type
  with affine dominant transition morphisms.
We let $\cX$ be the projective limit of $(X_i)_i$
 and assume that $\cX$ is Noetherian over $S$.
 \end{num}
 
\begin{prop} \label{prop:finite_cor&limit}
Let $Y$ be a $\Pmor$-scheme of finite type over $S$.
Then the canonical morphism
$$
\varphi:\ilim_{i \in I^{op}} \corrg S {X_i} Y
 \rightarrow c_0(\cX \times_S Y/\cX,\Rc).
$$
is an isomorphism.
\end{prop}
\begin{proof}
Note that according to \cite[IV, 8.3.8(i)]{SGA4},
 we can assume the conditions (2) of \ref{projective_system}
 is verified for $(X_i)_{i \in I}$. Thus conditions (1) to (4) of
 \emph{loc. cit.} are verified.
Then the surjectivity of $\varphi$ follows from \ref{prop:lifting_cycles}
 and the injectivity from \ref{lifting_ppty_special&universal}.
\end{proof}

\subsection{The fibred category of correspondences}

We can summarize the preceding constructions:
\begin{prop} \label{prop:pfibred_corr}
The $2$-functor
$$
\Pmorc:S \mapsto \Pmorcx S
$$
equipped with the pullback defined in \ref{df:pullback_cor}
 and with the tensor product of \ref{df:pcor_tensor}
 is a monoidal $\Pmor$-fibred category
\index{word}{fibred!monoidalPfibredcategory@monoidal $\Pmor$-fibred category!of finite correspondences}
  such that the functor
$$
\gamma:\Pmor \rightarrow \Pmorc
$$
 (see \eqref{eq:graph_functor}) is a morphism
 of monoidal $\Pmor$-fibred category.
\end{prop}
\begin{proof}
According to Lemma \ref{lm:forget_base_corr},
 for any $\Pmor$ morphisms $p$, $p^*$ admits a left adjoint
 $p_\sharp$.
We have checked that $\gamma$ is symmetric monoidal
 and commutes with $f^*$ and $p_\sharp$ (see respectively 
 \ref{num:gamma_monoidal}, \ref{lm:coherences_pullback_corr}
  and \ref{lm:forget_base_corr}).
But $\gamma$ is essentially surjective. Thus, to prove the properties
 \bc and \pf for the fibred category $\Pmorc$, we are reduced to the
 case of $\Pmor$ which is easy
 (see example \ref{ex:can_mon_weak_P-fibred&bc}).
This concludes.
\end{proof}

\begin{rem} Consider the definition above.
\begin{enumerate}
\item The category $\Pmorc$ is $\Rc$-linear.
 For any scheme $S$, $\Pmorcx S$ is additive.
For any finite family of schemes $(S_i)_{i \in I}$
 which admits a sum $S$ in $\sch$, the canonical map
$$
\Pmorcx S \rightarrow \bigoplus_{i \in I} \Pmorcx{S_i}
$$
is an isomorphism.
\item The functor $\gamma:\Pmor \rightarrow \Pmorc$
 is nothing else than the canonical geometric sections
 of $\Pmorc$ (see definition \ref{df:geometric_sections}).
 \end{enumerate}
\end{rem} 

We will apply these definitions in the particular cases $\Pmor=\sm$ 
 (resp. $\Pmor=\sft$)
the class of smooth separated (resp. separated) morphisms of finite type.
Note that we get a commutative square
$$
\xymatrix@C=30pt@R=15pt{
\sm\ar^-\gamma[r]\ar[d] & \smc\ar[d] \\
\sft\ar^-\gamma[r] & \sftc
}
$$
where the vertical maps are
 the obvious embeddings of monoidal $\sm$-fibred categories.

\begin{num}
Consider extensions of rings $\Rc \subset \Rc' \subset \QQ$.
The functors \eqref{eq:Pcor_change_of_rings} for various
 schemes $S$ in $\sch$ are compatible with the operations
 of a $\Pmor$-fibred category because it is just concerned
 with adding denominators in the coefficients of the
 finite correspondences considered.
 Thus they induce a morphism of monoidal $\Pmor$-fibred categories
 over $\sch$:
\begin{equation} \label{eq:Pcor_change_of_rings_bis}
\Pmorc \otimes_\Rc \Rc'
 \rightarrow 
 \renewcommand{\Rc}{\Lambda'} \Pmorc.
\end{equation} 
According to Proposition \ref{prop:Pcor_change_of_rings},
 we immediately get the following result:
\end{num}
\begin{prop}\label{prop:Pcor_change_of_rings_bis}
Then the morphism \eqref{eq:Pcor_change_of_rings_bis}
 is an equivalence of monoidal $\Pmor$-fibred categories.
\end{prop}

\begin{rem} \label{rem:Pcor&regular}
\renewcommand{\Rc}{\ZZ}
The restriction of the category $\Pmorc $ to the
 category of regular schemes was already defined
 in \cite{Deg7}. Indeed, one can check using the
 comparison of Suslin-Voevodsky's multiplicities
 with Serre's intersection multiplicities (using
 Tor-formulas ; \textit{cf.} \ref{th:mSV&tor}),
 that the operations $\tau^*$, $\tau_\sharp$, and $\otr$
 defined here coincide with that of \cite{Deg7}.
 This remark extends \ref{composition_reg=Tor}.
\end{rem}

%%% Local variables:
%%% mode: latex
%%% tex-main-file: "main"
%%% TeX-master: "main"
%%% ispell-local-dictionary: "english"
%%% end:

%%%%%%%%%%%%%%%%%%%%%%%%%%%%%%%%%%%%%%%%%%%%%%%%%%%%%%%%%%%%%%%%%%
% Transferts
%%%%%%%%%%%%%%%%%%%%%%%%%%%%%%%%%%%%%%%%%%%%%%%%%%%%%%%%%%%%%%%%%%

\section{Sheaves with transfers} \label{sec:fxtr}

%%%%%%%%%%%%%%%%%%%%%%%%%%%%%%%%%%%%%%%%%%%%%%%%%%%%%%%%%%%%%%%%%
%
% Définition préfaisceaux avec transferts
%
%
\renewcommand{\Rc}{\Lambda}
\begin{assumption} \label{num:fxtr_convention_class_P}
The category $\sch$ is the category of noetherian schemes
 of finite dimension.
We fix an admissible class $\Pmor$ of morphisms in $\sch$ 
 satisfying the following assumptions:
\begin{enumerate}
\item[(a)] Any morphism in $\Pmor$ is separated of finite type.
\item[(b)] Any \'etale separated morphism of finite type is in $\Pmor$.
\end{enumerate}
We fix a topology $t$ on $\S$ which is $\Pmor$-admissible
 and such that:
\begin{enumerate}
\item[(c)] For any scheme $S$,
 there is a class of covers $\mathcal E$ of the form $(p:W \rightarrow S)$
 with $p$ a $\Pmor$-morphism such that $t$ is the topology
 generated by $\mathcal E$ and the covers of the form
 $(U \rightarrow U \sqcup V,V \rightarrow U \sqcup V)$
 for any schemes $U$ and $V$ in $\sch$.
\end{enumerate}
We fix a ring of coefficients $\Rc$.
 Whenever we speak of $\Rc$-cycles 
 (or the premotivic category $\Pmorc$, etc...),
 we mean cycles with coefficients in the localization of $\ZZ$
 with respect to invertible primes in $\Rc$.
\end{assumption}

Note that in sections \ref{sec:examples_ftr} and \ref{sec:comparison_ftr},
 we will apply the conventions of
 section \ref{sec:premotivic_convention} by replacing the class
 of smooth morphisms of finite type (resp. morphisms of finite type)
 there by the class of smooth separated morphisms of finite type
 (resp. separated morphisms of finite type).

\subsection{Presheaves with transfers}

We consider the additive category $\Pmorcx S$ of definition
 \ref{df:Pmorc} and the graph functor 
 $\gamma:\Pmorx S \rightarrow \Pmorcx S$ of \eqref{eq:graph_functor}.
\begin{df}
\label{df:ptrx}
A \emph{presheaf with transfers}
\index{word}{presheaf!with transfers}
\index{word}{transfer|see{presheaf \emph{or} sheaf}}
 $F$ over $S$ is an additive presheaf of 
 $\Rc$-modules over $\Pmorcx S$. 
We denote by $\pshg{\Pmorcx S}$
\index{notat}{pshPcor@$\pshg{\Pmorcx S}$}
 the corresponding category.

If $X$ is a $\Pmor$-scheme over $S$,
 we denote by $\lrep S X$
 \index{notat}{LambdatrSX@$\lrep S X$}
  the presheaf with transfers represented by $X$.

We denote by $\hat \gamma_*$ the functor
\begin{equation} \label{graph_psh_tr}
\pshg {\Pmorcx S} \rightarrow \psh {\Pmorx S}, F \mapsto F \circ \gamma.
\end{equation}
\end{df}

Note that $\pshg{\Pmorcx S}$ is obviously a Grothendieck abelian category
 generated by the objects $\lrep S X$ for a $\Pmor$-scheme $X/S$.
Moreover, the following proposition is straightforward:
\begin{prop} \label{prop:pshtr&Ppremotivic}
There is an essentially unique
 Grothendieck abelian $\Pmor$-premotivic category $\pshg{\Pmorc}$
 which is geometrically generated (\textit{cf.} \ref{df:generating_twists}),
 whose fiber over a scheme $S$ is $\pshg{\Pmorcx S}$ and such
 that the functor $\lrepNP_S$ induces a morphism
 of additive monoidal $\Pmor$-fibred categories.
\begin{equation} \label{eq:Pmorc&pshtr}
\Pmorc \rightarrow \pshg{\Pmorc}.
\end{equation}
Moreover, the functor \eqref{graph_psh_tr} induces a morphism
 of abelian $\Pmor$-premotivic categories
$$
\hat \gamma^*:\psh {\Pmor} \rightleftarrows \pshg {\Pmorc}:\hat \gamma_*.
$$
\end{prop}
\begin{proof}
To help the reader, we recall the following consequence of Yoneda's lemma:
\begin{lm}
Let $F:(\Pmorcx S)^{op} \rightarrow \Mod \Rc$ be a presheaf with transfers.
Let $\cI$ be the category of representable presheaves with transfers
 over $F$. Then the canonical map
$$
\ilim_{\lrep S X \rightarrow F} \lrep S X\To F
$$
is an isomorphism.
The limit is taken in $\pshg{\Pmorcx S}$ and runs over $\cI$.
\end{lm}
This lemma allows us to define the structural left adjoint of $\pshg{\Pmorc}$
 (recall $f^*$, $p_\sharp$ for $p$ a $\Pmor$-morphism and the tensor product)
 because they are indeed determined by \eqref{eq:Pmorc&pshtr}.
 The existence of the structural right adjoints is formal.

The same lemma allows to get the adjunction $(\hat \gamma^*,\hat \gamma_*)$.
\end{proof}

\begin{rem} \label{eq:computing_left_adj_pshtr}
Note that for any presheaf with transfers $F$ over $S$,
 and any morphism $f:T \rightarrow S$ (resp. $\Pmor$-morphism $p:S \rightarrow S'$),
 we get as usual $f_*F=F \circ f^*$ (resp. $p^*F=F \circ p_\sharp$)
 where the functor $f^*$ (resp. $p_\sharp$) on the right hand side is
 taken with respect to the $\Pmor$-fibred category $\Pmorc$.
\end{rem}

Let us state the following lemma for future use.
\begin{lm}\label{lm:ptr_property(C)}
Let Let $(S_\alpha)_{\alpha\in A}$ be a projective system of schemes 
\index{word}{projective system, of schemes}
 in $\sch$, with dominant affine transition maps, and
such that $S=\varprojlim_{\alpha\in A}S_\alpha$ is representable in $\sch$.

Consider an index $\alpha_0 \in A$ and a presheaf with transfers $F$
 over $S_{\alpha_0}$. For any index $\alpha/\alpha_0$, we denote by
 $F_\alpha$ (resp. $F$) the pullback of $F_{\alpha_0}$ over $S_\alpha$
 (resp. $S$) in the sense of the premotivic structure on $\pshg{\Pmorc}$.

Then the canonical map:
$$
\varinjlim_{\alpha \in A/\alpha_0} F_\alpha(S_\alpha) \longrightarrow F(S)
$$
is an isomorphism.
\end{lm}
\begin{proof}
The presheaf $F_{\alpha_0}$ can be written as an inductive limit
 of representable sheaves of the form $\lrep{S_{\alpha_0}}{X_{\alpha_0}}$
 of a $\Pmor$-scheme $X_{\alpha_0}/S_{\alpha_0}$.
 As the global section functor on presheaves with transfers commute
 with inductive limit, we are reduced
 to the case where $F=\lrep{S_{\alpha_0}}{X_{\alpha_0}}$.
 In this case, the lemma follows directly from
 Proposition \ref{prop:finite_cor&limit}.
\end{proof}

\subsection{Sheaves with transfers}

\begin{df}
\label{df:ftrx}
A $t$-sheaf with transfers
\index{word}{sheaf!tsheafwithtr@$t$-sheaf with transfers}
 over $S$ is a presheaf with transfers $F$ such
that the functor $F \circ \gamma_S$ is a $t$-sheaf. 
We denote by $\shg t {\Pmorcx S}$
\index{notat}{shtPcor@$\shg t {\Pmorcx S}$}
 the full subcategory of $\psh {\Pmorcx S}$ of sheaves with transfers.
\end{df}
According to this definition, we get a canonical faithful functor
$$
\gamma_*:\shg t {\Pmorcx S} \rightarrow \sh t {\Pmorx S}, F \mapsto F \circ \gamma.
$$

\begin{ex}
A particularly important case for us is the case when $t=\nis$
 is the Nisnevich topology.
 According to the original definition of Voevodsky,
 a Nisnevich sheaf with transfers will be called simply
  a \emph{sheaf with transfers}.
\index{word}{sheaf!sheaf with transfers}
\end{ex}
%The notation of this definition is convenient for this section.
%After this section\footnote{See definition \ref{df:premotivic_ftr} 
%and what follow for the definitive conventions.}, 
%we will use the following simplified notation:
%\begin{align*}
% \ptr S&:=\pshg {\smcx S}, \qquad  & \uptr S:=\pshg {\sftcx S} \\
% \ftrx t S&:=\shg t {\smcx S},
%  & \uftrx t S:=\shg t {\sftcx S}.
%\end{align*}

\begin{rem}
Later on, in the case $\Pmor=\sft$, 
we will use the notation $\ulrep S X$
\index{notat}{LambdatrSX@$\ulrep S X$}
 to denote the presheaf on
the big site $\sftcx S$ represented by a separated $S$-scheme of
finite type.
\end{rem}

\begin{prop}
\label{prop:lrep+et,qfh}
Let $X$ be an $\Pmor$-scheme over $S$.
\begin{enumerate}
\item The presheaf $\lrep S X$ is an \'etale sheaf with transfers.
\index{word}{sheaf!etalesheafwithtr@\'etale sheaf with transfers}
\item If $\car X \subset \Rc^\times$,
 $\lrep S X$ is a $\qfh$-sheaf with transfers.
\end{enumerate}
\end{prop}
\begin{proof}
For point (1), we follow the proof of \cite[4.2.4]{Deg7}: the computation
of the pullback by an \'etale map is given in our context by 
point (3) of Proposition \ref{basic_propert_corr}. Moreover,
the property for a cycle $\alpha/Y$ to be $\Rc$-universal is \'etale-local on $Y$
according to \ref{cor:special_etale_local}.
For point (2), we refer to \cite[4.2.7]{SV1}.
\end{proof}

%\begin{ex}
%In this example and the next lemma,
% $\Pmor$ is the class of separated morphisms
% and $f:X \rightarrow S$ an \'etale separated morphism.
%
%If $f$ is an \'etale cover,
% the presheaf with transfer $\lrep S X$ 
% is the constant \'etale $\Rc$-sheaf represented by $X$
% in the big site $\Pmorx S$.
% 
%On the contrary, if $f(X)$ does not contain any connected component
% of $S$, then $\lrep S X=0$.
% 
%This is shown by the next lemma.
%\end{ex}

We can actually describe explicitly representable presheaves with transfers
 in the following case:
\begin{prop} \label{prop:rep_et_ftr}
Let $S$ be a scheme and $X$ be a finite \'etale $S$-scheme.
Then for any $\Pmor$-scheme $Y$ over $S$,
$$
\Gamma(Y,\lrep S X)=\pi_0(Y \times_S X).\Rc.
$$
\end{prop}
\noindent This readily follows from the following lemma:
\begin{lm}\label{lm:rel_cycle_on_etale_ext_0}
Let $f:X \rightarrow S$ be an \'etale separated morphism of finite type.
Let $\pi^{finite}_0(X/S)$ be the set
 of connected components $V$ of $X$
 such that $f(V)$ is equal to a connected component of $S$
 (\emph{i.e.} $f$ is finite over $V$).

Then $c_0(X/S,\Rc)=\pi_0^{finite}(X/S).\Rc$.
\end{lm}
\begin{proof}
We can assume that $S$ is reduced and connected.

We first treat the case where $X=S$.
Consider a finite $\Rc$-universal $S$-cycle
 $\alpha$ with domain $S$.
Write $\alpha =\sum_{i \in I} n_i.\acycl{Z_i}_S$ in standard form.
By definition, $Z_i$ dominates an irreducible component of
$S$ thus $Z_i$ is equal to that irreducible component. \\
Consider $S_0$ an irreducible component of $S$ 
and an index $i \in I$ such that 
$S_0 \cap Z_i$ is not empty.
Consider a point $s \in S_0 \cap Z_i$.
We have obviously $\alpha_s=n_i.\acycl{\spec{\kappa(s)}} \neq 0$. 
Thus there exists a component of $\alpha$ which dominates $S_0$ 
\emph{i.e.} $\exists j \in I$ such that $Z_j=S_0$.
Moreover, computing $\alpha_s$ using alternatively $Z_i$ and $Z_j$
gives $n_i=n_j$. \\
As $S$ is noetherian, we see inductively $\{Z_i | i\in I\}$ is
the set of irreducible components of $S$ and for any $i,j \in I$,
$n_i=n_j$. Thus $c_0(S/S,\Rc)=\ZZ$.

Consider now the case of an \'etale $S$-scheme $X$.
By additivity of $c_0$, we can assume that $X$ is connected.
Consider the following canonical map:
$$
c_0(X/S,\Rc) \rightarrow c_0(X \times_S X/X,\Rc),
 \alpha \mapsto \alpha \otimesf_S X.
$$
Note that considering the projection
 $p:X \times_S X \rightarrow X$, by definition,
  $\alpha \otimesf_S X=p^*(\alpha)$. \\
Consider the diagonal $\delta:X \rightarrow X \times_S X$ of $X/S$.
Because $X/S$ is \'etale and separated,
 $\delta$ is a direct factor of $X \times_S X$ and we can write
 $X \times_S X=X \sqcup U$.
Because $c_0$ is additive,
$$
c_0(X \times_S X/X,\Rc)=c_0(X/X,\Rc) \oplus c_0(U/X,\Rc).
$$
Moreover, the projection on the first factor is induced
by the map $\delta^*$ on cycles.
Because $\delta^*p^*=1$, we deduce that a cycke in $c_0(X/S,\Rc)$
 corresponds uniquely to a cycle in $c_0(X/X,\Rc)$.
According to the preceding case,
 this latter group is the free group generated
 by the cycle $\acycl X$.
This latter cycle is $\Rc$-universal over $S$,
 because $X/S$ is flat.
Thus, if $X/S$ is finite, it is an element of $c_0(X/S,\Rc)$
 so that $c_0(X/S,\Rc)=\Rc$.
Otherwise, not any of the $\Rc$-linear combination of $\acycl X$
 belongs to $c_0(X/S,\Rc)$ so that $c_0(X/S,\Rc)=0$.
\end{proof}

\subsection{Associated sheaf with transfers}

\begin{num}
Recall from \ref{hypothesissite} that we denote by $(\Pmorx S)^\amalg$
 the category of $I$-diagrams of objects in $\Pmorx S$ indexed
 by a discrete category $I$.
Given any simplicial object $\cX$ of $(\Pmorx S)^\amalg$,
 we will consider the complex $\lrep S \cX$ of $\pshg{\Pmorcx S}$
 applying the definition of \ref{num:motives_simplicial_sm}
 to the Grothendieck $\Pmor$-fibred category $\pshg{\Pmor}$.

Consider a $t$-cover $p:W \rightarrow X$ in $\Pmorx X$.
We denote by $W_X^n$ the $n$-fold product of $W$ over $X$
 (in the category $\Pmorx X$).
We denote by $\check S(W/X)$ the \v Cech simplicial object 
 of $\Pmorcx S$ such that $\check S_n(W/X)=W^{n+1}_X$.
 The canonical morphism
  $\check S(W/X) \rightarrow X$ is a $t$-hypercover according
   to definition \ref{hypothesissite}.
 We will call these particular type of $t$-hypercovers
  the \emph{\v Cech $t$-hypercovers}
\index{word}{hypercover!Cechthypercover@\v Cech $t$-hypercovers}
   of $X$.
\end{num}
%
%%We consider the functor
%%\begin{align*}
%%\lrepNP_S:(\Pmorcx X)^\sqcup & \rightarrow \psh{\Pmorcx S} \\
%% (U'_i)_{i \in I} &\mapsto \oplus_{i \in I} \lrep S {U'_i}.
%%\end{align*}
%Applying the functor $\lrepNP_S$ to the simplicial object $\check S_*(W/X)$
% we get a simplicial presheaf with transfers ;
% according to convention \ref{num:motives_simplicial_sm},
% we simply denote by $\lrep S {\check S_*(W/X)}$ the associated complex
% in the category $\pshg {\Pmorcx S}$.
%We thus get a morphism of complexes
%\begin{equation} \label{eq:Cech&cor}
%\lrep S {\check S_*(W/X)} \rightarrow \lrep S X
%\end{equation}
%where we consider $\lrep S X$ as a complex concentrated in degree $0$.
\begin{df} \label{df:top_comp_tr}
We will say that the admissible topology $t$ on $\Pmor$ is
 \emph{compatible with transfers} 
\index{word}{topology!compatible with transfers}
\index{word}{compatible with transfers|see{topology}}
 (resp. \emph{weakly compatible with transfers})
\index{word}{topology!weakly compatible with transfers}
 if for any scheme $S$ 
 and any $t$-hypercover (resp. any \v Cech $t$-hypercover)
 $\cX \rightarrow X$ in the site $\Pmorx S$,
 the canonical morphism of complexes
\begin{equation} \label{eq:simplicial&cor}
\lrep S {\cX} \rightarrow \lrep S X
\end{equation}
 induces a quasi-isomorphism 
 of the associated $t$-sheaves on $\Pmorx S$.
\end{df}
Obviously, if $t$ is compatible with transfers
 then it is weakly compatible with transfers.

Recall from \ref{prop:lrep+et,qfh} that,
 in the cases $t=\nis, \et$,
 \eqref{eq:simplicial&cor} is actually a morphism of complexes
 of $t$-sheaves with transfers.
The following proposition is a generalization of \cite[3.1.3]{V1}
 but its proof is in fact the same.
\begin{prop}
\label{lm:acyc}
The Nisnevich (resp \'etale) topology $t$ on $\Pmor$ is
 weakly compatible with transfers.
\end{prop}
\begin{proof}
We consider a $t$-cover $p:W \rightarrow X$,
 the associated \v Cech hypercover $\cX=\check S(W/X)$ of $X$
 and we prove that the map \eqref{eq:simplicial&cor} is a quasi-isomorphism
 of $t$-sheaves.
Recall that a point of $\Pmorx S$ for the topology $t$ is given by
 an essentially affine pro-object $(V_i)_{i \in I}$ of $\Pmorx S$.
 Moreover, its projective limit $\cV$ in the category of schemes
 is in particular a local henselian noetherian scheme. \\
It will be sufficient to check that the fiber of \eqref{eq:simplicial&cor}
 at the point $(V_i)_{i \in I}$ is a quasi-isomorphism.
 Thus, according to Proposition \ref{prop:finite_cor&limit},
 we can assume that $S=\cV$ is a local henselian scheme
 and we are to reduce to prove that the complex of abelian
 groups
$$
\hdots \rightarrow c_0(W \times_X W/S,\Rc)
 \rightarrow c_0(W/S,\Rc)
 \xrightarrow{p_*} c_0(X/S,\Rc)
 \rightarrow 0
$$
is acyclic. We denote this complex by $C$.

Recall that the abelian group $c_0(X/S)$
 is covariantly functorial in $X$ with respect to
 separated morphisms of finite type $f:X' \rightarrow X$
  (\textit{cf.} paragraph \ref{num:c_0&pushout}).
Moreover, if $f$ is an immersion, $f_*$ is obviously injective.

Let $\cF_0$ be the set
 of closed subschemes $Z$ of $X$ such that $Z/S$ is finite.
Given a closed subscheme $Z$ in $\cF_0$,
 we let $C_Z$ be the complex of abelian groups
\begin{equation} \label{eq:Cech&cor_bis}
\hdots \rightarrow c_0(W_Z \times_Z W_Z/S,\Rc)
 \rightarrow c_0(W_Z/S,\Rc)
 \xrightarrow{p_{Z*}} c_0(Z/S,\Rc)
 \rightarrow 0
\end{equation}
where $p_Z$ is the $t$-cover obtained by pullback along $Z \rightarrow X$.
From what we have just recalled,
 we can identify $C_Z$ with a subcomplex of $C$.
The set $\cF_0$ can be ordered by inclusion,
 and $C$ is the union of its subcomplexes $C_Z$.
If $\cF_0$ is empty, then $C=0$ and the proposition is clear.
Otherwise, $\cF_0$ is filtered and we can write:
$$
C=\ilim_{Z \in \cF_0} C_Z.
$$
Thus, it will be sufficient to prove that $C_Z$ is acyclic
 for any $Z \in \cF_0$.
Because $S$ is henselian and $Z$ is finite over $S$,
 $Z$ is indeed a finite sum of local henselian schemes.
 This implies that the $t$-cover $p_Z$,
 which is in particular \'etale surjective,
 admits a splitting $s:Z \rightarrow W_Z$.
 Then the complex \eqref{eq:Cech&cor_bis} is contractible
 with contracting homotopy defined by the family of morphisms
$$
(s \times_Z 1_{W_Z^n})_*:c_0(W_Z^{n}/S,\Rc)
 \rightarrow c_0(W_Z^{n+1}/S,\Rc).
$$
\end{proof}

\begin{num} Considering an additive abelian presheaf $G$ on $\Pmorx S$,
 the natural transformation
$$
X \mapsto \Hom_{\pshg{\Pmorx S}}(\hat \gamma_*\lrep S X,G)
$$
defines a presheaf with transfers over $S$.\footnote{Actually,
 this defines a right adjoint to the functor $\hat \gamma_*$.}
We will denote by $G_\tau$ its restriction to the site $\Pmorx S$.
Note that this definition can be applied in the case where $G$
 is a $t$-sheaf on $\Pmorx S$,
  because under the assumption \ref{num:fxtr_convention_class_P} on $t$, 
  it is in particular an additive presheaf.
\end{num}
\begin{df}
We will say that $t$ is mildly compatible with transfers
\index{word}{topology!mildly compatible with transfers}
 if for any scheme $S$ and any $t$-sheaf $F$ on $\Pmorx S$,
 $F_\tau$ is a $t$-sheaf on $\Pmorx S$.
\end{df}
If $t$ is weakly compatible with transfers
 then is it mildly compatible with transfers.

\begin{rem} \label{rem:mildly&epi_lrep}
Assume $t$ is mildly compatible with transfers.
Then for any scheme $S$,
 any $t$-cover $p:W \rightarrow X$ in $\Pmorx S$ induces a morphism
$$
p_*:\lrep S W \rightarrow \lrep S X
$$
which is an epimorphism of the associated $t$-sheaves on $\Pmorx S$.
This means that for any correspondence $\alpha \in \corr S Y X$,
there exists a $t$-cover $q:W' \rightarrow Y$ and a correspondence
$\alpha' \in \corr S {W'} W$ making the following diagram commutative:
\begin{equation} \label{eq:lift_corr}
\begin{split}
\xymatrix@=18pt{
W'\ar^-{\hat \alpha}[r]\ar|/-8pt/\bullet[r]\ar_q[d] & W\ar^p[d] \\
Y\ar^-\alpha[r]\ar|/-10pt/\bullet[r] & X
}
\end{split}
\end{equation}
\end{rem}

\begin{lm}
Assume $t$ is mildly compatible with transfers.

Let $S$ be a scheme and $P^{tr}$ be a presheaf with transfers over $S$.
We put $P=P^{tr} \circ \gamma$ as a presheaf on $\Pmorx S$.
We denote by $F$ the $t$-sheaf associated with $P$
 and by $\eta:P \rightarrow F$ the canonical natural transformation.

Then there exists a unique pair $(F^{tr},\eta^{tr})$
 such that:
\begin{enumerate}
\item $F^{tr}$ is a sheaf with transfers over $S$
 such that $F^{tr} \circ \gamma=F$.
\item $\eta^{tr}:P^{tr} \rightarrow F^{tr}$ is a natural transformation
 of presheaves with transfers such that the induced transformation
$$
P=(P^{tr} \circ \gamma) \rightarrow (F^{tr} \circ \gamma)=F
$$
coincides with $\eta$.
\end{enumerate}
\end{lm}
\begin{proof}
As a preliminary observation,
 we note that given a presheaf $G$ on $\Pmorx S$,
 the data of a presheaf with transfers $G^{tr}$ such that $G^{tr} \circ \gamma=G$
 is equivalent to the data for any $\Pmor$-schemes $X$ and $Y$ of a bilinear product
\begin{equation} \label{eq:product_tr}
G(X) \otimes_\ZZ \corr S Y X \rightarrow G(Y),
 \rho \otimes \alpha \mapsto \langle \rho,\alpha \rangle
\end{equation}
such that:
\begin{enumerate}
\item[(a)] For any morphism $f:Y' \rightarrow Y$ in $\Pmorx S$,
 $f^*\langle \rho,\alpha \rangle=\langle \rho,\alpha \circ f \rangle$.
\item[(b)]  For any morphism $f:X \rightarrow X'$ in $\Pmorx S$,
 if $\rho=f^*(\rho')$,
 $\langle \rho,\alpha \rangle=\langle \rho',f \circ \alpha \rangle$.
\item[(c)]  When $X=Y$, for any $\rho \in F(X)$, $\langle \rho,1_X \rangle=\rho$.
\item[(d)] For any finite $S$-correspondence $\beta \in \corr S Z Y$,
 $\langle \langle \rho,\alpha \rangle,\beta \rangle=
 \langle \rho,\alpha \circ \beta \rangle$.
\end{enumerate}
On the other hand, the data of products of the form \eqref{eq:product_tr}
 for any $\Pmor$-schemes $X$ and $Y$ over $S$ which satisfy 
 the conditions (a) and (b) above is equivalent to the data of a natural transformation
$$
\phi:G \rightarrow G_\tau
$$
by putting $\langle \rho,\alpha \rangle_\phi
=\lbrack\phi_X(\rho)\rbrack_Y.\alpha$.

Applying this to the presheaf with transfers $P^{tr}$,
 we obtain a canonical natural transformation
$$
\psi:P \rightarrow P_\tau.
$$
By assumption on $t$, $F_\tau$ is a $t$-sheaf.
Thus, there exists a unique natural transformation $\psi$
 such that the following diagram commutes:
$$
\xymatrix@=20pt{
P\ar^-{\psi}[r]\ar_\eta[d] & P_\tau\ar^{\eta_\tau}[d] \\
F\ar^-\phi[r] & F_\tau
}
$$
Thus we get products of the form \ref{eq:product_tr} associated with $\phi$
 which satisfies (a) and (b).
 The commutativity of the above diagram asserts
  they are compatible with the ones corresponding to $P^{tr}$
 and the unicity of the natural transformation $\phi$ implies
  the uniqueness statement of the lemma.

To finish the proof of the existence, we must show (c) and (d) for the product
 $\langle.,\rangle_\phi$.
Consider a couple $(\rho,\alpha) \in F(X) \times \corr S Y X$.
Because $F$ is the $t$-sheaf associated with $P$,
 there exists a $t$-cover $p:W \rightarrow X$ 
 and a section $\hat \rho \in P(W)$ such that $p^*(\rho)=\eta_W(\hat \rho)$.
Moreover, according to Remark \ref{rem:mildly&epi_lrep},
 we get a $t$-cover $q:W' \rightarrow Y$
 and a correspondence $\hat \alpha \in \corr S {W'} W$
  making the diagram \eqref{eq:lift_corr}
 commutative. Then we get using (a) and (b):
$$
q^*\langle \rho,\alpha \rangle_\phi
=\langle \rho,\alpha \circ q\rangle_\phi
=\langle \rho,p \circ \hat \alpha \rangle_\phi
=\langle p^*\rho,\hat \alpha \rangle_\phi
=\langle \eta_W(\hat \rho),\hat \alpha \rangle_\phi
=\langle \hat \rho,\hat \alpha \rangle_\psi.
$$
Because $q^*:F(X) \rightarrow F(W)$ is injective,
 we deduce easily from this principle the properties (c) and (d)
  and this concludes.
%\noindent \textit{Unicity}: Using the notation of \eqref{eq:product_tr},
% we show that for any $\alpha \in \corr S X Y$ and $\rho \in F_\sigma(Y)$,
% the element
%$$
%\langle \alpha,\rho \rangle_\phi \in F_\sigma(X)
%$$
%is uniquely determined by $\hat \phi$. 
%By definition of $F$, there exists a $t$-cover $q:V \rightarrow Y$ 
% and a section $\hat \rho \in \hat F(V)$ such that $q^*(\rho)=a_V(\hat \rho)$
% -- note that we have also used our assumption on $t$ and
%  the fact that $\hat F$ is additive.
%Moreover, because $t$ is weakly compatible with transfers, there exists
%a $t$-cover $U \rightarrow X$ and a correspondence $\beta \in \corr S U V$
%making the following diagram commutative:
%$$
%\xymatrix@=20pt{
%U\ar^\beta[r]\ar_p[d] & V\ar^q[d] \\
%X\ar^\alpha[r] & Y.
%}
%$$
%Thus,
%$$
%p^*\langle \alpha,\rho \rangle_\phi=\langle p^*(\alpha),\rho \rangle_\phi
%=\langle q_*(\beta),\rho \rangle_\phi=\langle \beta,q^*(\rho) \rangle_\phi
%=\langle \beta,\hat \rho \rangle_{\hat \phi}.
%$$
%This characterize uniquely $\phi$ as by assumption,
% $F_\sigma$ is a $t$-sheaf.
%
%\noindent \textit{Existence}: 
\end{proof}

\begin{num} Let us consider the canonical adjunction
$$
a^*_t:\psh{\Pmorx S} \rightleftarrows \sh t {\Pmorx S}:\mathcal O_t
$$
where $\mathcal O_t$ is the canonical forgetful functor.

We also denote by 
 $\mathcal O_t^{tr}:\shg t {\Pmorcx S} \rightarrow \pshg {\Pmorcx S}$
the obvious forgetful functor.
Trivially, the following relation holds:
\begin{equation} \label{eq:O_tr&gamma}
\hat \gamma_* \, a_{t,*}=a_{t,*}\,  \gamma_*.
\end{equation}
\end{num}
\begin{prop} \label{prop:exists_a^tr}
Using the notations above, the following condition on the admissible topology $t$
 are equivalent:
\begin{enumerate}
\item[(i)] $t$ is mildly compatible with transfers.
\item[(ii)] For any scheme $S$,
the functor $\cO_t^{tr}$ admits a left adjoint
$$a_t^{*}:\pshg {\Pmorcx S} \rightarrow \shg t {\Pmorcx S}$$
which is exact and
such that the exchange transformation 
\begin{equation} \label{eq:a_tr&gamma}
a^*_t \, \hat \gamma_* \rightarrow \, \gamma_* \, a^*_{t}
\end{equation}
induced by the identification \eqref{eq:O_tr&gamma}
is an isomorphism.
\end{enumerate}
Under these conditions, the following properties hold for any scheme $S$:
\begin{enumerate}
\item[(iii)] The category $\shg t {\Pmorcx S}$ is a Grothendieck abelian 
 category.
\item[(iv)] The functor $\gamma_*$ commutes with all limits and colimits.
\end{enumerate}
\end{prop}
\begin{proof}
The fact (i) implies (ii) follows from the preceding lemma
 as we can put $a_t^{tr}(F)=F^{tr}$ with the notation of the lemma.
 The fact this defines a functor, as well as the properties stated in (ii),
  follows from the uniqueness statement of \emph{loc. cit.}

Let us assume (ii). Then (iii) follows formally from (ii),
from the existence, adjunction property and exactness of $a_t^{*}$,
as $\pshg{\Pmorcx S}$ is a Grothendieck abelian category.
Moreover, we deduce from the isomorphism \eqref{eq:a_tr&gamma}
that $\gamma_*$ is exact:
indeed, $a^*_t$ and $\hat \gamma_*$ are exact. As $\gamma_*$
 commutes with arbitrary direct sums, we get (iv). \\
From this point, we deduce the existence of a right adjoint
 $\gamma^!$ to the functor $\gamma_*$.
 Using again the isomorphism \eqref{eq:a_tr&gamma},
 we obtain for any $t$-sheaves $F$ on $\Pmorx S$
 and any $\Pmor$-scheme $X/S$
 a canonical isomorphism $F_\tau(X)=\gamma^! F(X)$.
 This proves (i).
\end{proof}

\begin{num} Under the assumption of the previous proposition,
 given any $\Pmor$-scheme $X/S$, we will put
 $\lrep S X_t=a_t^{*} \lrep S X$. 
 The above proposition shows that the family $\lrep S X_t$
 for $\Pmor$-schemes $X/S$  is a generating family
 in $\shg t {\Pmorcx S}$. Moreover,
 we get easily the following corollary of the preceding proposition
 and Proposition \ref{prop:pshtr&Ppremotivic}:
\end{num}
\begin{cor} \label{cor:Ppremotivic_ftr}
Assume that $t$ is mildly compatible with transfers.
\index{word}{topology!mildly compatible with transfers}

Then there exists an essentially unique Grothendieck abelian
 $\Pmor$-premotivic category $\shg t {\Pmorc}$ which
 is geometrically generated, whose fiber over a scheme $S$
 is $\shg t {\Pmorcx S}$ and such that the $t$-sheafification functor
 induces an adjunction of abelian $\Pmor$-premotivic categories:
$$
a_t^{*}:\pshg{\Pmorc} \rightleftarrows \shg t {\Pmorc}:\cO_t^{tr}.
$$
Moreover, the functor $\gamma_*$ induces an adjunction
 of abelian $\Pmor$-premotivic categories:
\begin{equation} \label{eq:graph_premotivic_adj_ftr}
\gamma^*:\sh t {\Pmor} \rightleftarrows \shg t {\Pmorc}:\gamma_*.
\end{equation}
\end{cor}
\begin{rem} \label{rem:commutation_atr&gamma}
Notice moreover that $\gamma^* \, a^*_t=a_t^{*} \, \hat \gamma^*$.
\end{rem}
\begin{proof}
In fact, using the exactness of $a_t^{*}$,
given any sheaf $F$ with transfers $F$ over $S$,
 we get a canonical isomorphism
$$
F=\ilim_{\lrep S X_t \rightarrow F} \lrep S X_t
$$
where the limit is taken in $\shg t {\Pmorcx S}$
 and runs over the representable $t$-sheaves with transfers
 over $F$.
As in the proof of \ref{prop:pshtr&Ppremotivic},
 this allows defining uniquely the structural left adjoints
 of $\shg t {\Pmorc}$. The existence (and uniqueness)
 of the structural right adjoints then follows formally.
The same remark allows to construct the functor $\gamma^*$.
\end{proof}

\begin{rem} \label{rem:local&flasque_shtr}
Let us explicit the meaning of the preceding Corollary
 for a topology $t$ which is compatible with transfers.
 Given a complex $C$ with coefficients in the category $\shg t {\Pmorcx S}$,
 the following conditions are equivalent:
\begin{enumerate}
\item[(i)] $C$ is local (Definition \ref{df:basic_complexes&topology}),
\item[(i')] $\gamma_*(C)$ is local,
\item[(i'')] given any $\Pmor$-scheme $X/S$
 and any integer $n \in \ZZ$, the canonical map
$$
H^n(C(X)) \rightarrow H^n_t(X,\gamma_*(C))
$$
is an isomorphism,
\item[(ii)] $C$ is $t$-flasque (Definition \ref{df:basic_complexes&topology}),
\item[(ii')] $\gamma_*(C)$ is $t$-flasque,
\item[(ii'')] given any $t$-hypercover $p:\mathcal X \rightarrow X$
 in $\Pmor/S$ and any integer $n \in \ZZ$, the canonical map
$$
p^*:H^n(C(X)) \rightarrow H^n(C(\mathcal X))
$$
is an isomorphism.
\end{enumerate}
More precisely, the equivalence of (i) and (ii) is the preceding corollary,
 while the equivalence of (i) and (i') (resp. (ii) and (ii'))
 follows from the existence of the adjunction \eqref{eq:graph_premotivic_adj_ftr}
 and the fact $\gamma_*$ is exact.
 The equivalence between (i') and (i'') (resp. (ii') and (ii''))
 is a simple translation of Definition \ref{df:basic_complexes&topology}.
\end{rem}

\begin{num} Recall from Definition \ref{df:basic_complexes&topology}
 we say that the abelian $\Pmor$-premotivic category $\shg t {\Pmorc}$
 satisfies cohomological $t$-descent if for any scheme $S$,
 and any $t$-hypercover $\cX \rightarrow X$ in $\Pmorx S$,
 the induced morphism of complexes in $\shg t {\Pmorcx S}$
$$
\lrep S {\cX}_t \rightarrow \lrep S X_t
$$
is a quasi-isomorphism.
The preceding corollary thus gives the following one:
\end{num}
\begin{cor} \label{cor:compatibility_tr&top}
Assume $t$ is mildly compatible with transfers. \\
Then the following conditions are equivalent:
\begin{enumerate}
\item[(i)] The topology $t$ is compatible with transfers.
\item[(ii)] The abelian $\Pmor$-premotivic category
 $\shg t {\Pmorc}$ satisfies cohomological $t$-descent.
\item[(iii)] The abelian $\Pmor$-premotivic category
 $\shg t {\Pmorc}$ is compatible with $t$
 (see \ref{df:basic_complexes&topology}).
\end{enumerate}
%If we assume that any scheme $S$ in $\sch$
% is $t$-cohomologically bounded with respect to $t$-sheaves
% of $\Rc$-modules, then these conditions are equivalent
% to the following one:
%\begin{enumerate}
%\item[(iv)] The topology $t$ is weakly compatible with transfers.
%\end{enumerate}
\end{cor}
\begin{proof}
The equivalence of (i) and (ii) follows easily
 from the isomorphism \eqref{eq:a_tr&gamma}.
 The equivalence of (ii) and (iii) is Proposition \ref{adjunction&topology}
 applied to the adjunction \eqref{eq:graph_premotivic_adj_ftr},
 in view of \ref{prop:exists_a^tr}(iv).
%It remains only to prove that (iv) implies (ii). 
%Recall that the category $\Der(\sh t {\Pmorx S}$ has a canonical DG-structure
% (see for example \ref{num:canonical_DG-structure}).
%Point (ii) can be reformulated  by saying that for any complex $K$
% of $t$-sheaves on $\Pmorx S$, and any $t$-hypercover $\cX \rightarrow X$,
% the canonical map of $\Der(\Mod \Rc)$
%$$
%\derR \sHom_{\Der(\sh t {\Pmorx S})}(\gamma_*\lrep S X_t,K)
% \rightarrow \derR \sHom_{\Der(\sh t {\Pmorx S})}(\gamma_*\lrep S X_t,K)
%$$
%is an isomorphism.
%Recall also there is the injective model structure
% on $\Comp(\sh t {\Pmorx S})$ for which every object is cofibrant
% and with quasi-isomorphisms as weak equivalences
%  (recall are given in \cite[2.1]{CD1}).
%Replacing $K$ by a fibrant resolution for the injective model structure,
% we get for any simplicial objects $\cX$ of $\Pmorx S^\amalg$ that:
%$$
%\derR \sHom_{\Der(\sh t {\Pmorx S})}(\gamma_*\lrep S X_t,K)
%=\sHom_{\Der(\sh t {\Pmorx S})}(\gamma_*\lrep S X_t,K)
%$$
%and we have to prove that the presheaf
%$$
%E:\Pmorx S^{op} \rightarrow \Comp(\Rc-mod)
%X \mapsto \sHom_{\Der(\sh t {\Pmorx S})}(\gamma_*\lrep S X_t,K)
%$$
%is $t$-local in the category $\Comp(\sh t {\Pmorx S})$.
\end{proof}

\begin {num} Recall from Paragraph \ref{df:cd_structures}
 that a cd-structure $P$ on $\sch$ is the data of a family of commutative squares, called $P$-distinguished,
 of the form
\begin{equation} \label{eq:cd_dist_qr_tr}
\begin{split}
\xymatrix@=14pt{
B\ar^k[r]\ar_g[d]\ar@{}|Q[rd] & Y\ar^f[d] \\
A\ar_i[r] & X
}
\end{split}
\end{equation}
which is stable by isomorphisms. 
Further, we will consider the following assumptions on $P$:
\begin{enumerate}
\item[(a)] $P$ is complete, regular and bounded in the sense of \cite{voecd2}.
\item[(b)] Any $P$-distinguished square as above
 is made of $\Pmor$-morphisms and $k$ is an immersion.
\item[(c)] Any square as above which is cartesian and such that
 $X=A \sqcup Y$, $i$ and $f$ being the obvious immersions,
 is $P$-distinguished.
\end{enumerate}
Then the topology $t_P$ associated with $P$ (see \ref{df:cd_structures})
 is $\Pmor$-admissible and satisfy assumption \ref{num:fxtr_convention_class_P}(c).
Moreover, according to \cite[2.9]{voecd2}, we obtain the following properties:
\begin{enumerate}
\item[(d)] A  presheaf $F$ on $\Pmorx S$
 is a $t_P$-sheaf if and only if $F(\varnothing)=0$
 and for any $P$-distinguished square \eqref{eq:cd_dist_qr_tr}
 in $\Pmorx S$, the sequence
$$
0 \rightarrow F(X) \xrightarrow{f^*+e^*} F(Y) \oplus F(A)
 \xrightarrow{k^*-g^*} F(B)
$$
is exact.
\item[(e)] For any $P$-distinguished square \eqref{eq:cd_dist_qr_tr}
 the sequence of representable pre-sheaves on $\Pmorx S$
$$
0 \rightarrow \rep S B \xrightarrow{k_*-g_*}
 \rep S Y \oplus \rep S A \xrightarrow{f_*+e_*} \rep S X \rightarrow 0
$$
becomes exact on the associated $t_P$-sheaves.
\end{enumerate}
\end{num}
\begin{prop} \label{prop:compatible&cd_tr}
Consider a cd-structure
\index{word}{cdstructure@cd-structure}
 $P$ satisfying properties (a) and (b) above
 and assume that $t=t_P$ is the topology associated with $P$.
Then the following conditions are equivalent:
\begin{enumerate}
\item[(i)] The topology $t$ is compatible with transfers.
\index{word}{topology!compatible with transfers}
\item[(ii)] The topology $t$ is mildly compatible with transfers.
\index{word}{topology!mildly compatible with transfers}
\item[(iii)] For any scheme $S$
 and any $P$-distinguished square \eqref{eq:cd_dist_qr_tr} in $\Pmorx S$,
 the short sequence of presheaves with transfers over $S$
$$
0 \rightarrow \lrep S B \xrightarrow{k_*-g_*}
 \lrep S Y \oplus \lrep S A \xrightarrow{f_*+e_*} \lrep S X \rightarrow 0
$$
becomes exact on the associated $t$-sheaves on $\Pmorx S$.
\end{enumerate}
\end{prop}
\begin{proof}
The implication (i) $\Rightarrow$ (ii) is obvious.

The implication (ii) $\Rightarrow$ (iii) follows from
 point (e) above and the following facts: $\gamma^*$ is right exact
  (Corollary \ref{cor:Ppremotivic_ftr}), $\gamma^* a_t=a_t^{tr} \hat \gamma^*$
  (remark \ref{rem:commutation_atr&gamma}),
 $k_*:\lrep S B \rightarrow \lrep S Y$ is a monomorphism of
   presheaves with transfers (use \ref{basic_propert_corr}(2) and the
   fact $k$ is an immersion from assumption (b)).

Assume (iii).
Then we obtain (ii) as a direct consequence of the point (d) above.
Thus, to prove (i), we have only to prove that
 the abelian $\Pmor$-premotivic category $\shg t {\Pmorc}$ satisfies
 cohomological $t$-descent according to \ref{cor:compatibility_tr&top}.

Let $S$ be a scheme.
Recall that the category $\Der\big(\sh t {\Pmorx S}\big)$ has a canonical DG-structure
 (see for example \ref{num:canonical_DG-structure}).
The cohomological $t$-descent for $\shg t {\Pmorcx S}$
 can be reformulated  by saying that for any complex $K$
 of $t$-sheaves on $\Pmorx S$, and any $t$-hypercover $\cX \rightarrow X$,
 the canonical map of $\Der(\Mod \Rc)$
$$
\derR \Hom^\bullet_{\Der(\sh t {\Pmorx S})}(\gamma_*\lrep S X_t,K)
 \rightarrow \derR \Hom^\bullet_{\Der(\sh t {\Pmorx S})}(\gamma_*\lrep S \cX_t,K)
$$
is an isomorphism.
Recall also that there is the injective model structure
 on the category $\Comp(\sh t {\Pmorx S})$ for which every object is cofibrant,
 with quasi-isomorphisms as weak equivalences
  (see \cite[2.1]{CD1} for more details).
Replacing $K$ by a fibrant resolution for the injective model structure,
 we get for any simplicial objects $\cX$ of $\Pmorx S^\amalg$ that:
$$
\derR \Hom^\bullet_{\Der(\sh t {\Pmorx S})}(\gamma_*\lrep S \cX_t,K)
=\Hom^\bullet_{\Der(\sh t {\Pmorx S})}(\gamma_*\lrep S \cX_t,K).
$$
Thus it is sufficient to prove that the presheaf
$$
E:\Pmorx S^{op} \rightarrow \Comp(\Mod \Rc),
X \mapsto \Hom^\bullet_{\Der(\sh t {\Pmorx S})}(\gamma_*\lrep S X_t,K)
$$
satisfies $t$-descent in the sense of \ref{defdescente}.

We derive from (iii) that $E$ sends a $P$-distinguished square
 to a homotopy cartesian square in $\Der(\Mod \Rc)$.
 Thus the assertion follows from the arguments on $t$-descent 
 from \cite{voecd1,voecd2}.
\end{proof}

\begin{rem} \label{rem:cd_bounded_ftr}
It follows from Remark \ref{rem:local&flasque_shtr} that
 under the equivalent conditions (i), (ii), (iii) of the above
 corollary, the admissible topology $t=t_P$ is bounded in
 $\shg t {\Pmorc}$ in the sense of Definition \ref{df:hyper-bounded}.
 Over a scheme $S$,
 a bounded generating family is given by the following complexes
 of $\shg t {\Pmorcx S}$:
$$
\cdots \rightarrow 0 \rightarrow \lrep S B \xrightarrow{k_*-g_*}
 \lrep S Y \oplus \lrep S A \xrightarrow{f_*+e_*} \lrep S X \rightarrow 0
 \rightarrow \cdots
$$
 induced by a $P$-distinguished square of the form \eqref{eq:cd_dist_qr_tr}
  -- see also Example \ref{ex:Zar&Nis=bounded}.
\end{rem}

We end-up this section with a compatibility of certain sheaves
 with transfers with projective limits of schemes.
 This will be the key point to establish continuity
 for motivic complexes.
\begin{prop} \label{prop:ftr_property(C)}
Let $t$ be one of the topologies $\nis, \et, \cdh$.

Let $(S_\alpha)_{\alpha\in A}$ be a projective system of schemes 
\index{word}{projective system, of schemes}
 in $\sch$, with dominant affine transition maps, and
such that $S=\varprojlim_{\alpha\in A}S_\alpha$ is representable in $\sch$.

Consider an index $\alpha_0 \in A$ and a $t$-sheaf with transfers $F$
 over $\sftcx{S_0}$.
 For any index $\alpha/\alpha_0$, we denote by
 $F_\alpha$ (resp. $F$) the pullback of $F_{\alpha_0}$ over $S_\alpha$
 (resp. $S$) in the sense of the premotivic structure on
 $\shg t {\Pmorc}$ (obtained in Corollary \ref{cor:Ppremotivic_ftr}).

Then the canonical map:
$$
\varinjlim_{\alpha \in A/\alpha_0} F_\alpha(S_\alpha) \longrightarrow F(S)
$$
is an isomorphism.
\end{prop}
\begin{proof}
We consider the forgetful functor:
$\mathcal O_t^{tr}:\shg t {\sftc} \rightarrow \pshg {\sftc}$.
It is fully faithful and it commutes with the global section functor.
We want to prove the proposition by using Lemma \ref{lm:ptr_property(C)}.
Thus it is sufficient to prove that, for any morphism $f:X \rightarrow S$
 in $\sch$,
 the functor $\mathcal O_t^{tr}$ commutes with $f^*$.
 In other words, the pullback functor $\hat f^*$ for presheaves with
 transfers on $\sftc$ preserves $t$-sheaves with transfers:
 for any $t$-sheaf with transfers $F$ over $S$, $\hat f^*(F)$
 is a $t$-sheaf with transfers.

Let us first treat the case where $f$ is separated of finite type.
 Then $\hat f^*$ admits a left adjoint $\hat f_\sharp$ which preserves
 $t$-covers. Thus the property is clear.

In the general case, we write $f$ as a projective limit of morphisms
 of schemes $(f_\alpha:X_\alpha \rightarrow S)_{\alpha \in A}$
 such that the transition morphisms of the projective scheme
 $(X_\alpha)_{\alpha \in A}$ are affine and dominant
 and each $f_\alpha$ is separated of finite type.\footnote{Write
 the $\mathcal O_S$-algebra $f_*(\mathcal O_X)$ as the filtered
 union of its finite type sub-$\mathcal O_S$-algebras, ordered by
 inclusion.}
 To check that $\hat f^*(F)$ is a $t$-sheaf,
 we consider a $t$-cover $p:W \rightarrow X$
 of an $S$-scheme separated of finite type. 
 Because of our choice of topology $t$,
 there exists an index $\alpha_1/\alpha_0$
 such that $p:W \rightarrow X$ can be lifted as a $t$-cover
 $p_1:W_{\alpha_1} \rightarrow X_{\alpha_1}$ over $S_{\alpha_1}$.
 Using Lemma \ref{lm:ptr_property(C)} again,
  we now are reduced to prove that for any $\alpha/\alpha_1$,
  $\hat f_{\alpha_1}^*(F)$ satisfies the $t$-sheaf property
  with respect to the pullback of $p_1$ over $S_{\alpha}/S_{\alpha_1}$.
  This follows from the first case treated.
\end{proof}

\begin{rem}
The previous proposition generalizes \cite[Prop. 2.19]{Deg7}.
\end{rem}

\subsection{Examples} \label{sec:examples_ftr}

\begin{num} Assume that $t$ is the Nisnevich topology.
According to Lemma \ref{lm:acyc} and 
 Proposition \ref{prop:compatible&cd_tr},
 $t$ is then compatible with transfers.
 With the notation of Corollary \ref{cor:Ppremotivic_ftr},
 we get the following definition:
\begin{df} \label{df:premotivic_ftr}
We denote by \index{notat}{shtrLambda@$\ftr -$}
\index{notat}{shtrLambda@$\uftr -$}
$$
\ftr -, \qquad \uftr -
$$
the respective abelian premotivic 
 and generalized abelian premotivic categories defined 
 in Corollary \ref{cor:Ppremotivic_ftr}
 in the respective cases $\Pmor=\ssm$, $\Pmor=\ssft$.
\end{df}
From now on, the objects of $\ftr S$ (resp. $\uftr S$)
are called \emph{sheaves with transfers} over $S$
\index{word}{sheaf!sheaf with transfers}
(resp. \emph{generalized sheaves with transfers} over $S$).
\index{word}{sheaf!generalized sheaf with transfers}

Let $X$ be a separated $S$-scheme of finite type.
We let $\ulrep S X$ be the generalized sheaf with transfers
 represented by $X$ (\textit{cf.} \ref{prop:lrep+et,qfh}).
If $X$ is $S$-smooth,
 we denote by $\lrep S X$ its restriction to $\smcx S$
 -- \emph{i.e.} the sheaf with transfers over $S$ represented by $X$.
 
An important property of sheaves with transfers is that
 the abelian premotivic category $\ftr -$ (resp. $\uftr -$) is 
 compatible with the Nisnevich topology on $\ssm$ (resp. $\ssft$) 
 according to Proposition
 \ref{prop:compatible&cd_tr}.
Note moreover that it is compactly geometrically generated.
\end{num}
% \begin{rem}
% Consider an extension of rings $\ZZ \subset \Lambda \subset \Lambda' \subset \QQ$.
% Then, the canonical functor
% $$
% \psh {\Pmorcx S} \rightarrow
% \renewcommand{\Rc}{\Lambda'} \psh {\Pmorcx S}, F \mapsto F \otimes_\Lambda \Lambda'
% $$
% induces a canonical functor
% $$
% \ftr{S} \rightarrow
% \renewcommand{\Rc}{\Lambda'} \ftr{S}.
% $$
% This is evidently compatible with $f^*$, $p_\sharp$,
%  and tensor product.
% Thus we have defined a morphism of abelian premotivic categories:
% \begin{equation} \label{chg_coef_ftr}
% \ftr - \rightleftarrows
%  \renewcommand{\Rc}{\Lambda'} \ftr -.
% \end{equation}
% \end{rem}

\begin{num} \label{num:graph_functor_sh}
We also obtained an adjunction (resp. generalized adjunction)
 of premotivic abelian categories
\begin{align*}
\gamma^*:\sh {} \ssm & \rightleftarrows \ftr -:\gamma_* \\
\gamma^*:\sh {} \ssft & \rightleftarrows \uftr -:\gamma_*.
\end{align*}
Note that in each case $\gamma_*$ is conservative and exact
according to \ref{prop:exists_a^tr}(iv).
\end{num}

\begin{rem} \label{rem:comput_Ext_ftr}
An important application of the existence of the pair
 of adjoint functors $(\gamma^*,\gamma_*)$ is
 the following computation:
 given any complex $K$ of sheaves with transfers over $S$
 and any smooth $S$-scheme $X$,
\begin{align*}
\Hom_{\Der(\ftr S)}(\lrep S X,K[n])
&=\Hom_{\Der(\ftr S)}(\derL \gamma^* \rep S X,K[n]) \\
&=\Hom_{\Der(\sh {} \ssm)}(\rep S X,\gamma_*(K)[n]) \\
& =H^n_\nis(X,\gamma_*(K)).
\end{align*}
This is a generalization of \cite[chap. 5, 3.1.9]{FSV}
 to unbounded complexes and arbitrary bases.
\end{rem}

% Finally, the restriction of the canonical functor
% $a_\et^\tr:\psh {\Pmorcx S} \rightarrow \sh \et {\Pmorcx S}$
% to the category of $t$-sheaves with transfers yields an
% adjunction (resp. generalized adjunction) of premotivic abelian categories
% $$
% a^\tr_\et:\sh {\nis} {\Pmorc}
%  \rightleftarrows \sh {\et} {\Pmorc}:\cO_\et^\tr.
% $$

\begin{num} \label{num:restriction_uftr_ftr}
Let $S$ be a scheme.
Consider the inclusion functor
 $\varphi:\smcx S \rightarrow \sftcx S$.
It induces a functor
$$
\varphi^*:\uftrx {} S \rightarrow \ftrx {} S, F \mapsto F \circ \varphi
$$
which is obviously exact and commute with arbitrary direct sums.
In particular, it commutes with arbitrary colimits.
\end{num}
\begin{lm} \label{lm:enlargement_ftr_uftr}
With the notations above, the functor $\varphi^*$ admits a left adjoint
 $\varphi_!$ such that for any smooth $S$-scheme $X$,
  $\varphi_!(\lrep S X)=\ulrep S X$.
The functor $\varphi_!$ is fully faithful.
\end{lm}
In other words, we have defined an enlargement 
 of premotivic abelian categories (\textit{cf.} definition \ref{df:enlargement})
\begin{equation} \label{enlargement_ftr_uftr}
\varphi_!:\ftrx {} - \rightarrow \uftrx {} -:\varphi^*.
\end{equation}
\begin{proof}
Let $F$ be a sheaf with transfers. Let $\{X/F\}$ be the category
 of representable sheaf $\lrep S X$ over $F$ for a smooth $S$-scheme $X$.
We put 
$$
\varphi_!(F)=\ilim_{\{X/F\}} \ulrep S X.
$$
The adjunction property of $\varphi_!$ is immediate from the Yoneda lemma.
We prove that for any sheaf with transfers $F$, 
 the unit adjunction morphism $F \rightarrow \varphi^* \varphi_!(F)$ 
  is an isomorphism.
As already remarked, $\varphi^*$ commutes with colimits so that
we are restricted to the case where $F=\lrep S X$ which follows by definition.
\end{proof}

\begin{num} Assume now that $t=\cdh$ is the cdh-topology, and $\Pmor=\ssft$ is
 the class of separated morphisms of finite type.
 Recall the topology  $t$ is associated with the \emph{lower cd-structure}
 -- see Example \ref{ex:lower&upper_cd_structures}.
 Then the assumptions of Proposition \ref{prop:compatible&cd_tr}
  with respect to the lower cd-structure are fulfilled 
 according to \cite[4.3.3]{SV1} combined with \cite[4.2.9]{SV1}.
 Thus we get the following result:
\end{num}
\begin{prop}\label{prop:cdh&transfers}
The admissible topology $\cdh$ on $\ssft$ is compatible with transfers.
\end{prop}
As a corollary, we get a generalized premotivic abelian category
 whose fiber over a scheme $S$ is the category 
$\uftrx {\cdh} S$ of $\cdh$-sheaves with transfers on $\ssft$.
 It is compatible with the cdh-topology.
Moreover, the restriction of $a_{\cdh}$ to $\uftr S$ induces
a morphism of generalized premotivic categories; we get the
following commutative diagram of such morphisms:
$$
\xymatrix@=28pt{
\ush {} {-}\ar_{\gamma^*}[d]\ar^-{a^{*}_{\cdh}}[r]
 &  \ush {\cdh} {-}\ar^{\gamma^*_{\cdh}}[d] \\
\uftr -\ar^-{a^{*}_{\cdh}}[r] &  \uftrx {\cdh} -
}
$$

\subsection{Comparison results}  \label{sec:comparison_ftr}

\subsubsection{Change of coefficients}

\begin{num}
Assume the topology $t$ is mildly compatible with transfers
 and consider a localization $\Rc'$ of $\Rc$.

Then the morphism \eqref{eq:Pcor_change_of_rings_bis}
 of $\Pmor$-premotivic categories extends to an adjunction
 of abelian $\Pmor$-premotivic categories:
\begin{equation} \label{eq:fxtr_change_of_rings}
\shg t {\Pmorc} \otimes_\Rc \Rc' \rightleftarrows
\renewcommand{\Rc}{\Lambda'}
 \shg t {\Pmorc}
\end{equation}
Proposition \ref{prop:Pcor_change_of_rings_bis}
 immediately yields the following result:
\end{num}
\begin{prop}\label{prop:fxtr_change_of_rings}
Consider the above notations.
 Then the adjuction \eqref{eq:fxtr_change_of_rings}
 is an equivalence of $\Pmor$-premotivic categories.
\end{prop}

\begin{rem} \label{rem:fxtr&regular}
\renewcommand{\Rc}{\ZZ}
Remark \ref{rem:Pcor&regular} can be extended
 to sheaves with transfers:
 for any regular scheme $S$,
 the category $\ftr S=\shg {\nis} {\smcx S}$ defined here coincides
 with that defined in \cite{Deg7},
 as well as its operations of a $\Pmor$-premotivic
 category when restricted to regular schemes.
\end{rem}

\begin{rem}
In a previous version of this text,
 the preceding proposition was obtained under restrictive
 hypothesis. We have been able to remove these unnecessary assumptions
 thanks to point (3) of Remark \ref{rem:basic_ppties_corr}, which is a consequence
 of Proposition \ref{prop:multiplicity_denominators_bounded}
 (bound on the denominators of intersection multiplicities).
\end{rem}

\subsubsection{Representable $\qfh$-sheaves}

\begin{num}
 Let us denote by $\shqfh S$
 the category of $\qfh$-sheaves of $\Rc$-modules
\index{word}{sheaf!qfhsheaf@$\qfh$-sheaf}
  over $\sft/S$.
Remark that for an $S$-scheme $X$,
 the $\Rc$-presheaf represented by $X$ is not a sheaf
 for the $\qfh$-topology. We denote the associated sheaf
 by $\repx \qfh S X$.
We let $a_{\qfh}$ be the associated $\qfh$-sheaf functor.
Recall that for any $S$-scheme $X$,
 the graph functor (\ref{num:graph_functor_sh})
  induces a morphism of sheaves 
$$
\urep S X \xrightarrow{\gamma_{X/S}} \ulrep S X.
$$
We recall the following theorem of Suslin and Voevodsky
 (see \cite[4.2.7+4.2.12]{SV1}):
\end{num}
\begin{thm}\label{repsheaftransqfhsheaf}
Let $S$ be a scheme such that $\car S \subset \Rc^\times$.
Then, for any $S$-scheme $X$,
 the application of $a_\qfh$ to the map $\gamma_{X/S}$
 gives an isomorphism in $\shqfh S$:
$$
\repx \qfh S X \xrightarrow{\gamma^\qfh_{X/S}} \ulrep S X.
$$
\end{thm}

\begin{num} \label{num:restriction_qfh_uftr}
Assume $\car S \subset \Rc^\times$.
Using the previous theorem, we associate to any $\qfh$-sheaf
$F \in \shqfh S$ a presheaf with transfers
$$
\rho(F):X \mapsto \Hom_{\shqfh S}(\ulrep S X,F).
$$
We obviously get $\gamma^*\rho(F)=F$ as a presheaf over $\sft/S$
 so that $\rho(F)$ is a sheaf with transfers
and we have defined a functor
$$
\rho:\shqfh S \rightarrow \uftr S.
$$
For any $S$-scheme $X$, $\rho(\repx \qfh S X)=\ulrep S X$
 according to the previous proposition.
\end{num}
\begin{cor} \label{cor:inv_univ_homeo_ftr}
Assume $\car S \subset \Lambda^\times$.
Let $f:X' \rightarrow X$ be a morphism of $S$-schemes. \\
If $f$ is a universal homeomorphism, 
 then the map $f_*:\ulrep S {X'} \rightarrow \ulrep S X$
 is an isomorphism in $\uftr S$.
\end{cor}
\begin{proof}
Indeed, according to \cite[3.2.5]{V1}, $\repx \qfh S {X'} \rightarrow \repx \qfh S X$
is an isomorphism in $\shqfh S$ and we conclude by applying the functor 
$\rho$.
\end{proof}

\subsubsection{$\qfh$-sheaves and transfers}

\begin{prop}\label{qfhsheafwithtransfers}
Assume $\car S \subset \Lambda^\times$.
Any $\qfh$-sheaf of $\Rc$-modules over the category of $S$-schemes
of finite type is naturally endowed with a unique structure of
a sheaf with transfers, and any morphism of such $\qfh$-sheaves
is a morphism of sheaves with transfers.

In particular, the $\qfh$-sheafification functor defines a left exact
functor left adjoint to the forgetful functor $\rho:\shqfh S \rightarrow \uftr S$
introduced in \ref{num:restriction_qfh_uftr}.
\end{prop}

\begin{proof}
It follows from Theorem \ref{repsheaftransqfhsheaf} that
the category of $\Rc$-linear finite correspondences is
canonically equivalent to
the full subcategory of the category of $\qfh$-sheaves of
$\Rc$-modules spanned by the objects of shape $\repx \qfh S {X}$
for $X$ separated of finite type over $S$.
The first assertion is thus an immediate consequence
of Theorem \ref{repsheaftransqfhsheaf} and of the (additive)
Yoneda lemma. The fact that the $\qfh$-sheafification functor
defines a left adjoint to the restriction functor $\rho$ is then
obvious, while its left exactness is a consequence of the facts
that it is left exact (at the level of sheaves without transfers)
and that forgetting transfers defines a conservative and exact functor
from the category of Nisnevich sheaves with transfers to the
category of Nisnevich sheaves.
\end{proof}

Recall the following theorem:
\begin{thm} \label{H_etale&nisnevich_rational}
Assume $\Rc$ is a $\QQ$-algebra.
Let $F$ be an \'etale $\Rc$-sheaf on $\sft/S$.
Then for any $S$-scheme $X$, and any integer $i$,
the canonical map
$$
H^i_\nis(X,F) \rightarrow H^i_\et(X,F)
$$
is an isomorphism.
\end{thm}
\begin{proof}
Using the compatibility of \'etale cohomology with projective limits
of schemes, we are reduced to prove that $H^i_\et(X,F)=0$
whenever $X$ is henselian local and $i>0$.
Let $k$ be the residue field of $X$, $G$ its absolute Galois group
 and $F_0$ the restriction of $F$ to $\spec k$.
Then $F_0$ is a $G$-module and according to \cite[8.6]{SGA4}, 
$H^i_\et(X,F)=H^i(G,F_0)$. As $G$ is profinite, this group must be torsion
so that it vanishes by assumption.
\end{proof}

\begin{rem}
The preceding theorem also follows formally from Theorem \ref{ratdescent}.
\end{rem}

\begin{prop}
Assume $\Rc$ is a $\QQ$-algebra.
Let $S$ be an excellent scheme
 and $F$ be a $\qfh$-sheaf of $\Rc$-modules on $\sft/S$.
Then for any geometrically unibranch
\index{word}{scheme!geometrically unibranch}
 $S$-scheme $X$ of finite type,
 and any integer $i$, the canonical map
$$
H^i_\nis(X,F) \rightarrow H^i_\qfh(X,F)
$$
is an isomorphism.
\end{prop}
\begin{proof}
According to \ref{H_etale&nisnevich_rational},
$H^i_\nis(X,F)=H^i_\et(X,F)$.
Let $p:X' \rightarrow X$ be the normalization of $X$.
As $X$ is an excellent geometrically unibranch scheme,
 $p$ is a finite universal homeomorphism. 
It follows from \cite[VII, 1.1]{SGA4} that
 $H^i_\et(X,F)=H^i_\et(X',F)$
and from \cite[3.2.5]{V1} that $H^i_\qfh(X,F)=H^i_\qfh(X',F)$.
Thus we can assume that $X$ is normal, and the result is now exactly
 \cite[3.4.1]{V1}.
\end{proof}

\begin{cor} \label{rational_qfh_points}
Assume $\Rc$ is a $\QQ$-algebra.
Let $S$ be an excellent scheme.
\begin{enumerate}
\item Let $X$ be a geometrically unibranch $S$-scheme of finite type.
For any point $x$ of $X$,
 the local henselian scheme $X_x^h$ is a point for
  the category of sheaves $\shqfh S$ (\emph{i.e.}
  evaluating at $X_x^h$ defines an exact functor).
\item The family of points $X_x^h$ of the previous type
 is a conservative family for  $\shqfh S$.
\end{enumerate}
\end{cor}
\begin{proof}
The first point follows from the previous proposition.
For any excellent scheme $X$, the normalization morphism
 $X' \rightarrow X$ is a $\qfh$-cover. Thus the
 category $\shqfh S$ is equivalent
 to the category of $\qfh$-sheaves on the site made
 of geometrically unibranch $S$-schemes of finite type.
This implies the second assertion.
\end{proof}

\begin{num} Given any scheme $S$,
 we introduce the following composite functor 
 using the notations of \ref{num:restriction_qfh_uftr}
 and \ref{num:restriction_uftr_ftr}:
$$
\psi^*:\shqfh S \xrightarrow \rho \uftr S
 \xrightarrow{\varphi^*} \ftr S.
$$
\end{num}
\begin{thm}\label{enlargeabsheavestrqfhsheaves}
Assume $\Rc$ is a $\QQ$-algebra and
 let $S$ be a geometrically unibranch excellent scheme.
\index{word}{scheme!excellent}
\index{word}{scheme!geometrically unibranch}
Considering the above notation,
 the following conditions are true:
\begin{enumerate}
\item[(i)] For any $S$-scheme $X$ of finite type,
 $\psi^*\big(\repx \qfh S X\big)=\lrep S X$.
\item[(ii)] The functor $\psi^*$
 admits a left adjoint $\psi_!$.
\item[(iii)] For any smooth $S$-scheme $X$,
 $\psi_!\big(\lrep S X\big)=\repx \qfh S X$.
\item[(iv)] The functor $\psi^*$ is exact and preserves colimits.
\item[(v)] The functor $\psi_!$ is fully faithful.
\end{enumerate}
\end{thm}
According to property (iii), the functor $\psi_!$ commutes
 with pullbacks.
In other words, we have defined an enlargement of abelian premotivic categories
\index{word}{premotivic!enlargement of ---- category}
(\textit{cf.} definition \ref{df:enlargement})
 over the category of (noetherian) geometrically unibranch schemes:
\begin{equation} \label{qfh_enlargement}
\psi_!:\ftr - \rightleftarrows \shqfh -:\psi^*
\end{equation}
\begin{proof}
Point (i) follows from Theorem \ref{repsheaftransqfhsheaf}.
Recall the enlargement of \eqref{enlargement_ftr_uftr}:
$$
\varphi_!:\ftr - \rightarrow \uftr -:\varphi^*.
$$
We define the functor $\psi_!$ as the composite:
$$
\ftr S \xrightarrow{\varphi_!} \uftr S
 \xrightarrow{\gamma^*} \ush {} {S}
 \xrightarrow{a_{\qfh}} \shqfh S.
$$
According to the properties of the functors in this composite,
 $\psi_!$ is exact and preserves colimits.
Moreover, for any smooth $S$-scheme $X$, as $\lrep S X$
is a $\qfh$-sheaf over $\sft/S$ according to \ref{prop:lrep+et,qfh}, 
$\psi_!(\lrep S X)=\repx \qfh S X$ which proves (iii).
Property (ii) follows from (iii) and the fact $\psi_!$
commutes with colimits, while the sheaves
$\lrep S X$ for $X/S$ smooth generate $\ftr S$.

For any smooth $S$-scheme $X$, $\Gamma(X;\psi^*(F))=F(X)$.
Thus the exactness of $\psi^*$ follows from Corollary
\ref{rational_qfh_points} --- here we use the
 assumption that $S$ is geometrically unibranch and excellent,
 as it implies that $X$ satisfies the same properties, so that we can 
 apply the mentioned corollary.
 As $\psi^*$ obviously preserves direct sums, we get (iv).

To check that for any sheaf with transfers $F$ the unit map
$F \rightarrow \psi^* \psi_!(F)$ is an isomorphism, we thus are reduced
to the case where $F=\lrep S X$ for a smooth $S$-scheme $X$
 which follows from (i) and (iii).
\end{proof}
\section{Motivic complexes} \label{sec:voevodsky}

\begin{assumption} \label{num:voevodsky_convention_sch}
In this section, $\sch$ is the category
 of noetherian finite dimensional schemes.
 It is adequate in the sense of \ref{num:assumption1_sch}.
 Given a scheme $S$, we denote by $\ssm_S$
 the category smooth separated $S$-schemes of finite type.
 It is admissible in the sense of \ref{defclassP}.

\renewcommand{\Rc}{\Lambda}
We fix a ring of coefficients $\Rc$.
\end{assumption}

\subsection{Definition and basic properties}

\subsubsection{Premotivic categories}

\renewcommand{\Rc}{\Lambda}

According to Proposition \ref{prop:compatible&cd_tr}
 and Corollary \ref{cor:compatibility_tr&top},
 the abelian premotivic category $\ftr -$ constructed
 in \ref{df:premotivic_ftr} is compatible with Nisnevich
 topology. Thus we can apply to it the general definitions
  of section \ref{sec:fibred_derived}. This gives the following
  definition:
\begin{df} \label{df:Nis_DMe&DM}
We define the ($\Rc$-linear) category of \emph{motivic complexes}
\index{word}{motivic complex}
 (resp. \emph{stable motivic complexes} or simply \emph{motives})
\index{word}{motivic complex!stable}
\index{word}{motivic complex!stable|seealso{motive}}
\index{word}{motive}
following Definition \ref{df:effective_triangulated_premotives}
 (resp. Definition \ref{df:triangulated_premotives})
as
\begin{align*}
& \DMe=\DMue\left(\ftr -\right) \\
\text{resp. } & \DM=\DMu\left(\ftr -\right).
\end{align*}
\index{notat}{DMLambdaeff@$\DMe$}
\index{notat}{DMLambda@$\DM$}
Given a scheme $S$,
 we will put: $\DMex S=\DMe(S)$, $\DMx S=\DM(S)$.
\end{df}

\begin{num}
Let us unfold the preceding definition.
 Given a scheme $S$ in $\sch$,
  the triangulated category $\DMex S$ is equal to
  the $\AA^1$-localization of the derived category
  $\Der(\ftr S)$ of the category of sheaves with transfers
  over $S$.

Given a smooth scheme $S$-scheme $X$ of finite type,
 we have denoted by $\lrep S X$ the sheaf with transfers
 represented by $X$ over $S$.
 We will see this sheaf as an object of $\DMex S$,
  as a complex concentrated in degree $0$,
  and call it the effective motivic complex associated
  with $X/S$.

Recall the following operations as part of the 
 premotivic structure:
\begin{itemize}
\item Given any morphism $f:T \rightarrow S$ in $\sch$,
 there exists an adjunction of the form:
$$
\derL f^*:\DMex S \rightleftarrows \DMex T:\derR f_*\, .
$$
\item Given a separated smooth morphism of finite type
 $f:T \rightarrow S$ in $\sch$,
 there exists an adjunction of the form:
$$
\derL f_\sharp:\DMex S \rightleftarrows \DMex T:f^*=\derL f^*\, .
$$
\item Given any noetherian finite dimensional scheme $S$,
 the category $\DMex S$ is symmetric closed monoidal.
\end{itemize}
These operations are subject to the properties of a premotivic category:
 functoriality, smooth base change formula,
  smooth projection formula
   -- see section \ref{sec:pfibred} for more details.
 By construction, the triangulated premotivic category $\DMe$
  satisfies the homotopy property and the Nisnevich descent properties.

By construction  (cf. \eqref{eq:adj_DMue-DMu}), we get an adjunction
 of triangulated premotivic categories
\begin{equation} \label{eq:adjunction_DMeDM}
\sus:\DMe \rightleftarrows \DM:\lop.
\end{equation}
Considering the \emph{Tate motivic complex}
\index{word}{Tate!motivic complex}
\begin{equation}
\lrepNP_S(1):=\lrep S {\PP^1_S/\{1\}},
\end{equation}
the object $\sus(\lrepNP_S(1))$ is $\otimes$-invertible
 in $\DMx S$ and this property characterizes uniquely the homotopy
 category $\DMx S$ -- see Remark \ref{rem:stabilisation&universality}.
 Given a smooth separated $S$-scheme $X$ of finite type,
  we put: 
$$\mot S X:=\sus \lrep S X$$
\index{notat}{MSX@$\mot S X$}
 and simply call it the 
 motive associated with $X/S$.
Usually we denote by $\un_S$ the unit of the monoidal category $\DMx S$.

By construction, the premotivic category $\DM$ satisfies
 the homotopy, stability and Nisnevich descent properties
 (see Paragraph \ref{num:basic_DMu}).
\end{num}

\begin{ex} \label{ex:DM_field&regular}
\renewcommand{\Rc}{\ZZ}
\begin{itemize}
\item Let $k$ be a perfect field.
 Then $\DMex k$ contains as a full subcategory
  the category $\mathrm{DM}_-^{\mathit{eff}}(k)$ defined by Voevodsky
  (cf \cite[Chap. 5]{FSV}).
 This is the content of the proof of \cite[Chap. 5, Prop. 3.2.3]{FSV}.
 Indeed, recall from Paragraph \ref{num:AA^1-equivalences_weak&strong}
  that $\DMex k$ is equivalent to the full subcategory
   of $\Der(\ftr k)$
  made by the complexes which are $\AA^1$-local.
 Over a perfect field, Theorem 3.1.12 of \cite[Chap. 5]{FSV} implies that
 a complex of sheaves with transfers is $\AA^1$-local
 if and only if its homotopy sheaves are $\AA^1$-invariant.
\item Let $S$ be a regular scheme.
 The triangulated categories $\DMex S$ and $\DMx S$
 introduced here coincide with that constructed in \cite{CD1}.
 The same is true concerning the operations of premotivic triangulated
  categories (see Remark \ref{rem:fxtr&regular}).
\end{itemize}
\end{ex}

\begin{num}
Let $\Rc'$ be a localization of $\Rc$.
The premotivic adjunction
\begin{equation}
\ftr - \otimes_\Rc \Rc' \rightleftarrows
\renewcommand{\Rc}{\Lambda'}
 \ftr -
\end{equation}
obtained as a particular case of \eqref{eq:fxtr_change_of_rings}
 gives the following adjunctions of triangulated premotivic categories:
\begin{equation} \label{eq:changeofcoef_DM}
\begin{split}
\DM \otimes_\Lambda \Lambda' & \rightleftarrows
 \renewcommand{\Rc}{\Lambda'} \DM, \\
\DMe \otimes_\Lambda \Lambda' & \rightleftarrows
 \renewcommand{\Rc}{\Lambda'} \DMe.
\end{split}
\end{equation}
Proposition \ref{prop:fxtr_change_of_rings} immediately yields the following result:
\end{num}
\begin{prop} \label{prop:changeofcoef_DM}
The premotivic adjunctions \eqref{eq:changeofcoef_DM}
 are equivalences of triangulated premotivic categories.
\end{prop}
In other words, for any scheme $S$,
 the triangulated monoidal category
\renewcommand{\Rc}{\Lambda'}
 $\DMx S$ (resp. $\DMex S$)
 is the naive localization of the category 
\renewcommand{\Rc}{\Lambda}
 $\DMx S$ (resp. $\DMex S$) with respect
 to integers invertible in $\Lambda'$.

\bigskip

\subsubsection{Constructible and geometric motives}

\begin{num} \label{num:DM_generators}
The premotivic triangulated category $\DMe$ is geometrically generated:
 given any scheme $S$, the essentially small set $\mathcal G^{eff}_S$
 of motivic complexes of
 the form $\lrep S X$ for a smooth separated $S$-scheme $X$ of finite type
 form a set of generators in the triangulated category $\DMex S$.

Similarly, the premotivic triangulated category $\DM$ is $\ZZ$-generated
 where $\ZZ$ is the set of twists corresponding to the Tate twist:
 given any scheme $S$,
 the essentially small set $\mathcal G_S$
 of motives of the form $\mot S X(n)$
 for a smooth separated $S$-scheme $X$ of finite type
 and an integer $n \in \ZZ$
 form a set of generators in the triangulated category $\DMx S$.

Following the general conventions about premotivic triangulated
 category (Definition \ref{df:tau-geometric}), we define
 the notion of constructibility for motives as follows:
\end{num}
\begin{df} \label{df:contructible_V_motives}
Given any scheme $S$,
 we define the category of \emph{constructible motives}
\index{word}{constructible!motive}
\index{word}{motive!constructible|see{constructible}}
 (resp. \emph{constructible motivic complexes}) over $S$
\index{word}{constructible!motivic complex}
\index{word}{motivic complex!constructible|see{constructible}}
  as the thick triangulated subcategory
 of $\DMx S$ (resp. $\DMex S$)
 generated by $\mathcal G_S$ (resp. $\mathcal G^{eff}_S$).
 We denote it by $\DMcx S$ (resp. $\DMcex S$).
\end{df}

\begin{rem} 
Recall that $\DMc$ (resp. $\DMce$) is
 $\sm$-fibred monoidal subcategory of $\DM$ (resp. $\DMe$) over $\sch$.
In other words, constructible motives (resp. motivic complexes)
 are stable by the operations $f^*$, $p_\sharp$ for $p$ smooth
 and tensor product.
 This is obvious from definitions.
\end{rem}

\begin{num} Let $S$ be a scheme.
Consider the triangulated subcategory $\V_S$ of $\K^b(\smcx S)$ generated
by complexes of one the following forms~:
\begin{enumerate}
\item
for any distinguished square 
$\xymatrix@=10pt{W\ar^k[r]\ar_g[d] & V\ar^f[d] \\ U\ar^j[r] & X}$ 
of smooth $S$-schemes,
$$
[W]
 \xrightarrow{g_*-k_*} [U] \oplus [V]
 \xrightarrow{j^*+f^*} [X]
$$
\item for any smooth $S$-scheme $X$, 
$p:\AA^1_X \rightarrow X$ the canonical projection.
$$
[\AA^1_X] \xrightarrow{p_*} [X].
$$
\end{enumerate}
\end{num}

\begin{df} \label{df:gm_V_mot}
We define the category $\DMgmex S$
\index{notat}{DMgmeffS@$\DMgmex S$}
 of \emph{geometric effective motives}
\index{word}{motive!geometric effective}
over $S$ as the pseudo-abelian envelope of the triangulated category
$$
\K^b(\smcx S)/\V_S.
$$
We define the category $\DMgm(S)$
\index{notat}{DMgmS@$\DMgmx S$}
 of \emph{geometric motives}
\index{word}{motive!geometric}
 over $S$
as the triangulated category obtained from $\DMgmex S$ by
formally inverting the Tate complex
$$
[\PP^1_S] \rightarrow [S].
$$
\end{df}

\begin{rem}
The categories of geometric motives (resp. effective geometric motives)
 over an arbitrary base, as defined here, already appears
 in the work of Ivorra \cite[sec. 1.3]{Ivo}.
\end{rem}

\begin{num}
According to this definition,
 we can construct for any scheme $S$ a commutative diagram of functors: 
\begin{equation} \label{eq:embedding_gm_motives}
\begin{split}
\xymatrix@R=10pt@C=20pt{
\DMgmex S\ar[r]\ar[d] & \DMex S\ar^{\sus}[d] \\
\DMgmx S\ar[r] & \DMx S
}
\end{split}
\end{equation}
where the right vertical map is the left adjoint
 of \eqref{eq:adjunction_DMeDM}.

Recall from Remark \ref{rem:cd_bounded_ftr}
 that the Nisnevich topology is bounded in $\ftr -$.
 Thus, as a corollary of 
Proposition \ref{prop:compact_DMue},
 Corollary \ref{cor:DMue_compact&constructible}
 and Corollary \ref{cor:geometric_premotives}
 we get the following result:
\end{num}
\begin{thm} \label{thm:DM_constructible_geometric}
The horizontal functors of the square \eqref{eq:embedding_gm_motives}
 are fully faithful and their essential images consist of constructible
 objects in the sense of Definition \ref{df:contructible_V_motives}.

Given any motive (resp. motivic complex) $\mathcal M$ over $S$,
 the following conditions are equivalent:
\begin{enumerate}
\item[(i)] $\mathcal M$ is geometric (\emph{i.e.} in the image
 of the horizontal map of diagram \eqref{eq:embedding_gm_motives}),
\item[(ii)] $\mathcal M$ is constructible,
\index{word}{constructible!motive}
\item[(iii)] $\mathcal M$ is compact.
\index{word}{compact}
\end{enumerate}

The triangulated category $\DMx S$
 (resp. $\DMex S$)
 is compactly generated. More precisely,
 the objects of the set of generators $\mathcal G_S$
 (resp. $\mathcal G^{eff}_S$) defined in Paragraph \ref{num:DM_generators}
  are compact.
\end{thm}

\begin{rem} \label{rem:DMgm_field}
If $S=\spec k$ is the spectrum of a perfect field,
 then the categories $\DMgmx S$ and $\DMgmex S$
 coincide with the categories introduced by Voevodsky
 in \cite[chap. 5, Sec. 2.1]{FSV}.
The above theorem is a generalization of 
 \cite[chap. 5, Th. 3.2.6]{FSV} to an arbitrary base
 (and the non-effective case).
\end{rem}

\subsubsection{Enlargement, descent and continuity}

\begin{num} We can apply the definitions of section \ref{sec:fibred_derived}
 to the generalized abelian premotivic category $\uftr -$ constructed in \ref{df:premotivic_ftr}
 \end{num} 
\begin{df}
We define the ($\Rc$-linear) category of
 \emph{generalized motivic complexes}
\index{word}{motivic complex!generalized}
 (resp. \emph{generalized motives})
\index{word}{motive!generalized}
following definition \ref{df:triangulated_premotives}
 (resp. definition \ref{df:effective_triangulated_premotives})
as
\begin{align*}
&\uDMe=\DMue\left(\uftr -\right) \\
\text{resp. } & \uDM=\DMu\left(\uftr -\right).
\end{align*}
\index{notat}{DMLambdaeff@$\uDMe$}
\index{notat}{DMLambda@$\uDM$}
\end{df}

\begin{num}
The advantage of this definition
 is that any separated $S$-scheme $X$ of finite type
 defines a generalized motivic complex, given by the sheaf with transfers
 $\ulrep S X$ seen as a complex concentrated in degree $0$
  (see Definition \ref{df:premotivic_ftr}).

The category $\uDMe$, as a generalized premotivic category,
 admits the following operations:
\begin{itemize}
\item Given any morphism $f:T \rightarrow S$ in $\sch$,
 there exists an adjunction of the form:
$$
\derL f^*:\uDMex S \rightleftarrows \uDMex T:\derR f_*\, .
$$
\item Given a separated morphism $f:T \rightarrow S$
 of finite type in $\sch$ (non necessarily smooth),
 there exists an adjunction of the form:
$$
\derL f_\sharp:\uDMex S \rightleftarrows \uDMex T:f^*=\derL f^*\, .
$$
\item Given any noetherian finite dimensional scheme $S$,
 the category $\uDMex S$ is symmetric closed monoidal.
\end{itemize}
These operations satisfy the properties of a generalized premotivic
 category for which we refer the reader to
  section \ref{sec:premotivic_convention}.

As in the non generalized case,
 we get from the general construction (see \eqref{eq:adj_DMue-DMu})
 an adjunction of triangulated generalized premotivic categories
\begin{equation} \label{eq:adjunction_DMeDM_generalized}
\sus:\uDMe \rightleftarrows \uDM:\lop.
\end{equation}
To any separated $S$-scheme $X$ of finite type,
 we associate a generalized motive as:
$$\umot S X:=\sus \ulrep S X.$$
\index{notat}{MSX@$\umot S X$}

By construction,
 the generalized premotivic category $\uDMe$ (resp. $\uDM$)
 satisfies the homotopy property, Nisnevich descent property
 (resp. and stability property).
\end{num}

\begin{num}
By virtue of Remark \ref{rem:cd_bounded_ftr},
 the Nisnevich topology is bounded in $\uftr -$.
 Therefore, as a corollary of Proposition \ref{prop:compact_DMue}
  (resp. Corollary \ref{cor:DMue_compact&constructible}),
 we obtain in particular that $\uDMex S$ (resp. $\uDMx S$)
 is compactly generated, with the essentially small family
 of objects $\ulrep S X$ (resp. $\umot S X(n)$)
 for a separated $S$-scheme of finite type $X$
  (resp. and an integer $n \in \ZZ$) as compact generators.

Recall that for any scheme $S$, the obvious restriction functor
$$\varphi^*:\ftr S \rightarrow \uftr S$$
 admits a left adjoint $\varphi_!$ which is fully faithful
  (Lemma \ref{lm:enlargement_ftr_uftr}).
 Moreover, the adjoint pair $(\varphi_!,\varphi^*)$
  satisfies the assumption of Proposition \ref{prop:enlargement_eff}
  so that applying Corollary \ref{cor:derived_enlargement}
  gives the following proposition:
\end{num}
\begin{prop}\label{prop:enlargement_DM_Voe}
Given any scheme $S$,
 the adjoint pair $(\varphi_!,\varphi^*)$ can be derived and induces
 the following pair of adjoint functors
\begin{equation} \label{eq:enlargement_DM_Voe}
\begin{split}
\varphi_!:\DMx S &\rightleftarrows \uDMx S:\varphi^*, \\
\text{resp. } \varphi_!:\DMex S &\rightleftarrows \uDMex S:\varphi^*,
\end{split}
\end{equation}
such that $\varphi_!$ is fully faithful.

More generally, the family of these adjunctions for a noetherian
 finite dimensional scheme $S$ defines
 an enlargement of premotivic categories
\index{word}{premotivic!enlargement of ---- category}
  (Definition \ref{df:enlargement}).
\end{prop}
The abuse of notations is justified because of the following
 essentially commutative diagram of functors:
\begin{equation} \label{eq:enlargement_DM&eff_Voe}
\begin{split}
\xymatrix@R=10pt@C=20pt{
\DMe\ar^{\sus}[r]\ar_{\varphi_!}[d] & \DM\ar^{\varphi_!}[d] \\
\uDMe\ar^{\sus}[r] & \uDM
}
\end{split}
\end{equation}
Recall that, given a smooth separated $S$-scheme $X$,
 we have the relation:
$$\varphi_!(\mot S X)=\umot S X.$$

\begin{rem}
Beware that the functor $\varphi^*$ is far from being conservative.
The following example was suggested by V.~Vologodsky: let $Z$ be a nowhere dense closed subscheme
of $S$. Then $\varphi^*(\umot S Z)=0$.
In fact, one can see that $\DMx S$ is a localization of the category $\uDMx S$
 with respect to the objects $\uM$ such that $\varphi^*(\uM)=0$.
\end{rem}

\renewcommand{\Rc}{\QQ}
\begin{num}\label{adjunctionDM_DMqfh}
With rational coefficients, the preceding proposition can be
 refined.
 Recall that the $\qfh$-sheafification functor (\ref{qfhsheafwithtransfers})
 induces by \ref{num:functoriality_DMu} a premotivic adjunction
$$\underline{\alpha}^* : \uDM
 \rightleftarrows\uDMV_{\qfh,\Rc}: \underline{\alpha}_* \, .$$
\end{num}

\begin{thm}\label{embedDMDMqfh}
If $S$ is a geometrically unibranch excellent
\index{word}{scheme!geometrically unibranch}
\index{word}{scheme!excellent}
 noetherian scheme
 of finite dimension
 then the following composite functor
 
$$
\underline{\alpha}^* \varphi_!:\DMx S\To \uDMV_{\qfh,\Rc}(S)
$$
is fully faithful.
\end{thm}
\begin{proof}
Note that $\DMex S$ and $\DMue(\sh \qfh S)$
are compactly generated; see example \ref{ex:Zar&Nis=bounded}
 and Proposition \ref{prop:compact_DMue}. Hence
this corollary follows from Theorem \ref{enlargeabsheavestrqfhsheaves}
 and Proposition \ref{prop:enlargement}.
\end{proof}
\renewcommand{\Rc}{\Lambda}

\begin{rem}
Recall this theorem can be rephrased by saying that 
 motives over $S$ satisfies $\qfh$-descent
\index{word}{descent!qfhdescent@$\qfh$-descent}
 -- see Remark \ref{rem:descent&derived_P-premotivic}
 and more generally Section \ref{sec:theorie_descente}.
In the next section,
 we will give applications of this fact to motivic cohomology.
\end{rem}

\begin{thm} \label{thm:DM_continuity}
The following assertions hold:
\begin{enumerate}
\item The triangulated premotivic categories
 $\DMe$ and $\DM$ are weakly continuous (Definition \ref{df:continuous}). 
\item The generalized triangulated premotivic categories
 $\uDMe$ and $\uDM$ are weakly continuous.
\end{enumerate}
\end{thm}
\begin{proof}
Note that Proposition \ref{prop:ftr_property(C)}
 shows precisely that the generalized premotivic abelian
 category $\uftr -$ satisfies Property (wC) of Paragraph
 \ref{num:continuity_derived_premotivic}.
 Therefore, the assertion (2) follows from
 Propositions \ref{prop:continuity_gen-A^1_premotivic}
 and \ref{prop:continuity_gen-stable-A^1_premotivic}.

Moreover, the assertion (1) follows from
 Corollary \ref{cor:continuity_stable-A^1_premotivic}
 given the enlargement obtained in Proposition
 \ref{prop:enlargement_DM_Voe}.
\end{proof}

\begin{ex}
From the previous theorem and Proposition \ref{continuityconstructible},
 we obtain in particular that for any pro-scheme $(S_\alpha)_{\alpha \in A}$
 with affine and dominant transition map such that 
 $S=\varprojlim_{\alpha \in A} S_\alpha$ is noetherian finite dimensional,
 there exists canonical equivalences of categories:
\begin{align*}
2\mbox{-}\varinjlim_\alpha \big(\DMgme(S_\alpha)\big) & \To \DMgme(S), \\
2\mbox{-}\varinjlim_\alpha \big(\DMgm(S_\alpha)\big) & \To \DMgm(S).
\end{align*}
This result generalizes \cite[4.16]{Ivo}.
\end{ex}

\subsection{Motivic cohomology} \label{sec:motivic_coh}

\subsubsection{Definition and functoriality}

\begin{df}
 Let $S$ be a scheme and $(n,m) \in \ZZ^2$ be a couple of integers.
We define the \emph{motivic cohomology}
\index{word}{cohomology!motivic}
 of $S$ in degree $n$ and twist $m$ 
 with coefficients in $\Rc$ as the $\Rc$-module
$$
\HHm n m S=\Hom_{\DMx S}\big(\un_S,\un_S(m)[n]\big).
$$
\index{notat}{HMnmS@$\HHm n m S$}
Assuming $m \geq 0$, we define the \emph{effective motivic cohomology}
\index{word}{cohomology!effective motivic}
 of $S$
 in degree $n$ and twist $m$ 
 with coefficients in $\Rc$ as the $\Rc$-module
$$
\HHme n m S=\Hom_{\DMex S}\big(\lrepNP_S,\lrepNP_S(m)[n]\big).
$$
\index{notat}{HMeffnmS@$\HHme n m S$}
\end{df}
Motivic cohomology (resp. effective motivic cohomology)
 is contravariant with respect to morphisms of schemes
 and the monoidal structure on $\DM$ (resp. $\DMe$)
 defines a ring structure compatible with pullbacks:
 given two cohomology classes:
$$
\alpha:\un_S \rightarrow \un_S(m)[n],
 \alpha':\un_S \rightarrow \un_S(m')[n'],
$$
one simply put:
$$
\alpha.\alpha'=\alpha \otimes_S \alpha'.
$$
The link between motivic cohomology and effective motivic cohomology
 is provided by Proposition \ref{spanierwhitehead1}.
 Given any scheme $S$ and any couple of integers $(n,m) \in \ZZ^2$,
 one has a canonical identification:
$$
\HHm n m S
=\ilim_{r>>0} \Hom_{\DMex S}\big(\lrepNP_S(r),\lrepNP_S(m+r)[n]\big).
$$

\begin{num} \label{num:changeofcoef_motivic_coh}
Let $\Rc'$ be a localization of $\Rc$.
Then using the left adjoint of the premotivic adjunction
 \eqref{eq:changeofcoef_DM},
 we get a canonical morphism
$$
\HHm n m S \otimes_\Rc \Rc'
\renewcommand{\Rc}{\Lambda'} \rightarrow \HHm n m S.
$$
It is obviously compatible with pullbacks and the product structure.
According to Proposition \ref{prop:changeofcoef_DM}, this map is an isomorphism.
\end{num}

\begin{ex} \label{ex:motivic_coh&Chow}
\renewcommand{\Rc}{\ZZ}
Let $k$ be a perfect field.
Given any smooth separated $k$-scheme $S$ of finite type,
 with structural morphism $f$,
 and any pair of integers $(n,m) \in \ZZ^2$,
 motivic cohomology as defined in the previous definition
 coincide with motivic cohomology as defined by Voevodsky
 in \cite[chap. 5]{FSV} according to the following computation
 and Remark \ref{rem:DMgm_field}:
\begin{align*}
\HHm n m X&=\Hom_{\DMx X}(\un_X,\un_X(m)[n])
=\Hom_{\DMx X}(\un_X,f^*(\un_k)(m)[n]) \\
&=\Hom_{\DMx k}(\derL f_\sharp(\un_X),\un_k(m)[n])
=\Hom_{\DMx k}(M_k(X),\un_k(m)[n]) \\
&=\Hom_{\DMgmx k}(M_k(X),\un_k(m)[n]).
\end{align*}
In particular, it coincides with higher Chow groups
\index{word}{group!higher Chow}
\index{word}{cohomology!higher Chow group|see{group}}
 (cf \cite{allagree}) according to the following formula:
$$
\HHm n m X=CH^{m}(X,2m-n).
$$
Recall in particular the following computations:
$$
\HHm n m X=\begin{cases}
\ZZ^{\pi_0(X)} & \text{if } n=m=0, \\
\GG(X) & \text{if } n=m=1, \\
CH^m(X) & \text{if } n=2m, \\
0 & \text{if } m<0, n>\min(m+\dim(X),2m)
\end{cases}
$$
where $CH^m(X)$ is the usual Chow group
\index{word}{group!Chow}
\index{word}{cohomology!Chow group|see{group}}
 of $m$-codimensional cycles in $X$.

Note we will extend the identification of motivic cohomology
 as defined in the previous definition with the general
 version defined by Voevodsky -- \cite{V3} --
 in section \ref{sec:motivic_coh_spectrum}.
\end{ex}

\begin{num} \label{num:corr_cohm}
Consider a separated morphism $p:X \rightarrow S$ of finite type.
Recall from the $\ssft$-fibred structure of $\uDM$
 that $\umot S X=\derL p_\sharp p^*(\un_S)$.
Using the adjunction property of the pair $(\derL p_\sharp,p^*)$,
 we easily get:
\begin{equation} \label{eq:mot_coh&umot}
\begin{split}
\HHm n m X
&=\Hom_{\DMx X}\big(\un_X,\un_X(m)[n]\big)
=\Hom_{\uDMx X}\big(\un_X,\un_X(m)[n]\big) \\
&=\Hom_{\uDMx S}\big(\umot S X,\un_S(m)[n]\big).
\end{split}
\end{equation}
In particular, given any finite $S$-correspondence
$\alpha:X \doto Y$
between separated $S$-schemes of finite type, we obtain a pullback
$$
\alpha^*:\HHm n m Y \rightarrow \HHm n m X
$$
which is, among other properties, functorial with respect to composition
of finite $S$-correspondences and extends the natural contravariant functoriality
of motivic cohomology.

In particular, given any finite $\Rc$-universal morphism
 $f:Y \rightarrow X$, 
we obtain a pushout
$$
f_*:\HHm n m Y \rightarrow \HHm n m X
$$
by considering the transpose of the graph of $f$.
\end{num}
\begin{prop}
 Let $f:Y \rightarrow X$ be a finite $\Rc$-universal morphism
\index{word}{morphism!finite $\Rc$-universal}
  of schemes.
 Assume $X$ is connected and let $d>0$ be the degree of $f$
  (\textit{cf.} \ref{df:degree_corr}).
Then for any element $x \in \HHm n m X$,
 $f_*f^*(x)=d.x$.
\end{prop}
\noindent This is a simple application of Proposition \ref{prop:degree_corr}.
We left to the reader the exercise to state projection and base change formulas for this pushout.

\begin{ex}
Let $f:Y \rightarrow X$ be a finite morphism.
Recall that $f$ is $\Rc$-universal in the following
 particular cases:
\begin{itemize}
\item $f$ is flat (see Example \ref{ex:universal&flat});
\item $X$ is regular and $f$ sends the generic points of $Y$ 
 to generic points of $X$ (see Corollary \ref{cor:universal&regular}).
\end{itemize}
In particular, motivic cohomology is covariant with respect to this
 kind of finite morphisms.
\end{ex}

Another important application of the generalized motives
 is obtained using the Corollary \ref{cor:inv_univ_homeo_ftr}:
\begin{prop}
Let $f:X' \rightarrow X$ be a separated universal homeomorphism 
of finite type.
\index{word}{homeomorphism, universal}
Assume that $\car X \subset \Rc^{\times}$.
Then the pullback functor
$$
\HHm n m X \rightarrow \HHm n m {X'}
$$
is an isomorphism.
\end{prop}

\begin{rem}
The preceding considerations hold similarly
 for the effective motivic cohomology.
\end{rem}

\begin{ex}
In characteristic $0$,
 motivic cohomology (effective and non-effective)
 is invariant under semi-normalization (\cite{Swan}).

When restricted to excellent geometrically unibranch scheme $X$,
 motivic cohomology (effective and non-effective)
  is invariant under normalization.
 Indeed, the normalization $X' \rightarrow X$ of such a scheme
 is a universal homeomorphism (\cite[$IV_0$, 23.2.2]{EGA4})
 of finite type.
\end{ex}

\subsubsection{Effective motivic cohomology in weight $0$ and $1$}

\begin{num}
Let $S$ be a scheme and $X$ a smooth $S$-scheme.
For any subscheme $Y$ of $X$, we denote by $\lrep S {X/Y}$
the cokernel of the canonical morphism of sheaf with transfers
$\lrep S Y \rightarrow \lrep S X$.
As this morphism is a monomorphism, we obtain a canonical 
distinguished triangle in $\DMex S$
$$
\lrep S Y \rightarrow \lrep S X \rightarrow \lrep S {X/Y}
 \rightarrow \lrep S X[1].
$$
Using this notation and according to Definition \ref{df:Tate_twist},
 the Tate motivic complex is defined as:
 $\lrepNP_S(1)=\lrep S {\PP^1_S/\{\infty\}}[-2]$.

The following computation is classical:
$$
\lrepNP_S(1)=\lrep S {\PP^1_S/\AA^1_S}[-2]
 =\lrep S {\AA^1_S/\GG}[-2];
$$
the first identification follows from homotopy invariance
 and the second one by Nisnevich descent
  (cf. Prop. \ref{prop:BG_property_derived_premotivic}).
\end{num}
\begin{prop}
\label{computing_first_tate_twist}
Suppose $S$ is a normal scheme.

Then the sheaf on $\sm_S$ represented by $\GG$ admits
 a canonical structure of a sheaf with transfers
 and there is a canonical isomorphism in $\DMex S$:
$$
\GG \otimes_\ZZ \Rc\xrightarrow{\ \simeq\ } \lrepNP_S(1)[1].
$$
\end{prop}
\begin{proof}
Let $U$ be an open subscheme of $\AA^1_S$
and $X$ be a smooth $S$-scheme.
Note that $X$ is normal according to \cite[18.10.7]{EGA4}.
Consider a cycle 
$$
\alpha=\sum_i n_i.\acycl{Z_i}
$$
of $X \times_S U$ with $n_i \in \Rc$ and $Z_i$ irreducible
 finite and dominant over an irreducible component of $X$.
Then $Z_i$ is a divisor in $X \times_S U$
 and according to \cite[21.14.3]{EGA4}, it is flat over $X$.
In other words, $\alpha$ is a Hilbert cycle
 which implies it is $\Rc$-universal (Example \ref{ex:universal&flat}).
 As a consequence, we obtain the equality
$$
H^i\Gamma(X;\suslin\lrep S U)=H_{-i}^{sing}(X \times_S U/X) \otimes_\ZZ \Rc
$$
where the functor $\suslin$ is the associated Suslin singular complex
 (see \eqref{eq:Sulin_singular_complex})
 and the right-hand side is the Suslin homology of $X \times_S U/X$
(\textit{cf.} \cite{SV1}).

Suppose in addition that $X$ and $U$ are affine and let $Z=\PP^1_S-U$.
According to a theorem of Suslin and Voevodsky
 (\textit{cf.} \cite[th. 3.1]{SV1}),
$$
H_{-i}^{sing}(X \times_S U/X)=
\left\{ \begin{array}{ll}
\pic(X \times_S \PP^1_S,X \times_S Z) & \text{if } i=0 \\
0 & \text{otherwise;}
\end{array}\right.
$$
the group on the left-hand side is the
 \emph{relative Picard group}.
\index{word}{group!relative Picard}
In particular, the complex $\suslin\lrep S U$,
 seen as a complex of presheaves with transfers,
 is concentrated in cohomological degree $0$
 and its $0$-th  cohomology is the presheaf
 $X \mapsto \pic(X \times_S \PP^1_S,X \times_S Z) \otimes_\ZZ \Rc$.

Consider the following exact sequence of
 presheaves with transfers:
$$
0 \To \lrepNP_S(\GG)\To\lrepNP_S(\AA^1_{S})
 \To \plrepNP_S(\AA^1_S/\GG) \To 0.
$$
Applying the functor $\suslin$ to it, relatively to the
 category of complexes of presheaves with transfers,
 we obtain a distinguished triangle in $\Der(\ptr S)$:
$$
\suslin \lrepNP_S(\GG)
 \To \suslin \lrepNP_S(\AA^1_{S})
 \To \suslin \plrepNP_S(\AA^1_S/\GG) \To \suslin \lrepNP_S(\GG)[1].
$$
Taking the associated long exact sequence of
 cohomology presheaves,
 we obtain that the complex of presheaves with transfers
 $\suslin \plrepNP_S(\AA^1_S/\GG)$ is concentrated in cohomological
 degree $0$ and $-1$, and we get an exact sequence of presheaves:
\begin{align*}
0 \To \hat H^{-1}\lbrack \suslin \plrepNP_S(\AA^1_S/\GG) \rbrack
 \rightarrow \hat H^0\lbrack \suslin \lrepNP_S(\GG) \rbrack 
 \rightarrow & \hat H^0\lbrack \suslin \lrepNP_S(\AA^1_S) \rbrack\\ 
& \rightarrow \hat H^0\lbrack \suslin \plrepNP_S(\AA^1_S/\GG) \rbrack 
 \rightarrow 0.
\end{align*}
By definition of the relative Picard group,
 given any smooth (affine) scheme $X$, we get an exact sequence
 of abelian groups:
\begin{equation} \label{eq:GG_ftr}
0 \rightarrow \GG(X) \rightarrow
\pic(X \times_S \PP^1_S,X_0 \sqcup X_\infty)
 \rightarrow \pic(X \times_S \PP^1_S,X_0) \rightarrow 0.
\end{equation}
Thus we deduce that:
\begin{align*}
\hat H^0\lbrack \suslin \plrepNP_S(\AA^1_S/\GG) \rbrack&=0, \\
\hat H^{-1}\lbrack \suslin \plrepNP_S(\AA^1_S/\GG) \rbrack&=\GG \otimes_\ZZ \Rc.
\end{align*}
This gives in particular a canonical isomorphism:
$$
\suslin \plrepNP_S(\AA^1_S/\GG)[-1] \simeq \GG \otimes_\ZZ \Rc
$$
in $\Der(\ptr S)$.
Taking its image in $\DMex S$ we obtain a canonical isomorphism
 which can be written as:
$$
\suslin \lrepNP_S(\AA^1_S/\GG)[-1] \simeq \GG \otimes_\ZZ \Rc.
$$
Thus we can conclude because, according to Lemma \ref{suslinA1eqtriv2},
 the canonical map
$$
\lrep S {\AA^1_S/\GG} \rightarrow \suslin \lrepNP_S(\AA^1_S/\GG)
$$ 
is an isomorphism in $\DMex S$.
\end{proof}

\begin{rem}
In the course of the proof, a canonical structure of a sheaf with
 transfers over $S$ on $\GG$ has naturally appeared
 -- described by the exact sequence \eqref{eq:GG_ftr}.
 This structure is classical (see \cite[Ex. 2.4]{MVW}).
 One can describe it as follows.

Let $X$ and $Y$ be smooth $S$-schemes.
 Assume $X$ is connected (thus irreducible as it is normal).
 Let $Z$ be a closed integral subscheme $Z$ of $X \times_S Y$ which is 
 finite surjective over $X$. Then $Z/X$ corresponds to an extension
 of function fields $L/K$. The norm map of $L/K$ induces a morphism
 of abelian groups: $N_{Z/X}:\GG(Z) \rightarrow \GG(X)$.
 Then we associate with $Z$, seen as a finite correspondence from $X$ to
 $Y$, the following morphism:
$$
\GG(Y) \xrightarrow{p^*} \GG(Z) \xrightarrow{N_{Z/X}} \GG(X)
$$
where $p:Y \rightarrow Z$ is the natural projection.
\end{rem}

The following proposition is well-known
 to the expert. We include a proof for completeness.
\begin{prop} \label{prop:cohomology_G_m_regular_case}
For any regular scheme $X$ and any integer $i \geq 0$,
$$
H^i_\nis(X,\GG)=\begin{cases}
\mathcal O_X(X)^\times & \text{if } i=0, \\
\pic(X) & \text{if } i=1, \\
0 & \text{otherwise}
\end{cases}
$$
where $\pic(X)$ is the Picard group of $X$.
\index{word}{group!Picard}
\end{prop}
\begin{proof}
Let $Y$ be any \'etale scheme over $X$.
We let $C^0(V)$ be the abelian group
 made by the invertible rational functions on $V$
 and $C^1(V)$ be the group of $1$-codimensional cycles in $V$.
Classically, one associates with any rational function $f$ on $V$
 its Weil divisor
\index{word}{divisor!Weil}
  $\mathrm{div}(f) \in C^1(V)$.
 Recall, when $V$ is integral with function field $K$,
 $f \in K$, one puts:
$$
\mathrm{div}_V(f)=\sum_{x \in V^{(1)}} v_x(f).x;
$$
the sum runs over the points of codimension $1$ in $V$
 and $v_x(f)$ is the valuation of $f$ corresponding
 to the valuation ring $\mathcal O_{X,x}$.

According to this definition, we get a complex:
$$
0 \rightarrow \GG(V) \rightarrow C^0(V) \xrightarrow{\mathrm{div}_V} C^1(V).
$$
This sequence is functorial with respect
 to pullback of \'etale $X$-schemes.
 Thus we have defined a morphism of presheaves on $X_\et$:
$$
\pi:\GG \rightarrow C^*.
$$
Given any Nisnevich distinguished square $Q$
 (Example \ref{ex:lower&upper_cd_structures}),
 one can check easily that the image of $Q$ by $C^0$ (resp. $C^1$)
 is cocartesian. As a consequence $C^*$
 is a complex of Nisnevich sheaves which satisfies the Brown-Gersten property
 -- \emph{i.e.} it is Nisnevich flasque in the sense of
 Definition \ref{df:basic_complexes&topology}
 according to Proposition \ref{prop:BG_property_derived_premotivic}
 applied to the derived category of Nisnevich sheaves over $X$.

On the other hand, $\pi$ is a quasi-isomorphism of Nisnevich sheaves
 over $S$: indeed it is well-known that for any regular local ring $A$,
 the sequence
$$
0 \rightarrow A^\times \rightarrow \mathrm{Frac}(A)^\times
 \xrightarrow{\mathrm{div}_A} Z^1(A) \rightarrow 0
$$
is exact. This is an easy consequence of the fact $A$ is
 a unique factorization domain -- the classical Auslander-Buchsbaum theorem,
\index{word}{Auslander-Buchsbaum theorem}
 (e.g. \cite[20.3]{Mat}).

In particular, we get $H^i(X,\GG)=H^i(C^*(X))$ and this concludes.
\end{proof}

The following theorem is a generalization
 of a well-known computation of Voevodsky 
 for smooth schemes over a perfect field.
 The second case is a corollary of
 the two preceding propositions.
\begin{thm}
\label{thm:effective_motivic_cohomology_low_weight}
Let $S$ be a scheme and $n \in \ZZ$ an integer.
The following computation holds:
\begin{enumerate}
\item 
$$
\HHme n 0 S=
\Hom_{\DMV^{\eff}(S)}(\lrepNP_S,\lrepNP_S[n])=
\left\{ \begin{array}{ll}
\Rc^{\pi_0(S)} & \text{if } n=0 \\
0 & \text{otherwise;}
\end{array} \right.
$$
\item if $S$ is regular,
\index{word}{scheme!regular}
$$
\HHme n 1 S=
\Hom_{\DMV^{\eff}(S)}(\lrepNP_S,\lrepNP_S(1)[n])=
\left\{ \begin{array}{ll}
\mathcal O_S(S)^\times \otimes_\ZZ \Rc & \text{if } n=1 \\
\pic(S) \otimes_\ZZ \Rc & \text{if } n=2 \\
0 & \text{otherwise}
\end{array} \right.
$$
\end{enumerate}
\end{thm}
\begin{proof}
For the first case, according to Proposition \ref{prop:rep_et_ftr},
the sheaf $\lrepNP_S$ is Nisnevich local and $\AA^1$-local
as a complex of sheaves. It is obviously acyclic for the Nisnevich topology.
Thus, we conclude using again \ref{prop:rep_et_ftr} in the case $n=0$.

Consider the second case. According to Proposition
 \ref{prop:cohomology_G_m_regular_case},
 the sheaf $\GG$ on $\sm_S$ is $\AA^1$-local.
 Thus according to Proposition \ref{computing_first_tate_twist}
 $\GG \otimes \Rc[-1]$ is an $\AA^1$-resolution of $\lrepNP_S(1)$.
 In particular,
\begin{align*}
\Hom_{\DMV^{\eff}(S)}(\lrepNP_S,\lrepNP_S(1)[n])
&=\Hom_{\Der(\ftr S)}(\lrepNP_S,\GG  \otimes \Rc[n-1])\\
&=H^{n-1}_\nis(S,\GG) \otimes \Rc
\end{align*}
where the second identification follows from Remark \ref{rem:comput_Ext_ftr}.
The conclusion follows from another application of Proposition
 \ref{prop:cohomology_G_m_regular_case}.
\end{proof}

The following corollary is
 a (very) weak cancellation result in $\DMV^{\eff}(S)$~:
\begin{cor}
Let $S$ be a regular scheme. Then

$$\derR \uHom(\lrepNP_S(1),\lrepNP_S(1))=\lrepNP_S.$$

Moreover, if $m=0$ or $m=1$, for any integer $n>m$,

$$\derR \uHom(\lrepNP_S(n),\lrepNP_S(m))=0.$$
\end{cor}
\begin{proof}
We consider the first assertion.
Any smooth $S$-scheme is regular.
Hence, it is sufficient to prove that for any connected regular
scheme $S$, for any integer $n \in \ZZ$,
$$
\Hom_{\DMV^{\eff}(S)}(\lrepNP_S (1),\lrepNP_S(1)[n])=
\left\{ \begin{array}{ll}
\Rc & \text{if } n=0 \\
0 & \text{otherwise.}
\end{array} \right.
$$
Using the exact triangle 
\begin{equation} \label{eq:triangle_tate}
\lrep S {\GG}
 \rightarrow \lrep S {\AA^1}
 \rightarrow \lrepNP_S(1)[2] \xrightarrow{+1} 
\end{equation}
and the second case of the previous theorem,
we obtain the following long exact sequence
\begin{align*}
\cdots & \rightarrow \Hom(\lrep S {\AA^1},\lrepNP_S(1)[n])
 \rightarrow \Hom(\lrep S {\GG},\lrepNP_S(1)[n]) \\
 &\rightarrow \Hom(\lrepNP_{S} (1),\lrepNP_S(1)[n-1])
 \rightarrow  \Hom(\lrep S {\AA^1},\lrepNP_S(1)[n+1]) \rightarrow \cdots
\end{align*}
Then we conclude using the previous theorem and the fact 
$$\pic(\AA^1 \times S)=\pic(\GG \times S)$$
whenever $S$ is regular.

For the last assertion, we are reduced to prove 
that if $S$ is a regular scheme, for any integers $n>0$ and $i$,
$$
\Hom_{\DMV^{\eff}(S)}(\lrepNP_S(n),\lrepNP_S[i])=0.
$$
This is obviously implied by the case $n=1$.

Consider the distinguished triangle (\ref{eq:triangle_tate}) again.
Then the long exact sequence attached to the cohomological functor
$\Hom_{\DMex S}(- ,\lrepNP_S)$ and applied to this triangle 
together with the first case of the previous 
theorem allows us to conclude.
\end{proof}

\subsubsection{The motivic cohomology ring spectrum}
\label{sec:motivic_coh_spectrum}

\begin{num}\label{premotadj:DMtildeDMtr11}
According to definition \ref{df:premotivic_ftr}
 and paragraph \ref{num:graph_functor_sh}, we have an adjunction
 of abelian premotivic categories
$$
\gamma^*:\sh {}{-}\rightleftarrows\ftr -:\gamma_*
$$
such that $\gamma_*$ is conservative and exact.
According to Paragraph \ref{num:functoriality_DMu},
 it induces an adjunction of triangulated premotivic categories
\begin{equation}
\derL\gamma^*:\DMt \rightleftarrows \DM:\derR \gamma_*.
\end{equation}
Composing with the premotivic adjunction between the stable homotopy category
 and the $\AA^1$-derived homotopy category \eqref{eq:adj_SH_DMt},
  we finally get a canonical premotivic adjunction:
\begin{equation}  \label{eq:premot_adj_SH_DM}
\varphi^*:\SH \rightleftarrows \DM:\varphi_*.
\end{equation}
Recall that, because $\varphi^*$ is monoidal,
 $\varphi_*$ is weakly monoidal.
 In particular, for any scheme $S$,
 one gets canonical morphisms
$$
\un_S \rightarrow \varphi_*(\un_S)\ , \quad
 \varphi_*(\un_S) \wedge \varphi_*(\un_S) \rightarrow \varphi_*(\un_S)
$$
which gives a structure of a commutative monoid
 to the spectrum $\varphi_*(\un_S)$ \emph{i.e.} a ring spectrum.
\end{num}
\begin{df} \label{df:motivic_coh_spectrum}
Given any scheme $S$,
 one defines the \emph{motivic cohomology ring spectrum}
\index{word}{spectrum!motivic cohomology ring spectrum}
  over $S$ with coefficients in $\Rc$ as the commutative ring spectrum:
$$
\HH_{\cM,S}^\Rc:=\varphi_*(\un_S).
$$
\end{df}
The properties of the functor $\varphi_*$ immediately implies that
 the ring spectrum $\HH_{\cM,S}^\Rc$ represents motivic cohomology.
 One now easily checks that this ring spectrum coincides
 with the original one of Voevodsky (see \cite[sec. 6.1]{V3})
 in the case $\Rc=\ZZ$.
 Therefore, our definition of motivic cohomology
  (with $\ZZ$-coefficients) agrees with that given by
  Voevodsky in \emph{loc. cit.}

\begin{num}
Consider a localization $\Rc'$ of $\Rc$.
Then one gets an essentially commutative diagram of right adjoints
 of premotivic adjunctions:
$$
\xymatrix@R=8pt{
& \Der_{\AA^1}(S,\Rc)\otimes_\Rc \Rc'\ar[ld]
 & \DMV(S,\Rc)\otimes_\Rc \Rc'\ar[l] \\
\SH(S) &  & \\
& \Der_{\AA^1}(S,\Rc')\ar[lu]\ar_{(1)}[uu]
 & \DMV(S,\Rc')\ar_{(2)}[uu]\ar[l] \\
}
$$
where the map $(1)$ is the canonical equivalence
 (see Proposition \ref{prop:stable-A^1-derived_chg_coef})
 and the map (2) is the equivalence from \eqref{eq:changeofcoef_DM}.
 Note in particular that $(2)$ is monoidal (as its reciprocal equivalence
 is monoidal as the left adjoint of a premotivic adjunction).
 Thus this essentially commutative diagram defines a canonical
  morphism of ring spectra:
\begin{equation}\label{eq:changeofcoef_motivic_spectrum}
\HH_{\cM,S}^\Rc \otimes_\Rc \Rc' \rightarrow \HH_{\cM,S}^{\Rc'}.
\end{equation}
As a corollary of Proposition \ref{prop:changeofcoef_DM},
 we get the following result:
\end{num}
\begin{prop}
The map \eqref{eq:changeofcoef_motivic_spectrum} is an isomorphism
 of ring spectra.
\end{prop}

\begin{rem} \label{rem:no_warn_changeofcoef_motivic_spectrum}
In a previous version of this text, we only get the above result
 in particular cases. The main argument for the general case obtained above
 can be traced back to Proposition \ref{prop:multiplicity_denominators_bounded}.
\end{rem}

\begin{num}\label{num:conj_V}
Let $f:T \rightarrow S$ be a morphism of schemes.
Recall from the structure of the premotivic adjunction
 $(\varphi^*,\varphi_*)$ defined above that we get
 an exchange morphism:
$$
f^*\varphi_* \rightarrow \varphi_*f^*
$$
Applying this natural transformation to the unit object
 $\un_S$ of $\DMx S$, one gets a canonical morphism
 of ring spectra:
$$
\tau_f:f^*(\HH_{\cM,S}^\Rc) \rightarrow \HH_{\cM,T}^\Rc.
$$
Remark that this shows the collection $(\HH_{\cM,S}^\Rc)$
 is a section of the fibred category $\SH$.
 Recall also the following conjecture of Voevodsky
  (\cite[conj. 17]{V_OpenPB}):
\end{num}
\begin{conj}[Voevodsky]\label{conj:voevodsky}
For any morphism $f$ as above, the map $\tau_f$
 is an isomorphism.
\end{conj}
 
\begin{rem}
At least, Voevodsky formulated this conjecture in the
case where $\Rc=\ZZ$. According to the preceding proposition,
this implies the case of any coefficients ring $\Rc \subset \QQ$.
We will solve affirmatively a particular case
of this conjecture in \ref{compHbeilHVoev} when $\Lambda=\QQ$.
We will see below that this conjecture of Voevodsky is
strongly related to the behaviour of the six operations
in $\DM$; see
Proposition \ref{prop:reformulation of the conjecture of voevodsky}.
References for other known cases of variants
of the conjecture may be found in
Remark \ref{rem:voevodskysconjecture}.
\end{rem}
  
\subsection{Orientation and purity}

\begin{num}
For any scheme $S$, we let $\PP^{\infty}_S$ be the ind-scheme
$$
S \rightarrow  \PP^1_S \rightarrow \cdots
 \rightarrow \PP^n_S \rightarrow \PP^{n+1}_S \rightarrow
$$
made of the obvious closed immersions. This ind-scheme
has a comultiplication given by the Segre embedding
$$
\PP^\infty_S \times_S \PP^\infty_S \rightarrow \PP^\infty_S
$$

Define $\lrep S {\PP^\infty}=\ilim\lrep S {\PP^n}$.
Applying Theorem \ref{thm:effective_motivic_cohomology_low_weight}
in the case $S=\spec{\ZZ}$, we obtain a canonical isomorphism:
$$
\Hom_{\DMex {\spec \ZZ}}(\lrep{}{\PP^\infty},\lrepNP(1)[2])
 =\pic(\PP^\infty) \otimes_\ZZ \Rc,
$$
whose aim is a free $\Rc$-algebra of power series in one variable.
Considering the canonical dual invertible sheaf on $\PP^{\infty}$, 
 we obtain a canonical formal generator of this $\Rc$-algebra
 and thus a morphism $\DMex {\spec \ZZ}$:
$$
\mathfrak c_1:\lrep{}{\PP^\infty} \rightarrow \lrepNP(1)[2].
$$
For any scheme $S$,
 considering the canonical projection $f:S \rightarrow \spec \ZZ$,
 we obtain by pullback along $f$ a morphism of $\DMex S$
$$
\mathfrak c_{1,S}:\lrep S {\PP^\infty_S}
 \rightarrow \lrepNP_S(1)[2].
$$
Consider $\GG$ as a sheaf of groups over $\sm_S$.
Following \cite[part 4]{MV}, we introduce its classifying space
$B\GG$ as a simplicial sheaf over $\sm_S$. 
From proposition 1.16 of \emph{loc. cit.}, we get
$
\Hom_{\Hpt^s(S)}(S_+,B\GG)=\pic(S).
$
Moreover, in the homotopy category of pointed simplicial sheaves $\Hpt(S)$,
we have a canonical isomorphism
$B\GG\simeq\PP^\infty_S$ (\textit{cf.} \emph{loc. cit.}, prop. 3.7).
Thus, finally, we obtain a canonical map of pointed sets
\begin{align*}
\pic(S)&=\Hom_{\Hpt^s(S)}(S_+,B\GG)
 \rightarrow \Hom_{\Hpt(S)}(S_+,\PP^\infty) \\
 & \rightarrow \Hom_{\DMex S}(\lrepNP_S,\lrep S {\PP^\infty/*})
 \rightarrow \Hom_{\DMex S}(\lrepNP_S,\lrep S {\PP^\infty}).
\end{align*}
\end{num}
\begin{df} \label{df:orientation_DM}
Consider the above notations.
We define the first motivic Chern class
\index{word}{class!Chern}
as the following composite morphism
\begin{align*}
c_1:\pic(S)
& \longrightarrow \Hom_{\DMex S}(\lrepNP_S,\lrep S {\PP^\infty_S})
\xrightarrow{(\mathfrak c_{1,S})_*}
 \Hom_{\DMex S}(\lrepNP_S,\lrepNP_S(1)[2]) \\
& \longrightarrow  \Hom_{\DMx S}(\un_S,\un_S(1)[2])=\HHm 2 1 S
\end{align*}
\end{df}
The first motivic Chern class is evidently compatible with pullback.

\begin{rem} \label{rem:warn_motivic_Chern_additivity}
Beware that the map
$$
\pic(S) \rightarrow \Hom_{\DMex S}(\lrepNP_S,\lrep S {\PP^\infty_S})
$$
defined above is not necessarily a morphism of abelian groups.
However, the composite:
$$
\pic(S)
 \longrightarrow \Hom_{\DMex S}(\lrepNP_S,\lrep S {\PP^\infty_S})
\xrightarrow{(\mathfrak c_{1,S})_*}
 \Hom_{\DMex S}(\lrepNP_S,\lrepNP_S(1)[2])
$$
is the isomorphism of 
Theorem \ref{thm:effective_motivic_cohomology_low_weight} 
when $S$ is normal.
In particular, it is a morphism of abelian groups in this case.
We will give an argument below for the general case.
\end{rem}

\begin{num}\label{pptyChern}
The triangulated category $\DMx S$ thus satisfies all the 
axioms of \cite[\textsection 2.1]{Deg8} (see also Paragraph 2.3.1 of \emph{loc. cit.}
 in the regular case). In particular, we derive from the
main results of \emph{loc. cit.} the following facts:
\begin{enumerate}
\item Let $p:P \rightarrow S$ be a projective bundle of rank $n$.
Let $c:\un_S \rightarrow \un_S(1)[2]$ be the first Chern
 class of the canonical line bundle on $P$.
Then the map
$$
\mot S P \xrightarrow{\sum_i p \otimes c^i} \bigoplus_{i=0}^n \un_S(i)[2i]
$$
is an isomorphism. This gives the projective bundle theorem in motivic
cohomology for any base scheme.

One deduces using the method of Grothendieck that
 motivic cohomology possesses Chern classes of vector bundles
 which satisfies all the usual properties (see remark below
 for additivity).
\item Let $i:Z \rightarrow X$ be a closed immersion
 between smooth separated $S$-schemes of finite type.
 Assume $i$ has pure codimension $c$ and let $j$ be the complementary
 open immersion.
 Then there is a canonical \emph{purity isomorphism}:
\index{word}{purity!isomorphism (relative)}
$$
\mathfrak p_{X,Z}:\mot S {X/X-Z} \rightarrow \mot S Z(c)[2c].
$$
This defines in particular the \emph{Gysin triangle}
\index{word}{triangle!Gysin}
$$
\mot S {X-Z} \xrightarrow{j_*} \mot S X \xrightarrow{i^*} \mot S Z(c)[2c]
 \xrightarrow{\partial_{X,Z}} \mot S {X-Z}[1].
$$
\item Let $f:Y \rightarrow X$ be a projective morphism between
 smooth separated $S$-schemes of finite type.
 Assume $f$ has pure relative dimension $d$.
 Then there is an associated \emph{Gysin morphism}
\index{word}{morphism!Gysin}
$$
f^*:\mot S X \rightarrow \mot S Y(d)[2d]
$$
functorial in $f$. We refer the reader to \emph{loc. cit}
 for various formulas involving the Gysin morphism (projection formula,
 excess intersection,...)

Note in particular that we deduce from that Gysin morphism
 the following map in motivic cohomology:
$$
f_*:\HHm n i Y \rightarrow \HHm{n+2d} {i+d} X.
$$
\item For any smooth projective $S$-scheme $X$,
 the premotive $\mot S X$ admits a \emph{strong dual}.
\index{word}{dual, strong}
 If $X$ has pure relative dimension $d$ over $S$,
 the strong dual of $\mot S X$ is $\mot S X(-d)[-2d]$.
\end{enumerate}
\end{num}

\begin{rem}
According to \emph{loc. cit.},
 there exists for any scheme $S$ a formal group law
 $F_S(x,y)$ with coefficients in the graded ring $\HHm {2*} * S$.
 If one consider the Segre embedding
$$
\Sigma:\PP^\infty_S \rightarrow \PP^\infty_S \times_S \PP^\infty_S
$$
one has: $F_S(x,y)=\sigma^*(1)$ through the isomorphism:
$$
\HHm {2*} *{\PP^\infty_S \times_S \PP^\infty_S}
 \simeq \HHm {2*} * S[[x,y]]
$$
which results from the projective bundle formula
 in motivic cohomology.

According to Remark \ref{rem:warn_motivic_Chern_additivity},
 whenever $S$ is normal, one gets $F_S(x,y)=x+y$.
 In particular, $F_{\spec \ZZ}(x,y)=x+y$.
 On the other hand, according to the above definition of $F_S(x,y)$,
 $F_S(x,y)$ is compatible with pullback. Thus one deduces
 that $F_S(x,y)=x+y$ for any scheme $S$.
\end{rem}

\begin{num} \renewcommand{\Rc}{\ZZ}
According to the properties that we have previously proved,
 motivic cohomology, and in particular the bigraded part $\HHm {2n} n X$,
 possesses many of the desired property of a generalized Chow theory
 for regular schemes (see  \cite[XIV, \textsection 8]{SGA6}).

Note in particular that the existence of Chern classes
 allows to define a Chern character:
$$
K_0(X) \otimes_\ZZ \QQ \xrightarrow{ch} \HHm {2*} * X \otimes \QQ
 \simeq \HH_\cM^{2*,*}(X,\QQ)
$$
where the final isomorphism follows from Paragraph
 \ref{num:changeofcoef_motivic_coh}.
 In particular, we will prove in the next section (Corollary \ref{compHbeilHVoev})
 that, when $X$ is regular,
 this map is an isomorphism as expected.
\end{num}

\begin{rem} \renewcommand{\Rc}{\ZZ}
Among the good properties of motivic cohomology is the
 fact it is defined, with its ring structure and natural functoriality,
 other arbitrary schemes.
 On the other hand, even when $X$ is regular, 
 one cannot describe at the moment the cohomology group $\HHm {2n} n X$
 in terms of classes of $n$-codimensional cycles in $X$ modulo
 an appropriate equivalence relation.

Let us however mention the two following interesting facts:
\begin{enumerate}
\item Let $X$ be a scheme of finite type over $\spec \ZZ$
 and $X_p$ be its fiber over a primer $p$.
Then one has a pullback map:
$$
\HHm {2n} n X \rightarrow \HHm {2n} n {X_p}, \sigma \mapsto \sigma_p.
$$
When $X$ is an arithmetic scheme (regular and flat over $\ZZ$)
 with good reduction at $p$, the target is the Chow group
 of $n$-codimensional cycles (see Example \ref{ex:motivic_coh&Chow}).
 Then $\sigma_p$ should be thought as the specialization
\index{word}{specialization}
 of its generic fiber (which lies in $\HHm {2n} n {X_\QQ}=CH^n(X_\QQ)$
 according to the Example \ref{ex:motivic_coh&Chow}). 
 This construction should coincide with other specialization
 maps in the arithmetic case (see for example
  \cite[\textsection 20.3]{Ful}).
\item Let $X$ be a smooth $S$-scheme. Then any $n$-codimensional closed 
 subscheme $Z$ of $X$ which is smooth over $S$ defines using the Gysin morphism
 an element
$$
[Z]=i_*(1) \in \HHm {2n} n X
$$
which should be called the fundamental class of $X$.
One can extract from \cite{Deg8} some expected properties
 of these fundamental classes (relation to Chern classes, pullback
 properties such as compatibility with base change).

In particular, any $S$-point of $X$ defines an element of $\HHm {2d} d X$
 where $d$ is the dimension of $X$ (assumed of pure dimension).
 In particular, the group $\HHm {2d} d X$ is close to a group of
 cycles in $X$ of relative dimension $0$ over $S$.
\end{enumerate}
\end{rem}

\begin{num}
We end up this series of remarks on motivic cohomology
 with the following construction that the reader might enjoy.

Let $S$ be any scheme and $\mathscr P_S$ be the category
 of smooth projective $S$-schemes. Given any scheme $X$ and $Y$
 in $\mathscr P_S$, one can use the group
$$
\HHm {2d} d {X \times_S Y}
$$
where $d$ is the relative dimension of $Y$ as a group
 of correspondences using the properties obtained so far from motivic
 cohomology.
 In particular, one can mimic the construction of the category
 of Chow motives over a field $k$ using the category $\mathscr P_S$
 and these correspondences.
 One obtains an additive monoidal category $\mathrm{Chow}'(S,\Rc)$ of \emph{strong Chow motives}.
\index{word}{motive!Chow (strong)}

According to the duality property of motives
 (Paragraph \ref{pptyChern}, point 4)
 one also obtains a canonical isomorphism
$$
\Hom_{\DMx S}(\mot S X,\mot S Y)=\HHm {2d} d {X \times_S Y}.
$$
Thus one deduces a canonical full embedding
 of monoidal categories:
$$
\mathrm{Chow}'(S,\Rc) \rightarrow \DMgmx S
$$
which extends the well-known case when $S$ is a perfect field.
\end{num}

\begin{rem}
Beware that, with rational coefficients,
 a sharper notion of Chow motives -- in more precise terms,
 these are motives \emph{of weight zero} --
 have been introduced recently (see \cite{Hebert}, \cite{Bondarko}).
\end{rem}

\subsection{The six functors}

\begin{num} \label{adjunction:DMt_DM}
Recall that according to Definition \ref{df:premotivic_ftr}
 and Paragraph \ref{num:graph_functor_sh}, we have an adjunction
 of abelian premotivic categories
$$
\gamma^*:\sh {}{-}\rightleftarrows\ftr -:\gamma_*\,
$$
such that $\gamma_*$ is exact and conservative.
Moreover, for any scheme $S$, any smooth $S$-schemes $X$, $Y$
 and any open immersion $j:U \rightarrow X$,
 the canonical map:
$$
j_*:\corr S Y U \rightarrow \corr S Y X
$$
is obviously a monomorphism.
Thus the abelian premotivic category $\ftr -$ satisfies
 the assumptions (i)-(iv) of Paragraph \ref{assumption:nearly_nisnevich}.
In particular, we deduce from Corollaries \ref{cor:supp_nearly_nis}
 and \ref{smoothloc} the following theorem:
 \end{num}
\begin{prop}\label{prop:DMV_supp_wloc}
The premotivic triangulated category $\DM$ satisfies the support property.
Moreover, for any scheme $S$ and any closed immersion $i:Z \rightarrow X$
 between smooth $S$-schemes, $\DM$ satisfies the localization property
 with respect to $i$, \locx i.
\end{prop}
An important corollary of this proposition is that
 given any separated morphism $f:Y \rightarrow X$ of finite
 type, one can construct an adjunction of triangulated categories:
$$
f_!:\DMx Y \rightleftarrows \DMx X:f^!
$$
such that $f_!=f_*$ when $f$ is proper (see Section \ref{sec:Deligne}).
We will elaborate on this fact at the end of this section.

\begin{num}
Note that in particular, the premotivic category $\DM$
 satisfies the weak localization property \wloc.
According to the premotivic adjunction \eqref{eq:premot_adj_SH_DM}
 and the existence of the first Chern class in motivic cohomology
 (Definition \ref{df:orientation_DM}),
 one can apply Example \ref{ex:or&GysinII} to the premotivic triangulated
 category $\DM$ (which satisfies the Nisnevich separation property by
 construction).
 This implies in particular that $\DM$ is oriented as a premotivic triangulated
 category (Definition \ref{df:premotivic_orientation}).

In particular, one can apply Corollary \ref{cor:orientationsandpurity}
 to $\DM$ and get the following result:
\end{num}
\begin{prop} \label{prop:orientations&purity_DMV}
Any smooth projective morphism $f$ is $\DM$-pure:
 the canonical purity map \eqref{eq:relative_purity_iso_or}
$$
f_\sharp \rightarrow f_!(d)[2d],
$$
is an isomorphism where $d$ is the relative dimension of $f$.
\end{prop}

In particular, $\DM$ is weakly pure.
 The only property of the premotivic triangulated category $\DM$
 that we cannot prove is the localization
 property for general closed immersions.
However, the properties we have seen so far allows to construct
 the 6 operations and establish some of its properties
 that are already of interest. Let us summarize
 this formalism, from Theorem \ref{thm:support} together
 with Lemma \ref{lm:purity=>BC,supp,PF}:
\begin{thm}\label{thm:DM_six_functors}
For any separated morphism of finite type $f:Y \rightarrow X$,
 there exists an essentially unique pair of adjoint functors
$$
f_!:\DMx Y \rightleftarrows \DMx X:f^!
$$
such that:
\begin{enumerate}
\item There exists a structure of a covariant (resp. contravariant) 
 $2$-functor on $f \mapsto f_!$ (resp. $f \mapsto f^!$).
\item There exists a natural transformation $\alpha_f:f_! \rightarrow f_*$
 which is an isomorphism when $f$ is proper.
 Moreover, $\alpha$ is a morphism of $2$-functors.
\item For any smooth projective morphism $f:X \rightarrow S$ 
 of relative dimension $d$,
 there are canonical natural isomorphisms
\begin{align*}
\piso^\mathfrak t_f:f_\sharp & \longrightarrow f_!(d)[2d] \\
\piso^{\prime\mathfrak t}_f:f^* & \longrightarrow f^!(-d)[-2d]
\end{align*}
which are dual to each other.
\item For any cartesian square:
$$
\xymatrix@=16pt{
Y'\ar^{f'}[r]\ar_{g'}[d]\ar@{}|\Delta[rd] & X'\ar^g[d] \\
Y\ar_f[r] & X,
}
$$
such that $f$ is separated of finite type,
there exist natural transformations
\begin{align*}
g^*f_! \xrightarrow\sim f'_!{g'}^*\, , \\
g'_*{f'}^! \xrightarrow\sim  f^!g_*\, ,
\end{align*}
which are isomorphisms in the following cases:
\begin{itemize}
\item $g$ is smooth;
\item $f$ is projective and smooth.
\end{itemize}
\item For any smooth projective morphism $f:Y \rightarrow X$,
 there exist natural isomorphisms
\begin{align*}
Ex(f_!^*,\otimes):
(f_!K) \otimes_X L &\xrightarrow{\ \sim\ } f_!(K \otimes_Y f^*L)\, ,\ \\
  \uHom_X(f_!(L),K) & \xrightarrow{\ \sim\ } f_* \uHom_Y(L,f^!(K))\, ,\ \\
  f^! \uHom_X(L,M)& \xrightarrow{\ \sim\ } \uHom_Y(f^*(L),f^!(M))\, .
\end{align*}
\end{enumerate}
\end{thm}

\begin{rem}
As an example of application, let us recall the construction
 of the general trace map
\index{word}{map, trace}
  (from \cite{SGA4}) in the case
 of a smooth projective morphism $f:Y \rightarrow X$ of
 relative dimension $d$. It is the following composite map:
$$
f_*f^* \xrightarrow{\alpha_f^{-1}} f_!f^*
 \xrightarrow{\piso^{\prime\mathfrak t}_f} f_!f^!(d)[2d]
 \xrightarrow{ad'(f_!,f^!)} 1(d)[2d].
$$
This allows one to recover the Gysin map associated with $f$,
 already constructed in Paragraph \ref{pptyChern}, as well as
 the duality property for the motive $\mot X Y$.
\end{rem}

We will reformulate Voevodsky's conjecture \ref{conj:voevodsky}
in terms of the six operations
as follows.

\begin{prop}\label{prop:reformulation of the conjecture of voevodsky}
We fix a base scheme $S$ as well as a ring of coefficients $\Lambda$.
The following assertions are equivalent:
\begin{itemize}
\item[(i)] for any $S$-schemes $X$ and $Y$ and any morphism of finite type $f:X\to Y$,
the canonical map $\tau_f:f^*(\HH_{\cM,X}^\Rc) \rightarrow \HH_{\cM,Y}^\Rc$ is
invertible;
\item[(ii)] for any $S$-scheme $X$, the canonical functor
$$\ho(\Mod{\HH_{\cM,X}^\Rc})\to\DMx X$$
is an equivalence of categories, and
$\DMx -$ is a motivic category over $S$-schemes;
\item[(iii)] the premotivic category $\DMx -$ has the localization property
for $S$-schemes;
\item[(iv)] $\DMx -$ is a motivic category over $S$-schemes.
\end{itemize}
\end{prop}

\begin{proof}
The fact that properties (iii) and (iv) are equivalent
is obvious, since the only missing
property that is not known for $\DM$ to be a motivic category
is the localization property.
Condition (iv) is obviously a consequence of condition (ii).

Keeping track of notations introduced in
paragraph \ref{premotadj:DMtildeDMtr11},
we shall observe that the forgetful functor
$$\varphi_*:\DM\to\SH$$
commutes with the operator $j_\sharp$, for
any open immersion $j$, as follows. Since the forgetful
functor from $\DMt$ to $\SH$ is conservative
and commutes with $j_\sharp$ for
any open immersion $j$, it is sufficient to prove that
the functor $\derR\gamma_*:\DM\to\DMt$ has the
same property, which is precisely
Proposition \ref{comm_oubli_ouv_stable}.

Let us check that condition (i) (i.e. Voevodsky's
conjecture \ref{conj:voevodsky}) is a consequence of
condition (iv). Let us assume that (iv) holds true, and that
we have a morphism of finite type $f:X\to Y$.
The property that the canonical map
$$\tau_f:f^*(\HH_{\cM,X}^\Rc) \rightarrow \HH_{\cM,Y}^\Rc$$
is invertible is local for the Zariski topology on $X$ and on $Y$,
so that we may assume that $f$ is affine. Since the map
$\tau_f$ is invertible for $f$ smooth, we observe from there that it is
sufficient to prove that $\tau_f$ is invertible when $f$ is a closed
immersion. Let $j:U\to Y$ be the open immersion complement to $f$.
Assuming (iv), there is a homotopy cofiber sequence
of the form
$$j_\sharp\un_U\to\un_Y\to f_*\un_X$$
in $\DMx Y$, the image of which is isomorphic
to the homotopy cofiber sequence
$$j_\sharp  \HH_{\cM,U}^\Rc\to \HH_{\cM,Y}^\Rc
\to f_*\HH_{\cM,X}^\Rc$$
in $\SH(Y)$, since the functor $\varphi_*$ commutes with
$j_\sharp$ (as recalled above) and with $f_*$ (for obvious reasons).
But the localization property in $\SH$ implies that the
homotopy cofiber of the map $j_\sharp  \HH_{\cM,U}^\Rc\to \HH_{\cM,Y}^\Rc$
is $f_*f^*\HH_{\cM,Y}^\Rc$. Since the functor $f_*$ is
conservative in $\SH$ (being fully faithful), this shows that
the map $\tau_f$ is invertible.

Let us assume that condition (i) is true.
Since the forgetful functor $\varphi_*$ is conservative
and commutes with $i_*$ for any closed immersion $i$,
in order to prove that condition (iv) holds, i.e. that $\DM$ has
the localization property, it is
sufficient to prove that condition (ii) of
Corollary \ref{cor:premotivic&i_*_bis} is verified in $\DM$.
We observe furthermore that, for any smooth and
projective morphism of $S$-schemes $p:X\to Y$ every where
of relative dimension $d$,
the functor $\varphi_*$ commutes with $p_\sharp$.
Indeed, for any object $M$ in $\DMx X$, we have:
\begin{align*}
\varphi_* p_\sharp(M)
&\simeq\varphi_* p_*(M)(d)[2d]\\
&\simeq p_*\varphi_*(M)(d)[2d]\\
&\simeq p_!(\Th_X(T_f)\otimes \varphi_*(M))\\
&\simeq p_\sharp\varphi_*(M)
\end{align*}
(where the identification
$\Th_X(T_f)\otimes\varphi_*(M)\simeq\varphi_*(M)(d)[2d]$
comes from the orientation on $\varphi_*(M)$ induced by
its $\HH_{\cM,X}^\Rc$-module structure).
This implies that the functor $\varphi_*$
commutes with $f_\sharp$ for any smooth morphism
of $S$-schemes $f$. Indeed, this is a local condition with
respect to the Zariski topology both on the source and on the
target of $f$, so that it is sufficient to check the case where
$f$ is quasi-projective. Since the case where $f$ is an open
immersion is already known, and since we just discussed
the case where $f$ is a smooth and projective, this
proves our claim. Finally, we observe that,
given a closed immersion $i:Z\to X$
as well as a smooth morphism $f:Y\to X$,
the diagram
$$g_\sharp \HH_{\cM,f^{-1}(X-Z)}^\Rc\to
f_\sharp\HH_{\cM,Y}^\Rc\to i_*i^*\HH_{\cM,Y}^\Rc$$
is a homotopy cofiber sequence, where $g:f^{-1}(X-Z)\to X$
is the restriction of $f$. Since the
functor $\varphi_*$ is conservative and commutes with $f_\sharp$,
$g_\sharp$ and $i_*$, this proves that
$\DM$ has the localization property, by
Corollary \ref{cor:premotivic&i_*_bis}.

If condition (i) is true, then, by virtue of Proposition \ref{abstractmotivicmodules},
there is a morphism of premotivic categories
$$\alpha^*:\ho(\Mod{\HH_{\cM}^\Rc})\rightleftarrows\DM:\alpha_* \, .$$
Furthermore, Proposition \ref{hmtlinearproperties} implies that,
under condition (i),
$\ho(\Mod{\HH_{\cM}^\Rc})$ is a motivic category
(in particular, has the localization property).
We just saw that $\DM$ is a motivic category as well.
To prove that the functor $\alpha^*$ is an equivalence of
categories, by virtue of Corollary \ref{equivgeneratorscompactgentriang},
it is sufficient to prove that, for any smooth
morphism $f:X\to Y$, the unit map
$$f_\sharp\HH_{\cM,X}^\Rc\to
\alpha_*\alpha^*f_\sharp\HH_{\cM,X}^\Rc\simeq
\alpha_*f_\sharp\alpha^*\HH_{\cM,X}^\Rc$$
is invertible. Since the operators $\alpha_*$ and $f_\sharp$
commute (when we forget the $\HH_{\cM,X}^\Rc$-module
structure, $\alpha_*$ is just $\varphi_*$), it is sufficient
to check this property when $f$ is the identity. But
the map
$$\HH_{\cM,X}^\Rc\to
\alpha_*\alpha^*\HH_{\cM,X}^\Rc$$
is invertible (in fact the identity), by definition.
\end{proof}

\begin{rem}\label{rem:voevodskysconjecture}
A variant of Voevodsky's conjecture would be that the map
$\tau_f:f^*(\HH_{\cM,X}^\Rc) \rightarrow \HH_{\cM,Y}^\Rc$ is
invertible for \emph{regular} $S$-schemes.
We invite the reader to check that this version of
the conjecture may be reformulated as in
Proposition \ref{prop:reformulation of the conjecture of voevodsky}
(restricting ourselves to regular $S$-schemes, obviously),
essentially with the same proof. Evidence for this
weaker form of the conjecture is given by the fact that
over any field of exponent characteristic $p$, it is true
with $\Lambda=\ZZ[1/p]$; see \cite{CD5}.
A variant consists in replacing $\DM$ by its $\cdh$-local
version. In equal characteristic zero, this is proved
for possibly singular scheme in \cite{CD5} (in characteristic $p>0$,
this also holds up to $p$-torsion). The $\cdh$-local version of
$\HH_{\cM}^\Rc$ should be isomorphic to Spitzweck's
motivic cohomology spectrum \cite{motivicHZ}.
\end{rem}

\part{Beilinson motives and algebraic K-theory}
\markboth{Beilinson motives and algebraic K-theory}{}
\renewcommand{\theassumption}{12.0}
\begin{assumption}
In all this part, $\sch$ is assumed to be the category
 of noetherian schemes of finite dimension.
 \end{assumption}
\renewcommand{\theassumption}{\thesection.0}

\section{Stable homotopy theory of schemes} \label{sec:SH_recall}

\subsection{Ring spectra}

Consider a base scheme $S$.

Recall that a ring spectrum
\index{word}{spectrum!ring ----}
 $E$ over $S$ is a monoid object 
 in the monoidal category $\SH(S)$.
We say that $E$ is commutative if it is commutative as a monoid
 in the symmetric monoidal category $\SH(S)$.
In what follows,
 we will assume that all our ring spectra are commutative
  without mentioning it. \\
The premotivic category is $\ZZ^2$-graded where the first index
 refers to the simplicial sphere
\index{word}{sphere!simplicial}
  and the second one to the Tate twist.
\index{word}{Tate!twist}
 According to our general convention,
  a cohomology theory representable in $\SH$
\index{word}{cohomology!representable}
 is $\ZZ^2$-graded accordingly:
given such a ring spectrum $E$, for any smooth $S$-scheme $X$,
 and any integer $(i,n) \in \ZZ^2$, we get a bigraded ring:
$$
E^{n,i}(X)=\Hom_{\SH(S)}\big(\sus X_+,E(i)[n]\big).
$$
%%The product is graded commutative with respect to the first index.
When $X$ is a pointed smooth $S$-scheme, it defines a pointed sheaf
 of sets still denoted by $X$ and we denote by $\tilde E^{n,i}(X)$
 for the corresponding cohomology ring. \\
The \emph{coefficient ring} associated with $E$ is the cohomology
 of the base $E^{**}:=E^{**}(S)$. The ring $E^{**}(X)$ 
 (resp. $\tilde E^{**}(X)$) is in fact an $E^{**}$-algebra.

\begin{num} \label{num:strict_ring_sp}
We say $E$ is a \emph{strict ring spectrum}
\index{word}{spectrum!strict ring ----}
 if there exists a monoid in the category of symmetric Tate spectra $E'$ 
 and an isomorphism of ring spectra $E \simeq E'$ in $\SH(S)$.
In this case, a module $M$ over the monoid $E$ in the monoidal category $\SH(S)$
 will be said to be \emph{strict} if there exists an $E'$-module $M'$
 in the category of symmetric Tate spectra, as well as an isomorphism
 of $E$-modules $M\simeq M'$ in $\SH(S)$.
\end{num}
 
\subsection{Orientation}

\begin{num} Consider the infinite tower
$$
\PP^1_S \rightarrow \PP^2_S \rightarrow \cdots \rightarrow \PP^n_S \rightarrow  \cdots
$$
of schemes pointed by the infinity.
We denote by $\PP^\infty_S$ the limit of this tower
 as a pointed Nisnevich sheaf of sets
 and by $\iota:\PP^1_S \rightarrow \PP^\infty_S$ the induced inclusion.
Recall the following definition, classically translated from topology:
\end{num}

\begin{df} \label{df:orientation}
Let $E$ be a ring spectrum over $S$.
An \emph{orientation}
\index{word}{orientation!of a ring spectrum}
 of $E$ is a cohomology class $c$ in $\tilde E^{2,1}(\PP^\infty_S)$
 such that $\iota^*(c)$ is sent to the unit of the coefficient ring of $E$
 by the canonical isomorphism $\tilde E^{2,1}(\PP^1_S)=E^{0,0}$.

We then say that $(E,c)$ is an \emph{oriented ring spectrum}.
We shall say also that $E$ is \emph{orientable} if there exists an orientation $c$.
\end{df}
According to \cite[1.16 and 3.7]{MV},
 we get a canonical map for any smooth $S$-scheme $X$
$$
\pic(X)=H^1(X,\GG)
 \rightarrow \Hom_{\H_\bullet(S)}(X_+,\PP^\infty)
 \rightarrow \Hom_{\SH(S)}(\sus X_+,\sus \PP^\infty)
$$
(the first map is an isomorphism whenever $S$ is regular (or even geometrically unibranch)).
Given this map, an orientation $c$ of a ring spectrum $E$ defines a map of sets
$$
c_{1,X}:\pic(X) \rightarrow E^{2,1}(X)
$$
which is natural in $X$ (and from its construction in \cite{MV},
 one can check that $c=c_{1,\PP^\infty_S}(\mathcal O(1))$).
Usually, we put $c_1=c_{1,X}$.

\begin{ex} \label{ex:oriented_spectra}
\begin{enumerate}
\item The original example of an oriented ring spectrum is
 the algebraic cobordism spectrum
\index{word}{spectrum!algebraic cobordism}
 $\MGL$
\index{notat}{MGL@$\MGL$}
 introduced by Voevodsky (\textit{cf.} \cite{V3}).
\item According to Definition \ref{df:orientation_DM},
 the motivic cohomology ring spectrum $\HH_{\cM,S}^\Rc$
\index{word}{spectrum!motivic cohomology ring spectrum}
 defined in \ref{df:motivic_coh_spectrum}
 is an oriented ring spectrum.
\item Consider a triangulated premotivic category $\T$ which satisfies
 the weak localization property \wloc and such that
 there exists an adjunction of triangulated premotivic categories:
$$
\varphi^*:\SH \rightleftarrows \T:\varphi_*.
$$
Recall that $\varphi^*$ is symmetric monoidal. 
Thus, its right adjoint is
 weakly symmetric monoidal and for any the spectrum 
$$
\HH_{\T,S}:=\varphi_*(\un_S)
$$
 admits a (commutative) ring structure.

Then $\T$ is oriented in the sense of Definition \ref{df:premotivic_orientation}
 if and only if the ring spectrum $\HH_{\T,S}$ is oriented in the
 sense of Definition \ref{df:orientation} -- see Example \ref{ex:or&GysinII}.
 Moreover, an orientation of $\T$ is equivalent to the data of
 orientations $\HH_{\T,S}$ for any scheme $S$ which are stable
 by pullbacks (on cohomology).
\end{enumerate}
\end{ex}

\begin{rem} \label{rem:strict_ring_sp&GysinII}
When $E$ is a strict ring spectrum, the category $\Mod E$ satisfies
 the axioms of \cite[2.1]{Deg8} (see example 2.12 of \emph{loc.cit.}).
\end{rem}

Recall the following result, which first appeared in \cite{Vezzosi}:
\begin{prop}[Morel]
Let $(E,c)$ be an oriented ring spectrum. Then:
$$
E^{**}(\PP^{\infty}_S)=E^{**}[[c]] \quad\text{and}\quad
E^{**}(\PP^\infty_S \times \PP^\infty_S)=E^{**}[[x,y]]\,  ,
$$
where $x$ (resp. $y$) is the pullback of $c$
 along the first (resp. second) projection.
\end{prop}

\begin{rem} When $E$ is a strict ring spectrum, this is \cite[3.2]{Deg8}
 according to remark \ref{rem:strict_ring_sp&GysinII}.
The proof follows an argument of Morel (\cite[lemma 3.3]{Deg8})
 and the arguments of \emph{op.cit.}, p. 634, can be easily used
 to obtain the proposition arguing directly for the cohomology
 functor $X \mapsto E^{*,*}(X)$.
\end{rem}

\begin{num} Recall that the Segre embedding
\index{word}{embedding, Segre}
$$
\PP^n_S \times \PP^m_S \rightarrow \PP^{n+m+nm}_S
$$
define a map
$$
\sigma:\PP^\infty_S \times \PP^\infty_S \rightarrow \PP^\infty_S.
$$
It gives the structure of an $H$-group
\index{word}{group!Hgroup@$H$-group}
 to $\PP^\infty_S$ in the homotopy category $\H(S)$.
Consider the hypothesis of the previous proposition. 
Then the pullback along $\sigma$ in $E$-cohomology induces a map
$$
E^{**}[[c]] \xrightarrow{\sigma^*} E^{**}[[x,y]]
$$
and following Quillen, we check that the formal power series
$\sigma^*(c)$ defines a formal group law
\index{word}{law, formal group}
 over the ring $E^{**}$.
\end{num}

\begin{df}
Let $(E,c)$ be an oriented ring spectrum and consider the previous notation.

The formal group law $F_E(x,y):=\sigma^*(c)$ will be called
 the formal group law associated to $(E,c)$.
\end{df}
Recall the formal group law has the form:
$$
F_E(x,y)=x+y+\sum_{i+j>0} a_{ij}.x^iy^j
$$
with $a_{ij}=a_{ji}$ in $E^{-2i-2j,-i-j}$. \\
The coefficients $a_{ij}$ describe the failure of additivity
 of the first Chern class $c_1$. 
Indeed,
 assuming the previous definition, 
 we get the following result:
\begin{prop}
Let $X$ be a smooth $S$-scheme.
\begin{enumerate}
\item For any line bundle $L/X$,
 the class $c_1(L)$ is nilpotent
\index{word}{nilpotent}
  in $E^{**}(X)$.
\item Suppose $X$ admits an ample line bundle.
 For any line bundles $L,L'$ over $X$,
$$
c_1(L_1 \otimes L_2)=F_E(c_1(L_1),c_1(L_2)) \in E^{2,1}(X).
$$
\end{enumerate}
\end{prop}
We refer to \cite[3.8]{Deg8} in the case where $E$ is strict;
 the proof is the same in the general case.

Recall the following theorem of Vezzosi (\textit{cf.} \cite[4.3]{Vezzosi}):
\begin{thm}[Vezzosi] Let $(E,c)$ be an oriented
spectra over $S$, with formal group law $F_E$.
Then there exists a bijection between the following sets:
\begin{enumerate}
\item[(i)] Orientation classes $c'$ of $E$.
\item[(ii)] Morphisms of ring spectra $\MGL \rightarrow E$ in $\SH(S)$.
\item[(iii)] Couples $(F,\varphi)$ where $F$ is a formal group law over $E^{**}$
 and $\varphi$ is a power series over $E^{**}$ which defines an isomorphism
 of formal group law: $\varphi$ is invertible as a power series and 
 $F_E(\varphi(x),\varphi(y))=F(x,y)$.
\end{enumerate}\label{equivorientations}
\end{thm}

\subsection{Rational category}

{\renewcommand{\Rc}{\QQ}
In what follows, we shall use frequently
the equivalence of premotivic categories (see \ref{par:link_stablehomotopy_homology})
$$
\SH_\QQ \rightleftarrows \DMt \, ,
$$
and shall identify freely any rational spectrum
\index{word}{spectrum!rational}
 over a scheme $S$
 with an object of $\DMtx S$.}

\section{Algebraic K-theory} \label{sec:Ktheory}
\subsection{The K-theory spectrum}

We consider the spectrum $\BGL_S$
\index{notat}{KGLS@$\BGL_S$}
 which represents homotopy invariant K-theory
\index{word}{cohomology!K-theory|see{K-theory}}
\index{word}{K-theory!homotopy invariant}
  in $\SH(S)$ according to Voevodsky
   (see \cite{cis4}, \cite[6.2]{V3}, \cite[5.2]{RiouBGL}
    and \cite{PPR}). 
It is characterized by the following properties:
\begin{enumerate}
\item[(K1)] For any morphism of schemes $f:T \rightarrow S$,
there is an isomorphism $f^*\BGL_S\simeq \BGL_T$ in $\SH(T)$.
\item[(K2)] For any regular scheme $S$ and any integer $n$, 
 there exists an isomorphism
$$
\Hom_{\SH(S)}\left(\un_S[n],\BGL_S\right) \rightarrow K_n(S)
$$
(where the right-hand side is Quillen algebraic K-theory as defined
 by Thomason and Trobaugh, \cite{TT}, in the case where $S$ does not admit
 an ample family)
\index{word}{K-theory!Quillen}
such that, for any morphism $f:T \rightarrow S$ of regular schemes,
 the following diagram is commutative:
$$
\xymatrix@=16pt{
\Hom\left(\un_S[n],\BGL_S\right)\ar[r]\ar[d]
 & \Hom\left(f^*\un_S[n],f^*\BGL_S\right)\ar@{=}[r]
 & \Hom\left(\un_T[n],\BGL_T\right)\ar[d] \\
K_n(S)\ar^{f^*}[rr] && K_n(T)
}
$$
(where the lower horizontal map is the pullback in Quillen algebraic K-theory
 along the morphism $f$ and the upper horizontal
 map is obtained by using the functor $f^*:\SH(S) \rightarrow \SH(T)$
 and the identification (K1)).
%%\end{enumerate}
%%Moreover, according to \cite{PPR} and \cite{GS},
%% it is a strict ringed spectrum:
%%\begin{enumerate}
\item[(K3)] For any scheme $S$,
 there exists a unique structure of a commutative monoid on $\BGL_S$
 which is compatible with base change -- using the identification (K1) --
 and induces the canonical ring structure on $K_0(S)$.
%% Moreover, this structure is strict.
\end{enumerate}
Thus, according to (K1) and (K3), the collection of the ring spectrum $\BGL_S$
 for any scheme $S$ form an absolute ring spectrum. As usual, when no confusion can
 arise, we will not indicate the base in the notation $\BGL$. \\
Note that (K1) implies formally that the isomorphism of (K2) can be
 extended for any smooth $S$-scheme $X$ (with $S$ regular), giving a natural isomorphism:
$$
\Hom_{\SH(S)}\left(\sus X_+[n],\BGL\right) \rightarrow K_n(X)\, .
$$

\subsection{Periodicity}

\begin{num} Recall from the construction the following property of the spectrum
 $\BGL$:
\begin{enumerate}
\item[(K4)] the spectrum $\BGL$ is a $\PP^1$-periodic spectrum in
the sense that there exists a canonical isomorphism
$$
\BGL \xrightarrow \sim
 \derR \sHom\left(\sus \PP^1_S,\BGL\right)=\BGL(-1)[-2].
$$
As usual, $\PP^1_S$ is pointed by the infinite point.
\end{enumerate}
This isomorphism, classically called the Bott isomorphism,
\index{word}{Bott isomorphism}
 is characterized uniquely by the fact that its adjoint
isomorphism (obtained by tensoring with $\un_S(1)[2]$) is equal to the composite
\begin{equation} \label{eq:inv_Bott_iso}
\gamma_u:\BGL(1)[2]
 \xrightarrow{1 \otimes u} \BGL \wedge \BGL
 \xrightarrow{\mu} \BGL.
\end{equation}
where $u:\sus \PP^1 \rightarrow \BGL$ corresponds to
the class $[\mathcal O(1)]-1$ in $\tilde K_0(\PP^1)$ 
through the isomorphism (K2) and $\mu$ is the structural map
of monoid from (K3).

Using the isomorphism of (K4), the property (K1) can be extended
as follows: for any smooth $S$-scheme $X$ and any integers $(i,n) \in \ZZ^2$,
there is a canonical isomorphism:
\begin{equation} \label{eq:rep_K-theory}
\BGL^{n,i}(X) \xrightarrow{\ \sim\ } K_{2i-n}(X).
\end{equation}
\end{num}

\begin{rem} \label{orientation_KGL}
The element $u$ is invertible in the ring $\BGL^{*,*}(S)$.
Its inverse is the \emph{Bott element} $\beta \in \BGL^{2,1}(S)$.
If we chose as an orientation of the ring spectrum $\BGL$ (\textit{cf.} \ref{df:orientation})
 the class
$$
\beta.([\mathcal O(1)]-1) \in \BGL^{2,1}(\PP^\infty),
$$
the corresponding formal group law is
 the multiplicative formal group law:
$$
F(x,y)=x+y+\beta^{-1}.xy.
$$
\end{rem}

\subsection{Modules over algebraic K-theory}

\begin{thm}[{\O}stv{\ae}r, R\"ondigs, Spitzweck] The spectrum
$\BGL$ can be represented canonically by
a cartesian monoid $\BGL'$,
\index{notat}{KGLprime@$\BGL'$}
 as well as by a homotopy cartesian commutative monoid $\BGL^\beta$
\index{notat}{KGLbeta@$\BGL^\beta$}
  in the fibred model category of symmetric $\PP^1$-spectra,
in such a way that there exists a morphism of monoids
$\BGL'\To\BGL^\beta$ which is a termwise stable
$\AA^1$-equivalence.\label{BGLgoodformodules}
\end{thm}

\begin{proof}
For any noetherian scheme of finite dimension $S$,
one has a strict commutative ring spectrum $\BGL^\beta_S$
which is canonically isomorphic to $\BGL_S$ in $\SH(S)$
as ring spectra; see \cite{ROS}. One can check that the objects
$\BGL^\beta_S$ do form a commutative monoid over the
diagram of all noetherian schemes of finite dimension
(i.e. a commutative monoid in the category of sections of the fibred
category of $\PP^1$-spectra over the category of noetherian schemes
of finite dimension), either by hand, by following the
explicit construction of \emph{loc. cit.}, either by
modifying its construction very slightly as follows:
one can perform \emph{mutatis mutandis} the construction of \emph{loc. cit.} in the
$\PP^1$-stabilization of the $\AA^1$-localization
of the model category of Nisnevich simplicial sheaves over
(any essentially small adequate subcategory of) the category
of all noetherian schemes of finite dimension, and get an object $\BGL^\beta$,
whose restriction to each of the categories $\sm/S$ is the
object $\BGL^\beta_S$. From this point of view, we clearly
have canonical maps $f^*(\BGL^\beta_S)\to\BGL^\beta_T$
for any morphism of schemes $f:T\To S$. The object $\BGL^\beta$
is homotopy cartesian, as the composed map
$$\derL f^*(\BGL_S)\simeq\derL f^*(\BGL^\beta_S)\To  f^*(\BGL^\beta_S) \to\BGL^\beta_T\simeq \BGL_T$$
is an isomorphism in $\SH(T)$. Consider now a cofibrant resolution
$$\BGL'_{\spec \ZZ}\To \BGL^\beta_{\spec \ZZ}$$
in the model category of monoids of the category of symmetric $\PP^1$-spectra over $\spec \ZZ$;
see Theorem \ref{cmfmonoids}. Then, we define, for each noetherian
scheme of finite dimension $S$, the $\PP^1$-spectrum $\BGL'_S$
as the pullback of $\BGL'_{\spec \ZZ}$ along the map $f:S\To\spec \ZZ$.
As the functor $f^*$ is a left Quillen functor, the object $\BGL'_S$ is
cofibrant (both as a monoid and as a $\PP^1$-spectrum), so that we get,
by construction, a termwise cofibrant cartesian strict $\PP^1$-ring spectrum $\BGL'$,
as well as a morphism $\BGL'\To\BGL^\beta$ which is a termwise stable $\AA^1$-equivalence.
\end{proof}

\begin{paragr}
For each noetherian scheme of finite dimension $S$, one can consider
the model categories of modules over $\BGL'_S$ and $\BGL^\beta_S$
respectively; see \ref{abstractcmfmodules}.
The change of scalars functor along the stable $\AA^1$-equivalence
$\BGL'_S\To \BGL^\beta_S$
defines a left Quillen equivalence, whence an equivalence of
homotopy categories:
$$\ho(\Mod{\BGL'_S})\simeq \ho(\Mod{\BGL^\beta_S}).$$
\end{paragr}

\begin{df}\label{df:KGL-mod}
We define the premotivic triangulated category of $\BGL$-modules over $\base$
 \index{word}{modules!KGLmodules@$\BGL$-modules}
 as the fibred triangulated category whose fiber over a scheme $S$ in $\base$
 is defined as:
$$
\Mod\BGL(S):=\ho(\Mod{\BGL^\beta_S}).
$$
\end{df}

\begin{paragr}
By definition, for any smooth $S$-scheme $X$, we have
a canonical isomorphism
$$\Hom_{\SH(S)}(\sus(X_+),\BGL[n])
\simeq \Hom_\BGL( \BGL_S(X),\BGL[n])$$
(where $\BGL_S(X)=\BGL_S\wedge^\derL_S \sus (X_+)$, while
$\Hom_\BGL$ stands for $\Hom_{\Mod{\BGL}(S)}$).

According to (K1) and (K3), for any regular scheme $X$,
 we thus get a canonical isomorphism:
\begin{equation} \label{eq:K-theory&module}
\epsilon_S:\Hom_{\BGL}(\BGL_S[n],\BGL_S)
 \xrightarrow{\ \sim\ } K_n(S).
\end{equation}
Using Bott periodicity (K4), and the compatibility with
base change, this isomorphism can be extended for any
smooth $S$-scheme $X$ and any pair $(n,m) \in \ZZ^2$:
\begin{equation} \label{eq:K-theory&module_general}
\epsilon_{X/S}:\Hom_{\BGL}(\BGL_S(X),\BGL_S(m)[n])
 \xrightarrow{\ \sim\ } K_{2m-n}(X).
\end{equation}
\end{paragr}

\begin{cor} \label{cor:modules_Kth_algebrique}
The premotivic triangulated category $\Mod\BGL)$ form a motivic category,
 and the functors
$$\SH(S) \To\Mod\BGL(S)\ , \quad M\mapsto \BGL_S \wedge^\derL_S M$$
for a scheme $S$ in $\base$
define a morphism of motivic categories
$$\SH \To \Mod\BGL$$
over the category of noetherian schemes of finite dimension.
\end{cor}

\begin{proof}
This follows from the preceding theorem
and from \ref{abstractmotivicmodules}
and \ref{hmtlinearproperties}.
\end{proof}

\subsection{K-theory with support}

\begin{num} Consider a closed immersion $i:Z \rightarrow S$
 with complementary open immersion $j:U \rightarrow S$. 
 Assume $S$ is regular. \\
We use the definition of \cite[2.13]{Gillet} for the
K-theory of $S$ with support in $Z$
\index{word}{K-theory!with support}
 denoted by $K_*^Z(S)$. In other words, we define
$K^Z(S)$ as the homotopy fiber of the restriction map
$$\derR\Gamma(S,\BGL_S)=K(S)\To K(U)=\derR\Gamma(U,\BGL_U)\, ,$$
and put: $K_n^Z(S)=\pi_n(K^Z(S))$.

Applying the derived global section functor $\derR\Gamma(S,-)$
to the homotopy fiber sequence
\begin{equation} \label{eq:DMtr_localization00}
i_!\, i^!\BGL_S\To\BGL_S\To j_*\, j^*\BGL_S\, ,
\end{equation}
we get a homotopy fiber sequence
\begin{equation} \label{eq:DMtr_localization01}
\derR\Gamma(S,i_!\, i^!\BGL_S)\To
\derR\Gamma(S,\BGL_S)\To\derR\Gamma(U,\BGL_S)
\end{equation}
from which we deduce an isomorphism in the stable homotopy
category of $S^1$-spectra:
\begin{equation} \label{eq:DMtr_localization02}
\derR\Gamma(Z,i^!\BGL_S)=\derR\Gamma(S,i_!\, i^!\BGL_S)\simeq K^Z(S)\, .
\end{equation}
%% Recall that the localization triangle in $\SH(S)$
%%  associated with $i$:
%% \begin{equation} \label{eq:DMtr_localization}
%% j_!j^!(\un_S)
%%  \To \un_S
%% \To i_*i^*(\un_S)
%% \To j_!j^!(\un_S)[1]
%% \end{equation}
%% is canonically isomorphic to the infinite suspension of
%% the homotopy fiber sequence in $\H_\bullet(S)$:
%% $$
%% U_+\To S_+ \rightarrow S/U.
%% $$
%% Using \cite[\textsection 3.1]{RiouBGL} and this property,
We thus get the following property:
\begin{enumerate}
\item[(K6)] There is a canonical isomorphism
$$
\Hom_{\SH(S)}\left(\un_S[n],i_!i^!\BGL_S\right) \rightarrow K_n^Z(S)
$$
which satisfies the following compatiblities:\\
\noindent (K6a) the following diagram is commutative:
$$
\xymatrix@C=12pt@R=18pt{
\Hom\!\left(\un[n+1],j_*j^*\BGL_S\right)\ar[r]\ar[d]
 & \Hom\!\left(\un[n],i_!i^!\BGL_S\right)\ar[r]\ar[d]
 & \Hom\!\left(\un[n],\BGL_S\right)\ar[d] \\
K_{n+1}(U)\ar[r]
 & K_n^Z(S)\ar[r]
 & K_n(S)
}
$$
where the upper horizontal arrows are induced by the localization
 sequence \eqref{eq:DMtr_localization00},
and the lower one is the canonical sequence of K-theory with support.
The extreme left and right vertical maps are the isomorphisms of (K2);\\
\noindent (K6b) for any morphism $f:Y \rightarrow S$ of regular schemes,
 $k:T \rightarrow Y$ the pullback of $i$ along $f$,
 the following diagram is commutative:
$$
\xymatrix@C=6.5pt@R=18pt{
\Hom(\un[n],i_!i^!\BGL_S)\ar[r]\ar[d]
 & \Hom(f^*\un[n],f^*i_!i^!\BGL_S)\ar[r]
 & \Hom(\un[n],k_!k^!\BGL_Y)\ar[d] \\
K_n^Z(S)\ar^{f^*}[rr] 
 && K_n^T(Y)
}
$$
where the lower horizontal map is given by the functoriality
of relative K-theory (induced by the funtoriality of K-theory)
and the upper one is obtained using the functor $f^*$ of $\SH$, 
 the canonical exchange morphism $f^*i_!i^!\To k_!k^!f^*$
 and the identification (K1).
\end{enumerate}
This property can be extended to the motivic category $\ho(\Mod{\BGL})$
 and we get a canonical isomorphism
\begin{equation} \label{eq:K-theory_support}
\epsilon_i:\Hom_{\BGL}(\BGL_S[n],i_!i^!\BGL_S)
 \xrightarrow{\ \sim\ } K_n^Z(S)
\end{equation}
satisfying the analog of (K6a) and (K6b).
\end{num}

\subsection{Fundamental class}

\begin{num} \label{num:df_Tor_indep}
Consider a cartesian square of regular schemes
$$
\xymatrix{
Z'\ar^k[r]\ar_g[d] & S'\ar^f[d] \\
Z\ar^i[r] & S
}
$$
with $i$ a closed immersion. We will say that this square
is \emph{Tor-independant}
\index{word}{square!Tor-independant}
 if $Z$ and $S'$ are Tor-independent over $S$
in the sense of \cite[III, 1.5]{SGA6}: for any $i>0$,
 $\mathrm{Tor}_i^S(\mathcal O_Z,\mathcal O_{S'})=0$.\footnote{For example,
 when $i$ is a regular closed immersion of codimension $1$,
 this happens if and only if the above square is transversal.}

In this case, when we assume in addition that all the schemes
 in the previous square are regular and that $i$ is a closed immersion
 we get from \cite[3.18]{TT}\footnote{When all the schemes in the square
 admit ample line bundles, we can refer to \cite[2.11]{Quillen}.} the formula
$$
f^*i_*=k_*g^*:K_*(Z) \rightarrow K_*(S')
$$
in Quillen K-theory. An important point for us is that there exists
a \emph{canonical homotopy} between these morphisms at the level of the
Waldhausen spectra.\footnote{In the proof of Quillen, 
one can also trace back a canonical homotopy with the restriction
mentioned in the preceding footnote.}
According to the localization theorem of Quillen \cite[3.1]{Quillen},
 we get:
 \end{num}
 
\begin{thm}[Quillen] For any closed immersion
$i:Z \rightarrow S$ between regular schemes,
 there exists a canonical isomorphism
$$
\mathfrak q_i:K_n^Z(S) \rightarrow K_n(Z).
$$
Moreover, this isomorphism is functorial with respect to 
the Tor-independent squares as above, with $i$ a closed immersion
and all the schemes regular.\label{thm:loc_Quillen}
\end{thm}

\begin{rem} \label{rem:i_*&forget-support}
In the condition of this theorem,
 the following diagram is commutative by construction:
$$
\xymatrix@R=-2pt{
K_n^Z(S)\ar[rd]\ar_{\mathfrak q_i}[dd] &  \\
& K_n(S) \\
K_n(Z)\ar_{i_*}[ru] &
}
$$
where the non labeled map is the canonical one.
\end{rem}

\begin{df} \label{df:fund_class_BGL}
Let $i:Z \rightarrow S$ be a closed immersion between regular schemes.

We define the \emph{fundamental class}
\index{word}{class!fundamental}
 associated with $i$ as the morphism of $\BGL$-modules:
$$
\eta_i:i_*\BGL_Z \rightarrow \BGL_S
$$
defined by the image of the unit element $1$ through the following morphism:
$$
K_0(Z) \xrightarrow{\mathfrak q_i^{-1}} K_0^Z(S)
 \xrightarrow{\epsilon_i^{-1}}
 \Hom(\BGL_S,i_!i^!\BGL_S)
 =\Hom(i_*\BGL_Z,\BGL_S).
$$
We also denote by $\eta'_i:\BGL_Z \rightarrow i^!\BGL_S$ the morphism
 obtained by adjunction.
\end{df}

\begin{rem}\label{rem:fundamental_class&pullback}
The fundamental class has the following functoriality properties.
\begin{itemize}
\item[(1)] By definition, and applying remark \ref{rem:i_*&forget-support},
the composite map
$$
\BGL_S \To i_*i^*(\BGL_S)=i_*\BGL_Z \xrightarrow{\eta_i} \BGL_S
$$
corresponds via the isomorphism $\epsilon_S$ to $i_*(1) \in K_0(S)$.
According to \cite[Exp.~VII, 2.7]{SGA6}, this class is equal to $\lambda_{-1}(N_i)$
 where $N_i$ is the conormal sheaf of the regular immersion $i$.
\item[(2)] In the situation of a Tor-independent square as in \ref{num:df_Tor_indep},
remark that $f^*\eta_i=\eta_k$ through the canonical exchange isomorphism
$f^*i_*=k_*g^*$ --- apply the functoriality of $\epsilon_i$ from (K6b)
 and the one of $\mathfrak q_i$.
\item[(3)] Using the identification $i^!i_*=1$, we get $\eta'_i=i^!\eta_i$.
Consider a cartesian square as in \ref{num:df_Tor_indep} and assume $f$
 is smooth. Then the square is Tor-independent and we get 
 $g^*\eta'_i=\eta'_k$ using the exchange isomorphism $g^*i^!=k^!f^*$.
\end{itemize}
\end{rem}

\subsection{Absolute purity for K-theory}

\begin{prop} \label{prop:key_abs_purity}
For any closed immersion $i:Z \rightarrow S$ between regular schemes,
 the following diagram is commutative:
$$
\xymatrix{
\Hom_\BGL(\BGL_Z[n],\BGL_Z)\ar^{\eta'_i}[r]\ar_{\epsilon_Z}[d]\ar@{}|{(*)}[rd]
 & \Hom_\BGL(\BGL_Z[n],i^!\BGL_S)\ar^{\epsilon_i}[d] \\
K_n(Z)\ar^{\mathfrak q_i^{-1}}[r] & K_n^Z(S)
}
$$
\end{prop}
\begin{proof}
In this proof, we denote by $[-,-]$ the bifunctor $\Hom_{\BGL}(-,-)$. \\
\textit{Step}~1: We assume that $i:Z \rightarrow S$ admits a retraction $p:S \rightarrow Z$. \\
Consider a $\BGL$-linear map $\alpha:\BGL_Z[n] \rightarrow \BGL_Z$. Then, 
$\eta'_i(\alpha)$ corresponds by adjunction to the composition
$$
i_*\BGL_Z[n] \xrightarrow{i_*(\alpha)} i_*\BGL_Z \xrightarrow{\eta_i} \BGL_S.
$$
Applying the projection formula for the motivic category $\ho(\Mod{\BGL})$,
 we get:
$$
i_*(\alpha)=i_*(1 \otimes i^*p^*(\alpha))=i_*(1) \otimes p^*(\alpha).
$$
Here $1$ stands for the identity morphism of the $\BGL$-module $\BGL_Z$.
This shows that $\eta'_i(\alpha)$ corresponds by adjunction
 to the composite map:
$$
\eta_i \otimes p^*(\alpha):i_*\BGL_Z[n]=i_*\BGL_Z[n]\otimes\BGL_S\To\BGL_S\otimes\BGL_S=\BGL_S
$$
(the tensor product is the $\BGL$-linear one).
By assumption, $i_*:K_*(Z) \rightarrow K_*(S)$ admits a retraction which implies
the canonical map $\mathcal O_i:K_*^Z(S) \rightarrow K_*(S)$ admits a retraction
(\textit{cf.} remark \ref{rem:i_*&forget-support}). 
To check that the diagram $(*)$ is commutative,
 we can thus compose with $\mathcal O_i$. \\
Recall the first point of remark \ref{rem:fundamental_class&pullback}:
applying property (K6a) and the fact the isomorphism
 $\epsilon_S:[\BGL_S[n],\BGL_S] \rightarrow K_n(S)$ is compatible with
 the algebra structures, we are finally reduced to prove that
$$
i_*(\alpha)=i_*(1).p^*(\alpha) \in K_{n}(S).
$$
This follows from the projection formula in K-theory
 (see \cite[2.10]{Quillen} and \cite[3.17]{TT}). \\
\textit{Step}~2: We shall reduce the general case to Step~1.
We consider the following deformation 
to the normal cone diagram: let $D$ be the blow-up of $\AA^1_S$ in the closed
subscheme $\{0\} \times Z$, $P$ be the projective completion of the
normal bundle of $Z$ in $S$
 and $s$ be the canonical section of $P/Z$;
we get the following diagram of regular schemes:
\begin{equation} \label{eq:deformation}\begin{split}
\xymatrix@=25pt{
Z\ar^-{s_1}[r]\ar_i[d] & \AA^1_Z\ar[d] & Z\ar_-{s_0}[l]\ar^s[d] \\
S\ar[r] & D & P\ar[l]
}\end{split}
\end{equation}
where $s_0$ (resp. $s_1$) is the zero (resp. unit) section of $\AA^1_Z$ over $Z$.
These squares are cartesian and Tor-independent in the sense of \ref{num:df_Tor_indep}.
The maps $s_0$ and $s_1$ induce isomorphisms in K-theory because $Z$
is regular. Thus, the second point of remark \ref{rem:fundamental_class&pullback}
allows reducing to the case of the immersion $s$
which was done in Step 1.
\end{proof}

\begin{num} Consider a cartesian square
$$
\xymatrix{
T\ar^k[r]\ar_g[d] & X\ar^f[d] \\
Z\ar^i[r] & S
}
$$
such that $S$ and $Z$ are regular, $i$ is a closed immersion and $f$ is smooth.
In this case, the following diagram is commutative
$$
\xymatrix@R=12pt@C=30pt{
\Hom_\BGL(\BGL_Z(T)[n],\BGL_Z)\ar^-{\eta'_i}[r]\ar@{=}[d]
 & \Hom_\BGL(\BGL_Z(T)[n],i^!\BGL_S)\ar@{=}[d] \\
\Hom_\BGL(\BGL_T[n],\BGL_T)\ar^-{\eta'_k}[r]
 & \Hom_\BGL(\BGL_T[n],k^!\BGL_X)
}
$$
using the adjunction $(g_\sharp,g_*)$, the exchange isomorphism $g^*i^!\simeq k^!f^*$
(which uses relative purity for smooth morphisms)
and the third point of remark \ref{rem:fundamental_class&pullback}.
In particular, the preceding proposition has the following consequences:
\end{num}

\begin{thm}[Absolute purity]\index{word}{purity!absolute} For any
closed immersion $i:Z \rightarrow S$ between regular schemes,
the map 
$$
\eta'_i:\BGL_Z \rightarrow i^!\BGL_S
$$
is an isomorphism in the category
$\ho(\Mod{\BGL})(Z)$ (or in $\SH(Z)$).\label{abspurityBGL}
\end{thm}

\begin{cor} \label{cor:fund_class&pushout}
Given a cartesian square as above, for any pair $(n,m) \in \ZZ^2$,
 the following diagram is commutative:
$$
\xymatrix@R=18pt@C=26pt{
\Hom(\BGL_S(X),i_*\BGL_Z(m)[n])\ar^-{\eta_{i}}[r]\ar@{=}[d]
 & \Hom(\BGL_S(X),\BGL_S(m)[n])\ar^-{\epsilon_{X/S}}_-\sim[dd] \\
\Hom(\BGL_Z(T),\BGL_Z(m)[n])\ar_-{\epsilon_{T/Z}}^-\sim[d] & \\
K_{2m-n}(T)\ar^{k_*}[r]
 & K_{2m-n}(X)
}
$$
where the vertical maps are the isomorphisms \eqref{eq:K-theory&module_general}.
\end{cor}

\subsection{Trace maps}

\begin{num} Let $S$ be a regular scheme.
Let $Y$ be a smooth $S$-scheme.
The obvious map $\mathrm{Pic}(Y) \To K_0(Y)$
together with the canonical maps
$$
K_0(Y)
 \xrightarrow{\sim} \Hom_\BGL(\BGL_S(Y),\BGL_S)
 \xrightarrow{\beta_*} \Hom_\BGL(\BGL_S(Y),\BGL_S(1)[2])
$$
defines Chern classes\index{word}{class!Chern}
in the category $\ho(\Mod{\BGL})(S)$;
they corresponds to the orientation defined in
 remark \ref{orientation_KGL}. \\
Let $p:P \rightarrow S$ be a projective bundle of rank $n$.
Let $v=[\mathcal O(1)]-1$ in $K_0(P)$.
It corresponds to a map $\mathfrak v:\BGL_S(P) \rightarrow \BGL_S$.
According to \cite[3.2]{Deg8} and our choice of Chern classes,
 the following map
$$
\BGL_S(P)
 \xrightarrow{\sum_i \beta^i.\mathfrak v^i \boxtimes p_*}
  \bigoplus_{0 \leq i \leq n} \BGL_S(i)[2i]
$$
is an isomorphism. As $\beta$ is invertible,
 it follows that the map
\begin{equation} \label{eq:iso_deloope}
\varphi_{P/S}:\BGL_S(P)
 \xrightarrow{\sum_i \mathfrak v^i \boxtimes p_*}
  \bigoplus_{0 \leq i \leq n} \BGL_S
\end{equation}
is an isomorphism as well. 
Using this formula, the map $\Hom(\varphi_{P/S},\BGL_S)$ is equal
 to the isomorphism of Quillen's projective bundle theorem in K-theory
 (\textit{cf.} \cite[4.3]{Quillen}):
$$
f_{P/S}:\bigoplus_{i=0}^n K_*(S) \rightarrow K_*(P),
 (S_0,...,S_n) \mapsto \sum_i p^*(S_i).v^i.
$$

Let $p_*:K_*(P) \rightarrow K_*(S)$ be the pushout by the
 projective morphism $p$. According to the projection formula,
 it is $K_*(S)$-linear. In particular, 
 it is determined by the $n+1$-uple $(a_0,...,a_n)$ 
 where $a_i=p_*(v^i) \in K_0(S)$
 through the isomorphism $f_{P/S}$.
 Let $\mathfrak a_i:\BGL_S \rightarrow \BGL_S$ be the map corresponding
 to $a_i$.
\end{num} 
 
\begin{df}
Consider the previous notations. We define the \emph{trace map}
\index{word}{map, trace}
 associated with the projection $p:P \rightarrow S$ 
 as the morphism of $\BGL$-modules $\Tr^\BGL_p:p_*(\BGL_P)\To\BGL_S$
 determined as the composite
$$
p_*(\BGL_P)=\derR\sHom(\BGL_S(P),\BGL_S)
 \xrightarrow{(\varphi_{P/S}^*)^{-1}} 
 \bigoplus_{i=0}^n \BGL_S \xrightarrow{(\mathfrak a_0,\hdots,\mathfrak a_n)}
 \BGL_S.
$$
\index{notat}{TrpKGL@$\Tr^\BGL_p$}
\end{df}
From this definition, it follows that $\mathrm{Tr}_p$ represents
 the push-forward by $p$ in K-theory:
$$
\xymatrix@C=28pt@R=8pt{
\Hom_\BGL(\BGL_S[n],p_*\BGL_P)\ar^-{\Tr^\BGL_{p*}}[r]\ar@{=}[d]
 & \Hom_\BGL(\BGL_S[n],\BGL_S)\ar^{\epsilon_S}[dd] \\
\Hom_\BGL(\BGL_P[n],\BGL_P)\ar_{\epsilon_P}[d] &  \\
K_n(P)\ar^-{p_*}[r] & K_n(S)
}
$$
Consider moreover a cartesian square:
$$
\xymatrix{
Q\ar^q[r]\ar_g[d] & P\ar^p[d] \\
Y\ar^f[r] & S
}
$$
such that $f$ is smooth. From the projective base change theorem,
 we get $f^*p_*p^*=q_*q^*g^*$. Using this identification,
 we easily obtain that $f^*\Tr^\BGL_p=\Tr^\BGL_q$.
Thus, we conclude that the map
$$\Hom_\BGL(\BGL_S(Y)[n],p_*\BGL_P)\xrightarrow{\Tr^\BGL_p}
\Hom_\BGL(\BGL_S(Y)[n],\BGL_S)$$
represents the usual pushout map
$$q_*:K_n(Q) \rightarrow K_n(Y)$$
through the canonical isomorphisms \eqref{eq:K-theory&module_general}.

\begin{num} Consider a projective morphism $f:T \rightarrow S$ 
 between regular schemes and choose a factorization
$$
T \xrightarrow i P \xrightarrow p S
$$
where $i$ is a closed immersion and $p$ is the projection of a projective
bundle. Let us define a morphism
$$
\Tr^\BGL_{(p,i)}:f_*\BGL_T=p_*i_*\BGL_T
 \xrightarrow{p_*\eta_i} p_*\BGL_P
 \xrightarrow{\Tr^\BGL_p} \BGL_S.
$$
According to \ref{cor:fund_class&pushout} and the previous paragraph,
 for any cartesian square
$$
\xymatrix{
Y\ar^g[r]\ar_b[d] & X\ar^a[d] \\
T\ar^f[r] & S
}
$$
such that $a$ is smooth,
 the following diagram is commutative.
\begin{equation} \label{eq:traces&pushout}\begin{split}
\xymatrix@R=18pt@C=28pt{
\Hom_\BGL(\BGL_S(X),f_*\BGL_T(m)[n])\ar^{\Tr^\BGL_{(p,i)*}}[r]\ar@{=}[d]
 & \Hom_\BGL(\BGL_S(X),\BGL_S(m)[n])\ar^-{\epsilon_{X/S}}_-\simeq[dd] \\
\Hom_\BGL(\BGL_T(Y),\BGL_Z(m)[n])\ar_-{\epsilon_{Y/T}}^-\simeq[d] & \\
K_{2m-n}(Y)\ar^{g_*}[r]
 & K_{2m-n}(X)
}\end{split}
\end{equation}
%Because the class of $\BGL$-modules of the form $\BGL_S(X)(-m)[-n]$
% is a generating family in $\BGL-mod_S$, we deduce from the Yoneda lemma
% that $\Tr^\BGL_{(p,i)}$ does not depend on the chosen factorization:
\end{num}

\begin{df} \label{df:trace_BGL}
Considering the above notations, we define the \emph{trace map}
associated to $f$ as the morphism
$$
\Tr^\BGL_f=\Tr^\BGL_{(p,i)}:f_*f^*\BGL_S \rightarrow \BGL_S.
$$
\end{df}

\begin{rem}
By definition, the trace map $\Tr^\BGL_f$ is a morphism
of $\BGL$-modules. As a consequence, the map
obtained by adjunction
$$\eta'_f:\BGL_T\simeq f^*\BGL_S\To f^!\BGL_S $$
is also a morphism of $\BGL$-module.
This implies that the morphism $\eta'_f$
(and thus also $\Tr^\BGL_f$) is completely determined by the element
$$\eta'_f\in\Hom_\BGL(\BGL_T,f^!\BGL_S)\simeq\Hom_{\SH(T)}(\unit_T,f^!\BGL_S)\, .$$
Moreover, as $p$ is smooth, there is a canonical isomorphism
$p^!\BGL_S\simeq \BGL_P$ (by relative purity for $p$
and by periodicity; see \cite[lemma 6.1.3.3]{RiouBGL}).
From there, we deduce from Theorem \ref{abspurityBGL}
that we have a canonical isomorphism
$$f^!\BGL_S\simeq i^!\BGL_P\simeq \BGL_T\, .$$
This implies that we have an isomorphism:
$$\Hom_{\SH(T)}(\unit_T,f^!\BGL_S)\simeq K_0(T)\, .$$
Hence, the map $\eta'_f$ is completely determined by a class in $K_0(T)$.
The problem of the functoriality of trace maps in the motivic
category $\ho(\Mod \BGL)$
is thus a matter of functoriality of this
element $\eta'_f$ in $K_0$, which can be translated faithfully to the problem
of the functoriality of pushforwards for $K_0$.

However, the only property of trace maps we shall use here is the following.
\end{rem}

\begin{prop}\label{flattracedegree}
Let $f:T \rightarrow S$ be a finite flat morphism of regular schemes such that
 the $\mathcal O_S$-module $f_*\mathcal O_T$ is (globally) free of rank $d$.
 Then the following composite map
$$
\BGL_S
\To f_*f^* \BGL_S
\xrightarrow{\Tr^\BGL_f} \BGL_S
$$is equal to $d.1_{\BGL_S}$ in $\ho(\Mod{\BGL})(S)$ (whence in $\SH(S)$ as well).
\end{prop}
\begin{proof}
Let $\varphi$ be the composite map of $\ho(\Mod{\BGL})(S)$
$$
\BGL_S \To f_*f^*\BGL_S
 \xrightarrow{\Tr_f} \BGL_S\, .
$$
As $\varphi$ is $\BGL_S$-linear by construction, it corresponds
to an element
$$\varphi \in \Hom_\BGL(\BGL_S,\BGL_S)
\simeq\Hom_{\SH(S)}(\unit_S,\BGL_S)\simeq K_0(S)\, .$$
According to the commutative diagram \eqref{eq:traces&pushout},
 if we apply the global sections functor
$\Hom_{\SH(S)}(\unit_S,-)$ to $\varphi$, we obtain
through the evident canonical isomorphisms the composition
of the usual pullback and pushforward by $f$ in K-theory:
$$
K_{0}(S) \xrightarrow{f^*} K_{0}(T)
 \xrightarrow{f_*} K_{0}(S).
$$
With these notations, the element of $K_0(S)$ corresponding to $\varphi$ is the
pushforward of $1_T=f^*(1_S)$ by $f$,
while the element corresponding to the identity of $\BGL_S$ is of course $1_S$.
Under our assumptions on $f$, it is obvious that we have the
identity $f_*(1_T)=d.1_S\in K_0(S)$.
This means that $\varphi$ is $d$ times the
 identity of $\BGL_S$.
\end{proof}

\section{Beilinson motives} \label{sec:Beilinson_motives}

\renewcommand{\Rc}{\QQ}

\subsection{The {$\gamma$}-filtration}

\begin{num} \label{num:gamma_filtration}
We denote by $\BGLr$
\index{notat}{KGLQ@$\BGLr$}
 the $\QQ$-localization of the absolute ring spectrum $\BGL$,
 considered as a cartesian section of $\DMt$.
From \cite[5.3.10]{RiouBGL}, this spectrum has the following property:
\begin{enumerate}
\item[(K5)] For any scheme $S$, there exists a canonical decomposition,
 called the \emph{Adams decomposition}
\index{word}{decomposition, Adams}
$$
\BGLrx S \simeq \bigoplus_{i \in \ZZ} \BGLrn i_S
$$
\index{notat}{KGLSi@$\BGLrn i_S$}
compatible with base change and such that for any regular scheme $S$,
the isomorphism of (K2) induces an isomorphism:
$$
\Hom_{\DMtx S}\Big(\QQ_S(X)[n],\BGLrn i_{S}\Big)
 \simeq K^{(i)}_n(X):=Gr_\gamma^i K_n(X)_\QQ
$$
where the right-hand side is the $i$-th graded piece
of the $\gamma$-filtration
\index{word}{filtration, $\gamma$-filtration}
 on K-theory groups.
\end{enumerate}
We will denote by 
\begin{align*}
& \pi_i:\BGLrx S \rightarrow \BGLrn i_S, \\
\text{resp. } &\iota_i:\BGLrn i_S \rightarrow \BGLrx S
\end{align*} the projection (resp. inclusion) defined 
by the decomposition (K3) and we put $p_i=\iota_i \pi_i$
for the corresponding projector on $\BGLrx S$.
\end{num}

\begin{df}[Riou] We define
the \emph{Beilinson motivic cohomology spectrum}
\index{word}{spectrum!Beilinson motivic cohomology spectrum}
 as the rational Tate spectrum $\HBx S=\BGLrn 0_S$.\label{df:Beilinson_mot_coh}
\end{df}

\begin{rem}
Note that, by definition, for any morphism of schemes $f:T\To S$,
we have $f^*\HBx S\simeq \HBx T$.
\end{rem}

\begin{lm} \label{lm:gr_adams&twists}
The isomorphism $\gamma_u$ of \eqref{eq:inv_Bott_iso} is homogeneous of degree $+1$ 
with respect to the graduation (K5).
In other words, for any integer $i \in \ZZ$,
 the following composite map is an isomorphism
$$
\BGLrn i(1)[2] \xrightarrow{\iota_i} \BGLr(1)[2]
 \xrightarrow{\gamma_u} \BGLr \xrightarrow{\pi_i} \BGLrn {i+1}.
$$
\end{lm}
For any integer $i \in \ZZ$, we thus get a canonical isomorphism
\begin{equation} \label{eq:gr_adams&twists_iso}
\HB(i)[2i] \xrightarrow{\sim} \BGLrn i.
\end{equation}
\begin{proof}
It is sufficient to check that, for $j \neq i+1$,
$$
\begin{cases}
p_j \circ \gamma_u \circ p_i=0, & \\
p_j \circ \gamma_u^{-1} \circ p_i=0 &
\end{cases}
$$
in $\Hom_{\DMtx S}(\BGLr,\BGLr)$. But according to \cite[5.3.1 and 5.3.6]{RiouBGL},
we have only to check these equalities for the induced endomorphism of $K_0$
(seen as a presheaf on the category of smooth schemes over $\spec \ZZ$).
This follows then from the compatibility of the projective bundle isomorphism
with the $\gamma$-filtration; see \cite[Exp.~VI, 5.6]{SGA6}.
\end{proof}

\begin{num} \label{num:HB=LQ}
Recall from \cite{NOS} that $\BGLr$ is canonically
isomorphic (with respect to the orientation \ref{orientation_KGL})
to the universal oriented rational ring spectrum with
multiplicative formal group law introduced in \cite{NOS}.
The isomorphism of the preceding corollary shows
 in particular that $\HB$ is obtained from $\BGLr$
 by killing the elements $\beta^n$ for $n \neq 0$.
In particular, this shows that $\HB$ is canonically
isomorphic to the spectrum denoted by $\mathsf L\QQ$
in \cite{NOS}, which corresponds to the universal
additive formal group law over $\QQ$. This implies that
$\HB$ has a natural structure of a (commutative) ring spectrum.
\end{num}

\begin{prop}\label{cor:HB^2=HB}
The multiplication map
$$
\mu:\HB \otimes \HB \rightarrow \HB
$$
is an isomorphism.
\end{prop}
This trivially implies that the following map is an isomorphism:
\begin{equation}\label{HBHBlocal}
1 \otimes \eta:\HB \rightarrow \HB \otimes \HB.
\end{equation}
\begin{proof}
It is enough to treat the case $S=\spec \ZZ$.
We will proove that the projector
$$
\psi:\HB \otimes \HB \xrightarrow \mu \HB
 \xrightarrow{1 \otimes \eta} \HB \otimes \HB
$$
is an isomorphism (in which case it is in fact the identity).
We do that for the isomorphic ring spectrum $\mathsf L\QQ$.

Let $\mathsf H^{top}\QQ$ be the topological spectrum representing
 rational singular cohomology. In the terminology of \cite{NOS},
 $\mathsf L\QQ$ is a Tate spectrum representing the
 Landweber exact cohomology
\index{word}{cohomology!Landweber exact}
 which corresponds to the  Adams graded $MU_*$-algebra $\QQ$ 
 obtained by killing
 every generator of the Lazard ring $MU_*$. The corresponding
 topological spectrum is of course $\mathsf H^{top}\QQ$. \\
According to \cite[9.2]{NOS}, the spectrum
 $E=\mathsf L\QQ \otimes \mathsf L\QQ$ is a Landweber
 exact spectrum corresponding to the $MU_*$-algebra
 $\QQ \otimes_{MU_*} \QQ=\QQ$.
In particular, the corresponding topological spectrum
 is simply $\mathsf H^{top}\QQ$.
Thus, according to \cite[9.7]{NOS}, applied with
 $F=E=\mathsf L\QQ \otimes \mathsf L\QQ$,
 we get an isomorphism of $\QQ$-vector spaces
$$
\End(\mathsf L\QQ \otimes \mathsf L\QQ)
 =\Hom_\QQ(\QQ,E_{**})=\QQ.
$$
Thus $\psi=\lambda.Id$ for $\lambda \in \QQ$.
But $\lambda=0$ is excluded 
 because $\psi$ is a projector on a non-trivial factor, so that we can conclude.
\end{proof}

\subsection{Definition}

\begin{df} \label{df:Beilinson_motives}
Let $S$ be any scheme.

We say that an object $E$ of $\DMtx S$ is \emph{$\HB$-acyclic}
\index{word}{acyclic, $\HB$-acyclic}
if $\HB \otimes E=0$ in $\DMtx S$. A morphism of $\DMtx S$ is an \emph{$\HB$-equivalence}
\index{word}{equivalence!HBequivalence@$\HB$-equivalence}
if its cone is $\HB$-acyclic
 (or, equivalently, if its tensor product with $\HB$
is an isomorphism).

An object $M$ of $\DMtx S$ is \emph{$\HB$-local}
\index{word}{local!HBlocal@$\HB$-local}
 if, for any $\HB$-acyclic
object $E$, the group $\Hom(E,M)$ vanishes.

We denote by $\DMB(S)$
\index{notat}{DMBS@$\DMB(S)$}
 the Verdier quotient of $\DMtx S$ by
the localizing subcategory made of $\HB$-acyclic objects
(i.e. the localization of $\DMtx S$ by the class of $\HB$-equivalences).
 
The objects of $\DMB(S)$ are called the \emph{Beilinson motives}.
\index{word}{motive!Beilinson}
\end{df}

\begin{prop}
An object $E$ of $\DMtx S$ is $\HB$-acyclic if and only if we have $\BGLr\otimes E=0$.
\end{prop}

\begin{proof}
This follows immediately from property (K5) (see \ref{num:gamma_filtration})
 and Lemma \ref{lm:gr_adams&twists}.
\end{proof}

\begin{prop}\label{beilinsoncmf}
The localization functor $\DMtx S\To\DMB(S)$ admits a fully faithful right adjoint
whose essential image in $\DMtx S$ is the full subcategory spanned by
$\HB$-local objects. More precisely, there is a left Bousfield localization of the stable
model category of symmetric Tate spectra
$\Spt(S,\QQ)$ by a small set of maps whose homotopy category is
precisely $\DMB(S)$.
\end{prop}

\begin{proof}
For each smooth $S$-scheme $X$ and any integers $n,i\in \ZZ$, we have a functor
with values in the category of $\QQ$-vector spaces
$$F_{X,n,i}=\Hom_{\DMtx S}(\Sigma^\infty\QQ_S(X),\HB\otimes(-)(i)[n]):
\Spt(S,\QQ)\To\Mod \QQ$$
which preserves filtered colimits.
We define the class of $\HB$-weak equivalences as the class of maps of
$\Spt(S,\QQ)$ whose image by $F_{X,n,i}$ is an isomorphism for all $X$ and $n$, $i$.
By virtue of \cite[Prop.~1.15 and 1.18]{beke1}, we can apply
Smith's theorem \cite[Theorem 1.7]{beke1} (with the
class of cofibrations of $\Spt(S,\QQ)$), which implies
the proposition.
\end{proof}

\begin{rem}
We shall often make the abuse of considering $\DMB(S)$
as a full subcategory in $\DMt(S)$, with an implicit reference to
the preceding proposition.
\end{rem}

Note that $\HB$-acyclic objects are stable
 by the operations $f^*$, $f_\sharp$ and $\otimes$,
 so that applying Corollary \ref{cor:derived_premotivic_localization_exists},
 we obtain a premotivic category $\DMB$ together with a premotivic adjunction:
\begin{equation} \label{eq:premotivic_adj_DMt_DMB}
\beta^* :\DMt \rightleftarrows \DMB:\beta_*.
\end{equation}

\begin{prop}\label{HBHBloc}
The spectrum $\HBx S$ is $\HB$-local and the unit map
$\eta_{\HB}:\un \rightarrow \HBx S$ is an $\HB$-equivalence in $\DMtx S$.
\end{prop}
\begin{proof}
The unit map $\eta:\un_S \rightarrow \HBx S$
is an $\HB$-equivalence by \ref{cor:HB^2=HB}.

Consider a rational spectrum $E$ over $S$ such that $E \otimes \HB=0$
and a map $f:E \rightarrow \HB$. It follows trivially from the commutative diagram
$$
\xymatrix@R=20pt@C=30pt{
E\ar^f[r]\ar_{1 \otimes \eta}[d]
 & \HBx S\ar_{1 \otimes \eta}[d]\ar@{=}[rd] & \\
E \otimes \HBx S\ar^-{f \otimes 1}[r] & \HBx S \otimes \HBx S\ar^\mu[r] & \HBx S
}$$
that $f=0$, which shows that $\HBx S$ is $\HB$-local.
\end{proof}

\begin{cor}\label{strictifyHB}
The family of ring spectra $\HBx S$ comes from a cofibrant cartesian
commutative monoid (\ref{PfibredRMod}) of the
symmetric monoidal fibred model category of Tate spectra
over the category of schemes.
\end{cor}

\begin{proof}
By virtue of Proposition \ref{HBHBloc}
and of Corollary \ref{strictifyloccommalg}, there exists
a cofibrant commutative monoid in the model category
of symmetric Tate spectra over $\spec \ZZ$ which is canonically isomorphic
to $\HBx \ZZ$ in $\DMtx {\spec \ZZ}$ (as commutative ring spectrum).
For a morphism of schemes $f:S\To \spec \ZZ$,
we can then define $\HBx S$ as the pullback of $\HBx \ZZ$
(at the level of the model categories); using Proposition \ref{Qlinearhomotopyalgebras},
we see that this defines a cofibrant cartesian commutative monoid on the fibred category
of spectra which is isomorphic to $\HBx S$
as commutative ring spectra in $\DMtx S$.
\end{proof}

\begin{num} From now on, we shall assume that $\HB$
\index{notat}{HB@$\HB$}
 is given by a cofibrant cartesian commutative monoid
of the symmetric monoidal fibred model category of Tate spectra
over the category of schemes. By virtue of
propositions \ref{cmfPfibredRMod} and \ref{hmtlinearproperties}),
we get the motivic category $\ho(\Mod{\HB})$ of $\HB$-modules,
\index{word}{module!HBmodule@$\HB$-module}
together with a commutative diagram of morphisms of premotivic categories
$$
\xymatrix@R=12pt{
\DMt \ar[rd]_{\beta} \ar[rr]^{{\HB \otimes(-)}} & & \ho(\Mod{\HB}) \\
&{\DMB} \ar[ru]_{\varphi} & }
$$
(any $\HB$-acyclic object becomes null in the homotopy
category of $\HB$-modules by definition, so that $\HB\otimes(-)$
factors uniquely through $\DMB$ by the universal property of
localization).
\end{num}

\begin{prop}\label{HBmodulesfullyfaithful}
The forgetful functor $U:\ho(\Mod\HB)(S)\To \DMtx S$ is fully faithful.
\end{prop} 
 
\begin{proof}
We have to prove that, for any $\HBx S$-module $M$, the map
$$\HBx S\otimes M\To M$$
is an isomorphism in $\DMt(S)$. As this is a natural
transformation between exact functors which commute with
small sums, and as $\DMt$ is a compactly generated
triangulated category, it is sufficient to check this for
$M=\HBx S\otimes E$, with $E$ a (compact) object
of $\DMt(S)$ (see \ref{Modcompactgenerators}). In this case, this
follows immediately from the isomorphism \eqref{HBHBlocal}.
\end{proof}

\begin{thm} \label{thm:DMB&HB-mod}
The functor $\DMB(S) \To\ho(\Mod{\HBx S})$ is
 an equivalence of triangulated monoidal categories.
\index{word}{equivalence!of triangulated monoidal categories}
\end{thm}
\begin{proof}
This follows formally from the preceding proposition by definition
of $\DMB$ (see for instance \cite[Chap.~I, Prop.~1.3]{GZ}).
\end{proof}

\begin{rem}
The preceding theorem shows that the premotivic category of
$\HB$-modules $\ho(\Mod\HB)$
as well as the morphism $\DMt\To\ho(\Mod\HB)$
are completely independent of the choice of the strictification of the (commutative) monoid
structure on $\HB$ given by Corollary \ref{strictifyHB}.
\end{rem}

\begin{cor}\label{cor:DMBmotivic}
The premotivic category $\DMB\simeq\ho(\Mod{\HB})$
 is a $\QQ$-linear motivic category.
\end{cor}

\begin{proof}
It follows from Proposition \ref{hmtlinearproperties}
 and Theorem \ref{thm:DMB&HB-mod} that $\DMB$ 
 satisfies the homotopy, stability and localization properties
 (because this is true for $\DMt$ by \ref{cor:DM_tilde&univ_derived_motivic}).
 It is also well generated because it is a localization
 of $\DMt$.
 Thus we can apply Remark \ref{rem:compatly&motivic}
 to conclude.
\end{proof}

\begin{rem}
One can also prove that $\DMB$ is motivic
much more directly: this follows from the fact that $\DMt$ is motivic
and that the six Grothendieck operations preserve $\HB$-acyclic objects,
so that all the properties of $\DMt$ induce their analogs on $\DMB$
by the $2$-universal property of localization (we leave this
as an easy exercise for the reader).
\end{rem}

\begin{df}\label{df:beilinsonmotcoh}
For a scheme $X$, we define its Beilinson motivic cohomology
\index{word}{cohomology!Beilinson motivic}
  by the formula:
$$\HB^q(X,\QQ(p))=\Hom_{\DMB(X)}(\unit_X,\unit_X(p)[q])\, .$$
\index{notat}{HBqXQp@$\HB^q(X,\QQ(p))$}
\end{df}
In fact, according to the preceding corollary,
 the cohomology theory defined above is represented
 by the ring spectrum $\HB$.
 In particular,
 we can now justify the terminology of Beilinson motives:
\begin{cor}\label{computeHBinDMB}
For any regular scheme $X$, we have a canonical isomorphism
$$\HB^q(X,\QQ(p))\simeq Gr_\gamma^{p} \, K_{2p-q}(X)_\QQ\, .$$
\end{cor}

\begin{num}
Recall from Paragraph \ref{num:HB=LQ} that $\HBx S$
 is canonically oriented for any scheme $S$.
 Moreover, these orientations are compatible with pullbacks
 with respect to $S$.
 This means in particular that the motivic triangulated category
 $\DMB$ is oriented (see Example \ref{ex:oriented_spectra}).

In particular, the fibred category $\DMB$ satisfies
 the usual Grothendieck 6 functors formalism.
\index{word}{formalism, Grothendieck 6 functors}
  We refer the reader to Theorem \ref{thm:cor3_Ayoub}
   for the precise statement.
\end{num}

It was remarked in Paragraph \ref{num:HB=LQ} that
 $\HBx S$ is the universal oriented ring spectrum with
 additive formal group law over $S$.
\index{word}{spectrum!universal oriented ring ---- with
 additive formal group law}
 This property can be expressed by the following nice description
 of Beilinson motives:
\begin{cor} \label{cor:carac_KBL-local}
Let $E$ be a rational spectrum over $S$. The following conditions are equivalent:
\begin{enumerate}
\item[(i)] $E$ is a Beilinson motive (i.e. is in the essential image of the right
adjoint of the localization functor $\DMt\To\DMB$);
\item[(ii)] $E$ is $\HB$-local;
\item[(iii)] the map $\eta \otimes 1_E:E \rightarrow \HB \otimes E$ is an isomorphism;
\item[(iv)] $E$ is an $\HB$-module in $\DMt$;
\index{word}{module!HBmodule@$\HB$-module}
\item[(v)] $E$ is admits a strict $\HB$-module structure.
\index{word}{module!strictHBmodule@strict $\HB$-module}
\end{enumerate}
If, in addition, $E$ is a commutative ring spectrum, these conditions are equivalent to
 the following ones:
\begin{enumerate}
\item[(Ri)] $E$ is orientable;
\index{word}{orientable}
\item[(Rii)] $E$ is an $\HB$-algebra;
\index{word}{algebra!HBalgebra@$\HB$-algebra}
\item[(Riii)] $E$ admits a unique structure of $\HB$-algebra;
\end{enumerate}
And, if $E$ is a strict commutative ring spectrum, these conditions
are equivalent to the following conditions:
\begin{enumerate}
\item[(Riv)] there exists a morphism of commutative monoids $\HB\To E$ in the stable model
category of Tate spectra;
\item[(Rv)] there exists a unique morphism $\HB\To E$ in the homotopy category
of commutative monoids of the category of Tate spectra.
\end{enumerate}
\end{cor}

\begin{proof}
The equivalence between statements (i)--(v) follows immediately from \ref{thm:DMB&HB-mod}.
If $E$ is a ring spectrum,
the equivalence with (Ri), (Rii) and R(iii) is a consequence of \ref{equivorientations}
and of the fact that $\MGL_\QQ$ is $\HB$-local; see \cite[Cor. 10.6]{NOS}.
It remains to prove the equivalence with (Riv) and (Rv).
Then, $E$ is $\HB$-local if and only if the map $E\To \HB\otimes E$
is an isomorphism. But this map can be seen as a morphism of
strict commutative ring spectra (using the model structure of \ref{cmfQlinearcomm}
applied to the model category of Tate spectra)
whose target is clearly an $\HB$-algebra,
so that (Riv) is equivalent to (ii). It remains to check that there is at most one
strict $\HB$-algebra structure on $E$ (up to homotopy), which follows from
the fact that $\HB$ is the initial object in the
homotopy category of commutative monoids of the model
category given by Theorem \ref{cmfQlinearcomm}
applied to the model structure of Proposition \ref{beilinsoncmf}.
\end{proof}

%% \num \label{num:uDMB}
%% Consider the enlargement (\textit{cf.} \ref{})
%% $$
%% \varphi_!:\DMtr \leftrightarrows \uDMtr:\varphi^!.
%% $$
%% We can extend the defintion of $\BGL$-acyclic (resp. $\BGL$-local, $\HB$-local)
%%  objects to the generalized category $\uDMtr$.
%%  We let $\uDMB(S)$ be the localization of $\uDMtr$ with respect to
%%  category of $\BGL$-acyclic generalized spectra. From this construction,
%%  we thus obtain a commutative square
%% $$
%% \xymatrix@R=12pt@C=20pt{
%% \DMtr\ar^{F_\mathcyr B}[r]\ar_{\varphi_!}[d] & \DMB\ar^{\psi_!}[d] \\
%% \uDMtr\ar^{\underline{F}_\mathcyr B}[r] & \uDMB
%% }
%% $$
%% where the vertical maps are enlargement and the horizontal one are premotivic
%%  adjunctions. \\
%% Moreover, we can extend Theorem \ref{thm:DMB&HB-mod} to the generalized
%%  context and we get an equivalence of generalized premotivic categories
%% $$
%% \umod{\HB} \xrightarrow{\ \underline \varphi \ } \uDMB.
%% $$
\begin{cor} \label{cor:HB[t,t^-1]=KGL}
One has the following properties.
\begin{enumerate}
\item The ring structure on the spectrum $\HB$ is
given by the following structural maps (with the notations
of \ref{num:gamma_filtration}).
\begin{align*}
& \HB \otimes \HB
 \xrightarrow{\iota_0 \otimes \iota_0} \BGLr \otimes \BGLr
  \xrightarrow{\mu_{\BGL}} \BGLr \xrightarrow{\pi_0} \HB, \\
& \QQ \xrightarrow{\eta_{\BGL}} \BGLr \xrightarrow{\pi_0} \HB.
\end{align*}
\item The map $\imath_0:\HB \rightarrow \BGLr$ is compatible
 with the monoid structures.
\item Let $\HB[t,t^{-1}]=\bigoplus_{i\in \ZZ}\HB(i)[2i]$ be the free $\HB$-algebra
generated by one invertible generator $t$ of bidegree $(2,1)$.
Then the section $u:\QQ(1)[2] \rightarrow \BGLr$
 induces an isomorphism of $\HB$-algebras:
$$
\gamma'_u:\HB[t,t^{-1}] \rightarrow \BGLr.
$$
\end{enumerate}
\end{cor}

\begin{proof}
Property (1) follows from properties (2) and (3).
Property (2) is a trivial consequence of the previous
corollary. Using the isomorphisms \eqref{eq:gr_adams&twists_iso}
of Lemma \ref{lm:gr_adams&twists}, we get a canonical isomorphism 
$$
\HBx S[t,t^{-1}]
 \xrightarrow \sim \bigoplus_{i \in \ZZ} \BGLrn i.
$$
Through this isomorphism, the map $\gamma'_u$ corresponds
to the Adams decomposition
\index{word}{decomposition, Adams}
 (i.e. to the isomorphism (K5) of \ref{num:gamma_filtration})
 from which we deduce property (3).
\end{proof}

\begin{rem} \label{rm:strict_Kth}
One deduces easily, from the preceding proposition and from \ref{cor:HB^2=HB},
another proof of the fact that $\BGLr$ is a strict commutative
ring spectrum.

The isomorphism (3) is in fact compatible with the grading
of each term: the factor $\HB.t^i$ is sent to the factor 
$\BGLrn i$. Recall also the parameter $t$ corresponds
to the unit $\beta^{-1}$ in $\BGL^{*,*}$.
\end{rem}

%% In the same way, using an extension of the techniques of \cite{RiouBGL},
%%  we get from the compatibility of the Adams graduation on K-groups with the
%%  multiplicative structure the following lemma:
\begin{cor}
The Adams decomposition is compatible with the monoid structure on $\BGLr$:
for any integer $i,j,l$ such that $l \neq i+j$, 
 the following map is zero.
$$
\BGLrn i \otimes \BGLrn j
 \xrightarrow{\iota_i \otimes \, \iota_j} \BGLr \otimes \BGLr
 \xrightarrow{\ \mu \ } \BGLr
 \xrightarrow{\ \pi_l \ } \BGLrn l
$$
\end{cor}

\begin{num} \label{num:Beilinson_mot_coef}
\index{word}{coefficients, for Beilinson motives}
Let $R$ be a $\QQ$-algebra with structural morphism $\varphi$.
Recall from Paragraph \ref{num:stable-A^1-derived_chg_coef}
 that we get an adjunction of premotivic triangulated categories:
$$
\varphi^*:\DMt \rightarrow
\renewcommand{\Rc}{R}
 \DMt:\varphi_*.
$$
Moreover, for any object $M$ and $N$ of $\DMt(S)$,
 the canonical map
\begin{equation} \label{eq:DMt_chg_coef}
\Hom(M,N) \otimes_\QQ R \rightarrow
\Hom(\varphi^*(M),\varphi^*(N)).
\end{equation}
is an isomorphism provided $M$ is compact
 or $R$ is a finite $\QQ$-vector space.

In particular, the ring spectrum
$
\BGL_R:=\varphi^*(\BGLr)
$ 
represents Quillen algebraic K-theory with coefficients in $R$
over regular schemes. We can repeat Definition \ref{df:Beilinson_motives}
with $R$-coefficients and this gives the category
$\DMB(S,R)$ of Beilinson motives with $R$-coefficients together
with an adjunction:
$$
\varphi^*:\DMB \rightarrow \DMB(-,R):\varphi_*.
$$
Moreover, using the canonical map \eqref{eq:DMt_chg_coef}
 and the fact it is an isomorphism when $M$ is a constructible
 Beilinson motives, we immediately extend 
 all the properties proved so far
 from $\QQ$-coefficients to $R$-coefficients.
\end{num}

%%\section{The fundamental theorems}

%%For all the results in this section,
%%we assume that any scheme in $\sch$ is excellent.

\subsection{Motivic proper descent}

Recall from Definition \ref{df:continuous}
 we have defined the notion of continuity for
 a triangulated premotivic category which is the homotopy
 category of a premotivic model category,
 such as the triangulated motivic category $\DMB$
 --- in this case, the notion of continuity is relative
 to the Tate twist.
\begin{prop} \label{DMBcontinuous}
The motivic triangulated category $\DMB$ is continuous.
\end{prop}
\begin{proof}
We consider the adjunction \eqref{eq:premotivic_adj_DMt_DMB}.
 According to Theorem \ref{thm:DMB&HB-mod},
 the functor $\beta_*$ commutes with pullbacks by arbitrary
 morphisms.
Thus the continuity property for $\DMB$ follows from the
 continuity property for $\DMt$ which was established
 in Example \ref{DMtcontinuous}.
\end{proof}

We will give the main applications of continuity
 in the section on constructible Beilinson motives.
 Recall from \ref{Nispointscontinuity} the following
 corollary of the continuity property of the motivic
 category $\DMB$:
\begin{cor}\label{cor:DMBcontinuous}
Let $X$ be a scheme, and consider an $X$-scheme $Y$ of finite type.
Given a point $x \in X$, we denote by $X_x^h$ the spectrum
of the local henselian ring of $X$ at the point $x$.
Let $a_x:Y \times_X X_x^h \rightarrow Y$ be the canonical map.
Then the family of functors
$$
\DMB(Y) \To \DMB(Y \times_X X_x^h) \ , \quad E \mapsto a_x^*(E)
$$
is conservative.
\end{cor}

As the reader might expect,
 this proposition is very useful to reduce global properties 
 of the motivic category $\DMB$ to local properties.
 This is in particular illustrated by the following proposition.
\begin{thm}\label{thm:DMBseparated}
The motivic category $\DMB$ is separated
 (on the category of noetherian schemes of finite dimension).
\end{thm}
\begin{proof}
According to Proposition \ref{prop:redseptoetale},
 it is sufficient to check that,
 for any finite surjective morphism $f:T \rightarrow S$, the pullback functor
$$
f^*:\DMB(S) \rightarrow \DMB(T)
$$
is conservative.

We argue by induction on the dimension of $S$.

Let us first treat the case where $\dim(S)=0$.
Using the localization property,
 we can assume that $S$ and $T$ are reduced
 (\textit{cf.} \ref{prop:easy_csq_loc}).
 Then $S$ is a disjoint sum of spectra of fields.
 In particular, $f$ is not only finite surjective but also flat.
 Moreover, it is also globally free.
It will be  sufficient to prove that,
 for any Beilinson motive $E$ over $S$, the adjunction map
$$
E \rightarrow f_*f^*(E)
$$
is a monomorphism in $\DMB$.
Using the projection formula in $\DMB$
applied to the finite morphism $f$ (point (5) of Theorem \ref{thm:cor3_Ayoub}),
this latter map is isomorphic to
$$
\big(\HB \To f_*f^*(\HB)\big) \otimes 1_E.
$$
We are finally reduced to 
prove that the map $\HBx S \To f_*f^*\HBx S$
is a monomorphism in $\DMB$
 (any monomorphism of a triangulated category splits). 
As $\HBx S$ is a direct factor of $\BGLrx S$,
it is sufficient to find a retraction of the adjunction map
$$
\BGLrx S \To f_*f^*\BGLrx S\, ,
$$
and this follows from Proposition \ref{flattracedegree}.

Let us finally solve the induction process.
Applying the preceding proposition,
 we can assume that $S$ is local henselian.
 Let $s$ be the closed point of $S$ and $U$ be the open complement.
 Let $f_s$ (resp. $f_U$) be the pullback of $f$ above $s$
 (resp. $U$).
 Using the localization property of $\DMB$ and the base change
 isomorphisms (point (4) of Theorem  \ref{thm:cor3_Ayoub}),
 it is sufficient to treat the case of the finite morphisms $f_U$
 and $f_s$. The case of $f_U$ follows by the induction hypothesis
 while the case of $f_s$ follows from the case treated previously.
 This ends up the induction process.
\end{proof}

According respectively to Proposition \ref{prop:charseparated1}
 and Theorem \ref{charseparated3}, 
 we deduce from the preceding proposition the following result:
\begin{thm} \label{DMB_etale&h-descent}
\begin{enumerate}
\item The motivic category $\DMB$ satisfies \'etale descent.
\index{word}{descent!etale@\'etale}
\item The motivic category $\DMB$ satisfies $\h$-descent
\index{word}{descent!hdescent@$\h$-descent}
 when restricted to quasi-excellent schemes.
\index{word}{scheme!quasi-excellent}
\end{enumerate}
\end{thm}
Recall this means that for any \'etale hypercover
 (resp. $\h$-hypercover of a quasi-excellent scheme)
 $p:\X \To X$ and for any Beilinson motive $E$ over $X$,
 the map
$$
p^*:\derR \Gamma(X,E)\To \derR\Gamma(\X,E)
=\derR \varprojlim_n \derR \Gamma(\X_n,E)
$$
is an isomorphism in the derived category
 of the category of $\QQ$-vector spaces
 (see Corollary \ref{reductionelementarymodcatbasic}
 taking into account Definition \ref{defderivedabsoluteglobal sections}).

%% \begin{paragr}\label{plongement1}
%% Consider the the enlargement
%% $$\DMtr(S) \To \uDMtr(S)\, .$$
%% For $t=\qfh$ or $t=h$, we have an adjunction
%% $$\underline{a}^*:\uDMtr(S)\rightleftarrows\DM_{t,\QQ}(S):\underline{a}_*$$
%% (where $\underline{a}^*$ is the functor induced by the $t$-sheafification
%% functor). The functor $\underline{a}_*$ fully faithful and
%% we know that its essential image
%% consists precisely of the objects of $\uDMtr(S)$
%% which satisfy $t$-descent. Moreover, an object
%% of $\DMtr(S)$ satisfies $t$-descent is and only if
%% its image in $\uDMtr(S)$ have the same property.
%% We thus get the following result.
%% \end{paragr}
%% 
%% \begin{cor} \label{DMB&h-descent2}
%% For $t=\qfh$ or $t=h$, fhe functor
%% $$\DMB(S)\xrightarrow{\mathcal{O}_\mathcyr B}\DMtr(S)
%% \To \uDMtr(S)\xrightarrow{\underline{a}^*} \DM_{t,\QQ}(S)$$
%% is fully faithful.
%% \end{cor}

\subsection{Motivic absolute purity}

\begin{thm}[Absolute purity]\index{word}{purity!absolute} Let
$i:Z \rightarrow S$ be a closed immersion between regular schemes.
Assume $i$ has pure codimension $n$.
Then, considering the notations of \ref{num:gamma_filtration},
definition \ref{df:fund_class_BGL}, and the identification \eqref{eq:gr_adams&twists_iso},
the composed map
$$
\HBx Z \xrightarrow{\iota_0}
 \BGLrx Z \xrightarrow{\eta'_i} i^!\BGLrx S \xrightarrow{\pi_n} i^!\HBx S(n)[2n]
$$
is an isomorphism.\label{DMBpurity}
\end{thm}
This isomorphism, of equivalently the map obtained by adjunction:
$$
i_*(\HBx Z) \rightarrow \HBx S(n)[2n]
$$
is called the \emph{fundamental class}
\index{word}{class!fundamental}
 associated with $i$.
In fact, this is a canonical class in the Beilinson motivic cohomology
 of $X$ with support in $Z$ of bidegree $(2n,n)$.

\begin{rem}
It follows from Remark \ref{rem:fundamental_class&pullback}
 that the fundamental class in Beilinson motivic cohomology
 is compatible with pullback with respect to Tor-independent
 square.
\index{word}{pullback!of fundamental class}
\index{word}{square!Tor-independant}
\end{rem}

\begin{proof}
We have only to check that the above composition induces
 an isomorphism after applying the functor
 $\Hom(\QQ_S(X),-(a)[b])$ for a smooth $S$-scheme $X$
 and a couple of integers $(a,b) \in \ZZ^2$.
Using Remark \ref{rem:fundamental_class&pullback}(3),
 this composition is compatible with smooth base change, and we
 can assume $X=S$.
Let us consider the projector
$$
p_{a}:K_r^Z(S)_\QQ=K_r(S/S-Z)_\QQ \rightarrow K_r(S/S-Z)_\QQ
$$
induced by $\pi_a \circ \iota_a : \BGLr \To \BGLr$,
 and denote by $K_r^{(a)}(S/S-Z)$ (with $r=2a-b$) its image.
By virtue of Proposition \ref{prop:key_abs_purity},
we only have to check that the following composite is an isomorphism:
$$
\rho_i:K_r^{(a)}(Z) \xrightarrow{\iota_{a}} K_r(Z)_\QQ 
 \xrightarrow{\, \mathfrak q_i^{-1}} K_r(S/S-Z)_\QQ
 \xrightarrow{\pi_a} K_r^{(a+n)}(S/S-Z).
$$
From \ref{thm:loc_Quillen}, the morphism $\rho_i$ is functorial
 with respect to Tor-independent cartesian squares of regular schemes
 (\textit{cf.} \ref{num:df_Tor_indep}).
Thus, using again the deformation diagram \eqref{eq:deformation},
 we get a commutative diagram
$$
\xymatrix{
K_r^{(a)}(Z)\ar[r]\ar_{\rho_i}[d]
 & K_r^{(a)}(\AA^1_Z)\ar[d]
 & K_r^{(a)}(Z)\ar^{\rho_s}[d]\ar[l] \\
K_r^{(a+n)}(S/S-Z)\ar[r] & K_r^{(a+n)}(D/D-\AA^1_Z)
 & K_r^{(a+n)}(P/P-Z)\ar[l]
}
$$
in which any of the horizontal maps is an isomorphism
 (as a direct factor of an isomorphism).
Thus,
 we are reduced to the case of the closed immersion $s:Z \rightarrow P$,
 canonical section of the projectivization of a vector bundle $E$
 (where $E$ is the normal bundle
 of the closed immersion $i$). Moreover, as the assertion is local on $Z$,
 we may assume $E$ is a trivial vector bundle.

Let $p:P \rightarrow Z$ be the canonical projection,
 $j:P-Z \rightarrow P$ the obvious open immersion.
Considering the element $v':=\big([\cO(1)]-1\big)$ of $K_0(P)$,
 we let $v$ be its projection on the first graded part
 of the $\gamma$-filtration, $v \in K_0^{(1)}(P)$. \\
Recall that, according to the projective bundle formula, 
the horizontal lines in the following commutative diagram
 are split short exact sequences:
$$
\xymatrix@=14pt{
0\ar[r] & K_r(P/P-Z)_\QQ\ar^-\nu[r]\ar[d]
 & K_r(P)_\QQ\ar^{j^*}[r]\ar[d]
  & K_r(P-Z)_\QQ\ar[d]\ar[r] & 0 \\
0\ar[r] & K_r^{(a+n)}(P/P-Z)\ar^-{\nu'}[r] & K_r^{(a+n)}(P)\ar[r]
 & K_r^{(a+n)}(P-Z)\ar[r] & 0.
}
$$
By assumption on $E$, $v^n$ lies in the kernel of $j^*$
 and the diagram allows to identify the graded piece 
 $K_r^{(a+n)}(P/P-Z)$ with the submodule of $K_r^{(a+n)}(P)$
 of the form $K_r^{(a)}(Z).v^n$.

On the other hand, $j^*s_*=0$:
 there exists a unique element $\epsilon \in K_0(Z)$ such that
 $s_*(1)=p^*(\epsilon).v^n$ in $K_0(P)$.
From the relation $p_*s_*(1)=1$,
 we obtain that $\epsilon$ is a unit in $K_0(Z)$,
 with inverse the element $p_*(v^n)$.
 By virtue of \cite[Exp.~VI, Cor.~5.8]{SGA6}, $p_*(v^n)$ belongs to the $0$-th $\gamma$-graded
 part of $K_0(P)_\QQ$ so that the same holds for its inverse $\epsilon$.
In the end, for any element $z \in K_r(Z)$, 
 we get the following expression:
$$
s_*(z)=s_*(1.s^*p^*(z))=s_*(1).p^*(z)=p^*(\epsilon.z).v^n.
$$
Thus, the commutative diagram
$$
\xymatrix@R=18pt@C=26pt{
K_r^{(a)}(Z)\ar[r] & K_r(Z)_\QQ\ar^-{\mathfrak{q}_s^{-1}}[r]\ar_{s_*}[rd]
 & K_r(P/P-Z)_\QQ\ar^-\nu[d]\ar[r]
 & K_r^{(a+n)}(P/P-Z)\ar^-{\nu'}[d] \\
 & & K_r(P)_\QQ\ar[r] & K_r^{(n)}(P) & 
}
$$
implies that the isomorphism $\mathfrak{q}_s^{-1}$
preserves the $\gamma$-filtration (up to a shift by $n$).
Hence, it induces an isomorphism on the graded pieces by functoriality.
\end{proof}

\section{Constructible Beilinson motives} \label{sec:constructible_Beilinson}
\index{word}{constructible!Beilinson motive|(}

\subsection{Definition and basic properties}

In this section, we apply the general results
 of Section \ref{sec:constructible_motives} to the triangulated
 motivic category $\DMB$. Let us first recall
 the definition of constructibility (Def. \ref{deftauconst})
 which corresponds to the Tate twist.
\begin{df}
Given any scheme $S$, we define the category $\DMBc(S)$
\index{notat}{DMBcS@$\DMBc(S)$}
 of \emph{constructible Beilinson motives}
  over $S$
 as the thick triangulated subcategory of
 $\DMB(S)$ generated by the motives of the form
 $\mot S X(i)$ for a smooth $S$-scheme $X$ and an integer $i \in \ZZ$.
\end{df}

\begin{rem}
Constructible Beilinson motives plays towards Beilinson motives
 the same role as complexes of \'etale sheaves with bounded cohomology
 and constructible cohomology sheaves plays against
 complexes of \'etale sheaves (in the case of torsion coefficients
 prime to the residue characteristics).
This fact will be even more striking
 after Theorems \ref{thm:DMBc_6operations}
 and \ref{thm:DMBc_duality}.
\end{rem}

\begin{num}
Recall from Corollary \ref{cor:DM_tilde&univ_derived_motivic}
 that $\DMt$ is compactly generated by the Tate twist.
\index{word}{generated!compactly generated!triangulated $\Pmor$-fibred}
 According to Theorem \ref{thm:DMB&HB-mod},
 the same is true for the motivic category $\DMB$.
 Thus Proposition \ref{constructequivcompact} gives
 the following criterion of constructibility for Beilinson
 motives: 
\end{num}
\begin{prop}
Given any base scheme $S$, a Beilinson motive $\cM$ over $S$
 is constructible
  if and only if it is compact.  
\end{prop}

\begin{rem}
In the sequel, we will give several concrete descriptions
 of the category of constructible Beilinson motives
  (see Corollaries \ref{cor:DMBc&Voevodsky} and \ref{cor:DMBc&Morel}).
\end{rem}

Recall from Proposition \ref{DMBcontinuous} that
 $\DMB$ is continuous (with respect to the Tate twist).
 Proposition \ref{continuityconstructible} thus implies
 the following properties of constructible Beilinson motives:
\begin{prop} \label{DMBc_continuity}
Let $(S_\alpha)_{\alpha \in A}$ be a pro-object of noetherian
 finite dimensional schemes with affine transition maps
 and such that the scheme $S=\plim_{\alpha \in A} S_\alpha$
 is noetherian of finite dimension. Then the canonical functor:
\begin{equation}
2\mbox{-}\varinjlim_\alpha\DMBc(S_\alpha)\To \DMBc(S)
\end{equation}
is an equivalence of monoidal triangulated categories.
\end{prop}

\begin{ex}
Under the assumptions of the above proposition,
 for any couple of integers $(p,q)$, the canonical map
$$
\ilim_{\alpha} \HB^q(S_\alpha,\QQ(p)) \rightarrow \HB^q(S,\QQ(p))
$$
is an isomorphism.\footnote{This result is to be compared with \cite[Sec. 7, 2.2]{Quillen}
 --- it concerns homotopy invariant K-theory rather than K-theory.}
\end{ex}

\subsection{Grothendieck 6 functors formalism and duality}

The motivic triangulated category $\DMB$ is separated
 (\ref{thm:DMBseparated}) and weakly pure
  (see Definition \ref{df:weaklytaupure} ;
   this follows directly from Theorem \ref{DMBpurity}).
 Thus the abstract Theorem \ref{grothendieck6op}
 gives the finiteness theorem,
\index{word}{finiteness theorem}
  which we state here explicitly to help the reader:
\begin{thm}\label{thm:DMBc_6operations}
The triangulated subcategory $\DMBc$ of $\DMB$ is stable
 by the following operations:
\begin{enumerate}
\item $f^*$ for any morphism of schemes $f$.
\item $f_*$ for any morphism $f:Y \rightarrow X$ of finite type
 such that $X$ is quasi-excellent (resp. any proper morphism $f$).
\item $f_!$ for any separated morphism of finite type $f$.
\item $f^!$ for any morphism $f:Y \rightarrow X$ of finite type
 such that $X$ is quasi-excellent.
\item $\otimes_X$ for any scheme $X$.
\item $\uHom_X$ for any quasi-excellent scheme $X$.
\end{enumerate}
\end{thm}
To be more precise, point (1) and (5) are obvious,
 the non respe condition of point (2) is the hardest fact
  and follows from Theorem \ref{thmfinitness},
 point (3) as well as the respe condition of point (2)
  is Corollary \ref{corthmfinitnessproper},
 point (4) is Corollary \ref{constructexcepinvim} and
 point (6) is Corollary \ref{constructibleinternalHom}.

\begin{num}
Let $B$ be an excellent scheme such that $\dim(B) \leq 2$.
Recall that $B$ satisfies wide resolution of singularities
 up to quotient singularities
\index{word}{resolution of singularities!wide ---- up to quotient singularities}
  (see Def. \ref{df:ressingupquotient}
  and the result of De Jong recalled in \ref{cor:dejongdimleq2}).
 Thus according to Corollary \ref{properregulargenerators},
  we get the following description of constructible Beilinson
  motives:
\end{num}
\begin{prop}
Let $S$ be a separated $B$-scheme of finite type,
 and $T \subset S$ a closed subscheme.
 Then the triangulated category $\DMBc(S)$ is the smallest
  triangulated category of $\DMB(S)$ which contained
  motives of the form
$$
f_*(\un_X)(n)
$$
where $n$ is an integer and $f:X \rightarrow S$ is a projective
 morphism such that $X$ is regular connected
 and $f^{-1}(T)_{red}$ is either empty, either $X$ of the support
 of a strict normal crossing divisor.
\end{prop}

The main motivation to introduce the notion of constructibility
 is Grothendieck duality.
\index{word}{duality, Grothendieck}
 We obtain this duality from the theoretical result
  on motivic triangulated categories,
   more precisely Corollary \ref{cor:localduality2}:
\begin{thm} \label{thm:DMBc_duality}
Let $B$ be an excellent scheme such that $\dim(B) \leq 2$
 and $S$ be a regular separated $B$-scheme of finite type.

Then for any separated morphism $f:X \rightarrow S$ of finite type,
 the premotive $f^!(\un_S)$ is a dualizing object of $\DMBc(X)$.
 In fact,
  if we put $D_X(M):=\uHom_X(M,f^!(\un_S))$
   for any constructible Beilinson motives $M$, 
   the following properties hold:
\begin{itemize}
\item[(a)] For any separated $S$-scheme of finite type $X$,
the functor $D_X$ preserves constructible objects.
\item[(b)] For any separated $S$-scheme of finite type $X$, the natural map
$$M\To D_X(D_X(M))$$
is an isomorphism for any constructible Beilinson motive $M$.
\item[(c)] For any separated $S$-scheme of finite type $X$, and
for any Beilinson motive $M$ and $N$ over $X$,
 if $N$ is constructible then we have a canonical isomorphism
$$D_X(M\otimes_X D_X(N))\simeq \uHom_X(M,N)\, . $$
\item[(d)] For any morphism between separated $S$-schemes of finite type
$f:Y\To X$, we have natural isomorphisms
\begin{align*}
D_Y(f^*(M))& \simeq f^!(D_X(M))\\
f^*(D_X(M))& \simeq D_Y(f^!(M))\\
D_X(f_!(N))& \simeq f_*(D_Y(N))\\
f_!(D_Y(N))& \simeq D_X(f_*(N))\\
\end{align*}
where $M$ (resp. $N$) is a constructible Beilinson motive over $X$
 (resp. $Y$).
\end{itemize}
\end{thm}
\index{word}{constructible!Beilinson motive|)}

\begin{num} \label{num:constructible_Beilinson_mot_coef}
\index{word}{coefficients, for Beilinson motives}
Let $R$ be a $\QQ$-algebra.\footnote{The examples we have in mind are:
 $R=E$ is a number field, $R=\CC$, 
 $R=\QQ_l, \bar \QQ_l$ for a prime $l$.} \\
 We define the premotivic triangulated category 
 of constructible Beilinson motives
 with coefficients in $R$ as the category of constructible objects
 of the category $\DMB(-,R)$ defined in 
 Paragraph \ref{num:Beilinson_mot_coef}.

According to \emph{loc. cit.}, for any constructible Beilinson motives
 with coefficients in $\QQ$, we get an isomorphism:
$$
\Hom_{\DMBc(S)}(M,N) \otimes_\QQ R
 \longrightarrow
  \Hom_{\DMBc(S,R)}\big(\derL \varphi^*(M),\derL \varphi^*(N)\big).
$$
It is straightforward to see that this isomorphism
 allows to extend all the results proved so far for
 Beilinson motives with coefficient in $\QQ$ to the case
 of $R$-coefficients.
\end{num}
\section{Comparison theorems}\label{sec:comparisons_Beilison}

\subsection{Comparison with Voevodsky motives} \label{sec:compare_BV}
\renewcommand{\Rc}{\QQ}
\begin{num} We consider the premotivic adjunction of \ref{adjunction:DMt_DM}
\begin{equation}\label{premotivicadjcomparisonthm00}
\gamma^*:\DMt \rightleftarrows \DM:\gamma_*\, .
\end{equation}
For a scheme $S$, $\gamma_*(\un_S)$ is a (strict) commutative ring spectrum,
and, for any object $M$ of $\DM(S)$, $\gamma_*(M)$ is naturally
endowed with a structure of $\gamma_*(\un_S)$-module.
On the other hand, as we have the projective bundle formula in $\DM(S)$ (\ref{pptyChern}),
$\gamma_*(\un_S)$ is orientable (\ref{equivorientations}), which implies that, for any
object $M$ of $\DM(S)$, $\gamma_*(M)$ is
an $\HBx S$-module, whence is $\HB$-local (\ref{cor:carac_KBL-local}).
As consequence, we get a canonical factorization of \eqref{premotivicadjcomparisonthm00}:
\begin{equation}\label{premotivicadjcomparisonthm01}
\DMt \xrightarrow{\beta^*} \DMB \xrightarrow{\varphi^*} \DM.
\end{equation}
Consider the commutative diagram of premotivic categories
\begin{equation}\label{premotivicadjcomparisonthm02}\begin{split}
\xymatrix{
\DMt \ar[r]^{\gamma^*}\ar[d]_{\rho_\sharp}&\DM\ar[d]^{\psi_\sharp} \\
\uDMt \ar[r]^{\underline{\gamma}^*}&\uDMV_\QQ
}\end{split}
\end{equation}
in which the two vertical maps are the canonical
enlargements, and, in particular, are fully faithful (see \ref{prop:enlargement}).

Let $t$ denotes either the $\qfh$-topology
\index{word}{topology!qfhtopology@$\qfh$-topology}
 or the $\h$-topology.
 \index{word}{topology!htopology@$\h$-topology}
We also have the following commutative triangle
\begin{equation}\label{premotivicadjcomparisonthm03}\begin{split}
\xymatrix{
\uDMt \ar[r]^{\underline{\gamma}^*}\ar@/_1.3pc/[rr]_{\underline{a}^*}
&\uDMV_\QQ\ar[r]^{\underline{\alpha}^*}&
\uDMV_{t,\QQ}
}\end{split}
\end{equation}
in which both $\underline{a}^*$ and $\underline{\alpha}^*$
are induced by the $t$-sheafification functor;
see \ref{ex:stable_AA^1-derived_categories} and \ref{adjunctionDM_DMqfh}.
We obtain from \eqref{premotivicadjcomparisonthm01},
\eqref{premotivicadjcomparisonthm02}, and \eqref{premotivicadjcomparisonthm03}
the commutative diagram of premotivic categories below, in which
$\chi_\sharp=\varphi^*\underline{\alpha}^*\psi_\sharp$.
\begin{equation}\label{premotivicadjcomparisonthm04}\begin{split}
\xymatrix{
\DMt \ar[r]^{\beta^*}\ar[d]_{\rho_\sharp}&\DMB\ar[d]^{\chi_\sharp} \\
\uDMt \ar[r]^{\underline{a}^*}&\uDMV_{t,\QQ}
}\end{split}
\end{equation}
\end{num}

From now on, we shall fix an excellent noetherian scheme of finite dimension $S$.

\begin{thm}\label{plongement01}
We have canonical equivalences of categories
$$\DMB(S)\simeq \DMV_{\qfh,\QQ}(S)\simeq\DMV_{\h,\QQ}(S)$$
(recall that, for $t=\qfh,\h$, $\DMV_{t,\QQ}(S)$ stands for
the localizing subcategory of $\uDMV_{t,\QQ}(S)$, spanned by the objects of shape
$\Sigma^\infty \QQ_S(X)(n)$, where $X$ runs over the family of smooth $S$-schemes,
and $n\leq 0$ is an integer; see \ref{ex:stable_AA^1-derived_categories}).
\end{thm}

\begin{proof}
Let $t$ denote the $\qfh$-topology or the $\h$-topology.
We shall prove that the functor
$$\chi_\sharp:\DMB(S)\To \uDMV_{t,\QQ}(S)$$
is fully faithful, and that its essential image is precisely
$\DMV_{t,\QQ}$. The functor
$$\beta_*:\DMB\To \DMt(S)$$
is fully faithful, so that its composition with its left adjoint
$\beta^*$ is canonically isomorphic to the identity.
In particular, we get isomorphisms of functors:
$$\chi_\sharp\simeq \chi_\sharp\, \beta^* \, \beta_* \simeq
\underline{a}^*\, \rho_\sharp\,  \beta_* \, .$$
The right adjoint of $\underline{a}^*$ is fully faithful, and
its essential image consists of the objects of
$\uDMt (S)$ which satisfy $t$-descent (\ref{caracdescentbyhmtppties}).
On the other hand, the functor $\rho_\sharp$ is
fully faithful, and an object of $\DMt(S)$
satisfies $t$-descent if and only if its image by
$\rho_\sharp$ satisfies $t$-descent (\ref{enlargeddescent}).
By virtue of Theorem \ref{DMB_etale&h-descent},
this implies immediately that $\chi_\sharp$ is fully faithful.
Let $\DMV_{t,\QQ}(S)$ be the localizing
subcategory of $\mathit{DM}_{t,\QQ}(S)$ spanned by the objects of shape
$\Sigma^\infty \QQ(X)(n)$, where $X$ runs over the family of smooth $S$-schemes,
and $n\leq 0$ is an integer (\ref{ex:stable_AA^1-derived_categories}).
We know that $\DMV_{t,\QQ}(S)$
is compactly generated (see \ref{ex:Zar&Nis=bounded},
\ref{prop:compact_DMue} and \ref{spanierwhitehead2}), and that $\chi_\sharp$ is
a fully faithful exact functor which preserves small sums
as well as compact objects from $\DMB(S)$ to $\DMV_{t,\QQ}(S)$.
As, by construction, there exists a generating family of compact objects
of $\DMV_{t,\QQ}(S)$
in the essential image of $\chi_\sharp$, this implies that $\chi_\sharp$
induces an equivalence of triangulated categories
$\DMB(S)\simeq \DMV_{t,\QQ}(S)$ (see \ref{equivgeneratorscompactgentriang}).
\end{proof}

Let us underline the following result which completes
Corollary \ref{cor:carac_KBL-local}:
\begin{thm}\label{caracHbmodbyhdescent}
Let $E$ be an object of $\DMtx S$. The following conditions
 are equivalent:
\begin{enumerate}
\item[(i)] $E$ is a Beilinson motive;
\item[(ii)] $E$ satisfies $\h$-descent;
\item[(iii)] $E$ satisfies $\qfh$-descent;
\end{enumerate}
\end{thm}

\begin{proof}
We already know that condition (i) implies condition (ii)
(second point of Theorem \ref{DMB_etale&h-descent}), and condition (ii) implies
obviously condition (iii). It is thus sufficient to prove that
condition (iii) implies condition (i).
If $E$ satisfies $\qfh$-descent, then $\rho_\sharp(E)$
satisfies $\qfh$-descent in $\uDMV(S,\QQ)$ as well.
The commutativity of \eqref{premotivicadjcomparisonthm03}
implies then that $\rho_\sharp(E)$ belongs to the essential
image of $\underline{\gamma}_*$ (the right adjoint
of $\underline{\gamma}^*$).
As $\rho_\sharp$ is fully faithful,
the commutativity of \eqref{premotivicadjcomparisonthm02}
thus implies that $E$ itself belongs to the essential
image of $\gamma_*$ (the right adjoint to $\gamma^*$).
In particular, $E$ is then a module over the ring spectrum
$\gamma_*(\unit_S)$, which is itself an $\HB$-algebra.
We conclude by Corollary \ref{cor:carac_KBL-local}.
\end{proof}

\begin{thm}\label{comparisonDMBDMV}
If $S$ is excellent and geometrically unibranch,
\index{word}{scheme!geometrically unibranch}
 then the comparison functor
$$\varphi^*:\DMB(S)\To \DM(S)$$
is an equivalence of triangulated monoidal categories.
\index{word}{equivalence!of triangulated monoidal categories}
\end{thm}

\begin{proof}
If $S$ is geometrically unibranch, then we know that the
composed functor
$$\DM(S)\xrightarrow{\psi_\sharp}\uDMV_\QQ(S)
\xrightarrow{\underline{\alpha}^*} \uDMV_{\qfh,\QQ}(S)$$
is fully faithful (\ref{embedDMDMqfh}). The commutative diagram
$$\xymatrix{
\DMB(S) \ar[r]^{\varphi^*}\ar@/_1.3pc/[rr]_{\chi_\sharp}
&\DM(S)\ar[r]^{\underline{\alpha}^*\psi_\sharp}&
\uDMV_{\qfh,\QQ}(S)
}$$
and Theorem \ref{plongement01} imply that $\varphi^*$
is fully faithful. As $\varphi^*$ is exact, preserves small
sums as well as compact objects, and as $\DM(S)$
has a generating family of compact objects in the
essential image of $\varphi^*$, the functor $\varphi^*$ has to be
an equivalence of categories (\ref{equivgeneratorscompactgentriang}).
\end{proof}

\begin{rem}
Some version of the preceding theorem
(the one obtained by replacing $\DMB$ by $\ho(\Mod\HB)$)
was already known in the case where $S$ is the spectrum of a perfect field;
see \cite[theorem 68]{RNOST}.
The proof used de~Jong's resolution of singularities by alterations
and Poincar\'e duality in a crucial way. The proof of the preceding theorem
we gave here relies on proper descent
but does not use any kind of resolution of singularities.
\end{rem}

The preceding theorem allows to give the following description
 of constructible Beilinson motives over geometrically unibranch schemes:
\begin{cor} \label{cor:DMBc&Voevodsky}
For any geometrically unibranch scheme $S$, the functor $\varphi^*$
 induces an equivalence of triangulated monoidal categories:
$$
\DMBc(S) \xrightarrow {\ \sim \ } \DMV_{gm}(S,\QQ)
$$
where the right hand side is the $\QQ$-linear version of the category
 of geometric (Voevodsky) motives
\index{word}{motive!geometric}
  (Definition \ref{df:gm_V_mot}).
\end{cor}
Note that we also applied Proposition \ref{prop:changeofcoef_DM}
 to get this corollary.

We finally point out the following important fact about
Voevodsky's motivic cohomology spectrum $\Hmx S=\gamma_*(\unit_S)$ with rational coefficients:
\begin{cor}\label{compHbeilHVoev}
\begin{enumerate}
\item For any geometrically unibranch excellent scheme $S$, the canonical map
$$
\HBx S \rightarrow \Hmrx S
$$
is an isomorphism of ring spectra.
\item For any morphism $f:T \rightarrow S$ of
 excellent geometrically unibranch schemes,
 the canonical map
$$
 f^*\Hmrx S \rightarrow\Hmrx T
$$
is an isomorphism of ring spectra.
\end{enumerate}
\end{cor}
The second part is the last conjecture of Voevodsky's paper \cite{V_OpenPB}
 with rational coefficients (and geometrically unibranch schemes)
 -- see also Paragraph \ref{num:conj_V}.

\begin{proof}
The first part is a trivial consequence of the previous theorem,
and the second follows from the first, as the Beilinson motivic
cohomology spectrum is stable by pullbacks.
\end{proof}
%%\subsection{Duality and finiteness}

\subsection{Comparison with Morel motives}
\index{word}{motive!Morel|(}
\begin{paragr}\label{morel}
Let $S$ be a scheme.
The permutation isomorphism
\begin{equation}\label{morel1}
\tau : \QQ(1)[1]\otimes^\derL_{\QQ}\QQ(1)[1]\To\QQ(1)[1]\otimes^\derL_{\QQ}\QQ(1)[1]
\end{equation}
satisfies the equation $\tau^2=1$ in $\DMtx S$.
Hence it defines an element $\epsilon$ in $\mathrm{End}_{\DMtx S}(\QQ)$
which also satisfies the relation $\epsilon^2=1$.
We define two projectors
\begin{equation}\label{morel2}
e_{+}=\frac{1-\epsilon}{2}
\quad\text{and}\quad
e_{-}=\frac{1+\epsilon}{2}\, .
\end{equation}
As the triangulated category $\DMtx S$ is pseudo abelian,
we can define two objects by the formul{\ae}:
\begin{equation}\label{morel3}
\QQ_{+}=\mathrm{Im}\, e_{+}
\quad\text{and}\quad
\QQ_{-}=\mathrm{Im}\, e_{-}\, .
\end{equation}
Then for an object $M$ of $\DMtx S$, we set
\begin{equation}\label{morel4}
M_{+}=\QQ_{+}\otimes^\derL_{\QQ}M
\quad\text{and}\quad
M_{-}=\QQ_{-}\otimes^\derL_{\QQ}M\, .
\end{equation}
It is obvious that for any objects $M$ and $N$ of $\DMtx S$,
one has
\begin{equation}\label{morel5}
\Hom_{\DMtx S}(M_{i},N_{j})=0\quad\text{for $i,j\in\{+,-\}$ with $i\neq j$.}
\end{equation}
Denote by $\DMtx S_{+}$ (resp. $\DMtx S_{-}$) the full subcategory
of $\DMtx S$ made of objects which are isomorphic to some $M_{+}$
(resp. some $M_{-}$) for an object $M$ in $\DMtx S$. Then
\eqref{morel5} implies that the direct sum functor $(M_+,M_-)\mapsto M_+\oplus M_-$
induces an equivalence of triangulated categories
\begin{equation}\label{morel6}
(\DMtx S_{+})\times(\DMtx S_{-})\simeq\DMtx S\, .
\end{equation}
We shall call $\DMtx S_{+}$
\index{notat}{DA1SQplus@$\DMtx S_{+}$}
 the \emph{category of Morel motives over $S$}.
\index{word}{motive!Morel}
The aim of this section is to compare this category with $\DMB(S)$ (see Theorem \ref{thm:Morel}).
This will consists essentially of proving that $\QQ_{+}$ is nothing else than
Beilinson's motivic spectrum $\HB$ (which was announced by
Morel in {\cite{SHQ}}).
The main ingredients of the proof are the description of $\DMB(S)$ as full subcategory
of $\DMtx S$, the homotopy $t$-structure on $\DMtx S$,
and Morel's computation of the endomorphism ring
of the motivic sphere spectrum in terms of
Milnor-Witt K-theory \cite{dmtilde2,dmtilde,KthMW,A1brouwer}.
\end{paragr}

\begin{paragr}
For a little while, we shall assume that $S$ is the spectrum of a field $k$.

Recall that the \emph{algebraic Hopf fibration}
\index{word}{fibration!algebraic Hopf}
 is the map
$$\AA^2-\{0\}\To \PP^1\ , \quad (x,y)\mapsto [x,y]\, .$$
This defines, by desuspension, a morphism
$$\eta:\QQ(1)[1]\To \QQ$$
in $\DMtx S$; see \cite[6.2]{dmtilde2} (recall that we identify
$\DMtx S$ with $\SH_\QQ(S)$ and that, under this identification,
$\QQ(1)[1]$ corresponds to $\Sigma^\infty(\GG)$).
\end{paragr}

\begin{lm}\label{epsilonHopf}
We have $\eta=\epsilon\eta$ in $\Hom_{\DMtx S}(\QQ(1)[1],\QQ)$.
\end{lm}

\begin{proof}
See \cite[6.2.3]{dmtilde2}.
\end{proof}

\begin{paragr}
Recall the \emph{homotopy $t$-structure}
\index{word}{tstructure@$t$-structure, homotopy}
 on $\DMtx S$; see \cite[5.2]{dmtilde2}.
To remain close to the conventions of \emph{loc. cit.}, we shall adopt
homological notations, so that, for any object $M$ of $\DMtx S$, we have
the following truncation triangle
$$\tau_{> 0}M\To M \To \tau_{\leq 0}M\To \tau_{> 0}M[1] \, .$$
We whall write $H_0$ for the zeroth homology functor in the sense of this $t$-structure.
This $t$-structure can be described in terms of generators, as in
\cite[definition 2.2.41]{ayoub}: the category $\DMtx S_{\geq 0}$
is the smallest full subcategory of $\DMtx S$ which contains the objects of shape
$\QQ_S(X)(m)[m]$ for $X$ smooth over $S$, $m\in\ZZ$, and which
satisfies the following stability conditions:
\begin{itemize}
\item[(a)] $\DMtx S_{\geq 0}$ is stable under suspension; i.e. for any object $M$
in $\DMtx S_{\geq 0}$, $M[1]$ is in $\DMtx S_{\geq 0}$;
\item[(b)] $\DMtx S_{\geq 0}$ is closed under extensions: for any distinguished triangle
$$M'\To M \To M''\To M'[1]\, ,$$
if $M'$ and $M''$ are in $\DMtx S_{\geq 0}$, so is $M$;
\item[(c)] $\DMtx S_{\geq 0}$ is closed under small sums.
\end{itemize}
With this description, it is easy to see that $\DMtx S_{\geq 0}$ is also closed
under tensor product (because the class of generators has this property).
The category $\DMtx S_{\leq 0}$ is the full subcategory
of $\DMtx S$ which consists of objects $M$ such that
$$\Hom_{\DMtx S}(\QQ_S(X)(m)[m+n],M)\simeq 0$$
for $X/S$ smooth, $m\in\ZZ$, and $n>0$; see \cite[2.1.72]{ayoub}.

Note that the heart of the homotopy $t$-structure is symmetric monoidal, with
tensor product $\otimes^h$ defined by the formula:
$$F\otimes^h G=H_0(F\otimes^\derL_S G)$$
(the unit object is $H_0(\QQ)$).

We shall still write $\eta:H_0(\QQ(1)[1])\To H_0(\QQ)$ for the map
induced by the algebraic Hopf fibration.
\end{paragr}

\begin{prop}\label{Tateexact}
Tensoring by $\QQ(n)[n]$ defines a $t$-exact endofunctor 
\index{word}{functor!texact@$t$-exact endofunctor}
of $\DMtx S$ for any integer $n$.
\end{prop}

\begin{proof}
As tensoring by $\QQ(n)[n]$ is an equivalence of categories,
it is sufficient to prove this for $n\geq 0$.
This is then a particular case of \cite[2.2.51]{ayoub}.
\end{proof}

\begin{prop}\label{somevanishing}
For any smooth $S$-scheme $X$ of dimension $d$, and for any
object $M$ of $\DMtx S$, the map
$$\Hom(\QQ_S(X),M)\To \Hom(\QQ_S(X),M_{\leq n})$$
is an isomorphism for $n>d$.
\end{prop}

\begin{proof}
Using \cite[lemma 5.2.5]{dmtilde2}, it is sufficient
to prove the analog for the homotopy $t$-structure
on $\DMte(S)$, which follows from \cite[lemma 3.3.3]{dmtilde3}.
\end{proof}

\begin{prop}\label{homotopytstructnondeg}
The homotopy $t$-structure is non degenerated.
\index{word}{tstructure@$t$-structure!non degenerated}
 Even better,
for any object $M$ of $\DMtx S$, we have canonical isomorphisms
$$\derL \varinjlim_n \tau_{> n}M\simeq M
\quad\text{and}\quad
\derR \varprojlim_n \tau_{> n}M\simeq 0 \, ,$$
as well as isomorphisms
$$\derL \varinjlim_n \tau_{\leq  n}M\simeq 0
\quad\text{and}\quad
M\simeq\derR \varprojlim_n \tau_{\leq n}M \, .$$
\end{prop}

\begin{proof}
The first assertion is a direct consequence of
propositions \ref{Tateexact} and \ref{somevanishing}
(because the objects of shape $\QQ_S(X)(m)[i]$,
for $X/S$ smooth, and $m,i\in\ZZ$, form a generating family).
As the objects $\QQ_S(X)(m)[m+n]$ are compact in $\DMtx S$,
the category $\DMtx S_{\leq 0}$ is closed under small sums.
As $\DMtx S_{\geq 0}$ is also closed under small sums,
we deduce easily that the truncation functors $\tau_{>0}$ and $\tau_{\leq 0}$
preserve small sums, which implies that the homology functor $H_0$
has the same property. Moreover, if
$$C_0\To \cdots \To C_n\To C_{n+1} \To \cdots$$
is a sequence of maps in $\DMtx S$, then $C=\derL\varinjlim_n C_n$
fits in a distinguished triangle of shape
$$\bigoplus_n C_n \overset{1-s}{\To} \bigoplus_n C_n \To C \To \bigoplus_n C_n[1]\, ,$$
where $s$ is the map induced by the maps $C_n\To C_{n+1}$.
This implies that, for any integer $i$, we have
$$\varinjlim_n H_i(C_n)\simeq H_i(C)$$
(where the colimit is taken in the heart of the homotopy $t$-structure).
As the homotopy $t$-structure is non degenerated,
this proves the two formulas
$$\derL \varinjlim_n \tau_{> n}M\simeq M
\quad\text{and}\quad
\derL \varinjlim_n \tau_{\leq  n}M\simeq 0 \, .$$

Let $X$ be a smooth $S$-scheme of finite type, and $p$, $q$ be some integer.
To prove that the map
$$\Hom(\QQ_S(X)(m)[i],M)\To \Hom(\QQ_S(X)(m)[i],\derR \varprojlim_n \tau_{\leq n}M)$$
is bijective, we may assume that $m=0$ (replacing $M$ by $M(-m)[-m]$
and $i$ by $i-m$, and using Proposition \ref{Tateexact}).
Consider  the Milnor short exact sequence below, with $A=\QQ_S(X)[i]$
(in which the first map is injective, but we will not use it):
$$ {\varprojlim_n}^1 \Hom(A[1],\tau_{\leq n}M)\to
\Hom(A,\derR \varprojlim_n \tau_{\leq n}M)\to
\varprojlim_n\Hom(A,\tau_{\leq n}M)\to 0\, .$$
%% As both $A[1]$ and $A$ belong to $\DMtx S_{\geq i}$, it is isomorphic
%% to the following short exact sequence
%% $$0\To {\varprojlim_n}^1 \Hom(A[1],\tau_{\geq i}\tau_{\leq n}M)\To
%% \Hom(A,\tau_{\geq i} \derR \varprojlim_n \tau_{\leq n}M)\To
%% \varprojlim_n\Hom(A,\tau_{\geq i}\tau_{\leq n}M)\To 0\, .$$
Using Proposition \ref{somevanishing}, as $\varprojlim^1$ of a constant functor vanishes,
we get that the map
$$\Hom(A,M)\To \Hom(A,\derR \varprojlim_n \tau_{\leq n}M)$$
is an isomorphism. This gives the isomorphism
$$M\simeq\derR \varprojlim_n \tau_{\leq n}M \, .$$
Using the previous isomorphism, and
by contemplating the homotopy limit of the homotopy cofiber sequences
$$\tau_{> n}M\To M \To \tau_{\leq n}M\, ,$$
we deduce the isomorphism $\derR \varprojlim_n \tau_{>  n}M\simeq 0$.
\end{proof}

\begin{lm}\label{Motcohleqzero}
We have $\HB\in \DMtx S_{\geq 0}$, so that we have a canonical map
$$\HB \To H_0(\HB)$$
in $\DMtx S$. In particular,
 for any object $M$ in the heart of the homotopy $t$-structure,
\index{word}{tstructure@$t$-structure!heart}
 if $M$ is endowed with an action of the
monoid $H_0(\HB)$, then $M$ has a natural structure of $\HB$-module
in $\DMtx S$.
\end{lm}

\begin{proof}
As $\HB$ is isomorphic to the motivic cohomology spectrum
in the sense of Voevodsky (\ref{compHbeilHVoev}),
the first assertion is the first assertion of \cite[theorem 5.3.2]{dmtilde2}.
Therefore, the truncation triangle for $\HB$ gives a triangle
$$\tau_{> 0}\HB\To \HB \To H_0(\HB)\To \tau_{> 0}\HB[1]\, ,$$
which gives the second assertion. For the third assertion,
consider an object $M$ in the heart of the homotopy $t$-structure,
endowed with an action of $H_0(\HB)$. Note that
$\DMtx S_{\geq 0}$ is closed under tensor product, so that
$\HB\otimes^\derL_S M$ is in $\DMtx S_{\geq 0}$.
Hence we have natural maps
$$\HB\otimes^\derL_S M\To H_0(\HB\otimes^\derL_S M)
%%\To H_0(\HB)\otimes^\derL_S M
\To H_0(H_0(\HB)\otimes^\derL_S M)
= H_0(\HB)\otimes^hM\, .$$
Then the structural map $H_0(\HB)\otimes^hM\To M$ defines a map
$\HB\otimes^\derL_S M\To M$ which gives the expected action (observe
that, as we already know that $\HB$-modules do form a thick subcategory
of $\DMtx S$ (\ref{HBmodulesfullyfaithful}),
we don't even need to check all the axioms of an internal module:
it is sufficient to check that the unit $\QQ\To \HB$ induces
a section $M\To \HB\otimes^\derL_S M$ of the map constructed above).
\end{proof}

\begin{lm}\label{MilnorWittexactsequence}
We have the following exact sequence in the heart of the homotopy $t$-structure.
$$H_0(\QQ(1)[1])\xrightarrow{\eta} H_0(\QQ)\To H_0(\HB)\To 0$$
\end{lm}

\begin{proof}
Using the equivalence of categories from the heart of the homotopy $t$-structure
to the category of homotopy modules in the sense of \cite[definition 5.2.4]{dmtilde2},
by virtue of Corollary \ref{compHbeilHVoev} and \cite[theorem 5.3.2]{dmtilde2},
we know that $H_0(\HB)$ corresponds to the homotopy module $\underline{K}^M_*\otimes \QQ$
associated with Milnor K-theory\index{word}{K-theory!Milnor},
while $H_0(\QQ)$ corresponds to the homotopy module $\underline{K}^{MW}_*\otimes \QQ$
associated with Milnor-Witt K-theory
\index{word}{K-theory!Milnor-Witt}
 (which follows easily from \cite[theorems 2.11, 6.13 and 6.40]{A1brouwer}).
Considering $\underline{K}^M_*$ and $\underline{K}^{MW}_*$
 as unramified sheaves
in the sense of Morel~\cite{A1brouwer}, this lemma is then a reformulation
of the isomorphism
$$K^{MW}_*(F)/\eta\simeq K^M_*(F)$$
for any field $F$; see \cite[remark 2.2]{A1brouwer}.
\end{proof}

\begin{prop}\label{QplusHB}
We have ${\HB}_+\simeq\HB$, and the induced map $\QQ_+\To \HB$
gives a canonical isomorphism $H_0(\QQ_+)\simeq H_0(\HB)$.
\end{prop}

\begin{proof}
The map $\epsilon(1)[1]: \QQ(1)[1]\To \QQ(1)[1]$
can be described geometrically as the morphism
associated with the pointed morphism
$$\imath:\GG\To \GG\ , \quad t\mapsto t^{-1}$$
(see the second assertion of \cite[lemma 6.1.1]{dmtilde2}).
In the decomposition
$$K_1(\GG)\simeq k[t,t^{-1}]^\times\simeq k^\times\oplus \ZZ\, ,$$
the map $\imath$ induces multiplication by $-1$ on $\ZZ$.
Using the periodicity isomorphism $\BGL(1)[2]\simeq \BGL$,
we get the identifications:
$$K_1(\GG)\supset\Hom_{\SH(k)}(\Sigma^\infty({\GG})[1],\BGL)\simeq
\Hom_{\BGL}(\BGL,\BGL)\simeq K_0(k)\simeq \ZZ\, .$$
Therefore, $\epsilon$ acts as the multiplication by $-1$ on the spectrum
$\BGL_\QQ$, whence on $\HB$ as well. This means precisely that ${\HB}_+\simeq\HB$.
By Lemma \ref{epsilonHopf}, the class $2\eta$ vanishes in $\QQ_+$,
so that, appyling the ($t$-exact) functor $M\mapsto M_+$
to the exact sequence of Lemma \ref{MilnorWittexactsequence},
we get an isomorphism $H_0(\QQ_+)\simeq H_0({\HB}_+)\simeq H_0(\HB)$.
\end{proof}

\begin{cor}\label{cor:QplusHB}
For any object $M$ in the heart of the homotopy $t$-structure,
$M_+$ is a Beilinson motive.
\index{word}{motive!Beilinson}
\end{cor}

\begin{proof}
The object $M$ is an $H_0(\QQ)$-module, so that
$M_+$ is an $H_0(\QQ_+)$-module. By virtue of
Proposition \ref{QplusHB}, $M_+$ is then a module over $H_0(\HB)$,
so that, by Lemma \ref{Motcohleqzero}, $M_+$ is naturally
endowed with an action of $\HB$.
\end{proof}

\begin{rem}
Until now, we did not really use the fact we are in a $\QQ$-linear context
(replacing $\HB$ by Voevodsky's motivic spectrum,
we just needed $2$ to be invertible in the preceding corollary).
However, the following result really uses $\QQ$-linearity
(because, in the proof, we see $\DMB(S)$ as a full subcategory
of $\DMtx S$; see Proposition \ref{beilinsoncmf}).
\end{rem}

\begin{thm}\label{thm:Morel}
For any noetherian scheme of finite dimension $S$,
the map $\QQ_+\To \HB$ is an isomorphism in $\DMtx S$.
As a consequence, we have a canonical equivalence of
 triangulated monoidal categories
$$\DMtx S_+\simeq\DMB(S)\, .$$
\end{thm}
This theorem has already been proved by Morel when $S$ is the spectrum
 of a perfect field -- where the left hand side is the rational
 category of Voevodsky motives. Morel announced that the category
 $\DMtx S_+$ should be \emph{the} category of rational motives
 and this theorem confirm his insight.
\begin{proof}
Observe that, if ever $\QQ_+\simeq \HB$, we have
$\DMtx S_+\simeq\DMB(S)$: this follows from the fact that an object $M$ of $\DMtx S$
belongs to $\DMtx S_+$ (resp. to $\DMB(S)$) if and only if there exists an isomorphism
$M\simeq M_+$ (resp. $M\simeq \HB\otimes^\derL_S M$; see \ref{cor:carac_KBL-local}).
It is thus sufficient to prove the first assertion.

As both $\QQ_+$ and $\HB$ are stable by pullback, it is sufficient
to treat the case where $S=\spec \ZZ$.
Using Corollary \ref{cor:DMBcontinuous}, we may
replace $S$ by any of its henselisations,
so that, by the localization property, it is sufficient
to treat the case where $S$ is the spectrum of a (perfect) field $k$.

We shall prove directly that, for any object $M$ of $\DMtx S$,
$M_+$ is an $\HB$-module (or, equivalently, is $\HB$-local).
Note that $\DMB(S)$ is closed
under homotopy limits and homotopy colimits in $\DMtx S$:
indeed the inclusion functor $\DMB\To \DMt$
has a left adjoint which preserves a family of compact generators,
whence it also has a left adjoint (\ref{prop:exist_right_adjoint}).
By virtue of Proposition \ref{homotopytstructnondeg}, we may thus
assume that $M$ is bounded with respect to the homotopy $t$-structure.
As $\DMB(S)$ is certainly closed under extensions in $\DMtx S$, we may
even assume that $M$ belongs to the heart the homotopy
$t$-structure. We conclude with Corollary \ref{cor:QplusHB}.
\end{proof}

\begin{cor}\label{cor:Morel}
For any noetherian scheme of finite dimension $S$,
if $-1$ is a sum of squares in all the residue fields of $S$, then
$\QQ_-\simeq 0$ in $\DMtx S$, and we have a canonical equivalence of
 triangulated monoidal categories
$$\DMtx S\simeq \DMB(S)\, .$$
\end{cor}

\begin{proof}
It is sufficient to prove that, under this assumption, $\QQ_-\simeq 0$.
As in the preceding proof, we may replace $S$ by any of its henselisations
(\ref{Nispointscontinuity}),  so that, by the localization property (and
by induction on the dimension), it is sufficient
to treat the case where $S$ is the spectrum of a field $k$.
We have to check that, if $-1$ is a sum of squares in $k$, then
we have $\epsilon=-1$.
Using \cite[remark 6.3.5 and lemma 6.3.7]{dmtilde2},
we see that, if $k$ is of characteristic $2$, we always have $\epsilon=-1$,
while, if the characteristic of $k$ is distinct from $2$, we have
a morphism of rings
$$\mathit{GW}(k)\To \Hom_{\DMt({\spec{k}})}(\QQ,\QQ) \, ,$$
where  $\mathit{GW}(k)$ denotes 
 the Grothendieck-Witt ring\footnote{\emph{i.e.}
 the Grothendieck group of quadratic forms}
\index{word}{ring!Grothendieck-Witt}
 over $k$. This morphism sends
 the class of the quadratic form $-X^2$ to $-\epsilon$ and this proves the result.
(For a more precise version of this, with integral coefficients, see \cite[proposition 2.13]{A1brouwer}.)
\end{proof}

\begin{num}
Recall from Example \ref{ex:A^1-Nis_Et_compacity_stable}
that we can describe the category $\DMtgmx S$ of compact objects of
 $\DMtx S$ as the triangulated monoidal category
 obtained from
$$
\left(\K^b\left(\Rc(\sm/S)\right)/(BG_S \cup \cT_{\AA^1_S})\right)^\natural
$$
by formally inverting the Tate twist.
The operation $\epsilon$ still acts on this category
 and the decomposition in $+$ and $-$ part of a motive
  respects constructibility as this is a decomposition by direct
  factors.
  The preceding theorem gives the following description of
  constructible Beilinson motives:
\end{num}
\begin{cor}  \label{cor:DMBc&Morel}
For any noetherian scheme of finite dimension $S$,
 there is a canonical equivalence of triangulated
 monoidal categories
$$
\DMBc(S) \simeq \DMtgmx S_+
$$
When $-1$ is a sum of square in all the residue fields
 of $S$, this equivalence can be written:
$$
\DMBc(S) \simeq \DMtgmx S.
$$
\end{cor}

\begin{paragr}
Consider the $\QQ$-linear \emph{\'etale motivic category}
$\Der_{\AA^1,\et}(-,\QQ)$, defined by
$$\Der_{\AA^1,\et}(S,\QQ)=\Der_{\AA^1}(\sh{\et}{\sm/S})$$
(see \ref{ex:stable_AA^1-derived_categories}).
The \'etale sheafification functor induces a morphism of
motivic categories
\begin{equation}\label{A1NistoA1etalefunct}
\DMtx S\To \Der_{\AA^1,\et}(S,\QQ)\, .
\end{equation}
We shall prove the following result,
 as an application of Theorem \ref{thm:Morel}.
\end{paragr}

\begin{thm}\label{etalemotives}
For any noetherian scheme of finite dimension $S$, there is a
canonical equivalence of categories
$$\DMB(S)\simeq \Der_{\AA^1,\et}(S,\QQ)\, .$$
\end{thm}
As for Theorem \ref{thm:Morel},
 the idea of this result is from F.~Morel who already proved it
 at least in the case of a base field.

In order prove the above Theorem, we shall study
the behaviour of the decomposition \eqref{morel3} in $\Der_{\AA^1,\et}(S,\QQ)$:

\begin{lm}\label{Qminustrivial}
We have $\QQ_-\simeq 0$ in $\Der_{\AA^1,\et}(S,\QQ)$.
\end{lm}

\begin{proof}
Proceeding as in the proof of Theorem \ref{thm:Morel}, we may assume that
$S$ is the spectrum of a perfect field $k$. By \'etale descent, we see that we may
replace $k$ by any of its finite extension. In particular, we may
assume that $-1$ is a sum of squares in $k$.
But then, by virtue of Corollary \ref{cor:Morel},
$\QQ_-\simeq 0$ in $\DMtx S$, so that, by functoriality, $\QQ_-\simeq 0$
in $\Der_{\AA^1,\et}(S,\QQ)$.
\end{proof}

\begin{proof}[Proof of Theorem \ref{etalemotives}]
Note that the functor \eqref{A1NistoA1etalefunct}
has a fully faithful right adjoint, whose essential image consists of
objects of $\DMtx S$ which satisfy \'etale descent.
As any Beilinson motive satisfies \'etale descent
 (first point of \ref{DMB_etale&h-descent}),
$\DMB(S)$ can be seen naturally as a full subcategory of
$\Der_{\AA^1,\et}(S,\QQ)$. On the other hand,
by virtue of the preceding lemma,
any object of $\DMtx S$ which satisfies \'etale descent
belongs to $\DMtx S_+$. Hence, by Theorem \ref{thm:Morel},
any object of $\DMtx S$ which satisfies \'etale descent
is a Beilinson motive. This achieves the proof.
\end{proof}

\begin{rem}
If $S$ is excellent, and if all the residue fields of $S$ are of characteristic zero,
one can prove Theorem \ref{etalemotives} independently
of Morel's theorem: this follows then directly from a descent argument,
namely from Corollary \ref{charseparated4}
and from Theorem \ref{caracHbmodbyhdescent}.
\end{rem}

\begin{cor}\label{DMetaleKth}
For any regular noetherian scheme of finite dimension $S$, we have
canonical isomorphisms
$$\Hom_{\Der_{\AA^1,\et}(S,\QQ)}(\QQ_S,\QQ_S(p)[q])\simeq
\mathit{Gr}^p_\gamma K_{2p-q}(S)_\QQ\, .$$
\end{cor}

\begin{proof}
This follows immediately from Theorem \ref{etalemotives}, by definition of $\DMB$ (\ref{computeHBinDMB}).
\end{proof}

\begin{cor}\label{compDMtetaleDMB}
For any geometrically unibranch excellent noetherian scheme of finite dimension $S$,
there is a canonical equivalence of symmetric monoidal triangulated categories
$$\Der_{\AA^1,\et}(S,\QQ)\simeq \DMx S\, .$$
\end{cor}

\begin{proof}
This follows from theorems \ref{comparisonDMBDMV} and \ref{etalemotives}.
\end{proof}

\begin{rem}
The preceding corollary extends immediately to the case
 of coefficients in a $\QQ$-algebra $R$
  (cf. Example \ref{num:stable-A^1-derived_chg_coef} for the left hand side
   and Paragraph \ref{num:Beilinson_mot_coef} for the right hand side).
\end{rem}

\begin{cor}\label{etaleversush}
Let $S$ be an excellent noetherian scheme of finite dimension.
An object of $\DMtx S$ satisfies $\h$-descent
\index{word}{descent!hdescent@$\h$-descent}
 if and only if it satisfies \'etale descent.
\index{word}{descent!etale@\'etale}
\end{cor}

\begin{proof}
This follows from theorems \ref{caracHbmodbyhdescent} and \ref{etalemotives}.
\end{proof}
\index{word}{motive!Morel|)}
\section{Realizations} \label{sec:Weil}

\subsection{Tilting}

\begin{paragr}
Let $\M$ be a stable perfect symmetric monoidal $\sm$-fibred combinatorial model category
over an adequate category of $\base$-schemes $\sch$, such that $\ho(\M)$ is motivic,
with generating set of twists $\tau$.

Consider a homotopy cartesian commutative monoid $\mathcal{E}$ in $\M$.
Then $\Mod{\mathcal{E}}$ is an $\sm$-fibred model category, such that $\ho(\Mod{\mathcal{E}})$
is motivic, and we have a morphism of motivic categories
(see \ref{abstractmotivicmodules} and \ref{hmtlinearproperties})
$$\ho(\M)\To \ho(\Mod{\mathcal{E}})\ , \quad M\mapsto \mathcal{E}\otimes^\derL M\, .$$
In practice, all the realization functors are obtained in this way (at least over fields),
which can be formulated as follows (for simplicity, we shall work here in a $\QQ$-linear context,
but, if we are ready to consider higher categorical constructions, there is no
reason to make such an assumption).
\end{paragr}

\begin{paragr}\label{towardtiltings}
Consider a quasi-excellent noetherian scheme $S$ of finite dimension,
 as well as
  two stable symmetric monoidal $\sm$-fibred combinatorial model
 categories $\M$
and $\M'$ over the category of $S$-schemes of finite type such that $\ho(\M)$
and $\ho(\M')$ are motivic (as triangulated premotivic categories).
We also assume that both $\ho(\M)$
and $\ho(\M')$ are $\QQ$-linear and separated.
%%(in particular, $\ho(\M)$ and $\ho(\M')$ are pure (\ref{purityoverfields}),
%%and are stable under the six operations of Grothendieck (\ref{grothendieck6op})).

Consider a Quillen adjunction
\begin{equation}\label{towardtiltings1}
\varphi^*:\M\rightleftarrows\M':\varphi_*\, ,
\end{equation}
inducing a morphism of $\sm$-fibred categories
\begin{equation}\label{towardtiltings2}
\derL \varphi^*:\ho(\M)\To\ho(\M')\, .
\end{equation}
We consider both $\ho(\M)$ and $\ho(\M')$ as endowed with their
Tate twists, which defines two motivic subcategories of
constructible objects $\ho(\M)_c$ and $\ho(\M')_c$, respectively.
The functor $\derL\varphi^*$ preserves constructible
objects, and thus defines a morphism of premotivic categories
\begin{equation}\label{towardtiltings3}
\derL \varphi^*:\ho(\M)_c\To\ho(\M')_c\, .
\end{equation}
\end{paragr}

\begin{prop}\label{realizationoverfields1}
Under the assumptions of \ref{towardtiltings}, if, for any
regular $S$-scheme of finite type $X$, and for any integers $p$ and $q$,
the map
$$\Hom_{\ho(\M)(X)}(\unit_X,\unit_X(p)[q])\To
\Hom_{\ho(\M')(X)}(\unit_X,\unit_X(p)[q])$$
is bijective, then the morphism \eqref{towardtiltings3}
is an equivalence of premotivic categories. Moreover, if both
$\ho(\M)$ and $\ho(\M')$ are compactly generated by their Tate twists,
\index{word}{generated!compactly generated!triangulated $\Pmor$-fibred}
 then the morphism \eqref{towardtiltings2}
  is an equivalence of motivic categories.
\index{word}{equivalence!of motivic categories}
\end{prop}

\begin{proof}
Note first that, for any separated $S$-scheme of finite type $X$, 
and for any integers $p$ and $q$,
the map
$$\Hom_{\ho(\M)(X)}(\unit_X,\unit_X(p)[q])\To
\Hom_{\ho(\M')(X)}(\unit_X,\unit_X(p)[q])$$
is bijective. Indeed, it is equivalent to prove that the maps
$$\derR\Gamma(X,\unit_X(p))\to\derR\Gamma(X,\varphi^*(\unit_X)(p))$$
are isomorphisms in the derived category of $\QQ$-vector spaces:
by $\h$-descent (\ref{charseparated3}), and by virtue of Gabber's weak uniformization
Theorem \ref{Gabber1}, it is sufficient to treat the case where $X$ is regular,
which is done by assumption. Let $T$ be an $S$-scheme of finite type.
To prove that the functor
$$\derL\varphi^*:\ho(\M)_c(T)\To\ho(\M')_c(T)$$
is fully faithful, it is sufficient to choose
two small families $\mathfrak{A}$ and $\mathfrak{B}$ of objects of $\ho(\M)(T)$
such that the thick subcategory generated by $\mathfrak{A}$ (by $\mathfrak{B}$,
respectively) contains $\ho(\M)(T)$, and to check that the map
$$\Hom_{\ho(\M)(T)}(A,B)\to\Hom_{\ho(\M')(T)}(\varphi^*(A),\varphi^*(B))$$
are bijective, where $A$ and $B$ run over $\mathfrak{A}$ and $\mathfrak{B}$,
respectively. By virtue of Proposition
\ref{thmfinitnessproper2},
it is thus sufficient to prove that,
for any separated smooth morphism $f:X\to T$, for any projective morphism $g:Y\To T$,
and for any integers $p$ and $q$, the map
$$\Hom_{\ho(\M)}(\derL f_\sharp(\unit_X),\derR g_*(\unit_Y)(p)[q])\To
\Hom_{\ho(\M')}(\derL f_\sharp(\unit_X),\derR g_*(\unit_Y)(p)[q])$$
is an isomorphism. Consider the pullback square
$$\xymatrix{
X\times_T Y\ar[r]^{\mathit{pr}_2}\ar[d]_{\mathit{pr}_1}& Y\ar[d]^g\\
X\ar[r]_f & T
}$$
From Proposition \ref{prop:motivic_adj_6_operations},
 the functor $\varphi^*$ commutes with $f_!$
  when $f$ is a separated morphism of finite type.
 One then easily concludes using this fact and the isomorphisms
  (obtained by adjunction and smooth (or proper) base change)
$$\begin{aligned}
\Hom(\derL f_\sharp(\unit_X),\derR g_*(\unit_Y)(p)[q])
& \simeq \Hom(\unit_X,\derL f^*\, \derR g_*(\unit_X)(p)[q])\\
& \simeq \Hom(\unit_X,\derR \mathit{pr}_{1,*}\, \derL \mathit{pr}^*_2(\unit_X)(p)[q])\\
& \simeq \Hom(\unit_X,\derR \mathit{pr}_{1,*}(\unit_{X\times_T Y})(p)[q])\, ,\\
& \simeq \Hom(\unit_{X\times_T Y},\unit_{X\times_S Y}(p)[q])\, ,
\end{aligned}$$
that \eqref{towardtiltings3} is fully faithful and that
$\ho(\M')_c(T)$ is the thick subcategory generated by the image
by $\derL\varphi^*$ of constructible objects of $\ho(\M)(T)$.
In other words, the functor \eqref{towardtiltings3} is an
equivalence of categories.

If, moreover, both $\ho(\M)$ and $\ho(\M')$ are compactly generated by their Tate twists, then
the sum preserving exact functor
$$\derL\varphi^*:\ho(\M)(T)\To\ho(\M')(T)$$
is an equivalence at the level of compact objects, hence is an
equivalence of categories (\ref{equivgeneratorscompactgentriang}).
\end{proof}

\begin{paragr}\label{tiltingatlast}
Under the assumptions of \ref{towardtiltings}, assume that $\M$ and
$\M'$ are strongly $\QQ$-linear (\ref{def:stronglyQQlinear}), left proper,
tractable, satisfy the monoid axiom, and have cofibrant unit objects.
%%and perfect (\ref{defperfectsymmonidcmf}).
Let $\mathcal{E}'$ be a fibrant resolution of $\unit$ in $\M'({\spec k})$.
By virtue of Theorem \ref{cmfQlinearcomm}, we may assume that $\mathcal{E}'$
is a fibrant and cofibrant commutative monoid in $\M'$. Then $\derR \varphi_*(\unit)=\varphi_*(\mathcal{E}')$
is a commutative monoid in $\M$. Let $\mathcal{E}$ be a cofibrant resolution of $\varphi_*(\mathcal{E}')$
in $\M(\spec k)$. Using Theorem \ref{cmfQlinearcomm}, we may assume that $\mathcal{E}$
is a fibrant and cofibrant commutative monoid, and that the map
$$\mathcal{E}\To \derR \varphi_*(\mathcal{E}')$$
is a morphism of commutative monoids (and a weak equivalence by
construction). We can see $\mathcal{E}$ and $\mathcal{E}'$ as cartesian commutative monoids in
$\M$ and $\M'$ respectively (by considering their pullbacks along
morphisms of finite type $f:X\To \spec k$).
We obtain the essentially commutative diagram of left Quillen functors
below (in which the lower horizontal map is the functor induced by $\varphi^*$
and by the change of scalars functor along the map $\varphi^*(\mathcal{E})\To\mathcal{E}'$):
\begin{equation}\begin{split}
\xymatrix{
\M\ar[r]\ar[d]& \M' \ar[d]\\
{\Mod{\mathcal{E}}}\ar[r] & {\Mod{\mathcal{E}'}}
}\end{split}
\end{equation}
where $\Mod{\mathcal{E}}$ and $\Mod{\mathcal{E}'}$
 are respectively the model premotivic categories
  of $\mathcal E$-modules and $\mathcal E'$-modules
 (see Proposition \ref{cmfPfibredRMod}).
\index{word}{modules!over a homotopy cartesian commutative monoid}
 
Note furthermore that the right hand vertical left Quillen functor
is a Quillen equivalence by construction (identifying $\M'(X)$
with $\unit_X$-modules, and using the fact that the morphism of monoids
$\unit_X\To \mathcal{E}'_X$ is a weak equivalence in $\M'(X)$). 
\end{paragr}

\begin{thm}\label{realizationoverfields2}
Consider the assumptions of \ref{tiltingatlast}, with $S=\spec k$
the spectrum of a field $k$. We suppose furthermore
that one of the following conditions is verified.
\begin{itemize}
\item[(i)] The field $k$ is perfect.
\item[(ii)] The motivic categories $\ho(\M)$ and $\ho(\M')$ are continuous
and semi-separated.
\end{itemize}
Then the morphism
$$\ho(\Mod{\mathcal{E}})_c\To\ho(\Mod{\mathcal{E}'})_c\simeq\ho(\M')_c$$
is an equivalence of motivic categories. Under these identifications, the morphism \eqref{towardtiltings3}
corresponds to the change of scalar functor
$$\ho(\M)_c\To \ho(\M')_c \simeq\ho(\Mod{\mathcal{E}})_c\ , \quad M\mapsto \mathcal{E}\otimes^\derL M\, .$$
If moreover both $\ho(\M)$ and $\ho(\M')$ are
 compactly generated by their Tate twists,
\index{word}{generated!compactly generated!triangulated $\Pmor$-fibred}
 then these identifications extend to non-constructible objects, so that,
in particular, the morphism \eqref{towardtiltings2}
corresponds to the change of scalar functor
$$\ho(\M)\To \ho(\M') \simeq \ho(\Mod{\mathcal{E}})\ , \quad M\mapsto \mathcal{E}\otimes^\derL M\, .$$
\end{thm}
\begin{rem}
This theorem can be thought as (a part of) a \emph{tilting theory}
\index{word}{tilting}
for motivic (homotopy) categories.
Remark that the theorem above readily implies that the morphism of
motivic categories
$$\varphi^*:\ho(\M)_c\to\ho(\M')$$
commutes with the six operations (because the, by virtue of
Theorem \ref{thm:realcommute6op}, the functor $M\mapsto\mathcal{E}\otimes^\derL M$
has this property, as well as the inclusion $\ho(\M')_c\subset\ho(\M')$).
\end{rem}
\begin{proof}
For any regular $k$-scheme of finite type $X$, and for any integers $p$ and $q$,
the map
$$\Hom_{\ho(\M)(X)}(\unit_X,\mathcal{E}_X(p)[q])\To
\Hom_{\ho(\M')(X)}(\unit_X,\mathcal{E}'_X(p)[q])$$
is bijective: this is easy to check whenever $X$ is smooth over $k$,
which proves the assertion under condition (i), while, under condition (ii),
we see immediately from Proposition \ref{insepclosure} that we may assume condition (i).
The first assertion is then a special case of the first assertion of Proposition \ref{realizationoverfields1}.
Similarly, by Proposition \ref{Modcompactgenerators}, the
second assertion follows from the second
assertion of Proposition \ref{realizationoverfields1}.
\end{proof}

\begin{ex}\label{ex:Betti}
Let $\M$ be the stable $\sm$-fibred model category of Tate spectra,
so that $\ho(\M)=\DMt$, and write $\M_\mathcyr{B}$
for the left Bousfield localization of $\M$ by the class of $\HB$-equivalences (see \ref{beilinsoncmf}),
so that $\ho(\M_\mathcyr{B})=\DMB$.

Let $k$ be a field of characteristic zero, endowed with an embedding
$\sigma:k\To \CC$.
Given a complex analytic manifold $X$, let $\M_{an}(X)$
\index{notat}{ManX@$\M_{an}(X)$}
be the category of complexes of sheaves of $\QQ$-vector spaces
on the smooth analytic site of $X$ (i.e. on the category
of smooth analytic $X$-manifolds, endowed with the Grothendieck topology
corresponding to open coverings), endowed with its
local model structure (see \cite[4.4.16]{ayoub2} and \cite{ayoub3}).
We shall write $\M^{\mathit{eff}}_{\mathit{Betti}}(X)$
\index{notat}{MBettieffX@$\M^{\mathit{eff}}_{\mathit{Betti}}(X)$}
 for the stable
left Bousfield localization of $\M_{an}(X)$ by the maps of shape
$\QQ(U\times\mathbf{D}^1)\To\QQ(U)$ for any analytic smooth $X(\CC)$-manifold $U$
(where $\mathbf{D}^1$ denotes the closed unit disc).
We define at last $\M_{\mathit{Betti}}(X)$
\index{notat}{MBettiX@$\M_{\mathit{Betti}}(X)$}
 as the stable model
category of analytic $\QQ(1)[1]$-spectra in $\M^{\mathit{eff}}_{\mathit{Betti}}(X)$,
where $\QQ(1)[1]$ stands for the cokernel of the map $\QQ\To\QQ(\AA^{1,\mathit{an}}-\{0\})$
induced by $1\in\CC$; see \cite[section 1]{ayoub3}.

Given a $k$-scheme of finite type $X$, we shall write
\begin{equation}\label{eq:ex:Betti1}
\Der_{\mathit{Betti}}(X):=\ho(\M_{\mathit{Betti}}(X))
\end{equation}
\index{notat}{DBettiX@$\Der_{\mathit{Betti}}(X)$}
(where the topological space $X(\CC)$ is endowed with its canonical
analytic structure). According to \cite[1.8 and 1.10]{ayoub3},
there exists canonical equivalences of categories
\begin{equation}\label{eq:ex:Betti2}
\Der_{\mathit{Betti}}(X)\simeq \ho(\M^{\mathit{eff}}_{\mathit{Betti}}(X))\simeq \Der(X(\CC),\QQ)\, ,
\end{equation}
where $\Der(X(\CC),\QQ)$ stands for the (unbounded) derived category of
the abelian category of sheaves of $\QQ$-vector spaces on the small
site of $X(\CC)$. By virtue of \cite[section 2]{ayoub3}, there exists a symmetric monoidal
left Quillen morphism of monoidal $\sm$-fibred model categories over the category
of $k$-schemes of finite type
\begin{equation}\label{eq:ex:Betti3}
\mathit{An}^* : \M\To \M_{\mathit{Betti}}\, ,
\end{equation}
which induces a morphism of motivic categories over the category of $k$-schemes
of finite type. Hence $\derR\mathit{An}_*(\unit)$ is a ring spectrum in $\DMt(\spec k)$
which represents Betti cohomology
\index{word}{cohomology!Betti}
 of smooth $k$-schemes.
As $\Der_{\mathit{Betti}}$ satisfies \'etale descent, it follows
from Corollary \ref{charseparated4} that it satisfies $\h$-descent, from
which, by virtue of Theorem \ref{caracHbmodbyhdescent},
the morphism \eqref{eq:ex:Betti3} defines a left Quillen functor
\begin{equation}\label{eq:ex:Betti4}
\mathit{An}^*:\M_{\mathcyr{B}}\To \M_{\mathit{Betti}}\, ,
\end{equation}
hence gives rise to a morphism of motivic categories
\begin{equation}\label{eq:ex:Betti5}
\DMB\To \Der_{\mathit{Betti}} \, ,
\end{equation}
the \emph{Betti realization functor}
\index{word}{realization functor!Betti}
 of Beilinson motives.

Appyling Theorem \ref{realizationoverfields2} to \eqref{eq:ex:Betti4},
we obtain a commutative ring spectrum $\mathcal{E}_{\mathit{Betti}}=\derR\mathit{An}_*(\unit)$
which represents Betti cohomology of smooth $k$-schemes, such that the
restriction of the functor \eqref{eq:ex:Betti5} to constructible objects
corresponds to the change of scalars functors $M\mapsto {\mathcal{E}_{\mathit{Betti}}}\otimes^\derL M$:
\begin{equation}\label{eq:ex:Betti6}
{\DMB}_{,c}(X)\To \ho(\Mod{\mathcal{E}_{\mathit{Betti}}})_c(X)\simeq \Der^b_c(X(\CC),\QQ)\ .
\end{equation}
It should be pointed out that, here, $\Der^b_c(X(\CC),\QQ)$ means the derived
category of sheaves which are \emph{constructible of geometric origin}
(i.e. constructible in the algebraic sense, and not in the analytic sense).

In other words, once Betti cohomology of smooth $k$-schemes is known, one
can reconstruct canonically the bounded derived categories of constructible sheaves
of geometric origin on $X(\CC)$ for
any $k$-scheme of finite type $X$, from the theory of mixed motives.
We expect all the realization functors to be of this shape (which should
follow from (some variant of) Theorem \ref{realizationoverfields2}): the (absolute) cohomology
of smooth $k$-schemes with constant coefficients determines the derived
categories of constructible sheaves of geometric origin over any $k$-scheme of finite type.
For instance, the geometric part of the theory of variations of mixed Hodge structures should be obtained
from Deligne cohomology, seen as a ring spectrum in $\DMB(k)$
(or, more precisely, in $\M_{\mathcyr{B}}(k)$). Work in progress of Brad Drew \cite{drew1,drew2}
goes in this direction.
\end{ex}

\subsection{Mixed Weil cohomologies}

Let $S$ be an excellent (regular) noetherian scheme of finite dimension,
and $\mathbf{K}$ a field of characteristic zero, called the \emph{field of
coefficients}.

\begin{paragr} \label{introducing_substratum}
Let $E$ be a Nisnevich sheaf of commutative differential graded $\mathbf{K}$-algebras
(i.e. is a commutative monoid in the category of
sheaves of complexes of $\mathbf{K}$-vector spaces).
We shall write
$$H^n(X,E)=\Hom_{\DMte(X)}(\QQ_X,E[n])$$
for any smooth $S$-scheme of finite type $X$, and for any integer $n$
(note that, if $E$ satisfies Nisnevich descent and is $\AA^1$-homotopy
invariant, which we can always assume, using \ref{cmfQlinearcomm},
then $H^n(X,E)=H^n(E(X))$).

We introduce %for the cohomology theory $H^*(\, .\, ,E)$
the following axioms~:
\begin{enumerate}
\item[W1] \textit{Dimension}.--- $H^{i}(S,E)\simeq
\begin{cases}
\mathbf{K}&\text{if $i=0$,}\\
0&\text{otherwise.}
\end{cases}$
\item[W2] \textit{Stability}.---
$\mathrm{dim}_{{\mathbf{K}}} H^i(\GG,E)=
\begin{cases}
1&\text{if $i=0$ or $i=1$,}\\
0&\text{otherwise.}
\end{cases}$
\item[W3] \textit{K\"unneth formula}.---
For any smooth $S$-schemes $X$ and $Y$, the exterior cup product
induces an isomorphism
$$
\bigoplus_{p+q=n}H^p(X,E) \otimes_\mathbf{K} H^q(Y,E)
 \overset{\sim}{\To} H^{n}(X\times_{S}Y,E) \ .
$$
\item[W3$^\prime$] \textit{Weak K\"unneth formula}.---
For any smooth $S$-scheme $X$, the exterior cup product
induces an isomorphism
$$
\bigoplus_{p+q=n}H^p(X,E) \otimes_\mathbf{K} H^q(\GG,E)
 \overset{\sim}{\To} H^{n}(X\times_{S}\GG,E) \ .
$$
\end{enumerate}
\end{paragr}

\begin{paragr}\label{defstabstructure}
Under assumptions W1 and W2, we will call any non-zero element
$c \in H^1(\GG,E)$ a \textit{stability class}. % for $H^*(\, . \, ,E)$.
Note that such a class corresponds to a non-trivial map
$$c:\QQ_S(1)\To E$$
in $\DMte(S)$ (using the decomposition $\QQ(\GG)=\QQ\oplus\QQ(1)[1]$).
In particular, possibly after replacing $E$ by
a fibrant resolution (so that $E$ is homotopy invariant and satisfies Nisnevich
descent), such a stability class
can be lifted to an actual map of complexes of presheaves.
Such a lift will be called a \emph{stability structure} on $E$.
\end{paragr}

\begin{df}
A sheaf of commutative differential graded $\mathbf{K}$-algebras $E$
as above is a \emph{mixed Weil cohomology}
\index{word}{cohomology!mixed Weil}
 (resp. a \emph{stable cohomology})
\index{word}{cohomology!stable}
if it satisfies the properties W1, W2 and W3 (resp. W1, W2 and W3$^\prime$)
stated above.
\end{df}

\begin{prop}\label{representabilityWeil}
Let $E$ be a stable cohomology.
There exists a (commutative) ring spectrum
\index{word}{spectrum!ring ---- (associated with a stable cohomology)}
 $\mathcal{E}$ in $\DMB(S)$
with the following properties.
\begin{itemize}
\item[(i)] For any smooth $S$-scheme $X$, and any integer $i$,
there is a canonical isomorphism of $\mathbf{K}$-vector spaces
$$H^i(X,E)\simeq \Hom_{\DMB(S)}(M_S(X),\mathcal{E}[i])\, .$$
\item[(ii)] Any choice of a stability structure on $E$ defines a map
$\QQ(1)\To\mathcal{E}$ in $\DMB(S)$, which induces an
$\mathcal{E}$-linear isomorphism
$\mathcal{E}(1)\simeq\mathcal{E}$.
\end{itemize}
\end{prop}

\begin{proof}
One defines explicitly the commutative ring spectrum $\mathcal{E}$
as follows. First, by virtue of Theorem \ref{cmfQlinearcomm},
we may assume that $E$ is a Nisnevich sheaf of commutative
differential graded algebras and is fibrant for the $\AA^1$-local projective model structure:
for any smooth $S$-scheme $X$, the two maps
$$H^n(E(X))\To H^n_\nis(X,E)\To H^n_\nis(X\times\AA^1,E)$$
are isomorphisms for any $n\in \ZZ$.
Let $L$ be the constant Nisnevich sheaf of complexes
of $\mathbf{K}$-vector spaces associated
to the kernel of the map induced by $S=\{1\}\subset\mathbf{G}_m$:
$$L=\mathrm{ker}\big(E(\mathbf{G}_m)\xrightarrow{ 1^* } E(S)\big)\, .$$
We remark that $L$ is cofibrant, and one defines
$$\mathcal{E}_n=\uHom(L^{\otimes n},E)$$
this sheaf being endowed with an action of the symmetric
group on $n$ letters by permuting the factors on $L^{\otimes n}$.
We then have canonical pairings
$$\uHom(L^{\otimes m},E)\otimes_\QQ\uHom(L^{\otimes n},E)
\To\uHom(L^{\otimes m+n},E\otimes_\QQ E)\To \uHom(L^{\otimes m+n},E)$$
which turn the collection $\mathcal{E}=\{\mathcal{E}_n\}_{n\geq 0}$ into a commutative
monoid in the category of symmetric sequences of sheaves of complexes
of $\QQ$-vector spaces; see Definition \ref{df:symmtensorproduct}.
On the other hand, we remark that $L$ is the constant sheaf
associated to $\Gamma(S,\uHom(\QQ(1)[1],E))$, from which we deduce that
there is a natural map
$$L\To \uHom(\QQ(1)[1],E)$$
which can be transposed into a canonical map
$$\QQ(1)[1]\To \uHom(L,E)=\mathcal{E}_1\, .$$
This defines a canonical structure of commutative monoid in
the category symmetric $\QQ(1)[1]$-spectra on the symmetric
sequence $\mathcal{E}$ (see Remark
\ref{rem:caractcommringspectra})\footnote{Here, we work with
$\QQ(1)[1]$-spectra. However, the paper \cite{CD2} is written
in the language of symmetric $\QQ(1)$-spectra. We leave as an exercise
to the reader the task of the translation, which consists in checking that the functor
$\{\mathcal{E}_n\}_{n\geq 0}\mapsto\{\mathcal{E}_n[n]\}_{n\geq 0}$
is a symmetric monoidal left Quillen equivalence from symmetric $\QQ(1)[1]$-spectra
to symmetric $\QQ(1)$-spectra, which is also a right Quillen functor (and thus, in
particular, preserves and detects stable $\AA^1$-equivalences).}.

By virtue of \cite[Proposition 2.1.6]{CD2},
for any smooth $S$-scheme $X$, and any integer $i$,
there is a canonical isomorphism of $\mathbf{K}$-vector spaces
$$H^i(X,E)\simeq \Hom_{\DMt(S)}(M_S(X),\mathcal{E}[i])\, ,$$
and any choice of a stability structure on $E$ defines an isomorphism
$\mathcal{E}(1)\simeq\mathcal{E}$.
Moreover, \cite[corollary 2.2.8]{CD2} and Theorem \ref{equivorientations}
assert that this ring spectrum $\mathcal{E}$ is oriented, so that,
by Corollary \ref{cor:carac_KBL-local},
$\mathcal{E}$ is an $\HB$-module, i.e. belongs to $\DMB(S)$.
\end{proof}

\begin{paragr}\label{num:df:E-modules}
Given a stable cohomology $E$ and its associated ring spectrum $\mathcal{E}$,
we can see $\mathcal{E}$ as a cartesian commutative monoid:
we define, for an $S$-scheme $X$, with structural map $f:X\To S$:
$$\mathcal{E}_X=\derL f^*(\mathcal{E})$$
(which means that we take a cofibrant replacement $\mathcal{E}'$ of $\mathcal{E}$
in the model category of commutative monoids of the category
of Tate spectra, and define $\mathcal{E}_X=f^*(\mathcal{E}')$), and put
\begin{equation}\label{defEmodWeil1}
\Der(X,\mathcal{E}):=\ho(\Mod{\mathcal{E}})(X)=\ho(\Mod{\mathcal{E}_X})\, .
\end{equation}
\index{notat}{DXE@$\Der(X,\mathcal{E})$}
We thus have realization functors
\index{word}{realization functor!(associated with a stable cohomology)}
\begin{equation}\label{defEmodWeil2}
\DMB(X)\To \Der(X,\mathcal{E})\ , \quad M\mapsto \mathcal{E}_X\otimes^\derL_X M
\end{equation}
which commute with the six operations of Grothendieck if ever $S$
is the spectrum of a field (Theorem \ref{thm:realcommute6op}).
Furthermore, $\Der(-,\mathcal{E})$ is a motivic category which is $\QQ$-linear
(in fact $\mathbf{K}$-linear), separated, and continuous.
%Observe furthermore that, if $S$ is the spectrum
%of a field , then $\Der(-,\mathcal{E})$ is also pure (\ref{purityoverfields}),
%so that the six Grothendieck operations preserves constructible
%objects in $\Der(-,\mathcal{E})$ (\ref{grothendieck6op}).

For an $S$-scheme $X$, define
$$H^q(X,E(p))=\Hom_{\DMB(X)}(\QQ_X,\mathcal{E}(p)[q])\simeq
\Hom_{\Der(X,\mathcal{E})}(\mathcal{E}_X,\mathcal{E}_X(p)[q])$$
(this notation is compatible with \ref{introducing_substratum}
by virtue of Proposition \ref{representabilityWeil}).
\end{paragr}

\begin{cor}\label{cor:descentWeil}
Any stable cohomology (in particular, any mixed Weil cohomology)
extends naturally to $S$-schemes of finite type, and this
extension satisfies cohomological $\h$-descent
\index{word}{descent!cohomological $\h$-descent}
 (in particular, \'etale descent
as well as proper descent).
\end{cor}

\begin{proof}
This follows immediately from the construction above and from
Theorem \ref{DMB_etale&h-descent}.
\end{proof}

\begin{paragr}
We denote by $\Der^{\vee}(S,\mathcal{E})$ the localizing subcategory of
$\Der(S,\mathcal{E})$ generated by its rigid objects (i.e. by the objects which have
strong duals). For instance, 
for any smooth and proper $S$-scheme $X$,
$\mathcal{E}(X)=\mathcal{E}\otimes^\derL_S M_S(X)$ belongs to
$\Der^{\vee}(S,\mathcal{E})$; see \ref{prop:purity&duality}.

If we denote by $\Der(\mathbf{K})$ the (unbounded) derived category of the
abelian category of $\mathbf{K}$-vector spaces, we get the following
interpretation of the K\"unneth formula.
\end{paragr}

\begin{thm}\label{smoothproperrealization}
If $E$ is a mixed Weil cohomology, then the functor
$$\derR\Hom_{\mathcal{E}}(\mathcal{E},-):\Der^{\vee}(S,\mathcal{E})
\To \Der(\mathbf{K})$$
is an equivalence of symmetric monoidal triangulated categories.
\end{thm}

\begin{proof}
This is \cite[theorem 2.6.2]{CD2}.
\end{proof}

\begin{thm}\label{ressingWeil}
If $S$ is the spectrum of a field, then $\Der^{\vee}(S,\mathcal{E})=\Der(S,\mathcal{E})$.
\end{thm}

\begin{proof}
This follows then from Corollary \ref{poincarerigidfield}.
\end{proof}

\begin{rem}
It is not reasonable to expect the analog of Theorem \ref{ressingWeil} to hold
whenever $S$ is of dimension $>0$; see (the proof of) \cite[corollary 3.2.7]{CD2}.
Heuristically, for higher dimensional schemes $X$, the rigid objects of $\Der(X,\mathcal{E})$
are extensions of some kind of locally constant sheaves
(in the $\ell$-adic setting, these correspond to \emph{$\QQ_\ell$-faisceaux lisses}).
\end{rem}

\begin{cor}\label{tiltingmixedWeil}
If $E$ is a mixed Weil cohomology, and if $S$ is the spectrum of a field, then the functor
$$\derR\Hom_{\mathcal{E}}(\mathcal{E},-):\Der(S,\mathcal{E})
\To \Der(\mathbf{K})$$
is an equivalence of symmetric monoidal triangulated categories.
\end{cor}

\begin{rem}
This result can be thought as a \emph{tilting theory}
\index{word}{tilting}
 for the spectra associated with mixed Weil cohomologies.
\end{rem}

\begin{paragr}\label{Weiloverfields}
Assume that $E$ is a mixed Weil cohomology, and that $S$ is the spectrum
of a field $k$. For each $k$-scheme of finite type $X$, denote by $\Der_c(X,\mathcal{E})$
the category of constructible objects of $\Der(X,\mathcal{E})$: by definition, this is
the thick triangulated subcategory of $\Der(X,\mathcal{E})$ generated by objects of shape
$\mathcal{E}(Y)=\mathcal{E}\otimes^\derL_XM_X(Y)$ for $Y$ smooth over $X$
(we can drop Tate twists because of \ref{representabilityWeil}~(ii)).
The category $\Der_c(X,\mathcal{E})$ also coincides with the category of compact objects
in $\Der(X,\mathcal{E})$; see \ref{constructequivcompact}.
Write $\Der^b(\mathbf{K})$ for the bounded derived category of
the abelian category of finite dimensional $\mathbf{K}$-vector spaces.
Note that $\Der^b(\mathbf{K})$ is canonically equivalent to the homotopy
category of perfect complexes of $\mathbf{K}$-modules, i.e. to the
category of compact objects of $\Der(\mathbf{K})$.
\end{paragr}

\begin{cor}\label{Weiloverfields2}
Under the assumptions of \ref{Weiloverfields}, we have a canonical equivalence
of symmetric monoidal triangulated categories
$$\Der_c(\spec k,\mathcal{E})\simeq \Der^b(\mathbf{K})\, .$$
\end{cor}

\begin{proof}
This follows from \ref{tiltingmixedWeil} and from the fact that
equivalences of categories preserve compact objects.
\end{proof}

\begin{cor}\label{compWeil}
Under the assumptions of \ref{Weiloverfields}, if $E'$ is another $\mathbf{K}$-linear
stable cohomology with associated ring spectrum $\mathcal{E'}$, any
morphism of presheaves of commutative differential $\mathbf{K}$-algebras
$E\To E'$ inducing an isomorphism $H^1(\GG,E)\simeq H^1(\GG,E')$
gives a canonical isomorphism $\mathcal{E}\simeq\mathcal{E}'$ in the
homotopy category of commutative ring spectra. In particular, we get
canonical equivalences of categories
$$\Der(X,\mathcal{E})\simeq \Der(X,\mathcal{E}')$$
for any $k$-scheme of finite type $X$ (and these are compatible
with the six operations of Grothendieck, as well as with the realization
functors).
\end{cor}

\begin{proof}
This follows from Theorem \ref{ressingWeil} and from
\cite[Theorem 2.6.5]{CD2}.
\end{proof}

The preceding corollary can be stated in the following way:
if $\mathcal{E}$ and $\mathcal{E}'$ are two (strict) commutative
ring spectra associated to $\mathbf{K}$-linear mixed Weil cohomologies
defined on smooth $k$-schemes $E$ and $E'$, respectively, then any
morphism $\mathcal{E}\To \mathcal{E}'$
in the homotopy category of (commutative) monoids
in the model category of $\mathbf{K}$-linear Tate spectra
is invertible. Moreover, $\mathcal{E}$ is isomorphic to $\mathcal{E}'$
if and only if $E$ is isomorphic to $E'$ (in the
appropriate homotopy categories of commutative monoids).
To be more precise (and more general), this last
assertion follows immediately from Corollary \ref{compWeil}
and from the following result.

\begin{prop}\label{prop:comptablevsnonstableWeil}
Let $\mathbf{E}$ be a commutative monoid in the $\AA^1$-stable
model category of sheaves of complexes of symmetric $\QQ(1)[1]$-spectra over
the Nisnevich smooth site of $k$. Suppose that there exists
an isomorphism $\mathbf{E}(1)\simeq\mathbf{E}$ in the homotopy
category of $\mathbf{E}$-modules and that
$$H^n(\spec k,\mathbf{E})=\begin{cases}
\mathbf{K}&\text{if $n=0$,}\\
0&\text{otherwise.}
\end{cases}$$
Then $E=\derR\Gamma(-,\mathbf{E})$
is a stable cohomology theory, and the commutative ring spectrum $\mathcal{E}$
associated to $E$ by Proposition \ref{representabilityWeil}
is canonically isomorphic to $\mathbf{E}$ in the
homotopy category of (strict) commutative ring spectra.
\end{prop}

\begin{proof}
By virtue of Theorem \ref{cmfQlinearcomm},
we may assume that $\mathbf{E}$ is (cofibrant and) fibrant.
The ring spectrum $\mathbf{E}$ is defined by a symmetric sequence of
complexes of Nisnevich sheaves of $\mathbf{K}$-vector spaces $\mathbf{E}_n$, $n\geq 0$,
(endowed with an action of the symmetric group on $n$-letters), together
with maps $\sigma_n:\mathbf{E}_n(1)[1]\To\mathbf{E}_{n+1}$ inducing
quasi-isomorphisms
$$\mathbf{E}_n\xrightarrow{\sim}\uHom(\mathbf{K}(1)[1],\mathbf{E}_{n+1})\, $$
as well as pairings
$$\mathbf{E}_m\otimes_\mathbf{K}\mathbf{E}_n\To\mathbf{E}_{m+n}$$
satisfying a few compatibilities. In particular,
$$E=\derR\Gamma(-,\mathbf{E})=\mathbf{E}_0$$
is naturally endowed with a structure of Nisnevich
sheaf of commutative differential graded algebras which satisfies
Nisnevich descent and $\AA^1$-homotopy invariance.
Moreover, for any integer $n\geq 0$, the Nisnevich
sheaf of complexes of $\mathbf{K}$-vector spaces $\mathbf{E}_n$ also has
the properties of Nisnevich descent and of $\AA^1$-homotopy invariance,
and is naturally endowed with a structure of $E$-module.
It is clear that $E$ is a stable cohomology theory, so that
(the proof of) Proposition \ref{representabilityWeil}
provides a commutative ring spectrum
$\mathcal{E}$ associated to it. With the notations introduced
in the proof of Proposition \ref{representabilityWeil}, we know that
$\mathcal{E}$ is made of the symmetric sequence
$\{\mathcal{E}_n=\uHom(L^{\otimes n},E)\}_{n\geq 0}$, where
$L$ is the constant sheaf associated to $\Gamma(S,\uHom(\mathbf{K}(1)[1],E))$.
Let us define $\mathcal{L}=L(1)[1]$.
We define a new symmetric sequence $\underline{\mathbf{E}}$ by
the formula
$$\underline{\mathbf{E}}_n=\uHom(\mathcal{L}^{\otimes n},\mathbf{E}_n)\, , \quad n\geq 0\, ,$$
where the symmetric group acts through the diagonal
$\mathfrak{S}_n\To\mathfrak{S}_n\times\mathfrak{S}_n$
by permutation of the factors on $\mathcal{L}^{\otimes n}$ and by
the structural action on $\mathbf{E}_n$.
We see that $\underline{\mathbf{E}}$ is a commutative
monoid in the category of symmetric sequences with the pairings
defined by the tensor product map
$$\uHom(\mathcal{L}^{\otimes m},\mathbf{E}_m)
\otimes_\mathbf{K}\uHom(\mathcal{L}^{\otimes n},\mathbf{E}_n)
\To
\uHom(\mathcal{L}^{\otimes m+n},\mathbf{E}_m\otimes_\mathbf{K} \mathbf{E}_n)$$
composed with the multiplication of $\mathbf{E}$:
$$
\uHom(\mathcal{L}^{\otimes m+n},\mathbf{E}_m\otimes_\mathbf{K} \mathbf{E}_n)
\To \uHom(\mathcal{L}^{\otimes m+n},\mathbf{E}_{m+n})\, .$$
Finally, we can compose the transposition of the map $\sigma_1:E(1)[1]\To\mathcal{E}_1$,
with the structural map $\mathbf{K}(1)[1]\To \uHom(L,E)=\mathcal{E}_1$,
to obtain:
$$\mathbf{K}(1)[1]\To\uHom(L,E)\To\uHom(L,\uHom(\mathbf{K}(1)[1],\mathbf{E}_1)
\simeq\uHom(\mathcal{L},\mathbf{E}_1)=\underline{\mathbf{E}}_1\, .$$
This defines a structure of commutative ring spectrum on $\underline{\mathbf{E}}$.
Note that $L$ is chain homotopy equivalent to $\mathbf{K}[-1]$, so that
the functors $\uHom(L^{\otimes n},-)$ preserve quasi-isomorphisms
(more precisely, $L$ is concentrated in cohomological degree $1$, and its first
cohomology sheaf is the constant sheaf associated to the
$\mathbf{K}$-vector space of dimension one $H^1(\mathbf{G}_m,\mathbf{E})$). 
%
%any choice of a
%non zero map $\mathbf{K}(1)\To\mathbf{E}$ in $\DMB(k)$
%induces an isomorphism $\mathbf{E}(1)\simeq\mathbf{E}$
%in the homotopy category of $\mathbf{E}$-modules (because
%the $\mathbf{K}$-vector space $H^1(\mathbf{G}_m,\mathbf{E})$
%is of dimension $1$), which
%can be promoted into an homotopy equivalence (because $\mathbf{E}$
%is cofibrant and fibrant), and such choices define a chain
%homotopy equivalence relating $\mathcal{L}$ and $\mathbf{K}(1)$).
Therefore, one has a quasi-isomorphism of commutative monoids of $\mathbf{K}$-linear
Tate spectra $\mathcal{E}\To\underline{\mathbf{E}}$,
defined by the canonical maps
$$\uHom(L^{\otimes n},E)\To\uHom(L^{\otimes n},\uHom(\mathbf{K}(n)[n],\mathbf{E}_n)\simeq
\uHom(\mathcal{L}^{\otimes n},\mathbf{E}_n)\, .$$
It remains to produce a quasi-isomorphism of commutative monoids of Tate spectra
$\mathbf{E}\To\underline{\mathbf{E}}$. We have a structural map $\mathbf{K}(1)[1]\To\uHom(L,E)$
which can be transposed into a map
$$\mathcal{L}=L(1)[1]\To E=\mathbf{E}_0\, .$$
As $E$ is a commutative monoid and each $\mathbf{E}_n$ an $E$-module, we have
natural maps
$$\mathcal{L}^{\otimes n}\otimes_{\mathbf{K}}\mathbf{E}_n
\To E^{\otimes n}\otimes_{\mathbf{K}}\mathbf{E}_n
\To E\otimes_{\mathbf{K}}\mathbf{E}_n\To \mathbf{E}_n$$
which can be transposed into $\mathfrak{S}_n$-equivariant maps
$$\mathbf{E}_n\To \uHom(\mathcal{L}^{\otimes n},\mathbf{E}_n)=\underline{\mathbf{E}}_n\, .$$
These define a morphism of commutative monoids of $\mathbf{K}$-linear Tate spectra
$\mathbf{E}\To\underline{\mathbf{E}}$. It remains to check that the maps
$\mathbf{E}_n\To\underline{\mathbf{E}}_n$ are quasi-isomorphisms for each $n\geq 0$.
As $\uHom(\mathbf{K}(n)[n],E)\simeq\mathbf{E}_n$, we can replace
$\mathbf{E}_n$ by $\uHom(L^{\otimes n},E)$. The case $n=1$ is then a reformulation
of Proposition \ref{representabilityWeil}~(\emph{ii}), and the general case follows
by an obvious induction.
\end{proof}

\begin{thm}\label{finitenessmixedWeilsixop}
Under the assumptions of paragraph \ref{Weiloverfields},
the six operations of Grothendieck preserve constructibility
in the motivic category $\Der(-,\mathcal{E})$,
 as defined in Paragraph \ref{num:df:E-modules}.
\end{thm}

\begin{proof}
Observe that the motivic category $\Der(-,\mathcal{E})$ is $\QQ$-linear
 and separated (because $\DMB$ is so, see \ref{hmtlinearproperties}),
 as well as $\tau$-compatible
 (because by Proposition \ref{purityoverfields}, it is even $\tau$-dualizable
 which is stronger than $\tau$-compatible; see Definition \ref{deftaupure}).
We conclude with \ref{grothendieck6op}.
\end{proof}

\begin{paragr}
As a consequence, we have, for any $k$-scheme of finite type $X$, a
realization functor
\index{word}{realization functor!of construcible motives}
$$\DMV_{\mathcyr{B},c}(X)\To \Der_c(X,\mathcal{E})$$
and we deduce from Theorem \ref{thm:realcommute6op}
that it preserves all of Grothendieck six operations.
For $X=\spec k$, by virtue of Corollary \ref{Weiloverfields2},
this corresponds to a symmetric monoidal exact realization functor
$$R:\DMV_{\mathcyr{B},c}(\spec k)\To \Der^b(\mathbf{K})\, .$$
This leads to a finiteness result:
\end{paragr}

\begin{cor}\label{Weiloverfields3}
Under the assumptions of \ref{Weiloverfields}, for any
$k$-scheme of finite type $X$, and for any objects
$M$ and $N$ in $\Der_c(X,\mathcal{E})$, $\Hom_{\mathcal{E}}(M,N[n])$
is a finite dimensional $\mathbf{K}$-vector space, and it is trivial
for all but a finite number of values of $n$.
\end{cor}

\begin{proof}
Let $f:X\To \spec k$ be the structural map. By virtue of \ref{finitenessmixedWeilsixop},
as $M$ and $N$ are constructible, the object $\derR f_*\, \derR \sHom_X(M,N)$
is constructible as well, i.e. is a compact object of $\Der(\spec k,\mathcal{E})$.
But $\derR\Hom_{\mathcal{E}}(M,N)$ is nothing else than
the image of $\derR f_*\, \derR \sHom_X(M,N)$
by the equivalence of categories given by Corollary \ref{tiltingmixedWeil}.
Hence $\derR\Hom_{\mathcal{E}}(M,N)$ is a compact object of
$\Der(\mathbf{K})$, which means that it belongs to $ \Der^b(\mathbf{K})$.
\end{proof}

\begin{paragr}
For a $\mathbf{K}$-vector space $V$ and an integer $n$, define
$$V(n)=\begin{cases}
V\otimes_{\mathbf{K}}\Hom_{\mathbf{K}}(H^1(\GG,E)^{\otimes n},\mathbf{K})&\text{if $n> 0$,}\\
V\otimes_{\mathbf{K}}H^1(\GG,E)^{\otimes (-n)}&\text{if $n\leq 0$.}
\end{cases}$$
Any choice of a generator in $\mathbf{K}(-1)=H^1(\GG,E)\simeq H^2(\PP^1_k,E)$
defines a natural isomorphism $V(n)\simeq V$ for any integer $n$.
We have canonical isomorphisms
$$H^q(X,E(p))\simeq H^q(X,E)(p)$$
(using the fact that the equivalence of Corollary \ref{Weiloverfields2}
is monoidal). The realization functors \eqref{defEmodWeil2}
induce in particular cycle class maps
$$\mathit{cl}_X:\HB^q(X,\QQ(p))\To H^q(X,E)(p)$$
(and similarly for cohomology with compact support,
for homology, and for Borel-Moore homology).
\end{paragr}

\begin{ex}\label{ex:deRham}
Let $k$ be a field of characteristic zero. We then have a
mixed Weil cohomology $E_{\mathit{dR}}$ defined by the
algebraic de Rham complex
\index{word}{complex!algebraic De Rham}
$$E_{\mathit{dR}}(X)=\Omega^*_{A/k}$$
for any smooth affine $k$-scheme of finite type $X=\spec A$
(algebraic de Rham cohomology
\index{word}{cohomology!algebraic De Rham}
 of smooth $k$-schemes of finite type
is obtained by Zariski descent); see \cite[3.1.5]{CD2}.
We obtain from \ref{representabilityWeil} a commutative ring spectrum $\mathcal{E}_{\mathit{dR}}$,
and, for a $k$-scheme of finite type $X$, we define
$$\Der_{\mathit{dR}}(X)=\Der_c(X,\mathcal{E}_{\mathit{dR}})\, .$$
We thus get a motivic category $\Der_{\mathit{dR}}$, and
we have a natural definition of algebraic de Rham cohomology of
$k$-schemes of finite type, given by
$$H^n_{\mathit{dR}}(X)=
\Hom_{\Der_{\mathit{dR}}(X)}(\mathcal{E}_{\mathit{dR},X},\mathcal{E}_{\mathit{dR},X}[n])\, .$$
This definition coincides with the usual one: this is true by definition for
separated smooth $k$-schemes of finite type, while the general case follows from
$\h$-descent (\ref{cor:descentWeil}) and from de Jong's Theorem \ref{cor:dejongdimleq2}
(or resolution of singularities \`a la Hironaka). 
We have, by construction, a \emph{de Rham realization functor}
\index{word}{realization functor!de Rham}
$$R_{\mathit{dR}}:\DMV_{\mathcyr{B},c}(X)\To \Der_{\mathit{dR}}(X)$$
which preserves the six operations of Grothendieck (Theorem \ref{thm:realcommute6op}).
In particular, we have cycle class maps
$$\HB^q(X,\QQ(p))\To H^q_{\mathit{dR}}(X)(p)\, .$$
Note that, for any field extension $k'/k$, we have natural isomorphisms
$$H^n_{\mathit{dR}}(X)\otimes_k k'\simeq H^n_{\mathit{dR}}(X\times_{\spec k}\spec{k'})\, .$$
\end{ex}

\begin{ex}\label{ex:deRhamanalytic}
Let $k$ be a field of characteristic zero, which is algebraically closed and
complete with respect to some valuation (archimedian or not).
We can then define a stable cohomology $E_{\mathit{dR,an}}$
as analytic de Rham cohomology
\index{word}{cohomology!analytic De Rham}
 of $X^{\mathit{an}}$,
for any smooth $k$-scheme of finite type $X$; see \cite[3.1.7]{CD2}.
As above, we get a ring spectrum $\mathcal{E}_{\mathit{dR,an}}$,
and for any $k$-scheme of finite type, a category of coefficients
$$\Der_{\mathit{dR,an}}(X)=\Der_c(X,\mathcal{E}_{\mathit{dR,an}})\, ,$$
which allows to define the analytic de Rham cohomology of any $k$-scheme
of finite type $X$ by
$$H^n_{\mathit{dR,an}}(X)=
\Hom_{\Der_{\mathit{dR,an}}(X)}(\mathcal{E}_{\mathit{dR,an},X},\mathcal{E}_{\mathit{dR,an},X}[n])\, .$$
We also have a realization functor
$$R_{\mathit{dR,an}}:\DMV_{\mathcyr{B},c}(X)\To \Der_{\mathit{dR,an}}(X)$$
which preserves the six operations of Grothendieck.

We then have a morphism of stable cohomologies
$$E_{\mathit{dR}} \To E_{\mathit{dR,an}}$$
which happens to be a quasi-isomorphism locally
for the Nisnevich topology (this is Grothendieck's theorem
in the case where $K$ is archimedian, and Kiehl's theorem in the case
where $K$ is non-archimedian; anyway, one obtains this directly
from Corollary \ref{compWeil}). This induces a canonical isomorphism
$$\mathcal{E}_{\mathit{dR}} \simeq \mathcal{E}_{\mathit{dR,an}}$$
in the homotopy category of commutative ring spectra.
In particular, $E_{\mathit{dR,an}}$ is a mixed Weil cohomology, and,
for any $k$-scheme of finite type, we have natural equivalences of categories
$$\Der_{\mathit{dR}}(X) \To \Der_{\mathit{dR,an}}(X) \ , \quad
M\mapsto \mathcal{E}_{\mathit{dR,an}}\otimes^\derL_{\mathcal{E}_{\mathit{dR}}}M$$
which commute with the six operations
of Grothendieck and are compatible with the realization functors.

Note that, in the case $k=\CC$, $\mathcal{E}_{\mathit{dR,an}}$ coincides with
Betti cohomology (after tensorization by $\CC$), so that we have canonical
fully faithful functors
$$\Der_{\mathit{Betti},c}(X)\otimes_\QQ\CC\To \Der_{\mathit{dR,an}}(X)$$
which are compatible with the realization functors.
More precisely, it follows from Proposition \ref{prop:comptablevsnonstableWeil}
that the Betti spectrum $\mathcal{E}_{\mathit{Betti}}$, obtained by appyling
Theorem \ref{realizationoverfields2} to Ayoub's realization functor \eqref{eq:ex:Betti4},
is the spectrum associated to $\QQ$-linear Betti cohomology, seen as a mixed Weil
cohomology, from Proposition \ref{representabilityWeil}. Therefore,
the holomorphic Poincar\'e Lemma, together with Corollary \ref{compWeil},
provide an isomorphism
$$\mathcal{E}_{\mathit{Betti}}\otimes_\QQ\CC\simeq\mathcal{E}_{\mathit{dR,an}}$$
in the homotopy category of commutative monoids of the model category
of $\CC$-linear Tate spectra. We thus have triangulated equivalences of categories
$$\Der^b_c(X(\CC),\CC)\simeq \ho(\Mod{\mathcal{E}_{\mathit{Betti}}\otimes_\QQ\CC})_c(X)
\simeq \Der_{\mathit{dR,an}}(X)$$
which commute with the six operations as well as with the realization
functors. In particular, for $X$ smooth,
by the Riemann-Hilbert correspondence\index{word}{Riemann-Hilbert},
$\Der_{\mathit{dR,an}}(X)$ is equivalent to the bounded derived category of
analytic regular holonomic $\mathcal{D}$-modules on $X$
which are constructible of geometric origin.
But we can go backward: as proved by Brad Drew in his thesis \cite{drew1},
using Corollary \ref{compWeil}, one can actually show directly (i.e. motivically)
that, for $X$ smooth, $\Der_{\mathit{dR,an}}(X)$ is equivalent to the bounded derived category of
analytic regular holonomic $\mathcal{D}$-modules on $X$
which are constructible of geometric origin, and thus give a new algebraic
proof of the Riemann-Hilbert correspondence. 
\end{ex}

\begin{ex}\label{ex:rigidcoh}
Let $V$ be a complete discrete valuation ring of mixed characteristic
with perfect residue field $k$ and field of functions $K$.
The Monsky-Washnitzer complex
\index{word}{complex!Monsky-Washnitzer}
 defines a stable cohomology
$E_{MW}$ over smooth $V$-schemes of finite type, defined by
$$E_{MW}(X)=\Omega^*_{A^\dagger/V}\otimes_V K$$
for any affine smooth $V$-scheme $X=\spec A$ (the case of a smooth
$V$-scheme of finite type is obtained by Zariski descent); see \cite[3.2.3]{CD2}.
Let $\mathcal{E}_{MW}$ be the corresponding ring spectrum in
$\DMB(\spec V)$, and write $j:\spec K\To \spec V$ and $i:\spec k \To \spec V$
for the canonical immersions. As we obviously have $j^*\mathcal{E}_{MW}=0$
(the Monsky-Washnitzer cohomology
\index{word}{cohomology!Monsky-Washnitzer}
 of a smooth $V$-scheme with empty
special fiber vanishes), we have a canonical isomorphism
$$\mathcal{E}_{MW}\simeq \derR i_*\, \derL i^*\mathcal{E}_{MW}\, .$$
We define the rigid cohomology spectrum $\mathcal{E}_{\mathit{rig}}$
in $\DMB(\spec k)$ by the formula
$$\mathcal{E}_{\mathit{rig}}= \derL i^*\mathcal{E}_{MW}\, .$$
This is a ring spectrum associated to a $K$-linear mixed Weil cohomology:
cohomology with coefficients in $\mathcal{E}_{\mathit{rig}}$
coincides with rigid cohomology
\index{word}{cohomology!rigid}
 in the sense of Berthelot,
and the K\"unneth formula for rigid cohomology holds
for smooth and projective $k$-schemes (as rigid
cohomology coincides then with crystalline cohomology), from which
we deduce the K\"unneth formula for smooth $k$-schemes of finite type;
see \cite[3.2.10]{CD2}. As before, we define
$$\Der_{\mathit{rig}}(X)=\Der_c(X,\mathcal{E}_\mathit{rig})$$
for any $k$-scheme of finite type $X$, and put
$$H^n_\mathit{rig}(X)=
\Hom_{\Der_{\mathit{rig}}(X)}(\mathcal{E}_{\mathit{rig},X},\mathcal{E}_{\mathit{rig},X}[n])\, .$$
Here again, we have, by construction, \emph{rigid realization functors}
\index{word}{realization functor!rigid}
$$R_{\mathit{rig}}:\DMV_{\mathcyr{B},c}(X)\To \Der_{\mathit{rig}}(X)$$
which preserve the six operations of Grothendieck (Theorem \ref{thm:realcommute6op}),
as well as (higher) cycle class maps
$\HB^q(X,\QQ(p))\To H^q_{\mathit{rig}}(X)(p)$.
\end{ex}

\backmatter
\nocite{HCHA,Bondarko,Bondarko0,Bondarko1,cdhweights,BondDeg,Olsson1,Olsson2,pepin1,pepin2,
scholbach1,scholbach2,shtukas,Soergel,Soergel2,Steenrod,hoyois1,integraltransfers,Wildes1,Wildes2,
vaish1,vaish2,vaish3,vaish4,drew1,drew2,BRTV,Robalo,lize1,lize2,lize3,jinyang,motivicHZ,AnLePa,
sixlectures,obituary,intromot,isosuslin,CD4,CD5,balzin,Cis5,ToVe,ayoub5,kahn,BaVKa,ayoub4,
totaro1,totaro2,Coniveau,localization,localization2,DLORV,levine5,levine6,logmotives,
jannsen1,jannsen2,seattle,Bei2,Bei3,periodsNori,AHuber,CD4,bachmann,aoki,
roendigs,hanamura0,hanamura1,hanamura2,hanamura3,arndt,gark,gluemixedmotives,
Rydh}
\bibliographystyle{alpha}
\markboth{}{References}
\bibliography{baseqc}

\newcommand{\etalchar}[1]{$^{#1}$}
\begin{thebibliography}{DL{\O}{\etalchar{+}}07}

\bibitem[AGV73]{SGA4}
M.~Artin, A.~Grothendieck, and J.-L. Verdier.
\newblock {\em Th\'eorie des topos et cohomologie \'etale des sch\'emas},
  volume 269, 270, 305 of {\em Lecture Notes in Mathematics}.
\newblock Springer-Verlag, 1972--1973.
\newblock S\'eminaire de G\'eom\'etrie Alg\'ebrique du Bois--Marie 1963--64
  (SGA~4).

\bibitem[ALP17]{AnLePa}
A.~Ananyevskiy, M.~Levine, and I.~Panin.
\newblock Witt sheaves and the {$\eta$}-inverted sphere spectrum.
\newblock {\em J. Topol.}, 10(2):370--385, 2017.

\bibitem[And04]{intromot}
Y.~Andr\'{e}.
\newblock {\em Une introduction aux motifs (motifs purs, motifs mixtes,
  p\'{e}riodes)}, volume~17 of {\em Panoramas et Synth\`eses}.
\newblock Soci\'{e}t\'{e} Math\'{e}matique de France, Paris, 2004.

\bibitem[Aok19]{aoki}
K.~Aoki.
\newblock The weight complex functor is symmetric monoidal.
\newblock arXiv$:$1904.01384, 2019.

\bibitem[Arn16]{arndt}
P.~Arndt.
\newblock {\em Abstract motivic homotopy theory}.
\newblock PhD thesis, University of Osnabr{\"u}ck, Osnabr{\"u}ck, Germany,
  2016.

\bibitem[Ayo07a]{ayoub}
J.~Ayoub.
\newblock {\em Les six op\'erations de {G}rothendieck et le formalisme des
  cycles \'evanescents dans le monde motivique ({I})}, volume 314 of {\em
  Ast\'erisque}.
\newblock Soc. Math. France, 2007.

\bibitem[Ayo07b]{ayoub2}
J.~Ayoub.
\newblock {\em Les six op\'erations de {G}rothendieck et le formalisme des
  cycles \'evanescents dans le monde motivique ({II})}, volume 315 of {\em
  Ast\'erisque}.
\newblock Soc. Math. France, 2007.

\bibitem[Ayo10]{ayoub3}
J.~Ayoub.
\newblock Note sur les op\'erations de {G}rothendieck et la r\'ealisation de
  {B}etti.
\newblock {\em J. Inst. Math. Jussieu}, 9(2):225--263, 2010.

\bibitem[Ayo14]{ayoub5}
J.~Ayoub.
\newblock La r\'ealisation \'etale et les op\'erations de {G}rothendieck.
\newblock {\em Ann. Sci. \'Ecole Norm. Sup.}, 47(1):1--145, 2014.

\bibitem[AZ12]{ayoub4}
J.~Ayoub and S.~Zucker.
\newblock Relative {A}rtin motives and the reductive {B}orel-{S}erre
  compactification of a locally symmetric variety.
\newblock {\em Invent. Math.}, 188(2):277--427, 2012.

\bibitem[Bac18a]{bachmann}
T.~Bachmann.
\newblock Motivic and real {\'e}tale stable homotopy theory.
\newblock {\em Compos. Math.}, 154(5):883--917, 2018.

\bibitem[Bac18b]{bachmann2}
T.~Bachmann.
\newblock Rigidity in etale motivic stable homotopy theory.
\newblock arXiv$:$1810.08028, 2018.

\bibitem[Bal19]{balzin}
E.~Balzin.
\newblock Reedy model structures in families.
\newblock arXiv$:$1803.00681v2, 2019.

\bibitem[Bar10]{Bar}
C.~Barwick.
\newblock On left and right model categories and left and right {B}ousfield
  localizations.
\newblock {\em Homology, Homotopy Appl.}, 12(2):245--320, 2010.

\bibitem[BBD82]{BBD}
A.~Be{\u\i}linson, J.~Bernstein, and P.~Deligne.
\newblock Faisceaux pervers.
\newblock {\em Ast\'erisque}, 100:5--171, 1982.

\bibitem[BD17]{BondDeg}
M.~V. Bondarko and F.~D\'{e}glise.
\newblock Dimensional homotopy t-structures in motivic homotopy theory.
\newblock {\em Adv. Math.}, 311:91--189, 2017.

\bibitem[Be{\u\i}84]{Bei2}
A.~Be{\u\i}linson.
\newblock Higher regulators and values of {$L$}-functions.
\newblock In {\em Current problems in mathematics, {V}ol. 24}, Itogi Nauki i
  Tekhniki, pages 181--238. Akad. Nauk SSSR, Vsesoyuz. Inst. Nauchn. i Tekhn.
  Inform., Moscow, 1984.

\bibitem[Be{\u\i}87]{Bei}
A.~Be{\u\i}linson.
\newblock Height pairing between algebraic cycles.
\newblock In {\em Current trends in arithmetical algebraic geometry ({A}rcata,
  {C}alif., 1985)}, volume~67 of {\em Contemp. Math.}, pages 1--24. Amer. Math.
  Soc., 1987.

\bibitem[Be{\u\i}12]{Bei3}
A.~Be{\u\i}linson.
\newblock Remarks on {G}rothendieck's standard conjectures.
\newblock In {\em Regulators}, volume 571 of {\em Contemp. Math.}, pages
  25--32. Amer. Math. Soc., Providence, RI, 2012.

\bibitem[Bek00]{beke1}
T.~Beke.
\newblock Sheafifiable homotopy model categories.
\newblock {\em Math. Proc. Camb. Phil. Soc.}, 129:447--475, 2000.

\bibitem[BG73]{BG}
K.~S. Brown and S.~M. Gersten.
\newblock Algebraic {$K$}-theory as generalized sheaf cohomology.
\newblock In {\em Higher {K}-theories {I} (Proc. Conf., Battelle Memorial
  Inst., Seattle, Wash., 1972)}, volume 341 of {\em Lecture Notes in
  Mathematics}, pages 266--292. Springer-Verlag, 1973.

\bibitem[BGI71]{SGA6}
P.~Berthelot, A.~Grothendieck, and L.~Illusie.
\newblock {\em Th\'eorie des intersections et th\'eor\`eme de
  {R}iemann-{R}och}, volume 225 of {\em Lecture Notes in Mathematics}.
\newblock Springer-Verlag, 1971.
\newblock S\'eminaire de G\'eom\'etrie Alg\'ebrique du Bois--Marie 1966--67
  (SGA~6).

\bibitem[BI15]{cdhweights}
M.~V. Bondarko and M.~A. Ivanov.
\newblock On {C}how weight structures for {$cdh$}-motives with integral
  coefficients.
\newblock {\em Algebra i Analiz}, 27(6):14--40, 2015.
\newblock Reprinted in St. Petersburg Math. J. {{\bf{2}}7} (2016), no. 6,
  869--888.

\bibitem[Bla03]{BL}
B.~Blander.
\newblock Local projective model structures on simplicial presheaves.
\newblock {\em $K$-Theory}, 24(3):283--301, 2003.

\bibitem[BM03]{BerMoe1}
C.~Berger and I.~Moerdijk.
\newblock Axiomatic homotopy theory for operads.
\newblock {\em Comment. Math. Helv.}, 78:805--831, 2003.

\bibitem[BM09]{BerMoe2}
C.~Berger and I.~Moerdijk.
\newblock On the derived category of an algebra over an operad.
\newblock {\em Georgian Math. J.}, 16(1):13--28, 2009.

\bibitem[Bon09]{Bondarko0}
M.~V. Bondarko.
\newblock Differential graded motives: weight complex, weight filtrations and
  spectral sequences for realizations; {V}oevodsky versus {H}anamura.
\newblock {\em J. Inst. Math. Jussieu}, 8(1):39--97, 2009.

\bibitem[Bon14]{Bondarko}
M.~V. Bondarko.
\newblock Weights for relative motives: relation with mixed complexes of
  sheaves.
\newblock {\em Int. Math. Res. Not. IMRN}, (17):4715--4767, 2014.

\bibitem[Bon15]{Bondarko1}
M.~V. Bondarko.
\newblock Mixed motivic sheaves (and weights for them) exist if `ordinary'
  mixed motives do.
\newblock {\em Compos. Math.}, 151(5):917--956, 2015.

\bibitem[Bou93]{BourbakiAC8&9}
N.~Bourbaki.
\newblock {\em \'{E}l\'ements de math\'ematique. {A}lg\`ebre commutative.
  {C}hapitres 8 et 9}.
\newblock Masson, 193.
\newblock Reprint of the 1983 original.

\bibitem[Bou98]{BourbakiAC}
N.~Bourbaki.
\newblock {\em Commutative algebra. {C}hapters 1--7}.
\newblock Elements of Mathematics (Berlin). Springer-Verlag, Berlin, 1998.
\newblock Translated from the French, Reprint of the 1989 English translation.

\bibitem[Bro74]{Brown}
K.~S. Brown.
\newblock Abstract homotopy theory and generalized sheaf cohomology.
\newblock {\em Trans. Amer. Math. Soc.}, 186:419--458, 1974.

\bibitem[BRTV18]{BRTV}
A.~Blanc, M.~Robalo, B.~To\"{e}n, and G.~Vezzosi.
\newblock Motivic realizations of singularity categories and vanishing cycles.
\newblock {\em J. \'{E}c. polytech. Math.}, 5:651--747, 2018.

\bibitem[BVK16]{BaVKa}
L.~Barbieri-Viale and B.~Kahn.
\newblock {\em On the derived category of 1-motives}, volume 381 of {\em
  Ast\'{e}risque}.
\newblock Soc. Math. France, 2016.

\bibitem[CD09]{CD1}
D.-C. Cisinski and F.~D{\'e}glise.
\newblock Local and stable homological algebra in {G}rothendieck abelian
  categories.
\newblock {\em Homology, Homotopy and Applications}, 11(1):219--260, 2009.

\bibitem[CD12]{CD2}
D.-C. Cisinski and F.~D{\'e}glise.
\newblock Mixed {W}eil cohomologies.
\newblock {\em Adv. Math.}, 230(1):55--130, 2012.

\bibitem[CD15]{CD5}
D.-C. Cisinski and F.~D{\'e}glise.
\newblock Integral mixed motives in equal characteristics.
\newblock {\em Doc. Math.}, (Extra volume: Alexander S. Merkurjev's sixtieth
  birthday):145--194, 2015.

\bibitem[CD16]{CD4}
D.-C. Cisinski and F.~D{\'e}glise.
\newblock \'{E}tale motives.
\newblock {\em Compos. Math.}, 152(3):556--666, 2016.

\bibitem[Cis03]{Cis1}
D.-C. Cisinski.
\newblock Images directes cohomologiques dans les cat\'egories de mod\`eles.
\newblock {\em Annales Math\'ematiques Blaise Pascal}, 10:195--244, 2003.

\bibitem[Cis06]{Cis3}
D.-C. Cisinski.
\newblock {\em Les pr\'efaisceaux comme mod\`eles des types d'homotopie},
  volume 308 of {\em Ast\'erisque}.
\newblock Soc. Math. France, 2006.

\bibitem[Cis08]{Cis2}
D.-C. Cisinski.
\newblock Propri\'et\'es universelles et extensions de {K}an d\'eriv\'ees.
\newblock {\em Theory and Applications of Categories}, 20(17):605--649, 2008.

\bibitem[Cis13]{cis4}
D.-C. Cisinski.
\newblock Descente par \'eclatements en {K}-th\'eorie invariante par homotopie.
\newblock {\em Ann. of Math.}, 177, 2013.

\bibitem[Cis18]{Cis5}
D.-C. Cisinski.
\newblock Cohomological methods in intersection theory.
\newblock Lecture notes from the LMS-CMI Research School `Homotopy Theory and
  Arithmetic Geometry: Motivic and Diophantine Aspects', Imperial College
  London, 2018.

\bibitem[Cis19]{HCHA}
D.-C. Cisinski.
\newblock {\em Higher categories and homotopical algebra}, volume 180 of {\em
  Cambridge studies in advanced mathematics}.
\newblock Cambridge University Press, 2019.

\bibitem[Con07]{Conrad}
B.~Conrad.
\newblock Deligne's notes on {N}agata compactifications.
\newblock {\em J. Ramanujan Math. Soc.}, 22(3):205--257, 2007.

\bibitem[Cra95]{crans}
S.~E. Crans.
\newblock Quillen closed model structures for sheaves.
\newblock {\em J. Pure Appl. Algebra}, 101:35--57, 1995.

\bibitem[CT11]{CisTab}
D.-C. Cisinski and G.~Tabuada.
\newblock Non connective ${K}$-theory via universal invariants.
\newblock {\em Compos. Math.}, 147(4):1281--1320, 2011.

\bibitem[D{\'e}g07]{Deg7}
F.~D{\'e}glise.
\newblock Finite correspondences and transfers over a regular base.
\newblock In {\em Algebraic cycles and motives. Vol. 1}, volume 343 of {\em
  London Math. Soc. Lecture Note Ser.}, pages 138--205. Cambridge Univ. Press,
  Cambridge, 2007.

\bibitem[D{\'e}g08]{Deg8}
F.~D{\'e}glise.
\newblock Around the {G}ysin triangle {II}.
\newblock {\em Doc. Math.}, 13:613--675, 2008.

\bibitem[Del71]{DelH3}
P.~Deligne.
\newblock Th\'{e}orie de {H}odge. {II}.
\newblock {\em Inst. Hautes \'{E}tudes Sci. Publ. Math.}, (40):5--57, 1971.

\bibitem[Del74]{DelH2}
P.~Deligne.
\newblock Th\'{e}orie de {H}odge. {III}.
\newblock {\em Inst. Hautes \'{E}tudes Sci. Publ. Math.}, (44):5--77, 1974.

\bibitem[Del87]{Del_det}
P.~Deligne.
\newblock Le d\'eterminant de la cohomologie.
\newblock In {\em Current trends in arithmetical algebraic geometry ({A}rcata,
  {C}alif., 1985)}, volume~67 of {\em Contemp. Math.}, pages 93--177. Amer.
  Math. Soc., 1987.

\bibitem[Del01]{Delnotes}
P.~Deligne.
\newblock Voevodsky lectures on cross functors.
\newblock unpublished notes available at
  {\url{https://www.math.ias.edu/vladimir/node/94}}, Fall 2001.

\bibitem[DFJK19]{DFJK}
F.~D{\'e}glise, J.~Fasel, F.~Jin, and A.~Khan.
\newblock Borel isomorphism and absolute purity.
\newblock arXiv$:$1902.02055, 2019.

\bibitem[dJ96]{dejong}
A.~J. de~Jong.
\newblock Smoothness, semi-stability and alterations.
\newblock {\em Publ. Math. IHES}, 83:51--93, 1996.

\bibitem[dJ97]{dejong2}
A.~J. de~Jong.
\newblock Families of curves and alterations.
\newblock {\em Annales de l'institut Fourier}, 47(2):599--621, 1997.

\bibitem[DL{\O}{\etalchar{+}}07]{DLORV}
B.~I. Dundas, M.~Levine, P.~A. {\O}stv{\ae}r, O.~R\"{o}ndigs, and V.~Voevodsky.
\newblock {\em Motivic homotopy theory}.
\newblock Universitext. Springer-Verlag, Berlin, 2007.
\newblock Lectures from the Summer School held in Nordfjordeid, August 2002.

\bibitem[Dre13]{drew1}
B.~Drew.
\newblock {\em R\'ealisations tannakiennes des motifs mixtes triangul\'es}.
\newblock PhD thesis, Univ.~Paris~13, 2013.

\bibitem[Dre18]{drew2}
B.~Drew.
\newblock Motivic hodge modules.
\newblock arXiv$:$1801.10129, 2018.

\bibitem[Dug01]{dug0}
D.~Dugger.
\newblock Combinatorial model categories have presentations.
\newblock {\em Adv. Math}, 164(1):177--201, 2001.

\bibitem[Dug06]{dug1}
D.~Dugger.
\newblock Spectral enrichment of model categories.
\newblock {\em Homology, Homotopy, and Applications}, 8(1):1--30, 2006.

\bibitem[Eke90]{Ekedahl}
T.~Ekedahl.
\newblock On the adic formalism.
\newblock In {\em The {G}rothendieck {F}estschrift, {V}ol. {II}}, volume~87 of
  {\em Progr. Math.}, pages 197--218. Birkh\"{a}user Boston, Boston, MA, 1990.

\bibitem[Ful98]{Ful}
W.~Fulton.
\newblock {\em Intersection theory}.
\newblock Springer, second edition, 1998.

\bibitem[Gar19]{gark}
G.~Garkusha.
\newblock Reconstructing rational stable motivic homotopy theory.
\newblock arXiv$:$1705.01635v3, 2019.

\bibitem[GD60]{EGA1}
A.~Grothendieck and J.~Dieudonn\'e.
\newblock \'{E}l\'ements de g\'eom\'etrie alg\'ebrique. {I}. {L}e langage des
  sch\'emas.
\newblock {\em Publ. Math. IHES}, 4, 1960.

\bibitem[GD61]{EGA2}
A.~Grothendieck and J.~Dieudonn\'e.
\newblock \'{E}l\'ements de g\'eom\'etrie alg\'ebrique. {II}. \'{E}tude globale
  \'el\'ementaire de quelques classes de morphismes.
\newblock {\em Publ. Math. IHES}, 8, 1961.

\bibitem[GD67]{EGA4}
A.~Grothendieck and J.~Dieudonn\'e.
\newblock \'{E}l\'ements de g\'eom\'etrie alg\'ebrique. {IV}. \'{E}tude locale
  des sch\'emas et des morphismes de sch\'emas {IV}.
\newblock {\em Publ. Math. IHES}, 20, 24, 28, 32, 1964-1967.

\bibitem[GG16]{gulet1}
S.~Gorchinskiy and V.~Guletski\u{\i}.
\newblock Symmetric powers in abstract homotopy categories.
\newblock {\em Adv. Math.}, 292:707--754, 2016.

\bibitem[GG18]{gulet2}
S.~Gorchinskiy and V.~Guletski\u{\i}.
\newblock Positive model structures for abstract symmetric spectra.
\newblock {\em Appl. Categ. Structures}, 26(1):29--46, 2018.

\bibitem[Gil81]{Gillet}
H.~Gillet.
\newblock Riemann-{R}och theorems for higher algebraic {$K$}-theory.
\newblock {\em Adv. in Math.}, 40(3):203--289, 1981.

\bibitem[Gro57]{Gro1}
A.~Grothendieck.
\newblock Sur quelques points d'alg\`ebre homologique.
\newblock {\em T\^ohoku Math. J.}, 9(2):119--221, 1957.

\bibitem[Gro77]{SGA5}
A.~Grothendieck.
\newblock {\em Cohomologie l-adique et {F}onctions {L}}, volume 589 of {\em
  Lecture Notes in Mathematics}.
\newblock Springer-Verlag, 1977.
\newblock S\'eminaire de G\'eom\'etrie Alg\'ebrique du Bois--Marie 1965--66
  (SGA~5).

\bibitem[Gro03]{SGA1}
A.~Grothendieck.
\newblock {\em Rev\^etements \'etales et groupe fondamental}, volume~3 of {\em
  Documents Math\'ematiques}.
\newblock Soc. Math. France, 2003.
\newblock S\'eminaire de G\'eom\'etrie Alg\'ebrique du Bois--Marie 1960--61
  (SGA~1). \'Edition recompos\'ee et annot\'ee du LNM 224, Springer, 1971.

\bibitem[GZ67]{GZ}
P.~Gabriel and M.~Zisman.
\newblock {\em Calculus of fractions and homotopy theory}, volume~35 of {\em
  Ergebnisse der {M}athematik}.
\newblock Springer-Verlag, 1967.

\bibitem[Han95]{hanamura1}
M.~Hanamura.
\newblock Mixed motives and algebraic cycles {I}.
\newblock {\em Math. Res. Lett.}, 2(6):811--821, 1995.

\bibitem[Han99]{hanamura3}
M.~Hanamura.
\newblock Mixed motives and algebraic cycles {III}.
\newblock {\em Math. Res. Lett.}, 6(1):61--82, 1999.

\bibitem[Han00]{hanamura0}
M.~Hanamura.
\newblock Homological and cohomological motives of algebraic varieties.
\newblock {\em Invent. Math.}, 142(2):319--349, 2000.

\bibitem[Han04]{hanamura2}
M.~Hanamura.
\newblock Mixed motives and algebraic cycles {II}.
\newblock {\em Invent. Math.}, 158(1):105--179, 2004.

\bibitem[Har66]{Har}
R.~Hartshorne.
\newblock {\em Residues and duality}.
\newblock Lecture notes of a seminar on the work of A. Grothendieck, given at
  Harvard 1963/64. With an appendix by P. Deligne. Lecture Notes in
  Mathematics, No. 20. Springer-Verlag, Berlin, 1966.

\bibitem[H{\'e}b11]{Hebert}
D.~H{\'e}bert.
\newblock Structure de poids \`a la {B}ondarko sur les motifs de {B}eilinson.
\newblock {\em Compos. Math.}, 147(5):1447--1462, 2011.

\bibitem[HK{\O}17]{Steenrod}
M.~Hoyois, S.~Kelly, and P.~A. {\O}stv{\ae}r.
\newblock The motivic {S}teenrod algebra in positive characteristic.
\newblock {\em J. Eur. Math. Soc.}, 19(12):3813--3849, 2017.

\bibitem[HMS17]{periodsNori}
A.~Huber and S.~M\"{u}ller-Stach.
\newblock {\em Periods and {N}ori motives}, volume~65 of {\em Ergebnisse der
  Mathematik und ihrer Grenzgebiete. 3. Folge. A Series of Modern Surveys in
  Mathematics [Results in Mathematics and Related Areas. 3rd Series. A Series
  of Modern Surveys in Mathematics]}.
\newblock Springer, Cham, 2017.
\newblock With contributions by Benjamin Friedrich and Jonas von Wangenheim.

\bibitem[Hor13]{hornbostel}
J.~Hornbostel.
\newblock Preorientations of the derived motivic multiplicative group.
\newblock {\em Algebr. Geom. Topol.}, 13(5):2667--2712, 2013.

\bibitem[Hov99]{Hovey}
M.~Hovey.
\newblock {\em Model Categories}, volume~63 of {\em Math. surveys and
  monographs}.
\newblock Amer. Math. Soc., 1999.

\bibitem[Hov01]{Hov}
M.~Hovey.
\newblock Spectra and symmetric spectra in general model categories.
\newblock {\em J. Pure Appl. Algebra}, 165(1):63--127, 2001.

\bibitem[Hoy17]{hoyois1}
M.~Hoyois.
\newblock The six operations in equivariant motivic homotopy theory.
\newblock {\em Adv. Math.}, 305:197--279, 2017.
\newblock Corrigendum in Adv. Math. 333 (2018).

\bibitem[HS15]{scholbach1}
A.~Holmstrom and J.~Scholbach.
\newblock Arakelov motivic cohomology {I}.
\newblock {\em J. Algebraic Geom.}, 24(4):719--754, 2015.

\bibitem[Hub00]{AHuber}
Annette Huber.
\newblock Realization of {V}oevodsky's motives.
\newblock {\em J. Algebraic Geom.}, 9(4):755--799, 2000.
\newblock Corrigendum in J. Algebraic Geom. 13 (2004).

\bibitem[ILO14]{gabber3}
Luc Illusie, Yves Laszlo, and Fabrice Orgogozo, editors.
\newblock {\em Travaux de {G}abber sur l'uniformisation locale et la
  cohomologie \'{e}tale des sch\'{e}mas quasi-excellents}, volume 363--364 of
  {\em Ast\'{e}risque}.
\newblock Soc. Math. France, 2014.
\newblock S\'{e}minaire \`a l'\'{E}cole Polytechnique 2006--2008. [Seminar of
  the Polytechnic School 2006--2008], With the collaboration of
  Fr\'{e}d\'{e}ric D\'{e}glise, Alban Moreau, Vincent Pilloni, Michel Raynaud,
  Jo\"{e}l Riou, Beno\^{\i}t Stroh, Michael Temkin and Weizhe Zheng.

\bibitem[Ivo07]{Ivo}
F.~Ivorra.
\newblock R\'ealisation {$l$}-adique des motifs triangul\'es g\'eom\'etriques.
  {I}.
\newblock {\em Doc. Math.}, 12:607--671, 2007.

\bibitem[Jan90]{jannsen2}
U.~Jannsen.
\newblock {\em Mixed motives and algebraic {$K$}-theory}, volume 1400 of {\em
  Lecture Notes in Mathematics}.
\newblock Springer-Verlag, Berlin, 1990.
\newblock With appendices by S. Bloch and C. Schoen.

\bibitem[Jan92]{jannsen1}
U.~Jannsen.
\newblock Motives, numerical equivalence, and semi-simplicity.
\newblock {\em Invent. Math.}, 107(3):447--452, 1992.

\bibitem[Jar00]{Jar}
J.~F. Jardine.
\newblock Motivic symmetric spectra.
\newblock {\em Doc. Math.}, 5:445--553, 2000.

\bibitem[Jar15]{JarLoc}
J.~F. Jardine.
\newblock {\em Local homotopy theory}.
\newblock Springer Monographs in Mathematics. Springer, 2015.

\bibitem[JKS94]{seattle}
Uwe Jannsen, Steven Kleiman, and Jean-Pierre Serre, editors.
\newblock {\em Motives}, volume~55 of {\em Proceedings of Symposia in Pure
  Mathematics}. American Mathematical Society, Providence, RI, 1994.

\bibitem[JY19]{jinyang}
F.~Jin and E.~Yang.
\newblock K{\"u}nneth formulas for motives and additivity of traces.
\newblock arXiv$:$1812.06441v2, 2019.

\bibitem[Kah18]{kahn}
B.~Kahn.
\newblock {\em Fonctions z\^{e}ta et {$L$} de vari\'{e}t\'{e}s et de motifs}.
\newblock Nano. Calvage et Mounet, Paris, 2018.

\bibitem[Kel17]{integraltransfers}
S.~Kelly.
\newblock {\em Voevodsky motives and {$l$}dh-descent}, volume 391 of {\em
  Ast\'erisque}.
\newblock Soc. Math. France, 2017.

\bibitem[Kel18]{isosuslin}
S.~Kelly.
\newblock Un isomorphisme de {S}uslin.
\newblock {\em Bull. Soc. Math. France}, 146(4):633--647, 2018.

\bibitem[Kra01]{Krause}
H.~Krause.
\newblock On {N}eeman's well generated triangulated categories.
\newblock {\em Doc. Math.}, 6:121--126, 2001.

\bibitem[Leh17]{pepin2}
S.~Pepin Lehalleur.
\newblock Constructible $1$-motives and exactness.
\newblock arXiv$:$1712.01180, 2017.

\bibitem[Leh19]{pepin1}
S.~Pepin Lehalleur.
\newblock Triangulated categories of relative $1$-motives.
\newblock arXiv$:$1512.00266v3, to appear in Advances in Mathematics, 2019.

\bibitem[Lev98]{LevineMM}
M.~Levine.
\newblock {\em Mixed motives}, volume~57 of {\em Math. surveys and monographs}.
\newblock Amer. Math. Soc., 1998.

\bibitem[Lev01]{localization}
M.~Levine.
\newblock Techniques of localization in the theory of algebraic cycles.
\newblock {\em J. Alg. Geom.}, 10:299--363, 2001.

\bibitem[Lev06]{localization2}
M.~Levine.
\newblock Chow's moving lemma and the homotopy coniveau tower.
\newblock {\em $K$-Theory}, 37(1-2):129--209, 2006.

\bibitem[Lev08]{Coniveau}
M.~Levine.
\newblock The homotopy coniveau tower.
\newblock {\em J. Topol.}, 1(1):217--267, 2008.

\bibitem[Lev10]{levine5}
M.~Levine.
\newblock Slices and transfers.
\newblock {\em Doc. Math.}, (Extra vol.: Andrei A. Suslin sixtieth
  birthday):393--443, 2010.

\bibitem[Lev13a]{levine6}
M.~Levine.
\newblock Convergence of {V}oevodsky's slice tower.
\newblock {\em Doc. Math.}, 18:907--941, 2013.

\bibitem[Lev13b]{sixlectures}
M.~Levine.
\newblock Six lectures on motives.
\newblock In {\em Autour des motifs---\'{E}cole d'\'{e}t\'{e}
  {F}ranco-{A}siatique de {G}\'{e}om\'{e}trie {A}lg\'{e}brique et de
  {T}h\'{e}orie des {N}ombres/{A}sian-{F}rench {S}ummer {S}chool on {A}lgebraic
  {G}eometry and {N}umber {T}heory. {V}ol. {II}}, volume~41 of {\em Panor.
  Synth\`eses}, pages 1--141. Soc. Math. France, Paris, 2013.

\bibitem[Lev18]{obituary}
M.~Levine.
\newblock Vladimir {V}oevodsky---an appreciation.
\newblock {\em Bull. Amer. Math. Soc. (N.S.)}, 55(4):405--425, 2018.

\bibitem[Lip78]{lipman}
J.~Lipman.
\newblock Desingularization of two dimensional schemes.
\newblock {\em Annals of Math.}, 107:151--207, 1978.

\bibitem[Lur17]{DAG3}
J.~Lurie.
\newblock Higher algebra.
\newblock preprint, 2017.

\bibitem[LZ15]{lize2}
Y.~Liu and W.~Zheng.
\newblock Gluing restricted nerves of $\infty$-categories.
\newblock arXiv$:$1211.5294v4, 2015.

\bibitem[LZ16]{lize1}
Y.~Liu and W.~Zheng.
\newblock Gluing pseudo functors via $n$-fold categories.
\newblock arXiv$:$1211.1877v6, 2016.

\bibitem[LZ17]{lize3}
Y.~Liu and W.~Zheng.
\newblock Enhanced six operations and base change theorem for higher {A}rtin
  stacks.
\newblock arXiv$:$1211.5948v3, 2017.

\bibitem[Mac98]{MLa}
S.~MacLane.
\newblock {\em Categories for the working mathematician}, volume~5 of {\em
  Graduate Texts in Mathematics}.
\newblock Springer-Verlag, New York, second edition, 1998.

\bibitem[Mal05]{Mal}
G.~Maltsiniotis.
\newblock {\em La th\'eorie de l'homotopie de {G}rothendieck}, volume 301 of
  {\em Ast\'erisque}.
\newblock Soc. Math. France, 2005.

\bibitem[Mat70]{Mat}
H.~Matsumura.
\newblock {\em Commutative algebra}.
\newblock W. A. Benjamin, Inc., New York, 1970.

\bibitem[Mor03]{dmtilde2}
F.~Morel.
\newblock An introduction to {$\mathbf{A}^1$}-homotopy theory.
\newblock In {\em Contemporary Developments in {A}lgebraic {K}-theory},
  volume~15 of {\em ICTP Lecture notes}, pages 357--441. 2003.

\bibitem[Mor04a]{dmtilde}
F.~Morel.
\newblock On the motivic {$\pi\sb 0$} of the sphere spectrum.
\newblock In {\em Axiomatic, enriched and motivic homotopy theory}, volume 131
  of {\em NATO Sci. Ser. II Math. Phys. Chem.}, pages 219--260. Kluwer Acad.
  Publ., Dordrecht, 2004.

\bibitem[Mor04b]{KthMW}
F.~Morel.
\newblock Sur les puissances de l'id\'eal fondamental de l'anneau de {W}itt.
\newblock {\em Comment. Math. Helv.}, 79(4):689--703, 2004.

\bibitem[Mor05]{dmtilde3}
F.~Morel.
\newblock The stable $\mathbf{A}^1$-connectivity theorems.
\newblock {\em K-theory}, 35:1--68, 2005.

\bibitem[Mor06]{SHQ}
F.~Morel.
\newblock Rational stable splitting of {G}rassmanians and the rational motivic
  sphere spectrum.
\newblock statement of results, 2006.

\bibitem[Mor12]{A1brouwer}
F.~Morel.
\newblock {\em $\mathbf{{A}}^1$-algebraic topology over a field}, volume 2052
  of {\em Lecture Notes in Mathematics}.
\newblock Springer, 2012.

\bibitem[MV99]{MV}
F.~Morel and V.~Voevodsky.
\newblock {$\mathbf{A}^1$-homotopy} theory of schemes.
\newblock {\em Publ. Math. IHES}, 90:45--143, 1999.

\bibitem[MVW06]{MVW}
C.~Mazza, V.~Voevodsky, and C.~Weibel.
\newblock {\em Lecture notes on motivic cohomology}, volume~2 of {\em Clay
  Mathematics Monographs}.
\newblock American Mathematical Society, Providence, RI, 2006.

\bibitem[Nee92]{Nee4}
A.~Neeman.
\newblock The connection between the {$K$}-theory localization theorem of
  {T}homason, {T}robaugh and {Y}ao and the smashing subcategories of
  {B}ousfield and {R}avenel.
\newblock {\em Ann. Sci. {\'{E}}cole Norm. Sup. (4)}, 25(5):547--566, 1992.

\bibitem[Nee96]{Nee3}
A.~Neeman.
\newblock The {G}rothendieck duality theorem via {B}ousfield's techniques and
  {B}rown representability.
\newblock {\em J. Amer. Math. Soc.}, 9(1):205--236, 1996.

\bibitem[Nee01]{Nee1}
A.~Neeman.
\newblock {\em Triangulated categories}, volume 148 of {\em Annals of
  Mathematics Studies}.
\newblock Princeton University Press, Princeton, NJ, 2001.

\bibitem[NS{\O}09]{NOS}
N.~Nauman, M.~Spitzweck, and P.~A. {\O}stv{\ae}r.
\newblock Motivic {L}andweber exactness.
\newblock {\em Doc. Math.}, 14:551--593, 2009.

\bibitem[NV15]{vaish1}
A.N. Nair and V.~Vaish.
\newblock Weightless cohomology of algebraic varieties.
\newblock {\em J. Algebra}, 424:147--189, 2015.

\bibitem[Ols15]{Olsson1}
M.~Olsson.
\newblock Borel-{M}oore homology, {R}iemann-{R}och transformations, and local
  terms.
\newblock {\em Adv. Math.}, 273:56--123, 2015.

\bibitem[Ols16]{Olsson2}
M.~Olsson.
\newblock Motivic cohomology, localized {C}hern classes, and local terms.
\newblock {\em Manuscripta Math.}, 149(1-2):1--43, 2016.

\bibitem[Par17a]{gluemixedmotives}
D.~Park.
\newblock Construction of triangulated categories of motives using the
  localization property.
\newblock arXiv$:$1708.00597, 2017.

\bibitem[Par17b]{logmotives}
D.~Park.
\newblock Triangulated categories of motives over fs log schemes.
\newblock arXiv$:$1707.09435, 2017.

\bibitem[Par19]{dpark}
D.~Park.
\newblock Equivariant cd-structures and descent theory.
\newblock {\em J. Pure Appl. Algebra}, 2019.
\newblock to appear.

\bibitem[PPR09]{PPR}
I.~Panin, K.~Pimenov, and O.~R\"ondigs.
\newblock On {V}oevodsky's algebraic {K}-theory spectrum.
\newblock In {\em Algebraic Topology}, volume~4 of {\em Abel. Symp.}, pages
  279--330. Springer, 2009.

\bibitem[PS18]{positive1}
D.~Pavlov and J.~Scholbach.
\newblock Homotopy theory of symmetric powers.
\newblock {\em Homology Homotopy Appl.}, 20(1):359--397, 2018.

\bibitem[Qui73]{Quillen}
D.~Quillen.
\newblock Higher algebraic {K}-theory.
\newblock In {\em Higher {K}-theories {I} (Proc. Conf., Battelle Memorial
  Inst., Seattle, Wash., 1972)}, volume 341 of {\em Lecture Notes in
  Mathematics}, pages 85--147. Springer-Verlag, 1973.

\bibitem[Rio10]{RiouBGL}
J.~Riou.
\newblock Algebraic {K}-theory, {$\mathbf A\sp 1$}-homotopy and
  {R}iemann-{R}och theorems.
\newblock {\em J. Topol.}, 3(2):229--264, 2010.

\bibitem[R{\O}08a]{RNOST}
O.~R{\"o}ndigs and P.~A. {\O}stv{\ae}r.
\newblock Modules over motivic cohomology.
\newblock {\em Adv. Math.}, 219(2):689--727, 2008.

\bibitem[R{\O}08b]{RO}
O.~R\"{o}ndigs and P.~A. {\O}stv{\ae}r.
\newblock Rigidity in motivic homotopy theory.
\newblock {\em Math. Ann.}, 341(3):651--675, 2008.

\bibitem[Rob15]{Robalo}
M.~Robalo.
\newblock {$K$}-theory and the bridge from motives to noncommutative motives.
\newblock {\em Adv. Math.}, 269:399--550, 2015.

\bibitem[R{\"o}n05]{roendigs}
O.~R{\"o}ndigs.
\newblock Functoriality in motivic homotopy theory.
\newblock preprint available at
  \url{https://www.math.uni-bielefeld.de/~oroendig/functoriality.dvi}, 2005.

\bibitem[RS{\O}10]{ROS}
O.~R{\"o}ndigs, M.~Spitzweck, and P.~A. {\O}stv{\ae}r.
\newblock Motivic strict ring models for {K}-theory.
\newblock {\em Proc. Amer. Math. Soc.}, 138(10):3509--3520, 2010.

\bibitem[Ryd10]{Rydh}
D.~Rydh.
\newblock Submersions and effective descent of {\'e}tale morphisms.
\newblock {\em Bull. Soc. Math. France}, 138(2):181--230, 2010.

\bibitem[Sch15]{scholbach2}
J.~Scholbach.
\newblock Arakelov motivic cohomology {II}.
\newblock {\em J. Algebraic Geom.}, 24(4):755--786, 2015.

\bibitem[Ser75]{Ser}
J.-P. Serre.
\newblock {\em Alg\`ebre locale , multiplicit\'es}, volume~11 of {\em Lecture
  Notes in Mathematics}.
\newblock Springer-Verlag, 1975.
\newblock third edition.

\bibitem[Shi04]{shipley}
B.~Shipley.
\newblock A convenient model category for commutative ring spectra.
\newblock {\em Contemp. Math.}, 346:473--483, 2004.

\bibitem[Spi99]{spivak}
M.~Spivakovsky.
\newblock A new proof of {D}. {P}opescu's theorem on smoothing of ring
  homomorphisms.
\newblock {\em Journal of the Amer. Math. Soc.}, 12(2):381--444, 1999.

\bibitem[Spi01]{spit}
M.~Spitzweck.
\newblock Operads, {A}lgebras and {M}odules in {M}onoidal {M}odel {C}ategories
  and {M}otives.
\newblock PhD thesis, Mathematisch Naturwissenschatlichen Fakult{\"a}t,
  Rheinischen Friedrich Wilhelms Universit{\"a}t, Bonn, 2001.

\bibitem[Spi18]{motivicHZ}
M.~Spitzweck.
\newblock A commutative {$\mathbf{P}^1$}-spectrum representing motivic
  cohomology over {D}edekind domains.
\newblock {\em M\'{e}m. Soc. Math. Fr. (N.S.)}, (157), 2018.

\bibitem[SS00]{SS}
S.~Schwede and B.~Shipley.
\newblock Algebras and modules in monoidal model categories.
\newblock {\em Proc. London Math. Soc.}, 80:491--511, 2000.

\bibitem[SV96]{SV0}
A.~Suslin and V.~Voevodsky.
\newblock Singular homology of abstract algebraic varieties.
\newblock {\em Invent. Math.}, 123(1):61--94, 1996.

\bibitem[SV00a]{SV2}
A.~Suslin and V.~Voevodsky.
\newblock {\em {B}loch-{K}ato conjecture and motivic cohomology with finite
  coefficients}, volume 548 of {\em NATO Sciences Series, Series C:
  Mathematical and Physical Sciences}, pages 117--189.
\newblock Kluwer, 2000.

\bibitem[SV00b]{SV1}
A.~Suslin and V.~Voevodsky.
\newblock {\em Relative cycles and {C}how sheaves}, volume 143 of {\em Annals
  of Mathematics Studies}, chapter~2, pages 10--86.
\newblock Princeton University Press, 2000.

\bibitem[SVW18]{Soergel2}
W.~Soergel, R.~Virk, and M.~Wendt.
\newblock Equivariant motives and geometric representation theory.
\newblock with an appendix by F. H{\"o}rmann and M. Wendt, arXiv$:$1809.05480,
  2018.

\bibitem[SW18]{Soergel}
W.~Soergel and M.~Wendt.
\newblock Perverse motives and graded derived category {${\mathcal{O}}$}.
\newblock {\em J. Inst. Math. Jussieu}, 17(2):347--395, 2018.

\bibitem[Swa80]{Swan}
R.~G. Swan.
\newblock On seminormality.
\newblock {\em J. Algebra}, 67(1):210--229, 1980.

\bibitem[Tot16]{totaro1}
B.~Totaro.
\newblock The motive of a classifying space.
\newblock {\em Geom. Topol.}, 20(4):2079--2133, 2016.

\bibitem[Tot18]{totaro2}
B.~Totaro.
\newblock Adjoint functors on the derived category of motives.
\newblock {\em J. Inst. Math. Jussieu}, 17(3):489--507, 2018.

\bibitem[TR19]{shtukas}
J.~Scholbach T.~Richarz.
\newblock The intersection motive of the moduli stack of shtukas.
\newblock arXiv$:$1901.04919, 2019.

\bibitem[TT90]{TT}
R.~W. Thomason and Thomas Trobaugh.
\newblock Higher algebraic {$K$}-theory of schemes and of derived categories.
\newblock In {\em The {G}rothendieck {F}estschrift, {V}ol.\ {III}}, volume~88
  of {\em Progr. Math.}, pages 247--435. Birkh\"auser Boston, Boston, MA, 1990.

\bibitem[TV19]{ToVe}
B.~To\"{e}n and G.~Vezzosi.
\newblock Trace and {K}{\"u}nneth formulas for singularity categories and
  applications.
\newblock arXiv$:$1710.05902v2, 2019.

\bibitem[Vai16]{vaish2}
V.~Vaish.
\newblock Motivic weightless complex and the relative {A}rtin motive.
\newblock {\em Adv. Math.}, 292:316--373, 2016.

\bibitem[Vai17]{vaish4}
V.~Vaish.
\newblock Punctual gluing of $t$-structures and weight structures.
\newblock arXiv$:$1705.00790, 2017.

\bibitem[Vai18]{vaish3}
V.~Vaish.
\newblock On weight truncations in the motivic setting.
\newblock {\em Adv. Math.}, 338:339--400, 2018.

\bibitem[Vez01]{Vezzosi}
G.~Vezzosi.
\newblock Brown-{P}eterson spectra in stable {$\mathbf A\sp 1$}-homotopy
  theory.
\newblock {\em Rend. Sem. Mat. Univ. Padova}, 106:47--64, 2001.

\bibitem[Voe92]{VBletter}
V.~Voevodsky.
\newblock A letter to beilinson.
\newblock K-theory preprint archive 0033, 1992.

\bibitem[Voe96]{V1}
V.~Voevodsky.
\newblock Homology of schemes.
\newblock {\em Selecta Math. (N.S.)}, 2(1):111--153, 1996.

\bibitem[Voe98]{V3}
V.~Voevodsky.
\newblock {$\mathbf A^1$}-homotopy theory.
\newblock In {\em Proceedings of the International Congress of Mathematicians,
  Vol. I (Berlin, 1998)}, pages 579--604 (electronic), 1998.

\bibitem[Voe02a]{allagree}
V.~Voevodsky.
\newblock Motivic cohomology groups are isomorphic to higher chow groups in any
  characteristic.
\newblock {\em Int. Math. Res. Not.}, 7:351--355, 2002.

\bibitem[Voe02b]{V_OpenPB}
V.~Voevodsky.
\newblock Open problems in the motivic stable homotopy theory. {I}.
\newblock In {\em Motives, polylogarithms and {H}odge theory, {P}art {I}
  ({I}rvine, {CA}, 1998)}, volume~3 of {\em Int. Press Lect. Ser.}, pages
  3--34. Int. Press, Somerville, MA, 2002.

\bibitem[Voe10a]{cancel}
V.~Voevodsky.
\newblock Cancellation theorem.
\newblock {\em Doc.Math. Extra volume: Andrei A. Suslin sixtieth birthday},
  pages 671--685, 2010.

\bibitem[Voe10b]{voecd1}
V.~Voevodsky.
\newblock Homotopy theory of simplicial sheaves in completely decomposable
  topologies.
\newblock {\em J. Pure Appl. Algebra}, 214(8):1384--1398, 2010.

\bibitem[Voe10c]{voecd2}
V.~Voevodsky.
\newblock Unstable motivic homotopy categories in {N}isnevich and
  cdh-topologies.
\newblock {\em J. Pure Appl. Algebra}, 214(8):1399--1406, 2010.

\bibitem[VSF00]{FSV}
V.~Voevodsky, A.~Suslin, and E.~M. Friedlander.
\newblock {\em Cycles, Transfers and Motivic homology theories}, volume 143 of
  {\em Annals of Mathematics Studies}.
\newblock Princeton Univ. Press, 2000.

\bibitem[Wil17]{Wildes1}
J.~Wildeshaus.
\newblock Intermediate extension of {C}how motives of {A}belian type.
\newblock {\em Adv. Math.}, 305:515--600, 2017.

\bibitem[Wil18]{Wildes2}
J.~Wildeshaus.
\newblock Weights and conservativity.
\newblock {\em Algebr. Geom.}, 5(6):686--702, 2018.

\end{thebibliography}
\pagebreak
\renewcommand{\Rc}{\Lambda}
\markboth{}{Index}
\printindex{word}{Index}
\markboth{}{Notation}
\printindex{notat}{Notation}
\markboth{}{}
\section*{Index of properties of $\Pmor$-fibred triangulated categories}
\addcontentsline{toc}{section}{Index of properties of $\Pmor$-fibred triangulated categories}
\label{appendix:Pfibred}
\begin{center}

\begin{tabular}{|l|c|c|c|l|}
\hline
Name & Symbol & Def. & related & Remark \\
& & & result & \\
\hline
\hline
additive & & \ref{df:additive_triangulated} && \\
\hline
adjoint property & \adj & \ref{df:adj_BC_PF_prop_tri_premotivic}
 & \ref{thm:support} & \\
\hline
adjoint property for $f$ &
 \adjx f & \ref{df:adj_BC_PF_prop_tri_premotivic}
  && $f$ morphism of schemes\\
\hline
cotransversality property & & \ref{df:P-fibred_transversality}
 && defined for any $\Pmor$-fibred category \\
\hline
homotopy property & \htp & \ref{df:ppty:htp} && \\
\hline
localization property & \loc & \ref{df:localization}
 & \ref{thm:localizationandwpur} & \\
& & & \ref{smoothloc} & \\
\hline
localization property for $i$ & \locx i
 & \textsection \ref{num:general_notations_loc} & & 
 $i$ closed immersion \\
\hline
motivic & & \ref{df:motivic_cat} & \ref{thm:cor3_Ayoub} & 
for premotivic triangulated categories, \\
& & & \ref{cor:DMBmotivic} & means: \htp, \stab, \loc, \adj \\
\hline
oriented & & \ref{df:premotivic_orientation}
 & \ref{cor:orientationsandpurity}
  & for triangulated premotivic categories \\
 & & && satisfying \wloc \\
\hline
projection formula & \PF & \ref{df:adj_BC_PF_prop_tri_premotivic} && \\
\hline
projection formula for $f$ &
 \PFx f & \ref{df:adj_BC_PF_prop_tri_premotivic}
  & \ref{thm:localizationandwpur} & $f$ morphism of schemes \\
\hline
proper base change property
 & \BC & \ref{df:adj_BC_PF_prop_tri_premotivic}
  & \ref{thm:localizationandwpur} & \\
\hline
proper base change property for $f$ &
 \BCx f & \ref{df:adj_BC_PF_prop_tri_premotivic}
  && $f$ morphism of schemes \\
\hline
purity property & \pur & \ref{df:pure_proper_morphism}
 & \ref{thm:localizationandwpur} & \\
\hline
separated & \sep & \ref{df:ppty:sep_ssep} & \ref{thmfinitness} &  \\
 & & & \ref{thm:localduality} & \\
 &&& \ref{thm:DMBseparated} & \\
\hline
semi-separated & \ssep & \ref{df:ppty:sep_ssep}
 & \ref{prop:charseparated1} &  \\
\hline
stability property & \stab & \ref{df:stabilityppty} && \\
\hline
support property & \supp & \ref{defsupportppty}
 & \ref{prop:support_exists_f_!_plus} &  \\
 &&& \ref{thm:support} & \\
 &&& \ref{prop:DMV_supp_wloc} & \\
\hline
$\tau$-compatible & & \ref{df:weaklytaupure} & \ref{grothendieck6op}
 & $\tau$ set of twists \\
\hline
$\tau$-continuous & & \ref{df:continuous} & \ref{DMtcontinuous} &
 for homotopy $\Pmor$-fibred categories, \\
 &&& \ref{thm:DM_continuity} & $\tau$ set of twists \\
 &&& \ref{DMBcontinuous} & \\
\hline
$\tau$-dualizable & & \ref{deftaupure} & \ref{thm:localduality}
 & $\tau$ set of twists\\
\hline
$t$-descent property & & \ref{defdescente}
 && for homotopy $\Pmor$-fibred categories, \\
&&&& $t$ topology \\
\hline
transversality property & & \ref{df:P-fibred_transversality}
 && for any $\Pmor$-fibred category \\
\hline
$t$-separated & \sepx{$t$} & \ref{df:ppty:sepx t} && $t$ topology \\
\hline
weak localization property & \wloc
 & \ref{df:weak_localization_ppty} & \ref{prop:DMV_supp_wloc} & \\
\hline
weak purity property & \wpur & \ref{df:pure_proper_morphism}
 & \ref{thm:localizationandwpur} & \\
 &&& \ref{cor:orientationsandpurity} & \\
\hline
\end{tabular}
\end{center}
\end{document}